\documentclass[numbers=enddot,12pt,final,onecolumn,notitlepage]{scrartcl}%
\usepackage[headsepline,footsepline,manualmark]{scrlayer-scrpage}
\usepackage[all,cmtip]{xy}
\usepackage{amsfonts}
\usepackage{amssymb}
\usepackage{framed}
\usepackage{amsmath}
\usepackage{comment}
\usepackage{color}
\usepackage[breaklinks=true]{hyperref}
\usepackage[sc]{mathpazo}
\usepackage[T1]{fontenc}
\usepackage{amsthm}
\usepackage{needspace}
\usepackage{allrunes}
\providecommand{\U}[1]{\protect\rule{.1in}{.1in}}
\theoremstyle{definition}
\newtheorem{theo}{Theorem}[section]
\newenvironment{theorem}[1][]
{\begin{theo}[#1]\begin{leftbar}}
{\end{leftbar}\end{theo}}
\newtheorem{lem}[theo]{Lemma}
\newenvironment{lemma}[1][]
{\begin{lem}[#1]\begin{leftbar}}
{\end{leftbar}\end{lem}}
\newtheorem{prop}[theo]{Proposition}
\newenvironment{proposition}[1][]
{\begin{prop}[#1]\begin{leftbar}}
{\end{leftbar}\end{prop}}
\newtheorem{defi}[theo]{Definition}
\newenvironment{definition}[1][]
{\begin{defi}[#1]\begin{leftbar}}
{\end{leftbar}\end{defi}}
\newtheorem{remk}[theo]{Remark}
\newenvironment{remark}[1][]
{\begin{remk}[#1]\begin{leftbar}}
{\end{leftbar}\end{remk}}
\newtheorem{coro}[theo]{Corollary}
\newenvironment{corollary}[1][]
{\begin{coro}[#1]\begin{leftbar}}
{\end{leftbar}\end{coro}}
\newtheorem{conv}[theo]{Convention}
\newenvironment{convention}[1][]
{\begin{conv}[#1]\begin{leftbar}}
{\end{leftbar}\end{conv}}
\newtheorem{quest}[theo]{Question}
\newenvironment{question}[1][]
{\begin{quest}[#1]\begin{leftbar}}
{\end{leftbar}\end{quest}}
\newtheorem{exmp}[theo]{Example}
\newenvironment{example}[1][]
{\begin{exmp}[#1]\begin{leftbar}}
{\end{leftbar}\end{exmp}}
\newenvironment{statement}{\begin{quote}}{\end{quote}}
\newenvironment{verlong}{}{}
\newenvironment{vershort}{}{}
\newenvironment{noncompile}{}{}
\excludecomment{verlong}
\includecomment{vershort}
\excludecomment{noncompile}

\newcommand{\QSym}{\operatorname{QSym}}

\newcommand{\NN}{\mathbb{N}}
\newcommand{\ZZ}{\mathbb{Z}}
\newcommand{\xx}{\mathbf{x}}
\newcommand{\tvi}{\left. \textarm{\tvimadur} \right.}
\newcommand{\bel}{\left. \textarm{\belgthor} \right.}
\newcommand{\are}{\ar@{-}}
\newcommand{\arinj}{\ar@{_{(}->}}
\newcommand{\arsurj}{\ar@{->>}}
\newcommand\arxiv[1]{\href{http://www.arxiv.org/abs/#1}{\texttt{arXiv:#1}}}
\let\sumnonlimits\sum
\let\prodnonlimits\prod
\renewcommand{\sum}{\sumnonlimits\limits}
\renewcommand{\prod}{\prodnonlimits\limits}
\begin{document}

\title{Shuffle-compatible permutation statistics II: the exterior peak set}
\author{Darij Grinberg}
\date{31 August 2026}
\maketitle
\tableofcontents

\section*{***}

This paper is a continuation of the work \cite{part1} by Gessel and Zhuang
(but can be read independently of the latter). It is devoted to the study of
shuffle-compatibility of permutation statistics -- a concept introduced in
\cite{part1}, although various instances of it had appeared earlier across
the literature.

In Section \ref{sect.notations}, we introduce the notations that we will need
throughout this paper. In Section \ref{sect.Zenri}, we prove that the exterior
peak set statistic $\operatorname{Epk}$ is shuffle-compatible (Theorem
\ref{thm.Epk.sh-co-a}), as conjectured by Gessel and Zhuang in \cite{part1}.
In Section \ref{sect.LR}, we introduce the concept of an \textquotedblleft
LR-shuffle-compatible\textquotedblright\ statistic, which is stronger than
shuffle-compatibility. We give a sufficient criterion for it and use it to
show that $\operatorname{Epk}$ and some other statistics are
LR-shuffle-compatible.

\begin{vershort}
The last three sections relate all of this to quasisymmetric functions; these
sections are only brief summaries, and we refer to \cite{verlong} for the
details. In Section \ref{sect.Descent}, we recall the concept of descent
statistics introduced in \cite{part1} and its connection to quasisymmetric
functions. Motivated by this connection, in Section \ref{sect.kernel}, we
define the kernel of a descent statistic, and study this kernel for
$\operatorname{Epk}$, giving two explicit generating sets for this kernel. In
Section \ref{sect.dendri}, we extend the quasisymmetric functions connection
to the concept of LR-shuffle-compatible statistics, and relate it to
dendriform algebras.
\end{vershort}

\begin{verlong}
The last three sections relate all of this to quasisymmetric functions: In
Section \ref{sect.Descent}, we recall the concept of descent statistics
introduced in \cite{part1} and its connection to quasisymmetric functions.
Motivated by this connection, in Section \ref{sect.kernel}, we define the
kernel of a descent statistic, and study this kernel for
$\operatorname{Epk}$, giving two explicit generating sets for this kernel.
In Section \ref{sect.dendri},
we extend the quasisymmetric functions connection to the
concept of LR-shuffle-compatible statistics, and relate it to dendriform algebras.
\end{verlong}

\begin{noncompile}
Much of the following is a rough draft, and some proofs (particularly in
Section \ref{sect.dendri}) are merely outlined.
\end{noncompile}

\subsection*{Acknowledgments}

We thank Yan Zhuang, Ira Gessel and Sara Billey for helpful conversations and
corrections. An anonymous referee has also helpfully pointed out mistakes.
The SageMath computer algebra system \cite{SageMath} has been
used in finding some of the results below.

\subsection{Remark on alternative versions}

\begin{vershort}
This paper also has a detailed version \cite{verlong}, which includes some
proofs that have been omitted from the present version as well as more details
on some other proofs and further results in Sections \ref{sect.Descent},
\ref{sect.kernel} and \ref{sect.dendri}.
\end{vershort}

\begin{verlong}
You are reading the detailed version of this paper. For the standard version
(which is shorter by virtue of omitting some proofs and even some results),
see \cite{vershort}.
\end{verlong}

\section{\label{sect.notations}Notations and definitions}

Let us first introduce the definitions and notations that we will use in the
rest of this paper. Many of these definitions appear in \cite{part1} already;
we have tried to deviate from the notations of \cite{part1} as little as possible.

\subsection{Permutations and other basic concepts}

\begin{definition}
We let $\mathbb{N}=\left\{  0,1,2,3,\ldots\right\}  $ and $\mathbb{P}=\left\{
1,2,3,\ldots\right\}  $. Both of these sets are understood to be equipped with
their standard total order.
Elements of $\mathbb{P}$ will be called \textit{letters} (despite being
numbers).
\end{definition}

\begin{definition}
Let $n\in\mathbb{Z}$. We shall use the notation $\left[  n\right]  $ for the
totally ordered set $\left\{  1,2,\ldots,n\right\}  $ (with the usual order
relation inherited from $\mathbb{Z}$). Note that $\left[  n\right]
=\varnothing$ when $n\leq0$.
\end{definition}

\begin{definition}
\label{def.perm}Let $n\in\mathbb{N}$. An $n$\textit{-permutation} shall mean a
word with $n$ letters, which are distinct and belong to $\mathbb{P}$.
Equivalently, an $n$-permutation shall be regarded as an injective map
$\left[  n\right]  \rightarrow\mathbb{P}$ (the image of $i$ under this map
being the $i$-th letter of the word).
\end{definition}

For example, $\left(  3,6,4\right)  $ and $\left(  9,1,2\right)  $ are
$3$-permutations, but $\left(  2,1,2\right)  $ is not.

\begin{definition}
A \textit{permutation} is defined to be an $n$-permutation for some
$n\in\mathbb{N}$. If $\pi$ is an $n$-permutation for some $n\in\mathbb{N}$,
then the number $n$ is called the \textit{size} of the
permutation $\pi$ and is denoted by $\left\vert \pi\right\vert $. A
permutation is said to be \textit{nonempty} if it is nonempty as a word (i.e.,
if its size is $>0$).
\end{definition}

Note that the meaning of \textquotedblleft permutation\textquotedblright\ we
have just defined is unusual (most authors define a permutation to be a
bijection from a set to itself); we are following \cite{part1} in defining
permutations this way.

\begin{definition}
Let $n\in\mathbb{N}$. Two $n$-permutations $\alpha$ and $\beta$ are said to be
\textit{order-isomorphic} if they have the following property: For every two
integers $i,j\in\left[  n\right]  $, we have $\alpha\left(  i\right)
<\alpha\left(  j\right)  $ if and only if $\beta\left(  i\right)
<\beta\left(  j\right)  $.
\end{definition}

\begin{definition}
\textbf{(a)} A \textit{permutation statistic} is a map $\operatorname{st}$
from the set of all permutations to an arbitrary set that has the following
property: Whenever $\alpha$ and $\beta$ are two order-isomorphic permutations,
we have $\operatorname{st} \alpha = \operatorname{st} \beta$.

\textbf{(b)} Let $\operatorname{st}$ be a permutation statistic. Two
permutations $\alpha$ and $\beta$ are said to be
$\operatorname{st}$\textit{-equivalent} if they satisfy
$\left\vert \alpha\right\vert = \left\vert \beta\right\vert $ and
$\operatorname{st}\alpha = \operatorname{st}\beta$.
The relation
\textquotedblleft$\operatorname{st}$-equivalent\textquotedblright\ is
an equivalence relation; its equivalence
classes are called $\operatorname{st}$\textit{-equivalence classes}.
\end{definition}

\begin{remark}
Let $n\in\mathbb{N}$. Let us call an $n$-permutation $\pi$ \textit{standard}
if its letters are $1,2,\ldots,n$ (in some order). The standard $n$%
-permutations are in bijection with the $n!$ permutations of the set $\left\{
1,2,\ldots,n\right\}  $ in the usual sense of this word (i.e., the bijections
from this set to itself).

It is easy to see that for each $n$-permutation $\sigma$, there exists a
\textbf{unique} standard $n$-permutation $\pi$ order-isomorphic to $\sigma$.
Thus, a permutation statistic is uniquely determined by its values on standard
permutations. Consequently, we can view permutation statistics as statistics
defined on standard permutations, i.e., on permutations in the usual sense of
the word.
\end{remark}

The word \textquotedblleft permutation statistic\textquotedblright\ is often
abbreviated as \textquotedblleft statistic\textquotedblright.

\subsection{Some examples of permutation statistics}

\begin{definition}
\label{def.Des-et-al}Let $n\in\mathbb{N}$. Let $\pi=\left(  \pi_{1},\pi
_{2},\ldots,\pi_{n}\right)  $ be an $n$-permutation.

\textbf{(a)} The \textit{descents} of $\pi$ are the elements $i\in\left[
n-1\right]  $ satisfying $\pi_{i}>\pi_{i+1}$.

\textbf{(b)} The \textit{descent set} of $\pi$ is defined to be the set of all
descents of $\pi$. This set is denoted by $\operatorname{Des}\pi$, and is
always a subset of $\left[  n-1\right]  $.

\textbf{(c)} The \textit{peaks} of $\pi$ are the elements $i\in\left\{
2,3,\ldots,n-1\right\}  $ satisfying $\pi_{i-1}<\pi_{i}>\pi_{i+1}$.

\textbf{(d)} The \textit{peak set} of $\pi$ is defined to be the set of all
peaks of $\pi$. This set is denoted by $\operatorname*{Pk}\pi$, and is always
a subset of $\left\{  2,3,\ldots,n-1\right\}  $.

\textbf{(e)} The \textit{left peaks} of $\pi$ are the elements $i\in\left[
n-1\right]  $ satisfying $\pi_{i-1}<\pi_{i}>\pi_{i+1}$, where we set $\pi
_{0}=0$.

\textbf{(f)} The \textit{left peak set} of $\pi$ is defined to be the set of
all left peaks of $\pi$. This set is denoted by $\operatorname*{Lpk}\pi$, and
is always a subset of $\left[  n-1\right]  $. It is easy to see that (for
$n\geq2$) we have%
\[
\operatorname*{Lpk}\pi=\operatorname*{Pk}\pi\cup\left\{  1\ \mid\ \pi_{1}%
>\pi_{2}\right\}  .
\]
(The strange notation ``$\left\{  1\ \mid\ \pi_{1}>\pi_{2}\right\}  $''
means the set of all numbers $1$ satisfying $\pi_1 > \pi_2$.
In other words, it is the $1$-element
set $\left\{  1\right\}  $ if $\pi_{1}>\pi_{2}$, and the empty set
$\varnothing$ otherwise.)

\textbf{(g)} The \textit{right peaks} of $\pi$ are the elements $i\in\left\{
2,3,\ldots,n\right\}  $ satisfying $\pi_{i-1}<\pi_{i}>\pi_{i+1}$, where we set
$\pi_{n+1}=0$.

\textbf{(h)} The \textit{right peak set} of $\pi$ is defined to be the set of
all right peaks of $\pi$. This set is denoted by $\operatorname*{Rpk}\pi$, and
is always a subset of $\left\{  2,3,\ldots,n\right\}  $. It is easy to see
that (for $n\geq2$) we have%
\[
\operatorname*{Rpk}\pi=\operatorname*{Pk}\pi\cup\left\{  n\ \mid\ \pi
_{n-1}<\pi_{n}\right\}  .
\]

\textbf{(i)} The \textit{exterior peaks} of $\pi$ are the elements
$i\in\left[  n\right]  $ satisfying $\pi_{i-1}<\pi_{i}>\pi_{i+1}$, where we
set $\pi_{0}=0$ and $\pi_{n+1}=0$.

\textbf{(j)} The \textit{exterior peak set} of $\pi$ is defined to be the set
of all exterior peaks of $\pi$. This set is denoted by
$\operatorname{Epk} \pi$, and is always a subset of $\left[  n\right]  $.
It is easy to see that (for $n\geq2$) we have
\begin{align*}
\operatorname{Epk}\pi &  =\operatorname*{Pk}\pi\cup\left\{  1\ \mid\ \pi
_{1}>\pi_{2}\right\}  \cup\left\{  n\ \mid\ \pi_{n-1}<\pi_{n}\right\} \\
&  =\operatorname*{Lpk}\pi\cup\operatorname*{Rpk}\pi
\end{align*}
(where, again, $\left\{  n\ \mid\ \pi_{n-1}<\pi_{n}\right\}$ is the
$1$-element set $\left\{n\right\}$ if $\pi_{n-1}<\pi_{n}$,
and otherwise is the empty set).

(For $n=1$, we have $\operatorname{Epk}\pi=\left\{  1\right\}  $.)
\end{definition}

For example, the $6$-permutation $\pi=\left(  4,1,3,9,6,8\right)  $ has%
\begin{align*}
\operatorname{Des}\pi &  =\left\{  1,4\right\}
,\ \ \ \ \ \ \ \ \ \ \operatorname*{Pk}\pi=\left\{  4\right\}  ,\\
\operatorname*{Lpk}\pi &  =\left\{  1,4\right\}
,\ \ \ \ \ \ \ \ \ \ \operatorname*{Rpk}\pi=\left\{  4,6\right\}
,\ \ \ \ \ \ \ \ \ \ \operatorname{Epk}\pi=\left\{  1,4,6\right\}  .
\end{align*}
For another example, the $6$-permutation $\pi=\left(  1,4,3,2,9,8\right)  $
has%
\begin{align*}
\operatorname{Des}\pi &  =\left\{  2,3,5\right\}
,\ \ \ \ \ \ \ \ \ \ \operatorname*{Pk}\pi=\left\{  2,5\right\}  ,\\
\operatorname*{Lpk}\pi &  =\left\{  2,5\right\}
,\ \ \ \ \ \ \ \ \ \ \operatorname*{Rpk}\pi=\left\{  2,5\right\}
,\ \ \ \ \ \ \ \ \ \ \operatorname{Epk}\pi=\left\{  2,5\right\}  .
\end{align*}

Notice that Definition \ref{def.Des-et-al} actually defines several
permutation statistics. For example, Definition \ref{def.Des-et-al}
\textbf{(b)} defines the permutation statistic $\operatorname{Des}$, whose
codomain is the set of all subsets of $\mathbb{P}$. Likewise, Definition
\ref{def.Des-et-al} \textbf{(d)} defines the permutation statistic
$\operatorname*{Pk}$, and Definition \ref{def.Des-et-al} \textbf{(f)} defines
the permutation statistic $\operatorname*{Lpk}$, whereas Definition
\ref{def.Des-et-al} \textbf{(h)} defines the permutation statistic
$\operatorname*{Rpk}$. The main permutation statistic that we will study in
this paper is $\operatorname{Epk}$, which is defined in Definition
\ref{def.Des-et-al} \textbf{(j)}; its codomain is the set of all subsets of
$\mathbb{P}$.

The following simple fact expresses the set $\operatorname{Epk}\pi$
corresponding to an $n$-permutation $\pi$ in terms of $\operatorname{Des}\pi$:

\begin{proposition}
\label{prop.Epk.through-Des}Let $n$ be a positive integer. Let $\pi$ be an
$n$-permutation. Then,%
\[
\operatorname{Epk}\pi=\left(  \operatorname{Des}\pi\cup\left\{  n\right\}
\right)  \setminus\left(  \operatorname{Des}\pi+1\right)  ,
\]
where $\operatorname{Des}\pi+1$ denotes the set $\left\{  i+1\ \mid
\ i\in\operatorname{Des}\pi\right\}  $.
\end{proposition}

\begin{vershort}
\begin{proof}
[Proof of Proposition \ref{prop.Epk.through-Des}.]The rather easy proof can be
found in the detailed version \cite{verlong} of this paper.
\end{proof}
\end{vershort}

\begin{verlong}
\begin{proof}
[Proof of Proposition \ref{prop.Epk.through-Des}.]Write $\pi$ in the form
$\pi=\left(  \pi_{1},\pi_{2},\ldots,\pi_{n}\right)  $. Set $\pi_{0}=0$ and
$\pi_{n+1}=0$. Recall that $\operatorname{Des}\pi$ is defined as the set of
all descents of $\pi$. In other words,%
\[
\operatorname{Des}\pi=\left(  \text{the set of all descents of }\pi\right)
=\left\{  i\in\left[  n-1\right]  \ \mid\ \pi_{i}>\pi_{i+1}\right\}
\]
(because the descents of $\pi$ are defined to be the $i\in\left[  n-1\right]
$ satisfying $\pi_{i}>\pi_{i+1}$).

But $\left(  \pi_{1},\pi_{2},\ldots,\pi_{n}\right)  =\pi$ is an $n$%
-permutation, and thus has no equal entries. Hence, for each $i\in\left[
n-1\right]  $, we have $\pi_{i}\neq\pi_{i+1}$. Thus, for each $i\in\left[
n-1\right]  $, we have the equivalence $\left(  \pi_{i}\geq\pi_{i+1}\right)
\Longleftrightarrow\left(  \pi_{i}>\pi_{i+1}\right)  $. Therefore,%
\[
\left\{  i\in\left[  n-1\right]  \ \mid\ \pi_{i}\geq\pi_{i+1}\right\}
=\left\{  i\in\left[  n-1\right]  \ \mid\ \pi_{i}>\pi_{i+1}\right\}
=\operatorname{Des}\pi.
\]

On the other hand, $\pi_{n}\in\mathbb{P}$, so that $\pi_{n}%
>0=\pi_{n+1}$. Hence, $n$ is an element of the set $\left\{  i\in\left\{
n\right\}  \ \mid\ \pi_{i}>\pi_{i+1}\right\}  $. Clearly, this set cannot have
any other element (since it is a subset of $\left\{  n\right\}  $); thus,
$\left\{  i\in\left\{  n\right\}  \ \mid\ \pi_{i}>\pi_{i+1}\right\}  =\left\{
n\right\}  $.

But $\left[  n\right]  =\left[  n-1\right]  \cup\left\{  n\right\}  $, so
that
\begin{align}
&  \left\{  i\in\left[  n\right]  \ \mid\ \pi_{i}>\pi_{i+1}\right\}
\nonumber\\
&  =\left\{  i\in\left[  n-1\right]  \cup\left\{  n\right\}  \ \mid\ \pi
_{i}>\pi_{i+1}\right\} \nonumber\\
&  =\underbrace{\left\{  i\in\left[  n-1\right]  \ \mid\ \pi_{i}>\pi
_{i+1}\right\}  }_{=\operatorname{Des}\pi}\cup\underbrace{\left\{
i\in\left\{  n\right\}  \ \mid\ \pi_{i}>\pi_{i+1}\right\}  }_{=\left\{
n\right\}  }=\operatorname{Des}\pi\cup\left\{  n\right\}  .
\label{pf.prop.Epk.through-Des.Desun}%
\end{align}

On the other hand, $\pi_{1}\in\mathbb{P}$, so that $\pi_{1}>0=\pi_{0}%
$. Hence, we do not have $\pi_{0}\geq\pi_{1}$. Thus, $0$ is not an element of
the set $\left\{  i\in\left\{  0\right\}  \ \mid\ \pi_{i}\geq\pi
_{i+1}\right\}  $. Clearly, this set cannot have any other element (since it
is a subset of $\left\{  0\right\}  $); thus,
$\left\{  i\in\left\{  0\right\} \ \mid\ \pi_{i}\geq\pi_{i+1}\right\}
=\varnothing$.

But $\left\{  0,1,\ldots,n-1\right\}  =\left\{  0\right\}  \cup\left[
n-1\right]  $, so that
\begin{align*}
&  \left\{  i\in\left\{  0,1,\ldots,n-1\right\}  \ \mid\ \pi_{i}\geq\pi
_{i+1}\right\} \\
&  =\left\{  i\in\left\{  0\right\}  \cup\left[  n-1\right]  \ \mid\ \pi
_{i}\geq\pi_{i+1}\right\} \\
&  =\underbrace{\left\{  i\in\left\{  0\right\}  \ \mid\ \pi_{i}\geq\pi
_{i+1}\right\}  }_{=\varnothing}\cup\underbrace{\left\{  i\in\left[
n-1\right]  \ \mid\ \pi_{i}\geq\pi_{i+1}\right\}  }_{=\operatorname{Des}\pi
}=\varnothing\cup\operatorname{Des}\pi=\operatorname{Des}\pi.
\end{align*}
Hence,%
\[
\operatorname{Des}\pi=\left\{  i\in\left\{  0,1,\ldots,n-1\right\}
\ \mid\ \pi_{i}\geq\pi_{i+1}\right\}  ,
\]
so that%
\begin{align}
\operatorname{Des}\pi+1  &  =\left\{  i+1\ \mid\ i\in\left\{  0,1,\ldots
,n-1\right\}  \ \text{satisfies}\ \pi_{i}\geq\pi_{i+1}\right\} \nonumber\\
&  =\left\{  j\in\left[  n\right]  \ \mid\ \pi_{j-1}\geq\pi_{j}\right\}
\nonumber\\
&  =\left\{  i\in\left[  n\right]  \ \mid\ \pi_{i-1}\geq\pi_{i}\right\}
\label{pf.prop.Epk.through-Des.Des+1}%
\end{align}
(here, we have renamed the index $j$ as $i$).

But $\operatorname{Epk}\pi$ is the set of all exterior peaks of $\pi$ (by the
definition of $\operatorname{Epk}\pi$). Thus,%
\begin{align*}
\operatorname{Epk}\pi &  =\left(  \text{the set of all exterior peaks of }%
\pi\right) \\
&  =\left\{  i\in\left[  n\right]  \ \mid\ \underbrace{\pi_{i-1}<\pi_{i}%
>\pi_{i+1}}_{\Longleftrightarrow\ \left(  \pi_{i}>\pi_{i+1}\text{ and }%
\pi_{i-1}<\pi_{i}\right)  }\right\} \\
&  \ \ \ \ \ \ \ \ \ \ \left(  \text{by the definition of an \textquotedblleft
exterior peak\textquotedblright\ of }\pi\right) \\
&  =\left\{  i\in\left[  n\right]  \ \mid\ \pi_{i}>\pi_{i+1}\text{ and
}\underbrace{\pi_{i-1}<\pi_{i}}_{\Longleftrightarrow\ \left(  \text{not }%
\pi_{i-1}\geq\pi_{i}\right)  }\right\} \\
&  =\left\{  i\in\left[  n\right]  \ \mid\ \pi_{i}>\pi_{i+1}\text{ and not
}\pi_{i-1}\geq\pi_{i}\right\} \\
&  =\underbrace{\left\{  i\in\left[  n\right]  \ \mid\ \pi_{i}>\pi
_{i+1}\right\}  }_{\substack{=\operatorname{Des}\pi\cup\left\{  n\right\}
\\\text{(by (\ref{pf.prop.Epk.through-Des.Desun}))}}}\setminus
\underbrace{\left\{  i\in\left[  n\right]  \ \mid\ \pi_{i-1}\geq\pi
_{i}\right\}  }_{\substack{=\operatorname{Des}\pi+1\\\text{(by
(\ref{pf.prop.Epk.through-Des.Des+1}))}}}\\
&  =\left(  \operatorname{Des}\pi\cup\left\{  n\right\}  \right)
\setminus\left(  \operatorname{Des}\pi+1\right)  .
\end{align*}
This proves Proposition \ref{prop.Epk.through-Des}.
\end{proof}
\end{verlong}

\subsection{Shuffles and shuffle-compatibility}

\begin{definition}
\label{def.shuffles}Let $\pi$ and $\sigma$ be two permutations.

\textbf{(a)} We say that $\pi$ and $\sigma$ are \textit{disjoint} if no letter
appears in both $\pi$ and $\sigma$.

\textbf{(b)} Assume that $\pi$ and $\sigma$ are disjoint. Set $m=\left\vert
\pi\right\vert $ and $n=\left\vert \sigma\right\vert $. Let $\tau$ be an
$\left(  m+n\right)  $-permutation. Then, we say that $\tau$ is a
\textit{shuffle} of $\pi$ and $\sigma$ if both $\pi$ and $\sigma$ are
subsequences of $\tau$.

\textbf{(c)} We let $S\left(  \pi,\sigma\right)  $ be the set of all shuffles
of $\pi$ and $\sigma$.
\end{definition}

For example, the permutations $\left(  3,1\right)  $ and $\left(
6,2,9\right)  $ are disjoint, whereas the permutations $\left(  3,1,2\right)
$ and $\left(  6,2,9\right)  $ are not. The shuffles of the two disjoint
permutations $\left(  3,1\right)  $ and $\left(  2,6\right)  $ are%
\begin{align*}
&  \left(  3,1,2,6\right)  ,\ \ \ \ \ \ \ \ \ \ \left(  3,2,1,6\right)
,\ \ \ \ \ \ \ \ \ \ \left(  3,2,6,1\right)  ,\\
&  \left(  2,3,1,6\right)  ,\ \ \ \ \ \ \ \ \ \ \left(  2,3,6,1\right)
,\ \ \ \ \ \ \ \ \ \ \left(  2,6,3,1\right)  .
\end{align*}

If $\pi$ and $\sigma$ are two disjoint permutations, and if $\tau$ is a shuffle
of $\pi$ and $\sigma$, then each letter of $\tau$ must be either a letter of
$\pi$ or a letter of $\sigma$.
\begin{vershort}
 (This follows easily from the pigeonhole principle.)
\end{vershort}
\begin{verlong}
\ \ \ \ \footnote{\textit{Proof.} Let $\pi$ and $\sigma$ be two disjoint
permutations.
Let $\tau$ be a shuffle of $\pi$ and $\sigma$.
We must prove that each letter of $\tau$ must be either a letter of
$\pi$ or a letter of $\sigma$.

Let $P$ be the set of all letters of $\pi$.
Let $S$ be the set of all letters of $\sigma$.
Let $T$ be the set of all letters of $\tau$.

Let $m = \left|\pi\right|$ and $n = \left|\sigma\right|$.
Hence, $\tau$ is an $\left(m+n\right)$-permutation (since $\tau$
is a shuffle of $\pi$ and $\sigma$).
In other words, $\tau$ is a permutation with $m+n$ letters.

The word $\pi$ is a permutation with $m$ letters (because
$m = \left|\pi\right|$). Thus, this word $\pi$ has exactly $m$
distinct letters. In other words, $\left|P\right| = m$
(since $P$ is the set of all letters of $\pi$).
Similarly, $\left|S\right| = n$ and $\left|T\right| = m+n$.
The sets $P$ and $S$ are disjoint (since the permutations $\pi$
and $\sigma$ are disjoint); thus,
$\left|P \cup S\right| = \underbrace{\left|P\right|}_{=m}
                       + \underbrace{\left|S\right|}_{=n}
= m+n$.

But $\pi$ is a subsequence of $\tau$ (since $\tau$
is a shuffle of $\pi$ and $\sigma$). Thus, $P \subseteq T$.
Similarly, $S \subseteq T$.
Combining $P \subseteq T$ with $S \subseteq T$, we obtain
$P \cup S \subseteq T$. Since
$\left|P \cup S\right| = m+n = \left|T\right|$, we thus conclude
that $P \cup S$ is a subset of $T$ but has the same size as $T$.
By the pigeonhole principle, this entails that $P \cup S = T$.

Now, each letter of $\tau$ must be an element of $T$ (by the definition
of $T$), thus an element of $P \cup S$ (since $T = P \cup S$), and thus
either an element of $P$ or an element of $S$. In other words, each
letter of $\tau$ must be either a letter of $\pi$ or a letter of $\sigma$
(by the definitions of $P$ and $S$). Qed.
}
\end{verlong}

If $\pi$ and $\sigma$ are two disjoint permutations, then $S\left(  \pi
,\sigma\right)  =S\left(  \sigma,\pi\right)  $ is an $\dbinom{m+n}{m}$-element
set, where $m=\left\vert \pi\right\vert $ and $n=\left\vert \sigma\right\vert
$.

Definition \ref{def.shuffles} \textbf{(b)} is used, e.g., in \cite{Greene88}.
From the point of view of combinatorics on words, it is somewhat naive, as it
fails to properly generalize to the case when the words $\pi$ and $\sigma$ are
no longer disjoint\footnote{In this general case, it is best to define a
shuffle of two words $\pi=\left(  \pi_{1},\pi_{2},\ldots,\pi_{m}\right)  $ and
$\sigma=\left(  \sigma_{1},\sigma_{2},\ldots,\sigma_{n}\right)  $ as a word of
the form $\left(  \gamma_{\eta\left(  1\right)  },\gamma_{\eta\left(
2\right)  },\ldots,\gamma_{\eta\left(  m+n\right)  }\right)  $, where $\left(
\gamma_{1},\gamma_{2},\ldots,\gamma_{m+n}\right)  $ is the word $\left(
\pi_{1},\pi_{2},\ldots,\pi_{m},\sigma_{1},\sigma_{2},\ldots,\sigma_{n}\right)
$, and where $\eta$ is some permutation of the set $\left\{  1,2,\ldots
,m+n\right\}  $ (that is, a bijection from this set to itself) satisfying
$\eta^{-1}\left(  1\right)  <\eta^{-1}\left(  2\right)  <\cdots<\eta
^{-1}\left(  m\right)  $ (this causes the letters $\pi_{1},\pi_{2},\ldots
,\pi_{m}$ to appear in the word $\left(  \gamma_{\eta\left(  1\right)
},\gamma_{\eta\left(  2\right)  },\ldots,\gamma_{\eta\left(  m+n\right)
}\right)  $ in this order) and $\eta^{-1}\left(  m+1\right)  <\eta^{-1}\left(
m+2\right)  <\cdots<\eta^{-1}\left(  m+n\right)  $ (this causes the letters
$\sigma_{1},\sigma_{2},\ldots,\sigma_{n}$ to appear in the word $\left(
\gamma_{\eta\left(  1\right)  },\gamma_{\eta\left(  2\right)  },\ldots
,\gamma_{\eta\left(  m+n\right)  }\right)  $ in this order). Furthermore, the
proper generalization of $S\left(  \pi,\sigma\right)  $ to this case would be
a multiset, not a mere set.}. But we will not be considering this general
case, since our results do not seem to straightforwardly extend to it
(although we might have to look more closely); thus, Definition
\ref{def.shuffles} will suffice for us.

\begin{definition}
\textbf{(a)} If $a_{1},a_{2},\ldots,a_{k}$ are finitely many arbitrary
objects, then $\left\{  a_{1},a_{2},\ldots,a_{k}\right\}
_{\operatorname*{multi}}$ denotes the \textbf{multiset} whose elements are
$a_{1},a_{2},\ldots,a_{k}$ (each appearing with the multiplicity with which it
appears in the list $\left(  a_{1},a_{2},\ldots,a_{k}\right)  $).

\textbf{(b)} Let $\left(  a_{i}\right)  _{i\in I}$ be a finite family of
arbitrary objects. Then, $\left\{  a_{i}\ \mid\ i\in I\right\}
_{\operatorname*{multi}}$ denotes the \textbf{multiset} whose elements are the
elements of this family (each appearing with the multiplicity with which it
appears in the family).
\end{definition}

For example, $\left\{  k^{2}\ \mid\ k\in\left\{  -2,-1,0,1,2\right\}
\right\}  _{\operatorname*{multi}}$ is the multiset that contains the element
$4$ twice, the element $1$ twice, and the element $0$ once (and no other
elements). This multiset can also be written in the form $\left\{
4,1,0,1,4\right\}  _{\operatorname*{multi}}$, or in the form $\left\{
0,1,1,4,4\right\}  _{\operatorname*{multi}}$.

\begin{definition}
Let $\operatorname{st}$ be a permutation statistic. We say that
$\operatorname{st}$ is \textit{shuffle-compatible} if and only if it has the
following property: For any two disjoint permutations $\pi$ and $\sigma$, the
multiset%
\[
\left\{  \operatorname{st} \tau \ \mid\ \tau\in S\left(
\pi,\sigma\right)  \right\}  _{\operatorname*{multi}}%
\]
depends only on $\operatorname{st} \pi$, $\operatorname{st} \sigma$,
$\left\vert \pi\right\vert $ and $\left\vert \sigma\right\vert $.
\end{definition}

In other words, a permutation statistic $\operatorname{st}$ is
shuffle-compatible if and only if it has the following property:

\begin{itemize}
\item If $\pi$ and $\sigma$ are two disjoint permutations, and if $\pi
^{\prime}$ and $\sigma^{\prime}$ are two disjoint permutations, and if these
permutations satisfy
\begin{align*}
\operatorname{st} \pi    &  =\operatorname{st}\left(
\pi^{\prime}\right)  ,\ \ \ \ \ \ \ \ \ \ \operatorname{st}
\sigma =\operatorname{st}\left(  \sigma^{\prime}\right)  ,\\
\left\vert \pi\right\vert  &  =\left\vert \pi^{\prime}\right\vert
\ \ \ \ \ \ \ \ \ \ \text{and}\ \ \ \ \ \ \ \ \ \ \left\vert \sigma\right\vert
=\left\vert \sigma^{\prime}\right\vert ,
\end{align*}
then
\[
\left\{  \operatorname{st} \tau   \ \mid\ \tau\in S\left(
\pi,\sigma\right)  \right\}  _{\operatorname*{multi}}=\left\{
\operatorname{st} \tau   \ \mid\ \tau\in S\left(  \pi^{\prime
},\sigma^{\prime}\right)  \right\}  _{\operatorname*{multi}}.
\]

\end{itemize}

The notion of a shuffle-compatible permutation statistic was coined by Gessel
and Zhuang in \cite{part1}, where various statistics were analyzed for their
shuffle-compatibility. In particular, it was shown in \cite{part1} that the
statistics $\operatorname{Des}$, $\operatorname*{Pk}$, $\operatorname*{Lpk}$
and $\operatorname*{Rpk}$ are shuffle-compatible. Our next goal is to prove
the same for the statistic $\operatorname{Epk}$.

\section{\label{sect.Zenri}Extending enriched $P$-partitions and the exterior
peak set}

We are going to define \textit{$\mathcal{Z}$-enriched }$P$\textit{-partitions}%
, which are a straightforward generalization of the notions of
\textquotedblleft$P$-partitions\textquotedblright\ \cite{Stanle72},
\textquotedblleft enriched $P$-partitions\textquotedblright\ \cite[\S 2]{Stembr97}
and \textquotedblleft left enriched $P$-partitions\textquotedblright
\ \cite{Peters05}. We will then consider a new particular case
of this notion, which leads to a proof of the shuffle-compatibility of
$\operatorname{Epk}$ conjectured in \cite{part1}
(Theorem~\ref{thm.Epk.sh-co-a} below).

We remark that Bruce Sagan and Duff Baker-Jarvis are currently working on an
alternative, bijective approach to the shuffle-compatibility of permutation
statistics, which may lead to a different proof of this fact.

\subsection{Lacunar sets}

First, let us briefly study \textit{lacunar sets}, a class of subsets
of $\mathbb{Z}$ that are closely connected to exterior peaks.
We start with the definition:

\begin{definition}
A set $S$ of integers is said to be \textit{lacunar} if each $s\in S$
satisfies $s+1\notin S$.
\end{definition}

In other words, a set of integers is lacunar if and only if it contains no two
consecutive integers. For example, the set $\left\{  2, 5, 7\right\}  $ is
lacunar, while the set $\left\{  2, 5, 6\right\}  $ is not.

Lacunar sets of integers are also called \textit{sparse} sets in some of the
literature (though the latter word has several competing meanings).

\begin{definition}
\label{def.lac.Ln}
Let $n\in\mathbb{N}$. We define a set $\mathbf{L}_{n}$ of subsets of $\left[
n\right]  $ as follows:

\begin{itemize}
\item If $n$ is positive, then $\mathbf{L}_{n}$ shall mean the set of all
nonempty lacunar subsets of $\left[  n\right]  $.

\item If $n=0$, then $\mathbf{L}_{n}$ shall mean the set $\left\{
\varnothing\right\}  $.
\end{itemize}
\end{definition}

For example,%
\begin{align*}
\mathbf{L}_{0}  &  =\left\{  \varnothing\right\}
;\ \ \ \ \ \ \ \ \ \ \mathbf{L}_{1}=\left\{  \left\{  1\right\}  \right\}
;\ \ \ \ \ \ \ \ \ \ \mathbf{L}_{2}
=\left\{ \left\{ 1 \right\} , \left\{  2\right\}  \right\}  ;\\
\mathbf{L}_{3}  &  =\left\{  \left\{  1\right\}  ,\left\{  2\right\}
,\left\{  3\right\}  ,\left\{  1,3\right\}  \right\}  .
\end{align*}

\begin{proposition}
\label{prop.lac.fib}Let $\left(  f_{0},f_{1},f_{2},\ldots\right)  $ be the
Fibonacci sequence (defined by $f_{0}=0$ and $f_{1}=1$ and the recursive
relation $f_{m}=f_{m-1}+f_{m-2}$ for all $m\geq2$). Let $n$ be a positive
integer. Then, $\left\vert \mathbf{L}_{n}\right\vert =f_{n+2}-1$.
\end{proposition}

\begin{proof}
[Proof of Proposition \ref{prop.lac.fib}.]Recall that $\mathbf{L}_{n}$ is the
set of all nonempty lacunar subsets of $\left[  n\right]  $ (since $n$ is
positive). Thus, $\left\vert \mathbf{L}_{n}\right\vert $ is the number of all
lacunar subsets of $\left[  n\right]  $ minus $1$ (since the empty set
$\varnothing$, which is clearly a lacunar subset of $\left[  n\right]  $, is
withheld from the count). But a known fact (see, e.g., \cite[Exercise 1.35
\textbf{a.}]{Stanley-EC1}) says that the number of lacunar subsets of $\left[
n\right]  $ is $f_{n+2}$. Combining the preceding two sentences, we conclude
that $\left\vert \mathbf{L}_{n}\right\vert =f_{n+2}-1$. This proves
Proposition \ref{prop.lac.fib}.
\end{proof}

The following observation is easy:

\begin{proposition}
\label{prop.Epk-lac}Let $n\in\mathbb{N}$. Let $\pi$ be an $n$-permutation.
Then, $\operatorname{Epk}\pi\in\mathbf{L}_{n}$.
\end{proposition}

\begin{proof}
[Proof of Proposition \ref{prop.Epk-lac}.]If $n=0$, then the statement is
obvious (since in this case, we have $\operatorname{Epk}\pi=\varnothing
\in\mathbf{L}_{0}$). Thus, WLOG assume that $n\neq0$. Hence, $n$ is positive.
Hence, $\mathbf{L}_{n}$ is the set of all nonempty lacunar subsets of $\left[
n\right]  $ (by the definition of $\mathbf{L}_{n}$).

The set $\operatorname{Epk}\pi$ is lacunar (since two consecutive integers
cannot both be exterior peaks of $\pi$), and is also nonempty (because
if we choose $i \in \left[n\right]$ such that $\pi\left(i\right)$ is the
largest letter of $\pi$, then $i$ is an exterior peak of $\pi$).
\begin{verlong}

[An alternative reason for the nonemptiness of $\operatorname{Epk}\pi$ is the
following: If $\operatorname{Epk}\pi$ was empty, then $\pi$ would have no
peaks, so that the sequence $\left(  \pi\left(  1\right)  ,\pi\left(
2\right)  ,\ldots,\pi\left(  n\right)  \right)  $ would be strictly decreasing
up to a certain point and then strictly increasing from there on; but then,
either $n$ or $1$ would be an exterior peak of $\pi$, which would contradict
the emptiness of $\operatorname{Epk}\pi$.]

\end{verlong}
Therefore, $\operatorname{Epk}\pi$ is a nonempty lacunar subset of
$\left[  n\right]  $. In other words, $\operatorname{Epk}\pi\in\mathbf{L}%
_{n}$ (since $\mathbf{L}_{n}$ is the set of all nonempty lacunar subsets of
$\left[  n\right]  $). This proves Proposition \ref{prop.Epk-lac}.
\end{proof}

Proposition \ref{prop.Epk-lac} actually has a sort of converse:

\begin{proposition}
\label{prop.when-Epk}Let $n\in\mathbb{N}$. Let $\Lambda$ be a subset of
$\left[  n\right]  $. Then, there exists an $n$-permutation $\pi$ satisfying
$\Lambda=\operatorname{Epk}\pi$ if and only if $\Lambda\in\mathbf{L}_{n}$.
\end{proposition}

\begin{vershort}
\begin{proof}
[Proof of Proposition \ref{prop.when-Epk}.]Omitted; see \cite{verlong} for a proof.
\end{proof}
\end{vershort}

\begin{verlong}
\begin{proof}
[Proof of Proposition \ref{prop.when-Epk}.]$\Longrightarrow:$ We need to prove
that for any $n$-permutation $\pi$, we have $\operatorname{Epk}\pi
\in\mathbf{L}_{n}$. But this follows immediately from Proposition
\ref{prop.Epk-lac}. This proves the $\Longrightarrow$ direction of Proposition
\ref{prop.when-Epk}.

$\Longleftarrow:$ Assume that $\Lambda\in\mathbf{L}_{n}$. We must prove that
there exists an $n$-permutation $\pi$ satisfying
$\Lambda = \operatorname{Epk} \pi$.
Such an $n$-permutation $\pi$ can be constructed as follows:

\begin{itemize}
\item If $n=0$, then we simply set $\pi=\left(  {}\right)  $. Thus, for the
rest of this construction, we WLOG assume that $n\neq0$. Hence, $n$ is
positive. Thus, from $\Lambda\in\mathbf{L}_{n}$, we conclude that $\Lambda$ is
a nonempty lacunar subset of $\left[  n\right]  $.

\item Write the set $\Lambda$ in the form $\Lambda=\left\{  u_{1}<u_{2}%
<\cdots<u_{\ell}\right\}  $ (where $\ell=\left\vert \Lambda\right\vert $).
Thus, $\ell\geq1$ (since $\Lambda$ is nonempty), and we can represent the set
$\left[  n\right]  \setminus\Lambda$ as a union of disjoint intervals as
follows:
\begin{align*}
&  \left[  n\right]  \setminus\Lambda\\
&  =\left[  1,u_{1}-1\right]  \cup\left[  u_{1}+1,u_{2}-1\right]  \cup\left[
u_{2}+1,u_{3}-1\right]  \cup\cdots\cup\left[  u_{\ell-1}+1,u_{\ell}-1\right]
\\
&  \ \ \ \ \ \ \ \ \ \ \cup\left[  u_{\ell}+1,n\right]  .
\end{align*}

\item Let $\pi$ take the values $n,n-1,\ldots,n-\ell+1$ on the elements of
$\Lambda$. (For example, this can be achieved by setting $\pi\left(
u_{i}\right)  =n+1-i$ for each $i\in\left[  \ell\right]  $.)

\item Let $\pi$ take the values $1,2,\ldots,n-\ell$ on the elements of
$\left[  n\right]  \setminus\Lambda$ in such a way that:

\begin{itemize}
\item[(A)] on each of the intervals $\left[  1,u_{1}-1\right]  ,\left[
u_{1}+1,u_{2}-1\right]  ,\left[  u_{2}+1,u_{3}-1\right]  ,\ldots,$%
\newline$\left[  u_{\ell-1}+1,u_{\ell}-1\right]  ,\left[  u_{\ell}+1,n\right]
$, the map $\pi$ is either strictly increasing or strictly decreasing;

\item[(B)] if the interval $\left[  1,u_{1}-1\right]  $ is nonempty, then the
map $\pi$ is strictly increasing on this interval;

\item[(C)] if the interval $\left[  u_{\ell}+1,n\right]  $ is nonempty, then
the map $\pi$ is strictly decreasing on this interval.
\end{itemize}

(This is indeed possible, because if the two intervals $\left[  1,u_{1}%
-1\right]  $ and $\left[  u_{\ell}+1,n\right]  $ are both nonempty, then they
are distinct (since $\ell\geq1$).)
\end{itemize}

Any $n$-permutation $\pi$ constructed in this way will satisfy $\Lambda
=\operatorname{Epk}\pi$. Indeed, it is clear that $\pi$ satisfies%
\[
\pi\left(  u\right)  >\pi\left(  v\right)  \ \ \ \ \ \ \ \ \ \ \text{for all
}u\in\Lambda\text{ and }v\in\left[  n\right]  \setminus\Lambda.
\]
Hence, any element of $\Lambda$ is an exterior peak of $\pi$. Conversely, an
element of $\left[  n\right]  \setminus\Lambda$ cannot be an exterior peak of
$\pi$ (because our construction of $\pi$ guarantees that any $s\in\left[
n\right]  \setminus\Lambda$ satisfies either $\left(  s-1\in\left[  n\right]
\text{ and }\pi\left(  s-1\right)  >\pi\left(  s\right)  \right)  $ or
$\left(  s+1\in\left[  n\right]  \text{ and }\pi\left(  s+1\right)
>\pi\left(  s\right)  \right)  $). Thus, the exterior peaks of $\pi$ are
precisely the elements of $\Lambda$; in other words, we have $\Lambda
=\operatorname{Epk}\pi$. This proves the $\Longleftarrow$ direction of
Proposition \ref{prop.when-Epk}.
\end{proof}
\end{verlong}

Next, let us introduce a total order on the finite subsets of $\mathbb{Z}$:

\begin{definition}
\label{def.order-on-P}
\textbf{(a)} Let $\mathbf{P}$ be the set of all finite subsets of $\mathbb{Z}$.

\textbf{(b)} If $A$ and $B$ are any two sets, then $A\bigtriangleup B$ shall
denote the \textit{symmetric difference} of $A$ and $B$. This is the set
$\left(  A\cup B\right)  \setminus\left(  A\cap B\right)  =\left(  A\setminus
B\right)  \cup\left(  B\setminus A\right)  $. It is well-known that the binary
operation $\bigtriangleup$ on sets is associative.

If $A$ and $B$ are two
distinct sets, then the set $A\bigtriangleup B$ is nonempty. Also, if
$A\in\mathbf{P}$ and $B\in\mathbf{P}$, then $A\bigtriangleup B\in\mathbf{P}$.
Thus, if $A$ and $B$ are two distinct sets in $\mathbf{P}$, then $\min\left(
A\bigtriangleup B\right)  \in\mathbb{Z}$ is well-defined.

\textbf{(c)} We define a binary relation $<$ on $\mathbf{P}$ as follows: For
any $A\in\mathbf{P}$ and $B\in\mathbf{P}$, we let $A<B$ if and only if $A\neq
B$ and $\min\left(  A\bigtriangleup B\right)  \in A$. (This definition makes
sense, because the condition $A\neq B$ ensures that $\min\left(
A\bigtriangleup B\right)  $ is well-defined.)
\end{definition}

Note that this relation $<$ is similar to the relation $<$ in \cite[Lemma
4.3]{ABN-peaks}.

\begin{proposition}
\label{prop.lac.order}The relation $<$ on $\mathbf{P}$ is the smaller relation
of a total order on $\mathbf{P}$.
\end{proposition}

\begin{vershort}

\begin{proof}
[Proof of Proposition \ref{prop.lac.order}.]See \cite{verlong} for this
straightforward argument (or imitate \cite[proof of Lemma 4.3]{ABN-peaks}).
\end{proof}
\end{vershort}

\begin{verlong}

\begin{proof}
[Proof of Proposition \ref{prop.lac.order}.]A quick proof of Proposition
\ref{prop.lac.order} can be obtained by imitating \cite[proof of Lemma
4.3]{ABN-peaks} (with the obvious changes made, such as replacing $\max$ by
$\min$, and reversing the relation).
But let us instead give a completely pedestrian proof:

First, we claim the following fact: If $X$ and $Y$ are two finite nonempty
subsets of $\mathbb{Z}$ satisfying $\min X\neq\min Y$, then%
\begin{equation}
X\neq Y\text{ and }\min\left(  X\bigtriangleup Y\right)  \in\left\{  \min
X,\min Y\right\}  . \label{pf.prop.lac.order.XY}%
\end{equation}

[\textit{Proof of (\ref{pf.prop.lac.order.XY}):} Let $X$ and $Y$ be two finite
nonempty subsets of $\mathbb{Z}$ satisfying $\min X\neq\min Y$. From $\min
X\neq\min Y$, we obtain $X\neq Y$. Hence, $\min\left(  X\bigtriangleup
Y\right)  \in\mathbb{Z}$ is well-defined.

The definition of $X\bigtriangleup Y$ yields $X\bigtriangleup Y=\left(
X\setminus Y\right)  \cup\left(  Y\setminus X\right)  \supseteq X\setminus Y$,
so that $X\setminus Y\subseteq X\bigtriangleup Y$.

But our claim is symmetric with respect to $X$ and $Y$ (since $Y\bigtriangleup
X=X\bigtriangleup Y$ and $\left\{  \min Y,\min X\right\}  =\left\{  \min
X,\min Y\right\}  $). Hence, we can WLOG assume that $\min X\leq\min Y$ (since
otherwise, we can just swap $X$ with $Y$ to ensure this). Assume this.
Combining $\min X\leq\min Y$ with $\min X\neq\min Y$, we obtain $\min X<\min
Y$. If we had $\min X\in Y$, then we would have $\min X\geq\min Y$ (since
every $y\in Y$ satisfies $y\geq\min Y$), which would contradict $\min X<\min
Y$. Hence, we cannot have $\min X\in Y$. Thus, $\min X\notin Y$. Combining
this with $\min X\in X$, we obtain $\min X\in X\setminus Y\subseteq
X\bigtriangleup Y$. Hence, $\min X\geq\min\left(  X\bigtriangleup Y\right)  $
(since every $z\in X\bigtriangleup Y$ satisfies $z\geq\min\left(
X\bigtriangleup Y\right)  $).

On the other hand, the definition of $X\bigtriangleup Y$ yields
$X\bigtriangleup Y=\left(  X\cup Y\right)  \setminus\left(  X\cap Y\right)
\subseteq X\cup Y$ and thus%
\[
\min\left(  X\bigtriangleup Y\right)  \geq\min\left(  X\cup Y\right)
=\min\left\{  \min X,\min Y\right\}  =\min X
\]
(since $\min X<\min Y$).
Combining this inequality with $\min X\geq\min\left(  X\bigtriangleup
Y\right)  $, we obtain $\min\left(  X\bigtriangleup Y\right)  =\min
X\in\left\{  \min X,\min Y\right\}  $. Thus, we have shown that $X\neq Y$ and
$\min\left(  X\bigtriangleup Y\right)  \in\left\{  \min X,\min Y\right\}  $.
This proves (\ref{pf.prop.lac.order.XY}).]

The binary relation $<$ is transitive\footnote{\textit{Proof.} Let
$A\in\mathbf{P}$, $B\in\mathbf{P}$ and $C\in\mathbf{P}$ be such that $A<B$ and
$B<C$. We shall prove that $A<C$.
\par
The sets $A$, $B$ and $C$ are elements of $\mathbf{P}$, and thus are finite
subsets of $\mathbb{Z}$ (by the definition of $\mathbf{P}$).
\par
We have $A<B$. In other words, $A\neq B$ and $\min\left(  A\bigtriangleup
B\right)  \in A$ (by the definition of the relation $<$).
\par
We have $B<C$. In other words, $B\neq C$ and $\min\left(  B\bigtriangleup
C\right)  \in B$ (by the definition of the relation $<$).
\par
Let $X=A\bigtriangleup B$. Then, the set $X$ is nonempty (since $A\neq B$) and
is a finite subset of $\mathbb{Z}$ (since $A$ and $B$ are finite subsets of
$\mathbb{Z}$). Thus, $\min X\in\mathbb{Z}$ is well-defined. From
$X=A\bigtriangleup B$, we obtain $\min X=\min\left(  A\bigtriangleup B\right)
\in A$.
\par
Let $Y=B\bigtriangleup C$. Then, the set $Y$ is nonempty (since $B\neq C$) and
is a finite subset of $\mathbb{Z}$ (since $B$ and $C$ are finite subsets of
$\mathbb{Z}$). Thus, $\min Y\in\mathbb{Z}$ is well-defined. From
$Y=B\bigtriangleup C$, we obtain $\min Y=\min\left(  B\bigtriangleup C\right)
\in B$.
\par
From $X=A\bigtriangleup B$ and $Y=B\bigtriangleup C$, we obtain%
\begin{align*}
X\bigtriangleup Y  &  =\left(  A\bigtriangleup B\right)  \bigtriangleup\left(
B\bigtriangleup C\right)  =A\bigtriangleup\underbrace{B\bigtriangleup
B}_{=\varnothing}\bigtriangleup C\ \ \ \ \ \ \ \ \ \ \left(  \text{since the
operation }\bigtriangleup\text{ is associative}\right) \\
&  =\underbrace{A\bigtriangleup\varnothing}_{=A}\bigtriangleup
C=A\bigtriangleup C.
\end{align*}
\par
We have $\min X\in X=A\bigtriangleup B=\left(  A\cup B\right)  \setminus
\left(  A\cap B\right)  $, so that $\min X\notin A\cap B$. If we had $\min
X=\min Y$, then we would have $\min X\in A\cap B$ (since $\min X\in A$ and
$\min X=\min Y\in B$), which would contradict $\min X\notin A\cap B$. Thus, we
cannot have $\min X=\min Y$. Hence, $\min X\neq\min Y$. Thus,
(\ref{pf.prop.lac.order.XY}) yields $X\neq Y$ and $\min\left(  X\bigtriangleup
Y\right)  \in\left\{  \min X,\min Y\right\}  $.
\par
If we had $A=C$, then we would have $X=\underbrace{A}_{=C}\bigtriangleup
B=C\bigtriangleup B=B\bigtriangleup C=Y$, which would contradict $X\neq Y$.
Thus, we cannot have $A=C$. Hence, $A\neq C$. Hence, $\min\left(
A\bigtriangleup C\right)  $ is well-defined. Let $\mu=\min\left(
A\bigtriangleup C\right)  $. Thus, $\mu=\min\left(
\underbrace{A\bigtriangleup C}_{=X\bigtriangleup Y}\right)  =\min\left(
X\bigtriangleup Y\right)  \in\left\{  \min X,\min Y\right\}  $.
\par
Assume (for the sake of contradiction) that $\mu\in C\setminus A$. Thus,
$\mu\notin A$. Hence, $\mu\neq\min X$ (since $\min X\in A$ but $\mu\notin A$).
Combining this with $\mu\in\left\{  \min X,\min Y\right\}  $, we obtain
$\mu\in\left\{  \min X,\min Y\right\}  \setminus\left\{  \min X\right\}
\subseteq\left\{  \min Y\right\}  $. Hence, $\mu=\min Y\in B$. But also
$\mu\in C\setminus A\subseteq C$. Combining $\mu\in B$ with $\mu\in C$, we
obtain $\mu\in B\cap C$. But $\mu=\min Y\in Y=B\bigtriangleup C=\left(  B\cup
C\right)  \setminus\left(  B\cap C\right)  $ (by the definition of
$B\bigtriangleup C$). Thus, $\mu\notin B\cap C$. This contradicts $\mu\in
B\cap C$. This contradiction shows that our assumption (that $\mu\in
C\setminus A$) was wrong. Hence, we don't have $\mu\in C\setminus A$.
\par
We have $\mu=\min\left(  A\bigtriangleup C\right)  \in A\bigtriangleup
C=\left(  A\setminus C\right)  \cup\left(  C\setminus A\right)  $. Hence,
either $\mu\in A\setminus C$ or $\mu\in C\setminus A$. Thus, $\mu\in
A\setminus C$ (since we don't have $\mu\in C\setminus A$). Therefore,
$\min\left(  A\bigtriangleup C\right)  =\mu\in A\setminus C\subseteq A$. So we
have shown that $A\neq C$ and $\min\left(  A\bigtriangleup C\right)  \in A$.
In other words, $A<C$ (by the definition of the relation $<$).
\par
Now, forget that we fixed $A$, $B$ and $C$. We thus have shown that if
$A\in\mathbf{P}$, $B\in\mathbf{P}$ and $C\in\mathbf{P}$ are such that $A<B$
and $B<C$, then $A<C$. In other words, the relation $<$ is transitive.},
irreflexive\footnote{\textit{Proof.} Let $A\in\mathbf{P}$ be such that $A<A$.
We shall derive a contradiction.
\par
We have $A<A$. In other words, $A\neq A$ and $\min\left(  A\bigtriangleup
A\right)  \in A$ (by the definition of the relation $<$). But $A\neq A$ is
absurd. Hence, we have found a contradiction.
\par
Now, forget that we fixed $A$. We thus have found a contradiction for each
$A\in\mathbf{P}$ satisfying $A<A$. Thus, no $A\in\mathbf{P}$ satisfies $A<A$.
In other words, the relation $<$ is irreflexive.} and
asymmetric\footnote{\textit{Proof.} We know that the binary relation $<$ is
transitive and irreflexive. Hence, this relation $<$ is asymmetric (since
every transitive irreflexive relation is asymmetric).}. Hence, $<$ is the
smaller relation of a partial order on $\mathbf{P}$. It remains to prove that
this partial order is a total order. In other words, we need to show that if
$A\in\mathbf{P}$ and $B\in\mathbf{P}$ are distinct, then either $A<B$ or
$B<A$. But this is easy\footnote{\textit{Proof.} Let $A\in\mathbf{P}$ and
$B\in\mathbf{P}$ be distinct. We must show that either $A<B$ or $B<A$.
\par
The sets $A$ and $B$ are elements of $\mathbf{P}$, and thus are finite subsets
of $\mathbb{Z}$ (by the definition of $\mathbf{P}$).
\par
The set $A\bigtriangleup B$ is nonempty (since $A$ and $B$ are distinct) and
is a finite subset of $\mathbb{Z}$ (since $A$ and $B$ are finite subsets of
$\mathbb{Z}$). Hence, $\min\left(  A\bigtriangleup B\right)  $ is
well-defined. Clearly,%
\[
\min\left(  A\bigtriangleup B\right)  \in A\bigtriangleup B=\left(  A\cup
B\right)  \setminus\left(  A\cap B\right)  \subseteq A\cup B.
\]
In other words, either $\min\left(  A\bigtriangleup B\right)  \in A$ or
$\min\left(  A\bigtriangleup B\right)  \in B$. In other words, we are in one
of the following two cases:
\par
\textit{Case 1:} We have $\min\left(  A\bigtriangleup B\right)  \in A$.
\par
\textit{Case 2:} We have $\min\left(  A\bigtriangleup B\right)  \in B$.
\par
Let us first consider Case 1. In this case, we have $\min\left(
A\bigtriangleup B\right)  \in A$. Thus, $A\neq B$ (since $A$ and $B$ are
distinct) and $\min\left(  A\bigtriangleup B\right)  \in A$. In other words,
$A<B$ (by the definition of the relation $<$). Hence, either $A<B$ or $B<A$.
Thus, our claim (that either $A<B$ or $B<A$) is proven in Case 1.
\par
Let us now consider Case 2. In this case, we have $\min\left(
A\bigtriangleup B\right)  \in B$. Hence, $\min\left(
\underbrace{B\bigtriangleup A}_{=A\bigtriangleup B}\right)  =\min\left(
A\bigtriangleup B\right)  \in B$. Thus, $B\neq A$ (since $A$ and $B$ are
distinct) and $\min\left(  B\bigtriangleup A\right)  \in B$. In other words,
$B<A$ (by the definition of the relation $<$). Hence, either $A<B$ or $B<A$.
Thus, our claim (that either $A<B$ or $B<A$) is proven in Case 2.
\par
Now, our claim (that either $A<B$ or $B<A$) has been proven in both Cases 1
and 2. Hence, this claim always holds. In other words, we always have either
$A<B$ or $B<A$. Qed.}. Thus, the proof of Proposition \ref{prop.lac.order} is complete.
\end{proof}
\end{verlong}

In the following, we shall regard the set $\mathbf{P}$ as a totally ordered
set, equipped with the order from Proposition \ref{prop.lac.order}. Thus, for
example, two sets $A$ and $B$ in $\mathbf{P}$ satisfy $A\geq B$ if and only if
either $A=B$ or $B<A$.

\begin{definition}
\label{def.S+1}Let $S$ be a subset of $\mathbb{Z}$. Then, we define a new
subset $S+1$ of $\mathbb{Z}$ by setting
\[
S+1=\left\{  i+1\ \mid\ i\in S\right\}  =\left\{  j\in\mathbb{Z}\ \mid\ j-1\in
S\right\}  .
\]
Note that $S+1\in\mathbf{P}$ if $S\in\mathbf{P}$.
\end{definition}

For example, $\left\{  2,5\right\}  +1=\left\{  3,6\right\}  $. Note that a
subset $S$ of $\mathbb{Z}$ is lacunar if and only if $S\cap\left(  S+1\right)
=\varnothing$.

\begin{proposition}
\label{prop.lac.RL}Let $\Lambda\in\mathbf{P}$ and $R\in\mathbf{P}$ be such
that the set $R$ is lacunar and $R\subseteq\Lambda\cup\left(  \Lambda
+1\right)  $. Then, $R\geq\Lambda$ (with respect to the total order on
$\mathbf{P}$).
\end{proposition}

\begin{proof}
[Proof of Proposition \ref{prop.lac.RL}.]Assume the contrary. Thus,
$R<\Lambda$ (since $\mathbf{P}$ is totally ordered). In other words,
$R\neq\Lambda$ and $\min\left(  R\bigtriangleup\Lambda\right)  \in R$ (by the
definition of the relation $<$). Let $\mu=\min\left(  R\bigtriangleup
\Lambda\right)  $. Thus, $\mu=\min\left(  R\bigtriangleup\Lambda\right)  \in
R\subseteq\Lambda\cup\left(  \Lambda+1\right)  $.

We have $\mu=\min\left(  R\bigtriangleup\Lambda\right)  \in R\bigtriangleup
\Lambda=\left(  R\cup\Lambda\right)  \setminus\left(  R\cap\Lambda\right)  $.
Hence, $\mu\notin R\cap\Lambda$. If we had $\mu\in\Lambda$, then we would have
$\mu\in R\cap\Lambda$ (since $\mu\in R$ and $\mu\in\Lambda$), which would
contradict $\mu\notin R\cap\Lambda$. Thus, we cannot have $\mu\in\Lambda$.
Hence, $\mu\notin\Lambda$. Combining $\mu\in\Lambda\cup\left(  \Lambda
+1\right)  $ with $\mu\notin\Lambda$, we obtain $\mu\in\left(  \Lambda
\cup\left(  \Lambda+1\right)  \right)  \setminus\Lambda\subseteq\Lambda+1$. In
other words, $\mu-1\in\Lambda$.

Every $x\in R\bigtriangleup\Lambda$ satisfies $x\geq\min\left(
R\bigtriangleup\Lambda\right)  $. Hence, if we had $\mu-1\in R\bigtriangleup
\Lambda$, then we would have $\mu-1\geq\min\left(  R\bigtriangleup
\Lambda\right)  =\mu$, which would contradict $\mu-1<\mu$. Thus, we cannot
have $\mu-1\in R\bigtriangleup\Lambda$. Thus, $\mu-1\notin R\bigtriangleup
\Lambda$. Combining this with $\mu-1\in\Lambda$, we obtain $\mu-1\in
\Lambda\setminus\left(  R\bigtriangleup\Lambda\right)  =R\cap\Lambda$ (since
every two sets $X$ and $Y$ satisfy $Y\setminus\left(  X\bigtriangleup
Y\right)  =X\cap Y$). Thus, $\mu-1\in R\cap\Lambda\subseteq R$.

But the set $R$ is lacunar. In other words, each $s\in R$ satisfies $s+1\notin
R$ (by the definition of \textquotedblleft lacunar\textquotedblright).
Applying this to $s=\mu-1$, we obtain $\left(  \mu-1\right)  +1\notin R$
(since $\mu-1\in R$). This contradicts $\left(  \mu-1\right)  +1=\mu\in R$.
This contradiction shows that our assumption was wrong; hence, Proposition
\ref{prop.lac.RL} is proven.
\end{proof}

\begin{verlong}
\begin{corollary}
\label{cor.lac.L=L}Let $\Lambda_{1}$ and $\Lambda_{2}$ be two finite lacunar
subsets of $\mathbb{Z}$ such that $\Lambda_{1}\cup\left(  \Lambda
_{1}+1\right)  =\Lambda_{2}\cup\left(  \Lambda_{2}+1\right)  $. Then,
$\Lambda_{1}=\Lambda_{2}$.
\end{corollary}

% \begin{vershort}
% \begin{proof}
% [Proof of Corollary \ref{cor.lac.L=L}.] Both $\Lambda_{1}$ and $\Lambda_{2}$
% belong to $\mathbf{P}$. Also, $\Lambda_{1}\subseteq\Lambda_{1}\cup\left(
% \Lambda_{1}+1\right)  =\Lambda_{2}\cup\left(  \Lambda_{2}+1\right)  $. Thus,
% Proposition \ref{prop.lac.RL} (applied to $R=\Lambda_{1}$ and $\Lambda
% =\Lambda_{2}$) yields $\Lambda_{1}\geq\Lambda_{2}$. The same argument (with
% the roles of $\Lambda_{1}$ and $\Lambda_{2}$ interchanged) yields $\Lambda
% _{2}\geq\Lambda_{1}$. Combining $\Lambda_{1}\geq\Lambda_{2}$ with $\Lambda
% _{2}\geq\Lambda_{1}$, we obtain $\Lambda_{1}=\Lambda_{2}$. This proves
% Corollary \ref{cor.lac.L=L}.
% \end{proof}
% \end{vershort}

% \begin{verlong}
\begin{proof}
[Proof of Corollary \ref{cor.lac.L=L}.]Both $\Lambda_{1}$ and $\Lambda_{2}$
are finite subsets of $\mathbb{Z}$, and thus belong to $\mathbf{P}$ (by the
definition of $\mathbf{P}$). Also,
\[
\Lambda_{1}\subseteq\Lambda_{1}\cup\left(  \Lambda_{1}+1\right)  =\Lambda
_{2}\cup\left(  \Lambda_{2}+1\right)  .
\]
Hence, Proposition \ref{prop.lac.RL} (applied to $R=\Lambda_{1}$ and
$\Lambda=\Lambda_{2}$) yields $\Lambda_{1}\geq\Lambda_{2}$ (with respect to
the total order on $\mathbf{P}$). The same argument (with the roles of
$\Lambda_{1}$ and $\Lambda_{2}$ interchanged) yields $\Lambda_{2}\geq
\Lambda_{1}$. Combining $\Lambda_{1}\geq\Lambda_{2}$ with $\Lambda_{2}%
\geq\Lambda_{1}$, we obtain $\Lambda_{1}=\Lambda_{2}$. This proves Corollary
\ref{cor.lac.L=L}.
\end{proof}
\end{verlong}

\subsection{\label{subsect.Zenri.gen}$\mathcal{Z}$-enriched $\left(
P,\gamma\right)  $-partitions}

\begin{convention}
By abuse of notation, we will often use the same notation for a poset
$P=\left(  X,\leq\right)  $ and its ground set $X$ when there is no danger of
confusion. In particular, if $x$ is some object, then \textquotedblleft$x\in
P$\textquotedblright\ shall mean \textquotedblleft$x\in X$\textquotedblright.
\end{convention}

\begin{definition}
A \textit{labeled poset} means a pair $\left(  P,\gamma\right)  $ consisting
of a finite poset $P=\left(  X,\leq\right)  $ and an injective map
$\gamma:X\rightarrow A$ for some totally ordered set $A$. The injective map
$\gamma$ is called the \textit{labeling} of the labeled poset $\left(
P,\gamma\right)  $. The poset $P$ is called the \textit{ground poset} of the
labeled poset $\left(  P,\gamma\right)  $.
\end{definition}

\begin{convention}
Let $\mathcal{N}$ be a totally ordered set, whose (strict) order relation will
be denoted by $\prec$. Let $+$ and $-$ be two distinct symbols. Let
$\mathcal{Z}$ be a subset of the set $\mathcal{N}\times\left\{  +,-\right\}
$. For each $q=\left(  n,s\right)  \in\mathcal{Z}$, we denote the element
$n\in\mathcal{N}$ by $\left\vert q\right\vert $, and we call the element
$s\in\left\{  +,-\right\}  $ the \textit{sign} of $q$. If $n\in\mathcal{N}$,
then we will denote the two elements $\left(  n,+\right)  $ and $\left(
n,-\right)  $ of $\mathcal{N}\times\left\{  +,-\right\}  $ by $+n$ and $-n$, respectively.

We equip the set $\mathcal{Z}$ with a total order, whose (strict) order
relation $\prec$ is defined by%
\[
\left(  n,s\right)  \prec\left(  n^{\prime},s^{\prime}\right)  \text{ if and
only if either }n\prec n^{\prime}\text{ or }\left(  n=n^{\prime}\text{ and
}s=-\text{ and }s^{\prime}=+\right)  .
\]
Let $\operatorname*{Pow}\mathcal{N}$ be the ring of all formal power series
over $\mathbb{Q}$ in the indeterminates $x_{n}$ for $n\in\mathcal{N}$.

We fix $\mathcal{N}$ and $\mathcal{Z}$ throughout Subsection
\ref{subsect.Zenri.gen}. That is, any result in this subsection is tacitly
understood to begin with \textquotedblleft Let $\mathcal{N}$ be a totally
ordered set, whose (strict) order relation will be denoted by $\prec$, and let
$\mathcal{Z}$ be a subset of the set $\mathcal{N}\times\left\{  +,-\right\}
$\textquotedblright; and the notations of this convention shall always be
in place throughout this Subsection.

Whenever $\prec$ denotes some strict order, the corresponding weak order will
be denoted by $\preccurlyeq$. (Thus, $a \preccurlyeq b$ means ``$a \prec b$ or
$a = b$''.)
\end{convention}

\begin{definition}
\label{def.ambivPp}Let $\left(  P,\gamma\right)  $ be a labeled poset. A
\textit{$\mathcal{Z}$-enriched }$\left(  P,\gamma\right)  $\textit{-partition}
means a map $f:P\rightarrow\mathcal{Z}$ such that for all $x<y$ in $P$, the
following conditions hold:

\begin{enumerate}
\item[\textbf{(i)}] We have $f\left(  x\right)  \preccurlyeq f\left(
y\right)  $.

\item[\textbf{(ii)}] If $f\left(  x\right)  =f\left(  y\right)  =+n$ for some
$n\in\mathcal{N}$, then $\gamma\left(  x\right)  <\gamma\left(  y\right)  $.

\item[\textbf{(iii)}] If $f\left(  x\right)  =f\left(  y\right)  =-n$ for some
$n\in\mathcal{N}$, then $\gamma\left(  x\right)  >\gamma\left(  y\right)  $.
\end{enumerate}

(Of course, this concept depends on $\mathcal{N}$ and $\mathcal{Z}$, but these
will always be clear from the context.)
\end{definition}

\Needspace{4cm}

\begin{example}
\label{exa.ambivPp.diamond}Let $P$ be the poset with the following Hasse
diagram:%
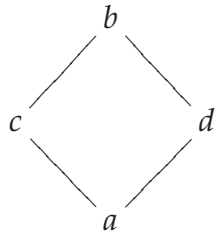
\[%
%TCIMACRO{\TeXButton{Hasse diagram}{\xymatrix{
%& b \are[dl] \are[dr] \\
%c \are[dr] & & d \are[dl] \\
%& a
%}}}%
%BeginExpansion
\xymatrix{
& b \are[dl] \are[dr] \\
c \are[dr] & & d \are[dl] \\
& a
}%
%EndExpansion
\]
(that is, the ground set of $P$ is $\left\{ a, b, c, d \right\}$, and
its order relation is given by $a < c < b$ and $a < d < b$).
Let $\gamma:P\rightarrow\mathbb{Z}$ be a map that satisfies
$\gamma\left(  a\right)  <\gamma\left(  b\right)  <\gamma\left(  c\right)
<\gamma\left(  d\right)  $ (for example, $\gamma$ could be the map that sends
$a,b,c,d$ to $2,3,5,7$, respectively).
Then, $\left(P, \gamma\right)$ is a labeled poset.
A $\mathcal{Z}$-enriched $\left(
P,\gamma\right)  $-partition is a map $f:P\rightarrow\mathcal{Z}$ satisfying
the following conditions:

\begin{enumerate}
\item[\textbf{(i)}] We have $f\left(  a\right)  \preccurlyeq f\left(
c\right)  \preccurlyeq f\left(  b\right)  $ and $f\left(  a\right)
\preccurlyeq f\left(  d\right)  \preccurlyeq f\left(  b\right)  $.

\item[\textbf{(ii)}] We cannot have $f\left(  c\right)  =f\left(  b\right)
=+n$ with $n\in\mathcal{N}$. \newline
We cannot have $f\left(  d\right)
=f\left(  b\right)  =+n$ with $n\in\mathcal{N}$.

\item[\textbf{(iii)}] We cannot have $f\left(  a\right)  =f\left(  c\right)
=-n$ with $n\in\mathcal{N}$. \newline
We cannot have $f\left(  a\right)
=f\left(  d\right)  =-n$ with $n\in\mathcal{N}$.
\end{enumerate}

For example, if $\mathcal{N}=\mathbb{P}$ (the totally ordered set of positive
integers, with its usual ordering) and $\mathcal{Z}=\mathcal{N}\times\left\{
+,-\right\}  $, then the map $f:P\rightarrow\mathcal{Z}$ sending $a,b,c,d$ to
$+2,-3,+2,-3$ (respectively) is a $\mathcal{Z}$-enriched $\left(
P,\gamma\right)  $-partition. Notice that the total ordering on $\mathcal{Z}$
in this case is given by%
\[
-1\prec+1\prec-2\prec+2\prec-3\prec+3\prec\cdots,
\]
rather than by the familiar total order on $\mathbb{Z}$.
\end{example}

The concept of a \textquotedblleft$\mathcal{Z}$-enriched $\left(
P,\gamma\right)  $-partition\textquotedblright\ generalizes three notions in
existing literature: that of a \textquotedblleft$\left(  P,\gamma\right)
$-partition\textquotedblright, that of an \textquotedblleft enriched $\left(
P,\gamma\right)  $-partition\textquotedblright, and that of a
\textquotedblleft left enriched $\left(  P,\gamma\right)  $%
-partition\textquotedblright\footnote{The ideas behind these three concepts
are due to Stanley \cite{Stanle72}, Stembridge \cite[\S 2]{Stembr97} and
Petersen \cite{Peters05}, respectively, but the precise definitions are not
standardized across the literature. We define a \textquotedblleft$\left(
P,\gamma\right)  $-partition\textquotedblright\ as in \cite[\S 1.1]{Stembr97};
this definition differs noticeably from Stanley's (in particular, Stanley
requires $f\left(  x\right)  \succcurlyeq f\left(  y\right)  $ instead of
$f\left(  x\right)  \preccurlyeq f\left(  y\right)  $, but the differences do
not end here). We define an \textquotedblleft enriched $\left(  P,\gamma
\right)  $-partition\textquotedblright\ as in \cite[\S 2]{Stembr97}. Finally,
we define a \textquotedblleft left enriched $\left(  P,\gamma\right)
$-partition\textquotedblright\ to be a $\mathcal{Z}$-enriched $\left(
P,\gamma\right)  $-partition where $\mathcal{N}=\mathbb{N}$ and $\mathcal{Z}%
=\left(  \mathcal{N}\times\left\{  +,-\right\}  \right)  \setminus\left\{
-0\right\}  $; this definition is equivalent to Petersen's
\cite[Definition 3.4.1]{Peters06} up to some differences of notation
(in particular, Petersen assumes that the ground set of $P$ is already
a subset of $\mathbb{P}$, and that the labeling $\gamma$ is the
canonical inclusion map $ P \to \mathbb{P}$; also, he identifies the elements
$+0, -1, +1, -2, +2, \ldots$ of
$\left(  \mathcal{N}\times\left\{  +,-\right\}  \right)  \setminus\left\{
-0\right\}  $ with the integers
$0, -1, +1, -2, +2, \ldots$, respectively).
Note that the definition Petersen gives in \cite[Definition 4.1]{Peters05}
is incorrect, and the one in \cite[Definition 3.4.1]{Peters06} is probably
his intent.}:

\begin{example}
\label{exa.ambivPp.abc}\textbf{(a)} If $\mathcal{N}=\mathbb{P}$ (the totally
ordered set of positive integers) and $\mathcal{Z}=\mathcal{N}\times\left\{
+\right\}  =\left\{  +n\ \mid\ n\in\mathcal{N}\right\}  $, then the
$\mathcal{Z}$-enriched $\left(  P,\gamma\right)  $-partitions are simply the
$\left(  P,\gamma\right)  $-partitions into $\mathcal{N}$, composed with the
canonical bijection $\mathcal{N}\rightarrow\mathcal{Z},\ n\mapsto\left(
+n\right)  $.

\textbf{(b)} If $\mathcal{N}=\mathbb{P}$ (the totally ordered set of positive
integers) and $\mathcal{Z}=\mathcal{N}\times\left\{  +,-\right\}  $, then the
$\mathcal{Z}$-enriched $\left(  P,\gamma\right)  $-partitions are the enriched
$\left(  P,\gamma\right)  $-partitions.

\textbf{(c)} If $\mathcal{N}=\mathbb{N}$ (the totally ordered set of
nonnegative integers) and $\mathcal{Z}=\left(  \mathcal{N}\times\left\{
+,-\right\}  \right)  \setminus\left\{  -0\right\}  $, then the $\mathcal{Z}%
$-enriched $\left(  P,\gamma\right)  $-partitions are the left enriched
$\left(  P,\gamma\right)  $-partitions. Note that $+0$ and $-0$ here stand for
the pairs $\left(  0,+\right)  $ and $\left(  0,-\right)  $; thus, they are
not equal.
\end{example}

\begin{definition}
If $\left(  P,\gamma\right)  $ is a labeled poset, then $\mathcal{E}\left(
P,\gamma\right)  $ shall denote the set of all $\mathcal{Z}$-enriched $\left(
P,\gamma\right)  $-partitions.
\end{definition}

\begin{definition}
Let $P$ be any finite poset. Then, $\mathcal{L}\left(  P\right)  $ shall
denote the set of all linear extensions of $P$. A linear extension of $P$
shall be understood simultaneously as a totally ordered set extending $P$ and
as a list $\left(  w_{1},w_{2},\ldots,w_{n}\right)  $ of all elements of $P$
such that no two integers $i<j$ satisfy $w_{i}\geq w_{j}$ in $P$.
\end{definition}

Let us prove some basic facts about $\mathcal{Z}$-enriched $\left(
P,\gamma\right)  $-partitions, straightforwardly generalizing classical
results proven by Stanley and Gessel (for the case of \textquotedblleft
plain\textquotedblright\ $\left(  P,\gamma\right)  $-partitions), Stembridge
\cite[Lemma 2.1]{Stembr97}
(for enriched $\left(  P,\gamma\right)  $-partitions) and
Petersen \cite[Lemma 3.4.1]{Peters06} (for left
enriched $\left(  P,\gamma\right)  $-partitions):

\begin{proposition}
\label{prop.fund-lem}For any labeled poset $\left(  P,\gamma\right)  $, we
have%
\[
\mathcal{E}\left(  P,\gamma\right)  =\bigsqcup_{w\in\mathcal{L}\left(
P\right)  }\mathcal{E}\left(  w,\gamma\right)  .
\]

\end{proposition}

\begin{vershort}
\begin{proof}
[Proof of Proposition \ref{prop.fund-lem}.]This is analogous to the proof of
\cite[Lemma 2.1]{Stembr97}.
See \cite{verlong} for details.
\end{proof}
\end{vershort}

\begin{verlong}
\begin{proof}
[Proof of Proposition \ref{prop.fund-lem}.] The following proof is a
straightforward generalization of \cite[proof of Lemma 2.1]{Stembr97}
(although we are
rewriting it along the lines of \cite[solution to Exercise 5.2.13]{HopfComb}).

If $z\in\mathcal{Z}$ is any element, then we say that $z$ is \textit{positive}
if the sign of $z$ is $+$ (that is, if $z=\left(  n,+\right)  $ for some
$n\in\mathcal{N}$), and we say that $z$ is \textit{negative} if the sign of
$z$ is $-$ (that is, if $z=\left(  n,-\right)  $ for some $n\in\mathcal{N}$).

Let $\left(P, \gamma\right)$ be a labeled poset.

Let $f\in\mathcal{E}\left(  P,\gamma\right)  $ be arbitrary. Thus, $f$ is a
$\mathcal{Z}$-enriched $\left(  P,\gamma\right)  $-partition (by the
definition of $\mathcal{E}\left(  P,\gamma\right)  $). We now define a binary
relation $<_{f}$ on the set $P$ as follows: For any two elements $x$ and $y$
of $P$, we shall have $x<_{f}y$ if and only if we have%
\begin{align*}
& \text{\textbf{either} }\left(  f\left(  x\right)  \prec f\left(  y\right)
\right)  \\
& \text{\textbf{or} }\left(  f\left(  x\right)  =f\left(  y\right)  \text{ is
positive, and }\gamma\left(  x\right)  <\gamma\left(  y\right)  \right)  \\
& \text{\textbf{or} }\left(  f\left(  x\right)  =f\left(  y\right)  \text{ is
negative, and }\gamma\left(  x\right)  >\gamma\left(  y\right)  \right)
\text{.}%
\end{align*}
It is easy to see that this relation $<_{f}$ is the smaller relation of a
total order on the set $P$
(indeed, it is transitive, irreflexive and asymmetric, and
every two distinct elements $x$ and $y$ satisfy either $x<_{f}y$ or
$y<_{f}x$).
Denote this total order by $\mathbf{w}_{f}$. Thus, $\mathbf{w}_{f}$ is a
total order on the set $P$.

This total order $\mathbf{w}_{f}$ is a linear extension of $P$.

[\textit{Proof:} Let $x$ and $y$ be two elements of $P$ such that $x<y$ in
$P$. We shall show that $x<y$ in $\mathbf{w}_{f}$.

Recall that $f$ is a $\mathcal{Z}$-enriched $\left(  P,\gamma\right)
$-partition. Thus, from $x<y$ in $P$, we obtain the following three facts:

\begin{enumerate}
\item[\textbf{(i)}] We have $f\left(  x\right)  \preccurlyeq f\left(
y\right)  $.

\item[\textbf{(ii)}] If $f\left(  x\right)  =f\left(  y\right)  =+n$ for some
$n\in\mathcal{N}$, then $\gamma\left(  x\right)  <\gamma\left(  y\right)  $.

\item[\textbf{(iii)}] If $f\left(  x\right)  =f\left(  y\right)  =-n$ for some
$n\in\mathcal{N}$, then $\gamma\left(  x\right)  >\gamma\left(  y\right)  $.
\end{enumerate}

From these three facts, we conclude that \textbf{either} $\left(  f\left(
x\right)  \prec f\left(  y\right)  \right)  $ \newline
\textbf{or} $\left(  f\left(
x\right)  =f\left(  y\right)  \text{ is positive, and }\gamma\left(  x\right)
<\gamma\left(  y\right)  \right)  $ \newline
\textbf{or} $\left(  f\left(  x\right)
=f\left(  y\right)  \text{ is negative, and }\gamma\left(  x\right)
>\gamma\left(  y\right)  \right)  $. In other words, $x<_{f}y$ (by the
definition of the relation $<_{f}$). In other words, $x<y$ in $\mathbf{w}_{f}$
(since $<_{f}$ is the smaller relation of $\mathbf{w}_{f}$).

Now, forget that we fixed $x$ and $y$. Thus, we have shown that if $x$ and $y$
are two elements of $P$ such that $x<y$ in $P$, then $x<y$ in $\mathbf{w}_{f}%
$. Hence, $\mathbf{w}_{f}$ is a linear extension of $P$ (since $\mathbf{w}%
_{f}$ is a total order).]

Furthermore, we have $f\in\mathcal{E}\left(  \mathbf{w}_{f},\gamma\right)  $.

[\textit{Proof:} Let $x$ and $y$ be two elements of $\mathbf{w}_{f}$ such that
$x<y$ in $\mathbf{w}_{f}$. We shall prove the following three statements:

\begin{enumerate}
\item[\textbf{(i)}] We have $f\left(  x\right)  \preccurlyeq f\left(
y\right)  $.

\item[\textbf{(ii)}] If $f\left(  x\right)  =f\left(  y\right)  =+n$ for some
$n\in\mathcal{N}$, then $\gamma\left(  x\right)  <\gamma\left(  y\right)  $.

\item[\textbf{(iii)}] If $f\left(  x\right)  =f\left(  y\right)  =-n$ for some
$n\in\mathcal{N}$, then $\gamma\left(  x\right)  >\gamma\left(  y\right)  $.
\end{enumerate}

Indeed, $x$ and $y$ are elements of $P$ (since $\mathbf{w}_{f}=P$ as sets).
Also, $x<y$ in $\mathbf{w}_{f}$. In other words, $x<_{f}y$ (since $<_{f}$ is
the smaller relation of $\mathbf{w}_{f}$). In other words, we have%
\begin{align*}
& \text{\textbf{either} }\left(  f\left(  x\right)  \prec f\left(  y\right)
\right)  \\
& \text{\textbf{or} }\left(  f\left(  x\right)  =f\left(  y\right)  \text{ is
positive, and }\gamma\left(  x\right)  <\gamma\left(  y\right)  \right)  \\
& \text{\textbf{or} }\left(  f\left(  x\right)  =f\left(  y\right)  \text{ is
negative, and }\gamma\left(  x\right)  >\gamma\left(  y\right)  \right)
\end{align*}
(by the definition of the relation $<_{f}$). In either of these three cases,
the three statements \textbf{(i)}, \textbf{(ii)} and \textbf{(iii)} above are true.

Now, forget that we fixed $x$ and $y$. We thus have shown that if $x$ and $y$
are two elements of $\mathbf{w}_{f}$ such that $x<y$ in $\mathbf{w}_{f}$, then
the three statements \textbf{(i)}, \textbf{(ii)} and \textbf{(iii)} above are
true. In other words, $f$ is a $\mathcal{Z}$-enriched $\left(  \mathbf{w}%
_{f},\gamma\right)  $-partition (by the definition of a $\mathcal{Z}%
$-enriched $\left(  \mathbf{w}_{f},\gamma\right)  $-partition). In other
words, $f\in\mathcal{E}\left(  \mathbf{w}_{f},\gamma\right)  $ (by the
definition of $\mathcal{E}\left(  \mathbf{w}_{f},\gamma\right)  $).]

Furthermore, the following holds:

\begin{statement}
\textit{Observation 1:} If $w\in\mathcal{L}\left(  P\right)  $ satisfies
$f\in\mathcal{E}\left(  w,\gamma\right)  $, then $w=\mathbf{w}_{f}$.
\end{statement}

[\textit{Proof of Observation 1:} Let $w\in\mathcal{L}\left(  P\right)  $ be
such that $f\in\mathcal{E}\left(  w,\gamma\right)  $. Thus, $f$ is a
$\mathcal{Z}$-enriched $\left(  w,\gamma\right)  $-partition (by the
definition of $\mathcal{E}\left(  w,\gamma\right)  $).

From $w \in \mathcal{L} \left( P \right)$, we conclude that $w = P$
as sets, and that $w$ is a total order.
Hence, $w = P = \mathbf{w}_f$ as sets.

Now, let $x$ and $y$ be two elements of $w$ such that $x<y$ in $w$. Since $f$
is a $\mathcal{Z}$-enriched $\left(  w,\gamma\right)  $-partition, we thus
conclude that the following three statements hold:

\begin{enumerate}
\item[\textbf{(i)}] We have $f\left(  x\right)  \preccurlyeq f\left(
y\right)  $.

\item[\textbf{(ii)}] If $f\left(  x\right)  =f\left(  y\right)  =+n$ for some
$n\in\mathcal{N}$, then $\gamma\left(  x\right)  <\gamma\left(  y\right)  $.

\item[\textbf{(iii)}] If $f\left(  x\right)  =f\left(  y\right)  =-n$ for some
$n\in\mathcal{N}$, then $\gamma\left(  x\right)  >\gamma\left(  y\right)  $.
\end{enumerate}

From these three statements, we conclude that \textbf{either} $\left(
f\left(  x\right)  \prec f\left(  y\right)  \right)  $ \newline
\textbf{or} $\left(
f\left(  x\right)  =f\left(  y\right)  \text{ is positive, and }\gamma\left(
x\right)  <\gamma\left(  y\right)  \right)  $ \newline
\textbf{or} $\left(  f\left(
x\right)  =f\left(  y\right)  \text{ is negative, and }\gamma\left(  x\right)
>\gamma\left(  y\right)  \right)  $. In other words, $x<_{f}y$ (by the
definition of the relation $<_{f}$). In other words, $x<y$ in $\mathbf{w}_{f}$
(since $<_{f}$ is the smaller relation of $\mathbf{w}_{f}$).

Now, forget that we fixed $x$ and $y$. We thus have shown that if $x$ and $y$
are two elements of $w$ such that $x<y$ in $w$, then $x<y$ in $\mathbf{w}_{f}%
$. In other words, $\mathbf{w}_{f}$ is a linear extension of $w$ (since
$\mathbf{w}_{f}$ is a total order). But the only linear extension of $w$ is
$w$ itself (since $w$ is a total order). Thus, we conclude that $\mathbf{w}%
_{f}=w$. This proves Observation 1.]

Now, forget that we fixed $f$. Thus, for each $f\in\mathcal{E}\left(
P,\gamma\right)  $, we have constructed a linear extension $\mathbf{w}_{f}$ of
$P$, and we have shown that it satisfies $f\in\mathcal{E}\left(
\mathbf{w}_{f},\gamma\right)  $ and Observation 1.

Hence, for each $f\in\mathcal{E}\left(  P,\gamma\right)  $, we have%
\[
f\in\mathcal{E}\left(  \mathbf{w}_{f},\gamma\right)  \subseteq\bigcup
_{w\in\mathcal{L}\left(  P\right)  }\mathcal{E}\left(  w,\gamma\right)
\]
(since $\mathbf{w}_{f}\in\mathcal{L}\left(  P\right)  $ (because
$\mathbf{w}_{f}$ is a linear extension of $P$)). In other words,
$\mathcal{E}\left(  P,\gamma\right)  \subseteq\bigcup_{w\in\mathcal{L}\left(
P\right)  }\mathcal{E}\left(  w,\gamma\right)  $. On the other hand, every
$w\in\mathcal{L}\left(  P\right)  $ satisfies $\mathcal{E}\left(
w,\gamma\right)  \subseteq\mathcal{E}\left(  P,\gamma\right)  $%
\ \ \ \ \footnote{\textit{Proof.} Let $w\in\mathcal{L}\left(  P\right)  $.
Thus, $w$ is a linear extension of $P$. Hence, $w = P$ as sets, and every
two elements $x$ and $y$ of $P$ satisfying $x < y$ in $P$ must also
satisfy $x < y$ in $w$.
Thus, every $\mathcal{Z}$-enriched
$\left(  w,\gamma\right)  $-partition is a $\mathcal{Z}$-enriched $\left(
P,\gamma\right)  $-partition (since the axioms for a $\left(  w,\gamma\right)
$-partition are at least as strong as those for a $\left(  P,\gamma\right)
$-partition). In other words, $\mathcal{E}\left(  w,\gamma\right)
\subseteq\mathcal{E}\left(  P,\gamma\right)  $, qed.}. Hence, $\bigcup
_{w\in\mathcal{L}\left(  P\right)  }\mathcal{E}\left(  w,\gamma\right)
\subseteq\mathcal{E}\left(  P,\gamma\right)  $. Combining this with
$\mathcal{E}\left(  P,\gamma\right)  \subseteq\bigcup_{w\in\mathcal{L}\left(
P\right)  }\mathcal{E}\left(  w,\gamma\right)  $, we obtain%
\begin{equation}
\mathcal{E}\left(  P,\gamma\right)  =\bigcup_{w\in\mathcal{L}\left(  P\right)
}\mathcal{E}\left(  w,\gamma\right)  .\label{pf.prop.fund-lem.union}%
\end{equation}

Finally, the sets $\mathcal{E}\left(  w,\gamma\right)  $ for distinct
$w\in\mathcal{L}\left(  P\right)  $ are disjoint.

[\textit{Proof:} Let $u$ and $v$ be two elements of $\mathcal{L}\left(
P\right)  $ such that $\mathcal{E}\left(  u,\gamma\right)  \cap\mathcal{E}%
\left(  v,\gamma\right)  \neq\varnothing$. We shall show that $u=v$.

We have $\mathcal{E}\left(  u,\gamma\right)  \cap\mathcal{E}\left(
v,\gamma\right)  \neq\varnothing$. Hence, there exists some $f\in
\mathcal{E}\left(  u,\gamma\right)  \cap\mathcal{E}\left(  v,\gamma\right)  $.
Consider this $f$.

We have%
\begin{align*}
f  & \in\mathcal{E}\left(  u,\gamma\right)  \cap\mathcal{E}\left(
v,\gamma\right)  \subseteq\mathcal{E}\left(  u,\gamma\right)  \subseteq
\bigcup_{w\in\mathcal{L}\left(  P\right)  }\mathcal{E}\left(  w,\gamma\right)
\ \ \ \ \ \ \ \ \ \ \left(  \text{since }u\in\mathcal{L}\left(  P\right)
\right)  \\
& =\mathcal{E}\left(  P,\gamma\right)  .
\end{align*}
Hence, Observation 1 (applied to $w=u$) yields $u=\mathbf{w}_{f}$. The same
argument (applied to $v$ instead of $u$) yields $v=\mathbf{w}_{f}$. Hence,
$u=\mathbf{w}_{f}=v$.

Now, forget that we fixed $u$ and $v$. We thus have shown that if $u$  and $v$
are two elements of $\mathcal{L}\left(  P\right)  $ such that $\mathcal{E}%
\left(  u,\gamma\right)  \cap\mathcal{E}\left(  v,\gamma\right)
\neq\varnothing$, then $u=v$. In other words, the sets $\mathcal{E}\left(
w,\gamma\right)  $ for distinct $w\in\mathcal{L}\left(  P\right)  $ are disjoint.]

Thus, the union $\bigcup_{w\in\mathcal{L}\left(  P\right)  }\mathcal{E}\left(
w,\gamma\right)  $ is a disjoint union. Hence, (\ref{pf.prop.fund-lem.union})
rewrites as $\mathcal{E}\left(  P,\gamma\right)  =\bigsqcup_{w\in
\mathcal{L}\left(  P\right)  }\mathcal{E}\left(  w,\gamma\right)  $. This
proves Proposition \ref{prop.fund-lem}.
\end{proof}
\end{verlong}

\begin{definition}
\label{def.GammaZ}Let $\left(  P,\gamma\right)  $ be a labeled poset. We
define a power series $\Gamma_{\mathcal{Z}}\left(  P,\gamma\right)
\in\operatorname*{Pow}\mathcal{N}$ by%
\[
\Gamma_{\mathcal{Z}}\left(  P,\gamma\right)  =\sum_{f\in\mathcal{E}\left(
P,\gamma\right)  }\prod_{p\in P}x_{\left\vert f\left(  p\right)  \right\vert
}.
\]
This is easily seen to be convergent in the usual topology on
$\operatorname*{Pow}\mathcal{N}$
(that is, the topology of coefficientwise convergence%
\footnote{This is the topology obtained by regarding
$\operatorname*{Pow}\mathcal{N}$ as a Cartesian product of
infinitely many copies of $\mathbb{Q}$ (each endowed with the
discrete topology), corresponding to the monomials.}).
(Indeed, for every monomial $\mathfrak{m}$ in $\operatorname*{Pow}\mathcal{N}$,
there exist at most $\left|P\right|! \cdot 2^{\left|P\right|}$ many
$f \in \mathcal{E} \left(P, \gamma\right)$ satisfying
$\prod_{p\in P}x_{\left\vert f\left(  p\right)  \right\vert } = \mathfrak{m}$.)
\end{definition}

\begin{corollary}
\label{cor.fund-lem}For any labeled poset $\left(  P,\gamma\right)  $, we have%
\[
\Gamma_{\mathcal{Z}}\left(  P,\gamma\right)  =\sum_{w\in\mathcal{L}\left(
P\right)  }\Gamma_{\mathcal{Z}}\left(  w,\gamma\right)  .
\]

\end{corollary}

\begin{proof}
[Proof of Corollary \ref{cor.fund-lem}.]Follows straight from Proposition
\ref{prop.fund-lem}.
\end{proof}

\begin{definition}
Let $P$ be any set. Let $A$ be a totally ordered set. Let $\gamma:P\rightarrow
A$ and $\delta:P\rightarrow A$ be two maps. We say that $\gamma$ and $\delta$
are \textit{order-isomorphic} if the following holds: For every pair $\left(
p,q\right)  \in P\times P$, we have $\gamma\left(  p\right)  \leq\gamma\left(
q\right)  $ if and only if $\delta\left(  p\right)  \leq\delta\left(
q\right)  $.
\end{definition}

\begin{lemma} \label{lem.ord-eq-Ppar}
Let $\left(P, \alpha\right)$ and $\left(P, \beta\right)$ be two labeled posets
with the same ground poset $P$.
Assume that the maps $\alpha$ and $\beta$ are order-isomorphic. Then:

\textbf{(a)} We have
$\mathcal{E}\left(P, \alpha\right) = \mathcal{E}\left(P, \beta\right)$.

\textbf{(b)} We have
$\Gamma_{\mathcal{Z}}\left(P, \alpha\right) = \Gamma_{\mathcal{Z}}\left(P, \beta\right)$.
\end{lemma}

\begin{proof}[Proof of Lemma~\ref{lem.ord-eq-Ppar}.]
\textbf{(a)} If $x$ and $y$ are two elements of $P$, then we have the following
equivalences:
\begin{align*}
\left( \alpha\left(x\right) \leq \alpha\left(y\right) \right)
  \ &\Longleftrightarrow \ %
  \left( \beta\left(x\right) \leq \beta\left(y\right) \right) ; \\
\left( \alpha\left(x\right) > \alpha\left(y\right) \right)
  \ &\Longleftrightarrow \ %
  \left( \beta\left(x\right) > \beta\left(y\right) \right) ; \\
\left( \alpha\left(x\right) < \alpha\left(y\right) \right)
  \ &\Longleftrightarrow \ %
  \left( \beta\left(x\right) < \beta\left(y\right) \right) .
\end{align*}
(Indeed, the first of these equivalences holds because $\alpha$ and $\beta$
are order-isomorphic; the second is the contrapositive of the first; the
third is obtained from the second by swapping $x$ with $y$.)

Hence, the conditions ``$\alpha\left(x\right) > \alpha\left(y\right)$''
and ``$\alpha\left(x\right) < \alpha\left(y\right)$'' in the definition
of a $\mathcal{Z}$-enriched $\left(P, \alpha\right)$-partition are
equivalent to the conditions
``$\beta\left(x\right) > \beta\left(y\right)$''
and ``$\beta\left(x\right) < \beta\left(y\right)$'' in the definition
of a $\mathcal{Z}$-enriched $\left(P, \beta\right)$-partition.
Therefore, the $\mathcal{Z}$-enriched $\left(P, \alpha\right)$-partitions
are precisely the
$\mathcal{Z}$-enriched $\left(P, \beta\right)$-partitions.
In other words,
$\mathcal{E}\left(P, \alpha\right) = \mathcal{E}\left(P, \beta\right)$.
This proves Lemma~\ref{lem.ord-eq-Ppar} \textbf{(a)}.

\textbf{(b)} Lemma~\ref{lem.ord-eq-Ppar} \textbf{(b)} follows from
Lemma~\ref{lem.ord-eq-Ppar} \textbf{(a)}.
\end{proof}

Let us recall the notion of the disjoint union of two posets:

\begin{definition}
\textbf{(a)} Let $P$ and $Q$ be two sets.
The \textit{disjoint union} of $P$ and $Q$ is the set
$\left( \left\{0\right\} \times P \right)
              \cup \left( \left\{1\right\} \times Q \right)$.
This set is denoted by $P \sqcup Q$,
and comes with two canonical injections
\begin{align*}
\iota_0 &: P \to P \sqcup Q, \qquad  p \mapsto \left(0, p\right), \qquad \text{and} \\
\iota_1 &: Q \to P \sqcup Q, \qquad  q \mapsto \left(1, q\right).
\end{align*}
The images of these two injections are disjoint, and their
union is $P \sqcup Q$.

If $f : P \sqcup Q \to X$ is any map, then the
\textit{restriction of $f$ to $P$} is understood to be the
map $f \circ \iota_0 : P \to X$, whereas the
\textit{restriction of $f$ to $Q$} is understood to be the
map $f \circ \iota_1 : Q \to X$.
(Of course, this notation is ambiguous when $P = Q$.)

When the sets $P$ and $Q$ are already disjoint, it is common
to identify their disjoint union $P \sqcup Q$ with their
union $P \cup Q$ via the map
\[
P \sqcup Q \to P \cup Q, \qquad \left(i, r\right) \mapsto r .
\]
Under this identification, the restriction of a map
$f : P \sqcup Q \to X$ to $P$ becomes identical with the
(literal) restriction $f\mid_P$ of the map $f : P \cup Q \to X$
(and similarly for the restrictions to $Q$).

\textbf{(b)} Let $P$ and $Q$ be two posets.
The \textit{disjoint union} of the posets $P$ and $Q$ is the poset
$P \sqcup Q$ whose ground set is the disjoint union $P \sqcup Q$,
and whose order relation is defined by the following rules:
\begin{itemize}
\item If $p$ and $p^{\prime}$ are two elements of $P$, then
      $\left(0, p\right) < \left(0, p^{\prime}\right)$ in $P \sqcup Q$
      if and only if $p < p^{\prime}$ in $P$.
\item If $q$ and $q^{\prime}$ are two elements of $Q$, then
      $\left(1, q\right) < \left(1, q^{\prime}\right)$ in $P \sqcup Q$
      if and only if $q < q^{\prime}$ in $Q$.
\item If $p \in P$ and $q \in Q$,
      then the elements $\left(0, p\right)$ and $\left(1, q\right)$
      of $P \sqcup Q$ are incomparable.
\end{itemize}
\end{definition}

\begin{proposition}
\label{prop.prod1}Let $\left(  P,\gamma\right)  $ and $\left(  Q,\delta
\right)  $ be two labeled posets. Let $\left(  P\sqcup Q,\varepsilon\right)  $
be a labeled poset whose ground poset $P\sqcup Q$ is the disjoint union of $P$
and $Q$, and whose labeling $\varepsilon$ is such that the restriction of
$\varepsilon$ to $P$ is order-isomorphic to $\gamma$ and such that the
restriction of $\varepsilon$ to $Q$ is order-isomorphic to $\delta$. Then,%
\[
\Gamma_{\mathcal{Z}}\left(  P,\gamma\right)  \Gamma_{\mathcal{Z}}\left(
Q,\delta\right)  =\Gamma_{\mathcal{Z}}\left(  P\sqcup Q,\varepsilon\right)  .
\]

\end{proposition}

\begin{proof}
[Proof of Proposition \ref{prop.prod1}.]We WLOG assume that the ground sets
$P$ and $Q$ are disjoint; thus, we can identify $P\sqcup Q$ with the union
$P\cup Q$.
Let us make this identification.

The restriction $\varepsilon\mid_P$ of $\varepsilon$ to $P$ is order-isomorphic
to $\gamma$. Hence, Lemma \ref{lem.ord-eq-Ppar} \textbf{(a)} (applied to
$\alpha = \varepsilon\mid_P$ and $\beta = \gamma$) yields
$\mathcal{E}\left(P, \varepsilon\mid_P\right)
= \mathcal{E}\left(P, \gamma\right)$.
Similarly,
$\mathcal{E}\left(Q, \varepsilon\mid_Q\right)
= \mathcal{E}\left(Q, \delta\right)$.

\begin{vershort}
It is easy to see that a map $f : P \sqcup Q \to \mathcal{Z}$
is a $\mathcal{Z}$-enriched $\left(P \sqcup Q, \varepsilon\right)$-partition
if and only if $f \mid_P$ is a $\mathcal{Z}$-enriched
$\left(  P,\varepsilon\mid_P\right)$-partition and
$f \mid_Q$ is a $\mathcal{Z}$-enriched
$\left(  Q,\varepsilon\mid_Q\right)$-partition.
\end{vershort}

\begin{verlong}
If $x$ and $y$ are two elements of $P \sqcup Q$, then $x < y$ in $P \sqcup Q$
holds if and only if
\begin{itemize}
\item \textbf{either} $x$ and $y$ both belong to $P$ and satisfy $x < y$ in $P$,
\item \textbf{or} $x$ and $y$ both belong to $Q$ and satisfy $x < y$ in $Q$.
\end{itemize}
Hence:
\begin{itemize}
\item If $f$ is a
$\mathcal{Z}$-enriched $\left(P \sqcup Q, \varepsilon\right)$-partition,
then $f \mid_P$ is a $\mathcal{Z}$-enriched
$\left(  P,\varepsilon\mid_P\right)$-partition and
$f \mid_Q$ is a $\mathcal{Z}$-enriched
$\left(  Q,\varepsilon\mid_Q\right)$-partition.
\item Conversely, if $g$ is a $\mathcal{Z}$-enriched
$\left(  P,\varepsilon\mid_P\right)$-partition, and if
$h$ is a $\mathcal{Z}$-enriched
$\left(  Q,\varepsilon\mid_Q\right)$-partition, then the unique map
$f : P \sqcup Q \to \mathcal{Z}$ satisfying
$\left(  f\mid_{P},f\mid_{Q}\right) = \left(g, h\right)$
(that is, the map $P \sqcup Q \to \mathcal{Z}$
that sends each $p \in P$ to $g \left( p \right)$
and sends each $q \in Q$ to $h \left( q \right)$) is
a $\mathcal{Z}$-enriched $\left(P \sqcup Q, \varepsilon\right)$-partition.
\end{itemize}
\end{verlong}

Therefore, the map
\begin{align*}
\mathcal{E}\left(  P\sqcup Q,\varepsilon\right)   &  \rightarrow
\mathcal{E}\left(  P,\varepsilon\mid_P\right)
\times\mathcal{E}\left(  Q,\varepsilon\mid_Q\right)
, \\
f  &  \mapsto\left(  f\mid_{P},f\mid_{Q}\right)
\end{align*}
is a bijection (this is easy to see). In other words, the map
\begin{align}
\mathcal{E}\left(  P\sqcup Q,\varepsilon\right)   &  \rightarrow
\mathcal{E}\left(  P,\gamma\right)  \times\mathcal{E}\left(  Q,\delta\right)
,\nonumber\\
f  &  \mapsto\left(  f\mid_{P},f\mid_{Q}\right)  \label{pf.prop.prod1.bij}%
\end{align}
is a bijection (since $\mathcal{E}\left(P, \varepsilon\mid_P\right)
= \mathcal{E}\left(P, \gamma\right)$
and
$\mathcal{E}\left(Q, \varepsilon\mid_Q\right)
= \mathcal{E}\left(Q, \delta\right)$).
Now, the definition of
$\Gamma_{\mathcal{Z}}\left(  P\sqcup Q,\varepsilon\right)$
yields
\begin{align*}
\Gamma_{\mathcal{Z}}\left(  P\sqcup Q,\varepsilon\right)   &  =\sum
_{f\in\mathcal{E}\left(  P\sqcup Q,\varepsilon\right)  }\underbrace{\prod
_{p\in P\sqcup Q}x_{\left\vert f\left(  p\right)  \right\vert }}_{=\left(
\prod_{p\in P}x_{\left\vert f\left(  p\right)  \right\vert }\right)  \left(
\prod_{p\in Q}x_{\left\vert f\left(  p\right)  \right\vert }\right)  }\\
&  = \sum_{f\in\mathcal{E}\left(  P\sqcup Q,\varepsilon\right)  }
\underbrace{
 \left( \prod_{p\in P}x_{\left\vert f\left(  p\right)  \right\vert }\right)
}_{
 =
 \prod_{p\in P}x_{\left\vert \left(f\mid_P\right) \left(  p\right)  \right\vert }
}
\underbrace{
 \left( \prod_{p\in Q}x_{\left\vert f\left(  p\right)  \right\vert }\right)
}_{
 =
 \prod_{p\in Q}x_{\left\vert \left(f\mid_Q\right) \left(  p\right)  \right\vert }
}
\\
&  = \sum_{f\in\mathcal{E}\left(  P\sqcup Q,\varepsilon\right)  }
 \left( \prod_{p\in P}x_{\left\vert \left(f\mid_P\right) \left(  p\right)  \right\vert }
 \right)
 \left(
 \prod_{p\in Q}x_{\left\vert \left(f\mid_Q\right) \left(  p\right)  \right\vert }
 \right)
\\
&  =
\sum_{\left(g, h\right) \in \mathcal{E}\left(  P,\gamma\right) 
            \times \mathcal{E}\left(  Q,\delta\right) }
\left( \prod_{p\in P}x_{\left\vert g\left(  p\right)  \right\vert }\right)
\left(\prod_{p\in Q}x_{\left\vert
h\left(  p\right)  \right\vert }\right) \\
&  \ \ \ \ \ \ \ \ \ \ \left(
\begin{array}
[c]{c}%
\text{here, we have substituted }\left(  g,h\right)  \text{ for }\left(
f\mid_{P},f\mid_{Q}\right)  \text{,}\\
\text{since the map (\ref{pf.prop.prod1.bij}) is a bijection}%
\end{array}
\right) \\
&  =\underbrace{\left(  \sum_{g\in\mathcal{E}\left(  P,\gamma\right)  }%
\prod_{p\in P}x_{\left\vert g\left(  p\right)  \right\vert }\right)
}_{= \sum_{f\in\mathcal{E}\left(  P,\gamma\right)  }%
\prod_{p\in P}x_{\left\vert f\left(  p\right)  \right\vert }
= \Gamma_{\mathcal{Z}}\left(  P,\gamma\right)  }
\ \cdot \ \underbrace{\left(
\sum_{h\in\mathcal{E}\left(  Q,\delta\right)  }\prod_{p\in Q}x_{\left\vert
h\left(  p\right)  \right\vert }\right)  }_{
= \sum_{f\in\mathcal{E}\left(  Q,\delta\right)  }\prod_{p\in Q}x_{\left\vert
f\left(  p\right)  \right\vert }
= \Gamma_{\mathcal{Z}}\left( Q,\delta\right)  }\\
&  =\Gamma_{\mathcal{Z}}\left(  P,\gamma\right)  \Gamma_{\mathcal{Z}}\left(
Q,\delta\right)  .
\end{align*}
This proves Proposition \ref{prop.prod1}.
\end{proof}

\begin{definition}
Let $n\in\mathbb{N}$. Let $\pi$ be any $n$-permutation. (Recall that we have
defined the concept of an \textquotedblleft$n$-permutation\textquotedblright%
\ in Definition \ref{def.perm}.) Then, $\left(  \left[  n\right]  ,\pi\right)
$ is a labeled poset (in fact, $\pi$ is an injective map $\left[  n\right]
\rightarrow\left\{  1,2,3,\ldots\right\}  $, and thus can be considered a
labeling). We define $\Gamma_{\mathcal{Z}}\left(  \pi\right)  $ to be the
power series $\Gamma_{\mathcal{Z}}\left(  \left[  n\right]  ,\pi\right)  $.
\end{definition}

Let us recall the concept of a ``poset homomorphism'':

\begin{definition}
Let $P$ and $Q$ be two posets.
A map $f : P \to Q$ is said to be a \textit{poset homomorphism} if
for any two elements $x$ and $y$ of $P$ satisfying $x \leq y$ in $P$,
we have $f\left(x\right) \leq f\left(y\right)$ in $Q$.
\end{definition}

It is well-known that if $U$ and $V$ are any two finite totally ordered sets
of the same size, then there is a unique poset isomorphism $U \to V$.
Thus, if $w$ is a finite totally ordered set with $n$ elements, then there
is a unique poset isomorphism $w \to \left[ n \right]$.
Now, we claim the following:

\begin{proposition}
\label{prop.Gamma=Gamma}Let $w$ be a finite totally ordered set with ground
set $W$. Let $n=\left\vert W\right\vert $. Let $\overline{w}$ be the unique
poset isomorphism $w\rightarrow\left[  n\right]  $. Let $\gamma:W\rightarrow
\left\{  1,2,3,\ldots\right\}  $ be any injective map. Then, $\Gamma
_{\mathcal{Z}}\left(  w,\gamma\right)  =\Gamma_{\mathcal{Z}}\left(
\gamma\circ\overline{w}^{-1}\right)  $.
\end{proposition}

\begin{proof}
[Proof of Proposition \ref{prop.Gamma=Gamma}.]Clearly,
$\left( w , \gamma \right)$ is a labeled poset (since $\gamma$ is injective).
The map $\gamma\circ\overline
{w}^{-1}:\left[  n\right]  \rightarrow\left\{  1,2,3,\ldots\right\}  $ is an
injective map, thus an $n$-permutation. Hence, $\Gamma_{\mathcal{Z}}\left(
\gamma\circ\overline{w}^{-1}\right)  $ is well-defined, and its definition
yields
$\Gamma_{\mathcal{Z}}\left(  \gamma\circ\overline{w}^{-1}\right)
=\Gamma_{\mathcal{Z}}\left(  \left[  n\right]  ,\gamma\circ\overline{w}%
^{-1}\right)  $.
But $\overline{w}$ is a poset isomorphism $w \to \left[ n \right]$, and
thus is an isomorphism of labeled posets\footnote{We define the notion of
an ``isomorphism of labeled posets'' in the obvious way:
If $\left(P, \alpha\right)$ and $\left(Q, \beta\right)$ are two labeled
posets, then a \textit{homomorphism of labeled posets} from
$\left(P, \alpha\right)$ to $\left(Q, \beta\right)$ means a poset
homomorphism $f : P \to Q$ satisfying $\alpha = \beta \circ f$.
A \textit{isomorphism of labeled posets} is an invertible homomorphism
of labeled posets whose inverse also is a homomorphism of labeled posets.
Note that this definition of an isomorphism is not equivalent to the
definition given in \cite[Section 1.1]{Stembr97}.}
from $\left(  w,\gamma\right)$ to $\left(  \left[
n\right]  ,\gamma\circ\overline{w}^{-1}\right)$.
Hence,
\begin{align*}
\mathcal{E}\left(  w,\gamma\right)   &  \rightarrow\mathcal{E}\left(  \left[
n\right]  ,\gamma\circ\overline{w}^{-1}\right)  ,\\
f  &  \mapsto f\circ\overline{w}^{-1}%
\end{align*}
is a bijection (since any isomorphism of labeled posets induces a bijection
between their $\mathcal{Z}$-enriched $\left(P,\gamma\right)$-partitions). %
\begin{vershort}
Furthermore, it satisfies $\prod_{p\in w}x_{\left\vert
f\left(  p\right)  \right\vert }=\prod_{p\in\left[  n\right]  }x_{\left\vert
\left(  f\circ\overline{w}^{-1}\right)  \left(  p\right)  \right\vert }$ for
each $f\in\mathcal{E}\left(  w,\gamma\right)  $. Hence, $\Gamma_{\mathcal{Z}%
}\left(  w,\gamma\right)  =\Gamma_{\mathcal{Z}}\left(  \left[  n\right]
,\gamma\circ\overline{w}^{-1}\right)  =\Gamma_{\mathcal{Z}}\left(  \gamma
\circ\overline{w}^{-1}\right)  $.
\end{vershort}
\begin{verlong}
Every $f\in\mathcal{E}\left(  w,\gamma\right)  $ satisfies%
\begin{align}
\prod_{p\in w}x_{\left\vert f\left(  p\right)  \right\vert }  & =\prod
_{p\in\left[  n\right]  }\underbrace{x_{\left\vert f\left(  \overline{w}%
^{-1}\left(  p\right)  \right)  \right\vert }}_{=x_{\left\vert \left(
f\circ\overline{w}^{-1}\right)  \left(  p\right)  \right\vert }}\nonumber\\
& \left(
\begin{array}
[c]{c}%
\text{here, we have substituted }\overline{w}^{-1}\left(  p\right)  \text{ for
}p\text{ in the product,}\\
\text{since }\overline{w}^{-1}:\left[  n\right]  \rightarrow w\text{ is a
bijection}%
\end{array}
\right)  \nonumber\\
& =\prod_{p\in\left[  n\right]  }x_{\left\vert \left(  f\circ\overline{w}%
^{-1}\right)  \left(  p\right)  \right\vert }%
.\label{pf.prop.Gamma=Gamma.prod=prod}%
\end{align}
But the definition of $\Gamma_{\mathcal{Z}}\left(  w,\gamma\right)  $ yields
\begin{align*}
\Gamma_{\mathcal{Z}}\left(  w,\gamma\right)    & =\sum_{f\in\mathcal{E}\left(
w,\gamma\right)  }\underbrace{\prod_{p\in w}x_{\left\vert f\left(  p\right)
\right\vert }}_{\substack{=\prod_{p\in\left[  n\right]  }x_{\left\vert \left(
f\circ\overline{w}^{-1}\right)  \left(  p\right)  \right\vert }\\\text{(by
(\ref{pf.prop.Gamma=Gamma.prod=prod}))}}}=\sum_{f\in\mathcal{E}\left(
w,\gamma\right)  }\prod_{p\in\left[  n\right]  }x_{\left\vert \left(
f\circ\overline{w}^{-1}\right)  \left(  p\right)  \right\vert }\\
& =\sum_{f\in\mathcal{E}\left(  \left[  n\right]  ,\gamma\circ\overline
{w}^{-1}\right)  }\prod_{p\in\left[  n\right]  }x_{\left\vert f\left(
p\right)  \right\vert }\\
& \ \ \ \ \ \ \ \ \ \ \left(
\begin{array}
[c]{c}%
\text{here, we have substituted }f\text{ for }f\circ\overline{w}^{-1}\text{ in
the sum, since}\\
\text{the map }\mathcal{E}\left(  w,\gamma\right)  \rightarrow\mathcal{E}%
\left(  \left[  n\right]  ,\gamma\circ\overline{w}^{-1}\right)  ,\ f\mapsto
f\circ\overline{w}^{-1}\\
\text{is a bijection}%
\end{array}
\right)  \\
& =\Gamma_{\mathcal{Z}}\left(  \left[  n\right]  ,\gamma\circ\overline{w}%
^{-1}\right)  \ \ \ \ \ \ \ \ \ \ \left(  \text{by the definition of }%
\Gamma_{\mathcal{Z}}\left(  \left[  n\right]  ,\gamma\circ\overline{w}%
^{-1}\right)  \right)  \\
& =\Gamma_{\mathcal{Z}}\left(  \gamma\circ\overline{w}^{-1}\right)  .
\end{align*}
This proves Proposition \ref{prop.Gamma=Gamma}.
\end{verlong}
\end{proof}

For the following corollary, let us recall that a bijective poset
homomorphism is not necessarily an isomorphism of posets (since its
inverse may or may not be a poset homomorphism).

\begin{corollary}
\label{cor.fund-lem2}Let $\left(  P,\gamma\right)  $ be a labeled poset. Let
$n=\left\vert P\right\vert $. Then,%
\[
\Gamma_{\mathcal{Z}}\left(  P,\gamma\right)  =\sum_{\substack{x:P\rightarrow
\left[  n\right]  \\\text{bijective poset}\\\text{homomorphism}}%
}\Gamma_{\mathcal{Z}}\left(  \gamma\circ x^{-1}\right)  .
\]

\end{corollary}

\begin{proof}
[Proof of Corollary \ref{cor.fund-lem2}.]For each totally ordered set $w$ with
ground set $P$, we let $\overline{w}$ be the unique poset isomorphism
$w\rightarrow\left[  n\right]  $. If $w$ is a linear extension of $P$, then
this map $\overline{w}$ is also a bijective poset homomorphism $P\rightarrow
\left[  n\right]  $ (since every poset homomorphism $w\rightarrow\left[
n\right]  $ is also a poset homomorphism $P\rightarrow\left[  n\right]  $).
Thus, for each $w\in\mathcal{L}\left(  P\right)  $, we have defined a
bijective poset homomorphism $\overline{w}:P\rightarrow\left[  n\right]  $. We
thus have defined a map%
\begin{align}
\mathcal{L}\left(  P\right)    & \rightarrow\left\{  \text{bijective poset
homomorphisms }P\rightarrow\left[  n\right]  \right\}  ,\nonumber\\
w  & \mapsto\overline{w}.\label{pf.cor.fund-lem2.bij}%
\end{align}
This map is injective (indeed, a linear extension $w\in\mathcal{L}\left(
P\right)  $ can be uniquely reconstructed from $\overline{w}$) and surjective
(because if $x$ is a bijective poset homomorphism $P\rightarrow\left[
n\right]  $, then the linear extension $w\in\mathcal{L}\left(  P\right)  $
defined (as a list) by $w=\left(  x^{-1}\left(  1\right)  ,x^{-1}\left(
2\right)  ,\ldots,x^{-1}\left(  n\right)  \right)  $ satisfies $x=\overline
{w}$). Hence, this map is a bijection. 

Corollary \ref{cor.fund-lem} yields%
\begin{equation}
\Gamma_{\mathcal{Z}}\left(  P,\gamma\right)  =\sum_{w\in\mathcal{L}\left(
P\right)  }\underbrace{\Gamma_{\mathcal{Z}}\left(  w,\gamma\right)
}_{\substack{=\Gamma_{\mathcal{Z}}\left(  \gamma\circ\overline{w}^{-1}\right)
\\\text{(by Proposition \ref{prop.Gamma=Gamma})}}}=\sum_{w\in\mathcal{L}%
\left(  P\right)  }\Gamma_{\mathcal{Z}}\left(  \gamma\circ\overline{w}%
^{-1}\right)  .\label{pf.cor.fund-lem2.1}%
\end{equation}
But recall that the map (\ref{pf.cor.fund-lem2.bij}) is a bijection. Thus, we
can substitute $x$ for $\overline{w}$ in the sum $\sum_{w\in\mathcal{L}\left(
P\right)  }\Gamma_{\mathcal{Z}}\left(  \gamma\circ\overline{w}^{-1}\right)  $,
obtaining%
\[
\sum_{w\in\mathcal{L}\left(  P\right)  }\Gamma_{\mathcal{Z}}\left(
\gamma\circ\overline{w}^{-1}\right)  =\sum_{\substack{x:P\rightarrow\left[
n\right]  \\\text{bijective poset}\\\text{homomorphism}}}\Gamma_{\mathcal{Z}%
}\left(  \gamma\circ x^{-1}\right)  .
\]
Hence, (\ref{pf.cor.fund-lem2.1}) becomes%
\[
\Gamma_{\mathcal{Z}}\left(  P,\gamma\right)  =\sum_{w\in\mathcal{L}\left(
P\right)  }\Gamma_{\mathcal{Z}}\left(  \gamma\circ\overline{w}^{-1}\right)
=\sum_{\substack{x:P\rightarrow\left[  n\right]  \\\text{bijective
poset}\\\text{homomorphism}}}\Gamma_{\mathcal{Z}}\left(  \gamma\circ
x^{-1}\right)  .
\]
This proves Corollary \ref{cor.fund-lem2}.
\end{proof}

\begin{corollary}
\label{cor.prod2}Let $n\in\mathbb{N}$ and $m\in\mathbb{N}$. Let $\pi$ be an
$n$-permutation and let $\sigma$ be an $m$-permutation such that $\pi$ and
$\sigma$ are disjoint. Then,%
\[
\Gamma_{\mathcal{Z}}\left(  \pi\right)  \Gamma_{\mathcal{Z}}\left(
\sigma\right)  =\sum_{\tau\in S\left(  \pi,\sigma\right)  }\Gamma
_{\mathcal{Z}}\left(  \tau\right)  .
\]

\end{corollary}

\begin{proof}
[Proof of Corollary \ref{cor.prod2}.]
\begin{verlong}
We first state a general fact about sets and maps:

\begin{statement}
\textit{Observation 0:} Let $U$, $V$ and $W$ be three sets. Let $\alpha
:U\rightarrow W$ and $\beta:V\rightarrow W$ be two injective maps such that
$\alpha\left(  U\right)  =\beta\left(  V\right)  $. Then, there exists a
unique bijective map $\lambda:U\rightarrow V$ satisfying $\alpha=\beta
\circ\lambda$.
\end{statement}

[\textit{Proof of Observation 0:} This is a simple exercise. (Roughly
speaking: If we replace the codomains of the maps $\alpha$ and $\beta$ by
$\alpha\left(  U\right)  =\beta\left(  V\right)  $, then the injective maps
$\alpha$ and $\beta$ become bijections, and the required map $\lambda$ becomes
$\beta^{-1}\circ\alpha$.)]

\end{verlong}
Consider the disjoint union $\left[  n\right]  \sqcup\left[  m\right]  $ of
the posets $\left[  n\right]  $ and $\left[  m\right]  $. (Note that this
disjoint union cannot be identified with the union $\left[  n\right]
\cup\left[  m\right]  $.) Let $\varepsilon$ be the map $\left[  n\right]
\sqcup\left[  m\right]  \rightarrow\left\{  1,2,3,\ldots\right\}  $ whose
restriction to $\left[  n\right]  $ is $\pi$ and whose restriction to $\left[
m\right]  $ is $\sigma$. This map $\varepsilon$ is injective, since $\pi$ and
$\sigma$ are disjoint permutations. Thus, $\left(  \left[  n\right]
\sqcup\left[  m\right]  ,\varepsilon\right)  $ is a labeled poset.

\begin{vershort}
The following two observations are easy to show (see \cite{verlong}
for detailed proofs):
\end{vershort}

\begin{verlong}
Let us make some observations:
\end{verlong}

\begin{statement}
\textit{Observation 1:} If $x$ is a bijective poset homomorphism $\left[
n\right]  \sqcup\left[  m\right]  \rightarrow\left[  n+m\right]  $, then
$\varepsilon\circ x^{-1}\in S\left(  \pi,\sigma\right)  $.
\end{statement}

\begin{verlong}
[\textit{Proof of Observation 1:} Let $x$ be a bijective poset homomorphism
$\left[  n\right]  \sqcup\left[  m\right]  \rightarrow\left[  n+m\right]  $.
Then, $x^{-1}$ is a well-defined bijective map $\left[  n+m\right]
\rightarrow\left[  n\right]  \sqcup\left[  m\right]  $ (since $x$ is
bijective). Hence, $\varepsilon\circ x^{-1}$ is an injective map $\left[
n+m\right]  \rightarrow\left\{  1,2,3,\ldots\right\}  $ (since $\varepsilon$
is injective), therefore an $\left(  n+m\right)  $-permutation. Moreover,
$\sigma$ is a subsequence of $\varepsilon\circ x^{-1}$%
\ \ \ \ \footnote{\textit{Proof.} Let us denote the elements $x\left(  \left(
1,1\right)  \right)  ,x\left(  \left(  1,2\right)  \right)  ,\ldots,x\left(
\left(  1,m\right)  \right)  $ of $\left[  n+m\right]  $ by $i_{1}%
,i_{2},\ldots,i_{m}$. Thus, $i_{q}=x\left(  \left(  1,q\right)  \right)  $ for
each $q\in\left[  m\right]  $. Hence, the elements $i_{1},i_{2},\ldots,i_{m}$
are the images of the distinct elements $\left(  1,1\right)  ,\left(
1,2\right)  ,\ldots,\left(  1,m\right)  $ under the map $x$. Therefore, these
elements $i_{1},i_{2},\ldots,i_{m}$ are distinct themselves (since $x$ is
injective).
\par
But $x$ is a poset homomorphism $\left[  n\right]  \sqcup\left[  m\right]
\rightarrow\left[  n+m\right]  $. Thus, $x\left(  \left(  1,1\right)  \right)
\leq x\left(  \left(  1,2\right)  \right)  \leq\cdots\leq x\left(  \left(
1,m\right)  \right)  $ (since $\left(  1,1\right)  \leq\left(  1,2\right)
\leq\cdots\leq\left(  1,m\right)  $ in $\left[  n\right]  \sqcup\left[
m\right]  $). In other words, $i_{1}\leq i_{2}\leq\cdots\leq i_{m}$ (since
$i_{q}=x\left(  \left(  1,q\right)  \right)  $ for each $q\in\left[  m\right]
$). Hence, $i_{1}<i_{2}<\cdots<i_{m}$ (since the elements $i_{1},i_{2}%
,\ldots,i_{m}$ are distinct). Hence,%
\[
\left(  \left(  \varepsilon\circ x^{-1}\right)  \left(  i_{1}\right)  ,\left(
\varepsilon\circ x^{-1}\right)  \left(  i_{2}\right)  ,\ldots,\left(
\varepsilon\circ x^{-1}\right)  \left(  i_{m}\right)  \right)
\]
is a subsequence of $\varepsilon\circ x^{-1}$. Since each $q\in\left[
m\right]  $ satisfies%
\begin{align*}
\left(  \varepsilon\circ x^{-1}\right)  \left(  \underbrace{i_{q}}_{=x\left(
\left(  1,q\right)  \right)  }\right)   &  =\left(  \varepsilon\circ
x^{-1}\right)  \left(  x\left(  \left(  1,q\right)  \right)  \right)
=\varepsilon\left(  \left(  1,q\right)  \right)  =\sigma\left(  q\right)  \\
&  \ \ \ \ \ \ \ \ \ \ \left(  \text{since the restriction of }\varepsilon
\text{ to }\left[  m\right]  \text{ is }\sigma\right)  ,
\end{align*}
we have%
\[
\left(  \left(  \varepsilon\circ x^{-1}\right)  \left(  i_{1}\right)  ,\left(
\varepsilon\circ x^{-1}\right)  \left(  i_{2}\right)  ,\ldots,\left(
\varepsilon\circ x^{-1}\right)  \left(  i_{m}\right)  \right)  =\left(
\sigma\left(  1\right)  ,\sigma\left(  2\right)  ,\ldots,\sigma\left(
m\right)  \right)  =\sigma.
\]
Thus, $\sigma$ is a subsequence of $\varepsilon\circ x^{-1}$ (since%
\[
\left(  \left(  \varepsilon\circ x^{-1}\right)  \left(  i_{1}\right)  ,\left(
\varepsilon\circ x^{-1}\right)  \left(  i_{2}\right)  ,\ldots,\left(
\varepsilon\circ x^{-1}\right)  \left(  i_{m}\right)  \right)
\]
is a subsequence of $\varepsilon\circ x^{-1}$). Qed.}. Similarly, $\pi$ is a
subsequence of $\varepsilon\circ x^{-1}$. Hence, $\varepsilon\circ x^{-1}$ is
an $\left(  n+m\right)  $-permutation having the property that both $\pi$ and
$\sigma$ are subsequences of $\varepsilon\circ x^{-1}$. In other words,
$\varepsilon\circ x^{-1}$ is a shuffle of $\pi$ and $\sigma$ (by the
definition of a shuffle). In other words, $\varepsilon\circ x^{-1}\in S\left(
\pi,\sigma\right)  $ (by the definition of $S\left(  \pi,\sigma\right)  $).
This proves Observation 1.]
\end{verlong}

\begin{statement}
\textit{Observation 2:} If $\tau\in S\left(  \pi,\sigma\right)  $, then there
exists a unique bijective poset homomorphism $x:\left[  n\right]
\sqcup\left[  m\right]  \rightarrow\left[  n+m\right]  $ satisfying
$\varepsilon\circ x^{-1}=\tau$.
\end{statement}

\begin{verlong}
[\textit{Proof of Observation 2:} Let $\tau\in S\left(  \pi,\sigma\right)  $.
Thus, $\tau$ is a shuffle of $\pi$ and $\sigma$ (by the definition of
$S\left(  \pi,\sigma\right)  $). In other words, $\tau$ is an $\left(
n+m\right)  $-permutation having the property that both $\pi$ and $\sigma$ are
subsequences of $\tau$ (by the definition of a shuffle). Since $\pi$ and
$\sigma$ are disjoint permutations, this entails that the letters of $\tau$
are the letters of $\pi$ and the letters of $\sigma$. Hence, the image
$\tau\left(  \left[  n+m\right]  \right)  $ of the map $\tau$ is the union
$\pi\left(  \left[  n\right]  \right)  \cup\sigma\left(  \left[  m\right]
\right)  $. But the image $\varepsilon\left(  \left[  n\right]  \sqcup\left[
m\right]  \right)  $ is the same union $\pi\left(  \left[  n\right]  \right)
\cup\sigma\left(  \left[  m\right]  \right)  $ (by the definition of
$\varepsilon$). Hence, $\varepsilon\left(  \left[  n\right]  \sqcup\left[
m\right]  \right)  =\tau\left(  \left[  n+m\right]  \right)  $. Also,
$\varepsilon:\left[  n\right]  \sqcup\left[  m\right]  \rightarrow\left\{
1,2,3,\ldots\right\}  $ and $\tau:\left[  n+m\right]  \rightarrow\left\{
1,2,3,\ldots\right\}  $ are injective maps (indeed, $\tau$ is injective since
$\tau$ is a permutation). Hence, Observation 0 (applied to $U=\left[
n\right]  \sqcup\left[  m\right]  $, $V=\left[  n+m\right]  $, $W=\left\{
1,2,3,\ldots\right\}  $, $\alpha=\varepsilon$ and $\beta=\tau$) shows that
there exists a unique bijective map $\lambda:\left[  n\right]  \sqcup\left[
m\right]  \rightarrow\left[  n+m\right]  $ satisfying $\varepsilon=\tau
\circ\lambda$. In other words, there exists a unique bijective map
$\lambda:\left[  n\right]  \sqcup\left[  m\right]  \rightarrow\left[
n+m\right]  $ satisfying $\varepsilon\circ\lambda^{-1}=\tau$ (since the
condition $\varepsilon=\tau\circ\lambda$ is equivalent to $\varepsilon
\circ\lambda^{-1}=\tau$). Consider this $\lambda$.

We have $\lambda\left(  \left(  0,1\right)  \right)  <\lambda\left(  \left(
0,2\right)  \right)  <\cdots<\lambda\left(  \left(  0,n\right)  \right)
$\ \ \ \ \footnote{\textit{Proof.} We know that $\pi$ is a subsequence of
$\tau$. In other words, there exist $n$ elements $i_{1}<i_{2}<\cdots<i_{n}$ of
$\left[  n+m\right]  $ satisfying $\pi=\left(  \tau\left(  i_{1}\right)
,\tau\left(  i_{2}\right)  ,\ldots,\tau\left(  i_{n}\right)  \right)  $.
Consider these $n$ elements.
\par
We have $\pi=\left(  \tau\left(  i_{1}\right)  ,\tau\left(  i_{2}\right)
,\ldots,\tau\left(  i_{n}\right)  \right)  $. In other words, $\pi\left(
p\right)  =\tau\left(  i_{p}\right)  $ for each $p\in\left[  n\right]  $.
\par
Now, let $p\in\left[  n\right]  $. Thus, $\pi\left(  p\right)  =\tau\left(
i_{p}\right)  $ (as we have just seen). But recall that $\varepsilon
\circ\lambda^{-1}=\tau$. In other words, $\varepsilon=\tau\circ\lambda$. On
the other hand, the restriction of $\varepsilon$ to $\left[  n\right]  $ is
$\pi$ (by the definition of $\varepsilon$). Hence, $\varepsilon\left(  \left(
0,p\right)  \right)  =\pi\left(  p\right)  =\tau\left(  i_{p}\right)  $. Thus,%
\[
\tau\left(  i_{p}\right)  =\underbrace{\varepsilon}_{=\tau\circ\lambda}\left(
\left(  0,p\right)  \right)  =\left(  \tau\circ\lambda\right)  \left(  \left(
0,p\right)  \right)  =\tau\left(  \lambda\left(  \left(  0,p\right)  \right)
\right)  .
\]
Since the map $\tau$ is injective (since $\tau$ is a permutation), we thus
obtain $i_{p}=\lambda\left(  \left(  0,p\right)  \right)  $.
\par
Now, forget that we fixed $p$. Thus, we have shown that $i_{p}=\lambda\left(
\left(  0,p\right)  \right)  $ for each $p\in\left[  n\right]  $. Hence, the
chain of inequalities $i_{1}<i_{2}<\cdots<i_{n}$ rewrites as $\lambda\left(
\left(  0,1\right)  \right)  <\lambda\left(  \left(  0,2\right)  \right)
<\cdots<\lambda\left(  \left(  0,n\right)  \right)  $. Qed.} and
$\lambda\left(  \left(  1,1\right)  \right)  <\lambda\left(  \left(
1,2\right)  \right)  <\cdots<\lambda\left(  \left(  1,m\right)  \right)
$\ \ \ \ \footnote{for similar reasons}. Thus, $\lambda$ is a poset
homomorphism $\left[  n\right]  \sqcup\left[  m\right]  \rightarrow\left[
n+m\right]  $. Hence, there exists at least one bijective poset homomorphism
$x:\left[  n\right]  \sqcup\left[  m\right]  \rightarrow\left[  n+m\right]  $
satisfying $\varepsilon\circ x^{-1}=\tau$ (namely, $x=\lambda$). Moreover,
this $x$ is unique (since $\lambda$ is the \textbf{only} bijective map
$\lambda:\left[  n\right]  \sqcup\left[  m\right]  \rightarrow\left[
n+m\right]  $ satisfying $\varepsilon\circ\lambda^{-1}=\tau$). Thus, there
exists a unique bijective poset homomorphism $x:\left[  n\right]
\sqcup\left[  m\right]  \rightarrow\left[  n+m\right]  $ satisfying
$\varepsilon\circ x^{-1}=\tau$. This proves Observation 2.]
\end{verlong}

Now, the map%
\begin{align*}
\left\{  \text{bijective poset homomorphisms }x:\left[  n\right]
\sqcup\left[  m\right]  \rightarrow\left[  n+m\right]  \right\}   &
\rightarrow S\left(  \pi,\sigma\right)  ,\\
x  &  \mapsto\varepsilon\circ x^{-1}%
\end{align*}
is well-defined (by Observation 1) and is a bijection (by Observation 2).
Hence, we can substitute $\varepsilon\circ x^{-1}$ for $\tau$ in the sum
$\sum_{\tau\in S\left(  \pi,\sigma\right)  }\Gamma_{\mathcal{Z}}\left(
\tau\right)  $. We thus obtain%
\begin{equation}
\sum_{\tau\in S\left(  \pi,\sigma\right)  }\Gamma_{\mathcal{Z}}\left(
\tau\right)  =\sum_{\substack{x:\left[  n\right]  \sqcup\left[  m\right]
\rightarrow\left[  n+m\right]  \\\text{bijective poset}\\\text{homomorphism}%
}}\Gamma_{\mathcal{Z}}\left(  \varepsilon\circ x^{-1}\right)  .
\label{pf.cor.prod2.sum=sum}%
\end{equation}

The definition of $\Gamma_{\mathcal{Z}}\left(  \pi\right)  $ yields
$\Gamma_{\mathcal{Z}}\left(  \pi\right)  =\Gamma_{\mathcal{Z}}\left(  \left[
n\right]  ,\pi\right)  $. The definition of $\Gamma_{\mathcal{Z}}\left(
\sigma\right)  $ yields $\Gamma_{\mathcal{Z}}\left(  \sigma\right)
=\Gamma_{\mathcal{Z}}\left(  \left[  m\right]  ,\sigma\right)  $. Multiplying
these two equalities, we obtain
\begin{align*}
\Gamma_{\mathcal{Z}}\left(  \pi\right)  \Gamma_{\mathcal{Z}}\left(
\sigma\right)   &  =\Gamma_{\mathcal{Z}}\left(  \left[  n\right]  ,\pi\right)
\Gamma_{\mathcal{Z}}\left(  \left[  m\right]  ,\sigma\right)  =\Gamma
_{\mathcal{Z}}\left(  \left[  n\right]  \sqcup\left[  m\right]  ,\varepsilon
\right) \\
&  \ \ \ \ \ \ \ \ \ \ \left(
\begin{array}
[c]{c}%
\text{by Proposition \ref{prop.prod1}, applied to }P=\left[  n\right]  \text{,
}\gamma=\pi\text{,}\\
Q=\left[  m\right]  \text{ and }\delta=\sigma
\end{array}
\right) \\
&  =\sum_{\substack{x:\left[  n\right]  \sqcup\left[  m\right]  \rightarrow
\left[  n+m\right]  \\\text{bijective poset}\\\text{homomorphism}}%
}\Gamma_{\mathcal{Z}}\left(  \varepsilon\circ x^{-1}\right) \\
&\ \ \ \ \ \ \ \ \ \ \left(
\begin{array}
[c]{c}%
\text{by Corollary \ref{cor.fund-lem2}, applied to}
\left[  n\right]  \sqcup\left[  m\right]  \text{, }\varepsilon \\
\text{ and
}n+m\text{ instead of }P\text{, }\gamma\text{ and }n
\end{array}
\right) \\
&  =\sum_{\tau\in S\left(  \pi,\sigma\right)  }\Gamma_{\mathcal{Z}}\left(
\tau\right)  \ \ \ \ \ \ \ \ \ \ \left(  \text{by (\ref{pf.cor.prod2.sum=sum}%
)}\right)  .
\end{align*}
This proves Corollary \ref{cor.prod2}.
\end{proof}

\subsection{Exterior peaks}

So far we have been doing general nonsense. Let us now specialize to a
situation that is connected to exterior peaks.

\begin{convention}
From now on, we set $\mathcal{N}=\left\{  0,1,2,\ldots\right\}  \cup\left\{
\infty\right\}  $, with total order given by $0\prec1\prec2\prec\cdots
\prec\infty$, and we set
\begin{align*}
\mathcal{Z}  &  =\left(  \mathcal{N}\times\left\{  +,-\right\}  \right)
\setminus\left\{  -0,+\infty\right\} \\
&  =\left\{  +0\right\}  \cup\left\{  +n\ \mid\ n\in\left\{  1,2,3,\ldots
\right\}  \right\}  \cup\left\{  -n\ \mid\ n\in\left\{  1,2,3,\ldots\right\}
\right\}  \cup\left\{  -\infty\right\}  .
\end{align*}
Recall that the total order on $\mathcal{Z}$ has%
\[
+0\prec-1\prec+1\prec-2\prec+2\prec\cdots\prec-\infty.
\]

\end{convention}

\begin{definition}
Let $S$ be a subset of $\mathbb{Z}$.
A map $\chi$ from $S$ to a totally ordered set $K$ is
said to be \textit{V-shaped} if there exists some $t\in S$ such that the map
$\chi\mid_{\left\{  s\in S\ \mid\ s\leq t\right\}  }$ is strictly decreasing
while the map $\chi\mid_{\left\{  s\in S\ \mid\ s\geq t\right\}  }$ is
strictly increasing. Notice that this $t\in S$ is uniquely determined in this
case; namely, it is the unique $k\in S$ that minimizes $\chi\left(  k\right)
$.
\end{definition}

Thus, roughly speaking, a map from a subset of $\ZZ$ to a totally
ordered set is \textit{V-shaped}
if and only if it is strictly decreasing up until a certain value of its
argument, and then strictly increasing afterwards.
For example, the $5$-permutation $\left(5, 1, 2, 3, 4\right)$ is V-shaped
(keep in mind that we regard $n$-permutations as injective maps
$\left[n\right] \to \mathbb{P}$), whereas the $4$-permutation
$\left(3, 1, 4, 2\right)$ is not.

\begin{verlong}
For later use, let us crystallize a simple criterion for V-shapedness:

\begin{lemma}
\label{lem.Vshaped.crit}Let $S$ be a finite nonempty subset of $\mathbb{Z}$.
Let $K$ be a totally ordered set. Let $\chi:S\rightarrow K$ be an injective
map. Let $I$ and $J$ be two subsets of $S$ such that $I\cup J=S$. Assume that
the map $\chi\mid_{I}$ is strictly decreasing, whereas the map $\chi\mid_{J}$
is strictly increasing. Assume further that each element of $I$ is smaller
than each element of $J$. Then, the map $\chi$ is V-shaped.
\end{lemma}

\begin{proof}
[Proof of Lemma \ref{lem.Vshaped.crit}.]The set $S$ is nonempty. Hence, there
exists some $k\in S$ that minimizes $\chi\left(  k\right)  $. Consider this
$k$. Thus,%
\begin{equation}
\chi\left(  k\right)  \leq\chi\left(  s\right)  \ \ \ \ \ \ \ \ \ \ \text{for
each }s\in S. \label{pf.lem.Vshaped.crit.min}%
\end{equation}

Now, we claim that if $a\in S$ and $b\in S$ satisfy $a<b$ and $\chi\left(
a\right)  <\chi\left(  b\right)  $, then%
\begin{equation}
b>k. \label{pf.lem.Vshaped.crit.bk}%
\end{equation}

[\textit{Proof of (\ref{pf.lem.Vshaped.crit.bk}):} Let $a\in S$ and $b\in S$
be such that $a<b$ and $\chi\left(  a\right)  <\chi\left(  b\right)  $. We
must prove that $b>k$.

Assume the contrary. Thus, $b\leq k$. But (\ref{pf.lem.Vshaped.crit.min})
(applied to $s=a$) yields $\chi\left(  k\right)  \leq\chi\left(  a\right)  $.
Hence, $\chi\left(  k\right)  \leq\chi\left(  a\right)  <\chi\left(  b\right)
$, so that $\chi\left(  k\right)  \neq\chi\left(  b\right)  $. Thus, $k\neq
b$. Combining this with $b\leq k$, we obtain $b<k$. Hence, $a<b<k$.

Assume (for the sake of contradiction) that $b\notin I$. Combining $b\in S$
with $b\notin I$, we obtain $b\in S\setminus I\subseteq J$ (since $I\cup
J=S$). Hence, $k\in J$\ \ \ \ \footnote{\textit{Proof.} Assume the contrary.
Thus, $k\notin J$. Combining $k\in S$ with $k\notin J$, we find $k\in
S\setminus J\subseteq I$ (since $I\cup J=S$). But recall that each element of
$I$ is smaller than each element of $J$. In other words, $i<j$ for each $i\in
I$ and each $j\in J$. Applying this to $i=k$ and $j=b$, we find $k<b$. This
contradicts $b<k$. This contradiction shows that our assumption was false,
qed.}. Now, the map $\chi\mid_{J}$ is strictly increasing. Thus, from $b<k$,
we obtain $\left(  \chi\mid_{J}\right)  \left(  b\right)  <\left(  \chi
\mid_{J}\right)  \left(  k\right)  $ (since $b\in J$ and $k\in J$). Thus,
$\chi\left(  b\right)  =\left(  \chi\mid_{J}\right)  \left(  b\right)
<\left(  \chi\mid_{J}\right)  \left(  k\right)  =\chi\left(  k\right)
\leq\chi\left(  b\right)  $ (by (\ref{pf.lem.Vshaped.crit.min}), applied to
$s=b$). This is absurd. This contradiction shows that our assumption (that
$b\notin I$) was false. Hence, we have $b\in I$.

The same argument (but using $a$ instead of $b$) shows that $a\in I$. Now,
recall that the map $\chi\mid_{I}$ is strictly decreasing. Hence, from $a<b$,
we obtain $\left(  \chi\mid_{I}\right)  \left(  a\right)  >\left(  \chi
\mid_{I}\right)  \left(  b\right)  $ (since $a\in I$ and $b\in I$). Thus,
$\chi\left(  a\right)  =\left(  \chi\mid_{I}\right)  \left(  a\right)
>\left(  \chi\mid_{I}\right)  \left(  b\right)  =\chi\left(  b\right)  $. This
contradicts $\chi\left(  a\right)  <\chi\left(  b\right)  $. This
contradiction shows that our assumption was wrong. Hence, $b>k$ is proven.
Thus, (\ref{pf.lem.Vshaped.crit.bk}) is proven.]

Next, we claim that if $a\in S$ and $b\in S$ satisfy $a<b$ and $\chi\left(
a\right)  >\chi\left(  b\right)  $, then%
\begin{equation}
a<k. \label{pf.lem.Vshaped.crit.ak}%
\end{equation}

[\textit{Proof of (\ref{pf.lem.Vshaped.crit.ak}):} Let $a\in S$ and $b\in S$
be such that $a<b$ and $\chi\left(  a\right)  >\chi\left(  b\right)  $. We
must prove that $a<k$.

Assume the contrary. Thus, $a\geq k$. But (\ref{pf.lem.Vshaped.crit.min})
(applied to $s=b$) yields $\chi\left(  k\right)  \leq\chi\left(  b\right)  $.
Hence, $\chi\left(  k\right)  \leq\chi\left(  b\right)  <\chi\left(  a\right)
$ (since $\chi\left(  a\right)  >\chi\left(  b\right)  $), so that
$\chi\left(  k\right)  \neq\chi\left(  a\right)  $. Thus, $k\neq a$. Combining
this with $a\geq k$, we obtain $a>k$. Thus, $k<a<b$.

Assume (for the sake of contradiction) that $a\notin J$. Combining $a\in S$
with $a\notin J$, we obtain $a\in S\setminus J\subseteq I$ (since $I\cup
J=S$). Hence, $k\in I$\ \ \ \ \footnote{\textit{Proof.} Assume the contrary.
Thus, $k\notin I$. Combining $k\in S$ with $k\notin I$, we find $k\in
S\setminus I\subseteq J$ (since $I\cup J=S$). But recall that each element of
$I$ is smaller than each element of $J$. In other words, $i<j$ for each $i\in
I$ and each $j\in J$. Applying this to $i=a$ and $j=k$, we find $a<k$. This
contradicts $k<a$. This contradiction shows that our assumption was false,
qed.}. Now, the map $\chi\mid_{I}$ is strictly decreasing. Thus, from $k<a$,
we obtain $\left(  \chi\mid_{I}\right)  \left(  k\right)  >\left(  \chi
\mid_{I}\right)  \left(  a\right)  $ (since $k\in I$ and $a\in I$). Thus,
$\chi\left(  k\right)  =\left(  \chi\mid_{I}\right)  \left(  k\right)
>\left(  \chi\mid_{I}\right)  \left(  a\right)  =\chi\left(  a\right)
\geq\chi\left(  k\right)  $ (since (\ref{pf.lem.Vshaped.crit.min}) (applied to
$s=a$) yields $\chi\left(  k\right)  \leq\chi\left(  a\right)  $). This is
absurd. This contradiction shows that our assumption (that $a\notin J$) was
false. Hence, we have $a\in J$.

The same argument (but using $b$ instead of $a$) shows that $b\in J$. Now,
recall that the map $\chi\mid_{J}$ is strictly increasing. Hence, from $a<b$,
we obtain $\left(  \chi\mid_{J}\right)  \left(  a\right)  <\left(  \chi
\mid_{J}\right)  \left(  b\right)  $ (since $a\in J$ and $b\in J$). Thus,
$\chi\left(  a\right)  =\left(  \chi\mid_{J}\right)  \left(  a\right)
<\left(  \chi\mid_{J}\right)  \left(  b\right)  =\chi\left(  b\right)  $. This
contradicts $\chi\left(  a\right)  >\chi\left(  b\right)  $. This
contradiction shows that our assumption was wrong. Hence, $a<k$ is proven.
Thus, (\ref{pf.lem.Vshaped.crit.ak}) is proven.]

Now, we claim that the map $\chi\mid_{\left\{  s\in S\ \mid\ s\leq k\right\}
}$ is strictly decreasing.

[\textit{Proof:} Let $a$ and $b$ be two elements of $\left\{  s\in
S\ \mid\ s\leq k\right\}  $ satisfying $a<b$. From $a\in\left\{  s\in
S\ \mid\ s\leq k\right\}  $, it follows that $a\in S$ and $a\leq k$.
Similarly, $b\in S$ and $b\leq k$.

We shall prove that $\chi\left(  a\right)  >\chi\left(  b\right)  $.

Indeed, assume the contrary. Thus, $\chi\left(  a\right)  \leq\chi\left(
b\right)  $. But $a<b$, so that $a\neq b$ and thus $\chi\left(  a\right)
\neq\chi\left(  b\right)  $ (since the map $\chi$ is injective). Combining
this with $\chi\left(  a\right)  \leq\chi\left(  b\right)  $, we obtain
$\chi\left(  a\right)  <\chi\left(  b\right)  $. Hence,
(\ref{pf.lem.Vshaped.crit.bk}) yields $b>k$. This contradicts $b\leq k$. This
contradiction shows that our assumption was false. Hence, $\chi\left(
a\right)  >\chi\left(  b\right)  $ is proven.

Now, forget that we fixed $a$ and $b$. We thus have shown that if $a$ and $b$
are two elements of $\left\{  s\in S\ \mid\ s\leq k\right\}  $ satisfying
$a<b$, then $\chi\left(  a\right)  >\chi\left(  b\right)  $. In other words,
the map $\chi\mid_{\left\{  s\in S\ \mid\ s\leq k\right\}  }$ is strictly decreasing.]

Next, we claim that the map $\chi\mid_{\left\{  s\in S\ \mid\ s\geq k\right\}
}$ is strictly increasing.

[\textit{Proof:} Let $a$ and $b$ be two elements of $\left\{  s\in
S\ \mid\ s\geq k\right\}  $ satisfying $a<b$. From $a\in\left\{  s\in
S\ \mid\ s\geq k\right\}  $, it follows that $a\in S$ and $a\geq k$.
Similarly, $b\in S$ and $b\geq k$.

We shall prove that $\chi\left(  a\right)  <\chi\left(  b\right)  $.

Indeed, assume the contrary. Thus, $\chi\left(  a\right)  \geq\chi\left(
b\right)  $. But $a<b$, so that $a\neq b$ and thus $\chi\left(  a\right)
\neq\chi\left(  b\right)  $ (since the map $\chi$ is injective). Combining
this with $\chi\left(  a\right)  \geq\chi\left(  b\right)  $, we obtain
$\chi\left(  a\right)  >\chi\left(  b\right)  $. Hence,
(\ref{pf.lem.Vshaped.crit.ak}) yields $a<k$. This contradicts $a\geq k$. This
contradiction shows that our assumption was false. Hence, $\chi\left(
a\right)  <\chi\left(  b\right)  $ is proven.

Now, forget that we fixed $a$ and $b$. We thus have shown that if $a$ and $b$
are two elements of $\left\{  s\in S\ \mid\ s\geq k\right\}  $ satisfying
$a<b$, then $\chi\left(  a\right)  <\chi\left(  b\right)  $. In other words,
the map $\chi\mid_{\left\{  s\in S\ \mid\ s\geq k\right\}  }$ is strictly increasing.]

Thus, we now know that the map $\chi\mid_{\left\{  s\in S\ \mid\ s\leq
k\right\}  }$ is strictly decreasing while the map $\chi\mid_{\left\{  s\in
S\ \mid\ s\geq k\right\}  }$ is strictly increasing. Hence, there exists some
$t\in S$ such that the map $\chi\mid_{\left\{  s\in S\ \mid\ s\leq t\right\}
}$ is strictly decreasing while the map $\chi\mid_{\left\{  s\in
S\ \mid\ s\geq t\right\}  }$ is strictly increasing (namely, $t=k$). In other
words, the map $\chi$ is V-shaped (by the definition of \textquotedblleft
V-shaped\textquotedblright). This proves Lemma \ref{lem.Vshaped.crit}.
\end{proof}
\end{verlong}

\begin{definition}
Let $n\in\mathbb{N}$.

\begin{enumerate}
\item[\textbf{(a)}] Let $f:\left[  n\right]  \rightarrow\mathcal{Z}$ be any
map. Then, $\left\vert f\right\vert $ shall denote the map $\left[  n\right]
\rightarrow\mathcal{N},\ i\mapsto\left\vert f\left(  i\right)  \right\vert $.

\item[\textbf{(b)}] Let $g:\left[  n\right]  \rightarrow\mathcal{N}$ be any
map. Then, we define a monomial $\mathbf{x}_{g}$ in $\operatorname*{Pow}%
\mathcal{N}$ by $\mathbf{x}_{g}=\prod_{i=1}^{n}x_{g\left(  i\right)  }$.
\end{enumerate}
\end{definition}

Using this definition, we can rewrite the definition of $\Gamma_{\mathcal{Z}%
}\left(  \pi\right)  $ as follows:

\begin{proposition}
\label{prop.Gamma.rewr}Let $n\in\mathbb{N}$. Let $\pi$ be any $n$-permutation.
Then,%
\begin{equation}
\Gamma_{\mathcal{Z}}\left(  \pi\right)  =\sum_{f\in\mathcal{E}\left(  \left[
n\right]  ,\pi\right)  }\prod_{p\in\left[  n\right]  }x_{\left\vert f\left(
p\right)  \right\vert }=\sum_{f\in\mathcal{E}\left(  \left[  n\right]
,\pi\right)  }\mathbf{x}_{\left\vert f\right\vert }.\label{eq.Gamma.rewr}%
\end{equation}

\end{proposition}

\begin{vershort}
\begin{proof}
[Proof of Proposition \ref{prop.Gamma.rewr}.] Easy consequence of the
definitions (see \cite{verlong} for details).
\end{proof}
\end{vershort}

\begin{verlong}
\begin{proof}
[Proof of Proposition \ref{prop.Gamma.rewr}.] The definition of $\Gamma
_{\mathcal{Z}}\left(  \pi\right)  $ yields%
\[
\Gamma_{\mathcal{Z}}\left(  \pi\right)  =\Gamma_{\mathcal{Z}}\left(  \left[
n\right]  ,\pi\right)  =\sum_{f\in\mathcal{E}\left(  \left[  n\right]
,\pi\right)  }\prod_{p\in\left[  n\right]  }x_{\left\vert f\left(  p\right)
\right\vert }\ \ \ \ \ \ \ \ \ \ \left(  \text{by the definition of }%
\Gamma_{\mathcal{Z}}\left(  \left[  n\right]  ,\pi\right)  \right)  .
\]
Also,%
\begin{align*}
\sum_{f\in\mathcal{E}\left(  \left[  n\right]  ,\pi\right)  }%
\underbrace{\mathbf{x}_{\left\vert f\right\vert }}_{\substack{=\prod_{i=1}%
^{n}x_{\left\vert f\right\vert \left(  i\right)  }\\\text{(by the definition
of }\mathbf{x}_{\left\vert f\right\vert }\text{)}}}  & =\sum_{f\in
\mathcal{E}\left(  \left[  n\right]  ,\pi\right)  }\underbrace{\prod_{i=1}%
^{n}}_{=\prod_{i\in\left[  n\right]  }}\underbrace{x_{\left\vert f\right\vert
\left(  i\right)  }}_{\substack{=x_{\left\vert f\left(  i\right)  \right\vert
}\\\text{(since }\left\vert f\right\vert \left(  i\right)  =\left\vert
f\left(  i\right)  \right\vert \text{)}}}=\sum_{f\in\mathcal{E}\left(  \left[
n\right]  ,\pi\right)  }\prod_{i\in\left[  n\right]  }x_{\left\vert f\left(
i\right)  \right\vert }\\
& =\sum_{f\in\mathcal{E}\left(  \left[  n\right]  ,\pi\right)  }\prod
_{p\in\left[  n\right]  }x_{\left\vert f\left(  p\right)  \right\vert }%
\end{align*}
(here, we have renamed the index $i$ as $p$ in the product). Combining these
two equalities, we obtain (\ref{eq.Gamma.rewr}). Thus, Proposition
\ref{prop.Gamma.rewr} is proven.
\end{proof}
\end{verlong}

\begin{definition}
\label{def.amen}
Let $n\in\mathbb{N}$. Let $g:\left[  n\right]  \rightarrow\mathcal{N}$ be any
map. Let $\pi$ be an $n$-permutation. We shall say that $g$ is $\pi
$\textit{-amenable} if it has the following properties:

\begin{enumerate}
\item[\textbf{(i')}] The map $\pi\mid_{g^{-1}\left(  0\right)  }$ is strictly
increasing. (This allows the case when $g^{-1}\left(  0\right)  =\varnothing$.)

\item[\textbf{(ii')}] For each $h\in g\left(  \left[  n\right]  \right)
\cap\left\{  1,2,3,\ldots\right\}  $, the map $\pi\mid_{g^{-1}\left(
h\right)  }$ is V-shaped.

\item[\textbf{(iii')}] The map $\pi\mid_{g^{-1}\left(  \infty\right)  }$ is
strictly decreasing. (This allows the case when $g^{-1}\left(  \infty\right)
=\varnothing$.)

\item[\textbf{(iv')}] The map $g$ is weakly increasing.
\end{enumerate}
\end{definition}

\begin{verlong}
\begin{proposition}
\label{prop.amen.equiv2}Let $n\in\mathbb{N}$. Let $\pi$ be any $n$%
-permutation. Let $f\in\mathcal{E}\left(  \left[  n\right]  ,\pi\right)  $.
Then, the map $\left\vert f\right\vert :\left[  n\right]  \rightarrow
\mathcal{N}$ is $\pi$-amenable.
\end{proposition}

\begin{proof}
[Proof of Proposition \ref{prop.amen.equiv2}.]We have $f\in\mathcal{E}\left(
\left[  n\right]  ,\pi\right)  $. Thus, $f$ is a $\mathcal{Z}$-enriched
$\left(  \left[  n\right]  ,\pi\right)  $-partition. In other words, $f$ is a
map $\left[  n\right]  \rightarrow\mathcal{Z}$ such that for all $x<y$ in
$\left[  n\right]  $, the following conditions hold:

\begin{enumerate}
\item[\textbf{(i)}] We have $f\left(  x\right)  \preccurlyeq f\left(
y\right)  $.

\item[\textbf{(ii)}] If $f\left(  x\right)  =f\left(  y\right)  =+h$ for some
$h\in\mathcal{N}$, then $\pi\left(  x\right)  <\pi\left(  y\right)  $.

\item[\textbf{(iii)}] If $f\left(  x\right)  =f\left(  y\right)  =-h$ for some
$h\in\mathcal{N}$, then $\pi\left(  x\right)  >\pi\left(  y\right)  $.
\end{enumerate}

(This follows from the definition of a $\mathcal{Z}$-enriched $\left(  \left[
n\right]  ,\pi\right)  $-partition.)

Condition \textbf{(i)} shows that the map $f$ is weakly increasing. Condition
\textbf{(ii)} shows that for each $h\in\mathcal{N}$, the map $\pi\mid
_{f^{-1}\left(  +h\right)  }$ is strictly increasing. Condition \textbf{(iii)}
shows that for each $h\in\mathcal{N}$, the map $\pi\mid_{f^{-1}\left(
-h\right)  }$ is strictly decreasing.

Now, set $g=\left\vert f\right\vert $. Thus, $g$ is a map $\left[  n\right]
\rightarrow\mathcal{N}$. We shall show that the map $g$ is $\pi$-amenable. In
order to prove this, we must show that the Properties \textbf{(i')},
\textbf{(ii')}, \textbf{(iii')} and \textbf{(iv')} in Definition
\ref{def.amen} hold. We shall now do exactly this.

If $x$ and $y$ are two elements of $\left[  n\right]  $ satisfying $x<y$, then
$g\left(  x\right)  \preccurlyeq g\left(  y\right)  $%
\ \ \ \ \footnote{\textit{Proof.} Let $x$ and $y$ be two elements of $\left[
n\right]  $ satisfying $x<y$. We must prove that $g\left(  x\right)
\preccurlyeq g\left(  y\right)  $.
\par
We have $f\left(  x\right)  \preccurlyeq f\left(  y\right)  $ (by Condition
\textbf{(i)}). Thus, $\left\vert f\left(  x\right)  \right\vert \preccurlyeq
\left\vert f\left(  y\right)  \right\vert $ (by the definition of the order on
$\mathcal{Z}$). But $g=\left\vert f\right\vert $. Hence, $g\left(  x\right)
=\left\vert f\right\vert \left(  x\right)  =\left\vert f\left(  x\right)
\right\vert $ (by the definition of $\left\vert f\right\vert $) and similarly
$g\left(  y\right)  =\left\vert f\left(  y\right)  \right\vert $. Thus,
$g\left(  x\right)  =\left\vert f\left(  x\right)  \right\vert \preccurlyeq
\left\vert f\left(  y\right)  \right\vert =g\left(  y\right)  $, qed.}. Thus,
the map $g$ is weakly increasing. In other words, Property \textbf{(iv')} holds.

Recall that $-0\notin\mathcal{Z}$. Thus, $g^{-1}\left(  0\right)
=f^{-1}\left(  +0\right)  $\ \ \ \ \footnote{\textit{Proof.} For each
$x\in\left[  n\right]  $, we have the following chain of logical equivalences:%
\begin{align*}
\left(  x\in g^{-1}\left(  0\right)  \right)  \  &  \Longleftrightarrow
\ \left(  g\left(  x\right)  =0\right)  \ \Longleftrightarrow\ \left(
\left\vert f\left(  x\right)  \right\vert =0\right)
\ \ \ \ \ \ \ \ \ \ \left(  \text{since }\underbrace{g}_{=\left\vert
f\right\vert }\left(  x\right)  =\left\vert f\right\vert \left(  x\right)
=\left\vert f\left(  x\right)  \right\vert \right) \\
&  \Longleftrightarrow\ \left(  f\left(  x\right)  =+0\text{ or }f\left(
x\right)  =-0\right) \\
&  \Longleftrightarrow\ \left(  f\left(  x\right)  =+0\right)
\ \ \ \ \ \ \ \ \ \ \left(
\begin{array}
[c]{c}%
\text{since }f\left(  x\right)  =-0\text{ cannot hold}\\
\text{(because }f\left(  x\right)  \in\mathcal{Z}\text{ but }-0\notin%
\mathcal{Z}\text{)}%
\end{array}
\right) \\
&  \Longleftrightarrow\ \left(  x\in f^{-1}\left(  +0\right)  \right)  .
\end{align*}
Thus, $g^{-1}\left(  0\right)  =f^{-1}\left(  +0\right)  $.}. But the map
$\pi\mid_{f^{-1}\left(  +0\right)  }$ is strictly increasing\footnote{because
for each $h\in\mathcal{N}$, the map $\pi\mid_{f^{-1}\left(  +h\right)  }$ is
strictly increasing}. In other words, the map $\pi\mid_{g^{-1}\left(
0\right)  }$ is strictly increasing (since $g^{-1}\left(  0\right)
=f^{-1}\left(  +0\right)  $). Hence, Property \textbf{(i')} holds. Similarly,
Property \textbf{(iii')} also holds\footnote{To prove this, we first need to
show that $g^{-1}\left(  \infty\right)  =f^{-1}\left(  -\infty\right)  $ (this
is similar to the proof of $g^{-1}\left(  0\right)  =f^{-1}\left(  +0\right)
$), and that the map $\pi\mid_{f^{-1}\left(  -\infty\right)  }$ is strictly
decreasing (because for each $h\in\mathcal{N}$, the map $\pi\mid
_{f^{-1}\left(  -h\right)  }$ is strictly decreasing).}.

Next, we shall show that Property \textbf{(ii')} holds. Indeed, fix $h\in
g\left(  \left[  n\right]  \right)  \cap\left\{  1,2,3,\ldots\right\}  $.
Thus, $h\in g\left(  \left[  n\right]  \right)  \cap\left\{  1,2,3,\ldots
\right\}  \subseteq g\left(  \left[  n\right]  \right)  $. Hence, there exists
some $y\in\left[  n\right]  $ such that $h=g\left(  y\right)  $. In other
words, the set $g^{-1}\left(  h\right)  $ is nonempty.

This set $g^{-1}\left(  h\right)  $ is an interval of $\left[  n\right]  $
(since the map $g$ is weakly increasing). The sets $f^{-1}\left(  -h\right)  $
and $f^{-1}\left(  +h\right)  $ are intervals of $\left[  n\right]  $ (since
the map $f$ is weakly increasing), and their union is
\[
f^{-1}\left(  -h\right)  \cup f^{-1}\left(  +h\right)  =g^{-1}\left(
h\right)
\]
\ \footnote{\textit{Proof.} For each $x\in\left[  n\right]  $, we have the
following chain of logical equivalences:%
\begin{align*}
\left(  x\in g^{-1}\left(  h\right)  \right)  \  &  \Longleftrightarrow
\ \left(  g\left(  x\right)  =h\right)  \ \Longleftrightarrow\ \left(
\left\vert f\left(  x\right)  \right\vert =h\right)
\ \ \ \ \ \ \ \ \ \ \left(  \text{since }\underbrace{g}_{=\left\vert
f\right\vert }\left(  x\right)  =\left\vert f\right\vert \left(  x\right)
=\left\vert f\left(  x\right)  \right\vert \right) \\
&  \Longleftrightarrow\ \left(  \underbrace{f\left(  x\right)  =+h}%
_{\Longleftrightarrow\ \left(  x\in f^{-1}\left(  +h\right)  \right)  }\text{
or }\underbrace{f\left(  x\right)  =-h}_{\Longleftrightarrow\ \left(  x\in
f^{-1}\left(  -h\right)  \right)  }\right) \\
&  \Longleftrightarrow\ \left(  x\in f^{-1}\left(  +h\right)  \text{ or }x\in
f^{-1}\left(  -h\right)  \right) \\
&  \Longleftrightarrow\ \left(  x\in f^{-1}\left(  +h\right)  \cup
f^{-1}\left(  -h\right)  \right)  .
\end{align*}
Thus, $g^{-1}\left(  h\right)  =f^{-1}\left(  +h\right)  \cup f^{-1}\left(
-h\right)  =f^{-1}\left(  -h\right)  \cup f^{-1}\left(  +h\right)  $, qed.}.
Moreover, each element of $f^{-1}\left(  -h\right)  $ is smaller than each
element of $f^{-1}\left(  +h\right)  $ (since $f$ is weakly increasing, and
$-h\prec+h$). Furthermore, recall that the map $\pi\mid_{f^{-1}\left(
-h\right)  }$ is strictly decreasing while the map $\pi\mid_{f^{-1}\left(
+h\right)  }$ is strictly increasing. In other words, the map $\left(  \pi
\mid_{g^{-1}\left(  h\right)  }\right)  \mid_{f^{-1}\left(  -h\right)  }$ is
strictly decreasing while the map $\left(  \pi\mid_{g^{-1}\left(  h\right)
}\right)  \mid_{f^{-1}\left(  +h\right)  }$ is strictly increasing (since
$\pi\mid_{f^{-1}\left(  -h\right)  }=\left(  \pi\mid_{g^{-1}\left(  h\right)
}\right)  \mid_{f^{-1}\left(  -h\right)  }$ and $\pi\mid_{f^{-1}\left(
+h\right)  }=\left(  \pi\mid_{g^{-1}\left(  h\right)  }\right)  \mid
_{f^{-1}\left(  +h\right)  }$). Hence, Lemma \ref{lem.Vshaped.crit} (applied
to $S=g^{-1}\left(  h\right)  $, $K=\mathbb{P}$, $\chi=\pi\mid_{g^{-1}\left(
h\right)  }$, $I=f^{-1}\left(  -h\right)  $ and $J=f^{-1}\left(  +h\right)  $)
shows that the map $\pi\mid_{g^{-1}\left(  h\right)  }$ is V-shaped. Thus,
Property \textbf{(ii')} holds. We have thus verified all four Properties
\textbf{(i')}, \textbf{(ii')}, \textbf{(iii')} and \textbf{(iv')}; thus, $g$
is $\pi$-amenable. In other words, $\left\vert f\right\vert $ is $\pi
$-amenable (since $g=\left\vert f\right\vert $). This proves Proposition
\ref{prop.amen.equiv2}.
\end{proof}

\begin{proposition}
\label{prop.amen.equiv3}Let $n\in\mathbb{N}$. Let $\pi$ be any $n$%
-permutation. Let $g:\left[  n\right]  \rightarrow\mathcal{N}$ be a $\pi
$-amenable map. Let $H$ be the set $g\left(  \left[  n\right]  \right)
\cap\left\{  1,2,3,\ldots\right\}  $. For each $h\in H$, we let $m_{h}$ be the
unique element $\mu$ of $g^{-1}\left(  h\right)  $ for which $\pi\left(
\mu\right)  $ is minimum.

\textbf{(a)} The elements $m_{h}$ for all $h\in H$ are well-defined and distinct.

\textbf{(b)} Let $f:\left[  n\right]  \rightarrow\mathcal{Z}$ be any map.
Then, $\left(  f\in\mathcal{E}\left(  \left[  n\right]  ,\pi\right)  \text{
and }\left\vert f\right\vert =g\right)  $ if and only if the following five
statements hold:

\begin{enumerate}
\item[\textbf{(x}$_{1}$\textbf{)}] For each $x\in g^{-1}\left(  0\right)  $,
we have $f\left(  x\right)  =+0$.

\item[\textbf{(x}$_{2}$\textbf{)}] For each $h\in H$ and each $x\in
g^{-1}\left(  h\right)  $ satisfying $x<m_{h}$, we have $f\left(  x\right)
=-h$.

\item[\textbf{(x}$_{3}$\textbf{)}] For each $h\in H$, we have $f\left(
m_{h}\right)  \in\left\{  -h,+h\right\}  $.

\item[\textbf{(x}$_{4}$\textbf{)}] For each $h\in H$ and each $x\in
g^{-1}\left(  h\right)  $ satisfying $x>m_{h}$, we have $f\left(  x\right)
=+h$.

\item[\textbf{(x}$_{5}$\textbf{)}] For each $x\in g^{-1}\left(  \infty\right)
$, we have $f\left(  x\right)  =-\infty$.
\end{enumerate}
\end{proposition}

\begin{proof}
[Proof of Proposition \ref{prop.amen.equiv3}.]The map $\pi:\left[  n\right]
\rightarrow\mathbb{P}$ is injective (since $\pi$ is an $n$-permutation). The
map $g$ is $\pi$-amenable, and thus satisfies the four Properties
\textbf{(i')}, \textbf{(ii')}, \textbf{(iii')} and \textbf{(iv')} from
Definition \ref{def.amen}. In particular, this map $g$ is weakly increasing
(because of Property \textbf{(iv')}).

\textbf{(a)} The elements $m_{h}$ for all $h\in H$ are
well-defined\footnote{\textit{Proof.} Let $h\in H$. We must prove that the
element $m_{h}$ is well-defined.
\par
We have $h\in H=g\left(  \left[  n\right]  \right)  \cap\left\{
1,2,3,\ldots\right\}  \subseteq g\left(  \left[  n\right]  \right)  $. Hence,
there exists some $y\in\left[  n\right]  $ such that $h=g\left(  y\right)  $.
Therefore, the preimage $g^{-1}\left(  h\right)  $ is nonempty. Therefore,
there exists an element $\mu$ of $g^{-1}\left(  h\right)  $ for which
$\pi\left(  \mu\right)  $ is minimum. Moreover, this element is unique (since
$\pi$ is injective). Thus, the unique element $\mu$ of $g^{-1}\left(
h\right)  $ for which $\pi\left(  \mu\right)  $ is minimum is well-defined. In
other words, $m_{h}$ is well-defined. Qed.} and
distinct\footnote{\textit{Proof.} Let $h_{1}\in H$ and $h_{2}\in H$ be such
that $m_{h_{1}}=m_{h_{2}}$. We must prove that $h_{1}=h_{2}$.
\par
We have $m_{h}\in g^{-1}\left(  h\right)  $ for each $h\in H$ (by the
definition of $m_{h}$). Hence, $m_{h_{1}}\in g^{-1}\left(  h_{1}\right)  $. In
other words, $g\left(  m_{h_{1}}\right)  =h_{1}$. Similarly, $g\left(
m_{h_{2}}\right)  =h_{2}$. Hence, $h_{1}=g\left(  \underbrace{m_{h_{1}}%
}_{=m_{h_{2}}}\right)  =g\left(  m_{h_{2}}\right)  =h_{2}$. Qed.}. This proves
Proposition \ref{prop.amen.equiv3} \textbf{(a)}.

\textbf{(b)} We must prove the following two claims:

\begin{statement}
\textit{Claim 1:} If $\left(  f\in\mathcal{E}\left(  \left[  n\right]
,\pi\right)  \text{ and }\left\vert f\right\vert =g\right)  $, then the five
statements \textbf{(x}$_{1}$\textbf{)}, \textbf{(x}$_{2}$\textbf{)},
\textbf{(x}$_{3}$\textbf{)}, \textbf{(x}$_{4}$\textbf{)} and \textbf{(x}$_{5}%
$\textbf{)} hold.
\end{statement}

\begin{statement}
\textit{Claim 2:} If the five statements \textbf{(x}$_{1}$\textbf{)},
\textbf{(x}$_{2}$\textbf{)}, \textbf{(x}$_{3}$\textbf{)}, \textbf{(x}$_{4}%
$\textbf{)} and \textbf{(x}$_{5}$\textbf{)} hold, then $\left(  f\in
\mathcal{E}\left(  \left[  n\right]  ,\pi\right)  \text{ and }\left\vert
f\right\vert =g\right)  $.
\end{statement}

[\textit{Proof of Claim 1:} Assume that $\left(  f\in\mathcal{E}\left(
\left[  n\right]  ,\pi\right)  \text{ and }\left\vert f\right\vert =g\right)
$. We must prove that the five statements \textbf{(x}$_{1}$\textbf{)},
\textbf{(x}$_{2}$\textbf{)}, \textbf{(x}$_{3}$\textbf{)}, \textbf{(x}$_{4}%
$\textbf{)} and \textbf{(x}$_{5}$\textbf{)} hold.

We have $f\in\mathcal{E}\left(  \left[  n\right]  ,\pi\right)  $. We have
$g=\left\vert f\right\vert $. Hence, everything that we said in the proof of
Proposition \ref{prop.amen.equiv2} is true in our current situation as well.
In particular, the Properties \textbf{(i')}, \textbf{(ii')}, \textbf{(iii')}
and \textbf{(iv')} in Definition \ref{def.amen} hold.

Now, statement \textbf{(x}$_{1}$\textbf{)} holds\footnote{\textit{Proof.} Let
$x\in g^{-1}\left(  0\right)  $. But $g^{-1}\left(  0\right)  =f^{-1}\left(
+0\right)  $ (as we have shown in the proof of Proposition
\ref{prop.amen.equiv2}). Hence, $x\in g^{-1}\left(  0\right)  =f^{-1}\left(
+0\right)  $, so that $f\left(  x\right)  =+0$.
\par
Now, forget that we fixed $x$. Thus, we have shown that for each $x\in
g^{-1}\left(  0\right)  $, we have $f\left(  x\right)  =+0$. In other words,
statement \textbf{(x}$_{1}$\textbf{)} holds.}, and statement \textbf{(x}$_{5}%
$\textbf{)} holds\footnote{\textit{Proof.} Let $x\in g^{-1}\left(
\infty\right)  $. But $g^{-1}\left(  \infty\right)  =f^{-1}\left(
-\infty\right)  $ (as we have shown in the proof of Proposition
\ref{prop.amen.equiv2}). Hence, $x\in g^{-1}\left(  \infty\right)
=f^{-1}\left(  -\infty\right)  $, so that $f\left(  x\right)  =-\infty$.
\par
Now, forget that we fixed $x$. Thus, we have shown that for each $x\in
g^{-1}\left(  \infty\right)  $, we have $f\left(  x\right)  =-\infty$. In
other words, statement \textbf{(x}$_{5}$\textbf{)} holds.}. Furthermore,
statement \textbf{(x}$_{3}$\textbf{)} holds\footnote{\textit{Proof.} Let $h\in
H$. Then, $m_{h}\in g^{-1}\left(  h\right)  $ (by the definition of $m_{h}$),
so that $g\left(  m_{h}\right)  =h$. Comparing this with $\underbrace{g}%
_{=\left\vert f\right\vert }\left(  m_{h}\right)  =\left\vert f\right\vert
\left(  m_{h}\right)  =\left\vert f\left(  m_{h}\right)  \right\vert $ (by the
definition of $\left\vert f\right\vert $), we obtain $\left\vert f\left(
m_{h}\right)  \right\vert =h$. In other words, either $f\left(  m_{h}\right)
=-h$ or $f\left(  m_{h}\right)  =+h$. In other words, $f\left(  m_{h}\right)
\in\left\{  -h,+h\right\}  $.
\par
Now, forget that we fixed $h$. Thus, we have shown that for each $h\in H$, we
have $f\left(  m_{h}\right)  \in\left\{  -h,+h\right\}  $. In other words,
statement \textbf{(x}$_{3}$\textbf{)} holds.}. Also, statement \textbf{(x}%
$_{2}$\textbf{)} holds\footnote{\textit{Proof.} Let $h\in H$ and $x\in
g^{-1}\left(  h\right)  $ be such that $x<m_{h}$. We shall show that $f\left(
x\right)  =-h$.
\par
Indeed, assume the contrary. Thus, $f\left(  x\right)  \neq-h$.
\par
We have $x\in g^{-1}\left(  h\right)  $, so that $g\left(  x\right)  =h$.
Comparing this with $\underbrace{g}_{=\left\vert f\right\vert }\left(
x\right)  =\left\vert f\right\vert \left(  x\right)  =\left\vert f\left(
x\right)  \right\vert $ (by the definition of $\left\vert f\right\vert $), we
obtain $\left\vert f\left(  x\right)  \right\vert =h$. In other words, either
$f\left(  x\right)  =-h$ or $f\left(  x\right)  =+h$. Since $f\left(
x\right)  \neq-h$, we thus must have $f\left(  x\right)  =+h$. But the map
$f:\left[  n\right]  \rightarrow\mathcal{Z}$ is weakly increasing (as we have
seen in the proof of Proposition \ref{prop.amen.equiv2}). Hence, from
$x<m_{h}$, we obtain $f\left(  x\right)  \preccurlyeq f\left(  m_{h}\right)
$. Hence, we cannot have $f\left(  m_{h}\right)  =-h$ (because if we had
$f\left(  m_{h}\right)  =-h$, then we would have $+h=f\left(  x\right)
\preccurlyeq f\left(  m_{h}\right)  =-h$, which would contradict the
definition of the order on $\mathcal{Z}$).
\par
But statement \textbf{(x}$_{3}$\textbf{)} yields $f\left(  m_{h}\right)
\in\left\{  -h,+h\right\}  $. In other words, either $f\left(  m_{h}\right)
=-h$ or $f\left(  m_{h}\right)  =+h$. Hence, $f\left(  m_{h}\right)  =+h$
(since we cannot have $f\left(  m_{h}\right)  =-h$). Combining $f\left(
x\right)  =+h$ and $f\left(  m_{h}\right)  =+h$, we obtain $f\left(  x\right)
=f\left(  m_{h}\right)  =+h$.
\par
But Condition \textbf{(ii)} from the proof of Proposition
\ref{prop.amen.equiv2} holds. Applying this condition to $y=m_{h}$, we obtain
$\pi\left(  x\right)  <\pi\left(  m_{h}\right)  $.
\par
But $m_{h}$ is the unique element $\mu$ of $g^{-1}\left(  h\right)  $ for
which $\pi\left(  \mu\right)  $ is minimum (by the definition of $m_{h}$).
Hence, $\pi\left(  m_{h}\right)  \leq\pi\left(  z\right)  $ for each $z\in
g^{-1}\left(  h\right)  $. Applying this to $z=x$, we obtain $\pi\left(
m_{h}\right)  \leq\pi\left(  x\right)  $. This contradicts $\pi\left(
x\right)  <\pi\left(  m_{h}\right)  $. This contradiction shows that our
assumption was false. Hence, $f\left(  x\right)  =-h$ is proven.
\par
Now, forget that we fixed $h$ and $x$. Thus, we have shown that for each $h\in
H$ and each $x\in g^{-1}\left(  h\right)  $ satisfying $x<m_{h}$, we have
$f\left(  x\right)  =-h$. In other words, statement \textbf{(x}$_{2}%
$\textbf{)} holds.}, and statement \textbf{(x}$_{4}$\textbf{)}
holds\footnote{\textit{Proof.} Let $h\in H$ and $x\in g^{-1}\left(  h\right)
$ be such that $x>m_{h}$. Thus, $m_{h}<x$. We shall show that $f\left(
x\right)  =+h$.
\par
Indeed, assume the contrary. Thus, $f\left(  x\right)  \neq+h$.
\par
We have $x\in g^{-1}\left(  h\right)  $, so that $g\left(  x\right)  =h$.
Comparing this with $\underbrace{g}_{=\left\vert f\right\vert }\left(
x\right)  =\left\vert f\right\vert \left(  x\right)  =\left\vert f\left(
x\right)  \right\vert $ (by the definition of $\left\vert f\right\vert $), we
obtain $\left\vert f\left(  x\right)  \right\vert =h$. In other words, either
$f\left(  x\right)  =-h$ or $f\left(  x\right)  =+h$. Since $f\left(
x\right)  \neq+h$, we thus must have $f\left(  x\right)  =-h$. But the map
$f:\left[  n\right]  \rightarrow\mathcal{Z}$ is weakly increasing (as we have
seen in the proof of Proposition \ref{prop.amen.equiv2}). Hence, from
$m_{h}<x$, we obtain $f\left(  m_{h}\right)  \preccurlyeq f\left(  x\right)
$. Hence, we cannot have $f\left(  m_{h}\right)  =+h$ (because if we had
$f\left(  m_{h}\right)  =+h$, then we would have $+h=f\left(  m_{h}\right)
\preccurlyeq f\left(  x\right)  =-h$, which would contradict the definition of
the order on $\mathcal{Z}$).
\par
But statement \textbf{(x}$_{3}$\textbf{)} yields $f\left(  m_{h}\right)
\in\left\{  -h,+h\right\}  $. In other words, either $f\left(  m_{h}\right)
=-h$ or $f\left(  m_{h}\right)  =+h$. Hence, $f\left(  m_{h}\right)  =-h$
(since we cannot have $f\left(  m_{h}\right)  =+h$). Combining $f\left(
x\right)  =-h$ and $f\left(  m_{h}\right)  =-h$, we obtain $f\left(
m_{h}\right)  =f\left(  x\right)  =-h$.
\par
But Condition \textbf{(iii)} from the proof of Proposition
\ref{prop.amen.equiv2} holds. Applying this condition to $m_{h}$ and $x$
instead of $x$ and $y$, we obtain $\pi\left(  m_{h}\right)  >\pi\left(
x\right)  $.
\par
But $m_{h}$ is the unique element $\mu$ of $g^{-1}\left(  h\right)  $ for
which $\pi\left(  \mu\right)  $ is minimum (by the definition of $m_{h}$).
Hence, $\pi\left(  m_{h}\right)  \leq\pi\left(  z\right)  $ for each $z\in
g^{-1}\left(  h\right)  $. Applying this to $z=x$, we obtain $\pi\left(
m_{h}\right)  \leq\pi\left(  x\right)  $. This contradicts $\pi\left(
m_{h}\right)  >\pi\left(  x\right)  $. This contradiction shows that our
assumption was false. Hence, $f\left(  x\right)  =+h$ is proven.
\par
Now, forget that we fixed $h$ and $x$. Thus, we have shown that for each $h\in
H$ and each $x\in g^{-1}\left(  h\right)  $ satisfying $x>m_{h}$, we have
$f\left(  x\right)  =+h$. In other words, statement \textbf{(x}$_{4}%
$\textbf{)} holds.}. Thus, the five statements \textbf{(x}$_{1}$\textbf{)},
\textbf{(x}$_{2}$\textbf{)}, \textbf{(x}$_{3}$\textbf{)}, \textbf{(x}$_{4}%
$\textbf{)} and \textbf{(x}$_{5}$\textbf{)} hold. This proves Claim 1.]

[\textit{Proof of Claim 2:} Assume that five statements \textbf{(x}$_{1}%
$\textbf{)}, \textbf{(x}$_{2}$\textbf{)}, \textbf{(x}$_{3}$\textbf{)},
\textbf{(x}$_{4}$\textbf{)} and \textbf{(x}$_{5}$\textbf{)} hold. We must
prove that $\left(  f\in\mathcal{E}\left(  \left[  n\right]  ,\pi\right)
\text{ and }\left\vert f\right\vert =g\right)  $.

For each $x\in\left[  n\right]  $, we have
\begin{equation}
\left\vert f\right\vert \left(  x\right)  =g\left(  x\right)  .
\label{pf.prop.amen.equiv3.c2.pf.1}%
\end{equation}

[\textit{Proof of (\ref{pf.prop.amen.equiv3.c2.pf.1}):} Let $x\in\left[
n\right]  $. We must prove (\ref{pf.prop.amen.equiv3.c2.pf.1}). The definition
of $\left\vert f\right\vert $ yields $\left\vert f\right\vert \left(
x\right)  =\left\vert f\left(  x\right)  \right\vert $.

The equality (\ref{pf.prop.amen.equiv3.c2.pf.1}) holds when $g\left(
x\right)  =0$\ \ \ \ \footnote{\textit{Proof.} Assume that $g\left(  x\right)
=0$. Thus, $x\in g^{-1}\left(  0\right)  $. Hence, statement \textbf{(x}$_{1}%
$\textbf{)} yields $f\left(  x\right)  =+0$ (since we know that statement
\textbf{(x}$_{1}$\textbf{)} holds). Thus, $\left\vert f\left(  x\right)
\right\vert =\left\vert +0\right\vert =0$. Hence, $\left\vert f\right\vert
\left(  x\right)  =\left\vert f\left(  x\right)  \right\vert =0=g\left(
x\right)  $. Thus, the equality (\ref{pf.prop.amen.equiv3.c2.pf.1}) holds,
qed.}. Hence, for the rest of this proof, we WLOG assume that $g\left(
x\right)  \neq0$.

The equality (\ref{pf.prop.amen.equiv3.c2.pf.1}) holds when $g\left(
x\right)  =\infty$\ \ \ \ \footnote{\textit{Proof.} Assume that $g\left(
x\right)  =\infty$. Thus, $x\in g^{-1}\left(  \infty\right)  $. Hence,
statement \textbf{(x}$_{5}$\textbf{)} yields $f\left(  x\right)  =-\infty$
(since we know that statement \textbf{(x}$_{5}$\textbf{)} holds). Thus,
$\left\vert f\left(  x\right)  \right\vert =\left\vert -\infty\right\vert
=\infty$. Hence, $\left\vert f\right\vert \left(  x\right)  =\left\vert
f\left(  x\right)  \right\vert =\infty=g\left(  x\right)  $. Thus, the
equality (\ref{pf.prop.amen.equiv3.c2.pf.1}) holds, qed.}. Hence, for the rest
of this proof, we WLOG assume that $g\left(  x\right)  \neq\infty$. Combining
this with $g\left(  x\right)  \neq0$, we obtain $g\left(  x\right)
\notin\left\{  0,\infty\right\}  $.

Set $h=g\left(  x\right)  $. We have
\begin{align*}
h  &  =g\left(  x\right)  \in\mathcal{N}\setminus\left\{  0,\infty\right\}
\ \ \ \ \ \ \ \ \ \ \left(  \text{since }g\left(  x\right)  \in\mathcal{N}%
\text{ and }g\left(  x\right)  \notin\left\{  0,\infty\right\}  \right) \\
&  =\left\{  1,2,3,\ldots\right\}  .
\end{align*}
Combining this with $h=g\left(  \underbrace{x}_{\in\left[  n\right]  }\right)
\in g\left(  \left[  n\right]  \right)  $, we obtain $h\in g\left(  \left[
n\right]  \right)  \cap\left\{  1,2,3,\ldots\right\}  =H$. Also, $x\in
g^{-1}\left(  h\right)  $ (since $g\left(  x\right)  =h$). Now, we are in one
of the following three cases:

\textit{Case 1:} We have $x<m_{h}$.

\textit{Case 2:} We have $x=m_{h}$.

\textit{Case 3:} We have $x>m_{h}$.

Let us first consider Case 1. In this case, we have $x<m_{h}$. Hence,
statement \textbf{(x}$_{2}$\textbf{)} yields $f\left(  x\right)  =-h$ (since
we know that statement \textbf{(x}$_{2}$\textbf{)} holds). Hence, $\left\vert
f\left(  x\right)  \right\vert =\left\vert -h\right\vert =h$. Hence,
$\left\vert f\right\vert \left(  x\right)  =\left\vert f\left(  x\right)
\right\vert =h=g\left(  x\right)  $. Thus, the equality
(\ref{pf.prop.amen.equiv3.c2.pf.1}) holds in Case 1.

Let us next consider Case 2. In this case, we have $x=m_{h}$. But statement
\textbf{(x}$_{3}$\textbf{)} yields $f\left(  m_{h}\right)  \in\left\{
-h,+h\right\}  $ (since we know that statement \textbf{(x}$_{3}$\textbf{)}
holds). Hence, $\left\vert f\left(  m_{h}\right)  \right\vert =h$. Hence,
$\left\vert f\right\vert \left(  x\right)  =\left\vert f\left(  \underbrace{x}%
_{=m_{h}}\right)  \right\vert =\left\vert f\left(  m_{h}\right)  \right\vert
=h=g\left(  x\right)  $. Thus, the equality (\ref{pf.prop.amen.equiv3.c2.pf.1}%
) holds in Case 2.

Let us finally consider Case 3. In this case, we have $x>m_{h}$. Hence,
statement \textbf{(x}$_{4}$\textbf{)} yields $f\left(  x\right)  =+h$ (since
we know that statement \textbf{(x}$_{4}$\textbf{)} holds). Hence, $\left\vert
f\left(  x\right)  \right\vert =\left\vert +h\right\vert =h$. Hence,
$\left\vert f\right\vert \left(  x\right)  =\left\vert f\left(  x\right)
\right\vert =h=g\left(  x\right)  $. Thus, the equality
(\ref{pf.prop.amen.equiv3.c2.pf.1}) holds in Case 3.

We have now shown that the equality (\ref{pf.prop.amen.equiv3.c2.pf.1}) holds
in each of the three Cases 1, 2 and 3. Thus,
(\ref{pf.prop.amen.equiv3.c2.pf.1}) always holds. This completes the proof of
(\ref{pf.prop.amen.equiv3.c2.pf.1}).]

The equality (\ref{pf.prop.amen.equiv3.c2.pf.1}) shows that $\left\vert
f\right\vert =g$. It remains to prove that $f\in\mathcal{E}\left(  \left[
n\right]  ,\pi\right)  $.

Consider the three Conditions \textbf{(i)}, \textbf{(ii)} and \textbf{(iii)}
from the proof of Proposition \ref{prop.amen.equiv2}. Now, let $x$ and $y$ be
two elements of $\left[  n\right]  $ satisfying $x<y$. We claim that the
following three conditions hold:

\begin{enumerate}
\item[\textbf{(i)}] We have $f\left(  x\right)  \preccurlyeq f\left(
y\right)  $.

\item[\textbf{(ii)}] If $f\left(  x\right)  =f\left(  y\right)  =+h$ for some
$h\in\mathcal{N}$, then $\pi\left(  x\right)  <\pi\left(  y\right)  $.

\item[\textbf{(iii)}] If $f\left(  x\right)  =f\left(  y\right)  =-h$ for some
$h\in\mathcal{N}$, then $\pi\left(  x\right)  >\pi\left(  y\right)  $.
\end{enumerate}

[\textit{Proof of Condition \textbf{(i)}:} Assume the contrary. Thus, we don't
have $f\left(  x\right)  \preccurlyeq f\left(  y\right)  $.

But $\left\vert f\right\vert =g$, so that $\left\vert f\right\vert \left(
x\right)  =g\left(  x\right)  $. But the definition of $\left\vert
f\right\vert $ yields $\left\vert f\right\vert \left(  x\right)  =\left\vert
f\left(  x\right)  \right\vert $. Hence, $\left\vert f\left(  x\right)
\right\vert =\left\vert f\right\vert \left(  x\right)  =g\left(  x\right)  $.
Similarly, $\left\vert f\left(  y\right)  \right\vert =g\left(  y\right)  $.
But the map $g$ is weakly increasing; thus, $g\left(  x\right)  \preccurlyeq
g\left(  y\right)  $ (since $x<y$). Hence, $\left\vert f\left(  x\right)
\right\vert =g\left(  x\right)  \preccurlyeq g\left(  y\right)  =\left\vert
f\left(  y\right)  \right\vert $.

A look back at the definition of the total order on $\mathcal{Z}$ reveals the
following: If two elements $\alpha$ and $\beta$ of $\mathcal{Z}$ satisfy
$\left\vert \alpha\right\vert \preccurlyeq\left\vert \beta\right\vert $ but
not $\alpha\preccurlyeq\beta$, then there must exist some $h\in\left\{
1,2,3,\ldots\right\}  $ such that $\alpha=+h$ and $\beta=-h$. Applying this to
$\alpha=f\left(  x\right)  $ and $\beta=f\left(  y\right)  $, we conclude that
there must exist some $h\in\left\{  1,2,3,\ldots\right\}  $ such that
$f\left(  x\right)  =+h$ and $f\left(  y\right)  =-h$ (because the elements
$f\left(  x\right)  $ and $f\left(  y\right)  $ of $\mathcal{Z}$ satisfy
$\left\vert f\left(  x\right)  \right\vert \preccurlyeq\left\vert f\left(
y\right)  \right\vert $ but not $f\left(  x\right)  \preccurlyeq f\left(
y\right)  $). Consider this $h$. We have $g\left(  x\right)  =\left\vert
\underbrace{f\left(  x\right)  }_{=+h}\right\vert =\left\vert +h\right\vert
=h$, so that $h=g\left(  \underbrace{x}_{\in\left[  n\right]  }\right)  \in
g\left(  \left[  n\right]  \right)  $. Combining this with $h\in\left\{
1,2,3,\ldots\right\}  $, we obtain $h\in g\left(  \left[  n\right]  \right)
\cap\left\{  1,2,3,\ldots\right\}  =H$. Also, $x\in g^{-1}\left(  h\right)  $
(since $g\left(  x\right)  =h$). If we had $x<m_{h}$, then statement
\textbf{(x}$_{2}$\textbf{)} would yield $f\left(  x\right)  =-h$, which would
contradict $f\left(  x\right)  =+h\neq-h$. Hence, we cannot have $x<m_{h}$.
Thus, we have $x\geq m_{h}$. But $x<y$, so that $y>x\geq m_{h}$. Also,
$g\left(  y\right)  =\left\vert \underbrace{f\left(  y\right)  }%
_{=-h}\right\vert =\left\vert -h\right\vert =h$, so that $y\in g^{-1}\left(
h\right)  $. Hence, statement \textbf{(x}$_{4}$\textbf{)} (applied to $y$
instead of $x$) yields $f\left(  y\right)  =+h$. This contradicts $f\left(
y\right)  =-h\neq+h$. This contradiction shows that our assumption was false.
Hence, Condition \textbf{(i)} is proven.]

[\textit{Proof of Condition \textbf{(ii)}:} Let $f\left(  x\right)  =f\left(
y\right)  =+h$ for some $h\in\mathcal{N}$. We must show that $\pi\left(
x\right)  <\pi\left(  y\right)  $.

We have $+h=f\left(  x\right)  \in\mathcal{Z}$. Thus, we cannot have
$h=\infty$ (since $h=\infty$ would lead to $+h=+\infty\notin\mathcal{Z}$,
which would contradict $+h\in\mathcal{Z}$). In other words, $h\neq\infty$.
Combining this with $h\in\mathcal{N}$, we obtain
$h\in\mathcal{N}\setminus\left\{ \infty\right\}  $.

The definition of $\left\vert f\right\vert $ yields $\left\vert f\right\vert
\left(  x\right)  =\left\vert f\left(  x\right)  \right\vert $. But
$g=\left\vert f\right\vert $. Hence, $g\left(  x\right)  =\left\vert
f\right\vert \left(  x\right)  =\left\vert \underbrace{f\left(  x\right)
}_{=+h}\right\vert =\left\vert +h\right\vert =h$, so that $x\in g^{-1}\left(
h\right)  $. Similarly, $y\in g^{-1}\left(  h\right)  $.

If $h=0$, then it is easy to see that $\pi\left(  x\right)  <\pi\left(
y\right)  $\ \ \ \ \footnote{\textit{Proof.} Assume that $h=0$. Now, recall
that $x\in g^{-1}\left(  h\right)  $. In other words, $x\in g^{-1}\left(
0\right)  $ (since $h=0$). Similarly, $y\in g^{-1}\left(  0\right)  $. But the
map $\pi\mid_{g^{-1}\left(  0\right)  }$ is strictly increasing (since
Condition \textbf{(i')} from Definition \ref{def.amen} holds). Hence, $\left(
\pi\mid_{g^{-1}\left(  0\right)  }\right)  \left(  x\right)  <\left(  \pi
\mid_{g^{-1}\left(  0\right)  }\right)  \left(  y\right)  $ (because $x<y$ and
because both $x$ and $y$ belong to $g^{-1}\left(  0\right)  $). Thus,
$\pi\left(  x\right)  =\left(  \pi\mid_{g^{-1}\left(  0\right)  }\right)
\left(  x\right)  <\left(  \pi\mid_{g^{-1}\left(  0\right)  }\right)  \left(
y\right)  =\pi\left(  y\right)  $, qed.}. Hence, for the rest of this proof of
$\pi\left(  x\right)  <\pi\left(  y\right)  $, we WLOG assume that we don't
have $h=0$. Thus, $h\neq0$. Combining this with $h\in\mathcal{N}%
\setminus\left\{  \infty\right\}  $, we obtain $h\in\left(  \mathcal{N}%
\setminus\left\{  \infty\right\}  \right)  \setminus\left\{  0\right\}
=\mathcal{N}\setminus\left\{  0,\infty\right\}  =\left\{  1,2,3,\ldots
\right\}  $.

Combining this with $h=g\left(  \underbrace{x}_{\in\left[  n\right]  }\right)
\in g\left(  \left[  n\right]  \right)  $, we obtain $h\in g\left(  \left[
n\right]  \right)  \cap\left\{  1,2,3,\ldots\right\}  =H$. Also, $x\in
g^{-1}\left(  h\right)  $. Hence, if we had $x<m_{h}$, then statement
\textbf{(x}$_{2}$\textbf{)} would yield $f\left(  x\right)  =-h$, which would
contradict $f\left(  x\right)  =+h\neq-h$. Therefore, we cannot have $x<m_{h}%
$. Hence, we have $x\geq m_{h}$. In other words, $x\in\left\{  s\in
g^{-1}\left(  h\right)  \ \mid\ s\geq m_{h}\right\}  $ (since $x\in
g^{-1}\left(  h\right)  $). Similarly, $y\in\left\{  s\in g^{-1}\left(
h\right)  \ \mid\ s\geq m_{h}\right\}  $.

But $h\in g\left(  \left[  n\right]  \right)  \cap\left\{  1,2,3,\ldots
\right\}  $. Thus, the map $\pi\mid_{g^{-1}\left(  h\right)  }$ is V-shaped
(since Condition \textbf{(ii')} from Definition \ref{def.amen} holds). In
other words, there exists some $t\in g^{-1}\left(  h\right)  $ such that the
map $\left(  \pi\mid_{g^{-1}\left(  h\right)  }\right)  \mid_{\left\{  s\in
g^{-1}\left(  h\right)  \ \mid\ s\leq t\right\}  }$ is strictly decreasing
while the map $\left(  \pi\mid_{g^{-1}\left(  h\right)  }\right)
\mid_{\left\{  s\in g^{-1}\left(  h\right)  \ \mid\ s\geq t\right\}  }$ is
strictly increasing. Clearly, this $t$ must be the unique element $\mu$ of
$g^{-1}\left(  h\right)  $ for which $\left(  \pi\mid_{g^{-1}\left(  h\right)
}\right)  \left(  \mu\right)  $ is minimum. In other words, this $t$ must be
the unique element $\mu$ of $g^{-1}\left(  h\right)  $ for which $\pi\left(
\mu\right)  $ is minimum (since $\left(  \pi\mid_{g^{-1}\left(  h\right)
}\right)  \left(  \mu\right)  =\pi\left(  \mu\right)  $ for each $\mu\in
g^{-1}\left(  h\right)  $). In other words, this $t$ must be $m_{h}$ (because
$m_{h}$ was defined to be the unique element $\mu$ of $g^{-1}\left(  h\right)
$ for which $\pi\left(  \mu\right)  $ is minimum). Thus, we conclude that the
map $\left(  \pi\mid_{g^{-1}\left(  h\right)  }\right)  \mid_{\left\{  s\in
g^{-1}\left(  h\right)  \ \mid\ s\leq m_{h}\right\}  }$ is strictly decreasing
while the map \newline$\left(  \pi\mid_{g^{-1}\left(  h\right)  }\right)
\mid_{\left\{  s\in g^{-1}\left(  h\right)  \ \mid\ s\geq m_{h}\right\}  }$ is
strictly increasing.

In particular, the map $\left(  \pi\mid_{g^{-1}\left(  h\right)  }\right)
\mid_{\left\{  s\in g^{-1}\left(  h\right)  \ \mid\ s\geq m_{h}\right\}  }$ is
strictly increasing. In other words, the map $\pi\mid_{\left\{  s\in
g^{-1}\left(  h\right)  \ \mid\ s\geq m_{h}\right\}  }$ is strictly increasing
(since \newline$\left(  \pi\mid_{g^{-1}\left(  h\right)  }\right)
\mid_{\left\{  s\in g^{-1}\left(  h\right)  \ \mid\ s\geq m_{h}\right\}  }%
=\pi\mid_{\left\{  s\in g^{-1}\left(  h\right)  \ \mid\ s\geq m_{h}\right\}
}$). Since both $x$ and $y$ belong to the set $\left\{  s\in g^{-1}\left(
h\right)  \ \mid\ s\geq m_{h}\right\}  $, we thus conclude that%
\[
\left(  \pi\mid_{\left\{  s\in g^{-1}\left(  h\right)  \ \mid\ s\geq
m_{h}\right\}  }\right)  \left(  x\right)  <\left(  \pi\mid_{\left\{  s\in
g^{-1}\left(  h\right)  \ \mid\ s\geq m_{h}\right\}  }\right)  \left(
y\right)
\]
(since $x<y$). Thus,%
\[
\pi\left(  x\right)  =\left(  \pi\mid_{\left\{  s\in g^{-1}\left(  h\right)
\ \mid\ s\geq m_{h}\right\}  }\right)  \left(  x\right)  <\left(  \pi
\mid_{\left\{  s\in g^{-1}\left(  h\right)  \ \mid\ s\geq m_{h}\right\}
}\right)  \left(  y\right)  =\pi\left(  y\right)  .
\]
This proves Condition \textbf{(ii)}.]

[\textit{Proof of Condition \textbf{(iii)}:} Let $f\left(  x\right)  =f\left(
y\right)  =-h$ for some $h\in\mathcal{N}$. We must show that $\pi\left(
x\right)  >\pi\left(  y\right)  $.

We have $-h=f\left(  x\right)  \in\mathcal{Z}$. Thus, we cannot have $h=0$
(since $h=0$ would lead to $-h=-0\notin\mathcal{Z}$, which would contradict
$-h\in\mathcal{Z}$). In other words, $h\neq0$.
Combining this with $h\in\mathcal{N}$, we obtain
$h\in\mathcal{N}\setminus\left\{  0\right\}  $.

The definition of $\left\vert f\right\vert $ yields $\left\vert f\right\vert
\left(  x\right)  =\left\vert f\left(  x\right)  \right\vert $. But
$g=\left\vert f\right\vert $. Hence, $g\left(  x\right)  =\left\vert
f\right\vert \left(  x\right)  =\left\vert \underbrace{f\left(  x\right)
}_{=-h}\right\vert =\left\vert -h\right\vert =h$, so that $x\in g^{-1}\left(
h\right)  $. Similarly, $y\in g^{-1}\left(  h\right)  $.

If $h=\infty$, then it is easy to see that $\pi\left(  x\right)  >\pi\left(
y\right)  $\ \ \ \ \footnote{\textit{Proof.} Assume that $h=\infty$. Now,
recall that $x\in g^{-1}\left(  h\right)  $. In other words, $x\in
g^{-1}\left(  \infty\right)  $ (since $h=\infty$). Similarly, $y\in
g^{-1}\left(  \infty\right)  $. But the map $\pi\mid_{g^{-1}\left(
\infty\right)  }$ is strictly decreasing (since Condition \textbf{(iii')} from
Definition \ref{def.amen} holds). Hence, $\left(  \pi\mid_{g^{-1}\left(
\infty\right)  }\right)  \left(  x\right)  >\left(  \pi\mid_{g^{-1}\left(
\infty\right)  }\right)  \left(  y\right)  $ (because $x<y$ and because both
$x$ and $y$ belong to $g^{-1}\left(  \infty\right)  $). Thus, $\pi\left(
x\right)  =\left(  \pi\mid_{g^{-1}\left(  \infty\right)  }\right)  \left(
x\right)  >\left(  \pi\mid_{g^{-1}\left(  \infty\right)  }\right)  \left(
y\right)  =\pi\left(  y\right)  $, qed.}. Hence, for the rest of this proof of
$\pi\left(  x\right)  >\pi\left(  y\right)  $, we WLOG assume that we don't
have $h=\infty$. Thus, $h\neq\infty$. Combining this with $h\in\mathcal{N}%
\setminus\left\{  0\right\}  $, we obtain $h\in\left(  \mathcal{N}%
\setminus\left\{  0\right\}  \right)  \setminus\left\{  \infty\right\}
=\mathcal{N}\setminus\left\{  0,\infty\right\}  =\left\{  1,2,3,\ldots
\right\}  $.

Combining this with $h=g\left(  \underbrace{x}_{\in\left[  n\right]  }\right)
\in g\left(  \left[  n\right]  \right)  $, we obtain $h\in g\left(  \left[
n\right]  \right)  \cap\left\{  1,2,3,\ldots\right\}  =H$. Also, $x\in
g^{-1}\left(  h\right)  $. Hence, if we had $x>m_{h}$, then statement
\textbf{(x}$_{4}$\textbf{)} would yield $f\left(  x\right)  =+h$, which would
contradict $f\left(  x\right)  =-h\neq+h$. Therefore, we cannot have $x>m_{h}%
$. Hence, we have $x\leq m_{h}$. In other words, $x\in\left\{  s\in
g^{-1}\left(  h\right)  \ \mid\ s\leq m_{h}\right\}  $ (since $x\in
g^{-1}\left(  h\right)  $). Similarly, $y\in\left\{  s\in g^{-1}\left(
h\right)  \ \mid\ s\leq m_{h}\right\}  $.

But $h\in g\left(  \left[  n\right]  \right)  \cap\left\{  1,2,3,\ldots
\right\}  $. Thus, the map $\pi\mid_{g^{-1}\left(  h\right)  }$ is V-shaped
(since Condition \textbf{(ii')} from Definition \ref{def.amen} holds). In
other words, there exists some $t\in g^{-1}\left(  h\right)  $ such that the
map $\left(  \pi\mid_{g^{-1}\left(  h\right)  }\right)  \mid_{\left\{  s\in
g^{-1}\left(  h\right)  \ \mid\ s\leq t\right\}  }$ is strictly decreasing
while the map $\left(  \pi\mid_{g^{-1}\left(  h\right)  }\right)
\mid_{\left\{  s\in g^{-1}\left(  h\right)  \ \mid\ s\geq t\right\}  }$ is
strictly increasing. Clearly, this $t$ must be the unique element $\mu$ of
$g^{-1}\left(  h\right)  $ for which $\left(  \pi\mid_{g^{-1}\left(  h\right)
}\right)  \left(  \mu\right)  $ is minimum. In other words, this $t$ must be
the unique element $\mu$ of $g^{-1}\left(  h\right)  $ for which $\pi\left(
\mu\right)  $ is minimum (since $\left(  \pi\mid_{g^{-1}\left(  h\right)
}\right)  \left(  \mu\right)  =\pi\left(  \mu\right)  $ for each $\mu\in
g^{-1}\left(  h\right)  $). In other words, this $t$ must be $m_{h}$ (because
$m_{h}$ was defined to be the unique element $\mu$ of $g^{-1}\left(  h\right)
$ for which $\pi\left(  \mu\right)  $ is minimum). Thus, we conclude that the
map $\left(  \pi\mid_{g^{-1}\left(  h\right)  }\right)  \mid_{\left\{  s\in
g^{-1}\left(  h\right)  \ \mid\ s\leq m_{h}\right\}  }$ is strictly decreasing
while the map \newline$\left(  \pi\mid_{g^{-1}\left(  h\right)  }\right)
\mid_{\left\{  s\in g^{-1}\left(  h\right)  \ \mid\ s\geq m_{h}\right\}  }$ is
strictly increasing.

In particular, the map $\left(  \pi\mid_{g^{-1}\left(  h\right)  }\right)
\mid_{\left\{  s\in g^{-1}\left(  h\right)  \ \mid\ s\leq m_{h}\right\}  }$ is
strictly decreasing. In other words, the map $\pi\mid_{\left\{  s\in
g^{-1}\left(  h\right)  \ \mid\ s\leq m_{h}\right\}  }$ is strictly decreasing
(since \newline$\left(  \pi\mid_{g^{-1}\left(  h\right)  }\right)
\mid_{\left\{  s\in g^{-1}\left(  h\right)  \ \mid\ s\leq m_{h}\right\}  }%
=\pi\mid_{\left\{  s\in g^{-1}\left(  h\right)  \ \mid\ s\leq m_{h}\right\}
}$). Since both $x$ and $y$ belong to the set $\left\{  s\in g^{-1}\left(
h\right)  \ \mid\ s\leq m_{h}\right\}  $, we thus conclude that%
\[
\left(  \pi\mid_{\left\{  s\in g^{-1}\left(  h\right)  \ \mid\ s\leq
m_{h}\right\}  }\right)  \left(  x\right)  >\left(  \pi\mid_{\left\{  s\in
g^{-1}\left(  h\right)  \ \mid\ s\leq m_{h}\right\}  }\right)  \left(
y\right)
\]
(since $x<y$). Thus,%
\[
\pi\left(  x\right)  =\left(  \pi\mid_{\left\{  s\in g^{-1}\left(  h\right)
\ \mid\ s\leq m_{h}\right\}  }\right)  \left(  x\right)  >\left(  \pi
\mid_{\left\{  s\in g^{-1}\left(  h\right)  \ \mid\ s\leq m_{h}\right\}
}\right)  \left(  y\right)  =\pi\left(  y\right)  .
\]
This proves Condition \textbf{(iii)}.]

Now, forget that we fixed $x$ and $y$. We thus have shown that for all $x<y$
in $\left[  n\right]  $, the conditions \textbf{(i)}, \textbf{(ii)} and
\textbf{(iii)} stated above hold. Thus, $f$ is a $\mathcal{Z}$-enriched
$\left(  \left[  n\right]  ,\pi\right)  $-partition (by the definition of a
$\mathcal{Z}$-enriched $\left(  \left[  n\right]  ,\pi\right)  $-partition).
In other words, $f\in\mathcal{E}\left(  \left[  n\right]  ,\pi\right)  $.
Hence, we have shown that $\left(  f\in\mathcal{E}\left(  \left[  n\right]
,\pi\right)  \text{ and }\left\vert f\right\vert =g\right)  $. This proves
Claim 2.]

Combining Claim 1 with Claim 2, we conclude that $\left(  f\in\mathcal{E}%
\left(  \left[  n\right]  ,\pi\right)  \text{ and }\left\vert f\right\vert
=g\right)  $ if and only if the five statements \textbf{(x}$_{1}$\textbf{)},
\textbf{(x}$_{2}$\textbf{)}, \textbf{(x}$_{3}$\textbf{)}, \textbf{(x}$_{4}%
$\textbf{)} and \textbf{(x}$_{5}$\textbf{)} hold. This proves Proposition
\ref{prop.amen.equiv3} \textbf{(b)}.
\end{proof}
\end{verlong}

\begin{proposition}
\label{prop.Epk-formula}Let $n\in\mathbb{N}$. Let $\pi$ be any $n$%
-permutation. Then,%
\[
\Gamma_{\mathcal{Z}}\left(  \pi\right)  =\sum_{\substack{g:\left[  n\right]
\rightarrow\mathcal{N}\\\text{is }\pi\text{-amenable}}}2^{\left\vert g\left(
\left[  n\right]  \right)  \cap\left\{  1,2,3,\ldots\right\}  \right\vert
}\mathbf{x}_{g}.
\]

\end{proposition}

\begin{vershort}
\begin{proof}
[Proof of Proposition \ref{prop.Epk-formula} (sketched).]
The claim will immediately follow from (\ref{eq.Gamma.rewr})
once we have shown the following two observations:

\begin{statement}
\textit{Observation 1:} If $f\in\mathcal{E}\left(  \left[  n\right]
,\pi\right)  $, then the map $\left\vert f\right\vert :\left[  n\right]
\rightarrow\mathcal{N}$ is $\pi$-amenable.
\end{statement}

\begin{statement}
\textit{Observation 2:} If $g:\left[  n\right]  \rightarrow\mathcal{N}$ is a
$\pi$-amenable map, then there exist precisely $2^{\left\vert g\left(  \left[
n\right]  \right)  \cap\left\{  1,2,3,\ldots\right\}  \right\vert }$ maps
$f\in\mathcal{E}\left(  \left[  n\right]  ,\pi\right)  $ satisfying
$\left\vert f\right\vert =g$.
\end{statement}

It thus remains to prove these two observations. Let us do this:

[\textit{Proof of Observation 1:} Let $f\in\mathcal{E}\left(  \left[
n\right]  ,\pi\right)  $. Thus, $f$ is a $\mathcal{Z}$-enriched $\left(
\left[  n\right]  ,\pi\right)  $-partition. In other words, $f$ is a map
$\left[  n\right]  \rightarrow\mathcal{Z}$ such that for all $x<y$ in $\left[
n\right]  $, the following conditions hold:

\begin{enumerate}
\item[\textbf{(i)}] We have $f\left(  x\right)  \preccurlyeq f\left(
y\right)  $.

\item[\textbf{(ii)}] If $f\left(  x\right)  =f\left(  y\right)  =+h$ for some
$h\in\mathcal{N}$, then $\pi\left(  x\right)  <\pi\left(  y\right)  $.

\item[\textbf{(iii)}] If $f\left(  x\right)  =f\left(  y\right)  =-h$ for some
$h\in\mathcal{N}$, then $\pi\left(  x\right)  >\pi\left(  y\right)  $.
\end{enumerate}

(This is due to the
definition of a $\mathcal{Z}$-enriched $\left(  \left[  n\right]  ,\pi\right)
$-partition.)

Condition \textbf{(i)} shows that the map $f$ is weakly increasing. Condition
\textbf{(ii)} shows that for each $h\in\mathcal{N}$, the map $\pi\mid
_{f^{-1}\left(  +h\right)  }$ is strictly increasing. Condition \textbf{(iii)}
shows that for each $h\in\mathcal{N}$, the map $\pi\mid_{f^{-1}\left(
-h\right)  }$ is strictly decreasing.

Now, set $g=\left\vert f\right\vert $. Then, $g^{-1}\left(  0\right)
=f^{-1}\left(  +0\right)  $ (since $-0\notin\mathcal{Z}$). But the map
$\pi\mid_{f^{-1}\left(  +0\right)  }$ is strictly increasing\footnote{because
for each $h\in\mathcal{N}$, the map $\pi\mid_{f^{-1}\left(  +h\right)  }$ is
strictly increasing}. Thus, the map $\pi\mid_{g^{-1}\left(  0\right)  }$ is
strictly increasing (since $g^{-1}\left(  0\right)  =f^{-1}\left(  +0\right)
$). Hence, Property \textbf{(i')} in Definition~\ref{def.amen} holds.
Similarly, Property \textbf{(iii')} in that definition also holds.

Now, fix $h\in g\left(  \left[  n\right]  \right)  \cap\left\{  1,2,3,\ldots
\right\}  $. Then, the set $g^{-1}\left(  h\right)  $ is nonempty (since $h\in
g\left(  \left[  n\right]  \right)  $), and can be written as the union of its
two disjoint subsets $f^{-1}\left(  +h\right)  $ and $f^{-1}\left(  -h\right)
$. Furthermore, each element of $f^{-1}\left(  -h\right)  $ is smaller than
each element of $f^{-1}\left(  +h\right)  $ (since $f$ is weakly increasing),
and we know that the map $\pi\mid_{f^{-1}\left(  -h\right)  }$ is strictly
decreasing while the map $\pi\mid_{f^{-1}\left(  +h\right)  }$ is strictly
increasing. Hence, the map $\pi\mid_{g^{-1}\left(  h\right)  }$ is strictly
decreasing up until some value of its argument (namely, either the largest
element of $f^{-1}\left(  -h\right)  $, or the smallest
element of $f^{-1}\left(  +h\right)  $, depending on which of these two
elements has the smaller image under $\pi$), and then strictly increasing
from there on. In other words, the map $\pi\mid_{g^{-1}\left(  h\right)  }$ is
V-shaped. Thus, Property \textbf{(ii')} in Definition~\ref{def.amen} holds.
Finally, Property \textbf{(iv')} in Definition~\ref{def.amen} holds because
$f$ is weakly increasing. We have
thus checked all four properties in Definition~\ref{def.amen};
thus, $g$ is $\pi$-amenable.
In other words, $\left| f \right|$ is $\pi$-amenable
(since $g = \left| f \right|$).
This proves Observation 1.]

[\textit{Proof of Observation 2:} Let $g:\left[  n\right]  \rightarrow
\mathcal{N}$ be a $\pi$-amenable map. Consider a map $f\in\mathcal{E}\left(
\left[  n\right]  ,\pi\right)  $ satisfying $\left\vert f\right\vert =g$. We
are wondering to what extent the map $f$ is determined by $g$ and $\pi$.

Everything that we said in the proof of Observation 1 still holds in
our situation (since $g = \left| f \right|$).

In order to determine the map $f$, it clearly suffices to determine the sets
$f^{-1}\left(  q\right)  $ for all $q\in\mathcal{Z}$. In other words, it
suffices to determine the set $f^{-1}\left(  +0\right)  $, the set
$f^{-1}\left(  -\infty\right)  $ and the sets $f^{-1}\left(  +h\right)  $ and
$f^{-1}\left(  -h\right)  $ for all $h\in\left\{  1,2,3,\ldots\right\}  $.

Recall from the proof of Observation 1 that $g^{-1}\left(  0\right)
=f^{-1}\left(  +0\right)  $. Thus, $f^{-1}\left(  +0\right)  $ is uniquely
determined by $g$. Similarly, $f^{-1}\left(  -\infty\right)  $ is uniquely
determined by $g$. Thus, we can focus on the remaining sets $f^{-1}\left(
+h\right)  $ and $f^{-1}\left(  -h\right)  $ for $h\in\left\{  1,2,3,\ldots
\right\}  $.

Fix $h\in\left\{  1,2,3,\ldots\right\}  $. Recall that the set $g^{-1}\left(
h\right)  $ is the union of its two disjoint subsets $f^{-1}\left(  +h\right)
$ and $f^{-1}\left(  -h\right)  $. Thus, $f^{-1}\left(  +h\right)  $ and
$f^{-1}\left(  -h\right)  $ are complementary subsets of $g^{-1}\left(
h\right)  $. If $g^{-1}\left(  h\right)  =\varnothing$, then this uniquely
determines $f^{-1}\left(  +h\right)  $ and $f^{-1}\left(  -h\right)  $. Thus,
we focus only on the case when $g^{-1}\left(  h\right)  \neq\varnothing$.

So assume that $g^{-1}\left(  h\right)  \neq\varnothing$. Hence, $h\in
g\left(  \left[  n\right]  \right)  $, so that $h\in g\left(  \left[
n\right]  \right)  \cap\left\{  1,2,3,\ldots\right\}  $. Since the map $g$ is
$\pi$-amenable, we thus conclude that the map $\pi\mid_{g^{-1}\left(
h\right)  }$ is V-shaped (by Property \textbf{(ii')} in
Definition~\ref{def.amen}).

The map $g$ is weakly increasing (by Property \textbf{(iv')} in
Definition~\ref{def.amen}). Hence,
$g^{-1}\left(  h\right)  $ is an interval of $\left[  n\right]  $. Let
$\alpha\in\mathbb{Z}$ and $\gamma\in\mathbb{Z}$ be such that $g^{-1}\left(
h\right)  =\left[  \alpha,\gamma\right]  $ (where $\left[  \alpha
,\gamma\right]  $ means the interval $\left\{  \alpha,\alpha+1,\ldots
,\gamma\right\}  $).

As in the proof of Observation 1, we can see that each element of
$f^{-1}\left(  -h\right)  $ is smaller than each element of $f^{-1}\left(
+h\right)  $. Since the union of $f^{-1}\left(  -h\right)  $ and
$f^{-1}\left(  +h\right)  $ is $g^{-1}\left(  h\right)  =\left[  \alpha
,\gamma\right]  $, we thus conclude that (roughly speaking) the
sets $f^{-1}\left(  -h\right)  $ and $f^{-1}\left(  +h\right)  $
partition the interval $g^{-1}\left(h\right)$ in two (possibly empty)
sub-intervals such that the interval $f^{-1}\left(  -h\right)  $ lies
completely to the left of $f^{-1}\left(  +h\right)  $.
Hence, there exists some $\beta\in\left[
\alpha-1,\gamma\right]  $ such that $f^{-1}\left(  -h\right)  =\left[
\alpha,\beta\right]  $ and $f^{-1}\left(  +h\right)  =\left[  \beta
+1,\gamma\right]  $. Consider this $\beta$. Clearly, $f^{-1}\left(  -h\right)
$ and $f^{-1}\left(  +h\right)  $ are uniquely determined by $\beta$; we just
need to find out which values $\beta$ can take.

As in the proof of Observation 1, we can see that the map $\pi\mid
_{f^{-1}\left(  -h\right)  }$ is strictly decreasing while the map $\pi
\mid_{f^{-1}\left(  +h\right)  }$ is strictly increasing. Let $k$ be the
element of $g^{-1}\left(  h\right)  $ minimizing $\pi\left(  k\right)  $.
Then, the map $\pi$ is strictly decreasing on the set $\left\{  u\in
g^{-1}\left(  h\right)  \ \mid\ u\leq k\right\}  $ and strictly increasing on
the set $\left\{  u\in g^{-1}\left(  h\right)  \ \mid\ u\geq k\right\}  $
(since the map $\pi\mid_{g^{-1}\left(  h\right)  }$ is V-shaped).

The map $\pi\mid_{f^{-1}\left(  -h\right)  }$ is strictly decreasing. In other
words, the map $\pi$ is strictly decreasing on the set $f^{-1}\left(
-h\right)  =\left[  \alpha,\beta\right]  $. On the other hand, the map $\pi$
is strictly increasing on the set $\left\{  u\in g^{-1}\left(  h\right)
\ \mid\ u\geq k\right\}  $. Hence, the two sets $\left[  \alpha,\beta\right]
$ and $\left\{  u\in g^{-1}\left(  h\right)  \ \mid\ u\geq k\right\}  $ cannot
have more than one element in common
(since $\pi$ is strictly decreasing on one
and strictly increasing on the other). Thus, $k\geq\beta$. A similar argument
shows that $k\leq\beta+1$. Combining these inequalities, we obtain
$k\in\left\{  \beta,\beta+1\right\}  $, so that $\beta\in\left\{
k,k-1\right\}  $. This shows that $\beta$ can take only two values: $k$ and
$k-1$.

Now, let us forget that we fixed $h$.
We have shown that for each $h\in g\left(
\left[  n\right]  \right)  \cap\left\{  1,2,3,\ldots\right\}  $, the sets
$f^{-1}\left(  +h\right)  $ and $f^{-1}\left(  -h\right)  $ are uniquely
determined once the integer $\beta$ is chosen, and that this integer $\beta$
can be chosen in two ways. (As we have seen, all other values of $h$ do not
matter.) Thus, in total, the map $f$ is uniquely determined up to $\left\vert
g\left(  \left[  n\right]  \right)  \cap\left\{  1,2,3,\ldots\right\}
\right\vert $ decisions, where each decision allows choosing from two values.
Thus, there are at most $2^{\left\vert g\left(  \left[  n\right]  \right)
\cap\left\{  1,2,3,\ldots\right\}  \right\vert }$ maps $f\in\mathcal{E}\left(
\left[  n\right]  ,\pi\right)  $ satisfying $\left\vert f\right\vert =g$.
Working the above argument backwards, we see that each way of making these
decisions
actually leads to a map $f\in\mathcal{E}\left(  \left[  n\right]  ,\pi\right)
$ satisfying $\left\vert f\right\vert =g$; thus, there are \textbf{exactly}
$2^{\left\vert g\left(  \left[  n\right]  \right)  \cap\left\{  1,2,3,\ldots
\right\}  \right\vert }$ maps $f\in\mathcal{E}\left(  \left[  n\right]
,\pi\right)  $ satisfying $\left\vert f\right\vert =g$. This proves
Observation 2.]
\end{proof}
\end{vershort}

\begin{verlong}
\begin{proof}
[Proof of Proposition \ref{prop.Epk-formula}.]We shall show the following observation:

\begin{statement}
\textit{Observation 1:} If $g:\left[  n\right]  \rightarrow\mathcal{N}$ is a
$\pi$-amenable map, then there exist precisely $2^{\left\vert g\left(  \left[
n\right]  \right)  \cap\left\{  1,2,3,\ldots\right\}  \right\vert }$ maps
$f\in\mathcal{E}\left(  \left[  n\right]  ,\pi\right)  $ satisfying
$\left\vert f\right\vert =g$.
\end{statement}

[\textit{Proof of Observation 1:} Let $g:\left[  n\right]  \rightarrow
\mathcal{N}$ be a $\pi$-amenable map. Let $H$ be the set $g\left(  \left[
n\right]  \right)  \cap\left\{  1,2,3,\ldots\right\}  $. For each $h\in H$, we
let $m_{h}$ be the unique element $\mu$ of $g^{-1}\left(  h\right)  $ for
which $\pi\left(  \mu\right)  $ is minimum. Proposition \ref{prop.amen.equiv3}
\textbf{(a)} shows that the elements $m_{h}$ for all $h\in H$ are well-defined
and distinct.

Proposition \ref{prop.amen.equiv3} \textbf{(b)} shows that a map $f:\left[
n\right]  \rightarrow\mathcal{Z}$ satisfies \newline
$\left(  f\in\mathcal{E}\left(
\left[  n\right]  ,\pi\right)  \text{ and }\left\vert f\right\vert =g\right)
$ if and only if the five statements \textbf{(x}$_{1}$\textbf{)},
\textbf{(x}$_{2}$\textbf{)}, \textbf{(x}$_{3}$\textbf{)}, \textbf{(x}$_{4}%
$\textbf{)} and \textbf{(x}$_{5}$\textbf{)} from Proposition
\ref{prop.amen.equiv3} \textbf{(b)} hold. Thus, we can construct every map
$f:\left[  n\right]  \rightarrow\mathcal{Z}$ that satisfies $\left(
f\in\mathcal{E}\left(  \left[  n\right]  ,\pi\right)  \text{ and }\left\vert
f\right\vert =g\right)  $ by the following algorithm:

\begin{itemize}
\item[\textbf{Step 1:}] For each $x\in g^{-1}\left(  0\right)  $, set
$f\left(  x\right)  =+0$. (This is the only option, because we want statement
\textbf{(x}$_{1}$\textbf{)} to hold.)

\item[\textbf{Step 2:}] For each $h\in H$ and each $x\in g^{-1}\left(
h\right)  $, set the value $f\left(  x\right)  $ as follows:

\begin{itemize}
\item If $x<m_{h}$, then set $f\left(  x\right)  =-h$. (This is the only
option, because we want statement \textbf{(x}$_{2}$\textbf{)} to hold.)

\item If $x=m_{h}$, then set $f\left(  x\right)  $ to be either $-h$ or $+h$.
(These two options are the only options, because we want statement
\textbf{(x}$_{3}$\textbf{)} to hold. Notice that $m_{h}\in g^{-1}\left(
h\right)  $ (by the definition of $m_{h}$); therefore, this step ensures that
$f\left(  m_{h}\right)  \in\left\{  -h,+h\right\}  $ for each $h\in H$, and
therefore statement \textbf{(x}$_{3}$\textbf{)} holds indeed.)

\item If $x>m_{h}$, then set $f\left(  x\right)  =+h$. (This is the only
option, because we want statement \textbf{(x}$_{4}$\textbf{)} to hold.)
\end{itemize}

\item[\textbf{Step 3:}] For each $x\in g^{-1}\left(  \infty\right)  $, set
$f\left(  x\right)  =-\infty$. (This is the only option, because we want
statement \textbf{(x}$_{5}$\textbf{)} to hold.)
\end{itemize}

This algorithm indeed constructs a well-defined map $f:\left[  n\right]
\rightarrow\mathcal{Z}$ (because for each $x\in\left[  n\right]  $, the value
$f\left(  x\right)  $ is set exactly once during the above
algorithm\footnote{\textit{Proof.} Let $x\in\left[  n\right]  $ be arbitrary.
Note that $H=g\left(  \left[  n\right]  \right)  \cap\left\{  1,2,3,\ldots
\right\}  \subseteq\left\{  1,2,3,\ldots\right\}  $.
\par
We have $g\left(  x\right)  \in\mathcal{N}=\left\{  0\right\}  \cup\left\{
\infty\right\}  \cup\left\{  1,2,3,\ldots\right\}  $. Hence, we have either
$g\left(  x\right)  =0$ or $g\left(  x\right)  =\infty$ or $g\left(  x\right)
\in\left\{  1,2,3,\ldots\right\}  $. In other words, we are in one of the
following cases:
\par
\textit{Case 1:} We have $g\left(  x\right)  =0$.
\par
\textit{Case 2:} We have $g\left(  x\right)  =\infty$.
\par
\textit{Case 3:} We have $g\left(  x\right)  \in\left\{  1,2,3,\ldots\right\}
$.
\par
Let us first consider Case 1. In this case, we have $g\left(  x\right)  =0$.
Thus, $x\in g^{-1}\left(  0\right)  $. Hence, $f\left(  x\right)  $ is set
during Step 1 of the above algorithm. If we had $0\in H$, then we would have
$0\in H\subseteq\left\{  1,2,3,\ldots\right\}  $, which would contradict
$0\notin\left\{  1,2,3,\ldots\right\}  $. Thus, we cannot have $0\in H$.
Hence, $0\notin H$. Hence, $x\notin g^{-1}\left(  h\right)  $ for all $h\in H$
(because if we had $x\in g^{-1}\left(  h\right)  $ for some $h\in H$, then
this $h$ would satisfy $g\left(  x\right)  =h\in H$, which would contradict
$g\left(  x\right)  =0\notin H$). Therefore, Step 2 of the above algorithm
does not set $f\left(  x\right)  $. Finally, $x\notin g^{-1}\left(
\infty\right)  $ (since $g\left(  x\right)  =0\neq\infty$). Hence, Step 3 of
the above algorithm does not set $f\left(  x\right)  $. Altogether, we have
now shown that $f\left(  x\right)  $ is set during Step 1 of the algorithm,
but neither Step 2 nor Step 3 sets $f\left(  x\right)  $. Thus, the value
$f\left(  x\right)  $ is set exactly once during the above algorithm. So we
have proven our claim (that the value $f\left(  x\right)  $ is set exactly
once during the above algorithm) in Case 1.
\par
Let us next consider Case 2. In this case, we have $g\left(  x\right)
=\infty$. Thus, $x\in g^{-1}\left(  \infty\right)  $. Hence, $f\left(
x\right)  $ is set during Step 3 of the above algorithm. If we had $\infty\in
H$, then we would have $\infty\in H\subseteq\left\{  1,2,3,\ldots\right\}  $,
which would contradict $\infty\notin\left\{  1,2,3,\ldots\right\}  $. Thus, we
cannot have $\infty\in H$. Hence, $\infty\notin H$. Hence, $x\notin
g^{-1}\left(  h\right)  $ for all $h\in H$ (because if we had $x\in
g^{-1}\left(  h\right)  $ for some $h\in H$, then this $h$ would satisfy
$g\left(  x\right)  =h\in H$, which would contradict $g\left(  x\right)
=\infty\notin H$). Therefore, Step 2 of the above algorithm does not set
$f\left(  x\right)  $. Finally, $x\notin g^{-1}\left(  0\right)  $ (since
$g\left(  x\right)  =\infty\neq0$). Hence, Step 1 of the above algorithm does
not set $f\left(  x\right)  $. Altogether, we have now shown that $f\left(
x\right)  $ is set during Step 3 of the algorithm, but neither Step 2 nor Step
1 sets $f\left(  x\right)  $. Thus, the value $f\left(  x\right)  $ is set
exactly once during the above algorithm. So we have proven our claim (that the
value $f\left(  x\right)  $ is set exactly once during the above algorithm) in
Case 2.
\par
Let us next consider Case 3. In this case, we have $g\left(  x\right)
\in\left\{  1,2,3,\ldots\right\}  $. Combining this with $g\left(
\underbrace{x}_{\in\left[  n\right]  }\right)  \in g\left(  \left[  n\right]
\right)  $, we obtain $g\left(  x\right)  \in g\left(  \left[  n\right]
\right)  \cap\left\{  1,2,3,\ldots\right\}  =H$. Hence, there exists an $h\in
H$ such that $x\in g^{-1}\left(  h\right)  $ (namely, $h=g\left(  x\right)
$). Of course, this $h$ is unique (because the requirement $x\in g^{-1}\left(
h\right)  $ entails $g\left(  x\right)  =h$, which uniquely determines $h$).
Therefore, the value $f\left(  x\right)  $ is set exactly once during Step 2
of the algorithm. Moreover, from $g\left(  x\right)  \in\left\{
1,2,3,\ldots\right\}  $, we obtain $g\left(  x\right)  \neq0$, so that
$x\notin g^{-1}\left(  0\right)  $. Hence, Step 1 of the above algorithm does
not set $f\left(  x\right)  $. Similarly, Step 3 of the above algorithm does
not set $f\left(  x\right)  $. Altogether, we have now shown that $f\left(
x\right)  $ is set exactly once during Step 2 of the algorithm, but neither
Step 1 nor Step 3 sets $f\left(  x\right)  $. Thus, the value $f\left(
x\right)  $ is set exactly once during the above algorithm. So we have proven
our claim (that the value $f\left(  x\right)  $ is set exactly once during the
above algorithm) in Case 3.
\par
We have now proven our claim (that the value $f\left(  x\right)  $ is set
exactly once during the above algorithm) in each of the three Cases 1, 2 and
3. Hence, this claim always holds. Qed.}). Moreover, this map $f$ that it
constructs satisfies $\left(  f\in\mathcal{E}\left(  \left[  n\right]
,\pi\right)  \text{ and }\left\vert f\right\vert =g\right)  $ (since the five
statements \textbf{(x}$_{1}$\textbf{)}, \textbf{(x}$_{2}$\textbf{)},
\textbf{(x}$_{3}$\textbf{)}, \textbf{(x}$_{4}$\textbf{)} and \textbf{(x}$_{5}%
$\textbf{)} from Proposition \ref{prop.amen.equiv3} \textbf{(b)} hold for this
map $f$). Therefore, the number of all maps $f:\left[  n\right]
\rightarrow\mathcal{Z}$ satisfying $\left(  f\in\mathcal{E}\left(  \left[
n\right]  ,\pi\right)  \text{ and }\left\vert f\right\vert =g\right)  $ equals
the number of ways we can perform the above algorithm. But the latter number
is easy to compute: The only choices we have in the algorithm are the choices
we make during Step 2 when $x=m_{h}$ (since in all the other decisions, we
have only one option); there are altogether $\left\vert H\right\vert $ many of
these choices (one for each $h\in H$), and each choice involves exactly $2$
options. Hence, the number of ways we can perform the above algorithm is
$2^{\left\vert H\right\vert }$. Thus, the number of all maps $f:\left[
n\right]  \rightarrow\mathcal{Z}$ satisfying $\left(  f\in\mathcal{E}\left(
\left[  n\right]  ,\pi\right)  \text{ and }\left\vert f\right\vert =g\right)
$ is $2^{\left\vert H\right\vert }$. In other words,%
\begin{align*}
&  \left(  \text{the number of all maps }f:\left[  n\right]  \rightarrow
\mathcal{Z}\text{ satisfying }\left(  f\in\mathcal{E}\left(  \left[  n\right]
,\pi\right)  \text{ and }\left\vert f\right\vert =g\right)  \right)  \\
&  =2^{\left\vert H\right\vert }.
\end{align*}

Now,%
\begin{align*}
&  \left(  \text{the number of all maps }f\in\mathcal{E}\left(  \left[
n\right]  ,\pi\right)  \text{ satisfying }\left\vert f\right\vert =g\right)
\\
&  =\left(  \text{the number of all maps }f:\left[  n\right]  \rightarrow
\mathcal{Z}\text{ satisfying }\left(  f\in\mathcal{E}\left(  \left[  n\right]
,\pi\right)  \text{ and }\left\vert f\right\vert =g\right)  \right)  \\
&  \ \ \ \ \ \ \ \ \ \ \left(  \text{since each }f\in\mathcal{E}\left(
\left[  n\right]  ,\pi\right)  \text{ is automatically a map }\left[
n\right]  \rightarrow\mathcal{Z}\right)  \\
&  =2^{\left\vert H\right\vert }=2^{\left\vert g\left(  \left[  n\right]
\right)  \cap\left\{  1,2,3,\ldots\right\}  \right\vert }%
\ \ \ \ \ \ \ \ \ \ \left(  \text{since }H=g\left(  \left[  n\right]  \right)
\cap\left\{  1,2,3,\ldots\right\}  \right)  .
\end{align*}
In other words, there exist precisely $2^{\left\vert g\left(  \left[
n\right]  \right)  \cap\left\{  1,2,3,\ldots\right\}  \right\vert }$ maps
$f\in\mathcal{E}\left(  \left[  n\right]  ,\pi\right)  $ satisfying
$\left\vert f\right\vert =g$. This proves Observation 1.]

Now, (\ref{eq.Gamma.rewr}) yields%
\begin{align*}
\Gamma_{\mathcal{Z}}\left(  \pi\right)   &  =\underbrace{\sum_{f\in
\mathcal{E}\left(  \left[  n\right]  ,\pi\right)  }}_{\substack{=\sum
_{\substack{g:\left[  n\right]  \rightarrow\mathcal{N}\\\text{is }%
\pi\text{-amenable}}}\sum_{\substack{f\in\mathcal{E}\left(  \left[  n\right]
,\pi\right)  ;\\\left\vert f\right\vert =g}}\\\text{(by Proposition
\ref{prop.amen.equiv2})}}}\mathbf{x}_{\left\vert f\right\vert }=\sum
_{\substack{g:\left[  n\right]  \rightarrow\mathcal{N}\\\text{is }%
\pi\text{-amenable}}}\sum_{\substack{f\in\mathcal{E}\left(  \left[  n\right]
,\pi\right)  ;\\\left\vert f\right\vert =g}}\underbrace{\mathbf{x}_{\left\vert
f\right\vert }}_{\substack{=\mathbf{x}_{g}\\\text{(since }\left\vert
f\right\vert =g\text{)}}}\\
&  =\sum_{\substack{g:\left[  n\right]  \rightarrow\mathcal{N}\\\text{is }%
\pi\text{-amenable}}}\underbrace{\sum_{\substack{f\in\mathcal{E}\left(
\left[  n\right]  ,\pi\right)  ;\\\left\vert f\right\vert =g}}\mathbf{x}_{g}%
}_{\substack{=2^{\left\vert g\left(  \left[  n\right]  \right)  \cap\left\{
1,2,3,\ldots\right\}  \right\vert }\mathbf{x}_{g}\\\text{(by Observation 1)}%
}}=\sum_{\substack{g:\left[  n\right]  \rightarrow\mathcal{N}\\\text{is }%
\pi\text{-amenable}}}2^{\left\vert g\left(  \left[  n\right]  \right)
\cap\left\{  1,2,3,\ldots\right\}  \right\vert }\mathbf{x}_{g}.
\end{align*}
This proves Proposition \ref{prop.Epk-formula}.
\end{proof}
\end{verlong}

Now, let us observe that if $g:\left[  n\right]  \rightarrow\mathcal{N}$ is a
weakly increasing map (for some $n\in\mathbb{N}$), then the fibers of $g$
(that is, the subsets $g^{-1}\left(  h\right)  $ of $\left[  n\right]  $ for
various $h\in\mathcal{N}$) are intervals of $\left[  n\right]  $ (possibly
empty). Of course, when these fibers are nonempty, they have smallest elements
and largest elements. We shall next study these elements more closely.

\begin{definition}
\label{def.fiberends}Let $n\in\mathbb{N}$. Let $g:\left[  n\right]
\rightarrow\mathcal{N}$ be any map. We define a subset $\operatorname*{FE}%
\left(  g\right)  $ of $\left[  n\right]  $ as follows:%
\begin{align*}
\operatorname*{FE}\left(  g\right)   &  =\left\{  \min\left(  g^{-1}\left(
h\right)  \right)  \ \mid\ h\in\left\{  1,2,3,\ldots,\infty\right\}  \text{
with }g^{-1}\left(  h\right)  \neq\varnothing\right\} \\
&  \ \ \ \ \ \ \ \ \ \ \cup\left\{  \max\left(  g^{-1}\left(  h\right)
\right)  \ \mid\ h\in\left\{  0,1,2,3,\ldots\right\}  \text{ with }%
g^{-1}\left(  h\right)  \neq\varnothing\right\}  .
\end{align*}
In other words, $\operatorname*{FE}\left(  g\right)  $ is the set comprising
the smallest elements of all nonempty fibers of $g$ except for $g^{-1}\left(
0\right)  $ as well as the largest elements of all nonempty fibers of $g$
except for $g^{-1}\left(  \infty\right)  $. We shall refer to the elements of
$\operatorname*{FE}\left(  g\right)  $ as the \textit{fiber-ends} of $g$.
\end{definition}

\begin{lemma}
\label{lem.FE.exist}Let $n\in\mathbb{N}$. Let $\Lambda\in\mathbf{L}_{n}$.
Then, there exists a weakly increasing map $g:\left[  n\right]  \rightarrow
\mathcal{N}$ such that $\operatorname*{FE}\left(  g\right)  =\left(
\Lambda\cup\left(  \Lambda+1\right)  \right)  \cap\left[  n\right]  $.
\end{lemma}

\begin{vershort}
\begin{proof}
[Proof of Lemma \ref{lem.FE.exist} (sketched).]If $n=0$, then Lemma
\ref{lem.FE.exist} holds for obvious reasons. Thus, WLOG assume that $n\neq0$.
Hence, $n$ is a positive integer. Thus, $\mathbf{L}_{n}$ is the set of all
nonempty lacunar subsets of $\left[  n\right]  $ (by the definition of
$\mathbf{L}_{n}$). Therefore, from $\Lambda\in\mathbf{L}_{n}$, we conclude
that $\Lambda$ is a nonempty lacunar subset of $\left[  n\right]  $. Write
this subset $\Lambda$ in the form $\Lambda=\left\{  j_{1}<j_{2}<\cdots
<j_{k}\right\}  $. Thus, $k\geq1$ (since $\Lambda$ is nonempty).

From $\left\{  j_{1}<j_{2}<\cdots<j_{k}\right\}  =\Lambda\subseteq\left[
n\right]  $, we obtain $0<j_{1}<j_{2}<\cdots<j_{k}\leq n$. Thus,
$1<j_{1}+1<j_{2}+1<\cdots<j_{k}+1\leq n+1$. Hence, the $k-1$ numbers
$j_{1}+1,j_{2}+1,\ldots,j_{k-1}+1$ all belong to the set $\left[  n\right]  $,
whereas the number $j_{k}+1$ only belongs to this set if $j_{k}<n$. Hence,%
\begin{align*}
& \left\{  j_{1}+1,j_{2}+1,\ldots,j_{k}+1\right\}  \cap\left[  n\right]  \\
& =%
\begin{cases}
\left\{  j_{1}+1,j_{2}+1,\ldots,j_{k-1}+1,j_{k}+1\right\}  , & \text{if }%
j_{k}<n;\\
\left\{  j_{1}+1,j_{2}+1,\ldots,j_{k-1}+1\right\}  , & \text{if }j_{k}=n
\end{cases}
.
\end{align*}
But $\Lambda=\left\{  j_{1}<j_{2}<\cdots<j_{k}\right\}  =\left\{  j_{1}%
,j_{2},\ldots,j_{k}\right\}  $ and thus \newline$\Lambda+1=\left\{
j_{1}+1,j_{2}+1,\ldots,j_{k}+1\right\}  $. Hence,%
\begin{align}
\left(  \Lambda+1\right)  \cap\left[  n\right]    & =\left\{  j_{1}%
+1,j_{2}+1,\ldots,j_{k}+1\right\}  \cap\left[  n\right]  \nonumber\\
& =%
\begin{cases}
\left\{  j_{1}+1,j_{2}+1,\ldots,j_{k-1}+1,j_{k}+1\right\}  , & \text{if }%
j_{k}<n;\\
\left\{  j_{1}+1,j_{2}+1,\ldots,j_{k-1}+1\right\}  , & \text{if }j_{k}=n
\end{cases}
.\label{pf.lem.FE.exist.short.2}%
\end{align}

Now, consider the map $g:\left[  n\right]  \rightarrow\mathcal{N}$ defined by%
\begin{align*}
g\left(  x\right)   &  =%
\begin{cases}
\left(  \text{the number of all }\lambda\in\Lambda\text{ such that }%
\lambda<x\right)  , & \text{if }x\leq j_{k};\\
\infty, & \text{if }x>j_{k}%
\end{cases}
\\
&  \ \ \ \ \ \ \ \ \ \ \text{for each }x\in\left[  n\right]  .
\end{align*}
Thus,%
\begin{align*}
&  \left(  g\left(  1\right)  ,g\left(  2\right)  ,\ldots,g\left(  n\right)
\right)  \\
&  =\left(  \underbrace{0,0,\ldots,0}_{j_{1}\text{ entries}}%
,\underbrace{1,1,\ldots,1}_{j_{2}-j_{1}\text{ entries}},\underbrace{2,2,\ldots
,2}_{j_{3}-j_{2}\text{ entries}},\ldots,\underbrace{k-1,k-1,\ldots,k-1}%
_{j_{k}-j_{k-1}\text{ entries}},\underbrace{\infty,\infty,\ldots,\infty
}_{n-j_{k}\text{ entries}}\right)  .
\end{align*}
The $n$-tuple on the right hand side of this equality consists of a block of
$0$'s, followed by a block of $1$'s, followed by a block of $2$'s, and so on,
all the way up to a block of $\left(  k-1\right)  $'s, which is then followed
by a block of $\infty$'s. The first $k$ of these blocks are nonempty (since
$0<j_{1}<j_{2}<\cdots<j_{k}$). The last block is nonempty if $j_{k}<n$, and
empty if $j_{k}=n$. Thus, the map $g$ is weakly increasing, and its nonempty
fibers are
\begin{align*}
g^{-1}\left(  0\right)   &  =\left\{  1,2,\ldots,j_{1}\right\}  ,\\
g^{-1}\left(  1\right)   &  =\left\{  j_{1}+1,j_{1}+2,\ldots,j_{2}\right\}
,\\
g^{-1}\left(  2\right)   &  =\left\{  j_{2}+1,j_{2}+2,\ldots,j_{3}\right\}
,\\
&  \vdots,\\
g^{-1}\left(  k-1\right)   &  =\left\{  j_{k-1}+1,j_{k-1}+2,\ldots
,j_{k}\right\}  ,\\
g^{-1}\left(  \infty\right)   &  =\left\{  j_{k}+1,j_{k}+2,\ldots,n\right\}  ,
\end{align*}
except that $g^{-1}\left(  \infty\right)  $ is empty when $j_{k}=n$. Hence,
the definition of $\operatorname*{FE}\left(  g\right)  $ yields%
\begin{align*}
&  \operatorname*{FE}\left(  g\right)  \\
&  =\underbrace{\left\{  \min\left(  g^{-1}\left(  h\right)  \right)
\ \mid\ h\in\left\{  1,2,3,\ldots,\infty\right\}  \text{ with }g^{-1}\left(
h\right)  \neq\varnothing\right\}  }_{=%
\begin{cases}
\left\{  \min\left(  g^{-1}\left(  1\right)  \right)  ,\min\left(
g^{-1}\left(  2\right)  \right)  ,\ldots,\min\left(  g^{-1}\left(  k-1\right)
\right)  ,\min\left(  g^{-1}\left(  \infty\right)  \right)  \right\}  , &
\text{if }j_{k}<n;\\
\left\{  \min\left(  g^{-1}\left(  1\right)  \right)  ,\min\left(
g^{-1}\left(  2\right)  \right)  ,\ldots,\min\left(  g^{-1}\left(  k-1\right)
\right)  \right\}  , & \text{if }j_{k}=n
\end{cases}
}\\
&  \ \ \ \ \ \ \ \ \ \ \cup\underbrace{\left\{  \max\left(  g^{-1}\left(
h\right)  \right)  \ \mid\ h\in\left\{  0,1,2,3,\ldots\right\}  \text{ with
}g^{-1}\left(  h\right)  \neq\varnothing\right\}  }_{=\left\{  \max\left(
g^{-1}\left(  0\right)  \right)  ,\max\left(  g^{-1}\left(  1\right)  \right)
,\max\left(  g^{-1}\left(  2\right)  \right)  ,\ldots,\max\left(
g^{-1}\left(  k-1\right)  \right)  \right\}  }\\
&  =\underbrace{%
\begin{cases}
\left\{  \min\left(  g^{-1}\left(  1\right)  \right)  ,\min\left(
g^{-1}\left(  2\right)  \right)  ,\ldots,\min\left(  g^{-1}\left(  k-1\right)
\right)  ,\min\left(  g^{-1}\left(  \infty\right)  \right)  \right\}  , &
\text{if }j_{k}<n;\\
\left\{  \min\left(  g^{-1}\left(  1\right)  \right)  ,\min\left(
g^{-1}\left(  2\right)  \right)  ,\ldots,\min\left(  g^{-1}\left(  k-1\right)
\right)  \right\}  , & \text{if }j_{k}=n
\end{cases}
}_{\substack{=%
\begin{cases}
\left\{  j_{1}+1,j_{2}+1,\ldots,j_{k-1}+1,j_{k}+1\right\}  , & \text{if }%
j_{k}<n;\\
\left\{  j_{1}+1,j_{2}+1,\ldots,j_{k-1}+1\right\}  , & \text{if }j_{k}=n
\end{cases}
\\\text{(since }\min\left(  g^{-1}\left(  h\right)  \right)  =j_{h}+1\text{
for each }h\in\left\{  1,2,\ldots,k-1\right\}  \text{,}\\\text{and since }%
\min\left(  g^{-1}\left(  \infty\right)  \right)  =j_{k}+1\text{ if }%
j_{k}<n\text{)}}}\\
&  \ \ \ \ \ \ \ \ \ \ \cup\underbrace{\left\{  \max\left(  g^{-1}\left(
0\right)  \right)  ,\max\left(  g^{-1}\left(  1\right)  \right)  ,\max\left(
g^{-1}\left(  2\right)  \right)  ,\ldots,\max\left(  g^{-1}\left(  k-1\right)
\right)  \right\}  }_{\substack{=\left\{  j_{1},j_{2},\ldots,j_{k}\right\}
\\\text{(since }\max\left(  g^{-1}\left(  h-1\right)  \right)  =j_{h}\text{
for each }h\in\left\{  1,2,\ldots,k\right\}  \text{)}}}\\
&  =\underbrace{%
\begin{cases}
\left\{  j_{1}+1,j_{2}+1,\ldots,j_{k-1}+1,j_{k}+1\right\}  , & \text{if }%
j_{k}<n;\\
\left\{  j_{1}+1,j_{2}+1,\ldots,j_{k-1}+1\right\}  , & \text{if }j_{k}=n
\end{cases}
}_{\substack{=\left(  \Lambda+1\right)  \cap\left[  n\right]  \\\text{(by
(\ref{pf.lem.FE.exist.short.2}))}}}\cup\underbrace{\left\{  j_{1},j_{2}%
,\ldots,j_{k}\right\}  }_{\substack{=\Lambda=\Lambda\cap\left[  n\right]
\\\text{(since }\Lambda\subseteq\left[  n\right]  \text{)}}}\\
&  =\left(  \left(  \Lambda+1\right)  \cap\left[  n\right]  \right)
\cup\left(  \Lambda\cap\left[  n\right]  \right)  =\left(  \Lambda\cap\left[
n\right]  \right)  \cup\left(  \left(  \Lambda+1\right)  \cap\left[  n\right]
\right)  =\left(  \Lambda\cup\left(  \Lambda+1\right)  \right)  \cap\left[
n\right]  .
\end{align*}
Altogether, we have now shown that our map $g:\left[  n\right]  \rightarrow
\mathcal{N}$ is weakly increasing and satisfies $\operatorname*{FE}\left(
g\right)  =\left(  \Lambda\cup\left(  \Lambda+1\right)  \right)  \cap\left[
n\right]  $. Hence, such a map $g$ exists. Thus, Lemma \ref{lem.FE.exist} is proven.
\end{proof}
\end{vershort}

\begin{verlong}
\begin{proof}
[Proof of Lemma \ref{lem.FE.exist}.]If $n=0$, then Lemma \ref{lem.FE.exist}
holds\footnote{\textit{Proof.} Assume that $n=0$. Thus, $\Lambda\in
\mathbf{L}_{n}=\mathbf{L}_{0}$ (since $n=0$), so that $\Lambda\in
\mathbf{L}_{0}=\left\{  \varnothing\right\}  $, so that $\Lambda=\varnothing$.
Now, there is only one map $g:\left[  n\right]  \rightarrow\mathcal{N}$ (since
$\left[  n\right]  =\left[  0\right]  =\varnothing$), and this map $g$ is
weakly increasing and satisfies $\operatorname*{FE}\left(  g\right)  =\left(
\Lambda\cup\left(  \Lambda+1\right)  \right)  \cap\left[  n\right]  $ (since
both $\operatorname*{FE}\left(  g\right)  $ and $\left(  \Lambda\cup\left(
\Lambda+1\right)  \right)  \cap\left[  n\right]  $ are the empty set). Hence,
there exists a weakly increasing map $g:\left[  n\right]  \rightarrow
\mathcal{N}$ such that $\operatorname*{FE}\left(  g\right)  =\left(
\Lambda\cup\left(  \Lambda+1\right)  \right)  \cap\left[  n\right]  $. In
other words, Lemma \ref{lem.FE.exist} holds, qed.}. Hence, for the rest of
this proof, we WLOG assume that we don't have $n=0$. Hence, $n$ is a positive
integer. Thus, $\mathbf{L}_{n}$ is the set of all nonempty lacunar subsets of
$\left[  n\right]  $ (by the definition of $\mathbf{L}_{n}$). Therefore, from
$\Lambda\in\mathbf{L}_{n}$, we conclude that $\Lambda$ is a nonempty lacunar
subset of $\left[  n\right]  $. Write this subset $\Lambda$ in the form
$\Lambda=\left\{  j_{1}<j_{2}<\cdots<j_{k}\right\}  $. Thus, $k\geq1$ (since
$\Lambda$ is nonempty).

Hence, $j_{1}\in\left\{  j_{1}<j_{2}<\cdots<j_{k}\right\}  =\Lambda
\subseteq\left[  n\right]  $, so that $j_{1}>0$. In other words, $0<j_{1}$.
Combining this with $j_{1}<j_{2}<\cdots<j_{k}$, we obtain $0<j_{1}%
<j_{2}<\cdots<j_{k}$.

From $\Lambda=\left\{  j_{1}<j_{2}<\cdots<j_{k}\right\}  $, we obtain
$\max\Lambda=j_{k}$. Thus, $j_{k}=\max\Lambda\leq n$ (since $\Lambda
\subseteq\left[  n\right]  $).

We have $\Lambda=\left\{  j_{1}<j_{2}<\cdots<j_{k}\right\}  =\left\{
j_{1},j_{2},\ldots,j_{k}\right\}  $ and thus%
\[
\Lambda+1=\left\{  j_{1}+1,j_{2}+1,\ldots,j_{k}+1\right\}
\]
(by the definition of $\Lambda+1$).

We are now in one of the following two cases:

\textit{Case 1:} We have $n\notin\Lambda$.

\textit{Case 2:} We have $n\in\Lambda$.

Let us consider Case 1. In this case, we have $n\notin\Lambda$. Hence,
$\Lambda\subseteq\left[  n-1\right]  $ (since $\Lambda\subseteq\left[
n\right]  $). Therefore, $\Lambda+1\subseteq\left\{  2,3,\ldots,n\right\}
\subseteq\left[  n\right]  $. Now,%
\[
\underbrace{\Lambda}_{\subseteq\left[  n\right]  }\cup\underbrace{\left(
\Lambda+1\right)  }_{\subseteq\left[  n\right]  }\subseteq\left[  n\right]
\cup\left[  n\right]  =\left[  n\right]  ,
\]
so that $\Lambda\cup\left(  \Lambda+1\right)  =\left(  \Lambda\cup\left(
\Lambda+1\right)  \right)  \cap\left[  n\right]  $.

But $j_{k}=\max\Lambda\leq n-1$ (since $\Lambda\subseteq\left[  n-1\right]
$). Thus, $j_{k}\leq n-1<n$. Combining this with $0<j_{1}<j_{2}<\cdots<j_{k}$,
we obtain $0<j_{1}<j_{2}<\cdots<j_{k}<n$.

Now, consider the map $g:\left[  n\right]  \rightarrow\mathcal{N}$ defined by%
\begin{align*}
g\left(  x\right)   &  =%
\begin{cases}
\left(  \text{the number of all }\lambda\in\Lambda\text{ such that }%
\lambda<x\right)  , & \text{if }x\leq j_{k};\\
\infty, & \text{if }x>j_{k}%
\end{cases}
\\
&  \ \ \ \ \ \ \ \ \ \ \text{for each }x\in\left[  n\right]  .
\end{align*}
Thus,%
\begin{align*}
&  \left(  g\left(  1\right)  ,g\left(  2\right)  ,\ldots,g\left(  n\right)
\right) \\
&  =\left(  \underbrace{0,0,\ldots,0}_{j_{1}\text{ entries}}%
,\underbrace{1,1,\ldots,1}_{j_{2}-j_{1}\text{ entries}},\underbrace{2,2,\ldots
,2}_{j_{3}-j_{2}\text{ entries}},\ldots,\underbrace{k-1,k-1,\ldots,k-1}%
_{j_{k}-j_{k-1}\text{ entries}},\underbrace{\infty,\infty,\ldots,\infty
}_{n-j_{k}\text{ entries}}\right)  .
\end{align*}
The $n$-tuple on the right hand side of this equality consists of a block of
$0$'s, followed by a block of $1$'s, followed by a block of $2$'s, and so on,
all the way up to a block of $\left(  k-1\right)  $'s, which is then followed
by a block of $\infty$'s. Each of these blocks is nonempty (since
$0<j_{1}<j_{2}<\cdots<j_{k}<n$).

Hence, the map $g$ is weakly increasing, and its nonempty fibers are
\begin{align*}
g^{-1}\left(  0\right)   &  =\left\{  1,2,\ldots,j_{1}\right\}  ,\\
g^{-1}\left(  1\right)   &  =\left\{  j_{1}+1,j_{1}+2,\ldots,j_{2}\right\}
,\\
g^{-1}\left(  2\right)   &  =\left\{  j_{2}+1,j_{2}+2,\ldots,j_{3}\right\}
,\\
&  \vdots,\\
g^{-1}\left(  k-1\right)   &  =\left\{  j_{k-1}+1,j_{k-1}+2,\ldots
,j_{k}\right\}  ,\\
g^{-1}\left(  \infty\right)   &  =\left\{  j_{k}+1,j_{k}+2,\ldots,n\right\}  .
\end{align*}
Thus, the definition of $\operatorname*{FE}\left(  g\right)  $ yields%
\begin{align*}
&  \operatorname*{FE}\left(  g\right) \\
&  =\underbrace{\left\{  \min\left(  g^{-1}\left(  h\right)  \right)
\ \mid\ h\in\left\{  1,2,3,\ldots,\infty\right\}  \text{ with }g^{-1}\left(
h\right)  \neq\varnothing\right\}  }_{=\left\{  \min\left(  g^{-1}\left(
1\right)  \right)  ,\min\left(  g^{-1}\left(  2\right)  \right)  ,\ldots
,\min\left(  g^{-1}\left(  k-1\right)  \right)  ,\min\left(  g^{-1}\left(
\infty\right)  \right)  \right\}  }\\
&  \ \ \ \ \ \ \ \ \ \ \cup\underbrace{\left\{  \max\left(  g^{-1}\left(
h\right)  \right)  \ \mid\ h\in\left\{  0,1,2,3,\ldots\right\}  \text{ with
}g^{-1}\left(  h\right)  \neq\varnothing\right\}  }_{=\left\{  \max\left(
g^{-1}\left(  0\right)  \right)  ,\max\left(  g^{-1}\left(  1\right)  \right)
,\max\left(  g^{-1}\left(  2\right)  \right)  ,\ldots,\max\left(
g^{-1}\left(  k-1\right)  \right)  \right\}  }\\
&  =\underbrace{\left\{  \min\left(  g^{-1}\left(  1\right)  \right)
,\min\left(  g^{-1}\left(  2\right)  \right)  ,\ldots,\min\left(
g^{-1}\left(  k-1\right)  \right)  ,\min\left(  g^{-1}\left(  \infty\right)
\right)  \right\}  }_{\substack{=\left\{  j_{1}+1,j_{2}+1,\ldots
,j_{k-1}+1,j_{k}+1\right\}  \\\text{(since }\min\left(  g^{-1}\left(
h\right)  \right)  =j_{h}+1\text{ for each }h\in\left\{  1,2,\ldots
,k-1\right\}  \text{,}\\\text{and since }\min\left(  g^{-1}\left(
\infty\right)  \right)  =j_{k}+1\text{)}}}\\
&  \ \ \ \ \ \ \ \ \ \ \cup\underbrace{\left\{  \max\left(  g^{-1}\left(
0\right)  \right)  ,\max\left(  g^{-1}\left(  1\right)  \right)  ,\max\left(
g^{-1}\left(  2\right)  \right)  ,\ldots,\max\left(  g^{-1}\left(  k-1\right)
\right)  \right\}  }_{\substack{=\left\{  j_{1},j_{2},\ldots,j_{k}\right\}
\\\text{(since }\max\left(  g^{-1}\left(  h-1\right)  \right)  =j_{h}\text{
for each }h\in\left\{  1,2,\ldots,k\right\}  \text{)}}}\\
&  =\underbrace{\left\{  j_{1}+1,j_{2}+1,\ldots,j_{k-1}+1,j_{k}+1\right\}
}_{\substack{=\left\{  j_{1}+1,j_{2}+1,\ldots,j_{k}+1\right\}  \\=\Lambda
+1}}\cup\underbrace{\left\{  j_{1},j_{2},\ldots,j_{k}\right\}  }_{=\Lambda
}=\left(  \Lambda+1\right)  \cup\Lambda\\
&  =\Lambda\cup\left(  \Lambda+1\right)  =\left(  \Lambda\cup\left(
\Lambda+1\right)  \right)  \cap\left[  n\right]  .
\end{align*}

Altogether, we have now shown that our map $g:\left[  n\right]  \rightarrow
\mathcal{N}$ is weakly increasing and satisfies $\operatorname*{FE}\left(
g\right)  =\left(  \Lambda\cup\left(  \Lambda+1\right)  \right)  \cap\left[
n\right]  $. Hence, such a map $g$ exists. Thus, Lemma \ref{lem.FE.exist} is
proven in Case 1.

Let us consider Case 2. In this case, we have $n\in\Lambda$. Hence, $n\leq
\max\Lambda=j_{k}$. Combined with $j_{k}\leq n$, this yields $n=j_{k}$.

We have $\Lambda\subseteq\left[  n\right]  $ and thus $\Lambda\cap\left[
n\right]  =\Lambda$. Also,
\begin{equation}
\left(  \Lambda+1\right)  \cap\left[  n\right]  =\left\{  j_{1}+1,j_{2}%
+1,\ldots,j_{k-1}+1\right\}  \label{pf.lem.FE.exist.c2.1}%
\end{equation}
\footnote{\textit{Proof.} We have $j_{1}<j_{2}<\cdots<j_{k}$ and thus
$j_{1}+1<j_{2}+1<\cdots<j_{k}+1$. In other words, $j_{1}+1<j_{2}%
+1<\cdots<j_{k-1}+1<j_{k}+1$. In view of $j_{k}=n$, this rewrites as
$j_{1}+1<j_{2}+1<\cdots<j_{k-1}+1<n+1$. In other words, $j_{1}+1<j_{2}%
+1<\cdots<j_{k-1}+1\leq n$ (since an integer is $<n+1$ if and only if it is
$\leq n$). Therefore, the $k-1$ numbers $j_{1}+1,j_{2}+1,\ldots,j_{k-1}+1$ all
belong to the set $\left[  n\right]  $. In other words, $\left\{
j_{1}+1,j_{2}+1,\ldots,j_{k-1}+1\right\}  \subseteq\left[  n\right]  $.
\par
On the other hand, $\underbrace{j_{k}}_{=n}+1=n+1\notin\left[  n\right]  $ and
thus $\left\{  j_{k}+1\right\}  \cap\left[  n\right]  =\varnothing$. Now,%
\begin{align*}
\underbrace{\left(  \Lambda+1\right)  }_{\substack{=\left\{  j_{1}%
+1,j_{2}+1,\ldots,j_{k}+1\right\}  \\=\left\{  j_{1}+1,j_{2}+1,\ldots
,j_{k-1}+1\right\}  \cup\left\{  j_{k}+1\right\}  }}\cap\left[  n\right]   &
=\left(  \left\{  j_{1}+1,j_{2}+1,\ldots,j_{k-1}+1\right\}  \cup\left\{
j_{k}+1\right\}  \right)  \cap\left[  n\right] \\
&  =\underbrace{\left(  \left\{  j_{1}+1,j_{2}+1,\ldots,j_{k-1}+1\right\}
\cap\left[  n\right]  \right)  }_{\substack{=\left\{  j_{1}+1,j_{2}%
+1,\ldots,j_{k-1}+1\right\}  \\\text{(since }\left\{  j_{1}+1,j_{2}%
+1,\ldots,j_{k-1}+1\right\}  \subseteq\left[  n\right]  \text{)}}%
}\cup\underbrace{\left(  \left\{  j_{k}+1\right\}  \cap\left[  n\right]
\right)  }_{=\varnothing}\\
&  =\left\{  j_{1}+1,j_{2}+1,\ldots,j_{k-1}+1\right\}  \cup\varnothing\\
&  =\left\{  j_{1}+1,j_{2}+1,\ldots,j_{k-1}+1\right\}  .
\end{align*}
}.

Now, consider the map $g:\left[  n\right]  \rightarrow\mathcal{N}$ defined by%
\begin{align*}
g\left(  x\right)   &  =\left(  \text{the number of all }\lambda\in
\Lambda\text{ such that }\lambda<x\right) \\
&  \ \ \ \ \ \ \ \ \ \ \text{for each }x\in\left[  n\right]  .
\end{align*}
Thus,%
\begin{align*}
&  \left(  g\left(  1\right)  ,g\left(  2\right)  ,\ldots,g\left(  n\right)
\right) \\
&  =\left(  \underbrace{0,0,\ldots,0}_{j_{1}\text{ entries}}%
,\underbrace{1,1,\ldots,1}_{j_{2}-j_{1}\text{ entries}},\underbrace{2,2,\ldots
,2}_{j_{3}-j_{2}\text{ entries}},\ldots,\underbrace{k-1,k-1,\ldots,k-1}%
_{j_{k}-j_{k-1}\text{ entries}}\right)  .
\end{align*}
\footnote{There are no entries beyond the $\underbrace{k-1,k-1,\ldots
,k-1}_{j_{k}-j_{k-1}\text{ entries}}$ block in this tuple, because $j_{k}=n$.}
The $n$-tuple on the right hand side of this equality consists of a block of
$0$'s, followed by a block of $1$'s, followed by a block of $2$'s, and so on,
all the way up to a block of $\left(  k-1\right)  $'s. Each of these blocks is
nonempty (since $0<j_{1}<j_{2}<\cdots<j_{k}$).

Hence, the map $g$ is weakly increasing, and its nonempty fibers are
\begin{align*}
g^{-1}\left(  0\right)   &  =\left\{  1,2,\ldots,j_{1}\right\}  ,\\
g^{-1}\left(  1\right)   &  =\left\{  j_{1}+1,j_{1}+2,\ldots,j_{2}\right\}
,\\
g^{-1}\left(  2\right)   &  =\left\{  j_{2}+1,j_{2}+2,\ldots,j_{3}\right\}
,\\
&  \vdots,\\
g^{-1}\left(  k-1\right)   &  =\left\{  j_{k-1}+1,j_{k-1}+2,\ldots
,j_{k}\right\}  .
\end{align*}
Thus, the definition of $\operatorname*{FE}\left(  g\right)  $ yields%
\begin{align*}
&  \operatorname*{FE}\left(  g\right) \\
&  =\underbrace{\left\{  \min\left(  g^{-1}\left(  h\right)  \right)
\ \mid\ h\in\left\{  1,2,3,\ldots,\infty\right\}  \text{ with }g^{-1}\left(
h\right)  \neq\varnothing\right\}  }_{=\left\{  \min\left(  g^{-1}\left(
1\right)  \right)  ,\min\left(  g^{-1}\left(  2\right)  \right)  ,\ldots
,\min\left(  g^{-1}\left(  k-1\right)  \right)  \right\}  }\\
&  \ \ \ \ \ \ \ \ \ \ \cup\underbrace{\left\{  \max\left(  g^{-1}\left(
h\right)  \right)  \ \mid\ h\in\left\{  0,1,2,3,\ldots\right\}  \text{ with
}g^{-1}\left(  h\right)  \neq\varnothing\right\}  }_{=\left\{  \max\left(
g^{-1}\left(  0\right)  \right)  ,\max\left(  g^{-1}\left(  1\right)  \right)
,\max\left(  g^{-1}\left(  2\right)  \right)  ,\ldots,\max\left(
g^{-1}\left(  k-1\right)  \right)  \right\}  }\\
&  =\underbrace{\left\{  \min\left(  g^{-1}\left(  1\right)  \right)
,\min\left(  g^{-1}\left(  2\right)  \right)  ,\ldots,\min\left(
g^{-1}\left(  k-1\right)  \right)  \right\}  }_{\substack{=\left\{
j_{1}+1,j_{2}+1,\ldots,j_{k-1}+1\right\}  \\\text{(since }\min\left(
g^{-1}\left(  h\right)  \right)  =j_{h}+1\text{ for each }h\in\left\{
1,2,\ldots,k-1\right\}  \text{)}}}\\
&  \ \ \ \ \ \ \ \ \ \ \cup\underbrace{\left\{  \max\left(  g^{-1}\left(
0\right)  \right)  ,\max\left(  g^{-1}\left(  1\right)  \right)  ,\max\left(
g^{-1}\left(  2\right)  \right)  ,\ldots,\max\left(  g^{-1}\left(  k-1\right)
\right)  \right\}  }_{\substack{=\left\{  j_{1},j_{2},\ldots,j_{k}\right\}
\\\text{(since }\max\left(  g^{-1}\left(  h-1\right)  \right)  =j_{h}\text{
for each }h\in\left\{  1,2,\ldots,k\right\}  \text{)}}}\\
&  =\underbrace{\left\{  j_{1}+1,j_{2}+1,\ldots,j_{k-1}+1\right\}
}_{\substack{=\left(  \Lambda+1\right)  \cap\left[  n\right]  \\\text{(by
(\ref{pf.lem.FE.exist.c2.1}))}}}\cup\underbrace{\left\{  j_{1},j_{2}%
,\ldots,j_{k}\right\}  }_{\substack{=\Lambda=\Lambda\cap\left[  n\right]
\\\text{(since }\Lambda\subseteq\left[  n\right]  \text{)}}}=\left(  \left(
\Lambda+1\right)  \cap\left[  n\right]  \right)  \cup\left(  \Lambda
\cap\left[  n\right]  \right) \\
&  =\left(  \Lambda\cap\left[  n\right]  \right)  \cup\left(  \left(
\Lambda+1\right)  \cap\left[  n\right]  \right)  =\left(  \Lambda\cup\left(
\Lambda+1\right)  \right)  \cap\left[  n\right]  .
\end{align*}

Altogether, we have now shown that our map $g:\left[  n\right]  \rightarrow
\mathcal{N}$ is weakly increasing and satisfies $\operatorname*{FE}\left(
g\right)  =\left(  \Lambda\cup\left(  \Lambda+1\right)  \right)  \cap\left[
n\right]  $. Hence, such a map $g$ exists. Thus, Lemma \ref{lem.FE.exist} is
proven in Case 2.

We have now proven Lemma \ref{lem.FE.exist} in each of the two Cases 1 and 2.
Since these two Cases cover all possibilities, we thus conclude that Lemma
\ref{lem.FE.exist} always holds.
\end{proof}
\end{verlong}

\begin{verlong}
\begin{lemma}
\label{lem.V-shaped.no-peaks}Let $n\in\mathbb{N}$. Let $\pi$ be an
$n$-permutation. Let $S$ be a nonempty interval of the totally ordered set
$\left[  n\right]  $ such that the map $\pi\mid_{S}$ is V-shaped. Then,
$S\cap\operatorname{Epk}\pi\subseteq\left\{  \min S,\max S\right\}  $.
\end{lemma}

\begin{proof}
[Proof of Lemma \ref{lem.V-shaped.no-peaks}.]Let $j\in S\cap
\operatorname{Epk}\pi$. We are going to show that $j\in\left\{  \min S,\max
S\right\}  $.

Indeed, assume the contrary. Thus, $j\notin\left\{  \min S,\max S\right\}  $.
Hence, $j\neq\min S$ and $j\neq\max S$.

But $j\in S\cap\operatorname{Epk}\pi\subseteq S$. Thus, $j\geq\min S$ and
$j\leq\max S$. Combining $j\geq\min S$ with $j\neq\min S$, we obtain $j>\min
S$. Combining $j\leq\max S$ with $j\neq\max S$, we obtain $j<\max S$.

Recall the following basic fact: If $T$ is a nonempty interval of $\left[
n\right]  $, and if $t\in T$ satisfies $t>\min T$, then $t-1\in T$. Applying
this to $T=S$ and $t=j$, we conclude that $j-1\in S$ (since $j>\min S$).

Recall the following basic fact: If $T$ is a nonempty interval of $\left[
n\right]  $, and if $t\in T$ satisfies $t<\max T$, then $t+1\in T$. Applying
this to $T=S$ and $t=j$, we conclude that $j+1\in S$ (since $j<\max S$).

Recall that the map $\pi\mid_{S}$ is V-shaped. According to the definition of
\textquotedblleft V-shaped\textquotedblright, this means the following: There
exists some $t\in S$ such that the map $\left(  \pi\mid_{S}\right)
\mid_{\left\{  s\in S\ \mid\ s\leq t\right\}  }$ is strictly decreasing while
the map $\left(  \pi\mid_{S}\right)  \mid_{\left\{  s\in S\ \mid\ s\geq
t\right\}  }$ is strictly increasing. Consider this $t$.

The map $\left(  \pi\mid_{S}\right)  \mid_{\left\{  s\in S\ \mid\ s\leq
t\right\}  }$ is strictly decreasing. In other words, the map $\pi
\mid_{\left\{  s\in S\ \mid\ s\leq t\right\}  }$ is strictly decreasing (since
$\left(  \pi\mid_{S}\right)  \mid_{\left\{  s\in S\ \mid\ s\leq t\right\}
}=\pi\mid_{\left\{  s\in S\ \mid\ s\leq t\right\}  }$).

The map $\left(  \pi\mid_{S}\right)  \mid_{\left\{  s\in S\ \mid\ s\geq
t\right\}  }$ is strictly increasing. In other words, the map $\pi
\mid_{\left\{  s\in S\ \mid\ s\geq t\right\}  }$ is strictly increasing (since
$\left(  \pi\mid_{S}\right)  \mid_{\left\{  s\in S\ \mid\ s\geq t\right\}
}=\pi\mid_{\left\{  s\in S\ \mid\ s\geq t\right\}  }$).

We have $j\in S\cap\operatorname{Epk}\pi\subseteq\operatorname{Epk}\pi$; in
other words, $j$ is an exterior peak of $\pi$ (by the definition of
$\operatorname{Epk}\pi$). In other words, $j$ is an element of $\left[
n\right]  $ and satisfies $\pi_{j-1}<\pi_{j}>\pi_{j+1}$ (by the definition of
an exterior peak).

We are in one of the following two cases:

\textit{Case 1:} We have $j<t$.

\textit{Case 2:} We have $j\geq t$.

Let us first consider Case 1. In this case, we have $j<t$. Hence, $j\leq t$
and $j-1\leq j\leq t$. Thus, $j-1\in\left\{  s\in S\ \mid\ s\leq t\right\}  $
(since $j-1\in S$). Also, $j\leq t$. Hence, $j\in\left\{  s\in S\ \mid\ s\leq
t\right\}  $ (since $j\in S$).

Now, we know that $j-1$ and $j$ are two elements of the set $\left\{  s\in
S\ \mid\ s\leq t\right\}  $, and satisfy $j-1<j$. Hence,
\[
\left(  \pi\mid_{\left\{  s\in S\ \mid\ s\leq t\right\}  }\right)  \left(
j-1\right)  >\left(  \pi\mid_{\left\{  s\in S\ \mid\ s\leq t\right\}
}\right)  \left(  j\right)
\]
(since the map $\pi\mid_{\left\{  s\in S\ \mid\ s\leq t\right\}  }$ is
strictly decreasing). Thus,%
\[
\pi_{j-1}=\pi\left(  j-1\right)  =\left(  \pi\mid_{\left\{  s\in
S\ \mid\ s\leq t\right\}  }\right)  \left(  j-1\right)  >\left(  \pi
\mid_{\left\{  s\in S\ \mid\ s\leq t\right\}  }\right)  \left(  j\right)
=\pi\left(  j\right)  =\pi_{j}.
\]
This contradicts $\pi_{j-1}<\pi_{j}$. Thus, we have obtained a contradiction
in Case 1.

Let us now consider Case 2. In this case, we have $j\geq t$. Hence, $j+1\geq
j\geq t$. Thus, $j+1\in\left\{  s\in S\ \mid\ s\geq t\right\}  $ (since
$j+1\in S$). Also, $j\geq t$. Hence, $j\in\left\{  s\in S\ \mid\ s\geq
t\right\}  $ (since $j\in S$).

Now, we know that $j$ and $j+1$ are two elements of the set $\left\{  s\in
S\ \mid\ s\geq t\right\}  $, and satisfy $j<j+1$. Hence,
\[
\left(  \pi\mid_{\left\{  s\in S\ \mid\ s\geq t\right\}  }\right)  \left(
j\right)  <\left(  \pi\mid_{\left\{  s\in S\ \mid\ s\geq t\right\}  }\right)
\left(  j+1\right)
\]
(since the map $\pi\mid_{\left\{  s\in S\ \mid\ s\geq t\right\}  }$ is
strictly increasing). Thus,%
\[
\pi_{j}=\pi\left(  j\right)  =\left(  \pi\mid_{\left\{  s\in S\ \mid\ s\geq
t\right\}  }\right)  \left(  j\right)  <\left(  \pi\mid_{\left\{  s\in
S\ \mid\ s\geq t\right\}  }\right)  \left(  j+1\right)  =\pi\left(
j+1\right)  =\pi_{j+1}.
\]
This contradicts $\pi_{j}>\pi_{j+1}$. Thus, we have obtained a contradiction
in Case 2.

We have now obtained a contradiction in each of the two Cases 1 and 2. Since
these two Cases cover all possibilities, we thus conclude that we always have
a contradiction. This contradiction shows that our assumption was wrong. Thus,
$j\in\left\{  \min S,\max S\right\}  $ is proven.

Now, forget that we fixed $j$. We thus have shown that $j\in\left\{  \min
S,\max S\right\}  $ for each $j\in S\cap\operatorname{Epk}\pi$. In other
words, $S\cap\operatorname{Epk}\pi\subseteq\left\{  \min S,\max S\right\}  $.
This proves Lemma \ref{lem.V-shaped.no-peaks}.
\end{proof}

\begin{lemma}
\label{lem.Epk.V-sh}Let $n\in\mathbb{N}$. Let $\pi$ be an $n$-permutation. Let
$S$ be a nonempty interval of the totally ordered set $\left[  n\right]  $.

\textbf{(a)} If $1\in S$ and $S\cap\operatorname{Epk}\pi\subseteq\left\{
\max S\right\}  $, then the map $\pi\mid_{S}$ is strictly increasing.

\textbf{(b)} If $n\in S$ and $S\cap\operatorname{Epk}\pi\subseteq\left\{
\min S\right\}  $, then the map $\pi\mid_{S}$ is strictly decreasing.

\textbf{(c)} If $S\cap\operatorname{Epk}\pi\subseteq\left\{  \min S,\max
S\right\}  $, then the map $\pi\mid_{S}$ is V-shaped.
\end{lemma}

\begin{proof}
[Proof of Lemma \ref{lem.Epk.V-sh}.]Set $\pi_{0}=0$ and $\pi_{n+1}=0$. Thus,
$\pi_{i}\in\mathbb{N}$ is well-defined for each $i\in\left\{  0,1,\ldots
,n+1\right\}  $. The exterior peaks of $\pi$ are the elements $i\in\left[
n\right]  $ satisfying $\pi_{i-1}<\pi_{i}>\pi_{i+1}$ (by the definition of an
\textquotedblleft exterior peak\textquotedblright).

Note that $\pi$ is an $n$-permutation; thus, the entries of the word $\pi$ are
distinct. (But of course, $\pi_{0}$ and $\pi_{n+1}$ don't count as entries of
$\pi$.)

We now claim the following:

\begin{statement}
\textit{Observation 1:} Let $i$ and $j$ be two elements of $\left[  n\right]
$ satisfying $i\leq j+1$ and $\pi_{i-1}\leq\pi_{i}$ and $\pi_{j}\geq\pi_{j+1}%
$. Then,%
\[
\left\{  i,i+1,\ldots,j\right\}  \cap\operatorname{Epk}\pi\neq\varnothing.
\]

\end{statement}

[\textit{Proof of Observation 1:} We have $\pi_{i-1}\neq\pi_{i}$%
\ \ \ \ \footnote{\textit{Proof.} Assume the contrary. Thus, $\pi_{i-1}%
=\pi_{i}$. If we have $i=1$, then $\pi_{i-1}=\pi_{1-1}=\pi_{0}=0\notin%
\mathbb{P}$, which contradicts $\pi_{i-1}=\pi_{i}\in\mathbb{P}$ (since
$i\in\left[  n\right]  $). Thus, we cannot have $i=1$. Hence, we have $i\neq1$
and therefore $i>1$ (since $i\in\left[  n\right]  $). Thus, $i-1\in\left[
n\right]  $. Hence, both $i-1$ and $i$ are elements of $\left[  n\right]  $.
\par
Now, $\pi_{i-1}$ and $\pi_{i}$ are two distinct entries of the word $\pi$
(since $i-1$ and $i$ are two distinct elements of $\left[  n\right]  $). Thus,
$\pi_{i-1}\neq\pi_{i}$ (since the entries of $\pi$ are distinct). This
contradicts $\pi_{i-1}=\pi_{i}$. This contradiction shows that our assumption
was false. Qed.}.

Combining $\pi_{i-1}\leq\pi_{i}$ with $\pi_{i-1}\neq\pi_{i}$, we obtain
$\pi_{i-1}<\pi_{i}$. If we had $i=j+1$, then we would have $j=i-1$ and thus
$\pi_{j}=\pi_{i-1}<\pi_{i}=\pi_{j+1}$ (since $i=j+1$); but this would
contradict $\pi_{j}\geq\pi_{j+1}$. Hence, we cannot have $i=j+1$. Thus, $i\neq
j+1$. Combining this with $i\leq j+1$, we obtain $i<j+1$. Thus, $i\leq j$
(since $i$ and $j$ are integers).

Now, $i\in\left\{  i,i+1,\ldots,j\right\}  $ (since $i\leq j$) and $\pi
_{i-1}<\pi_{i}$. Hence, $i$ is an element $k\in\left\{  i,i+1,\ldots
,j\right\}  $ satisfying $\pi_{k-1}<\pi_{k}$. Hence, there exists at least one
element $k\in\left\{  i,i+1,\ldots,j\right\}  $ satisfying $\pi_{k-1}<\pi_{k}$
(namely, $i$). Let $p$ be the \textbf{largest} such element. Thus, $p$ itself
is an element of $\left\{  i,i+1,\ldots,j\right\}  $ and satisfies $\pi
_{p-1}<\pi_{p}$, and furthermore, every element $k\in\left\{  i,i+1,\ldots
,j\right\}  $ satisfying $\pi_{k-1}<\pi_{k}$ must satisfy%
\begin{equation}
k\leq p. \label{pf.lem.Epk.V-sh.o1.pf.max}%
\end{equation}

We have $p\in\left\{  i,i+1,\ldots,j\right\}  \subseteq\left[  n\right]  $
(since $i$ and $j$ are elements of $\left[  n\right]  $), so that $\pi_{p}%
\in\mathbb{P}$. Thus, $\pi_{p}>0$.

Next, we claim that $\pi_{p}>\pi_{p+1}$. Indeed, assume the contrary. Hence,
$\pi_{p}\leq\pi_{p+1}$. But recall that $\pi_{p}>0$. Hence, $p\neq
n$\ \ \ \ \footnote{\textit{Proof.} Assume the contrary. Thus, $p=n$, so that
$\pi_{p+1}=\pi_{n+1}=0$. Hence, $\pi_{p}\leq\pi_{p+1}=0$. But this contradicts
$\pi_{p}>0$. This contradiction shows that our assumption was false. Qed.}.
Combining $p\in\left[  n\right]  $ with $p\neq n$, we obtain $p\in\left[
n\right]  \setminus\left\{  n\right\}  \subseteq\left[  n-1\right]  $;
therefore, $p+1\in\left[  n\right]  $. Now we know that $p$ and $p+1$ both are
elements of $\left[  n\right]  $. Hence, $\pi_{p}$ and $\pi_{p+1}$ are two
distinct entries of the word $\pi$ (since $p\neq p+1$). Hence, $\pi_{p}\neq
\pi_{p+1}$ (since the entries of the word $\pi$ are distinct). Combining this
with $\pi_{p}\leq\pi_{p+1}$, we obtain $\pi_{p}<\pi_{p+1}$. If we had $p=j$,
then we would have $\pi_{p}=\pi_{j}\geq\pi_{j+1}=\pi_{p+1}$ (since $j=p$),
which would contradict $\pi_{p}<\pi_{p+1}$. Thus, we cannot have $p=j$. In
other words, we have $p\neq j$. Combining $p\in\left\{  i,i+1,\ldots
,j\right\}  $ with $p\neq j$, we obtain $p\in\left\{  i,i+1,\ldots,j\right\}
\setminus\left\{  j\right\}  =\left\{  i,i+1,\ldots,j-1\right\}  $, so that
\[
p+1\in\left\{  i+1,i+2,\ldots,j\right\}  \subseteq\left\{  i,i+1,\ldots
,j\right\}  .
\]
Hence, $p+1$ is an element $k\in\left\{  i,i+1,\ldots,j\right\}  $ satisfying
$\pi_{k-1}<\pi_{k}$ (since $\pi_{\left(  p+1\right)  -1}=\pi_{p}<\pi_{p+1}$).
Therefore, (\ref{pf.lem.Epk.V-sh.o1.pf.max}) (applied to $k=p+1$) yields
$p+1\leq p$. But this contradicts $p<p+1$. This contradiction shows that our
assumption was false. Hence, $\pi_{p}>\pi_{p+1}$ is proven.

Now, we know that $p\in\left[  n\right]  $ and $\pi_{p-1}<\pi_{p}>\pi_{p+1}$.
In other words, $p$ is an exterior peak of $\pi$ (by the definition of an
exterior peak). In other words, $p\in\operatorname{Epk}\pi$. Combining this
with $p\in\left\{  i,i+1,\ldots,j\right\}  $, we obtain $p\in\left\{
i,i+1,\ldots,j\right\}  \cap\operatorname{Epk}\pi$. Hence, the set $\left\{
i,i+1,\ldots,j\right\}  \cap\operatorname{Epk}\pi$ has at least one element
(namely, $p$). Therefore, this set is nonempty. In other words, $\left\{
i,i+1,\ldots,j\right\}  \cap\operatorname{Epk}\pi\neq\varnothing$. This
proves Observation 1.]

\textbf{(c)} Assume that $S\cap\operatorname{Epk}\pi\subseteq\left\{  \min
S,\max S\right\}  $. We must prove that the map $\pi\mid_{S}$ is V-shaped.

Write the interval $S$ in the form $S=\left\{  a,a+1,\ldots,b\right\}  $ for
some elements $a$ and $b$ of $\left[  n\right]  $. Thus, $a\leq b$ (since the
interval $S$ is nonempty) and $\min S=a$ and $\max S=b$. Thus, $S\cap
\operatorname{Epk}\pi\subseteq\left\{  \underbrace{\min S}_{=a}%
,\underbrace{\max S}_{=b}\right\}  =\left\{  a,b\right\}  $.

The set $S$ is nonempty. Hence, there exists some $p\in S$ that minimizes
$\pi_{p}$. Consider this $p$. Thus,%
\begin{equation}
\pi_{p}\leq\pi_{s}\ \ \ \ \ \ \ \ \ \ \text{for each }s\in S.
\label{pf.lem.Epk.V-sh.c.min}%
\end{equation}
Note that $p\in S=\left\{  a,a+1,\ldots,b\right\}  $, so that $a\leq p\leq b$.
Also, $p\in S\subseteq\left[  n\right]  $.

From $S=\left\{  a,a+1,\ldots,b\right\}  $, we obtain%
\[
\left\{  s\in S\ \mid\ s\leq p\right\}  =\left\{  s\in\left\{  a,a+1,\ldots
,b\right\}  \ \mid\ s\leq p\right\}  =\left\{  a,a+1,\ldots,p\right\}
\]
(since $p\leq b$). Also, from $S=\left\{  a,a+1,\ldots,b\right\}  $, we obtain%
\[
\left\{  s\in S\ \mid\ s\geq p\right\}  =\left\{  s\in\left\{  a,a+1,\ldots
,b\right\}  \ \mid\ s\geq p\right\}  =\left\{  p,p+1,\ldots,b\right\}
\]
(since $p \geq a$).

Now, we claim that the map $\pi\mid_{\left\{  s\in S\ \mid\ s\leq p\right\}
}$ is strictly decreasing.

[\textit{Proof:} Let $i\in\left\{  a+1,a+2,\ldots,p\right\}  $. Thus,
\begin{align*}
i-1  &  \in\left\{  a,a+1,\ldots,p-1\right\}  \subseteq\left\{  a,a+1,\ldots
,b\right\}  \ \ \ \ \ \ \ \ \ \ \left(  \text{since }p-1<p\leq b\right) \\
&  =S\subseteq\left[  n\right]  .
\end{align*}
Hence, $\pi_{i-1}$ is well-defined. Also, $i\in\left\{  a+1,a+2,\ldots
,p\right\}  \subseteq\left\{  a,a+1,\ldots,b\right\}  $ (since $a+1>a\geq a$
and $p\leq b$), so that $i\in\left\{  a,a+1,\ldots,b\right\}  =S\subseteq
\left[  n\right]  $. Thus, $\pi_{i}$ is well-defined. Also, from $i\in\left\{
a,a+1,\ldots,b\right\}  $, we obtain $i\geq a$.

Now, assume (for the sake of contradiction) that $\pi_{i-1}\leq\pi_{i}$.

But
$i\in\left\{  a+1,a+2,\ldots,p\right\}  $ yields $a+1\leq i\leq p$; hence,
$a\leq p-1\leq b$ (since $p-1<p\leq b$). Thus, $p-1\in\left\{  a,a+1,\ldots
,b\right\}  =S$. Hence, (\ref{pf.lem.Epk.V-sh.c.min}) (applied to $s=p-1$)
yields $\pi_{p}\leq\pi_{p-1}$. Thus, $\pi_{p-1}\geq\pi_{p}=\pi_{\left(
p-1\right)  +1}$.

Also, $i\in\left[  n\right]  $ and $p-1\in S\subseteq\left[  n\right]  $ and
$i\leq p=\left(  p-1\right)  +1$. Hence, Observation 1 (applied to $j=p-1$)
yields $\left\{  i,i+1,\ldots,p-1\right\}  \cap\operatorname{Epk}\pi
\neq\varnothing$. In other words, the set $\left\{  i,i+1,\ldots,p-1\right\}
\cap\operatorname{Epk}\pi$ has at least one element. Consider such an
element, and denote it by $z$. Thus,%
\[
z\in\left\{  i,i+1,\ldots,p-1\right\}  \cap\operatorname{Epk}\pi.
\]

Now,%
\begin{align*}
z  &  \in\underbrace{\left\{  i,i+1,\ldots,p-1\right\}  }_{\substack{\subseteq
\left\{  a,a+1,\ldots,b\right\}  \\\text{(since }i\geq a\text{ and }p-1\leq
b\text{)}}}\cap\operatorname{Epk}\pi
\subseteq
\underbrace{\left\{  a,a+1,\ldots
,b\right\}  }_{=S}\cap\operatorname{Epk}\pi\\
&  =S\cap\operatorname{Epk}\pi
\subseteq \left\{  a,b\right\}  .
\end{align*}
But $z\in\left\{  i,i+1,\ldots,p-1\right\}  \cap\operatorname{Epk} \pi
\subseteq\left\{  i,i+1,\ldots,p-1\right\}  $, so that $z\leq p-1<p\leq b$.
Hence, $z\neq b$. Combining $z\in\left\{  a,b\right\}  $ with $z\neq b$, we
obtain $z\in\left\{  a,b\right\}  \setminus\left\{  b\right\}  \subseteq
\left\{  a\right\}  $. In other words, $z=a$. But from $z\in\left\{
i,i+1,\ldots,p-1\right\}  $, we also obtain $z\geq i\geq a+1$ (since
$i\in\left\{  a+1,a+2,\ldots,p\right\}  $). Hence, $z\geq a+1>a$; this
contradicts $z=a$.

This contradiction shows that our assumption (that $\pi_{i-1}\leq\pi_{i}$) was
wrong. Hence, we don't have $\pi_{i-1}\leq\pi_{i}$. In other words, we have
$\pi_{i-1}>\pi_{i}$.

Now, forget that we fixed $i$. We thus have shown that $\pi_{i-1}>\pi_{i}$ for
each $i\in\left\{  a+1,a+2,\ldots,p\right\}  $. In other words, $\pi_{a}%
>\pi_{a+1}>\pi_{a+2}>\cdots>\pi_{p}$. In other words, the map $\pi
\mid_{\left\{  a,a+1,\ldots,p\right\}  }$ is strictly decreasing. In other
words, the map $\pi\mid_{\left\{  s\in S\ \mid\ s\leq p\right\}  }$ is
strictly decreasing (since $\left\{  s\in S\ \mid\ s\leq p\right\}  =\left\{
a,a+1,\ldots,p\right\}  $).]

Next, we claim that the map $\pi\mid_{\left\{  s\in S\ \mid\ s\geq p\right\}
}$ is strictly increasing.

[\textit{Proof:} Let $j\in\left\{  p,p+1,\ldots,b-1\right\}  $. Thus,%
\begin{align*}
j+1  &  \in\left\{  p+1,p+2,\ldots,b\right\}  \subseteq\left\{  a,a+1,\ldots
,b\right\}  \ \ \ \ \ \ \ \ \ \ \left(  \text{since }p+1>p\geq a\right) \\
&  =S\subseteq\left[  n\right]  .
\end{align*}
Hence, $\pi_{j+1}$ is well-defined. Also, $j\in\left\{  p,p+1,\ldots
,b-1\right\}  \subseteq\left\{  a,a+1,\ldots,b\right\}  $ (since $p\geq a$ and
$b-1<b$), so that $j\in\left\{  a,a+1,\ldots,b\right\}  =S\subseteq\left[
n\right]  $. Thus, $\pi_{j}$ is well-defined. Also, from $j\in\left\{
a,a+1,\ldots,b\right\}  $, we obtain $j\geq a$ and $j\leq b$.

Now, assume (for the sake of contradiction) that $\pi_{j}\geq\pi_{j+1}$.

But
$j\in\left\{  p,p+1,\ldots,b-1\right\}  $ yields $p\leq j\leq b-1$; hence,
$p+1\leq b$. Combined with $a\leq p<p+1$, this yields $a\leq p+1\leq b$. Thus,
$p+1\in\left\{  a,a+1,\ldots,b\right\}  =S$. Hence,
(\ref{pf.lem.Epk.V-sh.c.min}) (applied to $s=p+1$) yields $\pi_{p}\leq
\pi_{p+1}$. In other words, $\pi_{\left(  p+1\right)  -1}\leq\pi_{p+1}$ (since
$p=\left(  p+1\right)  -1$).

Also, $p+1\in S\subseteq\left[  n\right]  $ and $j\in\left[  n\right]  $ and
$p+1\leq j+1$ (since $p\leq j$). Hence, Observation 1 (applied to $i=p+1$)
yields $\left\{  p+1,p+2,\ldots,j\right\}  \cap\operatorname{Epk}\pi
\neq\varnothing$. In other words, the set $\left\{  p+1,p+2,\ldots,j\right\}
\cap\operatorname{Epk}\pi$ has at least one element. Consider such an
element, and denote it by $z$. Thus,%
\[
z\in\left\{  p+1,p+2,\ldots,j\right\}  \cap\operatorname{Epk}\pi.
\]

Now,%
\begin{align*}
z  &  \in\underbrace{\left\{  p+1,p+2,\ldots,j\right\}  }_{\substack{\subseteq
\left\{  a,a+1,\ldots,b\right\}  \\\text{(since }p+1\geq a\text{ and }j\leq
b\text{)}}}\cap\operatorname{Epk}\pi
\subseteq
\underbrace{\left\{  a,a+1,\ldots
,b\right\}  }_{=S}\cap\operatorname{Epk}\pi\\
&  =S\cap\operatorname{Epk}\pi
\subseteq \left\{  a,b\right\}  .
\end{align*}
But $z\in\left\{  p+1,p+2,\ldots,j\right\}  \cap\operatorname{Epk} \pi
\subseteq\left\{  p+1,p+2,\ldots,j\right\}  $, so that $z\geq p+1>p\geq a$.
Hence, $z\neq a$. Combining $z\in\left\{  a,b\right\}  $ with $z\neq a$, we
obtain $z\in\left\{  a,b\right\}  \setminus\left\{  a\right\}  \subseteq
\left\{  b\right\}  $. In other words, $z=b$. But from $z\in\left\{
p+1,p+2,\ldots,j\right\}  $, we also obtain $z\leq j\leq b-1$ (since
$j\in\left\{  p,p+1,\ldots,b-1\right\}  $). Hence, $z\leq b-1<b$; this
contradicts $z=b$.

This contradiction shows that our assumption (that $\pi_{j}\geq\pi_{j+1}$) was
wrong. Hence, we don't have $\pi_{j}\geq\pi_{j+1}$. In other words, we have
$\pi_{j}<\pi_{j+1}$.

Now, forget that we fixed $j$. We thus have shown that $\pi_{j}<\pi_{j+1}$ for
each $j\in\left\{  p,p+1,\ldots,b-1\right\}  $. In other words, $\pi_{p}%
<\pi_{p+1}<\pi_{p+2}<\cdots<\pi_{b}$. In other words, the map $\pi
\mid_{\left\{  p,p+1,\ldots,b\right\}  }$ is strictly increasing. In other
words, the map $\pi\mid_{\left\{  s\in S\ \mid\ s\geq p\right\}  }$ is
strictly increasing (since $\left\{  s\in S\ \mid\ s\geq p\right\}  =\left\{
p,p+1,\ldots,b\right\}  $).]

Thus, we now know that the map $\pi\mid_{\left\{  s\in S\ \mid\ s\leq
p\right\}  }$ is strictly decreasing while the map $\pi\mid_{\left\{  s\in
S\ \mid\ s\geq p\right\}  }$ is strictly increasing. In view of $\pi
\mid_{\left\{  s\in S\ \mid\ s\leq p\right\}  }=\left(  \pi\mid_{S}\right)
\mid_{\left\{  s\in S\ \mid\ s\leq p\right\}  }$ and $\pi\mid_{\left\{  s\in
S\ \mid\ s\geq p\right\}  }=\left(  \pi\mid_{S}\right)  \mid_{\left\{  s\in
S\ \mid\ s\geq p\right\}  }$, this rewrites as follows: The map $\left(
\pi\mid_{S}\right)  \mid_{\left\{  s\in S\ \mid\ s\leq p\right\}  }$ is
strictly decreasing while the map $\left(  \pi\mid_{S}\right)  \mid_{\left\{
s\in S\ \mid\ s\geq p\right\}  }$ is strictly increasing. Hence, there exists
some $t\in S$ such that the map $\left(  \pi\mid_{S}\right)  \mid_{\left\{
s\in S\ \mid\ s\leq t\right\}  }$ is strictly decreasing while the map
$\left(  \pi\mid_{S}\right)  \mid_{\left\{  s\in S\ \mid\ s\geq t\right\}  }$
is strictly increasing (namely, $t=p$). In other words, the map $\pi\mid_{S}$
is V-shaped (by the definition of \textquotedblleft V-shaped\textquotedblright%
). This proves Lemma \ref{lem.Epk.V-sh} \textbf{(c)}.

\textbf{(a)} Assume that $1\in S$ and $S\cap\operatorname{Epk}\pi
\subseteq\left\{  \max S\right\}  $. We must prove that the map $\pi\mid_{S}$
is strictly increasing.

Write the interval $S$ in the form $S=\left\{  a,a+1,\ldots,b\right\}  $ for
some elements $a$ and $b$ of $\left[  n\right]  $. Thus, $a\leq b$ (since the
interval $S$ is nonempty) and $\min S=a$ and $\max S=b$. But $1\in S$ and thus
$\min S\leq1$. Hence, $a=\min S\leq1$. Combined with $a\geq1$ (since
$a\in\left[  n\right]  $), this yields $a=1$. Now,%
\[
S=\left\{  a,a+1,\ldots,b\right\}  =\left\{  1,2,\ldots,b\right\}
\ \ \ \ \ \ \ \ \ \ \left(  \text{since }a=1\right)  .
\]

But $S\cap\operatorname{Epk}\pi\subseteq\left\{  \max S\right\}  =\left\{
b\right\}  $ (since $\max S=b$).

Let $j\in\left\{  1,2,\ldots,b-1\right\}  $. Thus, $1\leq j\leq b-1$. Hence,
$j\leq b-1<b\leq n$ (since $b\in\left[  n\right]  $), so that $j\in\left[
n\right]  $ (since $1\leq j\leq n$).

Assume (for the sake of contradiction) that $\pi_{j}\geq\pi_{j+1}$.

We have $\pi_{1-1}=\pi_{0}=0\leq\pi_{1}$ (since $\pi_{1}\in\mathbb{N}$) and
$1\in S\subseteq\left[  n\right]  $ and $j\in\left[  n\right]  $ and $1\leq
j+1$ (since $j+1>j\geq1$) and $\pi_{j}\geq\pi_{j+1}$. Thus, Observation 1
(applied to $i=1$) yields
\[
\left\{  1,2,\ldots,j\right\}  \cap\operatorname{Epk}\pi\neq\varnothing.
\]
In other words, the set $\left\{  1,2,\ldots,j\right\}  \cap
\operatorname{Epk}\pi$ has at least one element. Consider such an element,
and denote it by $z$. Thus,%
\[
z\in\left\{  1,2,\ldots,j\right\}  \cap\operatorname{Epk}\pi.
\]

Now,%
\begin{align*}
z  &  \in\underbrace{\left\{  1,2,\ldots,j\right\}  }_{\substack{\subseteq
\left\{  a,a+1,\ldots,b\right\}  \\\text{(since }1=a\text{ and }j\leq b-1\leq
b\text{)}}}\cap\operatorname{Epk}\pi
\subseteq
\underbrace{\left\{  a,a+1,\ldots
,b\right\}  }_{=S}\cap\operatorname{Epk}\pi\\
&  =S\cap\operatorname{Epk}\pi
\subseteq \left\{  b\right\}  .
\end{align*}
In other words, $z=b$.

But $z\in\left\{  1,2,\ldots,j\right\}  \cap\operatorname{Epk}\pi
\subseteq\left\{  1,2,\ldots,j\right\}  $, so that $z\leq j\leq b-1<b$. Hence,
$z\neq b$. This contradicts $z=b$.

This contradiction shows that our assumption (that $\pi_{j}\geq\pi_{j+1}$) was
wrong. Hence, we don't have $\pi_{j}\geq\pi_{j+1}$. In other words, we have
$\pi_{j}<\pi_{j+1}$.

Now, forget that we fixed $j$. We thus have shown that $\pi_{j}<\pi_{j+1}$ for
each $j\in\left\{  1,2,\ldots,b-1\right\}  $. In other words, $\pi_{1}<\pi
_{2}<\cdots<\pi_{b}$. In other words, the map $\pi\mid_{\left\{
1,2,\ldots,b\right\}  }$ is strictly increasing. In other words, the map
$\pi\mid_{S}$ is strictly increasing (since $S=\left\{  1,2,\ldots,b\right\}
$). This proves Lemma \ref{lem.Epk.V-sh} \textbf{(a)}.

\textbf{(b)} Assume that $n\in S$ and $S\cap\operatorname{Epk}\pi
\subseteq\left\{  \min S\right\}  $. We must prove that the map $\pi\mid_{S}$
is strictly decreasing.

Write the interval $S$ in the form $S=\left\{  a,a+1,\ldots,b\right\}  $ for
some elements $a$ and $b$ of $\left[  n\right]  $. Thus, $a\leq b$ (since the
interval $S$ is nonempty) and $\min S=a$ and $\max S=b$. But $n\in S$ and thus
$\max S\geq n$. Hence, $b=\max S\geq n$. Combined with $b\leq n$ (since
$b\in\left[  n\right]  $), this yields $b=n$. Now,%
\[
S=\left\{  a,a+1,\ldots,b\right\}  =\left\{  a,a+1,\ldots,n\right\}
\ \ \ \ \ \ \ \ \ \ \left(  \text{since }b=n\right)  .
\]

But $S\cap\operatorname{Epk}\pi\subseteq\left\{  \min S\right\}  =\left\{
a\right\}  $ (since $\min S=a$).

Let $i\in\left\{  a+1,a+2,\ldots,n\right\}  $. Thus, $a+1\leq i\leq n$. Hence,
$i\geq a+1>a\geq1$ (since $a\in\left[  n\right]  $), so that $i\in\left[
n\right]  $ (since $1\leq i\leq n$).

Assume (for the sake of contradiction) that $\pi_{i-1}\leq\pi_{i}$.

We have $\pi_{n+1}=0\leq\pi_{n}$ (since $\pi_{n}\in\mathbb{N}$) and thus
$\pi_{n}\geq\pi_{n+1}$. Also, $n\in S\subseteq\left[  n\right]  $ and
$i\in\left[  n\right]  $ and $i\leq n+1$ (since $i\in\left[  n\right]
\subseteq\left[  n+1\right]  $) and $\pi_{i-1}\leq\pi_{i}$. Thus, Observation
1 (applied to $j=n$) yields
\[
\left\{  i,i+1,\ldots,n\right\}  \cap\operatorname{Epk}\pi\neq\varnothing.
\]
In other words, the set $\left\{  i,i+1,\ldots,n\right\}  \cap
\operatorname{Epk}\pi$ has at least one element. Consider such an element,
and denote it by $z$. Thus,%
\[
z\in\left\{  i,i+1,\ldots,n\right\}  \cap\operatorname{Epk}\pi.
\]

Now,%
\begin{align*}
z  &  \in\underbrace{\left\{  i,i+1,\ldots,n\right\}  }_{\substack{\subseteq
\left\{  a,a+1,\ldots,b\right\}  \\\text{(since }i\geq a\text{ and
}n=b\text{)}}}\cap\operatorname{Epk}\pi
\subseteq
\underbrace{\left\{  a,a+1,\ldots
,b\right\}  }_{=S}\cap\operatorname{Epk}\pi\\
&  =S\cap\operatorname{Epk}\pi
\subseteq \left\{  a\right\}  .
\end{align*}
In other words, $z=a$.

But $z\in\left\{  i,i+1,\ldots,n\right\}  \cap\operatorname{Epk}\pi
\subseteq\left\{  i,i+1,\ldots,n\right\}  $, so that $z\geq i\geq a+1>a$.
Hence, $z\neq a$. This contradicts $z=a$.

This contradiction shows that our assumption (that $\pi_{i-1}\leq\pi_{i}$) was
wrong. Hence, we don't have $\pi_{i-1}\leq\pi_{i}$. In other words, we have
$\pi_{i-1}>\pi_{i}$.

Now, forget that we fixed $i$. We thus have shown that $\pi_{i-1}>\pi_{i}$ for
each $i\in\left\{  a+1,a+2,\ldots,n\right\}  $. In other words, $\pi_{a}%
>\pi_{a+1}>\pi_{a+2}>\cdots>\pi_{n}$. In other words, the map $\pi
\mid_{\left\{  a,a+1,\ldots,n\right\}  }$ is strictly decreasing. In other
words, the map $\pi\mid_{S}$ is strictly decreasing (since $S=\left\{
a,a+1,\ldots,n\right\}  $). This proves Lemma \ref{lem.Epk.V-sh} \textbf{(b)}.
\end{proof}

\begin{lemma}
\label{lem.fiberends.cut}Let $n\in\mathbb{N}$. Let $g:\left[  n\right]
\rightarrow\mathcal{N}$ be any weakly increasing map. Let $u\in g\left(
\left[  n\right]  \right)  $.

\textbf{(a)} If $u=0$, then
\[
\operatorname*{FE}\left(  g\right)  \cap g^{-1}\left(  u\right)  =\left\{
\max\left(  g^{-1}\left(  u\right)  \right)  \right\}  .
\]

\textbf{(b)} If $u\in\left\{  1,2,3,\ldots\right\}  $, then%
\[
\operatorname*{FE}\left(  g\right)  \cap g^{-1}\left(  u\right)  =\left\{
\min\left(  g^{-1}\left(  u\right)  \right)  ,\max\left(  g^{-1}\left(
u\right)  \right)  \right\}  .
\]

\textbf{(c)} If $u=\infty$, then
\[
\operatorname*{FE}\left(  g\right)  \cap g^{-1}\left(  u\right)  =\left\{
\min\left(  g^{-1}\left(  u\right)  \right)  \right\}  .
\]

\end{lemma}

\begin{proof}
[Proof of Lemma \ref{lem.fiberends.cut}.]From $u\in g\left(  \left[  n\right]
\right)  $, we conclude that there exists some $j\in\left[  n\right]  $ such
that $u=g\left(  j\right)  $. In other words, the fiber $g^{-1}\left(
u\right)  $ is nonempty. In other words, $g^{-1}\left(  u\right)
\neq\varnothing$.

Recall that%
\begin{align}
\operatorname*{FE}\left(  g\right)   &  =\left\{  \min\left(  g^{-1}\left(
h\right)  \right)  \ \mid\ h\in\left\{  1,2,3,\ldots,\infty\right\}  \text{
with }g^{-1}\left(  h\right)  \neq\varnothing\right\} \nonumber\\
&  \ \ \ \ \ \ \ \ \ \ \cup\left\{  \max\left(  g^{-1}\left(  h\right)
\right)  \ \mid\ h\in\left\{  0,1,2,3,\ldots\right\}  \text{ with }%
g^{-1}\left(  h\right)  \neq\varnothing\right\}
\label{pf.lem.fiberends.cut.FEg=}%
\end{align}
(by the definition of $\operatorname*{FE}\left(  g\right)  $).

If $u\in\left\{  0,1,2,3,\ldots\right\}  $, then%
\begin{align}
\max\left(  g^{-1}\left(  u\right)  \right)   &  \in\left\{  \max\left(
g^{-1}\left(  h\right)  \right)  \ \mid\ h\in\left\{  0,1,2,3,\ldots\right\}
\text{ with }g^{-1}\left(  h\right)  \neq\varnothing\right\} \nonumber\\
&  \ \ \ \ \ \ \ \ \ \ \left(
\begin{array}
[c]{c}%
\text{since }u\text{ is an }h\in\left\{  0,1,2,3,\ldots\right\}  \text{
satisfying }g^{-1}\left(  h\right)  \neq\varnothing\\
\text{(since }u\in\left\{  0,1,2,3,\ldots\right\}  \text{ and }g^{-1}\left(
u\right)  \neq\varnothing\text{)}%
\end{array}
\right) \nonumber\\
&  \subseteq\left\{  \min\left(  g^{-1}\left(  h\right)  \right)  \ \mid
\ h\in\left\{  1,2,3,\ldots,\infty\right\}  \text{ with }g^{-1}\left(
h\right)  \neq\varnothing\right\} \nonumber\\
&  \ \ \ \ \ \ \ \ \ \ \cup\left\{  \max\left(  g^{-1}\left(  h\right)
\right)  \ \mid\ h\in\left\{  0,1,2,3,\ldots\right\}  \text{ with }%
g^{-1}\left(  h\right)  \neq\varnothing\right\} \nonumber\\
&  =\operatorname*{FE}\left(  g\right)  \ \ \ \ \ \ \ \ \ \ \left(  \text{by
(\ref{pf.lem.fiberends.cut.FEg=})}\right)  .
\label{pf.lem.fiberends.cut.maxin}%
\end{align}

If $u\in\left\{  1,2,3,\ldots,\infty\right\}  $, then%
\begin{align}
\min\left(  g^{-1}\left(  u\right)  \right)   &  \in\left\{  \min\left(
g^{-1}\left(  h\right)  \right)  \ \mid\ h\in\left\{  1,2,3,\ldots
,\infty\right\}  \text{ with }g^{-1}\left(  h\right)  \neq\varnothing\right\}
\nonumber\\
&  \ \ \ \ \ \ \ \ \ \ \left(
\begin{array}
[c]{c}%
\text{since }u\text{ is an }h\in\left\{  1,2,3,\ldots,\infty\right\}  \text{
satisfying }g^{-1}\left(  h\right)  \neq\varnothing\\
\text{(since }u\in\left\{  1,2,3,\ldots,\infty\right\}  \text{ and }%
g^{-1}\left(  u\right)  \neq\varnothing\text{)}%
\end{array}
\right) \nonumber\\
&  \subseteq\left\{  \min\left(  g^{-1}\left(  h\right)  \right)  \ \mid
\ h\in\left\{  1,2,3,\ldots,\infty\right\}  \text{ with }g^{-1}\left(
h\right)  \neq\varnothing\right\} \nonumber\\
&  \ \ \ \ \ \ \ \ \ \ \cup\left\{  \max\left(  g^{-1}\left(  h\right)
\right)  \ \mid\ h\in\left\{  0,1,2,3,\ldots\right\}  \text{ with }%
g^{-1}\left(  h\right)  \neq\varnothing\right\} \nonumber\\
&  =\operatorname*{FE}\left(  g\right)  \ \ \ \ \ \ \ \ \ \ \left(  \text{by
(\ref{pf.lem.fiberends.cut.FEg=})}\right)  .
\label{pf.lem.fiberends.cut.minin}%
\end{align}

\textbf{(a)} Assume that $u=0$. We must prove that $\operatorname*{FE}\left(
g\right)  \cap g^{-1}\left(  u\right)  =\left\{  \max\left(  g^{-1}\left(
u\right)  \right)  \right\}  $.

We have $u=0\in\left\{  0,1,2,3,\ldots\right\}  $. Thus,
(\ref{pf.lem.fiberends.cut.maxin}) shows that $\max\left(  g^{-1}\left(
u\right)  \right)  \in\operatorname*{FE}\left(  g\right)  $. Combining this
with $\max\left(  g^{-1}\left(  u\right)  \right)  \in g^{-1}\left(  u\right)
$ (which is obvious), we obtain $\max\left(  g^{-1}\left(  u\right)  \right)
\in\operatorname*{FE}\left(  g\right)  \cap g^{-1}\left(  u\right)  $. In
other words,%
\begin{equation}
\left\{  \max\left(  g^{-1}\left(  u\right)  \right)  \right\}  \subseteq
\operatorname*{FE}\left(  g\right)  \cap g^{-1}\left(  u\right)  .
\label{pf.lem.fiberends.cut.a.1}%
\end{equation}

On the other hand, let $j\in\operatorname*{FE}\left(  g\right)  \cap
g^{-1}\left(  u\right)  $ be arbitrary. We shall show that $j\in\left\{
\max\left(  g^{-1}\left(  u\right)  \right)  \right\}  $.

Indeed, we have $j\in\operatorname*{FE}\left(  g\right)  \cap g^{-1}\left(
u\right)  \subseteq g^{-1}\left(  u\right)  $, so that $g\left(  j\right)
=u$. On the other hand,%
\begin{align*}
j  &  \in\operatorname*{FE}\left(  g\right)  \cap g^{-1}\left(  u\right)
\subseteq\operatorname*{FE}\left(  g\right) \\
&  =\left\{  \min\left(  g^{-1}\left(  h\right)  \right)  \ \mid\ h\in\left\{
1,2,3,\ldots,\infty\right\}  \text{ with }g^{-1}\left(  h\right)
\neq\varnothing\right\} \\
&  \ \ \ \ \ \ \ \ \ \ \cup\left\{  \max\left(  g^{-1}\left(  h\right)
\right)  \ \mid\ h\in\left\{  0,1,2,3,\ldots\right\}  \text{ with }%
g^{-1}\left(  h\right)  \neq\varnothing\right\} \\
&  \ \ \ \ \ \ \ \ \ \ \left(  \text{by (\ref{pf.lem.fiberends.cut.FEg=}%
)}\right)  .
\end{align*}
In other words, either $j\in\left\{  \min\left(  g^{-1}\left(  h\right)
\right)  \ \mid\ h\in\left\{  1,2,3,\ldots,\infty\right\}  \text{ with }%
g^{-1}\left(  h\right)  \neq\varnothing\right\}  $ or $j\in\left\{
\max\left(  g^{-1}\left(  h\right)  \right)  \ \mid\ h\in\left\{
0,1,2,3,\ldots\right\}  \text{ with }g^{-1}\left(  h\right)  \neq
\varnothing\right\}  $. Hence, we are in one of the following two cases:

\textit{Case 1:} We have $j\in\left\{  \min\left(  g^{-1}\left(  h\right)
\right)  \ \mid\ h\in\left\{  1,2,3,\ldots,\infty\right\}  \text{ with }%
g^{-1}\left(  h\right)  \neq\varnothing\right\}  $.

\textit{Case 2:} We have $j\in\left\{  \max\left(  g^{-1}\left(  h\right)
\right)  \ \mid\ h\in\left\{  0,1,2,3,\ldots\right\}  \text{ with }%
g^{-1}\left(  h\right)  \neq\varnothing\right\}  $.

Let us first consider Case 1. In this case, we have \newline$j\in\left\{
\min\left(  g^{-1}\left(  h\right)  \right)  \ \mid\ h\in\left\{
1,2,3,\ldots,\infty\right\}  \text{ with }g^{-1}\left(  h\right)
\neq\varnothing\right\}  $. In other words, $j=\min\left(  g^{-1}\left(
h\right)  \right)  $ for some $h\in\left\{  1,2,3,\ldots,\infty\right\}  $
satisfying $g^{-1}\left(  h\right)  \neq\varnothing$. Consider this $h$. We
have $j=\min\left(  g^{-1}\left(  h\right)  \right)  \in g^{-1}\left(
h\right)  $, so that $g\left(  j\right)  =h$. Hence, $h=g\left(  j\right)
=u$. Hence, $u=h\in\left\{  1,2,3,\ldots,\infty\right\}  $. This contradicts
$u=0\notin\left\{  1,2,3,\ldots,\infty\right\}  $. Thus, $j\in\left\{
\max\left(  g^{-1}\left(  u\right)  \right)  \right\}  $ (because anything
follows from a contradiction). Hence, we have proven $j\in\left\{  \max\left(
g^{-1}\left(  u\right)  \right)  \right\}  $ in Case 1.

Let us now consider Case 2. In this case, we have \newline$j\in\left\{
\max\left(  g^{-1}\left(  h\right)  \right)  \ \mid\ h\in\left\{
0,1,2,3,\ldots\right\}  \text{ with }g^{-1}\left(  h\right)  \neq
\varnothing\right\}  $. In other words, $j=\max\left(  g^{-1}\left(  h\right)
\right)  $ for some $h\in\left\{  0,1,2,3,\ldots\right\}  $ satisfying
$g^{-1}\left(  h\right)  \neq\varnothing$. Consider this $h$. We have
$j=\max\left(  g^{-1}\left(  h\right)  \right)  \in g^{-1}\left(  h\right)  $,
so that $g\left(  j\right)  =h$. Hence, $h=g\left(  j\right)  =u$. Hence,
$j=\max\left(  g^{-1}\left(  \underbrace{h}_{=u}\right)  \right)  =\max\left(
g^{-1}\left(  u\right)  \right)  \in\left\{  \max\left(  g^{-1}\left(
u\right)  \right)  \right\}  $. Hence, we have proven $j\in\left\{
\max\left(  g^{-1}\left(  u\right)  \right)  \right\}  $ in Case 2.

We have now proven $j\in\left\{  \max\left(  g^{-1}\left(  u\right)  \right)
\right\}  $ in each of the two Cases 1 and 2. Thus, $j\in\left\{  \max\left(
g^{-1}\left(  u\right)  \right)  \right\}  $ always holds.

Now, forget that we fixed $j$. We thus have proven that $j\in\left\{
\max\left(  g^{-1}\left(  u\right)  \right)  \right\}  $ for each
$j\in\operatorname*{FE}\left(  g\right)  \cap g^{-1}\left(  u\right)  $. In
other words,%
\[
\operatorname*{FE}\left(  g\right)  \cap g^{-1}\left(  u\right)
\subseteq\left\{  \max\left(  g^{-1}\left(  u\right)  \right)  \right\}  .
\]
Combining this with (\ref{pf.lem.fiberends.cut.a.1}), we obtain
$\operatorname*{FE}\left(  g\right)  \cap g^{-1}\left(  u\right)  =\left\{
\max\left(  g^{-1}\left(  u\right)  \right)  \right\}  $. This proves Lemma
\ref{lem.fiberends.cut} \textbf{(a)}.

\textbf{(b)} Assume that $u\in\left\{  1,2,3,\ldots\right\}  $. We must prove
that $\operatorname*{FE}\left(  g\right)  \cap g^{-1}\left(  u\right)
=\left\{  \min\left(  g^{-1}\left(  u\right)  \right)  ,\max\left(
g^{-1}\left(  u\right)  \right)  \right\}  $.

We have $u\in\left\{  1,2,3,\ldots\right\}  \subseteq\left\{  0,1,2,3,\ldots
\right\}  $. Thus, (\ref{pf.lem.fiberends.cut.maxin}) shows that $\max\left(
g^{-1}\left(  u\right)  \right)  \in\operatorname*{FE}\left(  g\right)  $.
Combining this with $\max\left(  g^{-1}\left(  u\right)  \right)  \in
g^{-1}\left(  u\right)  $ (which is obvious), we obtain $\max\left(
g^{-1}\left(  u\right)  \right)  \in\operatorname*{FE}\left(  g\right)  \cap
g^{-1}\left(  u\right)  $.

We have $u\in\left\{  1,2,3,\ldots\right\}  \subseteq\left\{  1,2,3,\ldots
,\infty\right\}  $. Thus, (\ref{pf.lem.fiberends.cut.minin}) shows that
$\min\left(  g^{-1}\left(  u\right)  \right)  \in\operatorname*{FE}\left(
g\right)  $. Combining this with $\min\left(  g^{-1}\left(  u\right)  \right)
\in g^{-1}\left(  u\right)  $ (which is obvious), we obtain $\min\left(
g^{-1}\left(  u\right)  \right)  \in\operatorname*{FE}\left(  g\right)  \cap
g^{-1}\left(  u\right)  $.

Now, we have shown that both numbers $\min\left(  g^{-1}\left(  u\right)
\right)  $ and $\max\left(  g^{-1}\left(  u\right)  \right)  $ belong to
$\operatorname*{FE}\left(  g\right)  \cap g^{-1}\left(  u\right)  $. In other
words,%
\begin{equation}
\left\{  \min\left(  g^{-1}\left(  u\right)  \right)  ,\max\left(
g^{-1}\left(  u\right)  \right)  \right\}  \subseteq\operatorname*{FE}\left(
g\right)  \cap g^{-1}\left(  u\right)  . \label{pf.lem.fiberends.cut.b.1}%
\end{equation}

On the other hand, let $j\in\operatorname*{FE}\left(  g\right)  \cap
g^{-1}\left(  u\right)  $ be arbitrary. We shall show that $j\in\left\{
\min\left(  g^{-1}\left(  u\right)  \right)  ,\max\left(  g^{-1}\left(
u\right)  \right)  \right\}  $.

Indeed, we have $j\in\operatorname*{FE}\left(  g\right)  \cap g^{-1}\left(
u\right)  \subseteq g^{-1}\left(  u\right)  $, so that $g\left(  j\right)
=u$. On the other hand,%
\begin{align*}
j  &  \in\operatorname*{FE}\left(  g\right)  \cap g^{-1}\left(  u\right)
\subseteq\operatorname*{FE}\left(  g\right) \\
&  =\left\{  \min\left(  g^{-1}\left(  h\right)  \right)  \ \mid\ h\in\left\{
1,2,3,\ldots,\infty\right\}  \text{ with }g^{-1}\left(  h\right)
\neq\varnothing\right\} \\
&  \ \ \ \ \ \ \ \ \ \ \cup\left\{  \max\left(  g^{-1}\left(  h\right)
\right)  \ \mid\ h\in\left\{  0,1,2,3,\ldots\right\}  \text{ with }%
g^{-1}\left(  h\right)  \neq\varnothing\right\} \\
&  \ \ \ \ \ \ \ \ \ \ \left(  \text{by (\ref{pf.lem.fiberends.cut.FEg=}%
)}\right)  .
\end{align*}
In other words, either $j\in\left\{  \min\left(  g^{-1}\left(  h\right)
\right)  \ \mid\ h\in\left\{  1,2,3,\ldots,\infty\right\}  \text{ with }%
g^{-1}\left(  h\right)  \neq\varnothing\right\}  $ or $j\in\left\{
\max\left(  g^{-1}\left(  h\right)  \right)  \ \mid\ h\in\left\{
0,1,2,3,\ldots\right\}  \text{ with }g^{-1}\left(  h\right)  \neq
\varnothing\right\}  $. Hence, we are in one of the following two cases:

\textit{Case 1:} We have $j\in\left\{  \min\left(  g^{-1}\left(  h\right)
\right)  \ \mid\ h\in\left\{  1,2,3,\ldots,\infty\right\}  \text{ with }%
g^{-1}\left(  h\right)  \neq\varnothing\right\}  $.

\textit{Case 2:} We have $j\in\left\{  \max\left(  g^{-1}\left(  h\right)
\right)  \ \mid\ h\in\left\{  0,1,2,3,\ldots\right\}  \text{ with }%
g^{-1}\left(  h\right)  \neq\varnothing\right\}  $.

Let us first consider Case 1. In this case, we have \newline$j\in\left\{
\min\left(  g^{-1}\left(  h\right)  \right)  \ \mid\ h\in\left\{
1,2,3,\ldots,\infty\right\}  \text{ with }g^{-1}\left(  h\right)
\neq\varnothing\right\}  $. In other words, $j=\min\left(  g^{-1}\left(
h\right)  \right)  $ for some $h\in\left\{  1,2,3,\ldots,\infty\right\}  $
satisfying $g^{-1}\left(  h\right)  \neq\varnothing$. Consider this $h$. We
have $j=\min\left(  g^{-1}\left(  h\right)  \right)  \in g^{-1}\left(
h\right)  $, so that $g\left(  j\right)  =h$. Hence, $h=g\left(  j\right)
=u$. Hence, $j=\min\left(  g^{-1}\left(  \underbrace{h}_{=u}\right)  \right)
=\min\left(  g^{-1}\left(  u\right)  \right)  \in\left\{  \min\left(
g^{-1}\left(  u\right)  \right)  ,\max\left(  g^{-1}\left(  u\right)  \right)
\right\}  $. Hence, we have proven $j\in\left\{  \min\left(  g^{-1}\left(
u\right)  \right)  ,\max\left(  g^{-1}\left(  u\right)  \right)  \right\}  $
in Case 1.

Let us now consider Case 2. In this case, we have \newline$j\in\left\{
\max\left(  g^{-1}\left(  h\right)  \right)  \ \mid\ h\in\left\{
0,1,2,3,\ldots\right\}  \text{ with }g^{-1}\left(  h\right)  \neq
\varnothing\right\}  $. In other words, $j=\max\left(  g^{-1}\left(  h\right)
\right)  $ for some $h\in\left\{  0,1,2,3,\ldots\right\}  $ satisfying
$g^{-1}\left(  h\right)  \neq\varnothing$. Consider this $h$. We have
$j=\max\left(  g^{-1}\left(  h\right)  \right)  \in g^{-1}\left(  h\right)  $,
so that $g\left(  j\right)  =h$. Hence, $h=g\left(  j\right)  =u$. Hence,
$j=\max\left(  g^{-1}\left(  \underbrace{h}_{=u}\right)  \right)  =\max\left(
g^{-1}\left(  u\right)  \right)  \in\left\{  \min\left(  g^{-1}\left(
u\right)  \right)  ,\max\left(  g^{-1}\left(  u\right)  \right)  \right\}  $.
Hence, we have proven $j\in\left\{  \min\left(  g^{-1}\left(  u\right)
\right)  ,\max\left(  g^{-1}\left(  u\right)  \right)  \right\}  $ in Case 2.

We have now proven $j\in\left\{  \min\left(  g^{-1}\left(  u\right)  \right)
,\max\left(  g^{-1}\left(  u\right)  \right)  \right\}  $ in each of the two
Cases 1 and 2. Thus, $j\in\left\{  \min\left(  g^{-1}\left(  u\right)
\right)  ,\max\left(  g^{-1}\left(  u\right)  \right)  \right\}  $ always holds.

Now, forget that we fixed $j$. We thus have proven that \newline$j\in\left\{
\min\left(  g^{-1}\left(  u\right)  \right)  ,\max\left(  g^{-1}\left(
u\right)  \right)  \right\}  $ for each $j\in\operatorname*{FE}\left(
g\right)  \cap g^{-1}\left(  u\right)  $. In other words,%
\[
\operatorname*{FE}\left(  g\right)  \cap g^{-1}\left(  u\right)
\subseteq\left\{  \min\left(  g^{-1}\left(  u\right)  \right)  ,\max\left(
g^{-1}\left(  u\right)  \right)  \right\}  .
\]
Combining this with (\ref{pf.lem.fiberends.cut.b.1}), we obtain
$\operatorname*{FE}\left(  g\right)  \cap g^{-1}\left(  u\right)  =\left\{
\min\left(  g^{-1}\left(  u\right)  \right)  ,\max\left(  g^{-1}\left(
u\right)  \right)  \right\}  $. This proves Lemma \ref{lem.fiberends.cut}
\textbf{(b)}.

\textbf{(c)} Assume that $u=\infty$. We must prove that $\operatorname*{FE}%
\left(  g\right)  \cap g^{-1}\left(  u\right)  =\left\{  \min\left(
g^{-1}\left(  u\right)  \right)  \right\}  $.

We have $u=\infty\in\left\{  1,2,3,\ldots,\infty\right\}  $. Thus,
(\ref{pf.lem.fiberends.cut.minin}) shows that $\min\left(  g^{-1}\left(
u\right)  \right)  \in\operatorname*{FE}\left(  g\right)  $. Combining this
with $\min\left(  g^{-1}\left(  u\right)  \right)  \in g^{-1}\left(  u\right)
$ (which is obvious), we obtain $\min\left(  g^{-1}\left(  u\right)  \right)
\in\operatorname*{FE}\left(  g\right)  \cap g^{-1}\left(  u\right)  $. In
other words,%
\begin{equation}
\left\{  \min\left(  g^{-1}\left(  u\right)  \right)  \right\}  \subseteq
\operatorname*{FE}\left(  g\right)  \cap g^{-1}\left(  u\right)  .
\label{pf.lem.fiberends.cut.c.1}%
\end{equation}

On the other hand, let $j\in\operatorname*{FE}\left(  g\right)  \cap
g^{-1}\left(  u\right)  $ be arbitrary. We shall show that $j\in\left\{
\min\left(  g^{-1}\left(  u\right)  \right)  \right\}  $.

Indeed, we have $j\in\operatorname*{FE}\left(  g\right)  \cap g^{-1}\left(
u\right)  \subseteq g^{-1}\left(  u\right)  $, so that $g\left(  j\right)
=u$. On the other hand,%
\begin{align*}
j  &  \in\operatorname*{FE}\left(  g\right)  \cap g^{-1}\left(  u\right)
\subseteq\operatorname*{FE}\left(  g\right) \\
&  =\left\{  \min\left(  g^{-1}\left(  h\right)  \right)  \ \mid\ h\in\left\{
1,2,3,\ldots,\infty\right\}  \text{ with }g^{-1}\left(  h\right)
\neq\varnothing\right\} \\
&  \ \ \ \ \ \ \ \ \ \ \cup\left\{  \max\left(  g^{-1}\left(  h\right)
\right)  \ \mid\ h\in\left\{  0,1,2,3,\ldots\right\}  \text{ with }%
g^{-1}\left(  h\right)  \neq\varnothing\right\} \\
&  \ \ \ \ \ \ \ \ \ \ \left(  \text{by (\ref{pf.lem.fiberends.cut.FEg=}%
)}\right)  .
\end{align*}
In other words, either $j\in\left\{  \min\left(  g^{-1}\left(  h\right)
\right)  \ \mid\ h\in\left\{  1,2,3,\ldots,\infty\right\}  \text{ with }%
g^{-1}\left(  h\right)  \neq\varnothing\right\}  $ or $j\in\left\{
\max\left(  g^{-1}\left(  h\right)  \right)  \ \mid\ h\in\left\{
0,1,2,3,\ldots\right\}  \text{ with }g^{-1}\left(  h\right)  \neq
\varnothing\right\}  $. Hence, we are in one of the following two cases:

\textit{Case 1:} We have $j\in\left\{  \min\left(  g^{-1}\left(  h\right)
\right)  \ \mid\ h\in\left\{  1,2,3,\ldots,\infty\right\}  \text{ with }%
g^{-1}\left(  h\right)  \neq\varnothing\right\}  $.

\textit{Case 2:} We have $j\in\left\{  \max\left(  g^{-1}\left(  h\right)
\right)  \ \mid\ h\in\left\{  0,1,2,3,\ldots\right\}  \text{ with }%
g^{-1}\left(  h\right)  \neq\varnothing\right\}  $.

Let us first consider Case 1. In this case, we have \newline$j\in\left\{
\min\left(  g^{-1}\left(  h\right)  \right)  \ \mid\ h\in\left\{
1,2,3,\ldots,\infty\right\}  \text{ with }g^{-1}\left(  h\right)
\neq\varnothing\right\}  $. In other words, $j=\min\left(  g^{-1}\left(
h\right)  \right)  $ for some $h\in\left\{  1,2,3,\ldots,\infty\right\}  $
satisfying $g^{-1}\left(  h\right)  \neq\varnothing$. Consider this $h$. We
have $j=\min\left(  g^{-1}\left(  h\right)  \right)  \in g^{-1}\left(
h\right)  $, so that $g\left(  j\right)  =h$. Hence, $h=g\left(  j\right)
=u$. Hence, $j=\min\left(  g^{-1}\left(  \underbrace{h}_{=u}\right)  \right)
=\min\left(  g^{-1}\left(  u\right)  \right)  \in\left\{  \min\left(
g^{-1}\left(  u\right)  \right)  \right\}  $. Hence, we have proven
$j\in\left\{  \min\left(  g^{-1}\left(  u\right)  \right)  \right\}  $ in Case 1.

Let us now consider Case 2. In this case, we have \newline$j\in\left\{
\max\left(  g^{-1}\left(  h\right)  \right)  \ \mid\ h\in\left\{
0,1,2,3,\ldots\right\}  \text{ with }g^{-1}\left(  h\right)  \neq
\varnothing\right\}  $. In other words, $j=\max\left(  g^{-1}\left(  h\right)
\right)  $ for some $h\in\left\{  0,1,2,3,\ldots\right\}  $ satisfying
$g^{-1}\left(  h\right)  \neq\varnothing$. Consider this $h$. We have
$j=\max\left(  g^{-1}\left(  h\right)  \right)  \in g^{-1}\left(  h\right)  $,
so that $g\left(  j\right)  =h$. Hence, $h=g\left(  j\right)  =u$. Hence,
$u=h\in\left\{  0,1,2,3,\ldots\right\}  $. This contradicts $u=\infty
\notin\left\{  0,1,2,3,\ldots\right\}  $. Thus, $j\in\left\{  \min\left(
g^{-1}\left(  u\right)  \right)  \right\}  $ (because anything follows from a
contradiction). Hence, we have proven $j\in\left\{  \min\left(  g^{-1}\left(
u\right)  \right)  \right\}  $ in Case 2.

We have now proven $j\in\left\{  \min\left(  g^{-1}\left(  u\right)  \right)
\right\}  $ in each of the two Cases 1 and 2. Thus, $j\in\left\{  \min\left(
g^{-1}\left(  u\right)  \right)  \right\}  $ always holds.

Now, forget that we fixed $j$. We thus have proven that $j\in\left\{
\min\left(  g^{-1}\left(  u\right)  \right)  \right\}  $ for each
$j\in\operatorname*{FE}\left(  g\right)  \cap g^{-1}\left(  u\right)  $. In
other words,%
\[
\operatorname*{FE}\left(  g\right)  \cap g^{-1}\left(  u\right)
\subseteq\left\{  \min\left(  g^{-1}\left(  u\right)  \right)  \right\}  .
\]
Combining this with (\ref{pf.lem.fiberends.cut.c.1}), we obtain
$\operatorname*{FE}\left(  g\right)  \cap g^{-1}\left(  u\right)  =\left\{
\min\left(  g^{-1}\left(  u\right)  \right)  \right\}  $. This proves Lemma
\ref{lem.fiberends.cut} \textbf{(c)}.
\end{proof}
\end{verlong}

\begin{proposition}
\label{prop.Epk.fiberends}Let $n\in\mathbb{N}$.
Let $\pi$ be an $n$-permutation.
Let $g:\left[  n\right]
\rightarrow\mathcal{N}$ be any weakly increasing map. Then, the map $g$ is
$\pi$-amenable if and only if $\operatorname{Epk}\pi\subseteq
\operatorname*{FE}\left(  g\right)  $.
\end{proposition}

\begin{vershort}
\begin{proof}
[Proof of Proposition \ref{prop.Epk.fiberends}.] The map $g$ is weakly
increasing. Thus, all nonempty fibers $g^{-1}\left(  h\right)  $ of $g$ are
intervals of $\left[  n\right]  $. Recall that $g$ is $\pi$-amenable if and
only if the four Properties \textbf{(i')}, \textbf{(ii')}, \textbf{(iii')} and
\textbf{(iv')} in Definition \ref{def.amen} hold. Consider these four
properties. Since Property \textbf{(iv')} automatically holds (since we
assumed $g$ to be weakly increasing), we thus only need to discuss the other three:

\begin{itemize}
\item Property \textbf{(i')} is equivalent to the statement that every
exterior peak of $\pi$ that lies in the fiber $g^{-1}\left(  0\right)  $ must
be the largest element of this fiber. (Indeed, a strictly increasing map is
characterized by having no exterior peaks except for the largest element
of its domain.)

\item Property \textbf{(ii')} is equivalent to the statement that every peak
of $\pi$ that lies in a fiber $g^{-1}\left(  h\right)  $ with $h\in g\left(
\left[  n\right]  \right)  \cap\left\{  1,2,3,\ldots\right\}  $ must be either
the smallest or the largest element of this fiber (because the restriction
$\pi\mid_{g^{-1}\left(  h\right)  }$ is V-shaped if and only if no peak of
$\pi$ appears in the interior of this fiber\footnote{The \textit{interior} of
an interval $\left\{  a,a+1,\ldots,b\right\}  $ of $\left[  n\right]  $ is
defined to be the interval $\left\{  a+1,a+2,\ldots,b-1\right\}  $. (This is
an empty interval if $a+1>b-1$.)}). Moreover, we can replace the word
\textquotedblleft peak\textquotedblright\ by \textquotedblleft exterior
peak\textquotedblright\ in this sentence (since the exterior peaks $1$ and $n$
must automatically be the smallest and the largest element of whatever fibers
they belong to).

\item Property \textbf{(iii')} is equivalent to the statement that every
exterior peak of $\pi$ that lies in the fiber $g^{-1}\left(  \infty\right)  $
must be the smallest element of this fiber. (Indeed, a strictly decreasing map
is characterized by having no exterior peaks except for the smallest element
of its domain.)
\end{itemize}

Combining all of these insights, we conclude that the four Properties
\textbf{(i')}, \textbf{(ii')}, \textbf{(iii')} and \textbf{(iv')} hold if and
only if every exterior peak of $\pi$ is a fiber-end of $g$. In other words,
$g$ is $\pi$-amenable if and only if every exterior peak of $\pi$ is a
fiber-end of $g$ (since $g$ is $\pi$-amenable if and only if the four
Properties \textbf{(i')}, \textbf{(ii')}, \textbf{(iii')} and \textbf{(iv')}
hold). In other words, $g$ is $\pi$-amenable if and only if
$\operatorname{Epk}\pi\subseteq\operatorname*{FE}\left(  g\right)  $. This
proves Proposition \ref{prop.Epk.fiberends}.
\end{proof}
\end{vershort}

\begin{verlong}
\begin{proof}
[Proof of Proposition \ref{prop.Epk.fiberends}.]The map $g$ is weakly
increasing. Hence, each fiber $g^{-1}\left(  h\right)  $ of $g$ (with
$h\in\mathcal{N}$) is an interval of the totally ordered set $\left[
n\right]  $.

We shall prove the following two statements:

\begin{statement}
\textit{Observation 1:} If the map $g$ is $\pi$-amenable, then
$\operatorname{Epk}\pi\subseteq\operatorname*{FE}\left(  g\right)  $.
\end{statement}

\begin{statement}
\textit{Observation 2:} If $\operatorname{Epk}\pi\subseteq\operatorname*{FE}%
\left(  g\right)  $, then the map $g$ is $\pi$-amenable.
\end{statement}

[\textit{Proof of Observation 1:} Assume that the map $g$ is $\pi$-amenable.
We must show that $\operatorname{Epk}\pi\subseteq\operatorname*{FE}\left(
g\right)  $.

We know that the map $g$ is $\pi$-amenable. In other words, the four
properties \textbf{(i')}, \textbf{(ii')}, \textbf{(iii')} and \textbf{(iv')}
of Definition \ref{def.amen} hold (by the definition of \textquotedblleft$\pi
$-amenable\textquotedblright). We shall use these four properties in the following.

Let $j\in\operatorname{Epk}\pi$. We intend to show that $j\in
\operatorname*{FE}\left(  g\right)  $.

We have $j\in\operatorname{Epk}\pi$; in other words, $j$ is an exterior peak
of $\pi$ (by the definition of $\operatorname{Epk}\pi$). In other words, $j$
is an element of $\left[  n\right]  $ and satisfies $\pi_{j-1}<\pi_{j}%
>\pi_{j+1}$ (by the definition of an exterior peak). Define $u\in\mathcal{N}$
by $u=g\left(  j\right)  $. Thus, $j\in g^{-1}\left(  u\right)  $ and
$u=g\left(  \underbrace{j}_{\in\left[  n\right]  }\right)  \in g\left(
\left[  n\right]  \right)  $.

Recall that each fiber $g^{-1}\left(  h\right)  $ of $g$ (with $h\in
\mathcal{N}$) is an interval of the totally ordered set $\left[  n\right]  $.
Applying this to $h=u$, we conclude that $g^{-1}\left(  u\right)  $ is an
interval of the totally ordered set $\left[  n\right]  $. This interval is
furthermore nonempty (since $j\in g^{-1}\left(  u\right)  $).

We have $u\in\mathcal{N}=\left\{  0\right\}  \cup\left\{  \infty\right\}
\cup\left\{  1,2,3,\ldots\right\}  $. In other words, either $u=0$ or
$u=\infty$ or $u\in\left\{  1,2,3,\ldots\right\}  $. Hence, we are in one of
the following three cases:

\textit{Case 1:} We have $u=0$.

\textit{Case 2:} We have $u=\infty$.

\textit{Case 3:} We have $u\in\left\{  1,2,3,\ldots\right\}  $.

Let us first consider Case 1. In this case, we have $u=0$. Thus, Lemma
\ref{lem.fiberends.cut} \textbf{(a)} yields
\begin{equation}
\operatorname*{FE}\left(  g\right)  \cap g^{-1}\left(  u\right)  =\left\{
\max\left(  g^{-1}\left(  u\right)  \right)  \right\}  .
\label{pf.prop.Epk.fiberends.o1.pf.c1.1}%
\end{equation}

Recall that property \textbf{(i')} holds. In other words, the map $\pi
\mid_{g^{-1}\left(  0\right)  }$ is strictly increasing. Since $0=u$, this
rewrites as follows: The map $\pi\mid_{g^{-1}\left(  u\right)  }$ is strictly increasing.

But $j=\max\left(  g^{-1}\left(  u\right)  \right)  $%
\ \ \ \ \footnote{\textit{Proof.} Assume the contrary. Thus, $j\neq\max\left(
g^{-1}\left(  u\right)  \right)  $. Combining this with $j\leq\max\left(
g^{-1}\left(  u\right)  \right)  $ (which follows from $j\in g^{-1}\left(
u\right)  $), we obtain $j<\max\left(  g^{-1}\left(  u\right)  \right)  $.
\par
But recall the following basic fact: If $S$ is a nonempty interval of $\left[
n\right]  $, and if $s\in S$ satisfies $s<\max S$, then $s+1\in S$. Applying
this to $S=g^{-1}\left(  u\right)  $ and $s=j$, we conclude that $j+1\in
g^{-1}\left(  u\right)  $ (because $g^{-1}\left(  u\right)  $ is a nonempty
interval of $\left[  n\right]  $, and because $j\in g^{-1}\left(  u\right)  $
satisfies $j<\max\left(  g^{-1}\left(  u\right)  \right)  $). Hence, both $j$
and $j+1$ are elements of $g^{-1}\left(  u\right)  $, and satisfy $j<j+1$.
Hence, $\left(  \pi\mid_{g^{-1}\left(  u\right)  }\right)  \left(  j\right)
<\left(  \pi\mid_{g^{-1}\left(  u\right)  }\right)  \left(  j+1\right)  $
(since the map $\pi\mid_{g^{-1}\left(  u\right)  }$ is strictly increasing).
Thus,%
\[
\pi_{j}=\pi\left(  j\right)  =\left(  \pi\mid_{g^{-1}\left(  u\right)
}\right)  \left(  j\right)  <\left(  \pi\mid_{g^{-1}\left(  u\right)
}\right)  \left(  j+1\right)  =\pi\left(  j+1\right)  =\pi_{j+1}.
\]
This contradicts $\pi_{j}>\pi_{j+1}$. This contradiction shows that our
assumption was false; qed.}. Hence,%
\begin{align*}
j  &  =\max\left(  g^{-1}\left(  u\right)  \right)  \in\left\{  \max\left(
g^{-1}\left(  u\right)  \right)  \right\}  =\operatorname*{FE}\left(
g\right)  \cap g^{-1}\left(  u\right)  \ \ \ \ \ \ \ \ \ \ \left(  \text{by
(\ref{pf.prop.Epk.fiberends.o1.pf.c1.1})}\right) \\
&  \subseteq\operatorname*{FE}\left(  g\right)  .
\end{align*}
Thus, we have shown that $j\in\operatorname*{FE}\left(  g\right)  $ in Case 1.

Let us next consider Case 2. In this case, we have $u=\infty$. Thus, Lemma
\ref{lem.fiberends.cut} \textbf{(c)} yields
\begin{equation}
\operatorname*{FE}\left(  g\right)  \cap g^{-1}\left(  u\right)  =\left\{
\min\left(  g^{-1}\left(  u\right)  \right)  \right\}  .
\label{pf.prop.Epk.fiberends.o1.pf.c2.1}%
\end{equation}

Recall that property \textbf{(iii')} holds. In other words, the map $\pi
\mid_{g^{-1}\left(  \infty\right)  }$ is strictly decreasing. Since $\infty
=u$, this rewrites as follows: The map $\pi\mid_{g^{-1}\left(  u\right)  }$ is
strictly decreasing.

But $j=\min\left(  g^{-1}\left(  u\right)  \right)  $%
\ \ \ \ \footnote{\textit{Proof.} Assume the contrary. Thus, $j\neq\min\left(
g^{-1}\left(  u\right)  \right)  $. Combining this with $j\geq\min\left(
g^{-1}\left(  u\right)  \right)  $ (which follows from $j\in g^{-1}\left(
u\right)  $), we obtain $j>\min\left(  g^{-1}\left(  u\right)  \right)  $.
\par
But recall the following basic fact: If $S$ is a nonempty interval of $\left[
n\right]  $, and if $s\in S$ satisfies $s>\min S$, then $s-1\in S$. Applying
this to $S=g^{-1}\left(  u\right)  $ and $s=j$, we conclude that $j-1\in
g^{-1}\left(  u\right)  $ (because $g^{-1}\left(  u\right)  $ is a nonempty
interval of $\left[  n\right]  $, and because $j\in g^{-1}\left(  u\right)  $
satisfies $j>\min\left(  g^{-1}\left(  u\right)  \right)  $). Hence, both
$j-1$ and $j$ are elements of $g^{-1}\left(  u\right)  $, and satisfy $j-1<j$.
Hence, $\left(  \pi\mid_{g^{-1}\left(  u\right)  }\right)  \left(  j-1\right)
>\left(  \pi\mid_{g^{-1}\left(  u\right)  }\right)  \left(  j\right)  $ (since
the map $\pi\mid_{g^{-1}\left(  u\right)  }$ is strictly decreasing). Thus,%
\[
\pi_{j-1}=\pi\left(  j-1\right)  =\left(  \pi\mid_{g^{-1}\left(  u\right)
}\right)  \left(  j-1\right)  >\left(  \pi\mid_{g^{-1}\left(  u\right)
}\right)  \left(  j\right)  =\pi\left(  j\right)  =\pi_{j}.
\]
This contradicts $\pi_{j-1}<\pi_{j}$. This contradiction shows that our
assumption was false; qed.}. Hence,%
\begin{align*}
j  &  =\min\left(  g^{-1}\left(  u\right)  \right)  \in\left\{  \min\left(
g^{-1}\left(  u\right)  \right)  \right\}  =\operatorname*{FE}\left(
g\right)  \cap g^{-1}\left(  u\right)  \ \ \ \ \ \ \ \ \ \ \left(  \text{by
(\ref{pf.prop.Epk.fiberends.o1.pf.c2.1})}\right) \\
&  \subseteq\operatorname*{FE}\left(  g\right)  .
\end{align*}
Thus, we have shown that $j\in\operatorname*{FE}\left(  g\right)  $ in Case 2.

Let us next consider Case 3. In this case, we have $u\in\left\{
1,2,3,\ldots\right\}  $. Hence, Lemma \ref{lem.fiberends.cut} \textbf{(b)}
yields
\begin{equation}
\operatorname*{FE}\left(  g\right)  \cap g^{-1}\left(  u\right)  =\left\{
\min\left(  g^{-1}\left(  u\right)  \right)  ,\max\left(  g^{-1}\left(
u\right)  \right)  \right\}  . \label{pf.prop.Epk.fiberends.o1.pf.c3.1}%
\end{equation}

Combining $u\in g\left(  \left[  n\right]  \right)  $ with $u\in\left\{
1,2,3,\ldots\right\}  $, we obtain $u\in g\left(  \left[  n\right]  \right)
\cap\left\{  1,2,3,\ldots\right\}  $.

Recall that property \textbf{(ii')} holds. In other words, for each $h\in
g\left(  \left[  n\right]  \right)  \cap\left\{  1,2,3,\ldots\right\}  $, the
map $\pi\mid_{g^{-1}\left(  h\right)  }$ is V-shaped. Applying this to $h=u$,
we conclude that the map $\pi\mid_{g^{-1}\left(  u\right)  }$ is V-shaped.
Hence, Lemma \ref{lem.V-shaped.no-peaks} (applied to $S=g^{-1}\left(
u\right)  $) yields%
\[
g^{-1}\left(  u\right)  \cap\operatorname{Epk}\pi\subseteq\left\{
\min\left(  g^{-1}\left(  u\right)  \right)  ,\max\left(  g^{-1}\left(
u\right)  \right)  \right\}  .
\]
Combining $j\in g^{-1}\left(  u\right)  $ with $j\in\operatorname{Epk}\pi$,
we obtain%
\begin{align*}
j  &  \in g^{-1}\left(  u\right)  \cap\operatorname{Epk}\pi\subseteq\left\{
\min\left(  g^{-1}\left(  u\right)  \right)  ,\max\left(  g^{-1}\left(
u\right)  \right)  \right\} \\
&  =\operatorname*{FE}\left(  g\right)  \cap g^{-1}\left(  u\right)
\ \ \ \ \ \ \ \ \ \ \left(  \text{by (\ref{pf.prop.Epk.fiberends.o1.pf.c3.1}%
)}\right) \\
&  \subseteq\operatorname*{FE}\left(  g\right)  .
\end{align*}
Thus, we have shown that $j\in\operatorname*{FE}\left(  g\right)  $ in Case 3.

We have now proven that $j\in\operatorname*{FE}\left(  g\right)  $ in each of
the three Cases 1, 2 and 3. Hence, $j\in\operatorname*{FE}\left(  g\right)  $
always holds.

Now, forget that we fixed $j$. We thus have proven that $j\in
\operatorname*{FE}\left(  g\right)  $ for each $j\in\operatorname{Epk}\pi$.
In other words, $\operatorname{Epk}\pi\subseteq\operatorname*{FE}\left(
g\right)  $. This proves Observation 1.]

[\textit{Proof of Observation 2:} Assume that $\operatorname{Epk}\pi
\subseteq\operatorname*{FE}\left(  g\right)  $. We must prove that the map $g$
is $\pi$-amenable.

We are going to show that the four properties \textbf{(i')}, \textbf{(ii')},
\textbf{(iii')} and \textbf{(iv')} of Definition \ref{def.amen} hold. Clearly,
property \textbf{(iv')} holds (since the map $g$ is weakly increasing). Let us
now prove the remaining three properties:

[\textit{Proof of property \textbf{(i')}:} We must show that the map $\pi
\mid_{g^{-1}\left(  0\right)  }$ is strictly increasing. If $g^{-1}\left(
0\right)  =\varnothing$, then this is obvious. Thus, for the rest of this
proof, we WLOG assume that we don't have $g^{-1}\left(  0\right)
=\varnothing$. Hence, $g^{-1}\left(  0\right)  \neq\varnothing$. In other
words, there exists some $j\in\left[  n\right]  $ such that $g\left(
j\right)  =0$. Consider this $j$. Clearly, $j\geq1$ (since $j\in\left[
n\right]  $), so that $1\leq j$.

The map $g$ is weakly increasing. Hence, from $1\leq j$, we obtain $g\left(
1\right)  \preccurlyeq g\left(  j\right)  =0$. Thus, $g\left(  1\right)  =0$
(since $0$ is the smallest element of $\mathcal{N}$). Thus, $1\in
g^{-1}\left(  0\right)  $.

Recall that each fiber $g^{-1}\left(  h\right)  $ of $g$ (with $h\in
\mathcal{N}$) is an interval of the totally ordered set $\left[  n\right]  $.
Applying this to $h=0$, we conclude that $g^{-1}\left(  0\right)  $ is an
interval of the totally ordered set $\left[  n\right]  $. This interval
$g^{-1}\left(  0\right)  $ is nonempty (since $1\in g^{-1}\left(  0\right)  $).

We have $0=g\left(  \underbrace{j}_{\in\left[  n\right]  }\right)  \in
g\left(  \left[  n\right]  \right)  $. Hence, Lemma \ref{lem.fiberends.cut}
\textbf{(a)} (applied to $u=0$) yields $\operatorname*{FE}\left(  g\right)
\cap g^{-1}\left(  0\right)  =\left\{  \max\left(  g^{-1}\left(  0\right)
\right)  \right\}  $. Now,%
\[
g^{-1}\left(  0\right)  \cap\underbrace{\operatorname{Epk}\pi}_{\subseteq
\operatorname*{FE}\left(  g\right)  }\subseteq g^{-1}\left(  0\right)
\cap\operatorname*{FE}\left(  g\right)  =\operatorname*{FE}\left(  g\right)
\cap g^{-1}\left(  0\right)  =\left\{  \max\left(  g^{-1}\left(  0\right)
\right)  \right\}  .
\]
Hence, Lemma \ref{lem.Epk.V-sh} \textbf{(a)} (applied to $S=g^{-1}\left(
0\right)  $) shows that the map $\pi\mid_{g^{-1}\left(  0\right)  }$ is
strictly increasing (since $g^{-1}\left(  0\right)  $ is a nonempty interval
of the totally ordered set $\left[  n\right]  $ satisfying $1\in g^{-1}\left(
0\right)  $). This proves property \textbf{(i')}.]

[\textit{Proof of property \textbf{(ii')}:} We must show that for each $h\in
g\left(  \left[  n\right]  \right)  \cap\left\{  1,2,3,\ldots\right\}  $, the
map $\pi\mid_{g^{-1}\left(  h\right)  }$ is V-shaped. So let us fix $h\in
g\left(  \left[  n\right]  \right)  \cap\left\{  1,2,3,\ldots\right\}  $.

Then, $h\in g\left(  \left[  n\right]  \right)  \cap\left\{  1,2,3,\ldots
\right\}  \subseteq g\left(  \left[  n\right]  \right)  $. Thus, there exists
some $j\in\left[  n\right]  $ such that $h=g\left(  j\right)  $. In other
words, the fiber $g^{-1}\left(  h\right)  $ is nonempty. In other words,
$g^{-1}\left(  h\right)  \neq\varnothing$. Also, $g^{-1}\left(  h\right)  $ is
an interval of the totally ordered set $\left[  n\right]  $%
\ \ \ \ \footnote{Indeed, we have previously shown that each fiber
$g^{-1}\left(  h\right)  $ of $g$ (with $h\in\mathcal{N}$) is an interval of
the totally ordered set $\left[  n\right]  $.}.

We have $h\in g\left(  \left[  n\right]  \right)  $ and $h\in g\left(  \left[
n\right]  \right)  \cap\left\{  1,2,3,\ldots\right\}  \subseteq\left\{
1,2,3,\ldots\right\}  $. Hence, Lemma \ref{lem.fiberends.cut} \textbf{(b)}
(applied to $u=h$) yields \newline
$\operatorname*{FE}\left(  g\right)  \cap
g^{-1}\left(  h\right)  =\left\{  \min\left(  g^{-1}\left(  h\right)  \right)
,\max\left(  g^{-1}\left(  h\right)  \right)  \right\}  $. Now,%
\begin{align*}
g^{-1}\left(  h\right)  \cap\underbrace{\operatorname{Epk}\pi}_{\subseteq
\operatorname*{FE}\left(  g\right)  }
&\subseteq g^{-1}\left(  h\right)
\cap\operatorname*{FE}\left(  g\right)  =\operatorname*{FE}\left(  g\right)
\cap g^{-1}\left(  h\right) \\
& =\left\{  \min\left(  g^{-1}\left(  h\right)
\right)  ,\max\left(  g^{-1}\left(  h\right)  \right)  \right\}  .
\end{align*}
Hence, Lemma \ref{lem.Epk.V-sh} \textbf{(c)} (applied to $S=g^{-1}\left(
h\right)  $) shows that the map $\pi\mid_{g^{-1}\left(  h\right)  }$ is
V-shaped (since $g^{-1}\left(  h\right)  $ is a nonempty interval of the
totally ordered set $\left[  n\right]  $). This proves property \textbf{(ii')}.]

[\textit{Proof of property \textbf{(iii')}:} We must show that the map
$\pi\mid_{g^{-1}\left(  \infty\right)  }$ is strictly decreasing. If
$g^{-1}\left(  \infty\right)  =\varnothing$, then this is obvious. Thus, for
the rest of this proof, we WLOG assume that we don't have $g^{-1}\left(
\infty\right)  =\varnothing$. Hence, $g^{-1}\left(  \infty\right)
\neq\varnothing$. In other words, there exists some $j\in\left[  n\right]  $
such that $g\left(  j\right)  =\infty$. Consider this $j$. Clearly, $j\leq n$
(since $j\in\left[  n\right]  $).

The map $g$ is weakly increasing. Hence, from $j\leq n$, we obtain $g\left(
j\right)  \preccurlyeq g\left(  n\right)  $. In view of $g\left(  j\right)
=\infty$, this rewrites as $\infty\preccurlyeq g\left(  n\right)  $. Thus,
$g\left(  n\right)  =\infty$ (since $\infty$ is the largest element of
$\mathcal{N}$). Thus, $n\in g^{-1}\left(  \infty\right)  $.

Recall that each fiber $g^{-1}\left(  h\right)  $ of $g$ (with $h\in
\mathcal{N}$) is an interval of the totally ordered set $\left[  n\right]  $.
Applying this to $h=\infty$, we conclude that $g^{-1}\left(  \infty\right)  $
is an interval of the totally ordered set $\left[  n\right]  $. This interval
$g^{-1}\left(  \infty\right)  $ is nonempty (since $n\in g^{-1}\left(
\infty\right)  $).

We have $\infty=g\left(  \underbrace{j}_{\in\left[  n\right]  }\right)  \in
g\left(  \left[  n\right]  \right)  $. Hence, Lemma \ref{lem.fiberends.cut}
\textbf{(c)} (applied to $u=\infty$) yields $\operatorname*{FE}\left(
g\right)  \cap g^{-1}\left(  \infty\right)  =\left\{  \min\left(
g^{-1}\left(  \infty\right)  \right)  \right\}  $. Now,%
\[
g^{-1}\left(  \infty\right)  \cap\underbrace{\operatorname{Epk}\pi
}_{\subseteq\operatorname*{FE}\left(  g\right)  }\subseteq g^{-1}\left(
\infty\right)  \cap\operatorname*{FE}\left(  g\right)  =\operatorname*{FE}%
\left(  g\right)  \cap g^{-1}\left(  \infty\right)  =\left\{  \min\left(
g^{-1}\left(  \infty\right)  \right)  \right\}  .
\]
Hence, Lemma \ref{lem.Epk.V-sh} \textbf{(b)} (applied to $S=g^{-1}\left(
\infty\right)  $) shows that the map $\pi\mid_{g^{-1}\left(  \infty\right)  }$
is strictly decreasing (since $g^{-1}\left(  \infty\right)  $ is a nonempty
interval of the totally ordered set $\left[  n\right]  $ satisfying $n\in
g^{-1}\left(  \infty\right)  $). This proves property \textbf{(iii')}.]

We have now shown that the four properties \textbf{(i')}, \textbf{(ii')},
\textbf{(iii')} and \textbf{(iv')} of Definition \ref{def.amen} hold. In other
words, the map $g$ is $\pi$-amenable (by Definition \ref{def.amen}). This
proves Observation 2.]

Combining Observation 1 with Observation 2, we conclude that the map $g$ is
$\pi$-amenable if and only if $\operatorname{Epk}\pi\subseteq
\operatorname*{FE}\left(  g\right)  $. This proves Proposition
\ref{prop.Epk.fiberends}.
\end{proof}
\end{verlong}

We can rewrite Proposition \ref{prop.Epk-formula} as follows, exhibiting its
analogy with \cite[Proposition 2.2]{Stembr97}:

\begin{proposition}
\label{prop.Epk-formula2}Let $n\in\mathbb{N}$. Let $\pi$ be any $n$%
-permutation. Then,%
\[
\Gamma_{\mathcal{Z}}\left(  \pi\right)  =\sum_{\substack{g:\left[  n\right]
\rightarrow\mathcal{N}\text{ is}\\\text{weakly increasing;}%
\\\operatorname{Epk}\pi\subseteq\operatorname*{FE}\left(  g\right)
}}2^{\left\vert g\left(  \left[  n\right]  \right)  \cap\left\{
1,2,3,\ldots\right\}  \right\vert }\mathbf{x}_{g}.
\]

\end{proposition}

\begin{vershort}
\begin{proof}
[Proof of Proposition \ref{prop.Epk-formula2}.]Each $\pi$-amenable map
$g:\left[  n\right]  \rightarrow\mathcal{N}$ is weakly increasing (because of
Property \textbf{(iv')} in Definition \ref{def.amen}). Hence, Proposition
\ref{prop.Epk.fiberends} yields that the $\pi$-amenable maps $g:\left[
n\right]  \rightarrow\mathcal{N}$ are precisely the weakly increasing maps
$g:\left[  n\right]  \rightarrow\mathcal{N}$ satisfying
$\operatorname{Epk} \pi\subseteq\operatorname*{FE}\left(  g\right)  $.
Thus, Proposition \ref{prop.Epk-formula2} follows from
Proposition \ref{prop.Epk-formula}.
\end{proof}
\end{vershort}

\begin{verlong}
\begin{proof}
[Proof of Proposition \ref{prop.Epk-formula2}.]Each $\pi$-amenable map
$g:\left[  n\right]  \rightarrow\mathcal{N}$ is weakly increasing (because of
Property \textbf{(iv')} in Definition \ref{def.amen}). Thus, the $\pi
$-amenable maps $g:\left[  n\right]  \rightarrow\mathcal{N}$ are precisely the
weakly increasing maps $g:\left[  n\right]  \rightarrow\mathcal{N}$ that are
$\pi$-amenable. Hence, we have the following equality of summation signs:%
\[
\sum_{\substack{g:\left[  n\right]  \rightarrow\mathcal{N}\\\text{is }%
\pi\text{-amenable}}}=\sum_{\substack{g:\left[  n\right]  \rightarrow
\mathcal{N}\text{ is}\\\text{weakly increasing;}\\g\text{ is }\pi
\text{-amenable}}}=\sum_{\substack{g:\left[  n\right]  \rightarrow
\mathcal{N}\text{ is}\\\text{weakly increasing;}\\
\operatorname{Epk} \pi\subseteq\operatorname*{FE}\left(  g\right)  }}
\]
(because for every weakly increasing map $g:\left[  n\right]  \rightarrow
\mathcal{N}$, we have the logical equivalence%
\[
\left(  g\text{ is }\pi\text{-amenable}\right)  \ \Longleftrightarrow\ \left(
\operatorname{Epk}\pi\subseteq\operatorname*{FE}\left(  g\right)  \right)
\]
(by Proposition \ref{prop.Epk.fiberends})).

Now, Proposition \ref{prop.Epk-formula} yields%
\[
\Gamma_{\mathcal{Z}}\left(  \pi\right)  =\underbrace{\sum_{\substack{g:\left[
n\right]  \rightarrow\mathcal{N}\\\text{is }\pi\text{-amenable}}}}%
_{=\sum_{\substack{g:\left[  n\right]  \rightarrow\mathcal{N}\text{
is}\\\text{weakly increasing;}\\\operatorname{Epk}\pi\subseteq
\operatorname*{FE}\left(  g\right)  }}}2^{\left\vert g\left(  \left[
n\right]  \right)  \cap\left\{  1,2,3,\ldots\right\}  \right\vert }%
\mathbf{x}_{g}=\sum_{\substack{g:\left[  n\right]  \rightarrow\mathcal{N}%
\text{ is}\\\text{weakly increasing;}\\\operatorname{Epk}\pi\subseteq
\operatorname*{FE}\left(  g\right)  }}2^{\left\vert g\left(  \left[  n\right]
\right)  \cap\left\{  1,2,3,\ldots\right\}  \right\vert }\mathbf{x}_{g}.
\]
This proves Proposition \ref{prop.Epk-formula2}.
\end{proof}
\end{verlong}

\begin{definition}
\label{def.KnL}Let $n\in\mathbb{N}$. If $\Lambda$ is any subset of $\left[
n\right]  $, then we define a power series $K_{n,\Lambda}^{\mathcal{Z}}%
\in\operatorname*{Pow}\mathcal{N}$ by%
\begin{equation}
K_{n,\Lambda}^{\mathcal{Z}}=\sum_{\substack{g:\left[  n\right]  \rightarrow
\mathcal{N}\text{ is}\\\text{weakly increasing;}\\\Lambda\subseteq
\operatorname*{FE}\left(  g\right)  }}2^{\left\vert g\left(  \left[  n\right]
\right)  \cap\left\{  1,2,3,\ldots\right\}  \right\vert }\mathbf{x}_{g}.
\label{eq.def.KnL.1}%
\end{equation}
Thus, if $\pi$ is an $n$-permutation, then
Proposition \ref{prop.Epk-formula2} shows that
\begin{equation}
\Gamma_{\mathcal{Z}}\left(  \pi\right)
=K_{n, \operatorname{Epk} \pi}^{\mathcal{Z}}.
\label{eq.def.KnL.2}%
\end{equation}
\begin{verlong}
(Indeed, \eqref{eq.def.KnL.2} follows by comparing
\[
\Gamma_{\mathcal{Z}}\left(  \pi\right)  =\sum_{\substack{g:\left[  n\right]
\rightarrow\mathcal{N}\text{ is}\\\text{weakly increasing;}%
\\\operatorname{Epk}\pi\subseteq\operatorname*{FE}\left(  g\right)
}}2^{\left\vert g\left(  \left[  n\right]  \right)  \cap\left\{
1,2,3,\ldots\right\}  \right\vert }\mathbf{x}_{g}
\ \ \ \ \ \ \ \ \ \ \left(\text{by Proposition \ref{prop.Epk-formula2}}\right)
\]
with
\[
K_{n,\operatorname{Epk}\pi}^{\mathcal{Z}}=\sum_{\substack{g:\left[  n\right]  \rightarrow
\mathcal{N}\text{ is}\\\text{weakly increasing;}\\\operatorname{Epk}\pi\subseteq
\operatorname*{FE}\left(  g\right)  }}2^{\left\vert g\left(  \left[  n\right]
\right)  \cap\left\{  1,2,3,\ldots\right\}  \right\vert }\mathbf{x}_{g}.
\ \ \ \ \ \ \ \ \ \ \left(\text{by the definition of } K_{n,\operatorname{Epk}\pi}^{\mathcal{Z}}\right)
.
\]
)
\end{verlong}
\end{definition}

\begin{remark}
\label{rmk.KnL.concrete}Let $n\in\mathbb{N}$. Let $\Lambda$ be any subset of
$\left[  n\right]  $. It is easy to see that if $g:\left[  n\right]
\rightarrow\mathcal{N}$ is a weakly increasing map, and if $i\in\left[
n\right]  $, then $i\in\operatorname*{FE}\left(  g\right)  $ holds if and only
if we don't have $g\left(  i-1\right)  =g\left(  i\right)  =g\left(
i+1\right)  $, where we use the convention that $g\left(  0\right)  =0$ and
$g\left(  n+1\right)  =\infty$. Hence, a weakly increasing map $g:\left[
n\right]  \rightarrow\mathcal{N}$ satisfies $\Lambda\subseteq
\operatorname*{FE}\left(  g\right)  $ if and only if no $i\in\Lambda$
satisfies $g\left(  i-1\right)  =g\left(  i\right)  =g\left(  i+1\right)  $,
where we use the convention that $g\left(  0\right)  =0$ and $g\left(
n+1\right)  =\infty$. Thus, (\ref{eq.def.KnL.1}) can be rewritten as follows:%
\begin{align}
K_{n,\Lambda}^{\mathcal{Z}}  &  =\sum_{\substack{g:\left[  n\right]
\rightarrow\mathcal{N}\text{ is}\\\text{weakly increasing;}\\\text{no }%
i\in\Lambda\text{ satisfies }g\left(  i-1\right)  =g\left(  i\right)
=g\left(  i+1\right)  \\\text{(where we set }g\left(  0\right)  =0\text{ and
}g\left(  n+1\right)  =\infty\text{)}}}2^{\left\vert g\left(  \left[
n\right]  \right)  \cap\left\{  1,2,3,\ldots\right\}  \right\vert }%
\mathbf{x}_{g}\nonumber\\
&  =\sum_{\substack{\left(  g_{1},g_{2},\ldots,g_{n}\right)  \in
\mathcal{N}^{n};\\0\preccurlyeq g_{1}\preccurlyeq g_{2}\preccurlyeq
\cdots\preccurlyeq g_{n}\preccurlyeq\infty;\\\text{no }i\in\Lambda\text{
satisfies }g_{i-1}=g_{i}=g_{i+1}\\\text{(where we set }g_{0}=0\text{ and
}g_{n+1}=\infty\text{)}}}2^{\left\vert \left\{  g_{1},g_{2},\ldots
,g_{n}\right\}  \cap\left\{  1,2,3,\ldots\right\}  \right\vert }x_{g_{1}%
}x_{g_{2}}\cdots x_{g_{n}} \label{eq.rmk.KnL.concrete.2}%
\end{align}
(here, we have substituted $\left(  g_{1},g_{2},\ldots,g_{n}\right)  $ for
$\left(  g\left(  1\right)  ,g\left(  2\right)  ,\ldots,g\left(  n\right)
\right)  $ in the sum). For example,
\begin{align*}
K_{3,\left\{  1,3\right\}  }^{\mathcal{Z}}  &  =\sum_{\substack{\left(
g_{1},g_{2},g_{3}\right)  \in\mathcal{N}^{3};\\0\preccurlyeq g_{1}\preccurlyeq
g_{2}\preccurlyeq g_{3}\preccurlyeq\infty;\\\text{no }i\in\left\{
1,3\right\}  \text{ satisfies }g_{i-1}=g_{i}=g_{i+1}\\\text{(where we set
}g_{0}=0\text{ and }g_{4}=\infty\text{)}}}2^{\left\vert \left\{  g_{1}%
,g_{2},g_{3}\right\}  \cap\left\{  1,2,3,\ldots\right\}  \right\vert }%
x_{g_{1}}x_{g_{2}}x_{g_{3}}\\
&  =\sum_{\substack{\left(  g_{1},g_{2},g_{3}\right)  \in\mathcal{N}%
^{3};\\0\preccurlyeq g_{1}\preccurlyeq g_{2}\preccurlyeq g_{3}\preccurlyeq
\infty;\\\text{neither }0=g_{1}=g_{2}\text{ nor }g_{2}=g_{3}=\infty\text{
holds}}}2^{\left\vert \left\{  g_{1},g_{2},g_{3}\right\}  \cap\left\{
1,2,3,\ldots\right\}  \right\vert }x_{g_{1}}x_{g_{2}}x_{g_{3}}.
\end{align*}

As a consequence of (\ref{eq.rmk.KnL.concrete.2}), we see that if we
substitute $0$ for $x_{0}$ and for $x_{\infty}$, then $K_{n,\Lambda
}^{\mathcal{Z}}$ becomes the power series%
\begin{align*}
&  \sum_{\substack{\left(  g_{1},g_{2},\ldots,g_{n}\right)  \in\mathcal{N}%
^{n};\\0\preccurlyeq g_{1}\preccurlyeq g_{2}\preccurlyeq\cdots\preccurlyeq
g_{n}\preccurlyeq\infty;\\\text{no }i\in\Lambda\text{ satisfies }g_{i-1}%
=g_{i}=g_{i+1}\\\text{(where we set }g_{0}=0\text{ and }g_{n+1}=\infty
\text{);}\\\text{none of the }g_{i}\text{ equals }0\text{ or }\infty
}}2^{\left\vert \left\{  g_{1},g_{2},\ldots,g_{n}\right\}  \cap\left\{
1,2,3,\ldots\right\}  \right\vert }x_{g_{1}}x_{g_{2}}\cdots x_{g_{n}}\\
&  =\sum_{\substack{\left(  g_{1},g_{2},\ldots,g_{n}\right)  \in\left\{
1,2,3,\ldots\right\}  ^{n};\\g_{1}\preccurlyeq g_{2}\preccurlyeq
\cdots\preccurlyeq g_{n};\\\text{no }i\in\Lambda\setminus\left\{  1,n\right\}
\text{ satisfies }g_{i-1}=g_{i}=g_{i+1}}}2^{\left\vert \left\{  g_{1}%
,g_{2},\ldots,g_{n}\right\}  \right\vert }x_{g_{1}}x_{g_{2}}\cdots x_{g_{n}}%
\end{align*}
in the indeterminates $x_{1},x_{2},x_{3},\ldots$. This is called the
\textquotedblleft shifted quasi-symmetric function $\Theta_{\Lambda
\setminus\left\{  1,n\right\}  }^{n}\left(  X\right)  $\textquotedblright\ in
\cite[(3.2)]{BilHai95}.
\end{remark}

Corollary \ref{cor.prod2} now leads directly to the following multiplication
rule (an analogue of \cite[(3.1)]{Stembr97}):

\begin{corollary}
\label{cor.KnEpk-prodrule}Let $n\in\mathbb{N}$ and $m\in\mathbb{N}$. Let $\pi$
be an $n$-permutation. Let $\sigma$ be an $m$-permutation such that $\pi$
and $\sigma$ are disjoint. Then,
\[
K_{n,\operatorname{Epk} \pi}^{\mathcal{Z}}
\cdot K_{m,\operatorname{Epk} \sigma}^{\mathcal{Z}}
= \sum_{\tau\in S\left(  \pi,\sigma\right)  }
   K_{n+m,\operatorname{Epk}\tau}^{\mathcal{Z}}.
\]

\end{corollary}

\begin{example}
Applying Corollary \ref{cor.KnEpk-prodrule} to $n=2$, $m=1$, $\pi=\left(
1,2\right)  $ and $\sigma=\left(  3\right)  $, we obtain%
\[
K_{2,\operatorname{Epk}\left(  1,2\right)  }^{\mathcal{Z}}\cdot
K_{1,\operatorname{Epk}\left(  3\right)  }^{\mathcal{Z}}%
=K_{3,\operatorname{Epk}\left(  3,1,2\right)  }^{\mathcal{Z}}%
+K_{3,\operatorname{Epk}\left(  1,3,2\right)  }^{\mathcal{Z}}%
+K_{3,\operatorname{Epk}\left(  1,2,3\right)  }^{\mathcal{Z}}.
\]
In other words,%
\[
K_{2,\left\{  2\right\}  }^{\mathcal{Z}}\cdot K_{1,\left\{  1\right\}
}^{\mathcal{Z}}=K_{3,\left\{  1,3\right\}  }^{\mathcal{Z}}+K_{3,\left\{
2\right\}  }^{\mathcal{Z}}+K_{3,\left\{  3\right\}  }^{\mathcal{Z}}.
\]

\end{example}

\begin{proof}
[Proof of Corollary \ref{cor.KnEpk-prodrule}.]From (\ref{eq.def.KnL.2}), we
obtain
$\Gamma_{\mathcal{Z}}\left(  \pi\right)
=K_{n,\operatorname{Epk}\pi}^{\mathcal{Z}}$.
Similarly, $\Gamma_{\mathcal{Z}}\left(  \sigma\right)
=K_{m,\operatorname{Epk}\sigma}^{\mathcal{Z}}$. Multiplying these two
equalities, we obtain $\Gamma_{\mathcal{Z}}\left(  \pi\right)  \cdot
\Gamma_{\mathcal{Z}}\left(  \sigma\right)  =K_{n,\operatorname{Epk}\pi
}^{\mathcal{Z}}\cdot K_{m,\operatorname{Epk}\sigma}^{\mathcal{Z}}$. Hence,%
\begin{align*}
K_{n,\operatorname{Epk}\pi}^{\mathcal{Z}}
\cdot K_{m,\operatorname{Epk} \sigma}^{\mathcal{Z}}
&  =\Gamma_{\mathcal{Z}}\left(  \pi\right)  \cdot
\Gamma_{\mathcal{Z}}\left(  \sigma\right)
=\sum_{\tau\in S\left(  \pi
,\sigma\right)  }\underbrace{\Gamma_{\mathcal{Z}}\left(  \tau\right)
}_{\substack{=K_{n+m,\operatorname{Epk}\tau}^{\mathcal{Z}}\\\text{(by
(\ref{eq.def.KnL.2}))}}}\ \ \ \ \ \ \ \ \ \ \left(  \text{by Corollary
\ref{cor.prod2}}\right) \\
&  =\sum_{\tau\in S\left(  \pi,\sigma\right)  }K_{n+m,\operatorname{Epk}\tau
}^{\mathcal{Z}}.
\end{align*}
This proves Corollary \ref{cor.KnEpk-prodrule}.
\end{proof}

\begin{verlong}
The following lemma is a variant of the principle of inclusion and exclusion
tailored to our setting:

\begin{lemma}
\label{lem.KnL.iex}Let $n\in\mathbb{N}$. For each subset $\Lambda$ of $\left[
n\right]  $, define a power series $L_{n,\Lambda}^{\mathcal{Z}}\in
\operatorname*{Pow}\mathcal{N}$ by%
\begin{equation}
L_{n,\Lambda}^{\mathcal{Z}}=\sum_{\substack{g:\left[  n\right]  \rightarrow
\mathcal{N}\text{ is}\\\text{weakly increasing;}\\\Lambda\cap
\operatorname*{FE}\left(  g\right)  =\varnothing}}2^{\left\vert g\left(
\left[  n\right]  \right)  \cap\left\{  1,2,3,\ldots\right\}  \right\vert
}\mathbf{x}_{g}. \label{eq.lem.KnL.iex.defL}%
\end{equation}
Then:

\textbf{(a)} For each subset $\Lambda$ of $\left[  n\right]  $, we have%
\[
K_{n,\Lambda}^{\mathcal{Z}}=\sum_{Q\subseteq\Lambda}\left(  -1\right)
^{\left\vert Q\right\vert }L_{n,Q}^{\mathcal{Z}}.
\]

\textbf{(b)} For each subset $\Lambda$ of $\left[  n\right]  $, we have%
\[
L_{n,\Lambda}^{\mathcal{Z}}=\sum_{Q\subseteq\Lambda}\left(  -1\right)
^{\left\vert Q\right\vert }K_{n,Q}^{\mathcal{Z}}.
\]

\end{lemma}

\begin{proof}
[Proof of Lemma \ref{lem.KnL.iex}.]Recall the following known fact: If $R$ is
a finite set, then%
\begin{equation}
\sum_{Q\subseteq R}\left(  -1\right)  ^{\left\vert Q\right\vert }=%
\begin{cases}
1, & \text{if }R=\varnothing;\\
0, & \text{if }R\neq\varnothing
\end{cases}
. \label{pf.lem.KnL.iex.1}%
\end{equation}

Let $g:\left[  n\right]  \rightarrow \mathcal{N}$ be a weakly increasing map.
Let $\Lambda$ be a subset of $\left[  n\right]  $. The subsets $Q$ of
$\Lambda$ satisfying $Q\cap\operatorname*{FE}\left(  g\right)  =\varnothing$
are precisely the subsets of $\Lambda\setminus\operatorname*{FE}\left(
g\right)  $. Therefore,
\begin{align}
\sum_{\substack{Q\subseteq\Lambda;\\Q\cap\operatorname*{FE}\left(  g\right)
=\varnothing}}\left(  -1\right)  ^{\left\vert Q\right\vert }  &
=\sum_{Q\subseteq\Lambda\setminus\operatorname*{FE}\left(  g\right)  }\left(
-1\right)  ^{\left\vert Q\right\vert }=%
\begin{cases}
1, & \text{if }\Lambda\setminus\operatorname*{FE}\left(  g\right)
=\varnothing;\\
0, & \text{if }\Lambda\setminus\operatorname*{FE}\left(  g\right)
\neq\varnothing
\end{cases}
\ \ \ \ \ \ \ \ \ \ \left(  \text{by (\ref{pf.lem.KnL.iex.1})}\right)
\nonumber\\
&  =%
\begin{cases}
1, & \text{if }\Lambda\subseteq\operatorname*{FE}\left(  g\right)  ;\\
0, & \text{if }\Lambda\not \subseteq \operatorname*{FE}\left(  g\right)
\end{cases}
. \label{pf.lem.KnL.iex.2}%
\end{align}
Also, the subsets $Q$ of $\Lambda$ satisfying $Q\subseteq\operatorname*{FE}%
\left(  g\right)  $ are precisely the subsets of $\Lambda\cap
\operatorname*{FE}\left(  g\right)  $. Thus,%
\begin{equation}
\sum_{\substack{Q\subseteq\Lambda;\\Q\subseteq\operatorname*{FE}\left(
g\right)  }}\left(  -1\right)  ^{\left\vert Q\right\vert }=\sum_{Q\subseteq
\Lambda\cap\operatorname*{FE}\left(  g\right)  }\left(  -1\right)
^{\left\vert Q\right\vert }=%
\begin{cases}
1, & \text{if }\Lambda\cap\operatorname*{FE}\left(  g\right)  =\varnothing;\\
0, & \text{if }\Lambda\cap\operatorname*{FE}\left(  g\right)  \neq\varnothing
\end{cases}
\label{pf.lem.KnL.iex.3}%
\end{equation}
(by (\ref{pf.lem.KnL.iex.1})).

Now, forget that we fixed $g$ and $\Lambda$. We thus have proven
(\ref{pf.lem.KnL.iex.2}) and (\ref{pf.lem.KnL.iex.3}) for each
weakly increasing map $g:\left[  n\right]  \rightarrow \mathcal{N}$
and each subset $\Lambda$ of $\left[  n\right]  $.

\textbf{(a)} Let $\Lambda$ be a subset of $\left[  n\right]  $. Then,%
\begin{align*}
&  \sum_{Q\subseteq\Lambda}\left(  -1\right)  ^{\left\vert Q\right\vert
}\underbrace{L_{n,Q}^{\mathcal{Z}}}_{\substack{=\sum_{\substack{g:\left[
n\right]  \rightarrow\mathcal{N}\text{ is}\\\text{weakly increasing;}%
\\Q\cap\operatorname*{FE}\left(  g\right)  =\varnothing}}2^{\left\vert
g\left(  \left[  n\right]  \right)  \cap\left\{  1,2,3,\ldots\right\}
\right\vert }\mathbf{x}_{g}\\\text{(by the definition of }L_{n,Q}%
^{\mathcal{Z}}\text{)}}}\\
&  =\sum_{Q\subseteq\Lambda}\left(  -1\right)  ^{\left\vert Q\right\vert }%
\sum_{\substack{g:\left[  n\right]  \rightarrow\mathcal{N}\text{
is}\\\text{weakly increasing;}\\Q\cap\operatorname*{FE}\left(  g\right)
=\varnothing}}2^{\left\vert g\left(  \left[  n\right]  \right)  \cap\left\{
1,2,3,\ldots\right\}  \right\vert }\mathbf{x}_{g}\\
&  =\sum_{\substack{g:\left[  n\right]  \rightarrow\mathcal{N}\text{
is}\\\text{weakly increasing}}}\underbrace{\left(  \sum_{\substack{Q\subseteq
\Lambda;\\Q\cap\operatorname*{FE}\left(  g\right)  =\varnothing}}\left(
-1\right)  ^{\left\vert Q\right\vert }\right)  }_{\substack{=%
\begin{cases}
1, & \text{if }\Lambda\subseteq\operatorname*{FE}\left(  g\right)  ;\\
0, & \text{if }\Lambda\not \subseteq \operatorname*{FE}\left(  g\right)
\end{cases}
\\\text{(by (\ref{pf.lem.KnL.iex.2}))}}}2^{\left\vert g\left(  \left[
n\right]  \right)  \cap\left\{  1,2,3,\ldots\right\}  \right\vert }%
\mathbf{x}_{g}\\
&  =\sum_{\substack{g:\left[  n\right]  \rightarrow\mathcal{N}\text{
is}\\\text{weakly increasing}}}%
\begin{cases}
1, & \text{if }\Lambda\subseteq\operatorname*{FE}\left(  g\right)  ;\\
0, & \text{if }\Lambda\not \subseteq \operatorname*{FE}\left(  g\right)
\end{cases}
2^{\left\vert g\left(  \left[  n\right]  \right)  \cap\left\{  1,2,3,\ldots
\right\}  \right\vert }\mathbf{x}_{g}\\
&  =\sum_{\substack{g:\left[  n\right]  \rightarrow\mathcal{N}\text{
is}\\\text{weakly increasing;}\\\Lambda\subseteq\operatorname*{FE}\left(
g\right)  }}2^{\left\vert g\left(  \left[  n\right]  \right)  \cap\left\{
1,2,3,\ldots\right\}  \right\vert }\mathbf{x}_{g}=K_{n,\Lambda}^{\mathcal{Z}}%
\end{align*}
(by (\ref{eq.def.KnL.1})). This proves Lemma \ref{lem.KnL.iex} \textbf{(a)}.

% \begin{vershort}
% \textbf{(b)} The proof of Lemma \ref{lem.KnL.iex} \textbf{(b)} is analogous,
% except that the roles of $K_{n,Q}^{\mathcal{Z}}$ and $L_{n,Q}^{\mathcal{Z}}$
% are interchanged and we need to use (\ref{pf.lem.KnL.iex.3}) instead of
% (\ref{pf.lem.KnL.iex.2}).
% \end{vershort}
% 
% \begin{verlong}
\textbf{(b)} Let $\Lambda$ be a subset of $\left[  n\right]  $. Then,%
\begin{align*}
&  \sum_{Q\subseteq\Lambda}\left(  -1\right)  ^{\left\vert Q\right\vert
}\underbrace{K_{n,Q}^{\mathcal{Z}}}_{\substack{=\sum_{\substack{g:\left[
n\right]  \rightarrow\mathcal{N}\text{ is}\\\text{weakly increasing;}%
\\Q\subseteq\operatorname*{FE}\left(  g\right)  }}2^{\left\vert g\left(
\left[  n\right]  \right)  \cap\left\{  1,2,3,\ldots\right\}  \right\vert
}\mathbf{x}_{g}\\\text{(by the definition of }K_{n,Q}^{\mathcal{Z}}\text{)}%
}}\\
&  =\sum_{Q\subseteq\Lambda}\left(  -1\right)  ^{\left\vert Q\right\vert }%
\sum_{\substack{g:\left[  n\right]  \rightarrow\mathcal{N}\text{
is}\\\text{weakly increasing;}\\Q\subseteq\operatorname*{FE}\left(  g\right)
}}2^{\left\vert g\left(  \left[  n\right]  \right)  \cap\left\{
1,2,3,\ldots\right\}  \right\vert }\mathbf{x}_{g}\\
&  =\sum_{\substack{g:\left[  n\right]  \rightarrow\mathcal{N}\text{
is}\\\text{weakly increasing}}}\underbrace{\left(  \sum_{\substack{Q\subseteq
\Lambda;\\Q\subseteq\operatorname*{FE}\left(  g\right)  }}\left(  -1\right)
^{\left\vert Q\right\vert }\right)  }_{\substack{=%
\begin{cases}
1, & \text{if }\Lambda\cap\operatorname*{FE}\left(  g\right)  =\varnothing;\\
0, & \text{if }\Lambda\cap\operatorname*{FE}\left(  g\right)  \neq\varnothing
\end{cases}
\\\text{(by (\ref{pf.lem.KnL.iex.3}))}}}2^{\left\vert g\left(  \left[
n\right]  \right)  \cap\left\{  1,2,3,\ldots\right\}  \right\vert }%
\mathbf{x}_{g}\\
&  =\sum_{\substack{g:\left[  n\right]  \rightarrow\mathcal{N}\text{
is}\\\text{weakly increasing}}}%
\begin{cases}
1, & \text{if }\Lambda\cap\operatorname*{FE}\left(  g\right)  =\varnothing;\\
0, & \text{if }\Lambda\cap\operatorname*{FE}\left(  g\right)  \neq\varnothing
\end{cases}
\ \ 2^{\left\vert g\left(  \left[  n\right]  \right)  \cap\left\{  1,2,3,
\ldots\right\}  \right\vert }\mathbf{x}_{g}\\
&  =\sum_{\substack{g:\left[  n\right]  \rightarrow\mathcal{N}\text{
is}\\\text{weakly increasing;}\\\Lambda\cap\operatorname*{FE}\left(  g\right)
=\varnothing}}2^{\left\vert g\left(  \left[  n\right]  \right)  \cap\left\{
1,2,3,\ldots\right\}  \right\vert }\mathbf{x}_{g}=L_{n,\Lambda}^{\mathcal{Z}}%
\end{align*}
(by (\ref{eq.lem.KnL.iex.defL})). This proves Lemma \ref{lem.KnL.iex}
\textbf{(b)}.
% \end{verlong}
\end{proof}
\end{verlong}

Recall Definition~\ref{def.lac.Ln}.

\begin{proposition}
\label{prop.KnL.lindep}Let $n\in\mathbb{N}$. Then, the family%
\[
\left(  K_{n,\Lambda}^{\mathcal{Z}}\right)  _{\Lambda\in \mathbf{L}_n}%
\]
is $\mathbb{Q}$-linearly independent.
\end{proposition}

\begin{vershort}
Our proof of Proposition \ref{prop.KnL.lindep} requires
the following definition:
\end{vershort}

\begin{verlong}
We shall give two proofs of Proposition \ref{prop.KnL.lindep}. The first one
relies on studying the coefficients of $K_{n,\Lambda}^{\mathcal{Z}}$; it needs
the following definition:
\end{verlong}

\begin{definition}
Let $\mathfrak{m}$ be any monomial in $\operatorname*{Pow}\mathcal{N}$ (that
is, a formal commutative product of indeterminates $x_{h}$ with $h\in
\mathcal{N}$). Let $f\in\operatorname*{Pow}\mathcal{N}$. Then, $\left[
\mathfrak{m}\right]  \left(  f\right)  $ shall mean the coefficient of
$\mathfrak{m}$ in the power series $f$. (For example, $\left[  x_{0}^{2}%
x_{3}\right]  \left(  3+5x_{0}^{2}x_{3}+6x_{0}+9x_{\infty}\right)  =5$ and
$\left[  x_{0}^{2}x_{3}\right]  \left(  x_{1}-x_{\infty}\right)  =0$.)
\end{definition}

\begin{vershort}
\begin{lemma}
\label{lem.KnL.coeff}Let $n\in\mathbb{N}$.

\textbf{(a)} If $g$ and $h$ are two weakly increasing maps $\left[  n\right]
\rightarrow\mathcal{N}$, then
\[
\left(\text{we have $\xx_g = \xx_h$ if and only if $g = h$}\right) .
\]

\textbf{(b)} Let $R\in\mathbf{L}_{n}$. Let $h:\left[  n\right]  \rightarrow
\mathcal{N}$ be a weakly increasing map. Then,%
\[
\left[  \mathbf{x}_{h}\right]  \left(  K_{n,R}^{\mathcal{Z}}\right)  =%
\begin{cases}
2^{\left\vert h\left(  \left[  n\right]  \right)  \cap\left\{  1,2,3,\ldots
\right\}  \right\vert }, & \text{if }R\subseteq\operatorname*{FE}\left(
h\right)  ;\\
0, & \text{otherwise}%
\end{cases}
.
\]

\end{lemma}

\begin{proof}
[Proof of Lemma \ref{lem.KnL.coeff}.]\textbf{(a)} A weakly increasing map
$g:\left[  n\right]  \rightarrow\mathcal{N}$ can be uniquely reconstructed
from the \textbf{multiset} $\left\{  g\left(  1\right)  ,g\left(  2\right)
,\ldots,g\left(  n\right)  \right\}  _{\operatorname*{multi}}$ of its values
(because it is weakly increasing, so there is only one way in which these
values can be ordered). Hence, a weakly increasing map $g:\left[  n\right]
\rightarrow\mathcal{N}$ can be uniquely reconstructed from the monomial
$\mathbf{x}_{g}$ (since this monomial
$\mathbf{x}_{g}=x_{g\left(1\right)}x_{g\left(2\right)}\cdots
x_{g\left(n\right)}$
encodes the multiset $\left\{  g\left(  1\right)  ,g\left(
2\right)  ,\ldots,g\left(  n\right)  \right\}  _{\operatorname*{multi}}$). In
other words, if $g$ and $h$ are two weakly increasing maps $\left[  n\right]
\rightarrow\mathcal{N}$, then $\mathbf{x}_{g}=\mathbf{x}_{h}$ holds if and
only if $g=h$. This proves Lemma \ref{lem.KnL.coeff} \textbf{(a)}.

\textbf{(b)} The definition of $K_{n,R}^{\mathcal{Z}}$ yields%
\[
K_{n,R}^{\mathcal{Z}}=\sum_{\substack{g:\left[  n\right]  \rightarrow
\mathcal{N}\text{ is}\\\text{weakly increasing;}\\R\subseteq\operatorname*{FE}%
\left(  g\right)  }}2^{\left\vert g\left(  \left[  n\right]  \right)
\cap\left\{  1,2,3,\ldots\right\}  \right\vert }\mathbf{x}_{g}.
\]
Thus,%
\begin{align*}
\left[  \mathbf{x}_{h}\right]  \left(  K_{n,R}^{\mathcal{Z}}\right)   &
=\left[  \mathbf{x}_{h}\right]  \left(  \sum_{\substack{g:\left[  n\right]
\rightarrow\mathcal{N}\text{ is}\\\text{weakly increasing;}\\R\subseteq
\operatorname*{FE}\left(  g\right)  }}2^{\left\vert g\left(  \left[  n\right]
\right)  \cap\left\{  1,2,3,\ldots\right\}  \right\vert }\mathbf{x}_{g}\right)
\\
&  =\sum_{\substack{g:\left[  n\right]  \rightarrow\mathcal{N}\text{
is}\\\text{weakly increasing;}\\R\subseteq\operatorname*{FE}\left(  g\right)
}}2^{\left\vert g\left(  \left[  n\right]  \right)  \cap\left\{
1,2,3,\ldots\right\}  \right\vert }\underbrace{\left[  \mathbf{x}_{h}\right]
\left(  \mathbf{x}_{g}\right)  }_{\substack{=%
\begin{cases}
1, & \text{if }\mathbf{x}_{g}=\mathbf{x}_{h};\\
0, & \text{if }\mathbf{x}_{g}\neq\mathbf{x}_{h}%
\end{cases}
\\\text{(since }\mathbf{x}_{h}\text{ and }\mathbf{x}_{g}\text{ are two
monomials)}}}\\
&  =\sum_{\substack{g:\left[  n\right]  \rightarrow\mathcal{N}\text{
is}\\\text{weakly increasing;}\\R\subseteq\operatorname*{FE}\left(  g\right)
}}2^{\left\vert g\left(  \left[  n\right]  \right)  \cap\left\{
1,2,3,\ldots\right\}  \right\vert }\underbrace{%
\begin{cases}
1, & \text{if }\mathbf{x}_{g}=\mathbf{x}_{h};\\
0, & \text{if }\mathbf{x}_{g}\neq\mathbf{x}_{h}%
\end{cases}
}_{\substack{=%
\begin{cases}
1, & \text{if }g=h;\\
0, & \text{if }g\neq h
\end{cases}
\\\text{(by Lemma \ref{lem.KnL.coeff} \textbf{(a)})}}}\\
&  =\sum_{\substack{g:\left[  n\right]  \rightarrow\mathcal{N}\text{
is}\\\text{weakly increasing;}\\R\subseteq\operatorname*{FE}\left(  g\right)
}}2^{\left\vert g\left(  \left[  n\right]  \right)  \cap\left\{
1,2,3,\ldots\right\}  \right\vert }%
\begin{cases}
1, & \text{if }g=h;\\
0, & \text{if }g\neq h
\end{cases}
\\
&  =\sum_{\substack{g:\left[  n\right]  \rightarrow\mathcal{N}\text{
is}\\\text{weakly increasing;}\\R\subseteq\operatorname*{FE}\left(  g\right)
;\\g=h}}2^{\left\vert g\left(  \left[  n\right]  \right)  \cap\left\{
1,2,3,\ldots\right\}  \right\vert }.
\end{align*}
The sum on the right hand side of this equality has a unique addend (namely,
its addend for $g=h$, which is $2^{\left\vert h\left(  \left[  n\right]
\right)  \cap\left\{  1,2,3,\ldots\right\}  \right\vert }$) when
$R\subseteq\operatorname*{FE}\left(  h\right)  $; otherwise it is an empty
sum. Hence, this sum simplifies as follows:%
\[
\sum_{\substack{g:\left[  n\right]  \rightarrow\mathcal{N}\text{
is}\\\text{weakly increasing;}\\R\subseteq\operatorname*{FE}\left(  g\right)
;\\g=h}}2^{\left\vert g\left(  \left[  n\right]  \right)  \cap\left\{
1,2,3,\ldots\right\}  \right\vert }=%
\begin{cases}
2^{\left\vert h\left(  \left[  n\right]  \right)  \cap\left\{  1,2,3,\ldots
\right\}  \right\vert }, & \text{if }R\subseteq\operatorname*{FE}\left(
h\right)  ;\\
0, & \text{otherwise}%
\end{cases}
.
\]
Hence,%
\[
\left[  \mathbf{x}_{h}\right]  \left(  K_{n,R}^{\mathcal{Z}}\right)
=\sum_{\substack{g:\left[  n\right]  \rightarrow\mathcal{N}\text{
is}\\\text{weakly increasing;}\\R\subseteq\operatorname*{FE}\left(  g\right)
;\\g=h}}2^{\left\vert g\left(  \left[  n\right]  \right)  \cap\left\{
1,2,3,\ldots\right\}  \right\vert }=%
\begin{cases}
2^{\left\vert h\left(  \left[  n\right]  \right)  \cap\left\{  1,2,3,\ldots
\right\}  \right\vert }, & \text{if }R\subseteq\operatorname*{FE}\left(
h\right)  ;\\
0, & \text{otherwise}%
\end{cases}
.
\]
This proves Lemma \ref{lem.KnL.coeff} \textbf{(b)}.
\end{proof}
\end{vershort}

\begin{verlong}
\begin{lemma}
\label{lem.KnL.coeff}Let $n\in\mathbb{N}$.

\textbf{(a)} If $g$ and $h$ are two weakly increasing maps $\left[  n\right]
\rightarrow\mathcal{N}$, then
\[%
\begin{cases}
1, & \text{if }\mathbf{x}_{g}=\mathbf{x}_{h};\\
0, & \text{if }\mathbf{x}_{g}\neq\mathbf{x}_{h}%
\end{cases}
=%
\begin{cases}
1, & \text{if }g=h;\\
0, & \text{if }g\neq h
\end{cases}
.
\]

\textbf{(b)} Let $R\in\mathbf{L}_{n}$. Let $h:\left[  n\right]  \rightarrow
\mathcal{N}$ be a weakly increasing map. Then,%
\[
\left[  \mathbf{x}_{h}\right]  \left(  K_{n,R}^{\mathcal{Z}}\right)  =%
\begin{cases}
2^{\left\vert h\left(  \left[  n\right]  \right)  \cap\left\{  1,2,3,\ldots
\right\}  \right\vert }, & \text{if }R\subseteq\operatorname*{FE}\left(
h\right)  ;\\
0, & \text{otherwise}%
\end{cases}
.
\]

\end{lemma}

\begin{proof}
[Proof of Lemma \ref{lem.KnL.coeff}.]\textbf{(a)} A weakly increasing map
$g:\left[  n\right]  \rightarrow\mathcal{N}$ can be uniquely reconstructed
from the \textbf{multiset} $\left\{  g\left(  1\right)  ,g\left(  2\right)
,\ldots,g\left(  n\right)  \right\}  _{\operatorname*{multi}}$ of its values
(because it is weakly increasing, so there is only one way in which these
values can be ordered). Hence, a weakly increasing map $g:\left[  n\right]
\rightarrow\mathcal{N}$ can be uniquely reconstructed from the monomial
$\mathbf{x}_{g}$ (since this monomial
$\mathbf{x}_{g}=x_{g\left(1\right)}x_{g\left(2\right)}\cdots
x_{g\left(n\right)}$
encodes the multiset $\left\{  g\left(  1\right)  ,g\left(
2\right)  ,\ldots,g\left(  n\right)  \right\}  _{\operatorname*{multi}}$). In
other words, if $g$ and $h$ are two weakly increasing maps $\left[  n\right]
\rightarrow\mathcal{N}$, then $\mathbf{x}_{g}=\mathbf{x}_{h}$ holds if and
only if $g=h$. Hence, if $g$ and $h$ are two weakly increasing maps $\left[
n\right]  \rightarrow\mathcal{N}$, then
\[%
\begin{cases}
1, & \text{if }\mathbf{x}_{g}=\mathbf{x}_{h};\\
0, & \text{if }\mathbf{x}_{g}\neq\mathbf{x}_{h}%
\end{cases}
=%
\begin{cases}
1, & \text{if }g=h;\\
0, & \text{if }g\neq h
\end{cases}
.
\]
This proves Lemma \ref{lem.KnL.coeff} \textbf{(a)}.

\textbf{(b)} The definition of $K_{n,R}^{\mathcal{Z}}$ yields%
\[
K_{n,R}^{\mathcal{Z}}=\sum_{\substack{g:\left[  n\right]  \rightarrow
\mathcal{N}\text{ is}\\\text{weakly increasing;}\\R\subseteq\operatorname*{FE}%
\left(  g\right)  }}2^{\left\vert g\left(  \left[  n\right]  \right)
\cap\left\{  1,2,3,\ldots\right\}  \right\vert }\mathbf{x}_{g}.
\]
Thus,%
\begin{align*}
\left[  \mathbf{x}_{h}\right]  \left(  K_{n,R}^{\mathcal{Z}}\right)   &
=\left[  \mathbf{x}_{h}\right]  \left(  \sum_{\substack{g:\left[  n\right]
\rightarrow\mathcal{N}\text{ is}\\\text{weakly increasing;}\\R\subseteq
\operatorname*{FE}\left(  g\right)  }}2^{\left\vert g\left(  \left[  n\right]
\right)  \cap\left\{  1,2,3,\ldots\right\}  \right\vert }\mathbf{x}_{g}\right)
\\
&  =\sum_{\substack{g:\left[  n\right]  \rightarrow\mathcal{N}\text{
is}\\\text{weakly increasing;}\\R\subseteq\operatorname*{FE}\left(  g\right)
}}2^{\left\vert g\left(  \left[  n\right]  \right)  \cap\left\{
1,2,3,\ldots\right\}  \right\vert }\underbrace{\left[  \mathbf{x}_{h}\right]
\left(  \mathbf{x}_{g}\right)  }_{\substack{=%
\begin{cases}
1, & \text{if }\mathbf{x}_{g}=\mathbf{x}_{h};\\
0, & \text{if }\mathbf{x}_{g}\neq\mathbf{x}_{h}%
\end{cases}
\\\text{(since }\mathbf{x}_{h}\text{ and }\mathbf{x}_{g}\text{ are two
monomials)}}}\\
&  =\sum_{\substack{g:\left[  n\right]  \rightarrow\mathcal{N}\text{
is}\\\text{weakly increasing;}\\R\subseteq\operatorname*{FE}\left(  g\right)
}}2^{\left\vert g\left(  \left[  n\right]  \right)  \cap\left\{
1,2,3,\ldots\right\}  \right\vert }\underbrace{%
\begin{cases}
1, & \text{if }\mathbf{x}_{g}=\mathbf{x}_{h};\\
0, & \text{if }\mathbf{x}_{g}\neq\mathbf{x}_{h}%
\end{cases}
}_{\substack{=%
\begin{cases}
1, & \text{if }g=h;\\
0, & \text{if }g\neq h
\end{cases}
\\\text{(by Lemma \ref{lem.KnL.coeff} \textbf{(a)})}}}\\
&  =\sum_{\substack{g:\left[  n\right]  \rightarrow\mathcal{N}\text{
is}\\\text{weakly increasing;}\\R\subseteq\operatorname*{FE}\left(  g\right)
}}2^{\left\vert g\left(  \left[  n\right]  \right)  \cap\left\{
1,2,3,\ldots\right\}  \right\vert }%
\begin{cases}
1, & \text{if }g=h;\\
0, & \text{if }g\neq h
\end{cases}
\\
&  =\sum_{\substack{g:\left[  n\right]  \rightarrow\mathcal{N}\text{
is}\\\text{weakly increasing;}\\R\subseteq\operatorname*{FE}\left(  g\right)
;\\g=h}}2^{\left\vert g\left(  \left[  n\right]  \right)  \cap\left\{
1,2,3,\ldots\right\}  \right\vert }.
\end{align*}
The sum on the right hand side of this equality has a unique addend (namely,
its addend for $g=h$, which is $2^{\left\vert h\left(  \left[  n\right]
\right)  \cap\left\{  1,2,3,\ldots\right\}  \right\vert }$) when
$R\subseteq\operatorname*{FE}\left(  h\right)  $; otherwise it is an empty
sum. Hence, this sum simplifies as follows:%
\[
\sum_{\substack{g:\left[  n\right]  \rightarrow\mathcal{N}\text{
is}\\\text{weakly increasing;}\\R\subseteq\operatorname*{FE}\left(  g\right)
;\\g=h}}2^{\left\vert g\left(  \left[  n\right]  \right)  \cap\left\{
1,2,3,\ldots\right\}  \right\vert }=%
\begin{cases}
2^{\left\vert h\left(  \left[  n\right]  \right)  \cap\left\{  1,2,3,\ldots
\right\}  \right\vert }, & \text{if }R\subseteq\operatorname*{FE}\left(
h\right)  ;\\
0, & \text{otherwise}%
\end{cases}
.
\]
Hence,%
\[
\left[  \mathbf{x}_{h}\right]  \left(  K_{n,R}^{\mathcal{Z}}\right)
=\sum_{\substack{g:\left[  n\right]  \rightarrow\mathcal{N}\text{
is}\\\text{weakly increasing;}\\R\subseteq\operatorname*{FE}\left(  g\right)
;\\g=h}}2^{\left\vert g\left(  \left[  n\right]  \right)  \cap\left\{
1,2,3,\ldots\right\}  \right\vert }=%
\begin{cases}
2^{\left\vert h\left(  \left[  n\right]  \right)  \cap\left\{  1,2,3,\ldots
\right\}  \right\vert }, & \text{if }R\subseteq\operatorname*{FE}\left(
h\right)  ;\\
0, & \text{otherwise}%
\end{cases}
.
\]
This proves Lemma \ref{lem.KnL.coeff} \textbf{(b)}.
\end{proof}
\end{verlong}

\begin{proof}
[First proof of Proposition \ref{prop.KnL.lindep}.]Recall Definition
\ref{def.order-on-P}. In the following, we shall regard the set $\mathbf{P}$
as a totally ordered set, equipped with the order from Proposition
\ref{prop.lac.order}.

\begin{vershort}
Clearly, $\mathbf{L}_{n}\subseteq\mathbf{P}$. Hence, we consider
$\mathbf{L}_{n}$ as a totally ordered set, whose total order is inherited from
$\mathbf{P}$.
\end{vershort}

\begin{verlong}
Clearly, $\mathbf{L}_{n}$ is a set of subsets of $\left[  n\right]  $, and
thus a set of finite subsets of $\mathbb{Z}$. In other words, $\mathbf{L}_{n}$
is a subset of $\mathbf{P}$. Hence, we consider $\mathbf{L}_{n}$ as a totally
ordered set, whose total order is inherited from $\mathbf{P}$.
\end{verlong}

Let $\left(  a_{R}\right)  _{R\in\mathbf{L}_{n}}\in\mathbb{Q}^{\mathbf{L}_{n}%
}$ be a family of scalars (in $\mathbb{Q}$) such that $\sum_{R\in
\mathbf{L}_{n}}a_{R}K_{n,R}^{\mathcal{Z}}=0$. We are going to show that
$\left(  a_{R}\right)  _{R\in\mathbf{L}_{n}}=\left(  0\right)  _{R\in
\mathbf{L}_{n}}$.

Indeed, assume the contrary. Thus, $\left(  a_{R}\right)  _{R\in\mathbf{L}%
_{n}}\neq\left(  0\right)  _{R\in\mathbf{L}_{n}}$. Hence, there exists some
$R\in\mathbf{L}_{n}$ such that $a_{R}\neq0$. Let $\Lambda$ be the
\textbf{largest} such $R$ (with respect to the total order on $\mathbf{L}_{n}$
we have introduced above). Hence, $\Lambda$ is an element of $\mathbf{L}_{n}$
and satisfies $a_{\Lambda}\neq0$; but every element $R\in\mathbf{L}_{n}$
satisfying $R>\Lambda$ must satisfy%
\begin{equation}
a_{R}=0. \label{pf.prop.KnL.lindep.larger=0}%
\end{equation}

\begin{verlong}
We have $\Lambda \in \mathbf{L}_n$. Thus, $\Lambda$ is a subset of
$\left[ n\right]$ (since $\mathbf{L}_n$ is a set of subsets of
$\left[ n\right]$). In other words, $\Lambda \subseteq \left[n\right]$.
\end{verlong}

Lemma \ref{lem.FE.exist} shows that there exists a weakly increasing map
$g:\left[  n\right]  \rightarrow\mathcal{N}$ such that $\operatorname*{FE}%
\left(  g\right)  =\left(  \Lambda\cup\left(  \Lambda+1\right)  \right)
\cap\left[  n\right]  $. Consider this $g$. Combining $\Lambda\subseteq
\Lambda\cup\left(  \Lambda+1\right)  $ with $\Lambda\subseteq\left[  n\right]
$, we obtain
\[
\Lambda\subseteq\left(  \Lambda\cup\left(  \Lambda+1\right)  \right)
\cap\left[  n\right]  =\operatorname*{FE}\left(  g\right)  .
\]

For every $R\in\mathbf{L}_{n}$ satisfying $R\neq\Lambda$, we have%
\begin{equation}
\left[  \mathbf{x}_{g}\right]  \left(  a_{R}K_{n,R}^{\mathcal{Z}}\right)  =0.
\label{pf.prop.KnL.lindep.coeff=0}%
\end{equation}

[\textit{Proof of (\ref{pf.prop.KnL.lindep.coeff=0}):} Let $R\in\mathbf{L}%
_{n}$ be such that $R\neq\Lambda$. We must prove
(\ref{pf.prop.KnL.lindep.coeff=0}).

\begin{vershort}
Assume the contrary. Thus, $\left[  \mathbf{x}_{g}\right]  \left(
a_{R}K_{n,R}^{\mathcal{Z}}\right)  \neq0$. In other words, $a_{R}\left[
\mathbf{x}_{g}\right]  \left(  K_{n,R}^{\mathcal{Z}}\right)  \neq0$. Hence,
$a_{R}\neq0$ and $\left[  \mathbf{x}_{g}\right]  \left(  K_{n,R}^{\mathcal{Z}%
}\right)  \neq0$.

From the definition of $\mathbf{L}_{n}$, it follows easily that every element
of $\mathbf{L}_{n}$ is a lacunar subset of $\left[  n\right]  $. Hence, $R$ is
a lacunar subset of $\left[  n\right]  $ (since $R\in\mathbf{L}_{n}$).

But Lemma \ref{lem.KnL.coeff} \textbf{(b)} (applied to $h=g$) yields%
\[
\left[  \mathbf{x}_{g}\right]  \left(  K_{n,R}^{\mathcal{Z}}\right)  =%
\begin{cases}
2^{\left\vert g\left(  \left[  n\right]  \right)  \cap\left\{  1,2,3,\ldots
\right\}  \right\vert }, & \text{if }R\subseteq\operatorname*{FE}\left(
g\right)  ;\\
0, & \text{otherwise}%
\end{cases}
.
\]
Hence, $\left[  \mathbf{x}_{g}\right]  \left(  K_{n,R}^{\mathcal{Z}}\right)
=0$ if $R\not \subseteq \operatorname*{FE}\left(  g\right)  $. Thus, we cannot
have $R\not \subseteq \operatorname*{FE}\left(  g\right)  $ (since $\left[
\mathbf{x}_{g}\right]  \left(  K_{n,R}^{\mathcal{Z}}\right)  \neq0$).
Therefore, we have
\[
R\subseteq\operatorname*{FE}\left(  g\right)  =\left(  \Lambda\cup\left(
\Lambda+1\right)  \right)  \cap\left[  n\right]  \subseteq\Lambda\cup\left(
\Lambda+1\right)  .
\]
Thus, Proposition \ref{prop.lac.RL} yields that $R\geq\Lambda$. Combining this
with $R\neq\Lambda$, we obtain $R>\Lambda$. Hence,
(\ref{pf.prop.KnL.lindep.larger=0}) yields $a_{R}=0$. This contradicts
$a_{R}\neq0$. This contradiction shows that our assumption was wrong. Hence,
(\ref{pf.prop.KnL.lindep.coeff=0}) is proven.]
\end{vershort}

\begin{verlong}
Assume the contrary. Thus, $\left[  \mathbf{x}_{g}\right]  \left(
a_{R}K_{n,R}^{\mathcal{Z}}\right)  \neq0$. Therefore, $a_{R}\left[
\mathbf{x}_{g}\right]  \left(  K_{n,R}^{\mathcal{Z}}\right)  =\left[
\mathbf{x}_{g}\right]  \left(  a_{R}K_{n,R}^{\mathcal{Z}}\right)  \neq0$.
Hence, $a_{R}\neq0$ and $\left[  \mathbf{x}_{g}\right]  \left(  K_{n,R}%
^{\mathcal{Z}}\right)  \neq0$. Thus, we don't have $\left[  \mathbf{x}%
_{g}\right]  \left(  K_{n,R}^{\mathcal{Z}}\right)  =0$ (since $\left[
\mathbf{x}_{g}\right]  \left(  K_{n,R}^{\mathcal{Z}}\right)  \neq0$).

Every element of $\mathbf{L}_{n}$ is a lacunar subset of $\left[  n\right]
$\ \ \ \ \footnote{Indeed, this is obvious when $n$ is positive (since
$\mathbf{L}_{n}$ is defined to be the set of all nonempty lacunar subsets of
$\left[  n\right]  $ in this case), but also obvious when $n=0$ (since
$\mathbf{L}_{n}=\left\{  \varnothing\right\}  $ in this case, but the set
$\varnothing$ is a lacunar subset of $\left[  n\right]  $).}. Hence, $R$ is a
lacunar subset of $\left[  n\right]  $ (since $R\in\mathbf{L}_{n}$).

But Lemma \ref{lem.KnL.coeff} \textbf{(b)} (applied to $h=g$) yields%
\[
\left[  \mathbf{x}_{g}\right]  \left(  K_{n,R}^{\mathcal{Z}}\right)  =%
\begin{cases}
2^{\left\vert g\left(  \left[  n\right]  \right)  \cap\left\{  1,2,3,\ldots
\right\}  \right\vert }, & \text{if }R\subseteq\operatorname*{FE}\left(
g\right)  ;\\
0, & \text{otherwise}%
\end{cases}
.
\]
Hence, $\left[  \mathbf{x}_{g}\right]  \left(  K_{n,R}^{\mathcal{Z}}\right)
=0$ if $R\not \subseteq \operatorname*{FE}\left(  g\right)  $. Thus, we cannot
have $R\not \subseteq \operatorname*{FE}\left(  g\right)  $ (since we don't
have $\left[  \mathbf{x}_{g}\right]  \left(  K_{n,R}^{\mathcal{Z}}\right)
=0$). Therefore, we have
\[
R\subseteq\operatorname*{FE}\left(  g\right)  =\left(  \Lambda\cup\left(
\Lambda+1\right)  \right)  \cap\left[  n\right]  \subseteq\Lambda\cup\left(
\Lambda+1\right)  .
\]
Thus, Proposition \ref{prop.lac.RL} yields that $R\geq\Lambda$ (with respect
to the total order on $\mathbf{P}$). Combining this with $R\neq\Lambda$, we
obtain $R>\Lambda$. Hence, (\ref{pf.prop.KnL.lindep.larger=0}) yields
$a_{R}=0$. This contradicts $a_{R}\neq0$. This contradiction shows that our
assumption was wrong. Hence, (\ref{pf.prop.KnL.lindep.coeff=0}) is proven.]
\end{verlong}

On the other hand, Lemma \ref{lem.KnL.coeff} \textbf{(b)} (applied to $h=g$
and $R=\Lambda$) yields%
\[
\left[  \mathbf{x}_{g}\right]  \left(  K_{n,\Lambda}^{\mathcal{Z}}\right)  =%
\begin{cases}
2^{\left\vert g\left(  \left[  n\right]  \right)  \cap\left\{  1,2,3,\ldots
\right\}  \right\vert }, & \text{if }\Lambda\subseteq\operatorname*{FE}\left(
g\right)  ;\\
0, & \text{otherwise}%
\end{cases}
=2^{\left\vert g\left(  \left[  n\right]  \right)  \cap\left\{  1,2,3,\ldots
\right\}  \right\vert }%
\]
(since $\Lambda\subseteq\operatorname*{FE}\left(  g\right)  $).

\begin{vershort}
Now, recall that $\sum_{R\in\mathbf{L}_{n}}a_{R}K_{n,R}^{\mathcal{Z}}=0$.
Hence, $\left[  \mathbf{x}_{g}\right]  \left(  \sum_{R\in\mathbf{L}_{n}}%
a_{R}K_{n,R}^{\mathcal{Z}}\right)  =\left[  \mathbf{x}_{g}\right]  \left(
0\right)  =0$. Therefore,%
\begin{align*}
0 &  =\left[  \mathbf{x}_{g}\right]  \left(  \sum_{R\in\mathbf{L}_{n}}%
a_{R}K_{n,R}^{\mathcal{Z}}\right)  =\sum_{R\in\mathbf{L}_{n}}\left[
\mathbf{x}_{g}\right]  \left(  a_{R}K_{n,R}^{\mathcal{Z}}\right)  \\
&  =\left[  \mathbf{x}_{g}\right]  \left(  a_{\Lambda}K_{n,\Lambda
}^{\mathcal{Z}}\right)  +\sum_{\substack{R\in\mathbf{L}_{n};\\R\neq\Lambda
}}\underbrace{\left[  \mathbf{x}_{g}\right]  \left(  a_{R}K_{n,R}%
^{\mathcal{Z}}\right)  }_{\substack{=0\\\text{(by
(\ref{pf.prop.KnL.lindep.coeff=0}))}}}\ \ \ \ \ \ \ \ \ \ \left(  \text{since
}\Lambda\in\mathbf{L}_{n}\right)  \\
&  =\left[  \mathbf{x}_{g}\right]  \left(  a_{\Lambda}K_{n,\Lambda
}^{\mathcal{Z}}\right)  =a_{\Lambda}\underbrace{\left[  \mathbf{x}_{g}\right]
\left(  K_{n,\Lambda}^{\mathcal{Z}}\right)  }_{=2^{\left\vert g\left(  \left[
n\right]  \right)  \cap\left\{  1,2,3,\ldots\right\}  \right\vert }%
}=\underbrace{a_{\Lambda}}_{\neq0}\underbrace{2^{\left\vert g\left(  \left[
n\right]  \right)  \cap\left\{  1,2,3,\ldots\right\}  \right\vert }}_{\neq
0}\neq0.
\end{align*}
This contradiction shows that our assumption was false. Hence, $\left(
a_{R}\right)  _{R\in\mathbf{L}_{n}}=\left(  0\right)  _{R\in\mathbf{L}_{n}}$
is proven.
\end{vershort}

\begin{verlong}
Hence,%
\[
\left[  \mathbf{x}_{g}\right]  \left(  a_{\Lambda}K_{n,\Lambda}^{\mathcal{Z}%
}\right)  =a_{\Lambda}\underbrace{\left[  \mathbf{x}_{g}\right]  \left(
K_{n,\Lambda}^{\mathcal{Z}}\right)  }_{=2^{\left\vert g\left(  \left[
n\right]  \right)  \cap\left\{  1,2,3,\ldots\right\}  \right\vert }%
}=\underbrace{a_{\Lambda}}_{\neq0}\underbrace{2^{\left\vert g\left(  \left[
n\right]  \right)  \cap\left\{  1,2,3,\ldots\right\}  \right\vert }}_{\neq
0}\neq0.
\]

Now, recall that $\sum_{R\in\mathbf{L}_{n}}a_{R}K_{n,R}^{\mathcal{Z}}=0$.
Hence, $\left[  \mathbf{x}_{g}\right]  \left(  \sum_{R\in\mathbf{L}_{n}}%
a_{R}K_{n,R}^{\mathcal{Z}}\right)  =\left[  \mathbf{x}_{g}\right]  \left(
0\right)  =0$. Therefore,%
\begin{align*}
0 &  =\left[  \mathbf{x}_{g}\right]  \left(  \sum_{R\in\mathbf{L}_{n}}%
a_{R}K_{n,R}^{\mathcal{Z}}\right)  =\sum_{R\in\mathbf{L}_{n}}\left[
\mathbf{x}_{g}\right]  \left(  a_{R}K_{n,R}^{\mathcal{Z}}\right)  \\
&  =\left[  \mathbf{x}_{g}\right]  \left(  a_{\Lambda}K_{n,\Lambda
}^{\mathcal{Z}}\right)  +\sum_{\substack{R\in\mathbf{L}_{n};\\R\neq\Lambda
}}\underbrace{\left[  \mathbf{x}_{g}\right]  \left(  a_{R}K_{n,R}%
^{\mathcal{Z}}\right)  }_{\substack{=0\\\text{(by
(\ref{pf.prop.KnL.lindep.coeff=0}))}}}\\
&  \ \ \ \ \ \ \ \ \ \ \left(
\begin{array}
[c]{c}%
\text{here, we have split off the addend}\\
\text{for }R=\Lambda\text{ from the sum (since }\Lambda\in\mathbf{L}%
_{n}\text{)}%
\end{array}
\right)  \\
&  =\left[  \mathbf{x}_{g}\right]  \left(  a_{\Lambda}K_{n,\Lambda
}^{\mathcal{Z}}\right)  +\underbrace{\sum_{\substack{R\in\mathbf{L}%
_{n};\\R\neq\Lambda}}0}_{=0}=\left[  \mathbf{x}_{g}\right]  \left(
a_{\Lambda}K_{n,\Lambda}^{\mathcal{Z}}\right)  \neq0.
\end{align*}
This is absurd. This contradiction shows that our assumption was false. Hence,
$\left(  a_{R}\right)  _{R\in\mathbf{L}_{n}}=\left(  0\right)  _{R\in
\mathbf{L}_{n}}$ is proven.
\end{verlong}

Now, forget that we fixed $\left(  a_{R}\right)  _{R\in\mathbf{L}_{n}}$. We
thus have shown that if $\left(  a_{R}\right)  _{R\in\mathbf{L}_{n}}%
\in\mathbb{Q}^{\mathbf{L}_{n}}$ is a family of scalars (in $\mathbb{Q}$) such
that $\sum_{R\in\mathbf{L}_{n}}a_{R}K_{n,R}^{\mathcal{Z}}=0$, then $\left(
a_{R}\right)  _{R\in\mathbf{L}_{n}}=\left(  0\right)  _{R\in\mathbf{L}_{n}}$.
In other words, the family $\left(  K_{n,R}^{\mathcal{Z}}\right)
_{R\in\mathbf{L}_{n}}$ is $\mathbb{Q}$-linearly independent. In other words,
the family $\left(  K_{n,\Lambda}^{\mathcal{Z}}\right)  _{\Lambda\in
\mathbf{L}_{n}}$ is $\mathbb{Q}$-linearly independent. This proves Proposition
\ref{prop.KnL.lindep}.
\end{proof}

\begin{vershort}
A second proof of Proposition \ref{prop.KnL.lindep} can be found in
\cite{verlong}.
\end{vershort}

\begin{verlong}
\begin{proof}
[Second proof of Proposition \ref{prop.KnL.lindep} (sketched).]
Let $\Omega$ be the subset
\[
\left\{  1,3,5,\ldots\right\}  \cap\left[  n\right]  =\left\{  i\in\left[
n\right]  \ \mid\ i\text{ is odd}\right\}
\]
of $\left[  n\right]  $. This is
clearly a lacunar subset of $\left[  n\right]  $.

We are going to prove the following claim:

\begin{statement}
\textit{Claim 1:} \textbf{(a)} If $n$ is odd, then the only syzygy\footnote{If
$\left(  v_{h}\right)  _{h\in H}$ is a family of vectors in a vector space
over a field $\mathbb{F}$, then a \textit{syzygy} of this family $\left(
v_{h}\right)  _{h\in H}$ means a family $\left(  \lambda_{h}\right)  _{h\in
H}\in\mathbb{F}^{H}$ of scalars in $\mathbb{F}$ satisfying $\sum_{h\in
H}\lambda_{h}v_{h}=0$.
\par
Thus, a syzygy is what is commonly called a \textquotedblleft linear
dependence relation\textquotedblright\ (at least when the scalars $\lambda
_{h}$ are not all $0$). By abuse of notation, we shall speak of the
\textquotedblleft syzygy $\sum_{h\in H}\lambda_{h}v_{h}=0$\textquotedblright%
\ meaning not the equality $\sum_{h\in H}\lambda_{h}v_{h}=0$ but the family of
coefficients $\left(  \lambda_{h}\right)  _{h\in H}$.
\par
When we say \textquotedblleft the only syzygy\textquotedblright, we mean
\textquotedblleft the only nonzero syzygy up to scalar
multiples\textquotedblright.} of the family
$\left(  K_{n,\Lambda }^{\mathcal{Z}}\right) %
_{\Lambda\subseteq\left[  n\right]  \text{ is lacunar}}$
is $\sum_{\Lambda\subseteq\Omega}\left(  -1\right)  ^{\left\vert
\Lambda\right\vert }K_{n,\Lambda}^{\mathcal{Z}}=0$.

\textbf{(b)} If $n$ is even, then the family $\left(  K_{n,\Lambda
}^{\mathcal{Z}}\right)  _{\Lambda\subseteq\left[  n\right]  \text{ is
lacunar}}$ is $\mathbb{Q}$-linearly independent.
\end{statement}

Claim 1 (once proven) will clearly yield
Proposition \ref{prop.KnL.lindep}. Indeed, the family
$\left(  K_{n,\Lambda}^{\mathcal{Z}}\right)  _{\Lambda\in \mathbf{L}_n}$
is a subfamily of the family
$\left(  K_{n,\Lambda }^{\mathcal{Z}}\right) %
_{\Lambda\subseteq\left[  n\right]  \text{ is lacunar}}$,
which is obtained from the latter family by removing the element
$K_{n,\varnothing}^{\mathcal{Z}}$ when $n > 0$ (by the definition
of $\mathbf{L}_n$).
But Claim 1 shows that all nontrivial syzygies of the latter family
(if there are any to begin with) involve the element
$K_{n,\varnothing}^{\mathcal{Z}}$, and thus disappear when this
element is removed.
Hence, it lets us conclude that the former family is
$\mathbb{Q}$-linearly independent.
Thus, it remains to prove Claim 1.

For each subset $\Lambda$ of $\left[  n\right]  $, define a power series
$L_{n,\Lambda}^{\mathcal{Z}}\in\operatorname*{Pow}\mathcal{N}$ by
(\ref{eq.lem.KnL.iex.defL}). Then, for each lacunar subset $\Lambda$ of
$\left[  n\right]  $, we have%
\begin{align*}
K_{n,\Lambda}^{\mathcal{Z}}  &  =\sum_{Q\subseteq\Lambda}\left(  -1\right)
^{\left\vert Q\right\vert }L_{n,Q}^{\mathcal{Z}}\ \ \ \ \ \ \ \ \ \ \left(
\text{by Lemma \ref{lem.KnL.iex} \textbf{(a)}}\right)  \text{ and}\\
L_{n,\Lambda}^{\mathcal{Z}}  &  =\sum_{Q\subseteq\Lambda}\left(  -1\right)
^{\left\vert Q\right\vert }K_{n,Q}^{\mathcal{Z}}\ \ \ \ \ \ \ \ \ \ \left(
\text{by Lemma \ref{lem.KnL.iex} \textbf{(b)}}\right)  .
\end{align*}
Hence, the two families $\left(  K_{n,\Lambda}^{\mathcal{Z}}\right)
_{\Lambda\subseteq\left[  n\right]  \text{ is lacunar}}$ and $\left(
L_{n,\Lambda}^{\mathcal{Z}}\right)  _{\Lambda\subseteq\left[  n\right]  \text{
is lacunar}}$ can be obtained from each other by a unitriangular transition
matrix (unitriangular with respect to inclusion\footnote{We are using the fact
that a subset of a lacunar subset is lacunar.}). Thus, the syzygies of these
two families are in bijection with each other. Hence, in order to prove Claim
1, it suffices to prove the following claim:

\begin{statement}
\textit{Claim 2:} \textbf{(a)} If $n$ is odd, then the only syzygy of the
family $\left(  L_{n,\Lambda}^{\mathcal{Z}}\right)  _{\Lambda\subseteq\left[
n\right]  \text{ is lacunar}}$ is $L_{n,\Omega}^{\mathcal{Z}}=0$.

\textbf{(b)} If $n$ is even, then the family $\left(  L_{n,\Lambda
}^{\mathcal{Z}}\right)  _{\Lambda\subseteq\left[  n\right]  \text{ is
lacunar}}$ is $\mathbb{Q}$-linearly independent.
\end{statement}

Let $G$ be the set of all weakly increasing maps $g:\left[  n\right]
\rightarrow\mathcal{N}$. Let $R$ be the free $\mathbb{Q}$-vector space with
basis $G$; its standard basis will be denoted by $\left(  \left[  g\right]
\right)  _{g\in G}$. We define a $\mathbb{Q}$-linear map%
\begin{align*}
\Phi:R  &  \rightarrow\operatorname*{Pow}\mathcal{N},\\
\left[  g\right]   &  \mapsto2^{\left\vert g\left(  \left[  n\right]  \right)
\cap\left\{  1,2,3,\ldots\right\}  \right\vert }\mathbf{x}_{g}.
\end{align*}
This map $\Phi$ is easily seen to be injective (since the maps $g\in G$ are
weakly increasing, and thus can be uniquely recovered from the monomials
$\mathbf{x}_{g}$).

For each subset $\Lambda$ of $\left[  n\right]  $, we define an element
$\widetilde{L}_{\Lambda}$ of $R$ by%
\[
\widetilde{L}_{\Lambda}=\sum_{\substack{g\in G;\\\Lambda\cap\operatorname*{FE}%
\left(  g\right)  =\varnothing}}\left[  g\right]  .
\]
Then, each subset $\Lambda$ of $\left[  n\right]  $ satisfies%
\begin{align*}
L_{n,\Lambda}^{\mathcal{Z}}  &  =\underbrace{\sum_{\substack{g:\left[
n\right]  \rightarrow\mathcal{N}\text{ is}\\\text{weakly increasing;}%
\\\Lambda\cap\operatorname*{FE}\left(  g\right)  =\varnothing}}}%
_{\substack{=\sum_{\substack{g\in G;\\\Lambda\cap\operatorname*{FE}\left(
g\right)  =\varnothing}}\\\text{(by the definition of }G\text{)}%
}}\underbrace{2^{\left\vert g\left(  \left[  n\right]  \right)  \cap\left\{
1,2,3,\ldots\right\}  \right\vert }\mathbf{x}_{g}}_{\substack{=\Phi\left(
\left[  g\right]  \right)  \\\text{(by the definition of }\Phi\text{)}}%
}=\sum_{\substack{g\in G;\\\Lambda\cap\operatorname*{FE}\left(  g\right)
=\varnothing}}\Phi\left(  \left[  g\right]  \right) \\
&  =\Phi\left(  \underbrace{\sum_{\substack{g\in G;\\\Lambda\cap
\operatorname*{FE}\left(  g\right)  =\varnothing}}\left[  g\right]
}_{=\widetilde{L}_{\Lambda}}\right)  =\Phi\left(  \widetilde{L}_{\Lambda
}\right)  .
\end{align*}
Hence, the family $\left(  L_{n,\Lambda}^{\mathcal{Z}}\right)  _{\Lambda
\subseteq\left[  n\right]  \text{ is lacunar}}$ is the image of the family
$\left(  \widetilde{L}_{\Lambda}\right)  _{\Lambda\subseteq\left[  n\right]
\text{ is lacunar}}$ under the map $\Phi$. Thus, the syzygies of the two
families are in bijection (since the map $\Phi$ is injective\footnote{We are
here using the following obvious fact:
\par
Let $V$ and $W$ be two vector spaces over a field $\mathbb{F}$. Let $\left(
v_{h}\right)  _{h\in H}\in V^{H}$ be a family of vectors in $V$. Let
$\phi:V\rightarrow W$ be an injective $\mathbb{F}$-linear map. Then, the
syzygies of the families $\left(  v_{h}\right)  _{h\in H}\in V^{H}$ and
$\left(  \phi\left(  v_{h}\right)  \right)  _{h\in H}\in W^{H}$ are in bijection.
(Actually, these syzygies, when regarded as families of scalars, are literally
the same.) In particular, the family $\left(  v_{h}\right)  _{h\in H}$ is
$\mathbb{F}$-linearly independent if and only if the family $\left(  \phi\left(
v_{h}\right)  \right)  _{h\in H}$ is $\mathbb{F}$-linearly independent.}).
Hence, in order to prove Claim 2, it suffices to prove the following claim:

\begin{statement}
\textit{Claim 3:} \textbf{(a)} If $n$ is odd, then the only syzygy of the
family $\left(  \widetilde{L}_{\Lambda}\right)  _{\Lambda\subseteq\left[
n\right]  \text{ is lacunar}}$ is $\widetilde{L}_{\Omega}=0$.

\textbf{(b)} If $n$ is even, then the family $\left(  \widetilde{L}_{\Lambda
}\right)  _{\Lambda\subseteq\left[  n\right]  \text{ is lacunar}}$ is
$\mathbb{Q}$-linearly independent.
\end{statement}

Let us agree that if $g\in G$, then we will set $g\left(  0\right)  =0$ and
$g\left(  n+1\right)  =\infty$. Hence, $g\left(  i\right)  $ will be a
well-defined element of $\mathcal{N}$ for each $i\in\left\{  0,1,\ldots
,n+1\right\}  $.

If $g\in G$, then we let $\operatorname*{Stag}\left(  g\right)  $ be the
subset $\left\{  i\in\left[  n+1\right]  \ \mid\ g\left(  i\right)  =g\left(
i-1\right)  \right\}  $ of $\left[  n+1\right]  $. It is easy to see that the
family $\left(  \sum_{\substack{g\in G;\\\operatorname*{Stag}\left(  g\right)
=T}}\left[  g\right]  \right)  _{T\subseteq\left[  n+1\right]  ;\ T\neq\left[
n+1\right]  }$ of elements of $R$ is $\mathbb{Q}$-linearly
independent\footnote{\textit{Proof.} Clearly, any two elements of this family
are supported on different basis elements (i.e., any $\left[  g\right]  $
appearing in one of them cannot appear in any other). It thus remains to show
that these elements are $\neq0$. In other words, it remains to show that for
any proper subset $T$ of $\left[  n+1\right]  $, we have $\sum_{\substack{g\in
G;\\\operatorname*{Stag}\left(  g\right)  =T}}\left[  g\right]  \neq0$. But
this is easy: Just construct some $g\in G$ satisfying $\operatorname*{Stag}%
\left(  g\right)  =T$.}. Hence, the family $\left(  \sum_{\substack{g\in
G;\\\operatorname*{Stag}\left(  g\right)  \supseteq T}}\left[  g\right]
\right)  _{T\subseteq\left[  n+1\right]  ;\ T\neq\left[  n+1\right]  }$ of
elements of $R$ is $\mathbb{Q}$-linearly independent, too (because this family
is obtained from the previous family $\left(  \sum_{\substack{g\in
G;\\\operatorname*{Stag}\left(  g\right)  =T}}\left[  g\right]  \right)
_{T\subseteq\left[  n+1\right]  ;\ T\neq\left[  n+1\right]  }$ via a
unitriangular change-of-basis matrix\footnote{unitriangular with respect to
the reverse inclusion order (notice that $\sum_{\substack{g\in
G;\\\operatorname*{Stag}\left(  g\right)  =T}}\left[  g\right]  =0$ for
$T=\left[  n+1\right]  $, so the exclusion of $\left[  n+1\right]  $ makes
sense and does not mess up our computations)}). Therefore, the only syzygy of
the family \newline$\left(  \sum_{\substack{g\in G;\\\operatorname*{Stag}%
\left(  g\right)  \supseteq T}}\left[  g\right]  \right)  _{T\subseteq\left[
n+1\right]  }$ is $\sum_{\substack{g\in G;\\\operatorname*{Stag}\left(
g\right)  \supseteq\left[  n+1\right]  }}\left[  g\right]  =0$ (since it is
easy to see that no $g\in G$ satisfies $\operatorname*{Stag}\left(  g\right)
\supseteq\left[  n+1\right]  $, which is why $\sum_{\substack{g\in
G;\\\operatorname*{Stag}\left(  g\right)  \supseteq\left[  n+1\right]
}}\left[  g\right]  $ is indeed $0$).

But if $\Lambda$ is a lacunar subset of $\left[  n\right]  $, and if $g\in G$,
then we have the following logical equivalence:%
\begin{align*}
&  \ \left(  \Lambda\cap\operatorname*{FE}\left(  g\right)  =\varnothing
\right) \\
&  \Longleftrightarrow\ \left(  \text{no }i\in\Lambda\text{ satisfies }%
i\in\operatorname*{FE}\left(  g\right)  \right) \\
&  \Longleftrightarrow\ \left(  \text{each }i\in\Lambda\text{ satisfies
}i\notin\operatorname*{FE}\left(  g\right)  \right) \\
&  \Longleftrightarrow\ \left(  \text{each }i\in\Lambda\text{ satisfies
}\underbrace{i\neq\min\left(  g^{-1}\left(  h\right)  \right)  \text{ for all
}h\in\left\{  1,2,3,\ldots,\infty\right\}  }_{\substack{\Longleftrightarrow
\ \left(  g\left(  i\right)  =g\left(  i-1\right)  \right)  \\\text{(since
}g\text{ is weakly increasing, and }g\left(  0\right)  =0\text{)}}}\right. \\
&  \ \ \ \ \ \ \ \ \ \ \ \ \ \ \ \ \ \ \ \ \left.  \text{and }%
\underbrace{i\neq\max\left(  g^{-1}\left(  h\right)  \right)  \ \text{for
all}\ h\in\left\{  0,1,2,3,\ldots\right\}  }_{\substack{\Longleftrightarrow
\ \left(  g\left(  i\right)  =g\left(  i+1\right)  \right)  \\\text{(since
}g\text{ is weakly increasing, and }g\left(  n+1\right)  =\infty\text{)}%
}}\right) \\
&  \Longleftrightarrow\ \left(  \text{each }i\in\Lambda\text{ satisfies
}\underbrace{g\left(  i\right)  =g\left(  i-1\right)  }_{\Longleftrightarrow
\ \left(  i\in\operatorname*{Stag}\left(  g\right)  \right)  }\text{ and
}\underbrace{g\left(  i\right)  =g\left(  i+1\right)  }_{\Longleftrightarrow
\ \left(  i+1\in\operatorname*{Stag}\left(  g\right)  \right)  }\right) \\
&  \Longleftrightarrow\ \left(  \text{each }i\in\Lambda\text{ satisfies }%
i\in\operatorname*{Stag}\left(  g\right)  \text{ and }i+1\in
\operatorname*{Stag}\left(  g\right)  \right) \\
&  \Longleftrightarrow\ \left(  \Lambda\cup\left(  \Lambda+1\right)
\subseteq\operatorname*{Stag}\left(  g\right)  \right)  \ \Longleftrightarrow
\ \left(  \operatorname*{Stag}\left(  g\right)  \supseteq\Lambda\cup\left(
\Lambda+1\right)  \right)  .
\end{align*}
Hence, if $\Lambda$ is a lacunar subset of $\left[  n\right]  $, then the
condition $\Lambda\cap\operatorname*{FE}\left(  g\right)  =\varnothing$ is
equivalent to the condition $\operatorname*{Stag}\left(  g\right)
\supseteq\Lambda\cup\left(  \Lambda+1\right)  $. Thus, for each lacunar subset
$\Lambda$ of $\left[  n\right]  $, the definition of $\widetilde{L}_{\Lambda}$
becomes%
\[
\widetilde{L}_{\Lambda}=\sum_{\substack{g\in G;\\\Lambda\cap\operatorname*{FE}%
\left(  g\right)  =\varnothing}}\left[  g\right]  =\sum_{\substack{g\in
G;\\\operatorname*{Stag}\left(  g\right)  \supseteq\Lambda\cup\left(
\Lambda+1\right)  }}\left[  g\right]  .
\]
Hence, the family $\left(  \widetilde{L}_{\Lambda}\right)  _{\Lambda
\subseteq\left[  n\right]  \text{ is lacunar}}$ is a subfamily of the family
$\left(  \sum_{\substack{g\in G;\\\operatorname*{Stag}\left(  g\right)
\supseteq T}}\left[  g\right]  \right)  _{T\subseteq\left[  n+1\right]  }$
(because if $\Lambda$ is a lacunar subset of $\left[  n\right]  $, then
$\Lambda\cup\left(  \Lambda+1\right)  $ is a well-defined subset of $\left[
n+1\right]  $, and moreover $\Lambda$ can be uniquely recovered from
$\Lambda\cup\left(  \Lambda+1\right)  $\ \ \ \ \footnote{This takes a bit of
thought to check. You need to show that if $\Lambda_{1}$ and $\Lambda_{2}$ are
two lacunar subsets of $\left[  n\right]  $ satisfying $\Lambda_{1}\cup\left(
\Lambda_{1}+1\right)  =\Lambda_{2}\cup\left(  \Lambda_{2}+1\right)  $, then
$\Lambda_{1}=\Lambda_{2}$. This follows from
Corollary~\ref{cor.lac.L=L}.}). The further
argument depends on the parity of $n$:

\begin{itemize}
\item If $n$ is odd, then the vanishing element $\sum_{\substack{g\in
G;\\\operatorname*{Stag}\left(  g\right)  \supseteq\left[  n+1\right]
}}\left[  g\right]  $ does appear in the family $\left(  \widetilde{L}%
_{\Lambda}\right)  _{\Lambda\subseteq\left[  n\right]  \text{ is lacunar}}$,
because there exists a lacunar subset $\Lambda$ of $\left[  n\right]  $
satisfying $\Lambda\cup\left(  \Lambda+1\right)  =\left[  n+1\right]  $:
Namely, this $\Lambda$ is $\Omega$. Thus, the only syzygy of the family
$\left(  \widetilde{L}_{\Lambda}\right)  _{\Lambda\subseteq\left[  n\right]
\text{ is lacunar}}$ is $\widetilde{L}_{\Omega}=0$ (since the only syzygy of
the family $\left(  \sum_{\substack{g\in G;\\\operatorname*{Stag}\left(
g\right)  \supseteq T}}\left[  g\right]  \right)  _{T\subseteq\left[
n+1\right]  }$ is $\sum_{\substack{g\in G;\\\operatorname*{Stag}\left(
g\right)  \supseteq\left[  n+1\right]  }}\left[  g\right]  =0$).

\item If $n$ is even, then the vanishing element $\sum_{\substack{g\in
G;\\\operatorname*{Stag}\left(  g\right)  \supseteq\left[  n+1\right]
}}\left[  g\right]  $ does not appear in the family $\left(  \widetilde{L}%
_{\Lambda}\right)  _{\Lambda\subseteq\left[  n\right]  \text{ is lacunar}}$,
since no lacunar subset $\Lambda$ of $\left[  n\right]  $ satisfies
$\Lambda\cup\left(  \Lambda+1\right)  =\left[  n+1\right]  $ \ \ \ \ %
\footnote{The easiest way to see this is as follows: If $\Lambda$ is a
lacunar subset of $\left[n\right]$, then the sets $\Lambda$ and
$\Lambda+1$ are disjoint (since $\Lambda$ is lacunar), and thus their
union $\Lambda\cup\left(  \Lambda+1\right)$ has size
$\left|\Lambda\right| + \left|\Lambda+1\right| =
2 \left|\Lambda\right|$ (since $\left|\Lambda+1\right| =
\left|\Lambda\right|$), which is even. But the set $\left[n+1\right]$
has size $n+1$, which is odd (since $n$ is even). Thus,
$\Lambda\cup\left(  \Lambda+1\right)$ cannot equal $\left[n+1\right]$
in this case.}. Hence, the
syzygy $\sum_{\substack{g\in G;\\\operatorname*{Stag}\left(  g\right)
\supseteq\left[  n+1\right]  }}\left[  g\right]  =0$ of the family $\left(
\sum_{\substack{g\in G;\\\operatorname*{Stag}\left(  g\right)  \supseteq
T}}\left[  g\right]  \right)  _{T\subseteq\left[  n+1\right]  }$ disappears
when we pass to the subfamily $\left(  \widetilde{L}_{\Lambda}\right)
_{\Lambda\subseteq\left[  n\right]  \text{ is lacunar}}$. Consequently, the
subfamily $\left(  \widetilde{L}_{\Lambda}\right)  _{\Lambda\subseteq\left[
n\right]  \text{ is lacunar}}$ is $\mathbb{Q}$-linearly independent.
\end{itemize}

This proves Claim 3. As explained above, this yields Claim 2, hence also Claim
1, and thus completes the proof of Proposition \ref{prop.KnL.lindep}.
\end{proof}
\end{verlong}

\begin{corollary}
\label{cor.KnL.lindep-all}The family%
\[
\left(  K_{n,\Lambda}^{\mathcal{Z}}\right) %
_{n \in \NN; \ \Lambda \in \mathbf{L}_n}
\]
is $\mathbb{Q}$-linearly independent.
\end{corollary}

\begin{vershort}
\begin{proof}
[Proof of Corollary \ref{cor.KnL.lindep-all}.]Follows from Proposition
\ref{prop.KnL.lindep} using gradedness; see \cite{verlong} for details.
\end{proof}
\end{vershort}

\begin{verlong}
\begin{proof}
[Proof of Corollary \ref{cor.KnL.lindep-all}.]For each $n\in\mathbb{N}$, the
family
$\left(  K_{n,\Lambda}^{\mathcal{Z}}\right)  _{\Lambda\in \mathbf{L}_n}$
is $\mathbb{Q}$-linearly
independent (by Proposition \ref{prop.KnL.lindep}). Furthermore, these
families for varying $n$ live in linearly disjoint subspaces of
$\operatorname*{Pow}\mathcal{N}$ (because for each $n \in \NN$ and $\Lambda
\subseteq\left[  n\right]  $, the power series $K_{n,\Lambda}^{\mathcal{Z}}$
is homogeneous of degree $n$). Thus, the union $\left(  K_{n,\Lambda
}^{\mathcal{Z}}\right)  _{n \in \NN;\ \Lambda \in \mathbf{L}_n}$ of all
these families must also be $\mathbb{Q}%
$-linearly independent.
This proves Corollary \ref{cor.KnL.lindep-all}.
\end{proof}
\end{verlong}

We can now finally prove what we came here for:

\begin{theorem}
\label{thm.Epk.sh-co-a}The permutation statistic $\operatorname{Epk}$ is shuffle-compatible.
\end{theorem}

\begin{proof}
[Proof of Theorem \ref{thm.Epk.sh-co-a}.]We must prove that
$\operatorname{Epk}$ is shuffle-compatible. In other words, we must prove
that for any two disjoint permutations $\pi$ and $\sigma$, the multiset
$\left\{  \operatorname{Epk} \tau   \ \mid\ \tau\in S\left(
\pi,\sigma\right)  \right\}  _{\operatorname*{multi}}$ depends only on
$\operatorname{Epk} \pi$, $\operatorname{Epk} \sigma$,
$\left\vert \pi\right\vert $ and $\left\vert \sigma \right\vert $.
In other words, we must prove that if $\pi$ and $\sigma$ are
two disjoint permutations, and if $\pi^{\prime}$ and $\sigma^{\prime}$ are two
disjoint permutations satisfying $\operatorname{Epk} \pi 
=\operatorname{Epk}\left(  \pi^{\prime}\right)  $,
$\operatorname{Epk} \sigma
=\operatorname{Epk}\left(  \sigma^{\prime}\right)  $,
$\left\vert \pi\right\vert =\left\vert \pi^{\prime}\right\vert $ and
$\left\vert \sigma\right\vert =\left\vert \sigma^{\prime}\right\vert $, then
the multiset $\left\{  \operatorname{Epk} \tau   \ \mid\ \tau\in
S\left(  \pi,\sigma\right)  \right\}  _{\operatorname*{multi}}$ equals the
multiset $\left\{  \operatorname{Epk} \tau   \ \mid\ \tau\in
S\left(  \pi^{\prime},\sigma^{\prime}\right)  \right\}
_{\operatorname*{multi}}$.

So let $\pi$ and $\sigma$ be two disjoint permutations, and let $\pi^{\prime}$
and $\sigma^{\prime}$ be two disjoint permutations satisfying
$\operatorname{Epk} \pi   =\operatorname{Epk}\left(
\pi^{\prime}\right)  $, $\operatorname{Epk} \sigma 
=\operatorname{Epk}\left(  \sigma^{\prime}\right)  $, $\left\vert
\pi\right\vert =\left\vert \pi^{\prime}\right\vert $ and $\left\vert
\sigma\right\vert =\left\vert \sigma^{\prime}\right\vert $.

Define $n\in\mathbb{N}$ by $n=\left\vert \pi\right\vert =\left\vert
\pi^{\prime}\right\vert $ (this is well-defined, since $\left\vert
\pi\right\vert =\left\vert \pi^{\prime}\right\vert $). Likewise, define
$m\in\mathbb{N}$ by $m=\left\vert \sigma\right\vert =\left\vert \sigma
^{\prime}\right\vert $. Thus, $\pi$ is an $n$-permutation, while $\sigma$ is
an $m$-permutation. Hence, each $\tau\in S\left(  \pi,\sigma\right)  $ is an
$\left(  n+m\right)  $-permutation, and therefore satisfies
$\operatorname{Epk}\tau\in\mathbf{L}_{n+m}$ (by Proposition
\ref{prop.Epk-lac}, applied to $n+m$ and $\tau$ instead of $n$ and $\pi$).
Thus, the multiset $\left\{  \operatorname{Epk} \tau 
\ \mid\ \tau\in S\left(  \pi,\sigma\right)  \right\}  _{\operatorname*{multi}%
}$ consists of elements of $\mathbf{L}_{n+m}$. The same holds for the multiset
$\left\{  \operatorname{Epk} \tau   \ \mid\ \tau\in S\left(
\pi^{\prime},\sigma^{\prime}\right)  \right\}  _{\operatorname*{multi}}$ (for
similar reasons).

Corollary \ref{cor.KnEpk-prodrule} yields
\begin{align*}
&  K_{n,\operatorname{Epk} \pi   }^{\mathcal{Z}}\cdot
K_{m,\operatorname{Epk} \sigma   }^{\mathcal{Z}}\\
&  =\sum_{\tau\in S\left(  \pi,\sigma\right)  }K_{n+m,\operatorname{Epk}\tau
}^{\mathcal{Z}}=\sum_{\substack{\Lambda\in\mathbf{L}_{n+m}}}\underbrace{\sum
_{\substack{\tau\in S\left(  \pi,\sigma\right)  ;\\
\operatorname{Epk} \tau=\Lambda}} K_{n+m,\Lambda}^{\mathcal{Z}}
}_{=\left\vert \left\{  \tau\in
S\left(  \pi,\sigma\right)  \ \mid\ \operatorname{Epk}\tau=\Lambda\right\}
\right\vert K_{n+m,\Lambda}^{\mathcal{Z}}}\\
&  \ \ \ \ \ \ \ \ \ \ \left(  \text{because each }\tau\in S\left(  \pi
,\sigma\right)  \text{ satisfies }\operatorname{Epk}\tau\in\mathbf{L}%
_{n+m}\right)  \\
&  =\sum_{\substack{\Lambda\in\mathbf{L}_{n+m}}}\left\vert \left\{  \tau\in
S\left(  \pi,\sigma\right)  \ \mid\ \operatorname{Epk}\tau=\Lambda\right\}
\right\vert K_{n+m,\Lambda}^{\mathcal{Z}}.
\end{align*}
The same argument (but using $\pi^{\prime}$ and $\sigma^{\prime}$ instead
of $\pi$ and $\sigma$) yields
\[
K_{n,\operatorname{Epk}\left(  \pi^{\prime}\right)  }^{\mathcal{Z}}\cdot
K_{m,\operatorname{Epk}\left(  \sigma^{\prime}\right)  }^{\mathcal{Z}}%
=\sum_{\substack{\Lambda\in\mathbf{L}_{n+m}}}\left\vert \left\{  \tau\in
S\left(  \pi^{\prime},\sigma^{\prime}\right)
\ \mid\ \operatorname{Epk} \tau=\Lambda\right\}  \right\vert
K_{n+m,\Lambda}^{\mathcal{Z}}.
\]
The left-hand sides of these two equalities are identical (since
$\operatorname{Epk} \pi   =\operatorname{Epk}\left(
\pi^{\prime}\right)  $ and $\operatorname{Epk} \sigma 
=\operatorname{Epk}\left(  \sigma^{\prime}\right)  $). Thus, their right-hand
sides must also be identical. In other words, we have%
\begin{align*}
&  \sum_{\substack{\Lambda\in\mathbf{L}_{n+m}}}\left\vert \left\{  \tau\in
S\left(  \pi,\sigma\right)  \ \mid\ \operatorname{Epk}\tau=\Lambda\right\}
\right\vert K_{n+m,\Lambda}^{\mathcal{Z}}\\
&  =\sum_{\substack{\Lambda\in\mathbf{L}_{n+m}}}\left\vert \left\{  \tau\in
S\left(  \pi^{\prime},\sigma^{\prime}\right)
\ \mid\ \operatorname{Epk} \tau =\Lambda\right\}
\right\vert K_{n+m,\Lambda}^{\mathcal{Z}}.
\end{align*}
Since the family $\left(  K_{n+m,\Lambda}^{\mathcal{Z}}\right)  _{\Lambda
\in\mathbf{L}_{n+m}}$ is $\mathbb{Q}$-linearly independent (by Proposition
\ref{prop.KnL.lindep}), this shows that%
\[
\left\vert \left\{  \tau\in S\left(  \pi,\sigma\right)  \ \mid
\ \operatorname{Epk}\tau=\Lambda\right\}  \right\vert =\left\vert \left\{
\tau\in S\left(  \pi^{\prime},\sigma^{\prime}\right)  \ \mid
\ \operatorname{Epk}\tau=\Lambda\right\}  \right\vert
\]
for each $\Lambda\in\mathbf{L}_{n+m}$.
\begin{verlong}
In other words, for each $\Lambda\in\mathbf{L}_{n+m}$,
the multiplicity with which $\Lambda$ appears in the multiset
$\left\{
\operatorname{Epk} \tau   \ \mid\ \tau\in S\left(  \pi
,\sigma\right)  \right\}  _{\operatorname*{multi}}$ equals the
multiplicity with which $\Lambda$ appears in the multiset
$\left\{  \operatorname{Epk} \tau   \ \mid\ \tau\in S\left(
\pi^{\prime},\sigma^{\prime}\right)  \right\}  _{\operatorname*{multi}}$.
\end{verlong}
In other words, the multiset $\left\{
\operatorname{Epk} \tau   \ \mid\ \tau\in S\left(  \pi
,\sigma\right)  \right\}  _{\operatorname*{multi}}$ equals the multiset
\newline
$\left\{  \operatorname{Epk} \tau   \ \mid\ \tau\in S\left(
\pi^{\prime},\sigma^{\prime}\right)  \right\}  _{\operatorname*{multi}}$
(because both of these multisets consist of elements of $\mathbf{L}_{n+m}$,
and the previous sentence shows that each of these elements appears with
equal multiplicities in them). This completes our proof of Theorem
\ref{thm.Epk.sh-co-a}.
\end{proof}

We end this section with a tangential remark for readers of
\cite{part1}:

\begin{remark}
\label{rmk.Epk.equivalents}Let us use the notations of \cite{part1}
(specifically, the concept of ``equivalent'' statistics defined
in \cite[Section 3.1]{part1}; and various specific statistics
defined in \cite[Section 2.2]{part1}).
The permutation statistics $\left(
\operatorname*{Lpk},\operatorname{val}\right)  $, $\left(
\operatorname*{Lpk},\operatorname*{udr}\right)  $ and $\left(
\operatorname*{Pk},\operatorname*{udr}\right)  $ are equivalent to
$\operatorname{Epk}$, and therefore are shuffle-compatible.
% (See \cite{part1} for the definitions of all language left undefined here.)
\end{remark}

\begin{proof}
[Proof of Remark \ref{rmk.Epk.equivalents} (sketched).]If
$\operatorname{st}_1$ and $\operatorname{st}_2$ are two permutation
statistics, then we shall write
$\operatorname{st}_1 \sim\operatorname{st}_2$ to mean \textquotedblleft%
$\operatorname{st}_1$ is equivalent to
$\operatorname{st}_2$\textquotedblright.

The permutation statistic $\operatorname{val}$ is equivalent to
$\operatorname*{epk}$, because of \cite[Lemma 2.1 \textbf{(e)}]{part1}. In
other words, $\operatorname{val}\sim\operatorname*{epk}$. Hence, $\left(
\operatorname*{Lpk},\operatorname{val}\right)  \sim\left(
\operatorname*{Lpk},\operatorname*{epk}\right)  $. But if $\pi$ is an
$n$-permutation, then $\operatorname{Epk}\pi$ can be computed from the
knowledge of $\operatorname*{Lpk}\pi$ and $\operatorname*{epk}\pi$ (indeed,
$\operatorname{Epk}\pi$ differs from $\operatorname*{Lpk}\pi$ only in the
possible element $n$, so that
\[
\operatorname{Epk}\pi=%
\begin{cases}
\operatorname*{Lpk}\pi, & \text{if }\operatorname*{epk}\pi=\left\vert
\operatorname*{Lpk}\pi\right\vert ;\\
\operatorname*{Lpk}\pi\cup\left\{  n\right\}  , & \text{if }%
\operatorname*{epk}\pi\neq\left\vert \operatorname*{Lpk}\pi\right\vert
\end{cases}
\]
) and vice versa (since $\operatorname*{Lpk}\pi=\left(  \operatorname{Epk}
\pi\right)  \setminus\left\{  n\right\}  $ and $\operatorname{epk}
\pi=\left\vert \operatorname{Epk}\pi\right\vert $). Thus, $\left(
\operatorname{Lpk},\operatorname{epk}\right)  \sim\operatorname{Epk}$.
Hence, altogether, we obtain $\left(  \operatorname{Lpk},\operatorname{val}
\right)  \sim\left(  \operatorname*{Lpk},\operatorname*{epk}\right)
\sim\operatorname{Epk}$.
\begin{verlong}
In other words, $\left(  \operatorname*{Lpk}%
,\operatorname{val}\right)  $ is equivalent to $\operatorname{Epk}$.
\end{verlong}

Moreover, \cite[Lemma 2.2 \textbf{(a)}]{part1} shows that
for any permutation $\pi$, the
knowledge of $\operatorname*{Lpk}\pi$ allows us to compute
$\operatorname*{udr}\pi$ from $\operatorname{val}\pi$ and vice versa. Hence,
$\left(  \operatorname*{Lpk},\operatorname*{udr}\right)  \sim\left(
\operatorname*{Lpk},\operatorname{val}\right)  \sim\operatorname{Epk}$.
\begin{verlong}
In
other words, $\left(  \operatorname*{Lpk},\operatorname*{udr}\right)  $ is
equivalent to $\operatorname{Epk}$.
\end{verlong}

On the other hand, $\left(  \operatorname*{Pk},\operatorname*{lpk}\right)
\sim\operatorname*{Lpk}$. (This is proven similarly to our proof of $\left(
\operatorname*{Lpk},\operatorname*{epk}\right)  \sim\operatorname{Epk}$.)

Also, $\operatorname*{udr}\sim\left(  \operatorname*{lpk},\operatorname{val}
\right)  $ (indeed, \cite[Lemma 2.2 \textbf{(b)} and \textbf{(c)}]{part1} show
how the value $\left(  \operatorname*{lpk},\operatorname{val}\right)  \left(
\pi\right)  $ can be computed from $\operatorname*{udr}\pi$, whereas
\cite[Lemma 2.2 \textbf{(a)}]{part1} shows the opposite direction). Hence,
$\left(  \operatorname*{Pk},\operatorname*{udr}\right)  \sim\left(
\operatorname*{Pk},\operatorname*{lpk},\operatorname{val}\right)  \sim\left(
\operatorname*{Lpk},\operatorname{val}\right)  $ (since $\left(
\operatorname*{Pk},\operatorname*{lpk}\right)  \sim\operatorname*{Lpk}$).
Therefore, $\left(  \operatorname*{Pk},\operatorname*{udr}\right)  \sim\left(
\operatorname*{Lpk},\operatorname{val}\right)  \sim\operatorname{Epk}$.
\begin{verlong}
In
other words, $\left(  \operatorname*{Pk},\operatorname*{udr}\right)  $ is
equivalent to $\operatorname{Epk}$.
\end{verlong}

We have now shown that the statistics $\left(  \operatorname{Lpk}%
,\operatorname{val}\right)  $, $\left(  \operatorname{Lpk}%
,\operatorname{udr}\right)  $ and $\left(  \operatorname*{Pk}%
,\operatorname{udr}\right)  $ are equivalent to $\operatorname{Epk}$. Thus,
\cite[Theorem 3.2]{part1} shows that they are shuffle-compatible (since
$\operatorname{Epk}$ is shuffle-compatible). This proves Remark
\ref{rmk.Epk.equivalents}.
\end{proof}

\begin{question}
Our concept of a \textquotedblleft$\mathcal{Z}$-enriched $\left(
P,\gamma\right)  $-partition\textquotedblright\ generalizes the concept of an
\textquotedblleft enriched $\left(  P,\gamma\right)  $%
-partition\textquotedblright\ by restricting ourselves to a subset
$\mathcal{Z}$ of $\mathcal{N}\times\left\{  +,-\right\}  $. (This does not
sound like much of a generalization when stated like this, but as we have seen
the behavior of the power series $\Gamma_{\mathcal{Z}}\left(  P,\gamma\right)
$ depends strongly on what $\mathcal{Z}$ is, and is not all anticipated by the
$\mathcal{Z}=\mathcal{N}\times\left\{  +,-\right\}  $ case.) A different
generalization of enriched $\left(  P,\gamma\right)  $-partitions (introduced
by Hsiao and Petersen in \cite{HsiPet10}) are the \textit{colored }$\left(
P,\gamma\right)  $\textit{-partitions}, where the two-element set $\left\{
+,-\right\}  $ is replaced by the set $\left\{  1,\omega,\ldots,\omega
^{m-1}\right\}  $ of all $m$-th roots of unity (where $m$ is a chosen positive
integer, and $\omega$ is a fixed primitive $m$-th root of unity). We can play
various games with this concept. The most natural thing to do seems to be to
consider $m$ arbitrary total orders $<_{0},<_{1},\ldots,<_{m-1}$ on the
codomain $A$ of the labeling $\gamma$ (perhaps with some nice properties such
as all intervals being finite) and an arbitrary subset $\mathcal{Z}$ of
$\mathcal{N}\times\left\{  1,\omega,\ldots,\omega^{m-1}\right\}  $, and define
a $\mathcal{Z}$-enriched colored $\left(  P,\gamma\right)  $-partition to be a
map $f:P\rightarrow\mathcal{Z}$ such that every $x<y$ in $P$ satisfy the
following conditions:

\begin{enumerate}
\item[\textbf{(i)}] We have $f\left(  x\right)  \preccurlyeq f\left(
y\right)  $. (Here, the total order on $\mathcal{N}\times\left\{
1,\omega,\ldots,\omega^{m-1}\right\}  $ is defined by%
\[
\left(  n,\omega^{i}\right)  \prec\left(  n^{\prime},\omega^{i^{\prime}%
}\right)  \text{ if and only if either }n\prec n^{\prime}\text{ or }\left(
n=n^{\prime}\text{ and }i<i^{\prime}\right)
\]
(for $i,i^{\prime}\in\left\{  0,1,\ldots,m-1\right\}  $).)

\item[\textbf{(ii)}] If $f\left(  x\right)  =f\left(  y\right)  =\left(
n,\omega^{i}\right)  $ for some $n\in\mathcal{N}$ and $i\in\left\{
0,1,\ldots,m-1\right\}  $, then $\gamma\left(  x\right)  <_{i}\gamma\left(
y\right)  $.
\end{enumerate}

Is this a useful concept, and can it be used to study permutation statistics?
\end{question}

\begin{question}
Corollary \ref{cor.KnEpk-prodrule} provides a formula for rewriting a product
of the form $K_{n,\Lambda}^{\mathcal{Z}}\cdot K_{m,\Omega}^{\mathcal{Z}}$ as a
$\mathbb{Q}$-linear combination of $K_{n+m,\Xi}^{\mathcal{Z}}$'s when
$\Lambda\in\mathbf{L}_{n}$ and $\Omega\in\mathbf{L}_{m}$ (because any such
$\Lambda$ and $\Omega$ can be written as $\Lambda=\operatorname{Epk}\pi$ and
$\Omega=\operatorname{Epk}\sigma$ for appropriate permutations $\pi$ and
$\sigma$). Thus, in particular, any such product belongs to the $\mathbb{Q}%
$-linear span of the $K_{n+m,\Xi}^{\mathcal{Z}}$'s. Is this still true if
$\Lambda$ and $\Omega$ are arbitrary subsets of $\left[  n\right]  $ and
$\left[  m\right]  $ rather than having to belong to $\mathbf{L}_{n}$ and to
$\mathbf{L}_{m}$ ? Computations with SageMath suggest that the answer is
\textquotedblleft yes\textquotedblright. For example,%
\begin{align*}
K_{2,\left\{  1,2\right\}  }^{\mathcal{Z}}\cdot K_{1,\left\{  1\right\}
}^{\mathcal{Z}}  & =K_{3,\left\{  2\right\}  }^{\mathcal{Z}}+2\cdot
K_{3,\left\{  1,3\right\}  }^{\mathcal{Z}}\ \ \ \ \ \ \ \ \ \ \text{and}\\
K_{2,\varnothing}^{\mathcal{Z}}\cdot K_{1,\left\{  1\right\}  }^{\mathcal{Z}}
& =K_{3,\varnothing}^{\mathcal{Z}}+K_{3,\left\{  2\right\}  }^{\mathcal{Z}%
}+K_{3,\left\{  1,3\right\}  }^{\mathcal{Z}}=K_{3,\left\{  1\right\}
}^{\mathcal{Z}}+K_{3,\left\{  2\right\}  }^{\mathcal{Z}}+K_{3,\left\{
3\right\}  }^{\mathcal{Z}}.
\end{align*}
Note that the $\mathbb{Q}$-linear span of the $K_{n+m,\Xi}^{\mathcal{Z}}$'s
for all $\Xi\subseteq\left[  n+m\right]  $ is (generally) larger than that of
the $K_{n+m,\Xi}^{\mathcal{Z}}$'s with $\Xi\in\mathbf{L}_{n+m}$.
\end{question}

\section{\label{sect.LR}LR-shuffle-compatibility}

In this section, we shall introduce the concept of \textquotedblleft
LR-shuffle-compatibility\textquotedblright\ (short for \textquotedblleft
left-and-right-shuffle-compatibility\textquotedblright), which is stronger
than usual shuffle-compatibility. We shall prove that $\operatorname{Epk}$
still is LR-shuffle-compatible, and study some other statistics that are and
some that are not.

\subsection{Left and right shuffles}

We begin by introducing \textquotedblleft left shuffles\textquotedblright\ and
\textquotedblleft right shuffles\textquotedblright. There is a well-known
notion of left and right shuffles of words (see, e.g., the operations $\prec$
and $\succ$ in \cite[Example 1]{EbMaPa07}). Specialized to permutations, it
can be defined in the following simple way:

\begin{definition}
\label{def.LRshuf}Let $\pi$ and $\sigma$ be two disjoint permutations. Then:

\begin{itemize}
\item A \textit{left shuffle} of $\pi$ and $\sigma$ means a shuffle $\tau$ of
$\pi$ and $\sigma$ such that the first letter of $\tau$ is the first letter of
$\pi$. (This makes sense only when $\pi$ is nonempty. Otherwise, there are no
left shuffles of $\pi$ and $\sigma$.)

\item A \textit{right shuffle} of $\pi$ and $\sigma$ means a shuffle $\tau$ of
$\pi$ and $\sigma$ such that the first letter of $\tau$ is the first letter of
$\sigma$. (This makes sense only when $\sigma$ is nonempty. Otherwise, there
are no right shuffles of $\pi$ and $\sigma$.)

\item We let $S_{\prec}\left(  \pi,\sigma\right)  $ denote the set of all left
shuffles of $\pi$ and $\sigma$.

\item We let $S_{\succ}\left(  \pi,\sigma\right)  $ denote the set of all
right shuffles of $\pi$ and $\sigma$.
\end{itemize}
\end{definition}

For example, the left shuffles of the two disjoint permutations $\left(
3,1\right)  $ and $\left(  2,6\right)  $ are%
\[
\left(  3,1,2,6\right)  ,\ \ \ \ \ \ \ \ \ \ \left(  3,2,1,6\right)
,\ \ \ \ \ \ \ \ \ \ \left(  3,2,6,1\right)  ,
\]
whereas their right shuffles are%
\[
\left(  2,3,1,6\right)  ,\ \ \ \ \ \ \ \ \ \ \left(  2,3,6,1\right)
,\ \ \ \ \ \ \ \ \ \ \left(  2,6,3,1\right)  .
\]
The permutations $\left(  {}\right)  $ and $\left(  1,3\right)  $ have only
one right shuffle, which is $\left(  1,3\right)  $, and they have no left shuffles.

Clearly, if $\pi$ and $\sigma$ are two disjoint permutations such that at
least one of $\pi$ and $\sigma$ is nonempty, then the two sets $S_{\prec
}\left(  \pi,\sigma\right)  $ and $S_{\succ}\left(  \pi,\sigma\right)  $ are
disjoint and their union is $S\left(  \pi,\sigma\right)  $ (because every
shuffle of $\pi$ and $\sigma$ is either a left shuffle or a right shuffle, but
not both).

Left and right shuffles have a recursive structure that makes them amenable to
inductive arguments. To state it, we need one more definition:

\begin{definition}
\label{def.LR.pi1}Let $n\in\mathbb{N}$. Let $\pi$ be an $n$-permutation.

\textbf{(a)} For each $i\in\left\{  1,2,\ldots,n\right\}  $, we let $\pi_{i}$
denote the $i$-th entry of $\pi$. Thus, $\pi=\left(  \pi_{1},\pi_{2}%
,\ldots,\pi_{n}\right)  $.

\textbf{(b)} If $a$ is a positive integer that does not appear in $\pi$, then
$a:\pi$ denotes the $\left(  n+1\right)  $-permutation $\left(  a,\pi_{1}%
,\pi_{2},\ldots,\pi_{n}\right)  $.

\textbf{(c)} If $n>0$, then $\pi_{\sim1}$ denotes the $\left(  n-1\right)
$-permutation $\left(  \pi_{2},\pi_{3},\ldots,\pi_{n}\right)  $.
\end{definition}

\begin{proposition}
\label{prop.LR.rec}Let $\pi$ and $\sigma$ be two disjoint permutations.

\textbf{(a)} We have $S_{\prec}\left(  \pi,\sigma\right)  =S_{\succ}\left(
\sigma,\pi\right)  $.

\textbf{(b)} If $\pi$ is nonempty, then the permutations $\pi_{\sim1}$ and
$\pi_{1}:\sigma$ are well-defined and disjoint, and satisfy $S_{\prec}\left(
\pi,\sigma\right)  =S_{\succ}\left(  \pi_{\sim1},\pi_{1}:\sigma\right)  $.

\textbf{(c)} If $\sigma$ is nonempty, then the permutations $\sigma_{\sim1}$
and $\sigma_{1}:\pi$ are well-defined and disjoint, and satisfy $S_{\succ
}\left(  \pi,\sigma\right)  =S_{\prec}\left(  \sigma_{1}:\pi,\sigma_{\sim
1}\right)  $.
\end{proposition}

\begin{vershort}
\begin{proof}
[Proof of Proposition \ref{prop.LR.rec}.]The fairly simple proof is left to
the reader, who can also find it in \cite{verlong}.
\end{proof}
\end{vershort}

\begin{verlong}
\begin{proof}
[Proof of Proposition \ref{prop.LR.rec}.]\textbf{(a)} The definition of left
shuffles shows that the left shuffles of $\pi$ and $\sigma$ are the shuffles
$\tau$ of $\pi$ and $\sigma$ such that the first letter of $\tau$ is the first
letter of $\pi$. Meanwhile, the definition of right shuffles shows that the
right shuffles of $\sigma$ and $\pi$ are the shuffles $\tau$ of $\sigma$ and
$\pi$ such that the first letter of $\tau$ is the first letter of $\pi$.
Comparing these two descriptions, we conclude that the left shuffles of $\pi$
and $\sigma$ are the same as the right shuffles of $\sigma$ and $\pi$ (since
the shuffles of $\pi$ and $\sigma$ are the same as the shuffles of $\sigma$
and $\pi$). In other words, $S_{\prec}\left(  \pi,\sigma\right)  =S_{\succ
}\left(  \sigma,\pi\right)  $. This proves Proposition \ref{prop.LR.rec}
\textbf{(a)}.

\textbf{(b)} We first make a simple observation:

\begin{statement}
\textit{Claim 1:} Let $\alpha$ and $\beta$ be two permutations such that
$\beta$ is nonempty. Assume that $\alpha$ is a subsequence of $\beta$,
but does not contain the letter $\beta_{1}$. Then, the permutation $\beta
_{1}:\alpha$ also is a subsequence of $\beta$.
\end{statement}

[\textit{Proof of Claim 1:} The letter $\beta_{1}$ does not appear in $\alpha$
(since $\alpha$ does not contain $\beta_{1}$). Thus, $\beta_{1}:\alpha$ is a
well-defined permutation.

We have assumed that $\alpha$ is a subsequence of $\beta$. In other
words, $\alpha=\left(  \beta_{i_{1}},\beta_{i_{2}},\ldots,\beta_{i_{k}%
}\right)  $ for some $k\in\mathbb{N}$ and some positive integers $i_{1}%
,i_{2},\ldots,i_{k}$ satisfying $i_{1}<i_{2}<\cdots<i_{k}$. Consider these
$i_{1},i_{2},\ldots,i_{k}$. From $\alpha=\left(  \beta_{i_{1}},\beta_{i_{2}%
},\ldots,\beta_{i_{k}}\right)  $, we obtain $\beta_{1}:\alpha=\left(
\beta_{1},\beta_{i_{1}},\beta_{i_{2}},\ldots,\beta_{i_{k}}\right)  $.

Let $g\in\left\{  1,2,\ldots,k\right\}  $. Then, $\alpha=\left(  \beta_{i_{1}%
},\beta_{i_{2}},\ldots,\beta_{i_{k}}\right)  $ clearly contains the letter
$\beta_{i_{g}}$. If we had $i_{g}=1$, then this would yield that $\alpha$
contains the letter $\beta_{1}$ (since $i_{g}=1$); but this would contradict
the assumption that $\alpha$ does not contain the letter $\beta_{1}$. Hence,
we cannot have $i_{g}=1$. Thus, $i_{g}>1$, so that $1<i_{g}$.

Now, forget that we fixed $g$. We thus have shown that $1<i_{g}$ for each
$g\in\left\{  1,2,\ldots,k\right\}  $. Combining this with $i_{1}<i_{2}%
<\cdots<i_{k}$, we obtain $1<i_{1}<i_{2}<\cdots<i_{k}$. Hence, $\left(
\beta_{1},\beta_{i_{1}},\beta_{i_{2}},\ldots,\beta_{i_{k}}\right)  $ is a
subsequence of $\beta$. In view of $\beta_{1}:\alpha=\left(  \beta_{1}%
,\beta_{i_{1}},\beta_{i_{2}},\ldots,\beta_{i_{k}}\right)  $, this rewrites as
follows: $\beta_{1}:\alpha$ is a subsequence of $\beta$. This proves Claim 1.]

Assume that $\pi$ is nonempty. The first letter of $\pi$ does not appear in
$\sigma$ (since $\pi$ and $\sigma$ are disjoint). In other words, the letter
$\pi_{1}$ does not appear in $\sigma$. Thus, the permutation $\pi_{1}:\sigma$
is well-defined. The permutation $\pi_{\sim1}$ is clearly well-defined.
Furthermore, the permutations $\pi_{\sim1}$ and $\pi_{1}:\sigma$ are
disjoint\footnote{\textit{Proof.} Let $\ell$ be a letter that appears in both
$\pi_{\sim1}$ and $\pi_{1}:\sigma$. We shall derive a contradiction.
\par
The permutations $\pi$ and $\sigma$ are disjoint; thus, no letter appears in
both $\pi$ and $\sigma$. In other words, a letter that appears in $\pi$ cannot
appear in $\sigma$. But the letter $\ell$ appears in $\pi_{\sim1}$, thus in
$\pi$ (since $\pi_{\sim1}$ is a subsequence of $\pi$). Hence, $\ell$ does not
appear in $\sigma$ (since a letter that appears in $\pi$ cannot appear in
$\sigma$). But $\ell$ appears in $\pi_{1}:\sigma$. But the only letter that
appears in $\pi_{1}:\sigma$ but not in $\sigma$ is the letter $\pi_{1}$ (due
to the construction of $\pi_{1}:\sigma$). Thus, the letter $\ell$ must be
$\pi_{1}$ (since $\ell$ is a letter that appears in $\pi_{1}:\sigma$ but not
in $\sigma$). In other words, $\ell=\pi_{1}$. Hence, the letter $\pi_{1}$
appears in $\pi_{\sim1}$ (since the letter $\ell$ appears in $\pi_{\sim1}$).
\par
But $\pi$ is a permutation, and thus the letters of $\pi$ are distinct. Hence,
the letter $\pi_{1}$ does not appear in $\pi_{\sim1}$. This contradicts the
fact that the letter $\pi_{1}$ appears in $\pi_{\sim1}$.
\par
Now, forget that we fixed $\ell$. We thus have found a contradiction for each
letter $\ell$ that appears in both $\pi_{\sim1}$ and $\pi_{1}:\sigma$. Hence,
no letter appears in both $\pi_{\sim1}$ and $\pi_{1}:\sigma$. In other words,
the permutations $\pi_{\sim1}$ and $\pi_{1}:\sigma$ are disjoint.}. It thus
remains to show that $S_{\prec}\left(  \pi,\sigma\right)  =S_{\succ}\left(
\pi_{\sim1},\pi_{1}:\sigma\right)  $.

Set $m=\left\vert \pi\right\vert $ and $n=\left\vert \sigma\right\vert $.
Thus, $\pi_{\sim1}$ is an $\left(  m-1\right)  $-permutation, whereas $\pi
_{1}:\sigma$ is an $\left(  n+1\right)  $-permutation.

Now, we shall prove that%
\begin{equation}
S_{\prec}\left(  \pi,\sigma\right)  \subseteq S_{\succ}\left(  \pi_{\sim1}%
,\pi_{1}:\sigma\right)  \label{pf.prop.LR.rec.b.1}%
\end{equation}
and%
\begin{equation}
S_{\succ}\left(  \pi_{\sim1},\pi_{1}:\sigma\right)  \subseteq S_{\prec}\left(
\pi,\sigma\right)  . \label{pf.prop.LR.rec.b.2}%
\end{equation}

[\textit{Proof of (\ref{pf.prop.LR.rec.b.1}):} Let $\tau\in S_{\prec}\left(
\pi,\sigma\right)  $. We shall show that $\tau\in S_{\succ}\left(  \pi_{\sim
1},\pi_{1}:\sigma\right)  $.

We have $\tau\in S_{\prec}\left(  \pi,\sigma\right)  $. In other words, $\tau$
is a left shuffle of $\pi$ and $\sigma$ (by the definition of $S_{\prec
}\left(  \pi,\sigma\right)  $). In other words, $\tau$ is a shuffle of $\pi$
and $\sigma$ such that the first letter of $\tau$ is the first letter of $\pi$
(by the definition of a left shuffle).

So the first letter of $\tau$ is the first letter of $\pi$. In other words,
$\tau_{1}=\pi_{1}$. Thus, the permutation $\tau$ is nonempty.

Also, $\tau$ is a shuffle of $\pi$ and $\sigma$. In other words, $\tau$ is an
$\left(  m+n\right)  $-permutation such that both $\pi$ and $\sigma$ are
subsequences of $\tau$ (by the definition of a shuffle).

So $\pi$ is a subsequence of $\tau$. Thus, $\pi_{\sim1}$ is a
subsequence of $\tau$ as well (since $\pi_{\sim1}$ is a subsequence of $\pi$).

Also, $\sigma$ is a subsequence of $\tau$, but does not contain the
letter $\tau_{1}$ (because $\tau_{1}=\pi_{1}$, which does not appear in
$\sigma$). Therefore, Claim 1 (applied to $\alpha=\sigma$ and $\beta=\tau$)
yields that $\tau_{1}:\sigma$ also is a subsequence of $\tau$. In
other words, $\pi_{1}:\sigma$ is a subsequence of $\tau$ (since
$\tau_{1}=\pi_{1}$).

Finally, recall that $\pi_{\sim1}$ is an $\left(  m-1\right)  $-permutation,
whereas $\pi_{1}:\sigma$ is an $\left(  n+1\right)  $-permutation. But $\tau$
is an $\left(  m+n\right)  $-permutation, hence an $\left(  \left(
m-1\right)  +\left(  n+1\right)  \right)  $-permutation (since $m+n=\left(
m-1\right)  +\left(  n+1\right)  $). So we know that $\tau$ is an
\newline$\left(  \left(  m-1\right)  +\left(  n+1\right)  \right)
$-permutation such that both $\pi_{\sim1}$ and $\pi_{1}:\sigma$ are
subsequences of $\tau$. In other words, $\tau$ is a shuffle of $\pi_{\sim1}$
and $\pi_{1}:\sigma$ (by the definition of a shuffle). Since the first letter
of $\tau$ is the first letter of $\pi_{1}:\sigma$ (because the first letter of
$\tau$ is $\tau_{1}=\pi_{1}$, whereas the first letter of $\pi_{1}:\sigma$ is
$\pi_{1}$ as well), we thus conclude that $\tau$ is a right shuffle of
$\pi_{\sim1}$ and $\pi_{1}:\sigma$ (by the definition of a right shuffle). In
other words, $\tau\in S_{\succ}\left(  \pi_{\sim1},\pi_{1}:\sigma\right)  $
(by the definition of $S_{\succ}\left(  \pi_{\sim1},\pi_{1}:\sigma\right)  $).

Since we have proven this for all $\tau\in S_{\prec}\left(  \pi,\sigma\right)
$, we thus conclude that $S_{\prec}\left(  \pi,\sigma\right)  \subseteq
S_{\succ}\left(  \pi_{\sim1},\pi_{1}:\sigma\right)  $. This proves
(\ref{pf.prop.LR.rec.b.1}).]

[\textit{Proof of (\ref{pf.prop.LR.rec.b.2}):} Let $\tau\in S_{\succ}\left(
\pi_{\sim1},\pi_{1}:\sigma\right)  $. We shall show that $\tau\in S_{\prec
}\left(  \pi,\sigma\right)  $.

We have $\tau\in S_{\succ}\left(  \pi_{\sim1},\pi_{1}:\sigma\right)  $. In
other words, $\tau$ is a right shuffle of $\pi_{\sim1}$ and $\pi_{1}:\sigma$
(by the definition of $S_{\succ}\left(  \pi_{\sim1},\pi_{1}:\sigma\right)  $).
In other words, $\tau$ is a shuffle of $\pi_{\sim1}$ and $\pi_{1}:\sigma$ such
that the first letter of $\tau$ is the first letter of $\pi_{1}:\sigma$ (by
the definition of a right shuffle).

So the first letter of $\tau$ is the first letter of $\pi_{1}:\sigma$. In
other words, the first letter of $\tau$ is $\pi_{1}$ (since the first letter
of $\pi_{1}:\sigma$ is $\pi_{1}$). In other words, $\tau_{1}=\pi_{1}$. Note
that the entries of $\pi$ are distinct (since $\pi$ is a permutation); thus,
the letter $\pi_{1}$ does not appear in $\pi_{\sim1}$. In other words, the
letter $\tau_{1}$ does not appear in $\pi_{\sim1}$ (since $\tau_{1}=\pi_{1}$).
Also, the permutation $\tau$ is nonempty (since $\tau_{1}$ exists).

The definitions of $\pi_{1}$, $\pi_{\sim1}$ and $\pi_{1}:\pi_{\sim1}$ yield
$\pi_{1}:\pi_{\sim1}=\pi$. In view of $\tau_{1}=\pi_{1}$, this rewrites as
$\tau_{1}:\pi_{\sim1}=\pi$.

Also, $\tau$ is a shuffle of $\pi_{\sim1}$ and $\pi_{1}:\sigma$. In other
words, $\tau$ is an $\left(  \left(  m-1\right)  +\left(  n+1\right)  \right)
$-permutation such that both $\pi_{\sim1}$ and $\pi_{1}:\sigma$ are
subsequences of $\tau$ (by the definition of a shuffle).

So $\pi_{1}:\sigma$ is a subsequence of $\tau$. Thus, $\sigma$ is
a subsequence of $\tau$ as well (since $\sigma$ is a subsequence of
$\pi_{1}:\sigma$).

Also, $\pi_{\sim1}$ is a subsequence of $\tau$, but does not contain
the letter $\tau_{1}$ (since the letter $\tau_{1}$ does not appear in
$\pi_{\sim1}$). Therefore, Claim 1 (applied to $\alpha=\pi_{\sim1}$ and
$\beta=\tau$) yields that $\tau_{1}:\pi_{\sim1}$ also is a subsequence
of $\tau$. In other words, $\pi$ also is a subsequence of $\tau$
(since $\tau_{1}:\pi_{\sim1}=\pi$).

Finally, recall that $\tau$ is an $\left(  \left(  m-1\right)  +\left(
n+1\right)  \right)  $-permutation, hence an $\left(  m+n\right)
$-permutation (since $\left(  m-1\right)  +\left(  n+1\right)  =m+n$). So we
know that $\tau$ is an $\left(  m+n\right)  $-permutation such that both $\pi$
and $\sigma$ are subsequences of $\tau$. In other words, $\tau$ is a
shuffle of $\pi$ and $\sigma$ (by the definition of a shuffle). Since the
first letter of $\tau$ is the first letter of $\pi$ (because the first letter
of $\tau$ is $\tau_{1}=\pi_{1}$), we thus conclude that $\tau$ is a left
shuffle of $\pi$ and $\sigma$ (by the definition of a left shuffle). In other
words, $\tau\in S_{\prec}\left(  \pi,\sigma\right)  $ (by the definition of
$S_{\prec}\left(  \pi,\sigma\right)  $).

Since we have proven this for all $\tau\in S_{\succ}\left(  \pi_{\sim1}%
,\pi_{1}:\sigma\right)  $, we thus conclude that $S_{\succ}\left(  \pi_{\sim
1},\pi_{1}:\sigma\right)  \subseteq S_{\prec}\left(  \pi,\sigma\right)  $.
This proves (\ref{pf.prop.LR.rec.b.2}).]

Combining (\ref{pf.prop.LR.rec.b.1}) with (\ref{pf.prop.LR.rec.b.2}), we
obtain $S_{\prec}\left(  \pi,\sigma\right)  =S_{\succ}\left(  \pi_{\sim1}%
,\pi_{1}:\sigma\right)  $. This completes the proof of Proposition
\ref{prop.LR.rec} \textbf{(b)}.

\textbf{(c)} Assume that $\sigma$ is nonempty. Proposition \ref{prop.LR.rec}
\textbf{(a)} (applied to $\sigma$ and $\pi$ instead of $\pi$ and $\sigma$)
yields $S_{\prec}\left(  \sigma,\pi\right)  =S_{\succ}\left(  \pi
,\sigma\right)  $.

Proposition \ref{prop.LR.rec} \textbf{(b)} (applied to $\sigma$ and $\pi$
instead of $\pi$ and $\sigma$) yields that the permutations $\sigma_{\sim1}$
and $\sigma_{1}:\pi$ are well-defined and disjoint, and satisfy $S_{\prec
}\left(  \sigma,\pi\right)  =S_{\succ}\left(  \sigma_{\sim1},\sigma_{1}%
:\pi\right)  $. Thus, Proposition \ref{prop.LR.rec} \textbf{(a)} (applied to
$\sigma_{1}:\pi$ and $\sigma_{\sim1}$ instead of $\pi$ and $\sigma$) yields
$S_{\prec}\left(  \sigma_{1}:\pi,\sigma_{\sim1}\right)  =S_{\succ}\left(
\sigma_{\sim1},\sigma_{1}:\pi\right)  $. Combining all the equalities we have
now proven, we obtain%
\[
S_{\succ}\left(  \pi,\sigma\right)  =S_{\prec}\left(  \sigma,\pi\right)
=S_{\succ}\left(  \sigma_{\sim1},\sigma_{1}:\pi\right)  =S_{\prec}\left(
\sigma_{1}:\pi,\sigma_{\sim1}\right)  .
\]
This completes the proof of Proposition \ref{prop.LR.rec} \textbf{(c)}.
\end{proof}
\end{verlong}

\subsection{LR-shuffle-compatibility}

We shall use the so-called \textit{Iverson bracket notation} for truth values:

\begin{definition}
\label{def.iverson}If $\mathcal{A}$ is any logical statement, then we define
an integer $\left[  \mathcal{A}\right]  \in\left\{  0,1\right\}  $ by%
\[
\left[  \mathcal{A}\right]  =%
\begin{cases}
1, & \text{if }\mathcal{A}\text{ is true};\\
0, & \text{if }\mathcal{A}\text{ is false}%
\end{cases}
.
\]
This integer $\left[  \mathcal{A}\right]  $ is known as the \textit{truth
value} of $\mathcal{A}$.
\end{definition}

Thus, for example, $\left[  4>2\right]  =1$ whereas $\left[  2>4\right]  =0$.

We can now define a notion similar to shuffle-compatibility:

\begin{definition}
\label{def.LRcomp}Let $\operatorname{st}$ be a permutation statistic. We say
that $\operatorname{st}$ is \textit{LR-shuffle-compatible} if and only if it
has the following property: For any two disjoint nonempty permutations $\pi$
and $\sigma$, the multisets%
\[
\left\{  \operatorname{st} \tau   \ \mid\ \tau\in S_{\prec
}\left(  \pi,\sigma\right)  \right\}  _{\operatorname*{multi}}
\ \ \ \ \ \ \ \ \ \ \text{and}\ \ \ \ \ \ \ \ \ \ \left\{  \operatorname{st}
\tau \ \mid\ \tau\in S_{\succ}\left(  \pi,\sigma\right)
\right\}  _{\operatorname*{multi}}
\]
depend only on $\operatorname{st} \pi   $, $\operatorname{st} \sigma$,
$\left\vert \pi\right\vert $, $\left\vert
\sigma\right\vert $ and $\left[  \pi_{1}>\sigma_{1}\right]  $.
\end{definition}

In other words, a permutation statistic $\operatorname{st}$ is
LR-shuffle-compatible if and only if every two disjoint nonempty permutations
$\pi$ and $\sigma$ and every two disjoint nonempty permutations $\pi^{\prime}$
and $\sigma^{\prime}$ satisfying%
\begin{align*}
\operatorname{st} \pi    &  =\operatorname{st}\left(
\pi^{\prime}\right)  ,\ \ \ \ \ \ \ \ \ \ \operatorname{st} \sigma
=\operatorname{st}\left(  \sigma^{\prime}\right)  ,\\
\left\vert \pi\right\vert  &  =\left\vert \pi^{\prime}\right\vert
,\ \ \ \ \ \ \ \ \ \ \left\vert \sigma\right\vert =\left\vert \sigma^{\prime
}\right\vert \ \ \ \ \ \ \ \ \ \ \text{and}\ \ \ \ \ \ \ \ \ \ \left[  \pi
_{1}>\sigma_{1}\right]  =\left[  \pi_{1}^{\prime}>\sigma_{1}^{\prime}\right]
\end{align*}
satisfy
\begin{align*}
\left\{  \operatorname{st} \tau   \ \mid\ \tau\in S_{\prec
}\left(  \pi,\sigma\right)  \right\}  _{\operatorname*{multi}}  &  =\left\{
\operatorname{st} \tau   \ \mid\ \tau\in S_{\prec}\left(
\pi^{\prime},\sigma^{\prime}\right)  \right\}  _{\operatorname*{multi}%
}\ \ \ \ \ \ \ \ \ \ \text{and}\\
\left\{  \operatorname{st} \tau   \ \mid\ \tau\in S_{\succ
}\left(  \pi,\sigma\right)  \right\}  _{\operatorname*{multi}}  &  =\left\{
\operatorname{st} \tau   \ \mid\ \tau\in S_{\succ}\left(
\pi^{\prime},\sigma^{\prime}\right)  \right\}  _{\operatorname*{multi}}.
\end{align*}

For example, the permutation statistic $\operatorname*{Pk}$ is not
LR-shuffle-compatible. Indeed, if we take $\pi=\left(  4,2,3\right)  $,
$\sigma=\left(  1\right)  $, $\pi^{\prime}=\left(  2,3,4\right)  $ and
$\sigma^{\prime}=\left(  1\right)  $, then the equalities%
\begin{align*}
\operatorname*{Pk} \pi &  =\operatorname*{Pk}\left(
\pi^{\prime}\right)  ,\ \ \ \ \ \ \ \ \ \ \operatorname*{Pk} \sigma
=\operatorname*{Pk}\left(  \sigma^{\prime}\right)  ,\\
\left\vert \pi\right\vert  &  =\left\vert \pi^{\prime}\right\vert
,\ \ \ \ \ \ \ \ \ \ \left\vert \sigma\right\vert =\left\vert \sigma^{\prime
}\right\vert \ \ \ \ \ \ \ \ \ \ \text{and}\ \ \ \ \ \ \ \ \ \ \left[  \pi
_{1}>\sigma_{1}\right]  =\left[  \pi_{1}^{\prime}>\sigma_{1}^{\prime}\right]
\end{align*}
are all satisfied, but%
\[
\left\{  \operatorname*{Pk} \tau \ \mid\ \tau\in S_{\succ
}\left(  \pi,\sigma\right)  \right\}  _{\operatorname*{multi}}=\left\{
\underbrace{\operatorname*{Pk}\left(  1,4,2,3\right)  }_{=\left\{  2\right\}
}\right\}  _{\operatorname*{multi}}=\left\{  \left\{  2\right\}  \right\}
_{\operatorname*{multi}}%
\]
is not the same as%
\[
\left\{  \operatorname*{Pk} \tau \ \mid\ \tau\in S_{\succ
}\left(  \pi^{\prime},\sigma^{\prime}\right)  \right\}
_{\operatorname*{multi}}=\left\{  \underbrace{\operatorname*{Pk}\left(
1,2,3,4\right)  }_{=\varnothing}\right\}  _{\operatorname*{multi}}=\left\{
\varnothing\right\}  _{\operatorname*{multi}}.
\]
Similarly, the permutation statistic $\operatorname*{Rpk}$ is not
LR-shuffle-compatible. As we will see in Theorem \ref{thm.LRcomp.Pks} further
below, the three statistics $\operatorname{Des}$, $\operatorname*{Lpk}$ and
$\operatorname{Epk}$ are LR-shuffle-compatible.

\subsection{Head-graft-compatibility}

We shall now define another compatibility concept for a permutation statistic,
which will later prove a useful stepping stone for checking the
LR-shuffle-compatibility of this statistic.

\begin{definition}
\label{def.head-comp}Let $\operatorname{st}$ be a permutation statistic. We
say that $\operatorname{st}$ is \textit{head-graft-compatible} if and only if
it has the following property: For any nonempty permutation $\pi$ and any
letter $a$ that does not appear in $\pi$, the element $\operatorname{st}
\left(  a:\pi\right)  $ depends only on $\operatorname{st} \pi 
$, $\left\vert \pi\right\vert $ and $\left[  a>\pi_{1}\right]  $.
\end{definition}

In other words, a permutation statistic $\operatorname{st}$ is
head-graft-compatible if and only if every nonempty permutation $\pi$, every
letter $a$ that does not appear in $\pi$, every nonempty permutation
$\pi^{\prime}$ and every letter $a^{\prime}$ that does not appear in
$\pi^{\prime}$ satisfying%
\[
\operatorname{st} \pi   =\operatorname{st}\left(  \pi^{\prime
}\right)  ,\ \ \ \ \ \ \ \ \ \ \left\vert \pi\right\vert =\left\vert
\pi^{\prime}\right\vert \ \ \ \ \ \ \ \ \ \ \text{and}%
\ \ \ \ \ \ \ \ \ \ \left[  a>\pi_{1}\right]  =\left[  a^{\prime}>\pi
_{1}^{\prime}\right]
\]
satisfy $\operatorname{st}\left(  a:\pi\right)  =\operatorname{st}\left(
a^{\prime}:\pi^{\prime}\right)  $.

For example, the permutation statistic $\operatorname*{Pk}$ is not
head-graft-compatible, because if we take $\pi=\left(  3,1\right)  $, $a=2$,
$\pi^{\prime}=\left(  3,4\right)  $ and $a^{\prime}=2$, then we do have%
\[
\operatorname*{Pk} \pi =\operatorname*{Pk}\left(  \pi^{\prime
}\right)  ,\ \ \ \ \ \ \ \ \ \ \left\vert \pi\right\vert =\left\vert
\pi^{\prime}\right\vert \ \ \ \ \ \ \ \ \ \ \text{and}%
\ \ \ \ \ \ \ \ \ \ \left[  a>\pi_{1}\right]  =\left[  a^{\prime}>\pi
_{1}^{\prime}\right]
\]
but we don't have $\operatorname*{Pk}\left(  a:\pi\right)  =\operatorname*{Pk}%
\left(  a^{\prime}:\pi^{\prime}\right)  $ (in fact, $\operatorname*{Pk}\left(
a:\pi\right)  =\operatorname*{Pk}\left(  2,3,1\right)  =\left\{  2\right\}  $
whereas $\operatorname*{Pk}\left(  a^{\prime}:\pi^{\prime}\right)
=\operatorname*{Pk}\left(  2,3,4\right)  =\varnothing$). Similarly, it can be
shown that $\operatorname*{Rpk}$ is not head-graft-compatible. As we will see
below (in Proposition \ref{prop.head-comp.Pks}), the permutation statistics
$\operatorname{Des}$, $\operatorname*{Lpk}$ and $\operatorname{Epk}$ are
head-graft-compatible; we will analyze a few other statistics in Subsection
\ref{subsect.LR.others}.

\begin{remark}
Let $\operatorname{st}$ be a head-graft-compatible permutation statistic.
Then, it is easy to see that%
\[
\operatorname{st}\left(  3,1,2\right)  =\operatorname{st}\left(
2,1,3\right)  \ \ \ \ \ \ \ \ \ \ \text{and}%
\ \ \ \ \ \ \ \ \ \ \operatorname{st}\left(  2,3,1\right)
=\operatorname{st}\left(  1,3,2\right)  .
\]
Moreover, these are the only restrictions that head-graft-compatibility places
on the values of $\operatorname{st}$ at $3$-permutations.
The restrictions placed on the values of $\operatorname{st}$ at
permutations of length $n > 3$ are more complicated, and depend
on its values on shorter permutations.
\end{remark}

It is usually easy to check if a given permutation statistic is
head-graft-compatible. For example:

\begin{proposition}
\label{prop.head-comp.Pks}\textbf{(a)} The permutation statistic
$\operatorname{Des}$ is head-graft-compatible.

\textbf{(b)} The permutation statistic $\operatorname*{Lpk}$ is head-graft-compatible.

\textbf{(c)} The permutation statistic $\operatorname{Epk}$ is head-graft-compatible.
\end{proposition}

\begin{proof}
[Proof of Proposition \ref{prop.head-comp.Pks}.]In this proof, we shall use
the following notation: If $S$ is a set of integers, and $p$ is an integer,
then $S+p$ shall denote the set $\left\{  s+p\ \mid\ s\in S\right\}  $.

\textbf{(a)} Let $\pi$ be a nonempty permutation. Let $a$ be a letter that
does not appear in $\pi$. We shall express
the element $\operatorname{Des}\left(  a:\pi\right)  $ in terms of
$\operatorname{Des}\pi$, $\left|\pi\right|$
and $\left[  a>\pi_{1}\right]  $.

Let $n=\left\vert \pi\right\vert $. Thus, $\pi=\left(  \pi_{1},\pi
_{2},\ldots,\pi_{n}\right)  $. Therefore, $a:\pi=\left(  a,\pi_{1},\pi_{2}%
,\ldots,\pi_{n}\right)  $. Hence, the descents of $a:\pi$ are obtained as follows:

\begin{itemize}
\item The number $1$ is a descent of $a:\pi$ if and only if $a>\pi_{1}$.

\item Adding $1$ to each descent of $\pi$ yields a descent of $a:\pi$. (That
is, if $i$ is a descent of $\pi$, then $i+1$ is a descent of $a:\pi$.)
\end{itemize}

These are all the descents of $a:\pi$. Thus,%
\begin{equation}
\operatorname{Des}\left(  a:\pi\right)  =\left\{  1\ \mid\ a>\pi_{1}\right\}
\cup\left(  \operatorname{Des}\pi+1\right)  .
\label{pf.prop.head-comp.Pks.a.1}%
\end{equation}
(The strange notation \textquotedblleft$\left\{  1\ \mid\ a>\pi_{1}\right\}
$\textquotedblright\ means exactly what it says: It is the set of all numbers
$1$ satisfying $a>\pi_{1}$. In other words, it is $\left\{  1\right\}  $ if
$a>\pi_{1}$, and $\varnothing$ otherwise.)

The equality (\ref{pf.prop.head-comp.Pks.a.1}) shows that $\operatorname{Des}
\left(  a:\pi\right)  $ depends only on $\operatorname{Des}\pi$, $\left\vert
\pi\right\vert $ and $\left[  a>\pi_{1}\right]  $ (indeed, the truth value
$\left[  a>\pi_{1}\right]  $ determines whether $a>\pi_{1}$ is true). In other
words, $\operatorname{Des}$ is head-graft-compatible (by the definition of
\textquotedblleft head-graft-compatible\textquotedblright). This proves
Proposition \ref{prop.head-comp.Pks} \textbf{(a)}.

\textbf{(b)} Let $\pi$ be a nonempty permutation. Let $a$ be a letter that
does not appear in $\pi$. We shall express
the element $\operatorname*{Lpk}\left(  a:\pi\right)  $ in terms of
$\operatorname*{Lpk}\pi$, $\left|\pi\right|$
and $\left[  a>\pi_{1}\right]  $.

Notice first that $a\neq\pi_{1}$ (since $a$ does not appear in $\pi$). Thus,
$a<\pi_{1}$ is true if and only if $a>\pi_{1}$ is false.

Let $n=\left\vert \pi\right\vert $. Therefore, $\pi=\left(  \pi_{1},\pi
_{2},\ldots,\pi_{n}\right)  $. Thus, $a:\pi=\left(  a,\pi_{1},\pi_{2}%
,\ldots,\pi_{n}\right)  $. Hence, the left peaks of $a:\pi$ are obtained as follows:

\begin{itemize}
\item The number $1$ is a left peak of $a:\pi$ if and only if $a>\pi_{1}$.

\item Adding $1$ to each left peak $i$ of $\pi$ yields a left peak $i+1$ of
$a:\pi$, except for the case when $i=1$ (in which case $i+1=2$ is a left peak
of $a:\pi$ only if $a<\pi_{1}$).
\end{itemize}

These are all the left peaks of $a:\pi$. Thus,%
\begin{equation}
\operatorname*{Lpk}\left(  a:\pi\right)  =\left\{  1\ \mid\ a>\pi_{1}\right\}
\cup%
\begin{cases}
\operatorname*{Lpk}\pi+1, & \text{if }a<\pi_{1};\\
\left(  \operatorname*{Lpk}\pi+1\right)  \setminus\left\{  2\right\}  , &
\text{if not }a<\pi_{1}%
\end{cases}
. \label{pf.prop.head-comp.Pks.b.1}%
\end{equation}

This equality shows that $\operatorname*{Lpk}\left(  a:\pi\right)  $ depends
only on $\operatorname*{Lpk}\pi$, $\left\vert \pi\right\vert $ and $\left[
a>\pi_{1}\right]  $ (indeed, the truth value $\left[  a>\pi_{1}\right]  $
determines whether $a>\pi_{1}$ is true and also determines whether $a<\pi_{1}$
is true\footnote{Indeed, $a<\pi_{1}$ is true if and only if $a>\pi_{1}$ is
false.}). In other words, $\operatorname*{Lpk}$ is head-graft-compatible (by
the definition of \textquotedblleft head-graft-compatible\textquotedblright).
This proves Proposition \ref{prop.head-comp.Pks} \textbf{(b)}.

\textbf{(c)} To obtain a proof of Proposition \ref{prop.head-comp.Pks}
\textbf{(c)}, it suffices to take our above proof of Proposition
\ref{prop.head-comp.Pks} \textbf{(b)} and replace every appearance of
\textquotedblleft left peak\textquotedblright\ and \textquotedblleft%
$\operatorname*{Lpk}$\textquotedblright\ by \textquotedblleft exterior
peak\textquotedblright\ and
\textquotedblleft$\operatorname{Epk}$\textquotedblright.
\end{proof}

\subsection{Proving LR-shuffle-compatibility}

Let us now state a sufficient criterion for the LR-shuffle-compatibility of a statistic:

\begin{theorem}
\label{thm.head-comp.LRcomp}Let $\operatorname{st}$ be a permutation
statistic that is both shuffle-compatible and head-graft-compatible. Then,
$\operatorname{st}$ is LR-shuffle-compatible.
\end{theorem}

Before we prove this theorem, let us introduce some terminology and state an
almost-trivial fact:

\begin{definition}
\textbf{(a)} If $A$ is a finite multiset, and if $g$ is any object, then
$\left\vert A\right\vert _{g}$ means the multiplicity of $g$ in $A$.

\textbf{(b)} If $A$ and $B$ are two finite multisets, then we say that
$B\subseteq A$ if and only if each object $g$ satisfies $\left\vert
B\right\vert _{g}\leq\left\vert A\right\vert _{g}$.

\textbf{(c)} If $A$ and $B$ are two finite multisets satisfying $B\subseteq
A$, then $A-B$ shall denote the \textquotedblleft multiset
difference\textquotedblright\ of $A$ and $B$; this is the finite
multiset $C$ such
that each object $g$ satisfies $\left\vert C\right\vert _{g}=\left\vert
A\right\vert _{g}-\left\vert B\right\vert _{g}$.
\end{definition}

For example, $\left\{  2,3,3\right\}  _{\operatorname*{multi}}\subseteq
\left\{  1,2,2,3,3\right\}  _{\operatorname*{multi}}$ and $\left\{
1,2,2,3,3\right\}  _{\operatorname*{multi}}-\left\{  2,3,3\right\}
_{\operatorname*{multi}}=\left\{  1,2\right\}  _{\operatorname*{multi}}$.

\begin{lemma}
\label{lem.LR.difference}Let $\pi$ and $\sigma$ be two disjoint permutations
such that at least one of $\pi$ and $\sigma$ is nonempty. Let
$\operatorname{st}$ be any permutation statistic. Then:

\textbf{(a)} We have%
\begin{align*}
&  \left\{  \operatorname{st} \tau   \ \mid\ \tau\in S_{\prec
}\left(  \pi,\sigma\right)  \right\}  _{\operatorname*{multi}}\\
&  =\left\{  \operatorname{st} \tau   \ \mid\ \tau\in S\left(
\pi,\sigma\right)  \right\}  _{\operatorname*{multi}}-\left\{
\operatorname{st} \tau   \ \mid\ \tau\in S_{\succ}\left(
\pi,\sigma\right)  \right\}  _{\operatorname*{multi}}.
\end{align*}

\textbf{(b)} We have%
\begin{align*}
&  \left\{  \operatorname{st} \tau   \ \mid\ \tau\in S_{\succ
}\left(  \pi,\sigma\right)  \right\}  _{\operatorname*{multi}}\\
&  =\left\{  \operatorname{st} \tau   \ \mid\ \tau\in S\left(
\pi,\sigma\right)  \right\}  _{\operatorname*{multi}}-\left\{
\operatorname{st} \tau   \ \mid\ \tau\in S_{\prec}\left(
\pi,\sigma\right)  \right\}  _{\operatorname*{multi}}.
\end{align*}

\end{lemma}

\begin{proof}
[Proof of Lemma \ref{lem.LR.difference}.]Recall that the two sets $S_{\prec
}\left(  \pi,\sigma\right)  $ and $S_{\succ}\left(  \pi,\sigma\right)  $ are
disjoint and their union is $S\left(  \pi,\sigma\right)  $. Thus, $S_{\succ
}\left(  \pi,\sigma\right)  \subseteq S\left(  \pi,\sigma\right)  $ and
$S_{\prec}\left(  \pi,\sigma\right)  =S\left(  \pi,\sigma\right)  \setminus
S_{\succ}\left(  \pi,\sigma\right)  $. Hence,%
\begin{align*}
&  \left\{  \operatorname{st} \tau   \ \mid\ \tau\in S_{\prec
}\left(  \pi,\sigma\right)  \right\}  _{\operatorname*{multi}}\\
&  =\left\{  \operatorname{st} \tau   \ \mid\ \tau\in S\left(
\pi,\sigma\right)  \right\}  _{\operatorname*{multi}}-\left\{
\operatorname{st} \tau   \ \mid\ \tau\in S_{\succ}\left(
\pi,\sigma\right)  \right\}  _{\operatorname*{multi}}.
\end{align*}
This proves Lemma \ref{lem.LR.difference} \textbf{(a)}. The proof of Lemma
\ref{lem.LR.difference} \textbf{(b)} is analogous.
\end{proof}

\begin{proof}
[Proof of Theorem \ref{thm.head-comp.LRcomp}.]We shall first show the following:

\begin{statement}
\textit{Claim 1:} Let $\pi$, $\pi^{\prime}$ and $\sigma$ be three nonempty
permutations. Assume that $\pi$ and $\sigma$ are disjoint. Assume that
$\pi^{\prime}$ and $\sigma$ are disjoint. Assume furthermore that%
\[
\operatorname{st} \pi   =\operatorname{st}\left(  \pi^{\prime
}\right)  ,\ \ \ \ \ \ \ \ \ \ \left\vert \pi\right\vert =\left\vert
\pi^{\prime}\right\vert \ \ \ \ \ \ \ \ \ \ \text{and}%
\ \ \ \ \ \ \ \ \ \ \left[  \pi_{1}>\sigma_{1}\right]  =\left[  \pi
_{1}^{\prime}>\sigma_{1}\right]  .
\]
Then,
\begin{align}
&  \left\{  \operatorname{st} \tau   \ \mid\ \tau\in S_{\prec
}\left(  \pi,\sigma\right)  \right\}  _{\operatorname*{multi}}\nonumber\\
&  =\left\{  \operatorname{st} \tau   \ \mid\ \tau\in S_{\prec
}\left(  \pi^{\prime},\sigma\right)  \right\}  _{\operatorname*{multi}}
\label{pf.thm.head-comp.LRcomp.c1.eq3}%
\end{align}
and%
\begin{align}
&  \left\{  \operatorname{st} \tau   \ \mid\ \tau\in S_{\succ
}\left(  \pi,\sigma\right)  \right\}  _{\operatorname*{multi}}\nonumber\\
&  =\left\{  \operatorname{st} \tau   \ \mid\ \tau\in S_{\succ
}\left(  \pi^{\prime},\sigma\right)  \right\}  _{\operatorname*{multi}}.
\label{pf.thm.head-comp.LRcomp.c1.eq4}%
\end{align}

\end{statement}

[\textit{Proof of Claim 1:} We shall prove Claim 1 by induction on $\left\vert
\sigma\right\vert $:

\textit{Induction base:} The case $\left\vert \sigma\right\vert =0$ cannot
happen (because $\sigma$ is assumed to be nonempty). Thus, Claim 1 is true in
the case $\left\vert \sigma\right\vert =0$. This completes the induction base.

\textit{Induction step:} Let $N$ be a positive integer. Assume (as the
induction hypothesis) that Claim 1 holds when $\left\vert \sigma\right\vert
=N-1$. We must now prove that Claim 1 holds when $\left\vert \sigma\right\vert
=N$.

Indeed, let $\pi$, $\pi^{\prime}$ and $\sigma$ be as in Claim 1, and assume
that $\left\vert \sigma\right\vert =N$. We must prove
(\ref{pf.thm.head-comp.LRcomp.c1.eq3}) and
(\ref{pf.thm.head-comp.LRcomp.c1.eq4}).

Proposition \ref{prop.LR.rec} \textbf{(c)} yields that the permutations
$\sigma_{\sim1}$ and $\sigma_{1}:\pi$ are well-defined and disjoint, and
satisfy%
\begin{equation}
S_{\succ}\left(  \pi,\sigma\right)  =S_{\prec}\left(  \sigma_{1}:\pi
,\sigma_{\sim1}\right)  . \label{pf.thm.head-comp.LRcomp.c1.pf.1}%
\end{equation}
Furthermore, $\left\vert \sigma_{\sim1}\right\vert =\left\vert \sigma
\right\vert -1=N-1$ (since $\left\vert \sigma\right\vert =N$).

Proposition \ref{prop.LR.rec} \textbf{(c)} (applied to $\pi^{\prime}$ instead
of $\pi$) yields that the permutations $\sigma_{\sim1}$ and $\sigma_{1}%
:\pi^{\prime}$ are well-defined and disjoint, and satisfy%
\begin{equation}
S_{\succ}\left(  \pi^{\prime},\sigma\right)  =S_{\prec}\left(  \sigma_{1}%
:\pi^{\prime},\sigma_{\sim1}\right)  . \label{pf.thm.head-comp.LRcomp.c1.pf.2}%
\end{equation}

The letter $\sigma_{1}$ does not appear in the permutation $\pi$ (since $\pi$
and $\sigma$ are disjoint). Similarly, the letter $\sigma_{1}$ does not appear
in the permutation $\pi^{\prime}$. Also, $\left\vert \sigma_{1}:\pi\right\vert
=\underbrace{\left\vert \pi\right\vert }_{=\left\vert \pi^{\prime}\right\vert
}+1=\left\vert \pi^{\prime}\right\vert +1=\left\vert \sigma_{1}:\pi^{\prime
}\right\vert $.

We have $\sigma_{1}\neq\pi_{1}$ (since $\pi$ and $\sigma$ are disjoint). Thus,
the statement $\left(  \sigma_{1}>\pi_{1}\right)  $ is equivalent to $\left(
\text{not }\pi_{1}>\sigma_{1}\right)  $. Hence, $\left[  \sigma_{1}>\pi
_{1}\right]  =\left[  \text{not }\pi_{1}>\sigma_{1}\right]  =1-\left[  \pi
_{1}>\sigma_{1}\right]  $. Similarly, $\left[  \sigma_{1}>\pi_{1}^{\prime
}\right]  =1-\left[  \pi_{1}^{\prime}>\sigma_{1}\right]  $. Hence,%
\[
\left[  \sigma_{1}>\pi_{1}\right]  =1-\underbrace{\left[  \pi_{1}>\sigma
_{1}\right]  }_{=\left[  \pi_{1}^{\prime}>\sigma_{1}\right]  }=1-\left[
\pi_{1}^{\prime}>\sigma_{1}\right]  =\left[  \sigma_{1}>\pi_{1}^{\prime
}\right]  .
\]

Both permutations $\sigma_{1}:\pi$ and $\sigma_{1}:\pi^{\prime}$ begin with
the letter $\sigma_{1}$. Thus, both $\left(  \sigma_{1}:\pi\right)  _{1}$ and
$\left(  \sigma_{1}:\pi^{\prime}\right)  _{1}$ equal $\sigma_{1}$. Hence,
$\left(  \sigma_{1}:\pi\right)  _{1}=\left(  \sigma_{1}:\pi^{\prime}\right)
_{1}$.

The statistic $\operatorname{st}$ is head-graft-compatible. In other words,
for any nonempty permutation $\varphi$ and any letter $a$ that does not appear
in $\varphi$, the element $\operatorname{st}\left(  a:\varphi\right)  $
depends only on $\operatorname{st}\left(  \varphi\right)  $, $\left\vert
\varphi\right\vert $ and $\left[  a>\varphi_{1}\right]  $ (by the definition
of \textquotedblleft head-graft-compatible\textquotedblright). Hence,
if $\varphi$ and $\varphi^{\prime}$ are two nonempty permutations, and if $a$
is any letter that does not appear in $\varphi$ and does not appear in
$\varphi^{\prime}$, and if we have $\operatorname{st} \varphi
=\operatorname{st}\left(  \varphi^{\prime}\right)  $ and $\left\vert
\varphi\right\vert =\left\vert \varphi^{\prime}\right\vert $ and $\left[
a>\varphi_{1}\right]  =\left[  a>\varphi_{1}^{\prime}\right]  $, then
$\operatorname{st}\left(  a:\varphi\right)  =\operatorname{st}\left(
a:\varphi^{\prime}\right)  $. Applying this to $a=\sigma_{1}$, $\varphi=\pi$
and $\varphi^{\prime}=\pi^{\prime}$, we obtain
\[
\operatorname{st}\left(  \sigma_{1}:\pi\right)  =\operatorname{st}\left(
\sigma_{1}:\pi^{\prime}\right)
\]
(since $\operatorname{st} \pi   =\operatorname{st}\left(
\pi^{\prime}\right)  $ and $\left\vert \pi\right\vert =\left\vert \pi^{\prime
}\right\vert $ and $\left[  \sigma_{1}>\pi_{1}\right]  =\left[  \sigma_{1}%
>\pi_{1}^{\prime}\right]  $).

Next, we claim that%
\begin{align}
&\left\{  \operatorname{st} \tau   \ \mid\ \tau\in S_{\prec
}\left(  \sigma_{1}:\pi,\sigma_{\sim1}\right)  \right\}
_{\operatorname*{multi}}
\nonumber\\
&=\left\{
\operatorname{st} \tau   \ \mid\ \tau\in S_{\prec}\left(
\sigma_{1}:\pi^{\prime},\sigma_{\sim1}\right)  \right\}
_{\operatorname*{multi}}.
\label{pf.thm.head-comp.LRcomp.c1.pf.5}%
\end{align}

[\textit{Proof of (\ref{pf.thm.head-comp.LRcomp.c1.pf.5}):} The permutations
$\sigma_{1}:\pi$ and $\sigma_{1}:\pi^{\prime}$ are clearly nonempty. Hence, if
$\sigma_{\sim1}$ is the $0$-permutation $\left(  {}\right)  $, then
$S_{\prec}\left(  \sigma_{1}:\pi,\sigma_{\sim1}\right)  =\left\{  \sigma
_{1}:\pi\right\}  $ and $S_{\prec}\left(  \sigma_{1}:\pi^{\prime},\sigma
_{\sim1}\right)  =\left\{  \sigma_{1}:\pi^{\prime}\right\}  $. Thus, if
$\sigma_{\sim1}$ is the $0$-permutation $\left(  {}\right)  $, then
(\ref{pf.thm.head-comp.LRcomp.c1.pf.5}) follows from%
\begin{align*}
&  \left\{  \operatorname{st} \tau   \ \mid\ \tau\in
\underbrace{S_{\prec}\left(  \sigma_{1}:\pi,\sigma_{\sim1}\right)
}_{=\left\{  \sigma_{1}:\pi\right\}  }\right\}  _{\operatorname*{multi}}\\
&  =\left\{  \operatorname{st} \tau   \ \mid\ \tau\in\left\{
\sigma_{1}:\pi\right\}  \right\}  _{\operatorname*{multi}}=\left\{
\underbrace{\operatorname{st}\left(  \sigma_{1}:\pi\right)  }%
_{=\operatorname{st}\left(  \sigma_{1}:\pi^{\prime}\right)  }\right\}
_{\operatorname*{multi}}\\
&  =\left\{  \operatorname{st}\left(  \sigma_{1}:\pi^{\prime}\right)
\right\}  _{\operatorname*{multi}}=\left\{  \operatorname{st} \tau
\ \mid\ \tau\in\underbrace{\left\{  \sigma_{1}:\pi^{\prime
}\right\}  }_{=S_{\prec}\left(  \sigma_{1}:\pi^{\prime},\sigma_{\sim1}\right)
}\right\}  _{\operatorname*{multi}}\\
&  =\left\{  \operatorname{st} \tau   \ \mid\ \tau\in S_{\prec
}\left(  \sigma_{1}:\pi^{\prime},\sigma_{\sim1}\right)  \right\}
_{\operatorname*{multi}}.
\end{align*}
Thus, for the rest of our proof of (\ref{pf.thm.head-comp.LRcomp.c1.pf.5}), we
WLOG assume that $\sigma_{\sim1}$ is not the $0$-permutation $\left(
{}\right)  $. Thus, $\sigma_{\sim1}$ is nonempty.

But recall that $\left\vert \sigma_{\sim1}\right\vert =N-1$. Hence, the
induction hypothesis allows us to apply Claim 1 to $\sigma_{1}:\pi$,
$\sigma_{1}:\pi^{\prime}$ and $\sigma_{\sim1}$ instead of $\pi$, $\pi^{\prime
}$ and $\sigma$ (because we know that the permutations $\sigma_{\sim1}$ and
$\sigma_{1}:\pi$ are disjoint; that the permutations $\sigma_{\sim1}$ and
$\sigma_{1}:\pi^{\prime}$ are disjoint; that $\operatorname{st}\left(
\sigma_{1}:\pi\right)  =\operatorname{st}\left(  \sigma_{1}:\pi^{\prime
}\right)  $ and $\left\vert \sigma_{1}:\pi\right\vert =\left\vert \sigma
_{1}:\pi^{\prime}\right\vert $; and that $\left[  \underbrace{\left(
\sigma_{1}:\pi\right)  _{1}}_{=\left(  \sigma_{1}:\pi^{\prime}\right)  _{1}%
}>\left(  \sigma_{\sim1}\right)  _{1}\right]  =\left[  \left(  \sigma_{1}%
:\pi^{\prime}\right)  _{1}>\left(  \sigma_{\sim1}\right)  _{1}\right]  $). We
therefore obtain%
\[
\left\{  \operatorname{st} \tau   \ \mid\ \tau\in S_{\prec
}\left(  \sigma_{1}:\pi,\sigma_{\sim1}\right)  \right\}
_{\operatorname*{multi}}=\left\{  \operatorname{st} \tau 
\ \mid\ \tau\in S_{\prec}\left(  \sigma_{1}:\pi^{\prime},\sigma_{\sim
1}\right)  \right\}  _{\operatorname*{multi}}%
\]
and%
\[
\left\{  \operatorname{st} \tau   \ \mid\ \tau\in S_{\succ
}\left(  \sigma_{1}:\pi,\sigma_{\sim1}\right)  \right\}
_{\operatorname*{multi}}=\left\{  \operatorname{st} \tau 
\ \mid\ \tau\in S_{\succ}\left(  \sigma_{1}:\pi^{\prime},\sigma_{\sim
1}\right)  \right\}  _{\operatorname*{multi}}.
\]
The first of these two equalities is precisely
(\ref{pf.thm.head-comp.LRcomp.c1.pf.5}). Thus,
(\ref{pf.thm.head-comp.LRcomp.c1.pf.5}) is proven.]

Now,
\begin{align}
&  \left\{  \operatorname{st} \tau   \ \mid\ \tau\in
\underbrace{S_{\succ}\left(  \pi,\sigma\right)  }_{\substack{=S_{\prec}\left(
\sigma_{1}:\pi,\sigma_{\sim1}\right)  \\\text{(by
(\ref{pf.thm.head-comp.LRcomp.c1.pf.1}))}}}\right\}  _{\operatorname*{multi}%
}\nonumber\\
&  =\left\{  \operatorname{st} \tau   \ \mid\ \tau\in S_{\prec
}\left(  \sigma_{1}:\pi,\sigma_{\sim1}\right)  \right\}
_{\operatorname*{multi}}\nonumber\\
&  =\left\{  \operatorname{st} \tau   \ \mid\ \tau
\in\underbrace{S_{\prec}\left(  \sigma_{1}:\pi^{\prime},\sigma_{\sim1}\right)
}_{\substack{=S_{\succ}\left(  \pi^{\prime},\sigma\right)  \\\text{(by
(\ref{pf.thm.head-comp.LRcomp.c1.pf.2}))}}}\right\}  _{\operatorname*{multi}%
}\ \ \ \ \ \ \ \ \ \ \left(  \text{by (\ref{pf.thm.head-comp.LRcomp.c1.pf.5}%
)}\right) \nonumber\\
&  =\left\{  \operatorname{st} \tau   \ \mid\ \tau\in S_{\succ
}\left(  \pi^{\prime},\sigma\right)  \right\}  _{\operatorname*{multi}}.
\label{pf.thm.head-comp.LRcomp.c1.pf.6}%
\end{align}
This proves (\ref{pf.thm.head-comp.LRcomp.c1.eq4}). It remains to prove
(\ref{pf.thm.head-comp.LRcomp.c1.eq3}).

Lemma \ref{lem.LR.difference} \textbf{(a)} yields%
\begin{align}
&  \left\{  \operatorname{st} \tau   \ \mid\ \tau\in S_{\prec
}\left(  \pi,\sigma\right)  \right\}  _{\operatorname*{multi}}\nonumber\\
&  =\left\{  \operatorname{st} \tau   \ \mid\ \tau\in S\left(
\pi,\sigma\right)  \right\}  _{\operatorname*{multi}}-\left\{
\operatorname{st} \tau   \ \mid\ \tau\in S_{\succ}\left(
\pi,\sigma\right)  \right\}  _{\operatorname*{multi}}.
\label{pf.thm.head-comp.LRcomp.c1.pf.7}%
\end{align}
Lemma \ref{lem.LR.difference} \textbf{(a)} (applied to $\pi^{\prime}$ instead
of $\pi$) yields%
\begin{align}
&  \left\{  \operatorname{st} \tau   \ \mid\ \tau\in S_{\prec
}\left(  \pi^{\prime},\sigma\right)  \right\}  _{\operatorname*{multi}%
}\nonumber\\
&  =\left\{  \operatorname{st} \tau   \ \mid\ \tau\in S\left(
\pi^{\prime},\sigma\right)  \right\}  _{\operatorname*{multi}}-\left\{
\operatorname{st} \tau   \ \mid\ \tau\in S_{\succ}\left(
\pi^{\prime},\sigma\right)  \right\}  _{\operatorname*{multi}}.
\label{pf.thm.head-comp.LRcomp.c1.pf.8}%
\end{align}

But recall that the statistic $\operatorname{st}$ is shuffle-compatible. In
other words, for any two disjoint permutations $\alpha$ and $\beta$, the
multiset%
\[
\left\{  \operatorname{st} \tau   \ \mid\ \tau\in S\left(
\alpha,\beta\right)  \right\}  _{\operatorname*{multi}}%
\]
depends only on $\operatorname{st} \alpha$,
$\operatorname{st} \beta$, $\left\vert \alpha\right\vert $
and $\left\vert \beta\right\vert $ (by the definition of
shuffle-compatibility). In other words, if $\alpha$ and $\beta$ are two
disjoint permutations, and if $\alpha^{\prime}$ and $\beta^{\prime}$ are two
disjoint permutations, and if we have
\[
\operatorname{st} \alpha = \operatorname{st}\left(
\alpha^{\prime}\right)  ,\ \ \ \ \ \ \ \ \ \ \operatorname{st} \beta
=\operatorname{st}\left(  \beta^{\prime}\right)
,\ \ \ \ \ \ \ \ \ \ \left\vert \alpha\right\vert =\left\vert \alpha^{\prime
}\right\vert \ \ \ \ \ \ \ \ \ \ \text{and}\ \ \ \ \ \ \ \ \ \ \left\vert
\beta\right\vert =\left\vert \beta^{\prime}\right\vert ,
\]
then%
\[
\left\{  \operatorname{st} \tau   \ \mid\ \tau\in S\left(
\alpha,\beta\right)  \right\}  _{\operatorname*{multi}}=\left\{
\operatorname{st} \tau   \ \mid\ \tau\in S\left(  \alpha
^{\prime},\beta^{\prime}\right)  \right\}  _{\operatorname*{multi}}.
\]
Applying this to $\alpha=\pi$, $\beta=\sigma$, $\alpha^{\prime}=\pi^{\prime}$
and $\beta^{\prime}=\sigma$, we obtain
\begin{equation}
\left\{  \operatorname{st} \tau   \ \mid\ \tau\in S\left(
\pi,\sigma\right)  \right\}  _{\operatorname*{multi}}=\left\{
\operatorname{st} \tau   \ \mid\ \tau\in S\left(  \pi^{\prime
},\sigma\right)  \right\}  _{\operatorname*{multi}}
\label{pf.thm.head-comp.LRcomp.c1.pf.9}%
\end{equation}
(since $\operatorname{st} \pi   =\operatorname{st}\left(
\pi^{\prime}\right)  $, $\operatorname{st} \sigma 
=\operatorname{st} \sigma   $, $\left\vert \pi\right\vert
=\left\vert \pi^{\prime}\right\vert $ and $\left\vert \sigma\right\vert
=\left\vert \sigma\right\vert $). Now, (\ref{pf.thm.head-comp.LRcomp.c1.pf.7})
becomes%
\begin{align*}
&  \left\{  \operatorname{st} \tau   \ \mid\ \tau\in S_{\prec
}\left(  \pi,\sigma\right)  \right\}  _{\operatorname*{multi}}\\
&  =\underbrace{\left\{  \operatorname{st} \tau   \ \mid
\ \tau\in S\left(  \pi,\sigma\right)  \right\}  _{\operatorname*{multi}}%
}_{\substack{=\left\{  \operatorname{st} \tau   \ \mid\ \tau\in
S\left(  \pi^{\prime},\sigma\right)  \right\}  _{\operatorname*{multi}%
}\\\text{(by (\ref{pf.thm.head-comp.LRcomp.c1.pf.9}))}}}-\underbrace{\left\{
\operatorname{st} \tau   \ \mid\ \tau\in S_{\succ}\left(
\pi,\sigma\right)  \right\}  _{\operatorname*{multi}}}_{\substack{=\left\{
\operatorname{st} \tau   \ \mid\ \tau\in S_{\succ}\left(
\pi^{\prime},\sigma\right)  \right\}  _{\operatorname*{multi}}\\\text{(by
(\ref{pf.thm.head-comp.LRcomp.c1.pf.6}))}}}\\
&  =\left\{  \operatorname{st} \tau   \ \mid\ \tau\in S\left(
\pi^{\prime},\sigma\right)  \right\}  _{\operatorname*{multi}}-\left\{
\operatorname{st} \tau   \ \mid\ \tau\in S_{\succ}\left(
\pi^{\prime},\sigma\right)  \right\}  _{\operatorname*{multi}}\\
&  =\left\{  \operatorname{st} \tau   \ \mid\ \tau\in S_{\prec
}\left(  \pi^{\prime},\sigma\right)  \right\}  _{\operatorname*{multi}%
}\ \ \ \ \ \ \ \ \ \ \left(  \text{by (\ref{pf.thm.head-comp.LRcomp.c1.pf.8}%
)}\right)  .
\end{align*}
Thus, (\ref{pf.thm.head-comp.LRcomp.c1.eq3}) is proven. Hence, we have proven
both (\ref{pf.thm.head-comp.LRcomp.c1.eq3}) and
(\ref{pf.thm.head-comp.LRcomp.c1.eq4}). This shows that Claim 1 holds for our
$\pi$, $\pi^{\prime}$ and $\sigma$. This completes the induction step. Thus,
Claim 1 is proven by induction.]

We shall next derive a \textquotedblleft mirror version\textquotedblright\ of
Claim 1:

\begin{statement}
\textit{Claim 2:} Let $\pi$, $\sigma$ and $\sigma^{\prime}$ be three nonempty
permutations. Assume that $\pi$ and $\sigma$ are disjoint. Assume that $\pi$
and $\sigma^{\prime}$ are disjoint. Assume furthermore that%
\[
\operatorname{st} \sigma   =\operatorname{st}\left(
\sigma^{\prime}\right)  ,\ \ \ \ \ \ \ \ \ \ \left\vert \sigma\right\vert
=\left\vert \sigma^{\prime}\right\vert \ \ \ \ \ \ \ \ \ \ \text{and}%
\ \ \ \ \ \ \ \ \ \ \left[  \pi_{1}>\sigma_{1}\right]  =\left[  \pi_{1}%
>\sigma_{1}^{\prime}\right]  .
\]
Then,
\begin{align*}
&  \left\{  \operatorname{st} \tau   \ \mid\ \tau\in S_{\prec
}\left(  \pi,\sigma\right)  \right\}  _{\operatorname*{multi}}\\
&  =\left\{  \operatorname{st} \tau   \ \mid\ \tau\in S_{\prec
}\left(  \pi,\sigma^{\prime}\right)  \right\}  _{\operatorname*{multi}}%
\end{align*}
and%
\begin{align*}
&  \left\{  \operatorname{st} \tau   \ \mid\ \tau\in S_{\succ
}\left(  \pi,\sigma\right)  \right\}  _{\operatorname*{multi}}\\
&  =\left\{  \operatorname{st} \tau   \ \mid\ \tau\in S_{\succ
}\left(  \pi,\sigma^{\prime}\right)  \right\}  _{\operatorname*{multi}}.
\end{align*}

\end{statement}

[\textit{Proof of Claim 2:} We have $\sigma_{1}\neq\pi_{1}$ (since $\pi$ and
$\sigma$ are disjoint). Thus, the statement $\left(  \sigma_{1}>\pi
_{1}\right)  $ is equivalent to $\left(  \text{not }\pi_{1}>\sigma_{1}\right)
$. Hence, $\left[  \sigma_{1}>\pi_{1}\right]  =\left[  \text{not }\pi
_{1}>\sigma_{1}\right]  =1-\left[  \pi_{1}>\sigma_{1}\right]  $. Similarly,
$\left[  \sigma_{1}^{\prime}>\pi_{1}\right]  =1-\left[  \pi_{1}>\sigma
_{1}^{\prime}\right]  $. Hence,%
\[
\left[  \sigma_{1}>\pi_{1}\right]  =1-\underbrace{\left[  \pi_{1}>\sigma
_{1}\right]  }_{=\left[  \pi_{1}>\sigma_{1}^{\prime}\right]  }=1-\left[
\pi_{1}>\sigma_{1}^{\prime}\right]  =\left[  \sigma_{1}^{\prime}>\pi
_{1}\right]  .
\]
Hence, Claim 1 (applied to $\sigma$, $\sigma^{\prime}$ and $\pi$ instead of
$\pi$, $\pi^{\prime}$ and $\sigma$) shows that%
\begin{align*}
&  \left\{  \operatorname{st} \tau   \ \mid\ \tau\in S_{\prec
}\left(  \sigma,\pi\right)  \right\}  _{\operatorname*{multi}}\\
&  =\left\{  \operatorname{st} \tau   \ \mid\ \tau\in S_{\prec
}\left(  \sigma^{\prime},\pi\right)  \right\}  _{\operatorname*{multi}}%
\end{align*}
and%
\begin{align*}
&  \left\{  \operatorname{st} \tau   \ \mid\ \tau\in S_{\succ
}\left(  \sigma,\pi\right)  \right\}  _{\operatorname*{multi}}\\
&  =\left\{  \operatorname{st} \tau   \ \mid\ \tau\in S_{\succ
}\left(  \sigma^{\prime},\pi\right)  \right\}  _{\operatorname*{multi}}.
\end{align*}
But Proposition \ref{prop.LR.rec} \textbf{(a)} yields $S_{\prec}\left(
\pi,\sigma\right)  =S_{\succ}\left(  \sigma,\pi\right)  $. Similarly,
$S_{\prec}\left(  \pi,\sigma^{\prime}\right)  =S_{\succ}\left(  \sigma
^{\prime},\pi\right)  $. Also, Proposition \ref{prop.LR.rec} \textbf{(a)}
(applied to $\sigma$ and $\pi$ instead of $\pi$ and $\sigma$) yields
$S_{\prec}\left(  \sigma,\pi\right)  =S_{\succ}\left(  \pi,\sigma\right)  $.
Similarly, $S_{\prec}\left(  \sigma^{\prime},\pi\right)  =S_{\succ}\left(
\pi,\sigma^{\prime}\right)  $. Using all these equalities, we find%
\begin{align*}
&  \left\{  \operatorname{st} \tau   \ \mid\ \tau\in
\underbrace{S_{\prec}\left(  \pi,\sigma\right)  }_{=S_{\succ}\left(
\sigma,\pi\right)  }\right\}  _{\operatorname*{multi}}\\
&  =\left\{  \operatorname{st} \tau   \ \mid\ \tau\in S_{\succ
}\left(  \sigma,\pi\right)  \right\}  _{\operatorname*{multi}}=\left\{
\operatorname{st} \tau   \ \mid\ \tau\in\underbrace{S_{\succ
}\left(  \sigma^{\prime},\pi\right)  }_{=S_{\prec}\left(  \pi,\sigma^{\prime
}\right)  }\right\}  _{\operatorname*{multi}}\\
&  =\left\{  \operatorname{st} \tau   \ \mid\ \tau\in S_{\prec
}\left(  \pi,\sigma^{\prime}\right)  \right\}  _{\operatorname*{multi}}%
\end{align*}
and%
\begin{align*}
&  \left\{  \operatorname{st} \tau   \ \mid\ \tau\in
\underbrace{S_{\succ}\left(  \pi,\sigma\right)  }_{=S_{\prec}\left(
\sigma,\pi\right)  }\right\}  _{\operatorname*{multi}}\\
&  =\left\{  \operatorname{st} \tau   \ \mid\ \tau\in S_{\prec
}\left(  \sigma,\pi\right)  \right\}  _{\operatorname*{multi}}=\left\{
\operatorname{st} \tau   \ \mid\ \tau\in\underbrace{S_{\prec
}\left(  \sigma^{\prime},\pi\right)  }_{=S_{\succ}\left(  \pi,\sigma^{\prime
}\right)  }\right\}  _{\operatorname*{multi}}\\
&  =\left\{  \operatorname{st} \tau   \ \mid\ \tau\in S_{\succ
}\left(  \pi,\sigma^{\prime}\right)  \right\}  _{\operatorname*{multi}}.
\end{align*}
Thus, Claim 2 is proven.]

Finally, we are ready to take on the LR-shuffle-compatibility of
$\operatorname{st}$:

\begin{statement}
\textit{Claim 3:} Let $\pi$ and $\sigma$ be two disjoint nonempty
permutations. Let $\pi^{\prime}$ and $\sigma^{\prime}$ be two disjoint
nonempty permutations. Assume that%
\begin{align*}
\operatorname{st} \pi    &  =\operatorname{st}\left(
\pi^{\prime}\right)  ,\ \ \ \ \ \ \ \ \ \ \operatorname{st}
\sigma  =\operatorname{st}\left(  \sigma^{\prime}\right)  ,\\
\left\vert \pi\right\vert  &  =\left\vert \pi^{\prime}\right\vert
,\ \ \ \ \ \ \ \ \ \ \left\vert \sigma\right\vert =\left\vert \sigma^{\prime
}\right\vert \ \ \ \ \ \ \ \ \ \ \text{and}\ \ \ \ \ \ \ \ \ \ \left[  \pi
_{1}>\sigma_{1}\right]  =\left[  \pi_{1}^{\prime}>\sigma_{1}^{\prime}\right]
.
\end{align*}
Then,
\[
\left\{  \operatorname{st} \tau   \ \mid\ \tau\in S_{\prec
}\left(  \pi,\sigma\right)  \right\}  _{\operatorname*{multi}}=\left\{
\operatorname{st} \tau   \ \mid\ \tau\in S_{\prec}\left(
\pi^{\prime},\sigma^{\prime}\right)  \right\}  _{\operatorname*{multi}}%
\]
and%
\[
\left\{  \operatorname{st} \tau   \ \mid\ \tau\in S_{\succ
}\left(  \pi,\sigma\right)  \right\}  _{\operatorname*{multi}}=\left\{
\operatorname{st} \tau   \ \mid\ \tau\in S_{\succ}\left(
\pi^{\prime},\sigma^{\prime}\right)  \right\}  _{\operatorname*{multi}}.
\]

\end{statement}

[\textit{Proof of Claim 3:} We have $\left[  \pi_{1}>\sigma_{1}\right]
=\left[  \pi_{1}^{\prime}>\sigma_{1}^{\prime}\right]  $. Since $\left[
\pi_{1}>\sigma_{1}\right]  =\left[  \pi_{1}^{\prime}>\sigma_{1}^{\prime
}\right]  $ is either $1$ or $0$, we must therefore be in one of the following
two cases:

\textit{Case 1:} We have $\left[  \pi_{1}>\sigma_{1}\right]  =\left[  \pi
_{1}^{\prime}>\sigma_{1}^{\prime}\right]  =1$.

\textit{Case 2:} We have $\left[  \pi_{1}>\sigma_{1}\right]  =\left[  \pi
_{1}^{\prime}>\sigma_{1}^{\prime}\right]  =0$.

Let us first consider Case 1. In this case, we have $\left[  \pi_{1}%
>\sigma_{1}\right]  =\left[  \pi_{1}^{\prime}>\sigma_{1}^{\prime}\right]  =1$.

There clearly exists a positive integer $N$ that is larger than all entries of
$\sigma$ and larger than all entries of $\sigma^{\prime}$. Consider such an
$N$. Let $n=\left\vert \pi\right\vert $; thus, $\pi=\left(  \pi_{1},\pi
_{2},\ldots,\pi_{n}\right)  $. Let $\gamma$ be the permutation $\left(
\pi_{1}+N,\pi_{2}+N,\ldots,\pi_{n}+N\right)  $. This permutation $\gamma$ is
order-isomorphic to $\pi$, but is disjoint from $\sigma$ (since all its
entries are $>N$, while all the entries of $\sigma$ are $<N$) and disjoint
from $\sigma^{\prime}$ (for similar reasons). Also, $\gamma_{1}%
=\underbrace{\pi_{1}}_{>0}+N>N>\sigma_{1}$ (since $N$ is larger than all
entries of $\sigma$), so that $\left[  \gamma_{1}%
>\sigma_{1}\right]  =1$. Similarly, $\left[  \gamma_{1}>\sigma_{1}^{\prime
}\right]  =1$.

The permutation $\gamma$ is order-isomorphic to $\pi$. Thus,
$\operatorname{st} \gamma = \operatorname{st} \pi$ (since $\operatorname{st}$
is a permutation statistic) and
$\left\vert \gamma\right\vert =\left\vert \pi\right\vert $. The permutation
$\gamma$ is furthermore nonempty (since it is order-isomorphic to the nonempty
permutation $\pi$). Also, $\operatorname{st} \gamma
=\operatorname{st} \pi =\operatorname{st}\left(  \pi^{\prime
}\right)  $ and $\left\vert \gamma\right\vert =\left\vert \pi\right\vert
=\left\vert \pi^{\prime}\right\vert $. Moreover, $\left[  \pi_{1}>\sigma
_{1}\right]  =1=\left[  \gamma_{1}>\sigma_{1}\right]  $ and $\left[
\gamma_{1}>\sigma_{1}\right]  =1=\left[  \gamma_{1}>\sigma_{1}^{\prime
}\right]  $ and
$\left[  \gamma_{1}>\sigma^{\prime}_{1}\right]  =1=\left[  \pi
_{1}^{\prime}>\sigma_{1}^{\prime}\right]  $. Hence, Claim 1 (applied to
$\gamma$ instead of $\pi^{\prime}$) yields%
\[
\left\{  \operatorname{st} \tau   \ \mid\ \tau\in S_{\prec
}\left(  \pi,\sigma\right)  \right\}  _{\operatorname*{multi}}=\left\{
\operatorname{st} \tau   \ \mid\ \tau\in S_{\prec}\left(
\gamma,\sigma\right)  \right\}  _{\operatorname*{multi}}%
\]
and%
\begin{equation}
\left\{  \operatorname{st} \tau   \ \mid\ \tau\in S_{\succ
}\left(  \pi,\sigma\right)  \right\}  _{\operatorname*{multi}}=\left\{
\operatorname{st} \tau   \ \mid\ \tau\in S_{\succ}\left(
\gamma,\sigma\right)  \right\}  _{\operatorname*{multi}}.\nonumber
\end{equation}
Furthermore, Claim 2 (applied to $\gamma$ instead of $\pi$) yields%
\[
\left\{  \operatorname{st} \tau   \ \mid\ \tau\in S_{\prec
}\left(  \gamma,\sigma\right)  \right\}  _{\operatorname*{multi}}=\left\{
\operatorname{st} \tau   \ \mid\ \tau\in S_{\prec}\left(
\gamma,\sigma^{\prime}\right)  \right\}  _{\operatorname*{multi}}%
\]
and%
\[
\left\{  \operatorname{st} \tau   \ \mid\ \tau\in S_{\succ
}\left(  \gamma,\sigma\right)  \right\}  _{\operatorname*{multi}}=\left\{
\operatorname{st} \tau   \ \mid\ \tau\in S_{\succ}\left(
\gamma,\sigma^{\prime}\right)  \right\}  _{\operatorname*{multi}}.
\]
Finally, Claim 1 (applied to $\gamma$ and $\sigma^{\prime}$ instead of $\pi$
and $\sigma$) yields%
\[
\left\{  \operatorname{st} \tau   \ \mid\ \tau\in S_{\prec
}\left(  \gamma,\sigma^{\prime}\right)  \right\}  _{\operatorname*{multi}%
}=\left\{  \operatorname{st} \tau   \ \mid\ \tau\in S_{\prec
}\left(  \pi^{\prime},\sigma^{\prime}\right)  \right\}
_{\operatorname*{multi}}%
\]
and%
\begin{equation}
\left\{  \operatorname{st} \tau   \ \mid\ \tau\in S_{\succ
}\left(  \gamma,\sigma^{\prime}\right)  \right\}  _{\operatorname*{multi}%
}=\left\{  \operatorname{st} \tau   \ \mid\ \tau\in S_{\succ
}\left(  \pi^{\prime},\sigma^{\prime}\right)  \right\}
_{\operatorname*{multi}}.\nonumber
\end{equation}
Combining the equalities we have found, we obtain%
\begin{align*}
\left\{  \operatorname{st} \tau   \ \mid\ \tau\in S_{\prec
}\left(  \pi,\sigma\right)  \right\}  _{\operatorname*{multi}}  &  =\left\{
\operatorname{st} \tau   \ \mid\ \tau\in S_{\prec}\left(
\gamma,\sigma\right)  \right\}  _{\operatorname*{multi}}\\
&  =\left\{  \operatorname{st} \tau   \ \mid\ \tau\in S_{\prec
}\left(  \gamma,\sigma^{\prime}\right)  \right\}  _{\operatorname*{multi}}\\
&  =\left\{  \operatorname{st} \tau   \ \mid\ \tau\in S_{\prec
}\left(  \pi^{\prime},\sigma^{\prime}\right)  \right\}
_{\operatorname*{multi}}.
\end{align*}
The same argument (but with the symbols \textquotedblleft$S_{\prec}%
$\textquotedblright\ and \textquotedblleft$S_{\succ}$\textquotedblright%
\ interchanged) yields
\[
\left\{  \operatorname{st} \tau   \ \mid\ \tau\in S_{\succ
}\left(  \pi,\sigma\right)  \right\}  _{\operatorname*{multi}}=\left\{
\operatorname{st} \tau   \ \mid\ \tau\in S_{\succ}\left(
\pi^{\prime},\sigma^{\prime}\right)  \right\}  _{\operatorname*{multi}}.
\]
Thus, Claim 3 is proven in Case 1.

Let us now consider Case 2. In this case, we have $\left[  \pi_{1}>\sigma
_{1}\right]  =\left[  \pi_{1}^{\prime}>\sigma_{1}^{\prime}\right]  =0$.

There clearly exists a positive integer $N$ that is larger than all entries of
$\pi$ and larger than all entries of $\pi^{\prime}$. Consider such an $N$. Set
$m=\left\vert \sigma\right\vert $. Thus, $\sigma=\left(  \sigma_{1},\sigma
_{2},\ldots,\sigma_{m}\right)  $.
Let $\delta$ be the permutation $\left(  \sigma_{1}+N,\sigma_{2}%
+N,\ldots,\sigma_{m}+N\right)  $. This permutation $\delta$ is
order-isomorphic to $\sigma$, but is disjoint from $\pi$ (since all its
entries are $>N$, while all the entries of $\pi$ are $<N$) and disjoint from
$\pi^{\prime}$ (for similar reasons). Also, $\delta_{1}=\underbrace{\sigma
_{1}}_{>0}+N>N>\pi_{1}$ (since $N$ is larger than all entries of $\pi$),
so that we don't have $\pi_{1}>\delta_{1}$. Thus,
$\left[  \pi_{1}>\delta_{1}\right]  =0$. Similarly, $\left[  \pi_{1}^{\prime
}>\delta_{1}\right]  =0$.

The permutation $\delta$ is order-isomorphic to $\sigma$. Thus,
$\operatorname{st} \delta = \operatorname{st} \sigma$
(since $\operatorname{st}$ is a permutation statistic) and
$\left\vert \delta\right\vert =\left\vert \sigma\right\vert $. The permutation
$\delta$ is furthermore nonempty (since it is order-isomorphic to the nonempty
permutation $\sigma$). Also, $\operatorname{st} \delta
=\operatorname{st} \sigma   =\operatorname{st}\left(
\sigma^{\prime}\right)  $ and $\left\vert \delta\right\vert =\left\vert
\sigma\right\vert =\left\vert \sigma^{\prime}\right\vert $. Moreover, $\left[
\pi_{1}>\sigma_{1}\right]  =0=\left[  \pi_{1}>\delta_{1}\right]  $ and
$\left[  \pi_{1}>\delta_{1}\right]  =0=\left[  \pi_{1}^{\prime}>\delta
_{1}\right]  $ and $\left[  \pi_{1}^{\prime}>\delta_{1}\right]  =0=\left[
\pi_{1}^{\prime}>\sigma_{1}^{\prime}\right]  $. Hence, Claim 2 (applied to
$\delta$ instead of $\sigma^{\prime}$) yields%
\[
\left\{  \operatorname{st} \tau   \ \mid\ \tau\in S_{\prec
}\left(  \pi,\sigma\right)  \right\}  _{\operatorname*{multi}}=\left\{
\operatorname{st} \tau   \ \mid\ \tau\in S_{\prec}\left(
\pi,\delta\right)  \right\}  _{\operatorname*{multi}}%
\]
and%
\begin{equation}
\left\{  \operatorname{st} \tau   \ \mid\ \tau\in S_{\succ
}\left(  \pi,\sigma\right)  \right\}  _{\operatorname*{multi}}=\left\{
\operatorname{st} \tau   \ \mid\ \tau\in S_{\succ}\left(
\pi,\delta\right)  \right\}  _{\operatorname*{multi}}.\nonumber
\end{equation}
Furthermore, Claim 1 (applied to $\delta$ instead of $\sigma$) yields%
\[
\left\{  \operatorname{st} \tau   \ \mid\ \tau\in S_{\prec
}\left(  \pi,\delta\right)  \right\}  _{\operatorname*{multi}}=\left\{
\operatorname{st} \tau   \ \mid\ \tau\in S_{\prec}\left(
\pi^{\prime},\delta\right)  \right\}  _{\operatorname*{multi}}%
\]
and%
\[
\left\{  \operatorname{st} \tau   \ \mid\ \tau\in S_{\succ
}\left(  \pi,\delta\right)  \right\}  _{\operatorname*{multi}}=\left\{
\operatorname{st} \tau   \ \mid\ \tau\in S_{\succ}\left(
\pi^{\prime},\delta\right)  \right\}  _{\operatorname*{multi}}.
\]
Finally, Claim 2 (applied to $\pi^{\prime}$ and $\delta$ instead of $\pi$ and
$\sigma$) yields%
\[
\left\{  \operatorname{st} \tau   \ \mid\ \tau\in S_{\prec
}\left(  \pi^{\prime},\delta\right)  \right\}  _{\operatorname*{multi}%
}=\left\{  \operatorname{st} \tau   \ \mid\ \tau\in S_{\prec
}\left(  \pi^{\prime},\sigma^{\prime}\right)  \right\}
_{\operatorname*{multi}}%
\]
and%
\begin{equation}
\left\{  \operatorname{st} \tau   \ \mid\ \tau\in S_{\succ
}\left(  \pi^{\prime},\delta\right)  \right\}  _{\operatorname*{multi}%
}=\left\{  \operatorname{st} \tau   \ \mid\ \tau\in S_{\succ
}\left(  \pi^{\prime},\sigma^{\prime}\right)  \right\}
_{\operatorname*{multi}}.\nonumber
\end{equation}
Combining the equalities we have found, we obtain%
\begin{align*}
\left\{  \operatorname{st} \tau   \ \mid\ \tau\in S_{\prec
}\left(  \pi,\sigma\right)  \right\}  _{\operatorname*{multi}}  &  =\left\{
\operatorname{st} \tau   \ \mid\ \tau\in S_{\prec}\left(
\pi,\delta\right)  \right\}  _{\operatorname*{multi}}\\
&  =\left\{  \operatorname{st} \tau   \ \mid\ \tau\in S_{\prec
}\left(  \pi^{\prime},\delta\right)  \right\}  _{\operatorname*{multi}}\\
&  =\left\{  \operatorname{st} \tau   \ \mid\ \tau\in S_{\prec
}\left(  \pi^{\prime},\sigma^{\prime}\right)  \right\}
_{\operatorname*{multi}}.
\end{align*}
The same argument (but with the symbols \textquotedblleft$S_{\prec}%
$\textquotedblright\ and \textquotedblleft$S_{\succ}$\textquotedblright%
\ interchanged) yields
\[
\left\{  \operatorname{st} \tau   \ \mid\ \tau\in S_{\succ
}\left(  \pi,\sigma\right)  \right\}  _{\operatorname*{multi}}=\left\{
\operatorname{st} \tau   \ \mid\ \tau\in S_{\succ}\left(
\pi^{\prime},\sigma^{\prime}\right)  \right\}  _{\operatorname*{multi}}.
\]
Thus, Claim 3 is proven in Case 2.

We have now proven Claim 3 in each of the two Cases 1 and 2. Hence, Claim 3
always holds.]

Claim 3 says that for any two disjoint nonempty permutations $\pi$ and
$\sigma$, the multisets%
\[
\left\{  \operatorname{st} \tau   \ \mid\ \tau\in S_{\prec
}\left(  \pi,\sigma\right)  \right\}  _{\operatorname*{multi}}%
\ \ \ \ \ \ \ \ \ \ \text{and}\ \ \ \ \ \ \ \ \ \ \left\{  \operatorname{st}
\tau \ \mid\ \tau\in S_{\succ}\left(  \pi,\sigma\right)
\right\}  _{\operatorname*{multi}}
\]
depend only on $\operatorname{st} \pi   $, $\operatorname{st} \sigma$,
$\left\vert \pi\right\vert $, $\left\vert
\sigma\right\vert $ and $\left[  \pi_{1}>\sigma_{1}\right]  $. In other words,
the statistic $\operatorname{st}$ is LR-shuffle-compatible (by the definition
of \textquotedblleft LR-shuffle-compatible\textquotedblright). This proves
Theorem \ref{thm.head-comp.LRcomp}.
\end{proof}

Combining Theorem \ref{thm.head-comp.LRcomp} with Proposition
\ref{prop.head-comp.Pks}, we obtain the following:

\begin{theorem}
\label{thm.LRcomp.Pks}\textbf{(a)} The permutation statistic
$\operatorname{Des}$ is LR-shuffle-compatible.

\textbf{(b)} The permutation statistic $\operatorname*{Lpk}$ is LR-shuffle-compatible.

\textbf{(c)} The permutation statistic $\operatorname{Epk}$ is LR-shuffle-compatible.
\end{theorem}

\begin{proof}
[Proof of Theorem \ref{thm.LRcomp.Pks}.]\textbf{(a)} The permutation statistic
$\operatorname{Des}$ is shuffle-compatible (by \cite[\S 2.4]{part1}) and
head-graft-compatible (by Proposition \ref{prop.head-comp.Pks} \textbf{(a)}).
Thus, Theorem \ref{thm.head-comp.LRcomp} (applied to $\operatorname{st}
=\operatorname{Des}$) shows that the permutation statistic
$\operatorname{Des}$ is LR-shuffle-compatible. This proves Theorem
\ref{thm.LRcomp.Pks} \textbf{(a)}.

\textbf{(b)} The permutation statistic $\operatorname*{Lpk}$ is
shuffle-compatible (by \cite[Theorem 4.9 \textbf{(a)}]{part1}) and
head-graft-compatible (by Proposition \ref{prop.head-comp.Pks} \textbf{(b)}).
Thus, Theorem \ref{thm.head-comp.LRcomp} (applied to $\operatorname{st}
=\operatorname{Lpk}$) shows that the permutation statistic
$\operatorname{Lpk}$ is LR-shuffle-compatible. This proves Theorem
\ref{thm.LRcomp.Pks} \textbf{(b)}.

\textbf{(c)} The permutation statistic $\operatorname{Epk}$ is
shuffle-compatible (by Theorem \ref{thm.Epk.sh-co-a}) and
head-graft-compatible (by Proposition \ref{prop.head-comp.Pks} \textbf{(c)}).
Thus, Theorem \ref{thm.head-comp.LRcomp} (applied to $\operatorname{st}
=\operatorname{Epk}$) shows that the permutation statistic
$\operatorname{Epk}$ is LR-shuffle-compatible. This proves Theorem
\ref{thm.LRcomp.Pks} \textbf{(c)}.
\end{proof}

\subsection{\label{subsect.LR.others}Some other statistics}

The question of LR-shuffle-compatibility can be asked about any statistic. We
have so far answered it for $\operatorname{Des}$, $\operatorname*{Pk}$,
$\operatorname*{Lpk}$, $\operatorname*{Rpk}$ and $\operatorname{Epk}$. In
this section, we shall analyze it for some further statistics.

\subsubsection{The descent number $\operatorname{des}$}

The permutation statistic $\operatorname{des}$ (called the \textit{descent
number}) is defined as follows: For each permutation $\pi$, we set
$\operatorname{des}\pi=\left\vert \operatorname{Des}\pi\right\vert $ (that
is, $\operatorname{des}\pi$ is the number of all descents of $\pi$). It was
proven in \cite[Theorem 4.6 \textbf{(a)}]{part1} that this statistic
$\operatorname{des}$ is shuffle-compatible. We now claim the following:

\begin{proposition}
\label{prop.LRcomp.des}The permutation statistic $\operatorname{des}$ is
head-graft-compatible and LR-shuffle-compatible.
\end{proposition}

\begin{proof}
[Proof of Proposition \ref{prop.LRcomp.des}.]From
(\ref{pf.prop.head-comp.Pks.a.1}), we easily obtain the following: If $\pi$ is
a nonempty permutation, and if $a$ is a letter that does not appear in $\pi$,
then%
\[
\operatorname{des}\left(  a:\pi\right)  =\operatorname{des}\pi+\left[
a>\pi_{1}\right]  .
\]
Thus, $\operatorname{des}\left(  a:\pi\right)  $ depends only on
$\operatorname{des}\pi$, $\left\vert \pi\right\vert $ and $\left[  a>\pi
_{1}\right]  $. In other words, $\operatorname{des}$ is head-graft-compatible
(by the definition of \textquotedblleft
head-graft-compatible\textquotedblright). Hence,
Theorem \ref{thm.head-comp.LRcomp}
(applied to $\operatorname{st}=\operatorname{des}$) shows that
the permutation statistic $\operatorname{des}$ is
LR-shuffle-compatible. This proves Proposition \ref{prop.LRcomp.des}.
\end{proof}

\subsubsection{The major index $\operatorname*{maj}$}

The permutation statistic $\operatorname*{maj}$ (called the \textit{major
index}) is defined as follows: For each permutation $\pi$, we set
$\operatorname*{maj}\pi=\sum_{i\in\operatorname{Des}\pi}i$ (that is,
$\operatorname*{maj}\pi$ is the sum of all descents of $\pi$). It was proven
in \cite[Theorem 3.1 \textbf{(a)}]{part1} that this statistic
$\operatorname*{maj}$ is shuffle-compatible.

\begin{vershort}
However, $\operatorname*{maj}$ is neither head-graft-compatible nor
LR-shuffle-compatible. For example, if we take $\pi=\left(  5,4,2,3\right)  $,
$a=1$, $\pi^{\prime}=\left(  3,4,5,2\right)  $ and $a^{\prime}=1$, then we do
have%
\[
\operatorname*{maj} \pi =\operatorname*{maj}\left(  \pi
^{\prime}\right)  ,\ \ \ \ \ \ \ \ \ \ \left\vert \pi\right\vert =\left\vert
\pi^{\prime}\right\vert \ \ \ \ \ \ \ \ \ \ \text{and}%
\ \ \ \ \ \ \ \ \ \ \left[  a>\pi_{1}\right]  =\left[  a^{\prime}>\pi
_{1}^{\prime}\right]
\]
but we don't have $\operatorname*{maj}\left(  a:\pi\right)
=\operatorname*{maj}\left(  a^{\prime}:\pi^{\prime}\right)  $. Thus,
$\operatorname*{maj}$ is not head-graft-compatible. Using Corollary
\ref{cor.LRcomp.back} below, this entails that
$\operatorname*{maj}$ is not LR-shuffle-compatible.
\end{vershort}

\begin{verlong}
However, $\operatorname*{maj}$ is neither head-graft-compatible nor
LR-shuffle-compatible. For example, if we take $\pi=\left(  5,4,2,3\right)  $,
$a=1$, $\pi^{\prime}=\left(  3,4,5,2\right)  $ and $a^{\prime}=1$, then we do
have%
\[
\operatorname*{maj} \pi = \operatorname*{maj}\left(  \pi
^{\prime}\right)  ,\ \ \ \ \ \ \ \ \ \ \left\vert \pi\right\vert =\left\vert
\pi^{\prime}\right\vert \ \ \ \ \ \ \ \ \ \ \text{and}%
\ \ \ \ \ \ \ \ \ \ \left[  a>\pi_{1}\right]  =\left[  a^{\prime}>\pi
_{1}^{\prime}\right]
\]
but we don't have $\operatorname*{maj}\left(  a:\pi\right)
=\operatorname*{maj}\left(  a^{\prime}:\pi^{\prime}\right)  $. Thus,
$\operatorname*{maj}$ is not head-graft-compatible. Using Proposition
\ref{prop.LRcomp.head} below, this entails that $\operatorname*{maj}$ is not LR-shuffle-compatible.
\end{verlong}

\subsubsection{The joint statistic $\left(  \operatorname{des}
,\operatorname*{maj}\right)  $}

The next permutation statistic we shall study is the so-called joint statistic
$\left(  \operatorname{des},\operatorname*{maj}\right)  $. This statistic is
defined as the permutation statistic that sends each permutation $\pi$ to the
ordered pair $\left(  \operatorname{des}\pi,\operatorname*{maj}\pi\right)  $.
(Calling it $\left(  \operatorname{des},\operatorname*{maj}\right)  $ is thus
a slight abuse of notation.) It was proven in \cite[Theorem 4.5 \textbf{(a)}%
]{part1} that this statistic $\left(  \operatorname{des},\operatorname{maj}%
\right)  $ is shuffle-compatible. We now claim the following:

\begin{proposition}
\label{prop.LRcomp.desmaj}The permutation statistic $\left(
\operatorname{des},\operatorname*{maj}\right)  $ is head-graft-compatible and LR-shuffle-compatible.
\end{proposition}

\begin{proof}
[Proof of Proposition \ref{prop.LRcomp.desmaj}.]From
(\ref{pf.prop.head-comp.Pks.a.1}), we easily obtain the following: If $\pi$ is
a nonempty permutation, and if $a$ is a letter that does not appear in $\pi$,
then%
\begin{align*}
\operatorname{des}\left(  a:\pi\right)
&  =\operatorname{des}\pi+\left[ a>\pi_{1}\right]
\ \ \ \ \ \ \ \ \ \ \text{and}\\
\operatorname{maj}\left(  a:\pi\right)
&  =\operatorname{des} \pi + \operatorname{maj}\pi +\left[  a>\pi_{1}\right]  .
\end{align*}
Thus, $\left(  \operatorname{des},\operatorname*{maj}\right)  \left(
a:\pi\right)  $ depends only on $\left(  \operatorname{des}
,\operatorname*{maj}\right)  \left(  \pi\right)  $, $\left\vert \pi\right\vert
$ and $\left[  a>\pi_{1}\right]  $. In other words, $\left(
\operatorname{des},\operatorname*{maj}\right)  $ is head-graft-compatible (by
the definition of \textquotedblleft head-graft-compatible\textquotedblright).
Hence, Theorem \ref{thm.head-comp.LRcomp} (applied to $\operatorname{st}
=\left(  \operatorname{des},\operatorname*{maj}\right)  $) shows that the
permutation statistic $\left(  \operatorname{des},\operatorname*{maj}\right)
$ is LR-shuffle-compatible. This proves Proposition \ref{prop.LRcomp.desmaj}.
\end{proof}

\subsubsection{The comajor index $\operatorname*{comaj}$}

The permutation statistic $\operatorname*{comaj}$ (called the \textit{comajor
index}) is defined as follows: For each permutation $\pi$, we set
$\operatorname*{comaj}\pi=\sum_{k\in\operatorname{Des}\pi}\left(  n-k\right)
$, where $n=\left\vert \pi\right\vert $. It was proven in \cite[\S 3.2]{part1}
that this statistic $\operatorname*{comaj}$ is shuffle-compatible. We now
claim the following:

\begin{proposition}
\label{prop.LRcomp.comaj}The permutation statistic $\operatorname*{comaj}$ is
head-graft-compatible and LR-shuffle-compatible.
\end{proposition}

\begin{proof}
[Proof of Proposition \ref{prop.LRcomp.comaj}.]From
(\ref{pf.prop.head-comp.Pks.a.1}), we easily obtain the following: If $\pi$ is
a nonempty permutation, and if $a$ is a letter that does not appear in $\pi$,
then%
\[
\operatorname*{comaj}\left(  a:\pi\right)  =\operatorname*{comaj}\pi+\left[
a>\pi_{1}\right]  \cdot\left\vert \pi\right\vert .
\]
Thus, $\operatorname*{comaj}\left(  a:\pi\right)  $ depends only on
$\operatorname*{comaj}\pi$, $\left\vert \pi\right\vert $ and $\left[
a>\pi_{1}\right]  $. In other words, $\operatorname*{comaj}$ is
head-graft-compatible (by the definition of \textquotedblleft
head-graft-compatible\textquotedblright). Hence, Theorem
\ref{thm.head-comp.LRcomp} (applied to $\operatorname{st}
=\operatorname*{comaj}$) shows that the permutation statistic
$\operatorname*{comaj}$ is LR-shuffle-compatible. This proves Proposition
\ref{prop.LRcomp.comaj}.
\end{proof}

\subsection{Left- and right-shuffle-compatibility}

In this section, we shall study two notions closely related to LR-shuffle-compatibility:

\begin{definition}
\label{def.LR.left-right}Let $\operatorname{st}$ be a permutation statistic.

\textbf{(a)} We say that $\operatorname{st}$ is
\textit{left-shuffle-compatible} if for any two disjoint nonempty permutations
$\pi$ and $\sigma$ having the property that $\pi_{1}>\sigma_{1}$, the multiset
$\left\{  \operatorname{st} \tau   \ \mid\ \tau\in S_{\prec
}\left(  \pi,\sigma\right)  \right\}  _{\operatorname*{multi}}$ depends only
on $\operatorname{st} \pi$, $\operatorname{st} \sigma$,
$\left\vert \pi\right\vert $ and $\left\vert \sigma \right\vert $.

\textbf{(b)} We say that $\operatorname{st}$ is
\textit{right-shuffle-compatible} if for any two disjoint nonempty
permutations $\pi$ and $\sigma$ having the property that $\pi_{1}>\sigma_{1}$,
the multiset $\left\{  \operatorname{st} \tau   \ \mid\ \tau\in
S_{\succ}\left(  \pi,\sigma\right)  \right\}  _{\operatorname*{multi}}$
depends only on $\operatorname{st} \pi$, $\operatorname{st} \sigma$,
$\left\vert \pi\right\vert $ and $\left\vert \sigma\right\vert $.
\end{definition}

For a shuffle-compatible permutation statistic, these two notions are
equivalent to the notions of LR-shuffle-compatibility and
head-graft-compatibility, as the following proposition reveals:

\begin{proposition}
\label{prop.LRcomp.equivs}Let $\operatorname{st}$ be a shuffle-compatible
permutation statistic. Then, the following assertions are equivalent:

\begin{itemize}
\item \textit{Assertion }$\mathcal{A}_{1}$\textit{:} The statistic
$\operatorname{st}$ is LR-shuffle-compatible.

\item \textit{Assertion }$\mathcal{A}_{2}$\textit{:} The statistic
$\operatorname{st}$ is left-shuffle-compatible.

\item \textit{Assertion }$\mathcal{A}_{3}$\textit{:} The statistic
$\operatorname{st}$ is right-shuffle-compatible.

\item \textit{Assertion }$\mathcal{A}_{4}$\textit{:} The statistic
$\operatorname{st}$ is head-graft-compatible.
\end{itemize}
\end{proposition}

\begin{vershort}
\begin{proof}
[Proof of Proposition \ref{prop.LRcomp.equivs}.]Omitted; see \cite{verlong}.
\end{proof}
\end{vershort}

\begin{verlong}
\begin{proof}
[Proof of Proposition \ref{prop.LRcomp.equivs} (sketched).]We shall prove the
implications $\mathcal{A}_{1}\Longrightarrow\mathcal{A}_{2}$, $\mathcal{A}%
_{2}\Longrightarrow\mathcal{A}_{3}$, $\mathcal{A}_{3}\Longrightarrow
\mathcal{A}_{2}$, $\mathcal{A}_{3}\Longrightarrow\mathcal{A}_{4}$ and
$\mathcal{A}_{4}\Longrightarrow\mathcal{A}_{1}$.

The implication $\mathcal{A}_{4}\Longrightarrow\mathcal{A}_{1}$ follows from
Theorem \ref{thm.head-comp.LRcomp}.

\begin{noncompile}
Let us first recall that $\operatorname{st}$ is shuffle-compatible. In other
words, the following holds:

\begin{statement}
\textit{Claim 0:} Let $\pi$ and $\sigma$ be two disjoint permutations. Let
$\pi^{\prime}$ and $\sigma^{\prime}$ be two disjoint permutations. Assume that%
\begin{align*}
\operatorname{st} \pi    &  =\operatorname{st}\left(
\pi^{\prime}\right)  ,\ \ \ \ \ \ \ \ \ \ \operatorname{st} \sigma
=\operatorname{st}\left(  \sigma^{\prime}\right)  ,\\
\left\vert \pi\right\vert  &  =\left\vert \pi^{\prime}\right\vert
\ \ \ \ \ \ \ \ \ \ \text{and}\ \ \ \ \ \ \ \ \ \ \left\vert \sigma\right\vert
=\left\vert \sigma^{\prime}\right\vert .
\end{align*}
Then,
\[
\left\{  \operatorname{st} \tau   \ \mid\ \tau\in S\left(
\pi,\sigma\right)  \right\}  _{\operatorname*{multi}}=\left\{
\operatorname{st} \tau   \ \mid\ \tau\in S\left(  \pi^{\prime
},\sigma^{\prime}\right)  \right\}  _{\operatorname*{multi}}.
\]

\end{statement}
\end{noncompile}

\textit{Proof of the implication }$\mathcal{A}_{1}\Longrightarrow
\mathcal{A}_{2}$\textit{:} Assume that Assertion $\mathcal{A}_{1}$ holds. In
other words, the statistic $\operatorname{st}$ is LR-shuffle-compatible. We
shall show that Assertion $\mathcal{A}_{2}$ holds.

The statistic $\operatorname{st}$ is LR-shuffle-compatible. In other words,
for any two disjoint nonempty permutations $\pi$ and $\sigma$, the multisets%
\[
\left\{  \operatorname{st} \tau   \ \mid\ \tau\in S_{\prec
}\left(  \pi,\sigma\right)  \right\}  _{\operatorname*{multi}}%
\ \ \ \ \ \ \ \ \ \ \text{and}
\ \ \ \ \ \ \ \ \ \ \left\{  \operatorname{st} \tau
\ \mid\ \tau\in S_{\succ}\left(  \pi,\sigma\right)
\right\}  _{\operatorname*{multi}}%
\]
depend only on $\operatorname{st} \pi   $, $\operatorname{st} \sigma$,
$\left\vert \pi\right\vert $, $\left\vert \sigma\right\vert $ and
$\left[  \pi_{1}>\sigma_{1}\right]  $. Hence, for any
two disjoint nonempty permutations $\pi$ and $\sigma$, the multiset $\left\{
\operatorname{st} \tau   \ \mid\ \tau\in S_{\prec}\left(
\pi,\sigma\right)  \right\}  _{\operatorname*{multi}}$ depends only on
$\operatorname{st} \pi$, $\operatorname{st} \sigma$,
$\left\vert \pi\right\vert $, $\left\vert \sigma\right\vert$ and
$\left[  \pi_{1}>\sigma_{1}\right]  $. Hence, for any two disjoint
nonempty permutations $\pi$ and $\sigma$ having the property that $\pi
_{1}>\sigma_{1}$, the multiset \newline$\left\{  \operatorname{st}
\tau \ \mid\ \tau\in S_{\prec}\left(  \pi,\sigma\right)  \right\}
_{\operatorname*{multi}}$ depends only on $\operatorname{st} \pi$,
$\operatorname{st} \sigma$, $\left\vert \pi\right\vert $ and
$\left\vert \sigma\right\vert $ (indeed, it no longer
depends on $\left[  \pi_{1}>\sigma_{1}\right]  $, because our condition
$\pi_{1}>\sigma_{1}$ ensures that $\left[  \pi_{1}>\sigma_{1}\right]  =1$). In
other words, the statistic $\operatorname{st}$ is left-shuffle-compatible (by
the definition of \textquotedblleft left-shuffle-compatible\textquotedblright%
). In other words, Assertion $\mathcal{A}_{2}$ holds. This proves the
implication $\mathcal{A}_{1}\Longrightarrow\mathcal{A}_{2}$.

\textit{Proof of the implication }$\mathcal{A}_{2}\Longrightarrow
\mathcal{A}_{3}$\textit{:} Assume that Assertion $\mathcal{A}_{2}$ holds. In
other words, the statistic $\operatorname{st}$ is left-shuffle-compatible. We
shall show that Assertion $\mathcal{A}_{3}$ holds.

If $\pi$ and $\sigma$ are two disjoint nonempty permutations, then%
\begin{align*}
&  \left\{  \operatorname{st} \tau   \ \mid\ \tau\in S_{\succ
}\left(  \pi,\sigma\right)  \right\}  _{\operatorname*{multi}}\\
&  =\left\{  \operatorname{st} \tau   \ \mid\ \tau\in S\left(
\pi,\sigma\right)  \right\}  _{\operatorname*{multi}}-\left\{
\operatorname{st} \tau   \ \mid\ \tau\in S_{\prec}\left(
\pi,\sigma\right)  \right\}  _{\operatorname*{multi}}%
\end{align*}
(by Lemma \ref{lem.LR.difference} \textbf{(b)}). Hence, if $\pi$ and $\sigma$
are two disjoint nonempty permutations, then the multiset $\left\{
\operatorname{st} \tau   \ \mid\ \tau\in S_{\succ}\left(
\pi,\sigma\right)  \right\}  _{\operatorname*{multi}}$ is uniquely determined
by $\left\{  \operatorname{st} \tau   \ \mid\ \tau\in S\left(
\pi,\sigma\right)  \right\}  _{\operatorname*{multi}}$ and $\left\{
\operatorname{st} \tau   \ \mid\ \tau\in S_{\prec}\left(
\pi,\sigma\right)  \right\}  _{\operatorname*{multi}}$.

Thus, altogether, we know that for any two disjoint nonempty permutations
$\pi$ and $\sigma$ having the property that $\pi_{1}>\sigma_{1}$, the
following holds:

\begin{itemize}
\item The multiset $\left\{  \operatorname{st} \tau 
\ \mid\ \tau\in S_{\prec}\left(  \pi,\sigma\right)  \right\}
_{\operatorname*{multi}}$ depends only on $\operatorname{st} \pi$,
$\operatorname{st} \sigma$,
$\left\vert \pi\right\vert $ and $\left\vert \sigma\right\vert $
(since the statistic $\operatorname{st}$ is left-shuffle-compatible);

\item The multiset $\left\{  \operatorname{st} \tau 
\ \mid\ \tau\in S\left(  \pi,\sigma\right)  \right\}  _{\operatorname*{multi}%
}$ depends only on $\operatorname{st} \pi   $,
$\operatorname{st} \sigma   $, $\left\vert \pi\right\vert $ and
$\left\vert \sigma\right\vert $ (since the statistic $\operatorname{st}$ is shuffle-compatible);

\item The multiset $\left\{  \operatorname{st} \tau 
\ \mid\ \tau\in S_{\succ}\left(  \pi,\sigma\right)  \right\}
_{\operatorname*{multi}}$ is uniquely determined by \newline$\left\{
\operatorname{st} \tau   \ \mid\ \tau\in S\left(  \pi
,\sigma\right)  \right\}  _{\operatorname*{multi}}$ and $\left\{
\operatorname{st} \tau   \ \mid\ \tau\in S_{\prec}\left(
\pi,\sigma\right)  \right\}  _{\operatorname*{multi}}$.
\end{itemize}

Hence, the multiset $\left\{  \operatorname{st} \tau 
\ \mid\ \tau\in S_{\succ}\left(  \pi,\sigma\right)  \right\}
_{\operatorname*{multi}}$ also depends only on $\operatorname{st} \pi$,
$\operatorname{st} \sigma$,
$\left\vert \pi\right\vert $ and $\left\vert \sigma\right\vert $.
In other words, the
statistic $\operatorname{st}$ is right-shuffle-compatible (by the definition
of \textquotedblleft right-shuffle-compatible\textquotedblright). In other
words, Assertion $\mathcal{A}_{3}$ holds. This proves the implication
$\mathcal{A}_{2}\Longrightarrow\mathcal{A}_{3}$.

\textit{Proof of the implication }$\mathcal{A}_{3}\Longrightarrow
\mathcal{A}_{2}$\textit{:} The proof of the implication $\mathcal{A}%
_{3}\Longrightarrow\mathcal{A}_{2}$ is entirely analogous to the above proof
of the implication $\mathcal{A}_{2}\Longrightarrow\mathcal{A}_{3}$ (except
that we need to interchange \textquotedblleft left\textquotedblright\ and
\textquotedblleft right\textquotedblright\ and likewise interchange
\textquotedblleft$S_{\prec}$\textquotedblright\ and \textquotedblleft%
$S_{\succ}$\textquotedblright).

\textit{Proof of the implication }$\mathcal{A}_{3}\Longrightarrow
\mathcal{A}_{4}$\textit{:} Assume that Assertion $\mathcal{A}_{3}$ holds. We
shall show that Assertion $\mathcal{A}_{4}$ holds.

We have already proven the implication $\mathcal{A}_{3}\Longrightarrow
\mathcal{A}_{2}$. Thus, Assertion $\mathcal{A}_{2}$ holds (since
$\mathcal{A}_{3}$ holds). In other words, the statistic $\operatorname{st}$
is left-shuffle-compatible. In other words, the following claim holds:

\begin{statement}
\textit{Claim 1:} Let $\pi$ and $\sigma$ be two disjoint nonempty permutations
having the property that $\pi_{1}>\sigma_{1}$. Let $\pi^{\prime}$ and
$\sigma^{\prime}$ be two disjoint nonempty permutations having the property
that $\pi_{1}^{\prime}>\sigma_{1}^{\prime}$. Assume that%
\begin{align*}
\operatorname{st} \pi    &  =\operatorname{st}\left(
\pi^{\prime}\right)  ,\ \ \ \ \ \ \ \ \ \ \operatorname{st} \sigma
=\operatorname{st}\left(  \sigma^{\prime}\right)  ,\\
\left\vert \pi\right\vert  &  =\left\vert \pi^{\prime}\right\vert
\ \ \ \ \ \ \ \ \ \ \text{and}\ \ \ \ \ \ \ \ \ \ \left\vert \sigma\right\vert
=\left\vert \sigma^{\prime}\right\vert .
\end{align*}
Then,
\[
\left\{  \operatorname{st} \tau   \ \mid\ \tau\in S_{\prec
}\left(  \pi,\sigma\right)  \right\}  _{\operatorname*{multi}}=\left\{
\operatorname{st} \tau   \ \mid\ \tau\in S_{\prec}\left(
\pi^{\prime},\sigma^{\prime}\right)  \right\}  _{\operatorname*{multi}}.
\]

\end{statement}

Also, Assertion $\mathcal{A}_{3}$ holds. In other words, the statistic
$\operatorname{st}$ is right-shuffle-compatible. In other words, the
following claim holds:

\begin{statement}
\textit{Claim 2:} Let $\pi$ and $\sigma$ be two disjoint nonempty permutations
having the property that $\pi_{1}>\sigma_{1}$. Let $\pi^{\prime}$ and
$\sigma^{\prime}$ be two disjoint nonempty permutations having the property
that $\pi_{1}^{\prime}>\sigma_{1}^{\prime}$. Assume that%
\begin{align*}
\operatorname{st} \pi    &  =\operatorname{st}\left(
\pi^{\prime}\right)  ,\ \ \ \ \ \ \ \ \ \ \operatorname{st} \sigma
=\operatorname{st}\left(  \sigma^{\prime}\right)  ,\\
\left\vert \pi\right\vert  &  =\left\vert \pi^{\prime}\right\vert
\ \ \ \ \ \ \ \ \ \ \text{and}\ \ \ \ \ \ \ \ \ \ \left\vert \sigma\right\vert
=\left\vert \sigma^{\prime}\right\vert .
\end{align*}
Then,
\[
\left\{  \operatorname{st} \tau   \ \mid\ \tau\in S_{\succ
}\left(  \pi,\sigma\right)  \right\}  _{\operatorname*{multi}}=\left\{
\operatorname{st} \tau   \ \mid\ \tau\in S_{\succ}\left(
\pi^{\prime},\sigma^{\prime}\right)  \right\}  _{\operatorname*{multi}}.
\]

\end{statement}

We are now going to prove the following claim:

\begin{statement}
\textit{Claim 3:} Let $\pi$ be a nonempty permutation, and let $a$ be a letter
that does not appear in $\pi$. Let $\pi^{\prime}$ be a nonempty permutation,
and let $a^{\prime}$ be a letter that does not appear in $\pi^{\prime}$.
Assume that%
\[
\operatorname{st} \pi   =\operatorname{st}\left(  \pi^{\prime
}\right)  ,\ \ \ \ \ \ \ \ \ \ \left\vert \pi\right\vert =\left\vert
\pi^{\prime}\right\vert \ \ \ \ \ \ \ \ \ \ \text{and}%
\ \ \ \ \ \ \ \ \ \ \left[  a>\pi_{1}\right]  =\left[  a^{\prime}>\pi
_{1}^{\prime}\right]  .
\]
Then, $\operatorname{st}\left(  a:\pi\right)  =\operatorname{st}\left(
a^{\prime}:\pi^{\prime}\right)  $.
\end{statement}

[\textit{Proof of Claim 3:} The $1$-permutations $\left(  a\right)  $ and
$\left(  a^{\prime}\right)  $ are order-isomorphic. Thus, $\operatorname{st}
\left(  \left(  a\right)  \right)  =\operatorname{st}\left(  \left(
a^{\prime}\right)  \right)  $ (since $\operatorname{st}$ is a permutation
statistic). Also, $\left\vert \left(  a\right)  \right\vert =1=\left\vert
\left(  a^{\prime}\right)  \right\vert $.

We are in one of the following two cases:

\textit{Case 1:} We have $a>\pi_{1}$.

\textit{Case 2:} We have $a\leq\pi_{1}$.

Let us first consider Case 1. In this case, we have $a>\pi_{1}$. Thus,
$\left[  a>\pi_{1}\right]  =1$. Comparing this with $\left[  a>\pi_{1}\right]
=\left[  a^{\prime}>\pi_{1}^{\prime}\right]  $, we obtain $\left[  a^{\prime
}>\pi_{1}^{\prime}\right]  =1$, so that $a^{\prime}>\pi_{1}^{\prime}$.

Consider the $1$-permutations $\left(  a\right)  $ and $\left(  a^{\prime
}\right)  $. Then, the first entry of $\left(  a\right)  $ is $\left(
a\right)  _{1}=a>\pi_{1}$. Also, $\left(  a^{\prime}\right)  _{1}=a^{\prime
}>\pi_{1}^{\prime}$. Also, the permutations $\left(  a\right)  $ and $\pi$ are
disjoint (since $a$ does not appear in $\pi$). Similarly, the permutations
$\left(  a^{\prime}\right)  $ and $\pi^{\prime}$ are disjoint. Furthermore,
$\operatorname{st}\left(  \left(  a\right)  \right)  =\operatorname{st}
\left(  \left(  a^{\prime}\right)  \right)  $, $\operatorname{st}
\pi  =\operatorname{st}\left(  \pi^{\prime}\right)  $, $\left\vert
\left(  a\right)  \right\vert =\left\vert \left(  a^{\prime}\right)
\right\vert $ and $\left\vert \pi\right\vert =\left\vert \pi^{\prime
}\right\vert $. Hence, Claim 1 (applied to $\left(  a\right)  $, $\pi$,
$\left(  a^{\prime}\right)  $ and $\pi^{\prime}$ instead of $\pi$, $\sigma$,
$\pi^{\prime}$ and $\sigma^{\prime}$) yields%
\begin{equation}
\left\{  \operatorname{st} \tau   \ \mid\ \tau\in S_{\prec
}\left(  \left(  a\right)  ,\pi\right)  \right\}  _{\operatorname*{multi}%
}=\left\{  \operatorname{st} \tau   \ \mid\ \tau\in S_{\prec
}\left(  \left(  a^{\prime}\right)  ,\pi^{\prime}\right)  \right\}
_{\operatorname*{multi}}. \label{pf.prop.LRcomp.equivs.3to4.4}%
\end{equation}
But $S_{\prec}\left(  \left(  a\right)  ,\pi\right)  =\left\{  a:\pi\right\}
$. Hence,%
\[
\left\{  \operatorname{st} \tau   \ \mid\ \tau\in S_{\prec
}\left(  \left(  a\right)  ,\pi\right)  \right\}  _{\operatorname*{multi}%
}=\left\{  \operatorname{st} \tau   \ \mid\ \tau\in\left\{
a:\pi\right\}  \right\}  _{\operatorname*{multi}}=\left\{  \operatorname{st}
\left(  a:\pi\right)  \right\}  _{\operatorname*{multi}}.
\]
Similarly,%
\[
\left\{  \operatorname{st} \tau   \ \mid\ \tau\in S_{\prec
}\left(  \left(  a^{\prime}\right)  ,\pi^{\prime}\right)  \right\}
_{\operatorname*{multi}}=\left\{  \operatorname{st}\left(  a^{\prime}%
:\pi^{\prime}\right)  \right\}  _{\operatorname*{multi}}.
\]
In light of the last two equalities, the equality
(\ref{pf.prop.LRcomp.equivs.3to4.4}) rewrites as $\left\{  \operatorname{st}
\left(  a:\pi\right)  \right\}  _{\operatorname*{multi}}=\left\{
\operatorname{st}\left(  a^{\prime}:\pi^{\prime}\right)  \right\}
_{\operatorname*{multi}}$. In other words, $\operatorname{st}\left(
a:\pi\right)  =\operatorname{st}\left(  a^{\prime}:\pi^{\prime}\right)  $.
Thus, we have proven $\operatorname{st}\left(  a:\pi\right)
=\operatorname{st}\left(  a^{\prime}:\pi^{\prime}\right)  $ in Case 1.

Let us now consider Case 2. In this case, we have $a\leq\pi_{1}$. But $a$ does
not appear in $\pi$; therefore, $a\neq\pi_{1}$. Combined with $a\leq\pi_{1}$,
this yields $a<\pi_{1}$. In other words, $\pi_{1}>a$. Also, from $a\leq\pi
_{1}$, we obtain $\left[  a>\pi_{1}\right]  =0$. Comparing this with $\left[
a>\pi_{1}\right]  =\left[  a^{\prime}>\pi_{1}^{\prime}\right]  $, we obtain
$\left[  a^{\prime}>\pi_{1}^{\prime}\right]  =0$, so that $a^{\prime}\leq
\pi_{1}^{\prime}$. But $a^{\prime}$ does not appear in $\pi^{\prime}$; thus,
$a^{\prime}\neq\pi_{1}^{\prime}$. Combined with $a^{\prime}\leq\pi_{1}%
^{\prime}$, this yields $a^{\prime}<\pi_{1}^{\prime}$. In other words,
$\pi_{1}^{\prime}>a^{\prime}$.

Consider the $1$-permutations $\left(  a\right)  $ and $\left(  a^{\prime
}\right)  $. Then, the first entry of $\left(  a\right)  $ is $\left(
a\right)  _{1}=a$; similarly, $\left(  a^{\prime}\right)  _{1}=a^{\prime}$.
Now, $\pi_{1}>a=\left(  a\right)  _{1}$ and $\pi_{1}^{\prime}>a^{\prime
}=\left(  a^{\prime}\right)  _{1}$. Also, the permutations $\pi$ and $\left(
a\right)  $ are disjoint (since $a$ does not appear in $\pi$). Similarly, the
permutations $\pi^{\prime}$ and $\left(  a^{\prime}\right)  $ are disjoint.
Furthermore, $\operatorname{st}\left(  \left(  a\right)  \right)
=\operatorname{st}\left(  \left(  a^{\prime}\right)  \right)  $,
$\operatorname{st} \pi   =\operatorname{st}\left(  \pi^{\prime
}\right)  $, $\left\vert \left(  a\right)  \right\vert =\left\vert \left(
a^{\prime}\right)  \right\vert $ and $\left\vert \pi\right\vert =\left\vert
\pi^{\prime}\right\vert $. Hence, Claim 2 (applied to $\pi$, $\left(
a\right)  $, $\pi^{\prime}$ and $\left(  a^{\prime}\right)  $ instead of $\pi
$, $\sigma$, $\pi^{\prime}$ and $\sigma^{\prime}$) yields%
\begin{equation}
\left\{  \operatorname{st} \tau   \ \mid\ \tau\in S_{\succ
}\left(  \pi,\left(  a\right)  \right)  \right\}  _{\operatorname*{multi}%
}=\left\{  \operatorname{st} \tau   \ \mid\ \tau\in S_{\succ
}\left(  \pi^{\prime},\left(  a^{\prime}\right)  \right)  \right\}
_{\operatorname*{multi}}. \label{pf.prop.LRcomp.equivs.3to4.5}%
\end{equation}
But $S_{\succ}\left(  \pi,\left(  a\right)  \right)  =\left\{  a:\pi\right\}
$. Hence,%
\[
\left\{  \operatorname{st} \tau   \ \mid\ \tau\in S_{\succ
}\left(  \pi,\left(  a\right)  \right)  \right\}  _{\operatorname*{multi}%
}=\left\{  \operatorname{st} \tau   \ \mid\ \tau\in\left\{
a:\pi\right\}  \right\}  _{\operatorname*{multi}}=\left\{  \operatorname{st}
\left(  a:\pi\right)  \right\}  _{\operatorname*{multi}}.
\]
Similarly,%
\[
\left\{  \operatorname{st} \tau   \ \mid\ \tau\in S_{\succ
}\left(  \pi^{\prime},\left(  a^{\prime}\right)  \right)  \right\}
_{\operatorname*{multi}}=\left\{  \operatorname{st}\left(  a^{\prime}%
:\pi^{\prime}\right)  \right\}  _{\operatorname*{multi}}.
\]
In light of the last two equalities, the equality
(\ref{pf.prop.LRcomp.equivs.3to4.5}) rewrites as $\left\{  \operatorname{st}
\left(  a:\pi\right)  \right\}  _{\operatorname*{multi}}=\left\{
\operatorname{st}\left(  a^{\prime}:\pi^{\prime}\right)  \right\}
_{\operatorname*{multi}}$. In other words, $\operatorname{st}\left(
a:\pi\right)  =\operatorname{st}\left(  a^{\prime}:\pi^{\prime}\right)  $.
Thus, we have proven $\operatorname{st}\left(  a:\pi\right)
=\operatorname{st}\left(  a^{\prime}:\pi^{\prime}\right)  $ in Case 2.

We have now proven $\operatorname{st}\left(  a:\pi\right)
=\operatorname{st}\left(  a^{\prime}:\pi^{\prime}\right)  $ in both Cases 1
and 2. Therefore, $\operatorname{st}\left(  a:\pi\right)  =\operatorname{st}
\left(  a^{\prime}:\pi^{\prime}\right)  $ always holds. This proves Claim 3.]

Claim 3 shows that the statistic $\operatorname{st}$ is
head-graft-compatible. In other words, Assertion $\mathcal{A}_{4}$ holds. This
proves the implication $\mathcal{A}_{3}\Longrightarrow\mathcal{A}_{4}$.

We have now proven the implications $\mathcal{A}_{1}\Longrightarrow
\mathcal{A}_{2}$, $\mathcal{A}_{2}\Longrightarrow\mathcal{A}_{3}$,
$\mathcal{A}_{3}\Longrightarrow\mathcal{A}_{4}$ and $\mathcal{A}%
_{4}\Longrightarrow\mathcal{A}_{1}$. Thus, all four statements $\mathcal{A}%
_{1}$, $\mathcal{A}_{2}$, $\mathcal{A}_{3}$ and $\mathcal{A}_{4}$ are
equivalent. This proves Proposition \ref{prop.LRcomp.equivs}.
\end{proof}
\end{verlong}

Note that on their own, the properties of left-shuffle-compatibility and
right-shuffle-compatibility are not equivalent. For example, the permutation
statistic that sends each nonempty permutation $\pi$ to the truth value
$\left[  \pi_{1}>\pi_{i}\text{ for all }i>1\right]  $ (and the $0$-permutation
$\left(  {}\right)  $ to $0$) is right-shuffle-compatible (because in the
definition of right-shuffle-compatibility, all the $\operatorname{st}
\tau$ will be $0$), but not left-shuffle-compatible.

\subsection{Properties of compatible statistics}

\begin{vershort}
Let us state some more facts on compatibility properties. We refer to
\cite{verlong} for their proofs. We begin with a converse to Theorem
\ref{thm.head-comp.LRcomp}:
\end{vershort}

\begin{verlong}
The following converse to Theorem \ref{thm.head-comp.LRcomp} holds:

\begin{proposition}
\label{prop.LRcomp.head}Let $\operatorname{st}$ be a permutation statistic
that is LR-shuffle-compatible. Then, $\operatorname{st}$ is
head-graft-compatible and shuffle-compatible.
\end{proposition}

Before we can prove this, we need three lemmas:

\begin{lemma}
\label{lem.LRcomp.head-l1}Let $\operatorname{st}$ be a head-graft-compatible
permutation statistic. If $\alpha$ and $\beta$ are two permutations satisfying
$\operatorname{Des}\alpha=\operatorname{Des}\beta$ and $\left\vert
\alpha\right\vert =\left\vert \beta\right\vert $, then
$\operatorname{st} \alpha = \operatorname{st} \beta$.
\end{lemma}

\begin{proof}
[Proof of Lemma \ref{lem.LRcomp.head-l1} (sketched).]We must prove the
following claim:

\begin{statement}
\textit{Claim 1:} Let $\alpha$ and $\beta$ be two permutations satisfying
$\operatorname{Des}\alpha=\operatorname{Des}\beta$ and $\left\vert
\alpha\right\vert =\left\vert \beta\right\vert $. Then,
$\operatorname{st} \alpha = \operatorname{st} \beta$.
\end{statement}

[\textit{Proof of Claim 1:} We shall prove Claim 1 by induction on $\left\vert
\alpha\right\vert $:

\textit{Induction base:} If $\left\vert \alpha\right\vert =0$, then Claim 1 is
true (because if $\left\vert \alpha\right\vert =0$, then both $\alpha$ and
$\beta$ equal the $0$-permutation $\left(  {}\right)  $, and therefore
satisfy $\alpha=\beta$ and thus $\operatorname{st} \alpha
=\operatorname{st} \beta$). This completes the induction base.

\textit{Induction step:} Let $N$ be a positive integer. Assume (as the
induction hypothesis) that Claim 1 holds for $\left\vert \alpha\right\vert
=N-1$. We must now prove that Claim 1 holds for $\left\vert \alpha\right\vert
=N$.

Let $\alpha$ and $\beta$ be as in Claim 1, and assume that $\left\vert
\alpha\right\vert =N$. We must prove that
$\operatorname{st} \alpha = \operatorname{st} \beta$.

If $N=1$, then this holds\footnote{\textit{Proof.} Assume that $N=1$. Thus,
$\alpha$ and $\beta$ are $1$-permutations (since $\left\vert \alpha\right\vert
=N=1$ and $\left\vert \beta\right\vert =\left\vert \alpha\right\vert =N=1$)
and thus are order-isomorphic. Hence,
$\operatorname{st} \alpha =\operatorname{st} \beta$ (since
$\operatorname{st}$ is a permutation statistic).}. Hence, for the rest of
this proof, we WLOG assume that $N\neq1$. Hence, $N\geq2$ (since $N$ is a
positive integer).

The permutation $\alpha$ is nonempty (since $\left\vert \alpha\right\vert
=N>0$). Thus, $\alpha_{\sim1}$ and $\alpha_{1}$ are well-defined. Similarly,
$\beta_{\sim1}$ and $\beta_{1}$ are well-defined (since $\left\vert
\beta\right\vert =\left\vert \alpha\right\vert =N>0$). Clearly, the letter
$\alpha_{1}$ does not appear in $\alpha_{\sim1}$ (since the letters of the
permutation $\alpha$ are distinct). Similarly, the letter $\beta_{1}$ does not
appear in $\beta_{\sim1}$.

We have $\left\vert \alpha_{\sim1}\right\vert =\underbrace{\left\vert
\alpha\right\vert }_{=N}-1=N-1$ and similarly $\left\vert \beta_{\sim
1}\right\vert =N-1$. Thus, $\left\vert \alpha_{\sim1}\right\vert
=N-1=\left\vert \beta_{\sim1}\right\vert $. Also, $\left\vert \alpha_{\sim
1}\right\vert =N-1\geq1$ (since $N\geq2$); thus, the permutation $\alpha
_{\sim1}$ is nonempty. Similarly, $\beta_{\sim1}$ is nonempty.

It is easy to see that
\begin{align*}
\operatorname{Des}\left(  \alpha_{\sim1}\right)   &  =\left\{  i-1\ \mid
\ i\in\left(  \operatorname{Des}\alpha\right)  \setminus\left\{  1\right\}
\right\}  \ \ \ \ \ \ \ \ \ \ \text{and}\\
\operatorname{Des}\left(  \beta_{\sim1}\right)   &  =\left\{  i-1\ \mid
\ i\in\left(  \operatorname{Des}\beta\right)  \setminus\left\{  1\right\}
\right\}  .
\end{align*}
The right hand sides of these two equalities are equal (since
$\operatorname{Des}\alpha=\operatorname{Des}\beta$). Thus, their left hand
sides are equal as well. In other words, $\operatorname{Des}\left(
\alpha_{\sim1}\right)  =\operatorname{Des}\left(  \beta_{\sim1}\right)  $.

Moreover, $\left\vert \alpha_{\sim1}\right\vert =N-1$.
Hence, the induction hypothesis reveals that
we can apply Claim 1 to $\alpha_{\sim1}$ and $\beta_{\sim1}$ instead of
$\alpha$ and $\beta$. We thus obtain $\operatorname{st}\left(  \alpha_{\sim
1}\right)  =\operatorname{st}\left(  \beta_{\sim1}\right)  $.

Furthermore, $\left(  \alpha_{\sim1}\right)  _{1}=\alpha_{2}$, so that
\[
\left[  \alpha_{1}>\left(  \alpha_{\sim1}\right)  _{1}\right]  =\left[
\alpha_{1}>\alpha_{2}\right]  =\left[  1\in\operatorname{Des}\alpha\right]
.
\]
Similarly,%
\[
\left[  \beta_{1}>\left(  \beta_{\sim1}\right)  _{1}\right]  =\left[
1\in\operatorname{Des}\beta\right]  .
\]
The right hand sides of these two equalities are equal (since
$\operatorname{Des}\alpha=\operatorname{Des}\beta$). Thus, their left hand
sides are equal as well. In other words, $\left[  \alpha_{1}>\left(
\alpha_{\sim1}\right)  _{1}\right]  =\left[  \beta_{1}>\left(  \beta_{\sim
1}\right)  _{1}\right]  $.

But recall that $\operatorname{st}$ is head-graft-compatible. In other words,
every nonempty permutation $\pi$, every letter $a$ that does not appear in
$\pi$, every nonempty permutation $\pi^{\prime}$ and every letter $a^{\prime}$
that does not appear in $\pi^{\prime}$ satisfying%
\[
\operatorname{st} \pi   =\operatorname{st}\left(  \pi^{\prime
}\right)  ,\ \ \ \ \ \ \ \ \ \ \left\vert \pi\right\vert =\left\vert
\pi^{\prime}\right\vert \ \ \ \ \ \ \ \ \ \ \text{and}%
\ \ \ \ \ \ \ \ \ \ \left[  a>\pi_{1}\right]  =\left[  a^{\prime}>\pi
_{1}^{\prime}\right]
\]
satisfy $\operatorname{st}\left(  a:\pi\right)  =\operatorname{st}\left(
a^{\prime}:\pi^{\prime}\right)  $. Applying this to $\pi=\alpha_{\sim1}$,
$a=\alpha_{1}$, $\pi^{\prime}=\beta_{\sim1}$ and $a^{\prime}=\beta_{1}$, we
obtain $\operatorname{st}\left(  \alpha_{1}:\alpha_{\sim1}\right)
=\operatorname{st}\left(  \beta_{1}:\beta_{\sim1}\right)  $. In view of
$\alpha_{1}:\alpha_{\sim1}=\alpha$ and $\beta_{1}:\beta_{\sim1}=\beta$, this
rewrites as $\operatorname{st} \alpha = \operatorname{st} \beta$.

Thus, we have proven that $\operatorname{st} \alpha
= \operatorname{st} \beta$. Hence, Claim 1 holds for
$\left\vert \alpha\right\vert =N$. This completes the induction step. Thus,
Claim 1 is proven.]

Lemma \ref{lem.LRcomp.head-l1} follows immediately from Claim 1.
\end{proof}

\begin{lemma}
\label{lem.LRcomp.head-l2}Let $\operatorname{st}$ be a head-graft-compatible
permutation statistic. Let $X$ be the codomain of $\operatorname{st}$. Let
$n\in\mathbb{N}$. Then, there exists a map
\[
F_{n}:\left\{  \text{subsets of }\left[  n-1\right]  \right\}  \rightarrow X
\]
such that every $n$-permutation $\tau$ satisfies $\operatorname{st} \tau
= F_{n}\left(  \operatorname{Des}\tau\right)  $.
\end{lemma}

\begin{proof}
[Proof of Lemma \ref{lem.LRcomp.head-l2}.]We define the map $F_{n}$ as follows:

Let $Z$ be any subset of $\left[  n-1\right]  $. Then, it is well-known that
there exists some $n$-permutation $\tau$ satisfying $Z=\operatorname{Des}
\tau$. Pick any such $\tau$. Then, $\operatorname{st} \tau$
does not depend on the choice of $\tau$ (because if $\alpha$ and $\beta$ are
two different $n$-permutations $\tau$ satisfying $Z=\operatorname{Des}\tau$,
then Lemma \ref{lem.LRcomp.head-l1} yields
$\operatorname{st} \alpha = \operatorname{st} \beta$).
Hence, we can set
$F_{n}\left(  Z\right)  $ to be $\operatorname{st} \tau$.

Thus, we have defined $F_{n}\left(  Z\right)  $ for each subset $Z$ of
$\left[  n-1\right]  $. This completes the definition of $F_{n}$. This
definition shows that every $n$-permutation $\tau$ satisfies
$\operatorname{st} \tau   =F_{n}\left(  \operatorname{Des}
\tau\right)  $. Thus, Lemma \ref{lem.LRcomp.head-l2} is proven.
\end{proof}

\begin{lemma}
\label{lem.LRcomp.head-l3}Let $\pi$ and $\sigma$ be two disjoint permutations.
Let $\pi^{\prime}$ and $\sigma^{\prime}$ be two disjoint permutations. Assume
that%
\begin{align*}
\operatorname{Des}\pi &  =\operatorname{Des} \left( \pi^{\prime} \right)
,\ \ \ \ \ \ \ \ \ \ \operatorname{Des}\sigma = \operatorname{Des}
\left( \sigma^{\prime} \right),\\
\left\vert \pi\right\vert  &  =\left\vert \pi^{\prime}\right\vert
\ \ \ \ \ \ \ \ \ \ \text{and}
\ \ \ \ \ \ \ \ \ \ \left\vert \sigma\right\vert
=\left\vert \sigma^{\prime}\right\vert .
\end{align*}
Then,
\[
\left\{  \operatorname{Des}\tau\ \mid\ \tau\in S\left(  \pi,\sigma\right)
\right\}  _{\operatorname*{multi}}=\left\{  \operatorname{Des}\tau
\ \mid\ \tau\in S\left(  \pi^{\prime},\sigma^{\prime}\right)  \right\}
_{\operatorname*{multi}}.
\]

\end{lemma}

\begin{proof}
[Proof of Lemma \ref{lem.LRcomp.head-l3}.]Lemma \ref{lem.LRcomp.head-l3} is
simply the statement that the permutation statistic $\operatorname{Des}$ is
shuffle-compatible. But this has been proven in \cite[\S 2.4]{part1}.
\end{proof}

\begin{proof}
[Proof of Proposition \ref{prop.LRcomp.head} (sketched).]We know that
$\operatorname{st}$ is LR-shuffle-compatible. In other words, the following holds:

\begin{statement}
\textit{Claim 1:} Let $\pi$ and $\sigma$ be two disjoint nonempty
permutations. Let $\pi^{\prime}$ and $\sigma^{\prime}$ be two disjoint
nonempty permutations. Assume that%
\begin{align*}
\operatorname{st} \pi    &  =\operatorname{st}\left(
\pi^{\prime}\right)  ,\ \ \ \ \ \ \ \ \ \ \operatorname{st} \sigma
=\operatorname{st}\left(  \sigma^{\prime}\right)  ,\\
\left\vert \pi\right\vert  &  =\left\vert \pi^{\prime}\right\vert
,\ \ \ \ \ \ \ \ \ \ \left\vert \sigma\right\vert =\left\vert \sigma^{\prime
}\right\vert \ \ \ \ \ \ \ \ \ \ \text{and}\ \ \ \ \ \ \ \ \ \ \left[  \pi
_{1}>\sigma_{1}\right]  =\left[  \pi_{1}^{\prime}>\sigma_{1}^{\prime}\right]
.
\end{align*}
Then,
\[
\left\{  \operatorname{st} \tau   \ \mid\ \tau\in S_{\prec
}\left(  \pi,\sigma\right)  \right\}  _{\operatorname*{multi}}=\left\{
\operatorname{st} \tau   \ \mid\ \tau\in S_{\prec}\left(
\pi^{\prime},\sigma^{\prime}\right)  \right\}  _{\operatorname*{multi}}%
\]
and%
\[
\left\{  \operatorname{st} \tau   \ \mid\ \tau\in S_{\succ
}\left(  \pi,\sigma\right)  \right\}  _{\operatorname*{multi}}=\left\{
\operatorname{st} \tau   \ \mid\ \tau\in S_{\succ}\left(
\pi^{\prime},\sigma^{\prime}\right)  \right\}  _{\operatorname*{multi}}.
\]

\end{statement}

Next, we want to show that $\operatorname{st}$ is head-graft-compatible. In
other words, we want to show that the following holds:

\begin{statement}
\textit{Claim 2:} Let $\pi$ be a nonempty permutation, and let $a$ be a letter
that does not appear in $\pi$. Let $\pi^{\prime}$ be a nonempty permutation,
and let $a^{\prime}$ be a letter that does not appear in $\pi^{\prime}$.
Assume that%
\[
\operatorname{st} \pi   =\operatorname{st}\left(  \pi^{\prime
}\right)  ,\ \ \ \ \ \ \ \ \ \ \left\vert \pi\right\vert =\left\vert
\pi^{\prime}\right\vert \ \ \ \ \ \ \ \ \ \ \text{and}%
\ \ \ \ \ \ \ \ \ \ \left[  a>\pi_{1}\right]  =\left[  a^{\prime}>\pi
_{1}^{\prime}\right]  .
\]
Then, $\operatorname{st}\left(  a:\pi\right)  =\operatorname{st}\left(
a^{\prime}:\pi^{\prime}\right)  $.
\end{statement}

[\textit{Proof of Claim 2:} The $1$-permutations $\left(  a\right)  $ and
$\left(  a^{\prime}\right)  $ are order-isomorphic. Thus, $\operatorname{st}
\left(  \left(  a\right)  \right)  =\operatorname{st}\left(  \left(
a^{\prime}\right)  \right)  $ (since $\operatorname{st}$ is a permutation
statistic).

We have assumed that $a$ does not appear in $\pi$. Hence, in
particular, $a \neq \pi_1$. Thus, $\pi_1>a$ is true if and only if
$a>\pi_1$ is false (since the third possibility, $a = \pi_1$,
is precluded by the preceding sentence).
Thus, $\left[\pi_1>a\right] = \left[a>\pi_1 \text{ is false}\right]
= 1 - \left[a>\pi_1\right]$.
Similarly, $\left[\pi^\prime_1>a^\prime\right]
= 1 - \left[a^\prime>\pi^\prime_1\right]$.
Now,
\[
\left[\pi_1>a\right]
= 1 - \underbrace{\left[  a>\pi_{1}\right]}_{= \left[  a^{\prime}>\pi_{1}^{\prime}\right]}
= 1 - \left[  a^{\prime}>\pi_{1}^{\prime}\right]
= \left[\pi^\prime_1>a^\prime\right].
\]

The $1$-permutations $\left(  a\right)  $ and $\left(  a^{\prime}\right)  $
satisfy $\left(  a\right)  _{1}=a$ and $\left(  a^{\prime}\right)
_{1}=a^{\prime}$. Thus, the inequality $\left[\pi_1>a\right]
= \left[\pi^\prime_1>a^\prime\right]$ (which we have just shown) rewrites
as $\left[  \pi_1 > \left(  a\right)  _{1}\right]  =\left[  \pi^\prime_1
> \left(  a^{\prime}\right)  _{1}\right]  $. Also, $\operatorname{st}\left(
\left(  a\right)  \right)  =\operatorname{st}\left(  \left(  a^{\prime
}\right)  \right)  $ and $\left\vert \left(  a\right)  \right\vert
=1=\left\vert \left(  a^{\prime}\right)  \right\vert $. Also, the
$1$-permutation $\left(  a\right)  $ is disjoint from $\pi$ (since $a$ does
not appear in $\pi$). Similarly, the $1$-permutation $\left(  a^{\prime
}\right)  $ is disjoint from $\pi^{\prime}$. Thus, Claim 1 (applied to
$\sigma=\left(  a\right)  $ and $\sigma^{\prime}=\left(  a^{\prime}\right)  $)
yields
\[
\left\{  \operatorname{st} \tau   \ \mid\ \tau\in S_{\prec
}\left(  \pi,\left(  a\right)  \right)  \right\}  _{\operatorname*{multi}%
}=\left\{  \operatorname{st} \tau   \ \mid\ \tau\in S_{\prec
}\left(  \pi^{\prime},\left(  a^{\prime}\right)  \right)  \right\}
_{\operatorname*{multi}}%
\]
and%
\[
\left\{  \operatorname{st} \tau   \ \mid\ \tau\in S_{\succ
}\left(  \pi,\left(  a\right)  \right)  \right\}  _{\operatorname*{multi}%
}=\left\{  \operatorname{st} \tau   \ \mid\ \tau\in S_{\succ
}\left(  \pi^{\prime},\left(  a^{\prime}\right)  \right)  \right\}
_{\operatorname*{multi}}.
\]
But it is easily seen that $S_{\succ}\left(  \pi,\left(  a\right)  \right)
=\left\{  a:\pi\right\}  $. Hence,%
\[
\left\{  \operatorname{st} \tau   \ \mid\ \tau\in S_{\succ
}\left(  \pi,\left(  a\right)  \right)  \right\}  _{\operatorname*{multi}%
}=\left\{  \operatorname{st} \tau   \ \mid\ \tau\in\left\{
a:\pi\right\}  \right\}  _{\operatorname*{multi}}=\left\{  \operatorname{st}
\left(  a:\pi\right)  \right\}  _{\operatorname*{multi}}.
\]
Similarly,%
\[
\left\{  \operatorname{st} \tau   \ \mid\ \tau\in S_{\succ
}\left(  \pi^{\prime},\left(  a^{\prime}\right)  \right)  \right\}
_{\operatorname*{multi}}=\left\{  \operatorname{st}\left(  a^{\prime}%
:\pi^{\prime}\right)  \right\}  _{\operatorname*{multi}}.
\]
Thus,%
\begin{align*}
\left\{  \operatorname{st}\left(  a:\pi\right)  \right\}
_{\operatorname*{multi}}  &  =\left\{  \operatorname{st} \tau 
\ \mid\ \tau\in S_{\succ}\left(  \pi,\left(  a\right)  \right)  \right\}
_{\operatorname*{multi}}\\
&  =\left\{  \operatorname{st} \tau   \ \mid\ \tau\in S_{\succ
}\left(  \pi^{\prime},\left(  a^{\prime}\right)  \right)  \right\}
_{\operatorname*{multi}}=\left\{  \operatorname{st}\left(  a^{\prime}%
:\pi^{\prime}\right)  \right\}  _{\operatorname*{multi}},
\end{align*}
so that $\operatorname{st}\left(  a:\pi\right)  =\operatorname{st}\left(
a^{\prime}:\pi^{\prime}\right)  $. This proves Claim 2.]

Now, Claim 2 shows that $\operatorname{st}$ is head-graft-compatible. It
remains to show that $\operatorname{st}$ is shuffle-compatible. First, we
show an auxiliary statement:

\begin{statement}
\textit{Claim 3:} Let $\pi$, $\pi^{\prime}$ and $\sigma$ be three
permutations. Assume that $\pi$ and $\pi^{\prime}$ are order-isomorphic.
Assume that $\pi$ and $\sigma$ are disjoint. Assume that $\pi^{\prime}$ and
$\sigma$ are disjoint. Then,%
\[
\left\{  \operatorname{st} \tau   \ \mid\ \tau\in S\left(
\pi,\sigma\right)  \right\}  _{\operatorname*{multi}}=\left\{
\operatorname{st} \tau   \ \mid\ \tau\in S\left(  \pi^{\prime
},\sigma\right)  \right\}  _{\operatorname*{multi}}.
\]

\end{statement}

[\textit{Proof of Claim 3:} We have $\left\vert \pi\right\vert =\left\vert
\pi^{\prime}\right\vert $ (since $\pi$ and $\pi^{\prime}$ are
order-isomorphic). Define $n\in\mathbb{N}$ by $n=\left\vert \pi\right\vert
=\left\vert \pi^{\prime}\right\vert $. Let $X$ be the codomain of
$\operatorname{st}$. Lemma \ref{lem.LRcomp.head-l2} shows that there exists a
map
\[
F_{n}:\left\{  \text{subsets of }\left[  n-1\right]  \right\}  \rightarrow X
\]
such that every $n$-permutation $\tau$ satisfies
\begin{equation}
\operatorname{st} \tau   =F_{n}\left(  \operatorname{Des} \tau\right)
. \label{pf.lem.LRcomp.head-l2.c3.pf.1}
\end{equation}
Consider this map $F_{n}$.

We know that $\operatorname{Des}$ is a permutation statistic. Thus,
$\operatorname{Des}\pi=\operatorname{Des}\left(  \pi^{\prime}\right)  $
(since $\pi$ and $\pi^{\prime}$ are order-isomorphic). Also, $\left\vert
\pi\right\vert =\left\vert \pi^{\prime}\right\vert $, $\operatorname{Des}
\sigma=\operatorname{Des}\sigma$ and $\left\vert \sigma\right\vert
=\left\vert \sigma\right\vert $. Hence, Lemma \ref{lem.LRcomp.head-l3}
(applied to $\sigma^{\prime}=\sigma$) yields
\begin{equation}
\left\{  \operatorname{Des}\tau\ \mid\ \tau\in S\left(  \pi,\sigma\right)
\right\}  _{\operatorname*{multi}}=\left\{  \operatorname{Des}\tau
\ \mid\ \tau\in S\left(  \pi^{\prime},\sigma\right)  \right\}
_{\operatorname*{multi}}. \label{pf.lem.LRcomp.head-l2.c3.pf.2}%
\end{equation}
Now,
\begin{align*}
&  \left\{  \underbrace{\operatorname{st} \tau   }%
_{\substack{=F_{n}\left(  \operatorname{Des}\tau\right)  \\\text{(by
(\ref{pf.lem.LRcomp.head-l2.c3.pf.1}))}}}\ \mid\ \tau\in S\left(  \pi
,\sigma\right)  \right\}  _{\operatorname*{multi}}\\
&  =\left\{  F_{n}\left(  \operatorname{Des}\tau\right)  \ \mid\ \tau\in
S\left(  \pi,\sigma\right)  \right\}  _{\operatorname*{multi}}=F_{n}\left(
\left\{  \operatorname{Des}\tau\ \mid\ \tau\in S\left(  \pi,\sigma\right)
\right\}  _{\operatorname*{multi}}\right) \\
&  =F_{n}\left(  \left\{  \operatorname{Des}\tau\ \mid\ \tau\in S\left(
\pi^{\prime},\sigma\right)  \right\}  _{\operatorname*{multi}}\right)
\ \ \ \ \ \ \ \ \ \ \left(  \text{by (\ref{pf.lem.LRcomp.head-l2.c3.pf.2}%
)}\right) \\
&  =\left\{  \underbrace{F_{n}\left(  \operatorname{Des}\tau\right)
}_{\substack{=\operatorname{st} \tau   \\\text{(by
(\ref{pf.lem.LRcomp.head-l2.c3.pf.1}))}}}\ \mid\ \tau\in S\left(  \pi^{\prime
},\sigma\right)  \right\}  _{\operatorname*{multi}}=\left\{
\operatorname{st} \tau   \ \mid\ \tau\in S\left(  \pi^{\prime
},\sigma\right)  \right\}  _{\operatorname*{multi}}.
\end{align*}
This proves Claim 3.]

Let us finally prove that $\operatorname{st}$ is shuffle-compatible. In other
words, let us prove the following claim:

\begin{statement}
\textit{Claim 4:} Let $\pi$ and $\sigma$ be two disjoint permutations. Let
$\pi^{\prime}$ and $\sigma^{\prime}$ be two disjoint permutations. Assume that%
\begin{align*}
\operatorname{st} \pi    &  =\operatorname{st}\left(
\pi^{\prime}\right)  ,\ \ \ \ \ \ \ \ \ \ \operatorname{st} \sigma
=\operatorname{st}\left(  \sigma^{\prime}\right)  ,\\
\left\vert \pi\right\vert  &  =\left\vert \pi^{\prime}\right\vert
\ \ \ \ \ \ \ \ \ \ \text{and}\ \ \ \ \ \ \ \ \ \ \left\vert \sigma\right\vert
=\left\vert \sigma^{\prime}\right\vert .
\end{align*}
Then,
\[
\left\{  \operatorname{st} \tau   \ \mid\ \tau\in S\left(
\pi,\sigma\right)  \right\}  _{\operatorname*{multi}}=\left\{
\operatorname{st} \tau   \ \mid\ \tau\in S\left(  \pi^{\prime
},\sigma^{\prime}\right)  \right\}  _{\operatorname*{multi}}.
\]

\end{statement}

[\textit{Proof of Claim 4:}
Define $n\in\mathbb{N}$ by $n=\left\vert \pi\right\vert =\left\vert
\pi^{\prime}\right\vert $. Define $m\in\mathbb{N}$ by $m=\left\vert
\sigma\right\vert =\left\vert \sigma^{\prime}\right\vert $.

If $n=0$, then both $\pi$ and $\pi^{\prime}$ equal the $0$-permutation
$\left(  {}\right)  $ (since $n=\left\vert \pi\right\vert =\left\vert
\pi^{\prime}\right\vert $). Hence, if $n=0$, then Claim 4 is true (because in
this case, we have%
\begin{align*}
\left\{  \operatorname{st} \tau   \ \mid\ \tau\in S\left(
\underbrace{\pi}_{=\left(  {}\right)  },\sigma\right)  \right\}
_{\operatorname*{multi}}  &  =\left\{  \operatorname{st} \tau 
\ \mid\ \tau\in\underbrace{S\left(  \left(  {}\right)  ,\sigma\right)
}_{=\left\{  \sigma\right\}  }\right\}  _{\operatorname*{multi}}\\
&  =\left\{  \operatorname{st} \tau   \ \mid\ \tau\in\left\{
\sigma\right\}  \right\}  _{\operatorname*{multi}}
= \left\{  \operatorname{st} \sigma \right\}  _{\operatorname*{multi}}
\end{align*}
and similarly $\left\{  \operatorname{st} \tau   \ \mid\ \tau\in
S\left(  \pi^{\prime},\sigma^{\prime}\right)  \right\}
_{\operatorname*{multi}}=\left\{  \operatorname{st}\left(  \sigma^{\prime
}\right)  \right\}  _{\operatorname*{multi}}$, and therefore the assertion of
Claim 4 reduces to $\left\{  \operatorname{st} \sigma 
\right\}  _{\operatorname*{multi}}=\left\{  \operatorname{st}\left(
\sigma^{\prime}\right)  \right\}  _{\operatorname*{multi}}$, which follows
immediately from the assumption $\operatorname{st} \sigma 
=\operatorname{st}\left(  \sigma^{\prime}\right)  $). Thus, for the rest of
this proof, we WLOG assume that $n\neq0$. For similar reasons, we WLOG assume
that $m\neq0$.

The permutations $\pi$ and $\pi^{\prime}$ are nonempty (since $\left\vert
\pi\right\vert =\left\vert \pi^{\prime}\right\vert =n\neq0$). Similarly, the
permutations $\sigma$ and $\sigma^{\prime}$ are nonempty.

Recall that $S\left(  \pi,\sigma\right)  =S\left(
\sigma,\pi\right)  $ and $S\left(  \pi^{\prime},\sigma^{\prime}\right)
=S\left(  \sigma^{\prime},\pi^{\prime}\right)  $. Hence, Claim 4 does not
change if we swap $\pi$ with $\sigma$ while simultaneously swapping
$\pi^{\prime}$ with $\sigma^{\prime}$. Thus, we WLOG assume that $\pi_{1}%
\geq\sigma_{1}$ (since otherwise, we can ensure this by performing these
swaps). But $\pi_{1}\neq\sigma_{1}$ (since $\pi$ and $\sigma$ are disjoint).
Combining this with $\pi_{1}\geq\sigma_{1}$, we obtain $\pi_{1}>\sigma_{1}$.

There clearly exists a positive integer $N$ that is larger than all entries of
$\sigma$ and larger than all entries of $\sigma^{\prime}$. Consider such an
$N$. From $n=\left\vert \pi^{\prime}\right\vert $, we obtain $\pi^{\prime
}=\left(  \pi_{1}^{\prime},\pi_{2}^{\prime},\ldots,\pi_{n}^{\prime}\right)  $.
Let $\gamma$ be the permutation $\left(  \pi_{1}^{\prime}+N,\pi_{2}^{\prime
}+N,\ldots,\pi_{n}^{\prime}+N\right)  $. This permutation $\gamma$ is
order-isomorphic to $\pi^{\prime}$, but is disjoint from $\sigma$ (since all
its entries are $>N$, while all the entries of $\sigma$ are $<N$) and disjoint
from $\sigma^{\prime}$ (for similar reasons). Also, $\gamma_{1}%
=\underbrace{\pi_{1}^{\prime}}_{>0}+N>N>\sigma_{1}$, so that $\left[
\gamma_{1}>\sigma_{1}\right]  =1$. Similarly, $\left[  \gamma_{1}>\sigma
_{1}^{\prime}\right]  =1$.

The permutation $\gamma$ is order-isomorphic to $\pi^{\prime}$. Thus,
$\operatorname{st} \gamma = \operatorname{st}\left( \pi^{\prime}\right)  $
(since $\operatorname{st}$ is a permutation statistic)
and $\left\vert \gamma\right\vert =\left\vert \pi^{\prime}\right\vert $. The
permutation $\gamma$ is furthermore nonempty (since it is order-isomorphic to
the nonempty permutation $\pi^{\prime}$). Also, $\operatorname{st}
\gamma =\operatorname{st}\left(  \pi^{\prime}\right)
=\operatorname{st} \pi   $ and $\left\vert \gamma\right\vert
=\left\vert \pi^{\prime}\right\vert =\left\vert \pi\right\vert $. Moreover,
from $\pi_{1}>\sigma_{1}$, we obtain $\left[  \pi_{1}>\sigma_{1}\right]
=1=\left[  \gamma_{1}>\sigma_{1}^{\prime}\right]  $. Hence, Claim 1 (applied
to $\gamma$ instead of $\pi^{\prime}$) yields
\begin{equation}
\left\{  \operatorname{st} \tau   \ \mid\ \tau\in S_{\prec
}\left(  \pi,\sigma\right)  \right\}  _{\operatorname*{multi}}=\left\{
\operatorname{st} \tau   \ \mid\ \tau\in S_{\prec}\left(
\gamma,\sigma^{\prime}\right)  \right\}  _{\operatorname*{multi}}
\label{pf.lem.LRcomp.head-l2.c4.pf.5}%
\end{equation}
and%
\begin{equation}
\left\{  \operatorname{st} \tau   \ \mid\ \tau\in S_{\succ
}\left(  \pi,\sigma\right)  \right\}  _{\operatorname*{multi}}=\left\{
\operatorname{st} \tau   \ \mid\ \tau\in S_{\succ}\left(
\gamma,\sigma^{\prime}\right)  \right\}  _{\operatorname*{multi}}.
\label{pf.lem.LRcomp.head-l2.c4.pf.6}%
\end{equation}

Now, if $A$ and $B$ are two finite multisets, then $A+B$ shall denote the
multiset union of $A$ and $B$; this is the finite multiset $C$ such that
each object
$g$ satisfies $\left\vert C\right\vert _{g}=\left\vert A\right\vert
_{g}+\left\vert B\right\vert _{g}$. Then, it is easy to see that%
\begin{align*}
&  \left\{  \operatorname{st} \tau   \ \mid\ \tau\in S\left(
\pi,\sigma\right)  \right\}  _{\operatorname*{multi}}\\
&  =\left\{  \operatorname{st} \tau   \ \mid\ \tau\in S_{\prec
}\left(  \pi,\sigma\right)  \right\}  _{\operatorname*{multi}}+\left\{
\operatorname{st} \tau   \ \mid\ \tau\in S_{\succ}\left(
\pi,\sigma\right)  \right\}  _{\operatorname*{multi}}%
\end{align*}
and
\begin{align*}
&  \left\{  \operatorname{st} \tau   \ \mid\ \tau\in S\left(
\gamma,\sigma^{\prime}\right)  \right\}  _{\operatorname*{multi}}\\
&  =\left\{  \operatorname{st} \tau   \ \mid\ \tau\in S_{\prec
}\left(  \gamma,\sigma^{\prime}\right)  \right\}  _{\operatorname*{multi}%
}+\left\{  \operatorname{st} \tau   \ \mid\ \tau\in S_{\succ
}\left(  \gamma,\sigma^{\prime}\right)  \right\}  _{\operatorname*{multi}}.
\end{align*}
Hence, by adding together the equalities (\ref{pf.lem.LRcomp.head-l2.c4.pf.5})
and (\ref{pf.lem.LRcomp.head-l2.c4.pf.6}) (using the operation $+$ that we
have just defined), we obtain the equality%
\[
\left\{  \operatorname{st} \tau   \ \mid\ \tau\in S\left(
\pi,\sigma\right)  \right\}  _{\operatorname*{multi}}=\left\{
\operatorname{st} \tau   \ \mid\ \tau\in S\left(  \gamma
,\sigma^{\prime}\right)  \right\}  _{\operatorname*{multi}}.
\]

On the other hand, Claim 3 (applied to $\gamma$ and $\sigma^{\prime}$ instead
of $\pi$ and $\sigma$) yields%
\[
\left\{  \operatorname{st} \tau   \ \mid\ \tau\in S\left(
\gamma,\sigma^{\prime}\right)  \right\}  _{\operatorname*{multi}}=\left\{
\operatorname{st} \tau   \ \mid\ \tau\in S\left(  \pi^{\prime
},\sigma^{\prime}\right)  \right\}  _{\operatorname*{multi}}.
\]
Hence,%
\begin{align*}
\left\{  \operatorname{st} \tau   \ \mid\ \tau\in S\left(
\pi,\sigma\right)  \right\}  _{\operatorname*{multi}}  &  =\left\{
\operatorname{st} \tau   \ \mid\ \tau\in S\left(  \gamma
,\sigma^{\prime}\right)  \right\}  _{\operatorname*{multi}}\\
&  =\left\{  \operatorname{st} \tau   \ \mid\ \tau\in S\left(
\pi^{\prime},\sigma^{\prime}\right)  \right\}  _{\operatorname*{multi}}.
\end{align*}
This proves Claim 4.]

Claim 4 shows that $\operatorname{st}$ is shuffle-compatible. This completes
the proof of Proposition \ref{prop.LRcomp.head}.
\end{proof}
\end{verlong}

\begin{corollary}
\label{cor.LRcomp.back}Let $\operatorname{st}$ be a LR-shuffle-compatible
permutation statistic. Then, $\operatorname{st}$ is shuffle-compatible,
left-shuffle-compatible, right-shuffle-compatible and head-graft-compatible.
\end{corollary}

\begin{verlong}
\begin{proof}
[Proof of Corollary \ref{cor.LRcomp.back}.]Proposition \ref{prop.LRcomp.head}
yields that $\operatorname{st}$ is head-graft-compatible and
shuffle-compatible. Thus, Proposition \ref{prop.LRcomp.equivs} yields that
$\operatorname{st}$ is left-shuffle-compatible and right-shuffle-compatible
as well. This proves Corollary \ref{cor.LRcomp.back}.
\end{proof}
\end{verlong}

\begin{corollary}
\label{cor.LRcomp.two}Let $\operatorname{st}$ be a permutation statistic that
is left-shuffle-compatible and right-shuffle-compatible. Then,
$\operatorname{st}$ is LR-shuffle-compatible.
\end{corollary}

\begin{verlong}
\begin{proof}
[Proof of Corollary \ref{cor.LRcomp.two}.]We have assumed that
$\operatorname{st}$ is left-shuffle-compatible. In other words, the following
claim holds:

\begin{statement}
\textit{Claim 1:} Let $\pi$ and $\sigma$ be two disjoint nonempty permutations
having the property that $\pi_{1}>\sigma_{1}$. Let $\pi^{\prime}$ and
$\sigma^{\prime}$ be two disjoint nonempty permutations having the property
that $\pi_{1}^{\prime}>\sigma_{1}^{\prime}$. Assume that%
\begin{align*}
\operatorname{st} \pi    &  =\operatorname{st}\left(
\pi^{\prime}\right)  ,\ \ \ \ \ \ \ \ \ \ \operatorname{st} \sigma
=\operatorname{st}\left(  \sigma^{\prime}\right)  ,\\
\left\vert \pi\right\vert  &  =\left\vert \pi^{\prime}\right\vert
\ \ \ \ \ \ \ \ \ \ \text{and}\ \ \ \ \ \ \ \ \ \ \left\vert \sigma\right\vert
=\left\vert \sigma^{\prime}\right\vert .
\end{align*}
Then,
\[
\left\{  \operatorname{st} \tau   \ \mid\ \tau\in S_{\prec
}\left(  \pi,\sigma\right)  \right\}  _{\operatorname*{multi}}=\left\{
\operatorname{st} \tau   \ \mid\ \tau\in S_{\prec}\left(
\pi^{\prime},\sigma^{\prime}\right)  \right\}  _{\operatorname*{multi}}.
\]

\end{statement}

Also, the statistic $\operatorname{st}$ is right-shuffle-compatible. In other
words, the following claim holds:

\begin{statement}
\textit{Claim 2:} Let $\pi$ and $\sigma$ be two disjoint nonempty permutations
having the property that $\pi_{1}>\sigma_{1}$. Let $\pi^{\prime}$ and
$\sigma^{\prime}$ be two disjoint nonempty permutations having the property
that $\pi_{1}^{\prime}>\sigma_{1}^{\prime}$. Assume that%
\begin{align*}
\operatorname{st} \pi    &  =\operatorname{st}\left(
\pi^{\prime}\right)  ,\ \ \ \ \ \ \ \ \ \ \operatorname{st} \sigma
= \operatorname{st}\left(  \sigma^{\prime}\right)  ,\\
\left\vert \pi\right\vert  &  =\left\vert \pi^{\prime}\right\vert
\ \ \ \ \ \ \ \ \ \ \text{and}\ \ \ \ \ \ \ \ \ \ \left\vert \sigma\right\vert
=\left\vert \sigma^{\prime}\right\vert .
\end{align*}
Then,
\[
\left\{  \operatorname{st} \tau   \ \mid\ \tau\in S_{\succ
}\left(  \pi,\sigma\right)  \right\}  _{\operatorname*{multi}}=\left\{
\operatorname{st} \tau   \ \mid\ \tau\in S_{\succ}\left(
\pi^{\prime},\sigma^{\prime}\right)  \right\}  _{\operatorname*{multi}}.
\]

\end{statement}

On the other hand, we want to prove that $\operatorname{st}$ is
LR-shuffle-compatible. In other words, we want to prove the following claim:

\begin{statement}
\textit{Claim 3:} Let $\pi$ and $\sigma$ be two disjoint nonempty
permutations. Let $\pi^{\prime}$ and $\sigma^{\prime}$ be two disjoint
nonempty permutations. Assume that%
\begin{align*}
\operatorname{st} \pi    &  =\operatorname{st}\left(
\pi^{\prime}\right)  ,\ \ \ \ \ \ \ \ \ \ \operatorname{st} \sigma
= \operatorname{st}\left(  \sigma^{\prime}\right)  ,\\
\left\vert \pi\right\vert  &  =\left\vert \pi^{\prime}\right\vert
,\ \ \ \ \ \ \ \ \ \ \left\vert \sigma\right\vert =\left\vert \sigma^{\prime
}\right\vert \ \ \ \ \ \ \ \ \ \ \text{and}\ \ \ \ \ \ \ \ \ \ \left[  \pi
_{1}>\sigma_{1}\right]  =\left[  \pi_{1}^{\prime}>\sigma_{1}^{\prime}\right]
.
\end{align*}
Then,
\[
\left\{  \operatorname{st} \tau   \ \mid\ \tau\in S_{\prec
}\left(  \pi,\sigma\right)  \right\}  _{\operatorname*{multi}}=\left\{
\operatorname{st} \tau   \ \mid\ \tau\in S_{\prec}\left(
\pi^{\prime},\sigma^{\prime}\right)  \right\}  _{\operatorname*{multi}}%
\]
and%
\[
\left\{  \operatorname{st} \tau   \ \mid\ \tau\in S_{\succ
}\left(  \pi,\sigma\right)  \right\}  _{\operatorname*{multi}}=\left\{
\operatorname{st} \tau   \ \mid\ \tau\in S_{\succ}\left(
\pi^{\prime},\sigma^{\prime}\right)  \right\}  _{\operatorname*{multi}}.
\]

\end{statement}

[\textit{Proof of Claim 3:} We are in one of the following two cases:

\textit{Case 1:} We have $\pi_{1}>\sigma_{1}$.

\textit{Case 2:} We have $\pi_{1}\leq\sigma_{1}$.

Let us first consider Case 1. In this case, we have $\pi_{1}>\sigma_{1}$.
Hence, $\left[  \pi_{1}>\sigma_{1}\right]  =1$. Comparing this with $\left[
\pi_{1}>\sigma_{1}\right]  =\left[  \pi_{1}^{\prime}>\sigma_{1}^{\prime
}\right]  $, we find $\left[  \pi_{1}^{\prime}>\sigma_{1}^{\prime}\right]
=1$. Hence, $\pi_{1}^{\prime}>\sigma_{1}^{\prime}$. Recall also that $\pi
_{1}>\sigma_{1}$. Hence, Claim 1 yields%
\[
\left\{  \operatorname{st} \tau   \ \mid\ \tau\in S_{\prec
}\left(  \pi,\sigma\right)  \right\}  _{\operatorname*{multi}}=\left\{
\operatorname{st} \tau   \ \mid\ \tau\in S_{\prec}\left(
\pi^{\prime},\sigma^{\prime}\right)  \right\}  _{\operatorname*{multi}}.
\]
But Claim 2 yields%
\[
\left\{  \operatorname{st} \tau   \ \mid\ \tau\in S_{\succ
}\left(  \pi,\sigma\right)  \right\}  _{\operatorname*{multi}}=\left\{
\operatorname{st} \tau   \ \mid\ \tau\in S_{\succ}\left(
\pi^{\prime},\sigma^{\prime}\right)  \right\}  _{\operatorname*{multi}}.
\]
Thus, Claim 3 is proven in Case 1.

Let us next consider Case 2. In this case, we have $\pi_{1}\leq\sigma_{1}$.

Applying Proposition \ref{prop.LR.rec} \textbf{(a)} to various pairs of
disjoint permutations, we obtain $S_{\prec}\left(  \pi,\sigma\right)
=S_{\succ}\left(  \sigma,\pi\right)  $ and $S_{\prec}\left(  \pi^{\prime
},\sigma^{\prime}\right)  =S_{\succ}\left(  \sigma^{\prime},\pi^{\prime
}\right)  $ and $S_{\prec}\left(  \sigma,\pi\right)  =S_{\succ}\left(
\pi,\sigma\right)  $ and $S_{\prec}\left(  \sigma^{\prime},\pi^{\prime
}\right)  =S_{\succ}\left(  \pi^{\prime},\sigma^{\prime}\right)  $.

But $\pi_{1}\neq\sigma_{1}$ (since $\pi$ and $\sigma$ are disjoint). Combined
with $\pi_{1}\leq\sigma_{1}$, this yields $\pi_{1}<\sigma_{1}$; thus,
$\sigma_{1}>\pi_{1}$. Also, $\left[  \pi_{1}>\sigma_{1}\right]  =0$ (since
$\pi_{1}\leq\sigma_{1}$). Comparing this with $\left[  \pi_{1}>\sigma
_{1}\right]  =\left[  \pi_{1}^{\prime}>\sigma_{1}^{\prime}\right]  $, we find
$\left[  \pi_{1}^{\prime}>\sigma_{1}^{\prime}\right]  =0$. Hence, $\pi
_{1}^{\prime}\leq\sigma_{1}^{\prime}$. But $\pi_{1}^{\prime}\neq\sigma
_{1}^{\prime}$ (since $\pi^{\prime}$ and $\sigma^{\prime}$ are disjoint).
Combined with $\pi_{1}^{\prime}\leq\sigma_{1}^{\prime}$, this yields $\pi
_{1}^{\prime}<\sigma_{1}^{\prime}$; thus, $\sigma_{1}^{\prime}>\pi_{1}%
^{\prime}$. Hence, Claim 2 (applied to $\sigma$, $\pi$, $\sigma^{\prime}$ and
$\pi^{\prime}$ instead of $\pi$, $\sigma$, $\pi^{\prime}$ and $\sigma^{\prime
}$) yields%
\[
\left\{  \operatorname{st} \tau   \ \mid\ \tau\in S_{\succ
}\left(  \sigma,\pi\right)  \right\}  _{\operatorname*{multi}}=\left\{
\operatorname{st} \tau   \ \mid\ \tau\in S_{\succ}\left(
\sigma^{\prime},\pi^{\prime}\right)  \right\}  _{\operatorname*{multi}}.
\]
In light of $S_{\succ}\left(  \sigma,\pi\right)  =S_{\prec}\left(  \pi
,\sigma\right)  $ and $S_{\succ}\left(  \sigma^{\prime},\pi^{\prime}\right)
=S_{\prec}\left(  \pi^{\prime},\sigma^{\prime}\right)  $, this rewrites as
\[
\left\{  \operatorname{st} \tau   \ \mid\ \tau\in S_{\prec
}\left(  \pi,\sigma\right)  \right\}  _{\operatorname*{multi}}=\left\{
\operatorname{st} \tau   \ \mid\ \tau\in S_{\prec}\left(
\pi^{\prime},\sigma^{\prime}\right)  \right\}  _{\operatorname*{multi}}.
\]

Also, Claim 1 (applied to $\sigma$, $\pi$, $\sigma^{\prime}$ and $\pi^{\prime
}$ instead of $\pi$, $\sigma$, $\pi^{\prime}$ and $\sigma^{\prime}$) yields%
\[
\left\{  \operatorname{st} \tau   \ \mid\ \tau\in S_{\prec
}\left(  \sigma,\pi\right)  \right\}  _{\operatorname*{multi}}=\left\{
\operatorname{st} \tau   \ \mid\ \tau\in S_{\prec}\left(
\sigma^{\prime},\pi^{\prime}\right)  \right\}  _{\operatorname*{multi}}.
\]
In light of $S_{\prec}\left(  \sigma,\pi\right)  =S_{\succ}\left(  \pi
,\sigma\right)  $ and $S_{\prec}\left(  \sigma^{\prime},\pi^{\prime}\right)
=S_{\succ}\left(  \pi^{\prime},\sigma^{\prime}\right)  $, this rewrites as%
\[
\left\{  \operatorname{st} \tau   \ \mid\ \tau\in S_{\succ
}\left(  \pi,\sigma\right)  \right\}  _{\operatorname*{multi}}=\left\{
\operatorname{st} \tau   \ \mid\ \tau\in S_{\succ}\left(
\pi^{\prime},\sigma^{\prime}\right)  \right\}  _{\operatorname*{multi}}.
\]
Thus, Claim 3 is proven in Case 2.

We have now proven Claim 3 in both Cases 1 and 2. Thus, Claim 3 is always proven.]

Claim 3 shows that $\operatorname{st}$ is LR-shuffle-compatible. This proves
Corollary \ref{cor.LRcomp.two}.
\end{proof}
\end{verlong}

\section{\label{sect.Descent}Descent statistics and quasisymmetric functions}

In this section, we shall recall the concepts of descent statistics and their
shuffle algebras (introduced in \cite{part1}), and apply them to
$\operatorname{Epk}$.

\subsection{Compositions}

\begin{definition}
A \textit{composition} is a finite list of positive integers. If $I=\left(
i_{1},i_{2},\ldots,i_{n}\right)  $ is a composition, then the nonnegative
integer $i_{1}+i_{2}+\cdots+i_{n}$ is called the \textit{size} of $I$ and is
denoted by $\left\vert I\right\vert $; we furthermore say that $I$ is a
\textit{composition of }$\left\vert I\right\vert $.
\end{definition}

\begin{definition}
\label{def.comps-to-sets}Let $n\in\mathbb{N}$. For each composition $I=\left(
i_{1},i_{2},\ldots,i_{k}\right)  $ of $n$, we define a subset
$\operatorname{Des}I$ of $\left[  n-1\right]  $ by%
\begin{align*}
\operatorname{Des}I  &  =\left\{  i_{1},i_{1}+i_{2},i_{1}+i_{2}+i_{3}%
,\ldots,i_{1}+i_{2}+\cdots+i_{k-1}\right\} \\
&  =\left\{  i_{1}+i_{2}+\cdots+i_{s}\ \mid\ s\in\left[  k-1\right]  \right\}
.
\end{align*}

On the other hand, for each subset $A=\left\{  a_{1}<a_{2}<\cdots
<a_{k}\right\}  $ of $\left[  n-1\right]  $, we define a composition
$\operatorname*{Comp}A$ of $n$ by%
\[
\operatorname*{Comp}A=\left(  a_{1},a_{2}-a_{1},a_{3}-a_{2},\ldots
,a_{k}-a_{k-1},n-a_{k}\right)  .
\]
(The definition of $\operatorname*{Comp}A$ should be understood to give
$\operatorname*{Comp}A=\left(  n\right)  $ if $A=\varnothing$ and $n > 0$,
and to give
$\operatorname*{Comp}A=\left( \right)  $ if $A=\varnothing$ and $n = 0$.
Note that
$\operatorname*{Comp}A$ depends not only on the set $A$ itself, but also on
$n$. We hope that $n$ will always be clear from the context when we use this notation.)

We thus have defined a map $\operatorname{Des}$ (from the set of all
compositions of $n$ to the set of all subsets of $\left[  n-1\right]  $) and a
map $\operatorname*{Comp}$ (in the opposite direction). These two maps are
mutually inverse bijections.
\end{definition}

\begin{definition}
\label{def.CompDes}Let $n\in\mathbb{N}$. Let $\pi=\left(  \pi_{1},\pi
_{2},\ldots,\pi_{n}\right)  $ be an $n$-permutation. The \textit{descent
composition} of $\pi$ is defined to be the composition $\operatorname*{Comp}%
\left(  \operatorname{Des}\pi\right)  $ of $n$. This composition is denoted
by $\operatorname*{Comp}\pi$.
\end{definition}

For example, the $6$-permutation $\pi=\left(  4,1,3,9,6,8\right)  $ has
$\operatorname*{Comp}\pi=\left(  1,3,2\right)  $. For another example, the
$6$-permutation $\pi=\left(  1,4,3,2,9,8\right)  $ has $\operatorname*{Comp}%
\pi=\left(  2,1,2,1\right)  $.

Definition \ref{def.CompDes} defines the permutation statistic
$\operatorname*{Comp}$, whose codomain is the set of all compositions.

\subsection{Descent statistics}

\begin{definition}
Let $\operatorname{st}$ be a permutation statistic. We say that
$\operatorname{st}$ is a \textit{descent statistic} if and only if
$\operatorname{st}\pi$ (for $\pi$ a permutation) depends only on the descent
composition $\operatorname*{Comp}\pi$ of $\pi$. In other words,
$\operatorname{st}$ is a descent statistic if and only if every two
permutations $\pi$ and $\sigma$ satisfying $\operatorname*{Comp}%
\pi=\operatorname*{Comp}\sigma$ satisfy $\operatorname{st}\pi
=\operatorname{st}\sigma$.
\end{definition}

Equivalently, a permutation statistic $\operatorname{st}$ is a descent
statistic if and only if every two permutations $\pi$ and $\sigma$ satisfying
$\left\vert \pi\right\vert =\left\vert \sigma\right\vert $ and
$\operatorname{Des}\pi=\operatorname{Des}\sigma$ satisfy $\operatorname{st}
\pi=\operatorname{st}\sigma$. (This is indeed equivalent, because for two
permutations $\pi$ and $\sigma$, the condition $\left(  \left\vert
\pi\right\vert =\left\vert \sigma\right\vert \text{ and }\operatorname{Des}
\pi=\operatorname{Des}\sigma\right)  $ is equivalent to $\left(
\operatorname*{Comp}\pi=\operatorname*{Comp}\sigma\right)  $.)

For example, the permutation statistic $\operatorname{Des}$ is a descent
statistic, because each permutation $\pi$ satisfies $\operatorname{Des}
\pi=\operatorname{Des}\left(  \operatorname*{Comp}\pi\right)  $. Also,
$\operatorname*{Pk}$ is a descent statistic, since each permutation $\pi$
satisfies%
\[
\operatorname*{Pk}\pi=\left(  \operatorname{Des}\pi\right)  \setminus\left(
\left\{  1\right\}  \cup\left(  \operatorname{Des}\pi+1\right)  \right)  ,
\]
where $\operatorname{Des}\pi+1$ denotes the set $\left\{  i+1\ \mid
\ i\in\operatorname{Des}\pi\right\}  $ (and, as we have just said,
$\operatorname{Des}\pi$ can be recovered from $\operatorname*{Comp}\pi$).
Furthermore, $\operatorname{Epk}$ is a descent statistic, since each
$n$-permutation $\pi$ (for a positive integer $n$) satisfies%
\[
\operatorname{Epk}\pi=\left(  \operatorname{Des}\pi\cup\left\{  n\right\}
\right)  \setminus\left(  \operatorname{Des}\pi+1\right)
\]
(and both $\operatorname{Des}\pi$ and $n$ can be recovered from
$\operatorname*{Comp}\pi$). The permutation statistics $\operatorname*{Lpk}$
and $\operatorname*{Rpk}$ (and, of course, $\operatorname*{Comp}$) are descent
statistics as well, as one can easily check.

In \cite[Corollary 1.6]{Oguz18}, Ezgi Kantarc{\i} O\u{g}uz has demonstrated that
not every shuffle-compatible permutation statistic is a descent statistic.
However, this changes if we require LR-shuffle-compatibility, because of
Corollary \ref{cor.LRcomp.back} and of the following fact:

\begin{proposition}
\label{prop.des-stat.hgc}Every head-graft-compatible permutation statistic is
a descent statistic.
\end{proposition}

\begin{vershort}
\begin{proof}
[Proof of Proposition \ref{prop.des-stat.hgc}.]See \cite{verlong}.
\end{proof}
\end{vershort}

\begin{verlong}
\begin{proof}
[Proof of Proposition \ref{prop.des-stat.hgc}.]Let $\operatorname{st}$ be a
head-graft-compatible permutation statistic. Then, we must show that
$\operatorname{st}$ is a descent statistic. But this follows directly from
Lemma \ref{lem.LRcomp.head-l1}.
\end{proof}

Combining Proposition \ref{prop.des-stat.hgc} with Corollary
\ref{cor.LRcomp.back}, we conclude that every LR-shuffle-compatible
permutation statistic is a descent statistic.
\end{verlong}

\begin{definition}
\label{def.des-stat.stI}Let $\operatorname{st}$ be a descent statistic. Then,
we can regard $\operatorname{st}$ as a map from the set of all compositions
(rather than from the set of all permutations). Namely, for any composition
$I$, we define $\operatorname{st}I$ (an element of the codomain of
$\operatorname{st}$) by setting%
\[
\operatorname{st}I=\operatorname{st}\pi\ \ \ \ \ \ \ \ \ \ \text{for any
permutation }\pi\text{ satisfying }\operatorname*{Comp}\pi=I.
\]
This is well-defined (because for every composition $I$, there exists at least
one permutation $\pi$ satisfying $\operatorname*{Comp}\pi=I$, and all such
permutations $\pi$ have the same value of $\operatorname{st}\pi$). In the
following, we shall regard every descent statistic $\operatorname{st}$
simultaneously as a map from the set of all permutations and as a map from the
set of all compositions.
\end{definition}

Note that this definition leads to a new interpretation of
$\operatorname{Des}I$ for a composition $I$: It is now defined as
$\operatorname{Des}\pi$ for any permutation $\pi$ satisfying
$\operatorname*{Comp}\pi=I$. This could clash with the old meaning of
$\operatorname{Des}I$ introduced in Definition \ref{def.comps-to-sets}.
Fortunately, these two meanings of $\operatorname{Des}I$ are exactly the
same, so there is no conflict of notation.

However, Definition \ref{def.des-stat.stI} causes an ambiguity for expressions
like \textquotedblleft$\operatorname{Des}\left(  i_{1},i_{2},\ldots
,i_{n}\right)  $\textquotedblright: Here, the \textquotedblleft$\left(
i_{1},i_{2},\ldots,i_{n}\right)  $\textquotedblright\ might be understood
either as a permutation, or as a composition, and the resulting descent sets
$\operatorname{Des}\left(  i_{1},i_{2},\ldots,i_{n}\right)  $ are not the
same. A similar ambiguity occurs for any descent statistic
$\operatorname{st}$ instead of $\operatorname{Des}$.
We hope that this ambiguity will not arise
in this paper due to our explicit typecasting of permutations and
compositions; but the reader should be warned that it can arise if one takes
the notation too literally.

\begin{definition}
\label{def.des-stat.eq-comp}Let $\operatorname{st}$ be a descent statistic.

\textbf{(a)} Two compositions $J$ and $K$ are said to be
$\operatorname{st}$\textit{-equivalent} if and only if they
have the same size and satisfy
$\operatorname{st}J=\operatorname{st}K$. Equivalently, two compositions $J$
and $K$ are $\operatorname{st}$-equivalent if and only if there exist two
$\operatorname{st}$-equivalent permutations $\pi$ and $\sigma$ satisfying
$J=\operatorname*{Comp}\pi$ and $K=\operatorname*{Comp}\sigma$.

\textbf{(b)} The relation \textquotedblleft$\operatorname{st}$%
-equivalent\textquotedblright\ is an equivalence relation on compositions; its
equivalence classes are called $\operatorname{st}$\textit{-equivalence
classes of compositions}.
\end{definition}

\subsection{Quasisymmetric functions}

We now recall the definition of quasisymmetric functions; see \cite[Chapter
5]{HopfComb} (and various other modern textbooks) for more details about this:

\begin{definition}
\textbf{(a)}
Consider the ring of power series $\mathbb{Q}\left[  \left[  x_{1},x_{2}%
,x_{3},\ldots\right]  \right]  $ in infinitely many commuting indeterminates
over $\mathbb{Q}$.
A power series $f\in\mathbb{Q}\left[  \left[  x_{1}%
,x_{2},x_{3},\ldots\right]  \right]  $ is said to be \textit{quasisymmetric}
if it has the following property:
\begin{itemize}
\item For any positive integers $a_{1}%
,a_{2},\ldots,a_{k}$ and any two strictly increasing sequences $\left(
i_{1}<i_{2}<\cdots<i_{k}\right)  $ and $\left(  j_{1}<j_{2}<\cdots
<j_{k}\right)  $ of positive integers, the coefficient of $x_{i_{1}}^{a_{1}%
}x_{i_{2}}^{a_{2}}\cdots x_{i_{k}}^{a_{k}}$ in $f$ equals the coefficient of
$x_{j_{1}}^{a_{1}}x_{j_{2}}^{a_{2}}\cdots x_{j_{k}}^{a_{k}}$ in $f$.
\end{itemize}

\textbf{(b)}
A \textit{quasisymmetric function} is a quasisymmetric power series
$f\in\mathbb{Q}\left[  \left[  x_{1},x_{2},x_{3},\ldots\right]  \right]  $
that has bounded degree (i.e., there exists an $N\in\mathbb{N}$ such that each
monomial appearing in $f$ has degree $\leq N$).

\textbf{(c)} The quasisymmetric functions
form a $\mathbb{Q}$-subalgebra of $\mathbb{Q}\left[  \left[  x_{1},x_{2}%
,x_{3},\ldots\right]  \right]  $; this $\mathbb{Q}$-subalgebra is denoted by
$\operatorname*{QSym}$ and called the \textit{ring of quasisymmetric
functions} over $\mathbb{Q}$.
This $\mathbb{Q}$-algebra $\operatorname{QSym}$ is graded (in the obvious way,
i.e., by the degree of a monomial).
\end{definition}

The $\mathbb{Q}$-algebra $\operatorname*{QSym}$ has much interesting structure
(e.g., it is a Hopf algebra), some of which we will introduce later when we
need it. One simple yet crucial feature of $\operatorname*{QSym}$ that we will
immediately use is the \textit{fundamental basis} of $\operatorname*{QSym}$:

\begin{definition}
For any composition $\alpha$, we define the \textit{fundamental quasisymmetric
function} $F_{\alpha}$ to be the power series%
\[
\sum_{\substack{i_{1}\leq i_{2}\leq\cdots\leq i_{n};\\i_{j}<i_{j+1}\text{ for
each }j\in\operatorname{Des}\alpha}}x_{i_{1}}x_{i_{2}}\cdots x_{i_{n}}%
\in\operatorname*{QSym},
\]
where $n=\left\vert \alpha\right\vert $ is the size of $\alpha$. The family
$\left(  F_{\alpha}\right)  _{\alpha\text{ is a composition}}$ is a basis of
the $\mathbb{Q}$-vector space $\operatorname*{QSym}$; it is known as the
\textit{fundamental basis} of $\operatorname*{QSym}$.
\end{definition}

The fundamental quasisymmetric function $F_{\alpha}$ is denoted
by $L_{\alpha}$ in \cite[\S 5.2]{HopfComb}.

The multiplication of fundamental quasisymmetric functions is intimately
related to shuffles of permutations:

\begin{verlong}
\begin{theorem}
\label{thm.4.1}Let $\pi$ and $\sigma$ be two disjoint permutations. Let
$J=\operatorname*{Comp}\pi$ and $K=\operatorname*{Comp}\sigma$. For any
composition $L$, let $c_{J,K}^{L}$ be the number of permutations with descent
composition $L$ among the shuffles of $\pi$ and $\sigma$. Then,%
\[
F_{J}F_{K}=\sum_{L}c_{J,K}^{L}F_{L}%
\]
(where the sum is over all compositions $L$).
\end{theorem}

Theorem \ref{thm.4.1} is \cite[Theorem 4.1]{part1}; it can also be written in
the following form:
\end{verlong}

\begin{proposition}
\label{prop.4.1.rewr}Let $\pi$ and $\sigma$ be two disjoint permutations.
Then,%
\[
F_{\operatorname*{Comp}\pi}F_{\operatorname*{Comp}\sigma}=\sum_{\chi\in
S\left(  \pi,\sigma\right)  }F_{\operatorname*{Comp}\chi}.
\]

\end{proposition}

\begin{vershort}
Proposition \ref{prop.4.1.rewr} is a restatement of \cite[Theorem 4.1]{part1},
and is proven in \cite[(5.2.6)]{HopfComb} (which makes the additional
requirement that the letters of $\pi$ are $1,2,\ldots,\left\vert
\pi\right\vert $ and the letters of $\sigma$ are $\left\vert \pi\right\vert
+1,\left\vert \pi\right\vert +2,\ldots,\left\vert \pi\right\vert +\left\vert
\sigma\right\vert $; but this requirement is not used in the proof and thus
can be dropped).
\end{vershort}

\begin{verlong}
For a proof of Proposition \ref{prop.4.1.rewr} (and therefore also of the
equivalent Theorem \ref{thm.4.1}), we refer to \cite[(5.2.6)]{HopfComb} (which
makes the additional requirement that the letters of $\pi$ are $1,2,\ldots
,\left\vert \pi\right\vert $ and the letters of $\sigma$ are $\left\vert
\pi\right\vert +1,\left\vert \pi\right\vert +2,\ldots,\left\vert
\pi\right\vert +\left\vert \sigma\right\vert $; but this requirement is not
used in the proof and thus can be dropped).
\end{verlong}

\subsection{Shuffle algebras}

Any shuffle-compatible permutation statistic $\operatorname{st}$ gives rise
to a \textit{shuffle algebra} $\mathcal{A}_{\operatorname{st}}$, defined as follows:

\begin{definition}
\label{def.Ast}Let $\operatorname{st}$ be a shuffle-compatible permutation
statistic. For each permutation $\pi$, let $\left[  \pi\right]
_{\operatorname{st}}$ denote the $\operatorname{st}$-equivalence class of
$\pi$.

Let $\mathcal{A}_{\operatorname{st}}$ be the free $\mathbb{Q}$-vector space
whose basis is the set of all $\operatorname{st}$-equivalence classes of
permutations. We define a multiplication on $\mathcal{A}_{\operatorname{st}}$
by setting%
\[
\left[  \pi\right]  _{\operatorname{st}}\left[  \sigma\right]
_{\operatorname{st}}=\sum_{\tau\in S\left(  \pi,\sigma\right)  }\left[
\tau\right]  _{\operatorname{st}}%
\]
for any two disjoint permutations $\pi$ and $\sigma$. It is easy to see that
this multiplication is well-defined and associative, and turns $\mathcal{A}%
_{\operatorname{st}}$ into a $\mathbb{Q}$-algebra whose unity is the
$\operatorname{st}$-equivalence class of the $0$-permutation
$\left( \right)$.
\begin{verlong}
(In the
particular case when $\operatorname{st}$ is a descent statistic, this shall
be proven again in Proposition \ref{prop.Ast.alg} \textbf{(a)} below.)
\end{verlong}
This
$\mathbb{Q}$-algebra is denoted by $\mathcal{A}_{\operatorname{st}}$, and is
called the \textit{shuffle algebra} of $\operatorname{st}$. It is a graded
$\mathbb{Q}$-algebra; its $n$-th graded component (for each $n\in\mathbb{N}$)
is spanned by the $\operatorname{st}$-equivalence classes of all $n$-permutations.
\end{definition}

This definition originates in \cite[\S 3.1]{part1}.
The following fact is implicit in \cite{part1}:

\begin{proposition}
\label{prop.Ast.alg}Let $\operatorname{st}$ be a shuffle-compatible descent statistic.

\begin{verlong}
\textbf{(a)} The multiplication on $\mathcal{A}_{\operatorname{st}}$ defined
in Definition \ref{def.Ast} is well-defined and associative, and turns
$\mathcal{A}_{\operatorname{st}}$ into a $\mathbb{Q}$-algebra whose unity is
the $\operatorname{st}$-equivalence class of the $0$-permutation
$\left( \right)$.

\textbf{(b)}
\end{verlong}
There is a surjective $\mathbb{Q}$-algebra homomorphism
$p_{\operatorname{st}}:\operatorname*{QSym}\rightarrow\mathcal{A}%
_{\operatorname{st}}$ that satisfies
\[
p_{\operatorname{st}}\left(  F_{\operatorname*{Comp}\pi}\right)  =\left[
\pi\right]  _{\operatorname{st}}\ \ \ \ \ \ \ \ \ \ \text{for every
permutation }\pi.
\]

\end{proposition}

A central result, connecting shuffle-compatibility of a descent
statistic with $\operatorname*{QSym}$, is \cite[Theorem 4.3]{part1},
which we restate as follows:

\begin{theorem}
\label{thm.4.3}Let $\operatorname{st}$ be a descent statistic.

\textbf{(a)} The descent statistic $\operatorname{st}$ is shuffle-compatible
if and only if there exist a $\mathbb{Q}$-algebra $A$ with basis $\left(
u_{\alpha}\right)  $ (indexed by $\operatorname{st}$-equivalence classes
$\alpha$ of compositions) and a $\mathbb{Q}$-algebra homomorphism
$\phi_{\operatorname{st}}:\operatorname*{QSym}\rightarrow A$ with the
property that whenever $\alpha$ is an $\operatorname{st}$-equivalence class
of compositions, we have%
\[
\phi_{\operatorname{st}}\left(  F_{L}\right)  =u_{\alpha}%
\ \ \ \ \ \ \ \ \ \ \text{for each }L\in\alpha.
\]

\textbf{(b)} In this case, the $\mathbb{Q}$-linear map%
\[
\mathcal{A}_{\operatorname{st}}\rightarrow A,\ \ \ \ \ \ \ \ \ \ \left[
\pi\right]  _{\operatorname{st}}\mapsto u_{\alpha},
\]
where $\alpha$ is the $\operatorname{st}$-equivalence class of the
composition $\operatorname*{Comp}\pi$, is a $\mathbb{Q}$-algebra isomorphism
$\mathcal{A}_{\operatorname{st}}\rightarrow A$.
\end{theorem}

\begin{vershort}
Proofs of Proposition \ref{prop.Ast.alg} and Theorem \ref{thm.4.3}
(independent of \cite{part1}) can be found in \cite{verlong}.
\end{vershort}

\begin{verlong}
As we said, Theorem \ref{thm.4.3} is precisely \cite[Theorem 4.3]{part1}.
For the sake of completeness, we shall give proofs of Proposition
\ref{prop.Ast.alg} and Theorem \ref{thm.4.3} (independent of \cite{part1}) in
Subsection \ref{subsect.K.pfAsh}.
\end{verlong}

\begin{noncompile}
We shall perhaps also use some other notations that were introduced in
\cite{part1}. The reader should thus consult \cite{part1} for unexplained notations.
\end{noncompile}

\subsection{The shuffle algebra of $\operatorname{Epk}$}

Theorem \ref{thm.Epk.sh-co-a} yields that the permutation statistic
$\operatorname{Epk}$ is shuffle-compatible. Hence, the shuffle algebra
$\mathcal{A}_{\operatorname{Epk}}$ is well-defined. We have little to say
about it:

\begin{theorem}
\label{thm.Epk.AEpk}\textbf{(a)} The shuffle algebra $\mathcal{A}%
_{\operatorname{Epk}}$ is a graded quotient algebra of $\operatorname*{QSym}$.

\textbf{(b)} Define the Fibonacci sequence
$\left(  f_{0},f_{1},f_{2},\ldots\right)  $ as in
Proposition~\ref{prop.lac.fib}.
Let $n$ be a positive integer. The $n$-th graded component of
$\mathcal{A}_{\operatorname{Epk}}$ has dimension $f_{n+2}-1$.
\end{theorem}

\begin{vershort}
\begin{proof}
[Proof of Theorem \ref{thm.Epk.AEpk}.]See \cite{verlong}.
\end{proof}
\end{vershort}

\begin{verlong}
\begin{proof}
[Proof of Theorem \ref{thm.Epk.AEpk} (sketched).]\textbf{(a)} Proposition
\ref{prop.Ast.alg} \textbf{(b)} (applied to $\operatorname{st}
=\operatorname{Epk}$) yields a surjective $\mathbb{Q}$-algebra homomorphism
$p_{\operatorname{Epk}}:\operatorname*{QSym}\rightarrow\mathcal{A}%
_{\operatorname{Epk}}$. It is easy to see that this $p_{\operatorname{Epk}}$
is furthermore graded (i.e., degree-preserving). Thus, $\mathcal{A}%
_{\operatorname{Epk}}$ is isomorphic to a quotient of $\operatorname*{QSym}$
as a graded algebra. This proves Theorem \ref{thm.Epk.AEpk} \textbf{(a)}.

\textbf{(b)} The $n$-th graded component of
$\mathcal{A}_{\operatorname{Epk}}$ has a basis
indexed by $\operatorname{Epk}$-equivalence classes of
compositions of $n$. These latter classes are in bijection with
$\operatorname{Epk}$-equivalence classes of $n$-permutations. In turn, the
latter classes are in bijection with the images of $n$-permutations under the
map $\operatorname{Epk}$. Finally, the latter images are the elements of
$\mathbf{L}_n$ (according to Proposition \ref{prop.when-Epk}).
Hence, the number of these images is $\left|\mathbf{L}_n\right|
= f_{n+2} - 1$ (by Proposition \ref{prop.lac.fib}). Combining all of the
preceding sentences, we conclude that the dimension of the $n$-th graded
component of $\mathcal{A}_{\operatorname{Epk}}$ is $f_{n+2}-1$. This proves
Theorem \ref{thm.Epk.AEpk} \textbf{(b)}.
\end{proof}
\end{verlong}

We can describe $\mathcal{A}_{\operatorname{Epk}}$ using the notations
of Section \ref{sect.Zenri}:

\begin{definition}
Let $\Pi_{\mathcal{Z}}$ be the $\mathbb{Q}$-vector subspace of
$\operatorname*{Pow}\mathcal{N}$ spanned by the family $\left(  K_{n,\Lambda
}^{\mathcal{Z}}\right)  _{n \in \NN;\  \Lambda \in \mathbf{L}_n}  $.
Then, $\Pi_{\mathcal{Z}}$ is also the $\mathbb{Q}$-vector subspace of
$\operatorname*{Pow}\mathcal{N}$ spanned by the family $\left(
K_{n,\operatorname{Epk}\pi}^{\mathcal{Z}}\right)  _{n\in\mathbb{N}%
;\ \pi\text{ is an }n\text{-permutation}}$ (by Proposition \ref{prop.when-Epk}%
). In other words, $\Pi_{\mathcal{Z}}$ is also the $\mathbb{Q}$-vector
subspace of $\operatorname*{Pow}\mathcal{N}$ spanned by the family $\left(
\Gamma_{\mathcal{Z}}\left(  \pi\right)  \right)  _{n\in\mathbb{N};\ \pi\text{
is an }n\text{-permutation}}$ (because of (\ref{eq.def.KnL.2})). Hence,
Corollary \ref{cor.prod2} shows that $\Pi_{\mathcal{Z}}$ is closed under
multiplication. Since furthermore $\Gamma_{\mathcal{Z}}\left(
\left( \right) \right)  =1$ (for the $0$-permutation
$\left( \right)$), we can thus conclude that $\Pi
_{\mathcal{Z}}$ is a $\mathbb{Q}$-subalgebra of $\operatorname*{Pow}%
\mathcal{N}$.
\end{definition}

\begin{theorem}
\label{thm.Epk.sh-co}The $\mathbb{Q}$-linear map
\[
\mathcal{A}_{\operatorname{Epk}} \to \Pi_{\mathcal{Z}}, \qquad
\left[  \pi\right]  _{\operatorname{Epk}}\mapsto
K_{n,\operatorname{Epk} \pi}^{\mathcal{Z}}
\]
is a $\mathbb{Q}$-algebra isomorphism.
\end{theorem}

\begin{vershort}
\begin{proof}
[Proof of Theorem \ref{thm.Epk.sh-co}.]See \cite{verlong}.
\end{proof}
\end{vershort}

\begin{verlong}
In the process of proving Theorem \ref{thm.Epk.sh-co}, we will also prove
Theorem \ref{thm.Epk.sh-co-a} again.

\begin{proof}
[Proof of Theorem \ref{thm.Epk.sh-co} (sketched).]For each positive integer
$n$ and each $n$-permutation $\pi$, we have%
\[
\operatorname{Epk}\pi=\left(  \operatorname{Des}\pi\cup\left\{  n\right\}
\right)  \setminus\left(  \operatorname{Des}\pi+1\right)
\]
(by Proposition \ref{prop.Epk.through-Des}). Thus, $\operatorname{Epk}\pi$ is
uniquely determined by $\operatorname{Des}\pi$ and $n$. (Of course, this
holds for $n=0$ as well, because in this case only one $\pi$ exists.) Hence,
$\operatorname{Epk}$ is a descent statistic. Thus, for every composition $L$,
a set $\operatorname{Epk}L$ is defined (according to Definition
\ref{def.des-stat.stI}); explicitly,
$\operatorname{Epk}L=\operatorname{Epk}\pi$ whenever
$\pi$ is a permutation satisfying $\operatorname*{Comp}\pi=L$.

\begin{noncompile}
We can spell out this definition: For each positive integer $n$ and each
composition $L$ of $n$, the subset $\operatorname{Epk}L$ of $\left[
n\right]  $ is given by%
\[
\operatorname{Epk}L=\left(  \operatorname{Des}L\cup\left\{  n\right\}
\right)  \setminus\left(  \operatorname{Des}L+1\right)  .
\]

\end{noncompile}

Recall that $\left(  F_{L}\right)  _{L\text{ is a composition}}$ is a basis of
the $\mathbb{Q}$-vector space $\operatorname*{QSym}$. Let $\phi
_{\operatorname{Epk}}:\operatorname*{QSym}\rightarrow\Pi_{\mathcal{Z}}$ be
the $\mathbb{Q}$-linear map that sends each $F_{L}$ (for each $n\in\mathbb{N}$
and each composition $L$ of $n$) to $K_{n,\operatorname{Epk}L}^{\mathcal{Z}%
}\in\Pi_{\mathcal{Z}}$. This $\mathbb{Q}$-linear map $\phi
_{\operatorname{Epk}}$ respects multiplication\footnote{\textit{Proof.} Let
$n\in\mathbb{N}$ and $m\in\mathbb{N}$. Let $J$ be a composition of $n$, and
let $K$ be a composition of $m$. Fix an $n$-permutation $\pi$ satisfying
$\operatorname*{Comp}\pi=J$, and fix an $m$-permutation $\sigma$ satisfying
$\operatorname*{Comp}\sigma=K$ such that $\pi$ and $\sigma$ are disjoint.
Now,
\begin{align*}
&  \underbrace{\phi_{\operatorname{Epk}}\left(  F_{J}\right)  }%
_{\substack{=K_{n,\operatorname{Epk}J}^{\mathcal{Z}}
=K_{n,\operatorname{Epk}\pi}^{\mathcal{Z}}\\
\text{(since }\operatorname{Epk}J=\operatorname{Epk} \pi\text{)}}}
\cdot\underbrace{\phi_{\operatorname{Epk}}\left(  F_{K}\right)
}_{\substack{=K_{m,\operatorname{Epk}K}^{\mathcal{Z}}%
=K_{m,\operatorname{Epk}\sigma}^{\mathcal{Z}}\\\text{(since }%
\operatorname{Epk}K=\operatorname{Epk}\sigma\text{)}}}\\
&  =K_{n,\operatorname{Epk}\pi}^{\mathcal{Z}}\cdot
K_{m,\operatorname{Epk}\sigma}^{\mathcal{Z}}
=\sum_{\tau\in S\left(  \pi,\sigma\right)  }%
K_{n+m,\operatorname{Epk}\tau}^{\mathcal{Z}}\ \ \ \ \ \ \ \ \ \ \left(
\text{by Corollary \ref{cor.KnEpk-prodrule}}\right)  .
\end{align*}
Comparing this with%
\begin{align*}
\phi_{\operatorname{Epk}}\left(  \underbrace{F_{J}F_{K}}_{\substack{=\sum
_{\tau\in S\left(  \pi,\sigma\right)  }F_{\operatorname*{Comp}\tau}\\\text{(by
Proposition~\ref{prop.4.1.rewr})}}}\right)   &  =\phi_{\operatorname{Epk}}
\left(  \sum_{\tau\in S\left(  \pi,\sigma\right)  }F_{\operatorname*{Comp}
\tau}\right)  =\sum_{\tau\in S\left(  \pi,\sigma\right)  }\underbrace{\phi
_{\operatorname{Epk}}\left(  F_{\operatorname*{Comp}\tau}\right)
}_{\substack{=K_{n+m,\operatorname{Epk}\left(  \operatorname*{Comp}%
\tau\right)  }^{\mathcal{Z}}=K_{n+m,\operatorname{Epk}\tau}^{\mathcal{Z}%
}\\\text{(since }\operatorname{Epk}\left(  \operatorname*{Comp}\tau\right)
=\operatorname{Epk}\tau\text{)}}}\\
&  =\sum_{\tau\in S\left(  \pi,\sigma\right)  }K_{n+m,\operatorname{Epk}\tau
}^{\mathcal{Z}},
\end{align*}
we obtain $\phi_{\operatorname{Epk}}\left(  F_{J}\right)  \cdot
\phi_{\operatorname{Epk}}\left(  F_{K}\right)  =\phi_{\operatorname{Epk}}
\left(  F_{J}F_{K}\right)  $. Since the map $\phi_{\operatorname{Epk}}$ is
$\mathbb{Q}$-linear, this yields that $\phi_{\operatorname{Epk}}$ respects
multiplication (since $\left(  F_{L}\right)  _{L\text{ is a composition}}$ is
a basis of the $\mathbb{Q}$-vector space $\operatorname*{QSym}$).} and sends
$1\in\operatorname*{QSym}$ to $1\in\Pi_{\mathcal{Z}}$\ \ \ \ \footnote{This is
easy, since $1=F_{\left(  {}\right)  }$.}. Thus, $\phi_{\operatorname{Epk}}$
is a $\mathbb{Q}$-algebra homomorphism.

The family $\left(  K_{n,\Lambda}^{\mathcal{Z}}\right) %
_{n \in \NN;\  \Lambda \in \mathbf{L}_n}$ spans $\Pi_{\mathcal{Z}}$ (by the
definition of $\Pi_{\mathcal{Z}}$) and is $\mathbb{Q}$-linearly independent
(by Corollary \ref{cor.KnL.lindep-all}). Thus, it is a basis of $\Pi
_{\mathcal{Z}}$.

For each positive integer $n$, there is a canonical bijection between the
$\operatorname{Epk}$-equivalence classes of $n$-permutations and the
$\Lambda \in \mathbf{L}_n$ (indeed, the bijection sends
any equivalence class $\left[  \pi\right]  _{\operatorname{Epk}}$ to
$\operatorname{Epk}\pi$)\ \ \ \ \footnote{Proposition \ref{prop.when-Epk}
shows that this is indeed a bijection.}. Hence, Theorem~\ref{thm.4.3}
\textbf{(a)} (applied to $A=\Pi_{\mathcal{Z}}$, $\operatorname{st}
=\operatorname{Epk}$ and $u_{\alpha}=K_{n,\Lambda}^{\mathcal{Z}}$, where
$\alpha$ is an $\operatorname{Epk}$-equivalence class of $n$-permutations and
where $\Lambda$ is the corresponding element of $\mathbf{L}_n$)
shows that the descent statistic $\operatorname{Epk}$ is
shuffle-compatible. This proves Theorem \ref{thm.Epk.sh-co-a} again.
Theorem~\ref{thm.4.3} \textbf{(b)} then yields that the $\mathbb{Q}$-linear
map%
\[
\mathcal{A}_{\operatorname{Epk}}\rightarrow\Pi_{\mathcal{Z}}%
,\ \ \ \ \ \ \ \ \ \ \left[  \pi\right]  _{\operatorname{Epk}}\mapsto
K_{n,\operatorname{Epk}\pi}^{\mathcal{Z}}%
\]
is a $\mathbb{Q}$-algebra isomorphism from $\mathcal{A}_{\operatorname{Epk}}$
to $\Pi_{\mathcal{Z}}$. Hence, Theorem \ref{thm.Epk.sh-co} is proven.
\end{proof}
\end{verlong}

\section{\label{sect.kernel}The kernel of the map $\operatorname*{QSym}%
\rightarrow\mathcal{A}_{\operatorname{Epk}}$}

\subsection{The kernel of a descent statistic}

Now, we shall focus on a feature of shuffle-compatible descent statistics that
seems to have been overlooked so far: their kernels.

\begin{vershort}
All proofs in this section are omitted; they can be found in \cite{verlong}.
\end{vershort}

\begin{definition}
Let $\operatorname{st}$ be a descent statistic. Then, $\mathcal{K}%
_{\operatorname{st}}$ shall mean the $\mathbb{Q}$-vector subspace of
$\operatorname*{QSym}$ spanned by all elements of the form $F_{J}-F_{K}$,
where $J$ and $K$ are two $\operatorname{st}$-equivalent compositions. (See
Definition \ref{def.des-stat.eq-comp} \textbf{(a)} for the definition of
\textquotedblleft$\operatorname{st}$-equivalent
compositions\textquotedblright.) We shall refer to $\mathcal{K}%
_{\operatorname{st}}$ as the \textit{kernel} of $\operatorname{st}$.
\end{definition}

The following basic linear-algebraic lemma will be useful:

\begin{lemma}
\label{lem.K.fist}Let $\operatorname{st}$ be a descent statistic. Let $A$ be
a $\mathbb{Q}$-vector space with basis $\left(  u_{\alpha}\right)  $ indexed
by $\operatorname{st}$-equivalence classes $\alpha$ of compositions. Let
$\phi_{\operatorname{st}}:\operatorname*{QSym}\rightarrow A$ be a
$\mathbb{Q}$-linear map with the property that whenever $\alpha$ is an
$\operatorname{st}$-equivalence class of compositions, we have%
\begin{equation}
\phi_{\operatorname{st}}\left(  F_{L}\right)  =u_{\alpha}%
\ \ \ \ \ \ \ \ \ \ \text{for each }L\in\alpha. \label{pf.prop.K.ideal.dir1.1}%
\end{equation}

Then, $\operatorname*{Ker}\left(  \phi_{\operatorname{st}}\right)
=\mathcal{K}_{\operatorname{st}}$.
\end{lemma}

\begin{verlong}
\begin{proof}
[Proof of Lemma \ref{lem.K.fist}.]Let us first show that $\operatorname*{Ker}%
\left(  \phi_{\operatorname{st}}\right)  \subseteq\mathcal{K}%
_{\operatorname{st}}$.

Indeed, let $x\in\operatorname*{Ker}\left(  \phi_{\operatorname{st}}\right)
$ be arbitrary. Write $x\in\operatorname*{QSym}$ in the form $x=\sum_{L}%
x_{L}F_{L}$, where the sum ranges over all compositions $L$, and where the
$x_{L}$ are elements of $\mathbb{Q}$ (all but finitely many of which are
zero). (This can be done, since $\left(  F_{L}\right)  _{L\text{ is a
composition}}$ is a basis of the $\mathbb{Q}$-vector space
$\operatorname*{QSym}$.) Now,
$x\in\operatorname*{Ker}\left(  \phi_{\operatorname{st}}\right)  $,
so that $\phi_{\operatorname{st}}\left( x\right)  =0$. Thus,%
\begin{align*}
0  &  =\phi_{\operatorname{st}}\left(  x\right)  =\sum_{L}x_{L}%
\phi_{\operatorname{st}}\left(  F_{L}\right)  \ \ \ \ \ \ \ \ \ \ \left(
\text{since }x=\sum_{L}x_{L}F_{L}\right) \\
&  =\sum_{\alpha}\sum_{L\in\alpha}x_{L}\underbrace{\phi_{\operatorname{st}}
\left(  F_{L}\right)  }_{\substack{=u_{\alpha}\\\text{(by
(\ref{pf.prop.K.ideal.dir1.1}))}}}\ \ \ \ \ \ \ \ \ \ \left(
\begin{array}
[c]{c}%
\text{where the first sum is over}\\
\text{all }\operatorname{st}\text{-equivalence classes }\alpha\text{ of
compositions}%
\end{array}
\right) \\
&  =\sum_{\alpha}\sum_{L\in\alpha}x_{L}u_{\alpha}=\sum_{\alpha}\left(
\sum_{L\in\alpha}x_{L}\right)  u_{\alpha}.
\end{align*}
Since the family $\left(  u_{\alpha}\right)  $ is linearly independent
(because it is a basis of $A$), we thus conclude that
\begin{equation}
\sum_{L\in\alpha}x_{L}=0 \label{pf.prop.K.ideal.dir1.2}%
\end{equation}
for each $\operatorname{st}$-equivalence class $\alpha$ of compositions.

Now, for each $\operatorname{st}$-equivalence class $\alpha$ of compositions,
we fix an element $L_{\alpha}$ of $\alpha$. Then, for each
$\operatorname{st}$-equivalence class $\alpha$ of compositions and each
composition $L\in\alpha $, we have
\begin{equation}
F_{L}-F_{L_{\alpha}}\in\mathcal{K}_{\operatorname{st}}
\label{pf.prop.K.ideal.dir1.3}%
\end{equation}
(since the compositions $L$ and $L_{\alpha}$ are $\operatorname{st}%
$-equivalent\footnote{since both compositions $L$ and $L_{\alpha}$ lie in the
same $\operatorname{st}$-equivalence class $\alpha$}).

Now,%
\begin{align*}
x  &  =\sum_{L}x_{L}F_{L}\\
&  =\sum_{\alpha}\sum_{L\in\alpha}x_{L}\underbrace{F_{L}}_{=\left(
F_{L}-F_{L_{\alpha}}\right)  +F_{L_{\alpha}}}\ \ \ \ \ \ \ \ \ \ \left(
\begin{array}
[c]{c}%
\text{where the first sum is over}\\
\text{all }\operatorname{st}\text{-equivalence classes }\alpha\text{ of
compositions}%
\end{array}
\right) \\
&  =\sum_{\alpha}\sum_{L\in\alpha}x_{L}\left(  \left(  F_{L}-F_{L_{\alpha}%
}\right)  +F_{L_{\alpha}}\right) \\
&  =\sum_{\alpha}\sum_{L\in\alpha}x_{L}\underbrace{\left(  F_{L}-F_{L_{\alpha
}}\right)  }_{\substack{\in\mathcal{K}_{\operatorname{st}}\\\text{(by
(\ref{pf.prop.K.ideal.dir1.3}))}}}+\sum_{\alpha}\underbrace{\sum_{L\in\alpha
}x_{L}}_{\substack{=0\\\text{(by (\ref{pf.prop.K.ideal.dir1.2}))}%
}}F_{L_{\alpha}}\\
&  \in\underbrace{\sum_{\alpha}\sum_{L\in\alpha}x_{L}\mathcal{K}%
_{\operatorname{st}}}_{\subseteq\mathcal{K}_{\operatorname{st}}%
}+\underbrace{\sum_{\alpha}0F_{L_{\alpha}}}_{=0}\subseteq\mathcal{K}%
_{\operatorname{st}}.
\end{align*}

Now, forget that we fixed $x$. We thus have proven that $x\in\mathcal{K}%
_{\operatorname{st}}$ for each $x\in\operatorname*{Ker}\left(  \phi
_{\operatorname{st}}\right)  $. In other words, $\operatorname*{Ker}\left(
\phi_{\operatorname{st}}\right)  \subseteq\mathcal{K}_{\operatorname{st}}$.

Conversely, it is easy to see that $\mathcal{K}_{\operatorname{st}}%
\subseteq\operatorname*{Ker}\left(  \phi_{\operatorname{st}}\right)
$\ \ \ \ \footnote{\textit{Proof.} Recall that $\mathcal{K}%
_{\operatorname{st}}$ is the $\mathbb{Q}$-vector subspace of
$\operatorname*{QSym}$ spanned by all elements of the form $F_{J}-F_{K}$,
where $J$ and $K$ are two $\operatorname{st}$-equivalent compositions. Hence,
in order to prove that $\mathcal{K}_{\operatorname{st}}\subseteq
\operatorname*{Ker}\left(  \phi_{\operatorname{st}}\right)  $, it suffices to
show that
$F_{J}-F_{K}\in\operatorname*{Ker}\left(  \phi_{\operatorname{st}}\right)  $,
whenever $J$ and $K$ are two $\operatorname{st}$-equivalent
compositions. So let $J$ and $K$ be two $\operatorname{st}$-equivalent
compositions. We must show that $F_{J}-F_{K}\in\operatorname*{Ker}\left(
\phi_{\operatorname{st}}\right)  $.
\par
The two compositions $J$ and $K$ are $\operatorname{st}$-equivalent. Hence,
they lie in one and the same $\operatorname{st}$-equivalence class. Let
$\alpha$ be this $\operatorname{st}$-equivalence class. Then, $J\in\alpha$
and therefore $\phi_{\operatorname{st}}\left(  F_{J}\right)  =u_{\alpha}$ (by
(\ref{pf.prop.K.ideal.dir1.1}), applied to $L=J$). Similarly, $\phi
_{\operatorname{st}}\left(  F_{K}\right)  =u_{\alpha}$. Now, $\phi
_{\operatorname{st}}\left(  F_{J}-F_{K}\right)  =\underbrace{\phi
_{\operatorname{st}}\left(  F_{J}\right)  }_{=u_{\alpha}}-\underbrace{\phi
_{\operatorname{st}}\left(  F_{K}\right)  }_{=u_{\alpha}}=u_{\alpha
}-u_{\alpha}=0$. In other words, $F_{J}-F_{K}\in\operatorname*{Ker}\left(
\phi_{\operatorname{st}}\right)  $. This completes our proof.}. Combining
this with $\operatorname*{Ker}\left(  \phi_{\operatorname{st}}\right)
\subseteq\mathcal{K}_{\operatorname{st}}$, we obtain $\mathcal{K}%
_{\operatorname{st}}=\operatorname*{Ker}\left(  \phi_{\operatorname{st}%
}\right)  $. This proves Lemma \ref{lem.K.fist}.
\end{proof}
\end{verlong}

Theorem~\ref{thm.4.3} easily yields the following fact:

\begin{proposition}
\label{prop.K.ideal}Let $\operatorname{st}$ be a descent statistic. Then,
$\operatorname{st}$ is shuffle-compatible if and only if $\mathcal{K}%
_{\operatorname{st}}$ is an ideal of $\operatorname*{QSym}$. Furthermore, in
this case, $\mathcal{A}_{\operatorname{st}}\cong\operatorname*{QSym}%
/\mathcal{K}_{\operatorname{st}}$ as $\mathbb{Q}$-algebras.
\end{proposition}

\begin{verlong}
\begin{proof}
[Proof of Proposition \ref{prop.K.ideal}.]$\Longrightarrow:$ Assume that
$\operatorname{st}$ is shuffle-compatible. Thus, Theorem~\ref{thm.4.3}
\textbf{(a)} shows that there exist a $\mathbb{Q}$-algebra $A$ with basis
$\left(  u_{\alpha}\right)  $ indexed by $\operatorname{st}$-equivalence
classes $\alpha$ of compositions, and a $\mathbb{Q}$-algebra homomorphism
$\phi_{\operatorname{st}}:\operatorname*{QSym}\rightarrow A$ with the
property that whenever $\alpha$ is an $\operatorname{st}$-equivalence class
of compositions, we have%
\begin{equation}
\phi_{\operatorname{st}}\left(  F_{L}\right)  =u_{\alpha}%
\ \ \ \ \ \ \ \ \ \ \text{for each }L\in\alpha.
\label{pf.prop.K.ideal.dir1.13}%
\end{equation}
Consider this $A$ and this $\phi_{\operatorname{st}}$.

Lemma \ref{lem.K.fist} yields that $\mathcal{K}_{\operatorname{st}}
=\operatorname*{Ker}\left(  \phi_{\operatorname{st}}\right)  $. But the map
$\phi_{\operatorname{st}}$ is a $\mathbb{Q}$-algebra homomorphism. Thus, its
kernel $\operatorname*{Ker}\left(  \phi_{\operatorname{st}}\right)  $ is an
ideal of $\operatorname*{QSym}$. In other words, $\mathcal{K}%
_{\operatorname{st}}$ is an ideal of $\operatorname*{QSym}$ (since
$\mathcal{K}_{\operatorname{st}}=\operatorname*{Ker}\left(  \phi
_{\operatorname{st}}\right)  $).

It remains to show that $\mathcal{A}_{\operatorname{st}}\cong%
\operatorname*{QSym}/\mathcal{K}_{\operatorname{st}}$ as $\mathbb{Q}$-algebras.
This is easy: Each
element of the basis $\left(  u_{\alpha}\right)  $ of the $\mathbb{Q}$-vector
space $A$ is contained in the image of $\phi_{\operatorname{st}}$ (because of
(\ref{pf.prop.K.ideal.dir1.13})). Therefore, the homomorphism $\phi
_{\operatorname{st}}$ is surjective. Thus, $\phi_{\operatorname{st}}\left(
\operatorname*{QSym}\right)  =A$. Hence, $A=\phi_{\operatorname{st}}\left(
\operatorname*{QSym}\right)  \cong\operatorname*{QSym}/\operatorname*{Ker}%
\left(  \phi_{\operatorname{st}}\right)  $ (by the homomorphism theorem). But
Theorem \ref{thm.4.3} \textbf{(b)} shows that
$\mathcal{A}_{\operatorname{st}}\cong A$.
Thus, $\mathcal{A}_{\operatorname{st}}\cong A\cong
\operatorname*{QSym}/\underbrace{\operatorname*{Ker}\left(  \phi
_{\operatorname{st}}\right)  }_{=\mathcal{K}_{\operatorname{st}}%
}=\operatorname*{QSym}/\mathcal{K}_{\operatorname{st}}$. This finishes the
proof of the $\Longrightarrow$ direction of Proposition \ref{prop.K.ideal}.

$\Longleftarrow:$ Assume that $\mathcal{K}_{\operatorname{st}}$ is an ideal
of $\operatorname*{QSym}$. We must prove that $\operatorname{st}$ is shuffle-compatible.

We shall not use this direction of Proposition \ref{prop.K.ideal}, so let us
merely sketch the proof. Let $A$ be the $\mathbb{Q}$-algebra
$\operatorname*{QSym}/\mathcal{K}_{\operatorname{st}}$. Let $\phi
_{\operatorname{st}}$ be the canonical projection $\operatorname*{QSym}%
\rightarrow A$; this is clearly a $\mathbb{Q}$-algebra homomorphism.

For each $\operatorname{st}$-equivalence class $\alpha$ of compositions, we
define an element $u_{\alpha}$ of $A$ by requiring that%
\[
u_{\alpha}=\phi_{\operatorname{st}}\left(  F_{L}\right)
\ \ \ \ \ \ \ \ \ \ \text{whenever }L\in\alpha.
\]
This is easily seen to be well-defined, because the image $\phi
_{\operatorname{st}}\left(  F_{L}\right)  $ depends only on $\alpha$ but not
on $L$ (indeed, if $J$ and $K$ are two elements of $\alpha$, then $J$ and $K$
are $\operatorname{st}$-equivalent, whence $F_{J}-F_{K}\in\mathcal{K}%
_{\operatorname{st}}$, whence $F_{J}\equiv F_{K}\operatorname{mod}
\mathcal{K}_{\operatorname{st}}$ and therefore $\phi_{\operatorname{st}}
\left(  F_{J}\right)  =\phi_{\operatorname{st}}\left(  F_{K}\right)  $).

It is not hard to see that the family $\left(  u_{\alpha}\right)  $ (where
$\alpha$ ranges over all $\operatorname{st}$-equivalence classes of
compositions) is a basis of the $\mathbb{Q}$-algebra $A$. Hence, Theorem
\ref{thm.4.3} \textbf{(a)} yields that $\operatorname{st}$ is
shuffle-compatible. This proves the $\Longleftarrow$ direction of Proposition
\ref{prop.K.ideal}.
\end{proof}
\end{verlong}

\begin{corollary}
\label{cor.Epk.ideal}The kernel $\mathcal{K}_{\operatorname{Epk}}$ of the
descent statistic $\operatorname{Epk}$ is an ideal of $\operatorname*{QSym}$.
\end{corollary}

\begin{verlong}
\begin{proof}
[Proof of Corollary \ref{cor.Epk.ideal}.]This follows from Proposition
\ref{prop.K.ideal} (applied to $\operatorname{st}=\operatorname{Epk}$),
because of Theorem \ref{thm.Epk.sh-co-a}.
\end{proof}
\end{verlong}

We can study the kernel of any descent statistic; in particular, the case of
shuffle-compatible descent statistics appears interesting. Since
$\operatorname*{QSym}$ is isomorphic to a polynomial ring (as an algebra), it
has many ideals, which are rather hopeless to classify or tame; but the ones
obtained as kernels of shuffle-compatible descent statistics might be worth discussing.

\subsection{An F-generating set of $\mathcal{K}_{\operatorname{Epk}}$}

Let us now focus on $\mathcal{K}_{\operatorname{Epk}}$, the kernel of
$\operatorname{Epk}$.

\begin{proposition}
\label{prop.K.Epk.F}If $J=\left(  j_{1},j_{2},\ldots,j_{m}\right)  $ and $K$
are two compositions, then we shall write $J\rightarrow K$ if there exists an
$\ell\in\left\{  2,3,\ldots,m\right\}  $ such that $j_{\ell}>2$ and $K=\left(
j_{1},j_{2},\ldots,j_{\ell-1},1,j_{\ell}-1,j_{\ell+1},j_{\ell+2},\ldots
,j_{m}\right)  $. (In other words, we write $J\rightarrow K$ if $K$ can be
obtained from $J$ by \textquotedblleft splitting\textquotedblright\ some entry
$j_{\ell}>2$ into two consecutive entries\footnotemark\ $1$ and $j_{\ell}-1$,
provided that this entry was not the first entry -- i.e., we had $\ell>1$ --
and that this entry was greater than $2$.)

The ideal $\mathcal{K}_{\operatorname{Epk}}$ of $\operatorname*{QSym}$ is
spanned (as a $\mathbb{Q}$-vector space) by all differences of the form
$F_{J}-F_{K}$, where $J$ and $K$ are two compositions satisfying $J\rightarrow
K$.
\end{proposition}

\footnotetext{The word \textquotedblleft consecutive\textquotedblright\ here
means \textquotedblleft in consecutive positions of $J$\textquotedblright, not
\textquotedblleft consecutive integers\textquotedblright. So two consecutive
entries of $J$ are two entries of the form $j_{p}$ and $j_{p+1}$ for some
$p\in\left\{  1,2,\ldots,m-1\right\}  $.}

\begin{example}
We have $\left(  2,1,4,4\right)  \rightarrow\left(  2,1,1,3,4\right)  $, since
the composition $\left(  2,1,1,3,4\right)  $ is obtained from $\left(
2,1,4,4\right)  $ by splitting the third entry (which is $4>2$) into two
consecutive entries $1$ and $3$.

Similarly, $\left(  2,1,4,4\right)  \rightarrow\left(  2,1,4,1,3\right)  $.

But we do not have $\left(  3,1\right)  \rightarrow\left(  1,2,1\right)  $,
because splitting the first entry of the composition is not allowed in the
definition of the relation $\rightarrow$. Also, we do not have $\left(
1,2,1\right)  \rightarrow\left(  1,1,1,1\right)  $, because the entry we are
splitting must be $>2$.

Two compositions $J$ and $K$ satisfying $J\rightarrow K$ must necessarily
satisfy $\left\vert J\right\vert =\left\vert K\right\vert $.

Here are all relations $\rightarrow$ between compositions of size $4$:%
\[
\left(  1,3\right)  \rightarrow\left(  1,1,2\right)  .
\]

Here are all relations $\rightarrow$ between compositions of size $5$:%
\begin{align*}
\left(  1,4\right)   &  \rightarrow\left(  1,1,3\right)  ,\\
\left(  1,3,1\right)   &  \rightarrow\left(  1,1,2,1\right)  ,\\
\left(  1,1,3\right)   &  \rightarrow\left(  1,1,1,2\right)  ,\\
\left(  2,3\right)   &  \rightarrow\left(  2,1,2\right)  .
\end{align*}
There are no relations $\rightarrow$ between compositions of size $\leq3$.
\end{example}

\begin{verlong}
\begin{proof}
[Proof of Proposition \ref{prop.K.Epk.F}.]We begin by proving some simple claims.

\begin{statement}
\textit{Claim 1:} Let $n\in\mathbb{N}$. Let $J$ and $K$ be two compositions of
size $n$. Then, $J\rightarrow K$ if and only if there exists some $k\in\left[
n-1\right]  $ such that
\begin{align*}
\operatorname{Des}K  &  =\operatorname{Des}J\cup\left\{  k\right\}
,\ \ \ \ \ \ \ \ \ \ k\notin\operatorname{Des}J,\\
k-1  &  \in\operatorname{Des}J\ \ \ \ \ \ \ \ \ \ \text{and}%
\ \ \ \ \ \ \ \ \ \ k+1\notin\operatorname{Des}J\cup\left\{  n\right\}  .
\end{align*}

\end{statement}

[\textit{Proof of Claim 1:} This is straightforward to check:
\textquotedblleft Splitting\textquotedblright\ an entry of a composition $C$
into two consecutive entries (summing up to the original entry) is always
tantamount to adding a new element to $\operatorname{Des}C$. The rest is
translating conditions.]

If $n$ is a positive integer, and $L$ is any composition of $n$, then
\begin{equation}
\operatorname{Epk}L=\left(  \operatorname{Des}L\cup\left\{  n\right\}
\right)  \setminus\left(  \operatorname{Des}L+1\right)  .
\label{pf.prop.K.Epk.F.1}%
\end{equation}
(This is a consequence of Proposition \ref{prop.Epk.through-Des}, applied to
any $n$-permutation $\pi$ satisfying $L=\operatorname*{Comp}\pi$.)

\begin{noncompile}
Also, $\operatorname{Epk}L=\operatorname*{Lpk}L\cup\operatorname*{Rpk}L$ for
any integer $n\geq2$ and any composition $L$ of $n$.
\end{noncompile}

\begin{statement}
\textit{Claim 2:} Let $J$ and $K$ be two compositions satisfying $J\rightarrow
K$. Then, $\operatorname{Epk}J=\operatorname{Epk}K$.
\end{statement}

[\textit{Proof of Claim 2:} Easy consequence of Claim 1 and
(\ref{pf.prop.K.Epk.F.1}).]

For any two integers $a$ and $b$, we set $\left[  a,b\right]  =\left\{
a,a+1,\ldots,b\right\}  $. (This is an empty set if $a>b$.)

It is easy to see that every composition $J$ of size $n>0$ satisfies%
\begin{equation}
\left[  \max\left(  \operatorname{Epk}J\right)  ,n-1\right]  \subseteq
\operatorname{Des}J . \label{pf.prop.K.Epk.F.2}
\end{equation}

[\textit{Proof of (\ref{pf.prop.K.Epk.F.2}):} Let $J$ be a
composition of size $n>0$. We shall show that $\left[  \max\left(
\operatorname{Epk}J\right)  ,n-1\right]  \subseteq\operatorname{Des}J$.
\par
Indeed, assume the contrary. Thus, $\left[  \max\left(  \operatorname{Epk}
J\right)  ,n-1\right]  \not \subseteq \operatorname{Des}J$. Hence, there
exists some $q\in\left[  \max\left(  \operatorname{Epk}J\right)  ,n-1\right]
$ satisfying $q\notin\operatorname{Des}J$. Let $r$ be the \textbf{largest}
such $q$.
\par
Thus, $r\in\left[  \max\left(  \operatorname{Epk}J\right)  ,n-1\right]  $ but
$r\notin\operatorname{Des}J$. From $r\in\left[  \max\left(
\operatorname{Epk}J\right)  ,n-1\right]  \subseteq\left[  n-1\right]  $, we
obtain $r+1\in\left[  n\right]  $. Also, from $r\notin\operatorname{Des}J$,
we obtain $r+1\notin\operatorname{Des}J+1$.
\par
From $r\in\left[  \max\left(  \operatorname{Epk}J\right)  ,n-1\right]  $, we
obtain $r\geq\max\left(  \operatorname{Epk}J\right)  $, so that
$r+1>r\geq\max\left(  \operatorname{Epk}J\right)  $ and therefore
$r+1\notin\operatorname{Epk}J$ (since a number that is higher than
$\max\left(  \operatorname{Epk}J\right)  $ cannot belong to
$\operatorname{Epk}J$).
\par
From (\ref{pf.prop.K.Epk.F.1}), we obtain $\operatorname{Epk}J=\left(
\operatorname{Des}J\cup\left\{  n\right\}  \right)  \setminus\left(
\operatorname{Des}J+1\right)  $.
\par
If we had $r+1\in\operatorname{Des}J\cup\left\{  n\right\}  $, then we would
have $r+1\in\left(  \operatorname{Des}J\cup\left\{  n\right\}  \right)
\setminus\left(  \operatorname{Des}J+1\right)  $ (since $r+1\notin%
\operatorname{Des}J+1$). This would contradict $r+1\notin\operatorname{Epk}
J=\left(  \operatorname{Des}J\cup\left\{  n\right\}  \right)  \setminus
\left(  \operatorname{Des}J+1\right)  $. Thus, we cannot have $r+1\in
\operatorname{Des}J\cup\left\{  n\right\}  $. Therefore, $r+1\notin%
\operatorname{Des}J\cup\left\{  n\right\}  $.
\par
Hence, $r+1\neq n$ (since $r+1\notin\operatorname{Des}J\cup\left\{
n\right\}  $ but $n\in\left\{  n\right\}  \subseteq\operatorname{Des}
J\cup\left\{  n\right\}  $). Combined with $r+1\in\left[  n\right]  $, this
yields $r+1\in\left[  n\right]  \setminus\left\{  n\right\}  =\left[
n-1\right]  $. Combined with $r+1>\max\left(  \operatorname{Epk}J\right)  $,
this yields $r+1\in\left[  \max\left(  \operatorname{Epk}J\right)
,n-1\right]  $. Also, $r+1\notin\operatorname{Des}J$ (since $r+1\notin%
\operatorname{Des}J\cup\left\{  n\right\}  $).
\par
Thus, $r+1$ is a $q\in\left[  \max\left(  \operatorname{Epk}J\right)
,n-1\right]  $ satisfying $q\notin\operatorname{Des}J$. This contradicts the
fact that $r$ is the \textbf{largest} such $q$ (since $r+1$ is clearly larger
than $r$). This contradiction proves that our assumption was wrong; thus,
(\ref{pf.prop.K.Epk.F.2}) is proven.]

For each $n\in\mathbb{N}$ and each subset $S$ of $\left[  n-1\right]  $, we
define a subset $S_{n}^{\circ}$ of $\left[  n-1\right]  $ as follows:%
\[
S_{n}^{\circ}=\left\{  s\in S\ \mid\ s-1\notin S\text{ or }\left[
s,n-1\right]  \subseteq S\right\}  .
\]

Also, for each $n\in\mathbb{N}$ and each nonempty subset $T$ of $\left[
n\right]  $, we define a subset $\rho_{n}\left(  T\right)  $ of $\left[
n-1\right]  $ as follows:%
\[
\rho_{n}\left(  T\right)  =%
\begin{cases}
T\setminus\left\{  n\right\}  , & \text{if }n\in T;\\
T\cup\left[  \max T,n-1\right]  , & \text{if }n\notin T
\end{cases}
.
\]

\begin{statement}
\textit{Claim 3:} Let $n$ be a positive integer.
Let $J$ be a composition of size $n$.
Then, $\left(  \operatorname{Des}J\right)  _{n}^{\circ}
=\rho_{n}\left( \operatorname{Epk}J\right)  $.
\end{statement}

[\textit{Proof of Claim 3:}
Proposition \ref{prop.Epk-lac} yields
$\operatorname{Epk}J\in\mathbf{L}_{n}$ (since
$\operatorname{Epk}J = \operatorname{Epk}\pi$
for an appropriate $n$-permutation $\pi$).
Since $n$ is positive, this entails that
$\operatorname{Epk}J$ is nonempty (by the definition of
$\mathbf{L}_{n}$). Thus, $\rho_{n}\left( \operatorname{Epk}J\right)$
is well-defined.

Let $g\in\left(  \operatorname{Des}J\right)
_{n}^{\circ}$. We shall show that $g\in\rho_{n}\left(  \operatorname{Epk}
J\right)  $.

We have $g\in\left(  \operatorname{Des}J\right)  _{n}^{\circ}\subseteq
\operatorname{Des}J$ (since $S_{n}^{\circ}\subseteq S$ for each subset $S$ of
$\left[  n-1\right]  $) and therefore $\operatorname{Des}J\neq\varnothing$.
Hence, $J$ is not the empty composition. In other words, $n>0$.

From (\ref{pf.prop.K.Epk.F.1}), we obtain $\operatorname{Epk}J=\left(
\operatorname{Des}J\cup\left\{  n\right\}  \right)  \setminus\left(
\operatorname{Des}J+1\right)  $. Thus, the set $\operatorname{Epk}J$ is
disjoint from $\operatorname{Des}J+1$. Furthermore, the set
$\operatorname{Epk}J$ is nonempty\footnote{Indeed, it contains at least the
smallest element of the set $\operatorname{Des}J\cup\left\{  n\right\}  $
(since $\operatorname{Epk}J=\left(  \operatorname{Des}J\cup\left\{
n\right\}  \right)  \setminus\left(  \operatorname{Des}J+1\right)  $).}.

We have $g\in\left(  \operatorname{Des}J\right)  _{n}^{\circ}$. Thus, $g$ is
an element of $\operatorname{Des}J$ satisfying $g-1\notin\operatorname{Des}
J$ or $\left[  g,n-1\right]  \subseteq\operatorname{Des}J$ (by the definition
of $\left(  \operatorname{Des}J\right)  _{n}^{\circ}$). We are thus in one of
the following two cases:

\textit{Case 1:} We have $g-1\notin\operatorname{Des}J$.

\textit{Case 2:} We have $\left[  g,n-1\right]  \subseteq\operatorname{Des}J$.

Let us first consider Case 1. In this case, we have $g-1\notin%
\operatorname{Des}J$. In other words, $g\notin\operatorname{Des}J+1$.
Combined with $g\in\operatorname{Des}J\subseteq\operatorname{Des}
J\cup\left\{  n\right\}  $, this yields $g\in\left(  \operatorname{Des}
J\cup\left\{  n\right\}  \right)  \setminus\left(  \operatorname{Des}
J+1\right)  =\operatorname{Epk}J$. Moreover, $g\neq n$ (since $g\in
\operatorname{Des}J\subseteq\left[  n-1\right]  $) and thus $g\in\left(
\operatorname{Epk}J\right)  \setminus\left\{  n\right\}  $ (since
$g\in\operatorname{Epk}J$). But each nonempty subset $T$ of $\left[
n\right]  $ satisfies $T\setminus\left\{  n\right\}  \subseteq\rho_{n}\left(
T\right)  $ (by the definition of $\rho_{n}\left(  T\right)  $). Applying this
to $T=\operatorname{Epk}J$, we obtain $\left(  \operatorname{Epk}J\right)
\setminus\left\{  n\right\}  \subseteq\rho_{n}\left(  \operatorname{Epk}
J\right)  $. Hence, $g\in\left(  \operatorname{Epk}J\right)  \setminus
\left\{  n\right\}  \subseteq\rho_{n}\left(  \operatorname{Epk}J\right)  $.
Thus, $g\in\rho_{n}\left(  \operatorname{Epk}J\right)  $ is proven in Case 1.

Let us now consider Case 2. In this case, we have $\left[  g,n-1\right]
\subseteq\operatorname{Des}J$. Hence, each of the elements $g,g+1,\ldots,n-1$
belongs to $\operatorname{Des}J$. In other words, each of the elements
$g+1,g+2,\ldots,n$ belongs to $\operatorname{Des}J+1$. Hence, none of the
elements $g+1,g+2,\ldots,n$ belongs to $\operatorname{Epk}J$ (since the set
$\operatorname{Epk}J$ is disjoint from $\operatorname{Des}J+1$). Thus,
$\max\left(  \operatorname{Epk}J\right)  \leq g$. Therefore, $g\in\left[
\max\left(  \operatorname{Epk}J\right)  ,n-1\right]  $ (since $g\in
\operatorname{Des}J\subseteq\left[  n-1\right]  $).

Also, $n\notin\operatorname{Epk}J$\ \ \ \ \footnote{\textit{Proof.} Assume
the contrary. Thus, $n\in\operatorname{Epk}J$. But none of the elements
$g+1,g+2,\ldots,n$ belongs to $\operatorname{Epk}J$. Hence, $n$ is not among
the elements $g+1,g+2,\ldots,n$. Therefore, $g\geq n$, so that $g=n$. This
contradicts $g\in\operatorname{Des}J\subseteq\left[  n-1\right]  $. This
contradiction shows that our assumption was wrong, qed.}. Hence, the
definition of $\rho_{n}\left(  \operatorname{Epk}J\right)  $ yields $\rho
_{n}\left(  \operatorname{Epk}J\right)  =\operatorname{Epk}J\cup\left[
\max\left(  \operatorname{Epk}J\right)  ,n-1\right]  $. Now,%
\[
g\in\left[  \max\left(  \operatorname{Epk}J\right)  ,n-1\right]
\subseteq\operatorname{Epk}J\cup\left[  \max\left(  \operatorname{Epk}
J\right)  ,n-1\right]  =\rho_{n}\left(  \operatorname{Epk}J\right)  .
\]
Hence, $g\in\rho_{n}\left(  \operatorname{Epk}J\right)  $ is proven in Case 2.

Thus, $g\in\rho_{n}\left(  \operatorname{Epk}J\right)  $ is proven in both
Cases 1 and 2. This shows that $g\in\rho_{n}\left(  \operatorname{Epk}
J\right)  $ always holds.

Forget that we fixed $g$. We thus have proven that $g\in\rho_{n}\left(
\operatorname{Epk}J\right)  $ for each $g\in\left(  \operatorname{Des}
J\right)  _{n}^{\circ}$. In other words, $\left(  \operatorname{Des}J\right)
_{n}^{\circ}\subseteq\rho_{n}\left(  \operatorname{Epk}J\right)  $.

Now, let $h\in\rho_{n}\left(  \operatorname{Epk}J\right)  $ be arbitrary. We
shall prove that $h\in\left(  \operatorname{Des}J\right)  _{n}^{\circ}$.

We are in one of the following two cases:

\textit{Case 1:} We have $n\in\operatorname{Epk}J$.

\textit{Case 2:} We have $n\notin\operatorname{Epk}J$.

Let us first consider Case 1. In this case, we have $n\in\operatorname{Epk}
J$, and thus $\rho_{n}\left(  \operatorname{Epk}J\right)
=\operatorname{Epk}J\setminus\left\{  n\right\}  $ (by the definition of
$\rho_{n}\left(  \operatorname{Epk}J\right)  $). Hence,
\[
h\in\rho_{n}\left(  \operatorname{Epk}J\right)  =\operatorname{Epk}
J\setminus\left\{  n\right\}  \subseteq\operatorname{Epk}J=\left(
\operatorname{Des}J\cup\left\{  n\right\}  \right)  \setminus\left(
\operatorname{Des}J+1\right)  .
\]
In other words, $h\in\operatorname{Des}J\cup\left\{  n\right\}  $ and
$h\notin\operatorname{Des}J+1$. Since $h\in\operatorname{Des}J\cup\left\{
n\right\}  $ and $h\neq n$ (because $h\in\rho_{n}\left(  \operatorname{Epk}
J\right)  \subseteq\left[  n-1\right]  $), we obtain $h\in\left(
\operatorname{Des}J\cup\left\{  n\right\}  \right)  \setminus\left\{
n\right\}  \subseteq\operatorname{Des}J$. From $h\notin\operatorname{Des}
J+1$, we obtain $h-1\notin\operatorname{Des}J$. Thus, $h$ is an element of
$\operatorname{Des}J$ satisfying $h-1\notin\operatorname{Des}J$ or $\left[
h,n-1\right]  \subseteq\operatorname{Des}J$ (in fact, $h-1\notin%
\operatorname{Des}J$ holds). Thus, $h\in\left(  \operatorname{Des}J\right)
_{n}^{\circ}$ (by the definition of $\left(  \operatorname{Des}J\right)
_{n}^{\circ}$). Thus, $h\in\left(  \operatorname{Des}J\right)  _{n}^{\circ}$
is proven in Case 1.

Let us now consider Case 2. In this case, we have $n\notin\operatorname{Epk}
J$. Hence, the definition of $\rho_{n}\left(  \operatorname{Epk}J\right)  $
yields $\rho_{n}\left(  \operatorname{Epk}J\right)  =\left(
\operatorname{Epk}J\right)  \cup\left[  \max\left(  \operatorname{Epk}
J\right)  ,n-1\right]  $. Thus, $h\in\rho_{n}\left(  \operatorname{Epk}
J\right)  =\left(  \operatorname{Epk}J\right)  \cup\left[  \max\left(
\operatorname{Epk}J\right)  ,n-1\right]  $.

If $h\in\operatorname{Epk}J$, then we can prove $h\in\left(
\operatorname{Des}J\right)  _{n}^{\circ}$ just as in Case 1. Hence, let us
WLOG assume that we don't have $h\in\operatorname{Epk}J$. Thus,
$h\notin\operatorname{Epk}J$. Combined with $h\in\left(  \operatorname{Epk}
J\right)  \cup\left[  \max\left(  \operatorname{Epk}J\right)  ,n-1\right]  $,
this yields%
\begin{align*}
h  &  \in\left(  \left(  \operatorname{Epk}J\right)  \cup\left[  \max\left(
\operatorname{Epk}J\right)  ,n-1\right]  \right)  \setminus\left(
\operatorname{Epk}J\right)  =\left[  \max\left(  \operatorname{Epk}J\right)
,n-1\right]  \setminus\left(  \operatorname{Epk}J\right) \\
&  \subseteq\left[  \max\left(  \operatorname{Epk}J\right)  ,n-1\right]
\subseteq\operatorname{Des}J\ \ \ \ \ \ \ \ \ \ \left(  \text{by
(\ref{pf.prop.K.Epk.F.2})}\right)  .
\end{align*}
Moreover, from $h\in\left[  \max\left(  \operatorname{Epk}J\right)
,n-1\right]  $, we obtain $h\geq\max\left(  \operatorname{Epk}J\right)  $, so
that%
\[
\left[  h,n-1\right]  \subseteq\left[  \max\left(  \operatorname{Epk}%
J\right)  ,n-1\right]  \subseteq\operatorname{Des}%
J\ \ \ \ \ \ \ \ \ \ \left(  \text{by (\ref{pf.prop.K.Epk.F.2})}\right)  .
\]
Hence, $h$ is an element of $\operatorname{Des}J$ satisfying $h-1\notin%
\operatorname{Des}J$ or $\left[  h,n-1\right]  \subseteq\operatorname{Des}J$
(namely, $\left[  h,n-1\right]  \subseteq\operatorname{Des}J$). In other
words, $h\in\left(  \operatorname{Des}J\right)  _{n}^{\circ}$ (by the
definition of $\left(  \operatorname{Des}J\right)  _{n}^{\circ}$). Thus,
$h\in\left(  \operatorname{Des}J\right)  _{n}^{\circ}$ is proven in Case 2.

We have now proven $h\in\left(  \operatorname{Des}J\right)  _{n}^{\circ}$ in
both Cases 1 and 2. Hence, $h\in\left(  \operatorname{Des}J\right)
_{n}^{\circ}$ always holds.

Forget that we fixed $h$. We thus have shown that $h\in\left(
\operatorname{Des}J\right)  _{n}^{\circ}$ for each $h\in\rho_{n}\left(
\operatorname{Epk}J\right)  $. In other words, $\rho_{n}\left(
\operatorname{Epk}J\right)  \subseteq\left(  \operatorname{Des}J\right)
_{n}^{\circ}$. Combining this with $\left(  \operatorname{Des}J\right)
_{n}^{\circ}\subseteq\rho_{n}\left(  \operatorname{Epk}J\right)  $, we obtain
$\left(  \operatorname{Des}J\right)  _{n}^{\circ}=\rho_{n}\left(
\operatorname{Epk}J\right)  $. This proves Claim 3.]

\begin{statement}
\textit{Claim 4:} Let $n\in\mathbb{N}$. Let $J$ and $K$ be two compositions of
size $n$ satisfying $\operatorname{Epk}J=\operatorname{Epk}K$. Then,
$\left(  \operatorname{Des}J\right)  _{n}^{\circ}=\left(  \operatorname{Des}%
K\right)  _{n}^{\circ}$.
\end{statement}

[\textit{Proof of Claim 4:} If $n=0$, then $J=K$ (since there is only
one composition of $0$), and thus Claim 4 is obvious in this case.
Hence, for the rest of this proof, we WLOG assume that $n\neq0$.
Hence, $n$ is a positive integer.

Thus, Claim 3 yields $\left(  \operatorname{Des}%
J\right)  _{n}^{\circ}=\rho_{n}\left(  \operatorname{Epk}J\right)  $ and
similarly $\left(  \operatorname{Des}K\right)  _{n}^{\circ}=\rho_{n}\left(
\operatorname{Epk}K\right)  $. Hence,%
\[
\left(  \operatorname{Des}J\right)  _{n}^{\circ}=\rho_{n}\left(
\underbrace{\operatorname{Epk}J}_{=\operatorname{Epk}K}\right)  =\rho
_{n}\left(  \operatorname{Epk}K\right)  =\left(  \operatorname{Des}K\right)
_{n}^{\circ}.
\]
This proves Claim 4.]

We let $\overset{\ast}{\rightarrow}$ be the transitive-and-reflexive closure
of the relation $\rightarrow$. Thus, two compositions $J$ and $K$ satisfy
$J\overset{\ast}{\rightarrow}K$ if and only if there exists a sequence
$\left(  L_{0},L_{1},\ldots,L_{\ell}\right)  $ of compositions satisfying
$L_{0}=J$ and $L_{\ell}=K$ and $L_{0}\rightarrow L_{1}\rightarrow
\cdots\rightarrow L_{\ell}$.

\begin{statement}
\textit{Claim 5:} Let $n\in\mathbb{N}$. Let $K$ be a composition of size $n$.
Then, $\operatorname*{Comp}\left(  \left(  \operatorname{Des}K\right)
_{n}^{\circ}\right)  \overset{\ast}{\rightarrow}K$.
\end{statement}

[\textit{Proof of Claim 5:} We shall prove Claim 5 by strong induction on
$\left\vert \left(  \operatorname{Des}K\right)  \setminus\left(
\operatorname{Des}K\right)  _{n}^{\circ}\right\vert $. Thus, we fix an
$n\in\mathbb{N}$ and a composition $K$ of size $n$, and we assume (as the
induction hypothesis) that each composition $J$ of size $n$ satisfying
$\left\vert \left(  \operatorname{Des}J\right)  \setminus\left(
\operatorname{Des}J\right)  _{n}^{\circ}\right\vert <\left\vert \left(
\operatorname{Des}K\right)  \setminus\left(  \operatorname{Des}K\right)
_{n}^{\circ}\right\vert $ satisfies $\operatorname*{Comp}\left(  \left(
\operatorname{Des}J\right)  _{n}^{\circ}\right)  \overset{\ast}{\rightarrow
}J$. Our goal is to prove that $\operatorname*{Comp}\left(  \left(
\operatorname{Des}K\right)  _{n}^{\circ}\right)  \overset{\ast}{\rightarrow
}K$.

Let $A=\operatorname{Des}K$. Thus, $K=\operatorname*{Comp}A$ and
$A\subseteq\left[  n-1\right]  $.

Applying (\ref{pf.prop.K.Epk.F.1}) to $L=K$, we obtain $\operatorname{Epk}%
K=\left(  \operatorname{Des}K\cup\left\{  n\right\}  \right)  \setminus
\left(  \operatorname{Des}K+1\right)  =\left(  A\cup\left\{  n\right\}
\right)  \setminus\left(  A+1\right)  $ (since $\operatorname{Des}K=A$).

Also, $A_{n}^{\circ}\subseteq A$ (since $S_{n}^{\circ}\subseteq S$ for any
subset $S$ of $\left[  n-1\right]  $). If $A_{n}^{\circ}=A$, then we are done
(because if $A_{n}^{\circ}=A$, then $\operatorname*{Comp}\left(  \left(
\underbrace{\operatorname{Des}K}_{=A}\right)  _{n}^{\circ}\right)
=\operatorname*{Comp}\left(  \underbrace{A_{n}^{\circ}}_{=A}\right)
=\operatorname*{Comp}A=K$, and therefore the reflexivity of $\overset{\ast
}{\rightarrow}$ shows that $\operatorname*{Comp}\left(  \left(
\operatorname{Des}K\right)  _{n}^{\circ}\right)  \overset{\ast}{\rightarrow
}K$). Hence, we WLOG assume that $A_{n}^{\circ}\neq A$. Thus, $A_{n}^{\circ}$
is a proper subset of $A$ (since $A_{n}^{\circ}\subseteq A$). Therefore, there
exists a $q\in A$ satisfying $q\notin A_{n}^{\circ}$. Let $k$ be the
\textbf{largest} such $q$. Thus, $k\in A$ and $k\notin A_{n}^{\circ}$. Hence,
$k\in A\setminus A_{n}^{\circ}$. Also, $k\in A\subseteq\left[  n-1\right]  $.

Let $J=\operatorname*{Comp}\left(  A\setminus\left\{  k\right\}  \right)  $.
Thus, $\operatorname{Des}J=A\setminus\left\{  k\right\}  $, so that
$A=\operatorname{Des}J\cup\left\{  k\right\}  $ (since $k\in A$). Hence,
$\operatorname{Des}K=A=\operatorname{Des}J\cup\left\{  k\right\}  $. Also,
$k\notin A\setminus\left\{  k\right\}  =\operatorname{Des}J$.

Furthermore, $k-1\in\operatorname{Des}J$\ \ \ \ \footnote{\textit{Proof.}
Assume the contrary. Thus, $k-1\notin\operatorname{Des}J$. Hence,
$k-1\notin\operatorname{Des}J\cup\left\{  k\right\}  $ as well (since
$k-1\neq k$). In other words, $k-1\notin A$ (since $\operatorname{Des}%
J\cup\left\{  k\right\}  =A$). Therefore, the element $k$ of $A$ satisfies
$k-1\notin A$ or $\left[  k,n-1\right]  \subseteq A$. Thus, the definition of
$A_{n}^{\circ}$ yields $k\in A_{n}^{\circ}$. This contradicts $k\notin
A_{n}^{\circ}$. This contradiction shows that our assumption was false; qed.}
and $k+1\notin\operatorname{Des}J\cup\left\{  n\right\}  $%
\ \ \ \ \footnote{\textit{Proof.} Assume the contrary. Thus, $k+1\in
\operatorname{Des}J\cup\left\{  n\right\}  $. In other words, we have
$k+1\in\operatorname{Des}J$ or $k+1=n$. In other words, we are in one of the
following two cases:
\par
\textit{Case 1:} We have $k+1\in\operatorname{Des}J$.
\par
\textit{Case 2:} We have $k+1=n$.
\par
Let us first consider Case 1. In this case, we have $k+1\in\operatorname{Des}%
J$. Hence, $k+1\in\operatorname{Des}J\subseteq\operatorname{Des}%
J\cup\left\{  k\right\}  =A$. If we had $k+1\notin A_{n}^{\circ}$, then $k+1$
would be a $q\in A$ satisfying $q\notin A_{n}^{\circ}$; this would contradict
the fact that $k$ is the \textbf{largest} such $q$ (since $k+1$ is larger than
$k$). Hence, we cannot have $k+1\notin A_{n}^{\circ}$. Thus, we must have
$k+1\in A_{n}^{\circ}$. In other words, $k+1$ is an element of $A$ satisfying
$\left(  k+1\right)  -1\notin A$ or $\left[  k+1,n-1\right]  \subseteq A$ (by
the definition of $A_{n}^{\circ}$). Since $\left(  k+1\right)  -1\notin A$ is
impossible (because $\left(  k+1\right)  -1=k\in A$), we thus have $\left[
k+1,n-1\right]  \subseteq A$. Now, $\left[  k,n-1\right]
=\underbrace{\left\{  k\right\}  }_{\substack{\subseteq A\\\text{(since }k\in
A\text{)}}}\cup\underbrace{\left[  k+1,n-1\right]  }_{\subseteq A}\subseteq
A\cup A=A$. Thus, the element $k$ of $A$ satisfies $k-1\notin A$ or $\left[
k,n-1\right]  \subseteq A$. In other words, $k\in A_{n}^{\circ}$ (by the
definition of $A_{n}^{\circ}$). This contradicts $k\notin A_{n}^{\circ}$.
Thus, we have found a contradiction in Case 1.
\par
Let us now consider Case 2. In this case, we have $k+1=n$. Hence, $k=n-1$, so
that $\left[  k,n-1\right]  =\left\{  k\right\}  \subseteq A$ (since $k\in
A$). Thus, the element $k$ of $A$ satisfies $k-1\notin A$ or $\left[
k,n-1\right]  \subseteq A$. In other words, $k\in A_{n}^{\circ}$ (by the
definition of $A_{n}^{\circ}$). This contradicts $k\notin A_{n}^{\circ}$.
Thus, we have found a contradiction in Case 2.
\par
We have therefore found a contradiction in each of the two Cases 1 and 2.
Thus, we always get a contradiction, so our assumption must have been wrong.
Qed.}. Hence, we have found a $k\in\left[  n-1\right]  $ satisfying%
\begin{align*}
\operatorname{Des}K  &  =\operatorname{Des}J\cup\left\{  k\right\}
,\ \ \ \ \ \ \ \ \ \ k\notin\operatorname{Des}J,\\
k-1  &  \in\operatorname{Des}J\ \ \ \ \ \ \ \ \ \ \text{and}%
\ \ \ \ \ \ \ \ \ \ k+1\notin\operatorname{Des}J\cup\left\{  n\right\}  .
\end{align*}
Therefore, Claim 1 yields $J\rightarrow K$. Thus, Claim 2 yields
$\operatorname{Epk}J=\operatorname{Epk}K$. Claim 4 therefore yields $\left(
\operatorname{Des}J\right)  _{n}^{\circ}=\left(  \operatorname{Des}K\right)
_{n}^{\circ}=A_{n}^{\circ}$ (since $\operatorname{Des}K=A$). Thus,%
\begin{align*}
\left\vert \underbrace{\left(  \operatorname{Des}J\right)  }_{=A\setminus
\left\{  k\right\}  }\setminus\underbrace{\left(  \operatorname{Des}J\right)
_{n}^{\circ}}_{=A_{n}^{\circ}}\right\vert  &  =\left\vert \underbrace{\left(
A\setminus\left\{  k\right\}  \right)  \setminus A_{n}^{\circ}}_{=\left(
A\setminus A_{n}^{\circ}\right)  \setminus\left\{  k\right\}  }\right\vert
=\left\vert \left(  A\setminus A_{n}^{\circ}\right)  \setminus\left\{
k\right\}  \right\vert \\
&  =\left\vert A\setminus A_{n}^{\circ}\right\vert
-1\ \ \ \ \ \ \ \ \ \ \left(  \text{since }k\in A\setminus A_{n}^{\circ
}\right) \\
&  <\left\vert A\setminus A_{n}^{\circ}\right\vert =\left\vert \left(
\operatorname{Des}K\right)  \setminus\left(  \operatorname{Des}K\right)
_{n}^{\circ}\right\vert
\end{align*}
(since $A=\operatorname{Des}K$). Thus, the induction hypothesis shows that
$\operatorname*{Comp}\left(  \left(  \operatorname{Des}J\right)  _{n}^{\circ
}\right)  \overset{\ast}{\rightarrow}J$. Combining this with $J\rightarrow K$,
we obtain $\operatorname*{Comp}\left(  \left(  \operatorname{Des}J\right)
_{n}^{\circ}\right)  \overset{\ast}{\rightarrow}K$ (since $\overset{\ast
}{\rightarrow}$ is the transitive-and-reflexive closure of the relation
$\rightarrow$). In light of $\left(  \operatorname{Des}J\right)  _{n}^{\circ
}=\left(  \operatorname{Des}K\right)  _{n}^{\circ}$, this rewrites as
$\operatorname*{Comp}\left(  \left(  \operatorname{Des}K\right)  _{n}^{\circ
}\right)  \overset{\ast}{\rightarrow}K$. Thus, Claim 5 is proven by induction.]

Now, let $\mathcal{K}^{\prime}$ be the $\mathbb{Q}$-vector subspace of
$\operatorname*{QSym}$ spanned by all differences of the form $F_{J}-F_{K}$,
where $J$ and $K$ are two compositions satisfying $J\rightarrow K$.

\begin{statement}
\textit{Claim 6:} Let $J$ and $K$ be two compositions such that
$J\overset{\ast}{\rightarrow}K$. Then, $F_{J}-F_{K}\in\mathcal{K}^{\prime}$.
\end{statement}

[\textit{Proof of Claim 6:} We have $J\overset{\ast}{\rightarrow}K$. By the
definition of the relation $\overset{\ast}{\rightarrow}$, this means that
there exists a sequence $\left(  L_{0},L_{1},\ldots,L_{\ell}\right)  $ of
compositions satisfying $L_{0}=J$ and $L_{\ell}=K$ and $L_{0}\rightarrow
L_{1}\rightarrow\cdots\rightarrow L_{\ell}$. Consider this sequence. For each
$i\in\left\{  0,1,\ldots,\ell-1\right\}  $, we have $L_{i}\rightarrow L_{i+1}$
and thus $F_{L_{i}}-F_{L_{i+1}}\in\mathcal{K}^{\prime}$ (by the definition of
$\mathcal{K}^{\prime}$). Therefore, $\sum_{i=0}^{\ell-1}\left(  F_{L_{i}%
}-F_{L_{i+1}}\right)  \in\mathcal{K}^{\prime}$. In light of%
\begin{align*}
\sum_{i=0}^{\ell-1}\left(  F_{L_{i}}-F_{L_{i+1}}\right)   &  =F_{L_{0}%
}-F_{L_{\ell}}\ \ \ \ \ \ \ \ \ \ \left(  \text{by the telescope
principle}\right) \\
&  =F_{J}-F_{K}\ \ \ \ \ \ \ \ \ \ \left(  \text{since }L_{0}=J\text{ and
}L_{\ell}=K\right)  ,
\end{align*}
this rewrites as $F_{J}-F_{K}\in\mathcal{K}^{\prime}$. This proves Claim 6.]

\begin{statement}
\textit{Claim 7:} We have $\mathcal{K}_{\operatorname{Epk}}\subseteq
\mathcal{K}^{\prime}$.
\end{statement}

[\textit{Proof of Claim 7:} Recall that $\mathcal{K}_{\operatorname{Epk}}$ is
the $\mathbb{Q}$-vector subspace of $\operatorname*{QSym}$ spanned by all
elements of the form $F_{J}-F_{K}$, where $J$ and $K$ are two
$\operatorname{Epk}$-equivalent compositions. Thus, it suffices to show that
if $J$ and $K$ are two $\operatorname{Epk}$-equivalent compositions, then
$F_{J}-F_{K}\in\mathcal{K}^{\prime}$.

So let $J$ and $K$ be two $\operatorname{Epk}$-equivalent compositions. We
must prove that $F_{J}-F_{K}\in\mathcal{K}^{\prime}$.

The compositions $J$ and $K$ are $\operatorname{Epk}$-equivalent; in other
words, they have the same size and satisfy $\operatorname{Epk}%
J=\operatorname{Epk}K$. Let $n=\left\vert J\right\vert =\left\vert
K\right\vert $. (This is well-defined, since the compositions $J$ and $K$ have
the same size.)

Claim 4 yields $\left(  \operatorname{Des}J\right)  _{n}^{\circ}=\left(
\operatorname{Des}K\right)  _{n}^{\circ}$. But Claim 5 yields
$\operatorname*{Comp}\left(  \left(  \operatorname{Des}K\right)  _{n}^{\circ
}\right)  \overset{\ast}{\rightarrow}K$. Hence, Claim 6 (applied to
$\operatorname*{Comp}\left(  \left(  \operatorname{Des}K\right)  _{n}^{\circ
}\right)  $ instead of $J$) shows that \newline
$F_{\operatorname*{Comp}\left(  \left(
\operatorname{Des}K\right)  _{n}^{\circ}\right)  }-F_{K}\in\mathcal{K}%
^{\prime}$. The same argument (applied to $J$ instead of $K$) shows that
$F_{\operatorname*{Comp}\left(  \left(  \operatorname{Des}J\right)
_{n}^{\circ}\right)  }-F_{J}\in\mathcal{K}^{\prime}$. Now,%
\begin{align*}
&  \left(  F_{\operatorname*{Comp}\left(  \left(  \operatorname{Des}K\right)
_{n}^{\circ}\right)  }-F_{K}\right)  -\left(  F_{\operatorname*{Comp}\left(
\left(  \operatorname{Des}J\right)  _{n}^{\circ}\right)  }-F_{J}\right) \\
&  =\left(  F_{\operatorname*{Comp}\left(  \left(  \operatorname{Des}%
K\right)  _{n}^{\circ}\right)  }-\underbrace{F_{\operatorname*{Comp}\left(
\left(  \operatorname{Des}J\right)  _{n}^{\circ}\right)  }}%
_{\substack{=F_{\operatorname*{Comp}\left(  \left(  \operatorname{Des}%
K\right)  _{n}^{\circ}\right)  }\\\text{(since }\left(  \operatorname{Des}%
J\right)  _{n}^{\circ}=\left(  \operatorname{Des}K\right)  _{n}^{\circ
}\text{)}}}\right)  +F_{J}-F_{K}\\
&  =\underbrace{\left(  F_{\operatorname*{Comp}\left(  \left(
\operatorname{Des}K\right)  _{n}^{\circ}\right)  }-F_{\operatorname*{Comp}%
\left(  \left(  \operatorname{Des}K\right)  _{n}^{\circ}\right)  }\right)
}_{=0}+F_{J}-F_{K}=F_{J}-F_{K},
\end{align*}
so that%
\[
F_{J}-F_{K}=\underbrace{\left(  F_{\operatorname*{Comp}\left(  \left(
\operatorname{Des}K\right)  _{n}^{\circ}\right)  }-F_{K}\right)  }%
_{\in\mathcal{K}^{\prime}}-\underbrace{\left(  F_{\operatorname*{Comp}\left(
\left(  \operatorname{Des}J\right)  _{n}^{\circ}\right)  }-F_{J}\right)
}_{\in\mathcal{K}^{\prime}}\in\mathcal{K}^{\prime}-\mathcal{K}^{\prime
}\subseteq\mathcal{K}^{\prime}.
\]
This proves Claim 7.]

\begin{statement}
\textit{Claim 8:} We have $\mathcal{K}^{\prime}\subseteq\mathcal{K}%
_{\operatorname{Epk}}$.
\end{statement}

[\textit{Proof of Claim 8:} Recall that $\mathcal{K}^{\prime}$ is the
$\mathbb{Q}$-vector subspace of $\operatorname*{QSym}$ spanned by all
differences of the form $F_{J}-F_{K}$, where $J$ and $K$ are two compositions
satisfying $J\rightarrow K$. Thus, it suffices to show that if $J$ and $K$ are
two compositions satisfying $J\rightarrow K$, then $F_{J}-F_{K}\in
\mathcal{K}_{\operatorname{Epk}}$.

So let $J$ and $K$ be two compositions satisfying $J\rightarrow K$. We must
prove that $F_{J}-F_{K}\in\mathcal{K}_{\operatorname{Epk}}$.

We have $J\rightarrow K$. Therefore, $\left|J\right| = \left|K\right|$ and
$\operatorname{Epk} J = \operatorname{Epk} K$ (by Claim 2).
Hence, the compositions $J$ and $K$ are $\operatorname{Epk}$-equivalent.
Thus, the definition of $\mathcal{K}_{\operatorname{Epk}}$ shows that
$F_{J}-F_{K} \in \mathcal{K}_{\operatorname{Epk}}$.
This proves Claim 8.]

Combining Claim 7 and Claim 8, we obtain $\mathcal{K}_{\operatorname{Epk}}
=\mathcal{K}^{\prime}$. Recalling the definition of $\mathcal{K}^{\prime}$,
we can rewrite this as follows: $\mathcal{K}_{\operatorname{Epk}}$ is the
$\mathbb{Q}$-vector subspace of $\operatorname*{QSym}$ spanned by all
differences of the form $F_{J}-F_{K}$, where $J$ and $K$ are two compositions
satisfying $J\rightarrow K$. This proves Proposition \ref{prop.K.Epk.F}.
\end{proof}
\end{verlong}

\subsection{An M-generating set of $\mathcal{K}_{\operatorname{Epk}}$}

Another characterization of the ideal $\mathcal{K}_{\operatorname{Epk}}$ of
$\operatorname*{QSym}$ can be obtained using the monomial basis of
$\operatorname*{QSym}$. Let us first recall how said basis is defined:

For any composition $\alpha=\left(  \alpha_{1},\alpha_{2},\ldots,\alpha_{\ell
}\right)  $, we let%
\[
M_{\alpha}=\sum_{i_{1}<i_{2}<\cdots<i_{\ell}}x_{i_{1}}^{\alpha_{1}}x_{i_{2}%
}^{\alpha_{2}}\cdots x_{i_{\ell}}^{\alpha_{\ell}}%
\]
(where the sum is over all strictly increasing $\ell$-tuples $\left(
i_{1},i_{2},\ldots,i_{\ell}\right)  $ of positive integers). This power series
$M_{\alpha}$ belongs to $\operatorname*{QSym}$. The family $\left(  M_{\alpha
}\right)  _{\alpha\text{ is a composition}}$ is a basis of the $\mathbb{Q}%
$-vector space $\operatorname*{QSym}$; it is called the \textit{monomial
basis} of $\operatorname*{QSym}$.

\begin{proposition}
\label{prop.K.Epk.M}If $J=\left(  j_{1},j_{2},\ldots,j_{m}\right)  $ and $K$
are two compositions, then we shall write $J\underset{M}{\rightarrow}K$ if
there exists an $\ell\in\left\{  2,3,\ldots,m\right\}  $ such that $j_{\ell
}>2$ and $K=\left(  j_{1},j_{2},\ldots,j_{\ell-1},2,j_{\ell}-2,j_{\ell
+1},j_{\ell+2},\ldots,j_{m}\right)  $. (In other words, we write
$J\underset{M}{\rightarrow}K$ if $K$ can be obtained from $J$ by
\textquotedblleft splitting\textquotedblright\ some entry $j_{\ell}>2$ into
two consecutive entries $2$ and $j_{\ell}-2$, provided that this entry was not
the first entry -- i.e., we had $\ell>1$ -- and that this entry was greater
than $2$.)

The ideal $\mathcal{K}_{\operatorname{Epk}}$ of $\operatorname*{QSym}$ is
spanned (as a $\mathbb{Q}$-vector space) by all sums of the form $M_{J}+M_{K}%
$, where $J$ and $K$ are two compositions satisfying
$J\underset{M}{\rightarrow}K$.
\end{proposition}

\begin{example}
We have $\left(  2,1,4,4\right)  \underset{M}{\rightarrow}\left(
2,1,2,2,4\right)  $, since the composition $\left(  2,1,2,2,4\right)  $ is
obtained from $\left(  2,1,4,4\right)  $ by splitting the third entry (which
is $4>2$) into two consecutive entries $2$ and $2$.

Similarly, $\left(  2,1,4,4\right)  \underset{M}{\rightarrow}\left(
2,1,4,2,2\right)  $ and $\left(  2,1,5,4\right)  \underset{M}{\rightarrow
}\left(  2,1,2,3,4\right)  $.

But we do not have $\left(  3,1\right)  \underset{M}{\rightarrow}\left(
2,1,1\right)  $, because splitting the first entry of the composition is not
allowed in the definition of the relation $\underset{M}{\rightarrow}$.

Two compositions $J$ and $K$ satisfying $J\underset{M}{\rightarrow}K$ must
necessarily satisfy $\left\vert J\right\vert =\left\vert K\right\vert $.

Here are all relations $\underset{M}{\rightarrow}$ between compositions of
size $4$:%
\[
\left(  1,3\right)  \underset{M}{\rightarrow}\left(  1,2,1\right)  .
\]

Here are all relations $\underset{M}{\rightarrow}$ between compositions of
size $5$:%
\begin{align*}
\left(  1,4\right)   &  \underset{M}{\rightarrow}\left(  1,2,2\right)  ,\\
\left(  1,3,1\right)   &  \underset{M}{\rightarrow}\left(  1,2,1,1\right)  ,\\
\left(  1,1,3\right)   &  \underset{M}{\rightarrow}\left(  1,1,2,1\right)  ,\\
\left(  2,3\right)   &  \underset{M}{\rightarrow}\left(  2,2,1\right)  .
\end{align*}
There are no relations $\underset{M}{\rightarrow}$ between compositions of
size $\leq3$.
\end{example}

\begin{verlong}
Before we start proving Proposition \ref{prop.K.Epk.M}, let us recall a basic
formula (\cite[(5.2.2)]{HopfComb}) that connects the monomial quasisymmetric
functions with the fundamental quasisymmetric functions:

\begin{proposition}
\label{prop.M-through-F.1}Let $n\in\mathbb{N}$. Let $\alpha$ be any
composition of $n$. Then,%
\[
M_{\alpha}=\sum_{\substack{\beta\text{ is a composition}\\\text{of }n\text{
that refines }\alpha}}\left(  -1\right)  ^{\ell\left(  \beta\right)
-\ell\left(  \alpha\right)  }F_{\beta}.
\]
Here, if $\gamma$ is any composition, then $\ell\left(  \gamma\right)  $
denotes the \textit{length} of $\gamma$ (that is, the number of entries of
$\gamma$).
\end{proposition}

\begin{proof}
[Proof of Proposition \ref{prop.M-through-F.1}.]This is
precisely \cite[(5.2.2)]{HopfComb}.
%(If the numbering shifts: This is the formula in Proposition 5.2.8.)
\end{proof}

\begin{proposition}
\label{prop.M-through-F.2}Let $n$ be a positive integer. Let $C$ be a subset
of $\left[  n-1\right]  $.

\textbf{(a)} Then,%
\[
M_{\operatorname*{Comp}C}=\sum_{B\supseteq C}\left(  -1\right)  ^{\left\vert
B\setminus C\right\vert }F_{\operatorname*{Comp}B}.
\]
(The bound variable $B$ in this sum and any similar sums is supposed to be a
subset of $\left[  n-1\right]  $; thus, the above sum ranges over all subsets
$B$ of $\left[  n-1\right]  $ satisfying $B\supseteq C$.)

\textbf{(b)} Let $k\in\left[  n-1\right]  $ be such that $k\notin C$. Then,%
\[
M_{\operatorname*{Comp}C}+M_{\operatorname*{Comp}\left(  C\cup\left\{
k\right\}  \right)  }=\sum_{\substack{B\supseteq C;\\k\notin B}}\left(
-1\right)  ^{\left\vert B\setminus C\right\vert }F_{\operatorname*{Comp}B}.
\]

\textbf{(c)} Let $k\in\left[  n-1\right]  $ be such that $k\notin C$ and
$k-1\notin C\cup\left\{  0\right\}  $. Then,%
\[
M_{\operatorname*{Comp}C}+M_{\operatorname*{Comp}\left(  C\cup\left\{
k\right\}  \right)  }=\sum_{\substack{B\supseteq C;\\k\notin B;\\k-1\notin
B}}\left(  -1\right)  ^{\left\vert B\setminus C\right\vert }\left(
F_{\operatorname*{Comp}B}-F_{\operatorname*{Comp}\left(  B\cup\left\{
k-1\right\}  \right)  }\right)  .
\]

\end{proposition}

\begin{proof}
[Proof of Proposition \ref{prop.M-through-F.2}.]
\textbf{(a)} If $U$ is any subset of $\left[  n-1\right]  $, then $\ell\left(
\operatorname*{Comp}U\right)  =%
\begin{cases}
\left\vert U\right\vert +1, & \text{if }n>0;\\
0, & \text{if }n=0
\end{cases}
$. This easy fact yields that every subset $B$ of $\left[  n-1\right]  $
satisfies%
\begin{align}
\ell\left(  \operatorname*{Comp}B\right)  -\ell\left(  \operatorname*{Comp}%
C\right)    & =%
\begin{cases}
\left\vert B\right\vert +1, & \text{if }n>0;\\
0, & \text{if }n=0
\end{cases}
-%
\begin{cases}
\left\vert C\right\vert +1, & \text{if }n>0;\\
0, & \text{if }n=0
\end{cases}
\nonumber\\
& =%
\begin{cases}
\left\vert B\right\vert -\left\vert C\right\vert , & \text{if }n>0;\\
0, & \text{if }n=0
\end{cases}
=\left\vert B\right\vert -\left\vert C\right\vert
\label{pf.prop.M-through-F.2.a.l-l}%
\end{align}
(indeed, if $n=0$, then both $B$ and $C$ must be the empty set, and thus
$0=\left\vert B\right\vert -\left\vert C\right\vert $ holds in this case).

Proposition \ref{prop.M-through-F.1} (applied to $\alpha=\operatorname*{Comp}%
C$) yields%
\begin{align*}
& M_{\operatorname*{Comp}C}\\
& =\sum_{\substack{\beta\text{ is a composition}\\\text{of }n\text{ that
refines }\operatorname*{Comp}C}}\left(  -1\right)  ^{\ell\left(  \beta\right)
-\ell\left(  \operatorname*{Comp}C\right)  }F_{\beta}\\
& =\underbrace{\sum_{\substack{B\subseteq\left[  n-1\right]
;\\\operatorname*{Comp}B\text{ refines }\operatorname*{Comp}C}}}%
_{\substack{=\sum_{\substack{B\subseteq\left[  n-1\right]  ;\\B\supseteq
C}}\\\text{(because for a subset }B\text{ of }\left[  n-1\right]
\text{,}\\\text{we have }\left(  \operatorname*{Comp}B\text{ refines
}\operatorname*{Comp}C\right)  \\\text{if and only if }B\supseteq C\text{)}%
}}\underbrace{\left(  -1\right)  ^{\ell\left(  \operatorname*{Comp}B\right)
-\ell\left(  \operatorname*{Comp}C\right)  }}_{\substack{=\left(  -1\right)
^{\left\vert B\right\vert -\left\vert C\right\vert }\\\text{(by
(\ref{pf.prop.M-through-F.2.a.l-l}))}}}F_{\operatorname*{Comp}B}\\
& \ \ \ \ \ \ \ \ \ \ \left(
\begin{array}
[c]{c}%
\text{here, we have substituted }\operatorname*{Comp}B\text{ for }\beta\text{,
since}\\
\text{the map }\operatorname*{Comp}:\left\{  \text{subsets of }\left[
n-1\right]  \right\}  \rightarrow\left\{  \text{compositions of }n\right\}  \\
\text{is a bijection}%
\end{array}
\right)  \\
& =\underbrace{\sum_{\substack{B\subseteq\left[  n-1\right]  ;\\B\supseteq
C}}}_{=\sum_{B\supseteq C}}\underbrace{\left(  -1\right)  ^{\left\vert
B\right\vert -\left\vert C\right\vert }}_{\substack{=\left(  -1\right)
^{\left\vert B\setminus C\right\vert }\\\text{(since }B\supseteq C\text{)}%
}}F_{\operatorname*{Comp}B}=\sum_{B\supseteq C}\left(  -1\right)  ^{\left\vert
B\setminus C\right\vert }F_{\operatorname*{Comp}B}.
\end{align*}
This proves Proposition \ref{prop.M-through-F.2} \textbf{(a)}.

\textbf{(b)} Proposition \ref{prop.M-through-F.2} \textbf{(a)} (applied to
$C\cup\left\{  k\right\}  $ instead of $C$) yields%
\begin{align}
M_{\operatorname*{Comp}\left(  C\cup\left\{  k\right\}  \right)  }  &
=\underbrace{\sum_{B\supseteq C\cup\left\{  k\right\}  }}_{=\sum
_{\substack{B\supseteq C;\\k\in B}}}\underbrace{\left(  -1\right)
^{\left\vert B\setminus\left(  C\cup\left\{  k\right\}  \right)  \right\vert
}}_{\substack{=\left(  -1\right)  ^{\left\vert B\setminus C\right\vert
-1}\\\text{(since }\left\vert B\setminus\left(  C\cup\left\{  k\right\}
\right)  \right\vert =\left\vert B\setminus C\right\vert -1\\\text{(since
}k\notin C\text{))}}}F_{\operatorname*{Comp}B}=\sum_{\substack{B\supseteq
C;\\k\in B}}\underbrace{\left(  -1\right)  ^{\left\vert B\setminus
C\right\vert -1}}_{=-\left(  -1\right)  ^{\left\vert B\setminus C\right\vert
}}F_{\operatorname*{Comp}B}\nonumber\\
&  =-\sum_{\substack{B\supseteq C;\\k\in B}}\left(  -1\right)  ^{\left\vert
B\setminus C\right\vert }F_{\operatorname*{Comp}B}.
\label{pf.prop.M-through-F.2.b.1}%
\end{align}
But Proposition \ref{prop.M-through-F.2} \textbf{(a)} yields
\begin{align*}
M_{\operatorname*{Comp}C}  &  =\sum_{B\supseteq C}\left(  -1\right)
^{\left\vert B\setminus C\right\vert }F_{\operatorname*{Comp}B}\\
&  =\sum_{\substack{B\supseteq C;\\k\in B}}\left(  -1\right)  ^{\left\vert
B\setminus C\right\vert }F_{\operatorname*{Comp}B}+\sum_{\substack{B\supseteq
C;\\k\notin B}}\left(  -1\right)  ^{\left\vert B\setminus C\right\vert
}F_{\operatorname*{Comp}B}.
\end{align*}
Adding (\ref{pf.prop.M-through-F.2.b.1}) to this equality, we obtain%
\begin{align*}
&  M_{\operatorname*{Comp}C}+M_{\operatorname*{Comp}\left(  C\cup\left\{
k\right\}  \right)  }\\
&  =\sum_{\substack{B\supseteq C;\\k\in B}}\left(  -1\right)  ^{\left\vert
B\setminus C\right\vert }F_{\operatorname*{Comp}B}+\sum_{\substack{B\supseteq
C;\\k\notin B}}\left(  -1\right)  ^{\left\vert B\setminus C\right\vert
}F_{\operatorname*{Comp}B}+\left(  -\sum_{\substack{B\supseteq C;\\k\in
B}}\left(  -1\right)  ^{\left\vert B\setminus C\right\vert }%
F_{\operatorname*{Comp}B}\right) \\
&  =\sum_{\substack{B\supseteq C;\\k\notin B}}\left(  -1\right)  ^{\left\vert
B\setminus C\right\vert }F_{\operatorname*{Comp}B}.
\end{align*}
This proves Proposition \ref{prop.M-through-F.2} \textbf{(b)}.

\textbf{(c)} We have $k-1\notin C\cup\left\{  0\right\}  $. Thus, $k-1\notin
C$ and $k-1\neq0$. From $k-1\neq0$, we obtain $k-1\in\left[  n-1\right]  $.

The map%
\begin{align}
&  \left\{  B\subseteq\left[  n-1\right]  \ \mid\ B\supseteq C\text{ and
}k\notin B\text{ and }k-1\notin B\right\} \nonumber\\
&  \rightarrow\left\{  B\subseteq\left[  n-1\right]  \ \mid\ B\supseteq
C\text{ and }k\notin B\text{ and }k-1\in B\right\} \nonumber
\end{align}
sending each $B$ to $B\cup\left\{  k-1\right\}  $ is a bijection (this is easy
to check using the facts that $k-1\notin C$ and $k-1\in\left[  n-1\right]  $).
We shall denote this map by $\Phi$.

Proposition \ref{prop.M-through-F.2} \textbf{(b)} yields
\begin{align*}
&  M_{\operatorname*{Comp}C}+M_{\operatorname*{Comp}\left(  C\cup\left\{
k\right\}  \right)  }\\
&  =\sum_{\substack{B\supseteq C;\\k\notin B}}\left(  -1\right)  ^{\left\vert
B\setminus C\right\vert }F_{\operatorname*{Comp}B}\\
&  =\sum_{\substack{B\supseteq C;\\k\notin B;\\k-1\notin B}}\left(  -1\right)
^{\left\vert B\setminus C\right\vert }F_{\operatorname*{Comp}B}%
+\underbrace{\sum_{\substack{B\supseteq C;\\k\notin B;\\k-1\in B}}\left(
-1\right)  ^{\left\vert B\setminus C\right\vert }F_{\operatorname*{Comp}B}%
}_{\substack{=\sum_{\substack{B\supseteq C;\\k\notin B;\\k-1\notin B}}\left(
-1\right)  ^{\left\vert \left(  B\cup\left\{  k-1\right\}  \right)  \setminus
C\right\vert }F_{\operatorname*{Comp}\left(  B\cup\left\{  k-1\right\}
\right)  }\\\text{(here, we have substituted }B\cup\left\{  k-1\right\}
\text{ for }B\text{ in the sum}\\\text{(since the map }\Phi\text{ is a
bijection))}}}\\
&  =\sum_{\substack{B\supseteq C;\\k\notin B;\\k-1\notin B}}\left(  -1\right)
^{\left\vert B\setminus C\right\vert }F_{\operatorname*{Comp}B}+\sum
_{\substack{B\supseteq C;\\k\notin B;\\k-1\notin B}}\underbrace{\left(
-1\right)  ^{\left\vert \left(  B\cup\left\{  k-1\right\}  \right)  \setminus
C\right\vert }}_{\substack{=\left(  -1\right)  ^{\left\vert B\setminus
C\right\vert +1}\\\text{(since }\left\vert \left(  B\cup\left\{  k-1\right\}
\right)  \setminus C\right\vert =\left\vert B\setminus C\right\vert
+1\\\text{(since }k-1\notin B\text{ and }k-1\notin C\text{))}}%
}F_{\operatorname*{Comp}\left(  B\cup\left\{  k-1\right\}  \right)  }\\
&  =\sum_{\substack{B\supseteq C;\\k\notin B;\\k-1\notin B}}\left(  -1\right)
^{\left\vert B\setminus C\right\vert }F_{\operatorname*{Comp}B}+\sum
_{\substack{B\supseteq C;\\k\notin B;\\k-1\notin B}}\underbrace{\left(
-1\right)  ^{\left\vert B\setminus C\right\vert +1}}_{=-\left(  -1\right)
^{\left\vert B\setminus C\right\vert }}F_{\operatorname*{Comp}\left(
B\cup\left\{  k-1\right\}  \right)  }\\
&  =\sum_{\substack{B\supseteq C;\\k\notin B;\\k-1\notin B}}\left(  -1\right)
^{\left\vert B\setminus C\right\vert }F_{\operatorname*{Comp}B}-\sum
_{\substack{B\supseteq C;\\k\notin B;\\k-1\notin B}}\left(  -1\right)
^{\left\vert B\setminus C\right\vert }F_{\operatorname*{Comp}\left(
B\cup\left\{  k-1\right\}  \right)  }\\
&  =\sum_{\substack{B\supseteq C;\\k\notin B;\\k-1\notin B}}\left(  -1\right)
^{\left\vert B\setminus C\right\vert }\left(  F_{\operatorname*{Comp}%
B}-F_{\operatorname*{Comp}\left(  B\cup\left\{  k-1\right\}  \right)
}\right)  .
\end{align*}
This proves Proposition \ref{prop.M-through-F.2} \textbf{(c)}.
\end{proof}
\end{verlong}

\begin{verlong}
\begin{proof}
[Proof of Proposition \ref{prop.K.Epk.M}.]We shall use the notation
$\left\langle f_{i}\ \mid\ i\in I\right\rangle $ for the $\mathbb{Q}$-linear
span of a family $\left(  f_{i}\right)  _{i\in I}$ of elements of a
$\mathbb{Q}$-vector space.

Define a $\mathbb{Q}$-vector subspace $\mathcal{M}$ of $\operatorname*{QSym}$
by%
\[
\mathcal{M}=\left\langle M_{J}+M_{K}\ \mid\ J\text{ and }K\text{ are
compositions satisfying }J\underset{M}{\rightarrow}K\right\rangle .
\]
Then, our goal is to prove that $\mathcal{K}_{\operatorname{Epk}}%
=\mathcal{M}$.

We have%
\begin{align*}
\mathcal{M}  &  =\left\langle M_{J}+M_{K}\ \mid\ J\text{ and }K\text{ are
compositions satisfying }J\underset{M}{\rightarrow}K\right\rangle \\
&  =\sum_{n\in\mathbb{N}}\left\langle M_{J}+M_{K}\ \mid\ J\text{ and }K\text{
are compositions of }n\text{ satisfying }J\underset{M}{\rightarrow
}K\right\rangle
\end{align*}
(because if $J$ and $K$ are two compositions satisfying
$J\underset{M}{\rightarrow}K$, then $J$ and $K$ have the same size).

Consider the binary relation $\rightarrow$ defined in Proposition
\ref{prop.K.Epk.F}. Then, Proposition \ref{prop.K.Epk.F} yields%
\begin{align*}
\mathcal{K}_{\operatorname{Epk}}  &  =\left\langle F_{J}-F_{K}\ \mid\ J\text{
and }K\text{ are compositions satisfying }J\rightarrow K\right\rangle \\
&  =\sum_{n\in\mathbb{N}}\left\langle F_{J}-F_{K}\ \mid\ J\text{ and }K\text{
are compositions of }n\text{ satisfying }J\rightarrow K\right\rangle
\end{align*}
(because if $J$ and $K$ are two compositions satisfying $J\rightarrow K$, then
$J$ and $K$ have the same size).

Now, fix $n\in\mathbb{N}$. Let $\Omega$ be the set of all pairs $\left(
C,k\right)  $ in which $C$ is a subset of $\left[  n-1\right]  $ and $k$ is an
element of $\left[  n-1\right]  $ satisfying $k\notin C$, $k-1\in C$ and
$k+1\notin C\cup\left\{  n\right\}  $.

For every $\left(  C,k\right)  \in\Omega$, we define two elements
$\mathbf{m}_{C,k}$ and $\mathbf{f}_{C,k}$ of $\operatorname*{QSym}$ by%
\begin{align}
\mathbf{m}_{C,k}  &  =M_{\operatorname*{Comp}C}+M_{\operatorname*{Comp}\left(
C\cup\left\{  k+1\right\}  \right)  }\ \ \ \ \ \ \ \ \ \ \text{and}%
\label{pf.prop.K.Epk.M.mCk=}\\
\mathbf{f}_{C,k}  &  =F_{\operatorname*{Comp}C}-F_{\operatorname*{Comp}\left(
C\cup\left\{  k\right\}  \right)  }. \label{pf.prop.K.Epk.M.fCk=}%
\end{align}
\footnote{These two elements are well-defined, because both $C\cup\left\{
k\right\}  $ and $C\cup\left\{  k+1\right\}  $ are subsets of $\left[
n-1\right]  $ (since $k+1\notin C\cup\left\{  n\right\}  $ shows that $k+1\neq
n$).}

We have the following:

\begin{statement}
\textit{Claim 1:} We have%
\begin{align*}
&  \left\langle M_{J}+M_{K}\ \mid\ J\text{ and }K\text{ are compositions of
}n\text{ satisfying }J\underset{M}{\rightarrow}K\right\rangle \\
&  =\left\langle \mathbf{m}_{C,k}\ \mid\ \left(  C,k\right)  \in
\Omega\right\rangle .
\end{align*}

\end{statement}

[\textit{Proof of Claim 1:} It is easy to see that two subsets $C$ and $D$ of
$\left[  n-1\right]  $ satisfy $\operatorname*{Comp}C\underset{M}{\rightarrow
}\operatorname*{Comp}D$ if and only if there exists some $k\in\left[
n-1\right]  $ satisfying $D=C\cup\left\{  k+1\right\}  $, $k\notin C$, $k-1\in
C$ and $k+1\notin C\cup\left\{  n\right\}  $.\ \ \ \ \footnote{To prove this,
recall that \textquotedblleft splitting\textquotedblright\ an entry of a
composition $J$ into two consecutive entries (summing up to the original
entry) is always tantamount to adding a new element to $\operatorname{Des}J$.
It suffices to show that the conditions under which an entry of a composition
$J$ can be split in the definition of the relation $\underset{M}{\rightarrow}$
are precisely the conditions $k\notin C$, $k-1\in C$ and $k+1\notin
C\cup\left\{  n\right\}  $ on $C=\operatorname{Des}J$, where $k+1$ denotes
the new element that we are adding to $\operatorname{Des}J$. This is
straightforward.}. Thus,%
\begin{align*}
&  \left\langle M_{\operatorname*{Comp}C}+M_{\operatorname*{Comp}D}%
\ \mid\ C\text{ and }D\text{ are subsets of }\left[  n-1\right]  \right. \\
&  \ \ \ \ \ \ \ \ \ \ \ \ \ \ \ \ \ \ \ \ \left.  \text{satisfying
}\operatorname*{Comp}C\underset{M}{\rightarrow}\operatorname*{Comp}%
D\right\rangle \\
&  =\left\langle M_{\operatorname*{Comp}C}+M_{\operatorname*{Comp}D}%
\ \mid\ C\text{ and }D\text{ are subsets of }\left[  n-1\right]  \right. \\
&  \ \ \ \ \ \ \ \ \ \ \ \ \ \ \ \ \ \ \ \ \left.  \text{such that there
exists some }k\in\left[  n-1\right]  \right. \\
&  \ \ \ \ \ \ \ \ \ \ \ \ \ \ \ \ \ \ \ \ \left.  \text{satisfying }%
D=C\cup\left\{  k+1\right\}  \text{, }k\notin C\text{, }k-1\in C\text{ and
}k+1\notin C\cup\left\{  n\right\}  \right\rangle \\
&  =\left\langle M_{\operatorname*{Comp}C}+M_{\operatorname*{Comp}\left(
C\cup\left\{  k+1\right\}  \right)  }\ \mid\ C\subseteq\left[  n-1\right]
\text{ and }k\in\left[  n-1\right]  \right. \\
&  \ \ \ \ \ \ \ \ \ \ \ \ \ \ \ \ \ \ \ \ \left.  \text{are such that
}k\notin C\text{, }k-1\in C\text{ and }k+1\notin C\cup\left\{  n\right\}
\right\rangle \\
&  =\left\langle \underbrace{M_{\operatorname*{Comp}C}+M_{\operatorname*{Comp}%
\left(  C\cup\left\{  k+1\right\}  \right)  }}_{\substack{=\mathbf{m}%
_{C,k}\\\text{(by (\ref{pf.prop.K.Epk.M.mCk=}))}}}\ \mid\ \left(  C,k\right)
\in\Omega\right\rangle \\
&  \ \ \ \ \ \ \ \ \ \ \left(  \text{by the definition of }\Omega\right) \\
&  =\left\langle \mathbf{m}_{C,k}\ \mid\ \left(  C,k\right)  \in
\Omega\right\rangle .
\end{align*}
Now, recall that $\operatorname*{Comp}$ is a bijection between the subsets of
$\left[  n-1\right]  $ and the compositions of $n$. Hence,%
\begin{align*}
&  \left\langle M_{J}+M_{K}\ \mid\ J\text{ and }K\text{ are compositions of
}n\text{ satisfying }J\underset{M}{\rightarrow}K\right\rangle \\
&  =\left\langle M_{\operatorname*{Comp}C}+M_{\operatorname*{Comp}D}%
\ \mid\ C\text{ and }D\text{ are subsets of }\left[  n-1\right]  \right. \\
&  \ \ \ \ \ \ \ \ \ \ \ \ \ \ \ \ \ \ \ \ \left.  \text{satisfying
}\operatorname*{Comp}C\underset{M}{\rightarrow}\operatorname*{Comp}%
D\right\rangle \\
&  =\left\langle \mathbf{m}_{C,k}\ \mid\ \left(  C,k\right)  \in
\Omega\right\rangle .
\end{align*}
This proves Claim 1.]

\begin{statement}
\textit{Claim 2:} We have%
\begin{align*}
&  \left\langle F_{J}-F_{K}\ \mid\ J\text{ and }K\text{ are compositions of
}n\text{ satisfying }J\rightarrow K\right\rangle \\
&  =\left\langle \mathbf{f}_{C,k}\ \mid\ \left(  C,k\right)  \in
\Omega\right\rangle .
\end{align*}

\end{statement}

[\textit{Proof of Claim 2:} It is easy to see that two subsets $C$ and $D$ of
$\left[  n-1\right]  $ satisfy $\operatorname*{Comp}C\rightarrow
\operatorname*{Comp}D$ if and only if there exists some $k\in\left[
n-1\right]  $ satisfying $D=C\cup\left\{  k\right\}  $, $k\notin C$, $k-1\in
C$ and $k+1\notin C\cup\left\{  n\right\}  $.\ \ \ \ \footnote{To prove this,
recall that \textquotedblleft splitting\textquotedblright\ an entry of a
composition $J$ into two consecutive entries (summing up to the original
entry) is always tantamount to adding a new element to $\operatorname{Des}J$.
It suffices to show that the conditions under which an entry of a composition
$J$ can be split in the definition of the relation $\rightarrow$ are precisely
the conditions $k\notin C$, $k-1\in C$ and $k+1\notin C\cup\left\{  n\right\}
$ on $C=\operatorname{Des}J$, where $k$ denotes
the new element that we are adding to $\operatorname{Des}J$.
This is straightforward.}. Thus,%
\begin{align*}
&  \left\langle F_{\operatorname*{Comp}C}-F_{\operatorname*{Comp}D}%
\ \mid\ C\text{ and }D\text{ are subsets of }\left[  n-1\right]  \right. \\
&  \ \ \ \ \ \ \ \ \ \ \ \ \ \ \ \ \ \ \ \ \left.  \text{satisfying
}\operatorname*{Comp}C\rightarrow\operatorname*{Comp}D\right\rangle \\
&  =\left\langle F_{\operatorname*{Comp}C}-F_{\operatorname*{Comp}D}%
\ \mid\ C\text{ and }D\text{ are subsets of }\left[  n-1\right]  \right. \\
&  \ \ \ \ \ \ \ \ \ \ \ \ \ \ \ \ \ \ \ \ \left.  \text{such that there
exists some }k\in\left[  n-1\right]  \right. \\
&  \ \ \ \ \ \ \ \ \ \ \ \ \ \ \ \ \ \ \ \ \left.  \text{satisfying }%
D=C\cup\left\{  k\right\}  \text{, }k\notin C\text{, }k-1\in C\text{ and
}k+1\notin C\cup\left\{  n\right\}  \right\rangle \\
&  =\left\langle F_{\operatorname*{Comp}C}-F_{\operatorname*{Comp}\left(
C\cup\left\{  k\right\}  \right)  }\ \mid\ C\subseteq\left[  n-1\right]
\text{ and }k\in\left[  n-1\right]  \right. \\
&  \ \ \ \ \ \ \ \ \ \ \ \ \ \ \ \ \ \ \ \ \left.  \text{are such that
}k\notin C\text{, }k-1\in C\text{ and }k+1\notin C\cup\left\{  n\right\}
\right\rangle \\
&  =\left\langle \underbrace{F_{\operatorname*{Comp}C}-F_{\operatorname*{Comp}%
\left(  C\cup\left\{  k\right\}  \right)  }}_{\substack{=\mathbf{f}%
_{C,k}\\\text{(by (\ref{pf.prop.K.Epk.M.fCk=}))}}}\ \mid\ \left(  C,k\right)
\in\Omega\right\rangle \\
&  \ \ \ \ \ \ \ \ \ \ \left(  \text{by the definition of }\Omega\right) \\
&  =\left\langle \mathbf{f}_{C,k}\ \mid\ \left(  C,k\right)  \in
\Omega\right\rangle .
\end{align*}
Now, recall that $\operatorname*{Comp}$ is a bijection between the subsets of
$\left[  n-1\right]  $ and the compositions of $n$. Hence,%
\begin{align*}
&  \left\langle F_{J}-F_{K}\ \mid\ J\text{ and }K\text{ are compositions of
}n\text{ satisfying }J\rightarrow K\right\rangle \\
&  =\left\langle F_{\operatorname*{Comp}C}-F_{\operatorname*{Comp}D}%
\ \mid\ C\text{ and }D\text{ are subsets of }\left[  n-1\right]  \right. \\
&  \ \ \ \ \ \ \ \ \ \ \ \ \ \ \ \ \ \ \ \ \left.  \text{satisfying
}\operatorname*{Comp}C\rightarrow\operatorname*{Comp}D\right\rangle \\
&  =\left\langle \mathbf{f}_{C,k}\ \mid\ \left(  C,k\right)  \in
\Omega\right\rangle .
\end{align*}
This proves Claim 2.]

We define a partial order on the set $\Omega$ by setting%
\[
\left(  B,k\right)  \geq\left(  C,\ell\right)  \ \ \ \ \ \ \ \ \ \ \text{if
and only if}\ \ \ \ \ \ \ \ \ \ \left(  k=\ell\text{ and }B\supseteq C\right)
.
\]
Thus, $\Omega$ is a finite poset.

\begin{statement}
\textit{Claim 3:} For every $\left(  C,\ell\right)  \in\Omega$, we have%
\[
\mathbf{m}_{C,\ell}=\sum_{\substack{\left(  B,k\right)  \in\Omega;\\\left(
B,k\right)  \geq\left(  C,\ell\right)  }}\left(  -1\right)  ^{\left\vert
B\setminus C\right\vert }\mathbf{f}_{B,k}.
\]

\end{statement}

[\textit{Proof of Claim 3:} Let $\left(  C,\ell\right)  \in\Omega$. Thus, $C$
is a subset of $\left[  n-1\right]  $ and $\ell$ is an element of $\left[
n-1\right]  $ satisfying $\ell\notin C$, $\ell-1\in C$ and $\ell+1\notin
C\cup\left\{  n\right\}  $. From $\ell+1\notin C\cup\left\{  n\right\}  $, we
obtain $\ell+1\notin C$ and $\ell+1\neq n$. From $\ell+1\neq n$, we obtain
$\ell+1\in\left[  n-1\right]  $. Also, $\left(  \ell+1\right)  -1=\ell\notin
C\cup\left\{  0\right\}  $ (since $\ell\notin C$ and $\ell\neq0$). Thus,
Proposition \ref{prop.M-through-F.2} \textbf{(c)} (applied to $k=\ell+1$)
yields\footnote{Here and in the following, the bound variable $B$ in a sum
always is understood to be a subset of $\left[  n-1\right]  $.}%
\[
M_{\operatorname*{Comp}C}+M_{\operatorname*{Comp}\left(  C\cup\left\{
\ell+1\right\}  \right)  }=\sum_{\substack{B\supseteq C;\\\ell+1\notin
B;\\\ell\notin B}}\left(  -1\right)  ^{\left\vert B\setminus C\right\vert
}\left(  F_{\operatorname*{Comp}B}-F_{\operatorname*{Comp}\left(
B\cup\left\{  \ell\right\}  \right)  }\right)  .
\]

But every $B\subseteq\left[  n-1\right]  $ satisfying $B\supseteq C$ must
satisfy $\ell-1\in B$ (since $\ell-1\in C\subseteq B$). Hence, we can
manipulate summation signs as follows:%
\begin{align}
\sum_{\substack{B\supseteq C;\\\ell+1\notin B;\\\ell\notin B}}  &
=\sum_{\substack{B\supseteq C;\\\ell+1\notin B;\\\ell\notin B;\\\ell-1\in
B}}=\sum_{\substack{B\supseteq C;\\\ell+1\notin B\cup\left\{  n\right\}
;\\\ell\notin B;\\\ell-1\in B}}\ \ \ \ \ \ \ \ \ \ \left(
\begin{array}
[c]{c}%
\text{since }\ell+1\notin B\text{ is equivalent to }\ell+1\notin B\cup\left\{
n\right\} \\
\text{(because }\ell+1\neq n\text{)}%
\end{array}
\right) \nonumber\\
&  =\sum_{\substack{B\supseteq C;\\\ell\notin B;\\\ell-1\in B;\\\ell+1\notin
B\cup\left\{  n\right\}  }}\nonumber\\
&  =\sum_{\substack{B\supseteq C;\\\left(  B,\ell\right)  \in\Omega
}}\ \ \ \ \ \ \ \ \ \ \left(
\begin{array}
[c]{c}%
\text{since the}\\
\text{condition }\left(  \ell\notin B\text{, }\ell-1\in B\text{ and }%
\ell+1\notin B\cup\left\{  n\right\}  \right) \\
\text{on a subset }B\text{ of }\left[  n-1\right]  \text{ is equivalent to
}\left(  B,\ell\right)  \in\Omega\\
\text{(by the definition of }\Omega\text{)}%
\end{array}
\right) \nonumber\\
&  =\sum_{\substack{\left(  B,\ell\right)  \in\Omega;\\B\supseteq C}%
}=\sum_{\substack{\left(  B,k\right)  \in\Omega;\\k=\ell;\\B\supseteq C}%
}=\sum_{\substack{\left(  B,k\right)  \in\Omega;\\\left(  B,k\right)
\geq\left(  C,\ell\right)  }} \label{pf.prop.K.Epk.M.sum-man}%
\end{align}
(since the condition $\left(  k=\ell\text{ and }B\supseteq C\right)  $ on a
$\left(  B,k\right)  \in\Omega$ is equivalent to $\left(  B,k\right)
\geq\left(  C,\ell\right)  $ (by the definition of the partial order on
$\Omega$)).

Now, the definition of $\mathbf{m}_{C,\ell}$ yields%
\begin{align*}
\mathbf{m}_{C,\ell}  &  =M_{\operatorname*{Comp}C}+M_{\operatorname*{Comp}%
\left(  C\cup\left\{  \ell+1\right\}  \right)  }\\
&  =\underbrace{\sum_{\substack{B\supseteq C;\\\ell+1\notin B;\\\ell\notin
B}}}_{\substack{=\sum_{\substack{\left(  B,k\right)  \in\Omega;\\\left(
B,k\right)  \geq\left(  C,\ell\right)  }}\\\text{(by
(\ref{pf.prop.K.Epk.M.sum-man}))}}}\left(  -1\right)  ^{\left\vert B\setminus
C\right\vert }\left(  F_{\operatorname*{Comp}B}-F_{\operatorname*{Comp}\left(
B\cup\left\{  \ell\right\}  \right)  }\right) \\
&  =\sum_{\substack{\left(  B,k\right)  \in\Omega;\\\left(  B,k\right)
\geq\left(  C,\ell\right)  }}\left(  -1\right)  ^{\left\vert B\setminus
C\right\vert }\left(  F_{\operatorname*{Comp}B}%
-\underbrace{F_{\operatorname*{Comp}\left(  B\cup\left\{  \ell\right\}
\right)  }}_{\substack{=F_{\operatorname*{Comp}\left(  B\cup\left\{
k\right\}  \right)  }\\\text{(since }\ell=k\\\text{(since }\left(  B,k\right)
\geq\left(  C,\ell\right)  \text{ and thus }k=\ell\text{))}}}\right) \\
&  =\sum_{\substack{\left(  B,k\right)  \in\Omega;\\\left(  B,k\right)
\geq\left(  C,\ell\right)  }}\left(  -1\right)  ^{\left\vert B\setminus
C\right\vert }\underbrace{\left(  F_{\operatorname*{Comp}B}%
-F_{\operatorname*{Comp}\left(  B\cup\left\{  k\right\}  \right)  }\right)
}_{\substack{=\mathbf{f}_{B,k}\\\text{(since (\ref{pf.prop.K.Epk.M.fCk=})
yields}\\\mathbf{f}_{B,k}=F_{\operatorname*{Comp}B}-F_{\operatorname*{Comp}%
\left(  B\cup\left\{  k\right\}  \right)  }\text{)}}}\\
&  =\sum_{\substack{\left(  B,k\right)  \in\Omega;\\\left(  B,k\right)
\geq\left(  C,\ell\right)  }}\left(  -1\right)  ^{\left\vert B\setminus
C\right\vert }\mathbf{f}_{B,k}.
\end{align*}
This proves Claim 3.]

Now, Claim 3 shows that the family $\left(  \mathbf{m}_{C,k}\right)  _{\left(
C,k\right)  \in\Omega}$ expands triangularly with respect to the family
$\left(  \mathbf{f}_{C,k}\right)  _{\left(  C,k\right)  \in\Omega}$ with
respect to the poset structure on $\Omega$. Moreover, the expansion is
\textbf{uni}triangular (because if $\left(  B,k\right)  =\left(
C,\ell\right)  $, then $B=C$ and thus $\left(  -1\right)  ^{\left\vert
B\setminus C\right\vert }=\left(  -1\right)  ^{\left\vert C\setminus
C\right\vert }=\left(  -1\right)  ^{0}=1$) and thus invertibly triangular
(this means that the diagonal entries are invertible). Therefore, by a
standard fact from linear algebra (see, e.g., \cite[Corollary 11.1.19
\textbf{(b)}]{HopfComb}), we conclude that the span of the family $\left(
\mathbf{m}_{C,k}\right)  _{\left(  C,k\right)  \in\Omega}$ equals the span of
the family $\left(  \mathbf{f}_{C,k}\right)  _{\left(  C,k\right)  \in\Omega}%
$. In other words,%
\[
\left\langle \mathbf{m}_{C,k}\ \mid\ \left(  C,k\right)  \in\Omega
\right\rangle =\left\langle \mathbf{f}_{C,k}\ \mid\ \left(  C,k\right)
\in\Omega\right\rangle .
\]
Now, Claim 1 yields%
\begin{align*}
&  \left\langle M_{J}+M_{K}\ \mid\ J\text{ and }K\text{ are compositions of
}n\text{ satisfying }J\underset{M}{\rightarrow}K\right\rangle \\
&  =\left\langle \mathbf{m}_{C,k}\ \mid\ \left(  C,k\right)  \in
\Omega\right\rangle =\left\langle \mathbf{f}_{C,k}\ \mid\ \left(  C,k\right)
\in\Omega\right\rangle \\
&  =\left\langle F_{J}-F_{K}\ \mid\ J\text{ and }K\text{ are compositions of
}n\text{ satisfying }J\rightarrow K\right\rangle
\end{align*}
(by Claim 2).

Now, forget that we fixed $n$. We thus have proven that
\begin{align*}
&  \left\langle M_{J}+M_{K}\ \mid\ J\text{ and }K\text{ are compositions of
}n\text{ satisfying }J\underset{M}{\rightarrow}K\right\rangle \\
&  =\left\langle F_{J}-F_{K}\ \mid\ J\text{ and }K\text{ are compositions of
}n\text{ satisfying }J\rightarrow K\right\rangle
\end{align*}
for each $n\in\mathbb{N}$. Thus,%
\begin{align*}
&  \sum_{n\in\mathbb{N}}\left\langle M_{J}+M_{K}\ \mid\ J\text{ and }K\text{
are compositions of }n\text{ satisfying }J\underset{M}{\rightarrow
}K\right\rangle \\
&  =\sum_{n\in\mathbb{N}}\left\langle F_{J}-F_{K}\ \mid\ J\text{ and }K\text{
are compositions of }n\text{ satisfying }J\rightarrow K\right\rangle .
\end{align*}
In light of%
\[
\mathcal{K}_{\operatorname{Epk}}=\sum_{n\in\mathbb{N}}\left\langle
F_{J}-F_{K}\ \mid\ J\text{ and }K\text{ are compositions of }n\text{
satisfying }J\rightarrow K\right\rangle
\]
and%
\[
\mathcal{M}=\sum_{n\in\mathbb{N}}\left\langle M_{J}+M_{K}\ \mid\ J\text{ and
}K\text{ are compositions of }n\text{ satisfying }J\underset{M}{\rightarrow
}K\right\rangle ,
\]
this rewrites as $\mathcal{M}=\mathcal{K}_{\operatorname{Epk}}$. In other
words, $\mathcal{K}_{\operatorname{Epk}}=\mathcal{M}$. This proves
Proposition \ref{prop.K.Epk.M}.
\end{proof}
\end{verlong}

\begin{question}
It is worth analyzing the kernels of other known descent statistics
(shuffle-compatible or not).
Let us say that a descent statistic $\operatorname{st}$ is
\textit{M-binomial}
if its kernel $\mathcal{K}_{\operatorname{st}}$ can be spanned by
elements of the form $\lambda M_{J}+\mu M_{K}$ with
$\lambda, \mu \in \mathbb{Q}$ and compositions $J,K$.
Then, Proposition \ref{prop.K.Epk.M} yields that $\operatorname{Epk}$
is M-binomial.
It is easy to see that the statistics $\operatorname{Des}$ and
$\operatorname{des}$ are M-binomial as well.
Computations using SageMath suggest that the statistics
$\operatorname{Lpk}$, $\operatorname{Rpk}$, $\operatorname*{Pk}$,
$\operatorname{Val}$,
$\operatorname{pk}$, $\operatorname{lpk}$, $\operatorname{rpk}$
and $\operatorname{val}$
(see \cite{part1} for some of their definitions) are M-binomial, too
(at least for compositions of size $\leq 9$);
this would be nice to prove.
On the other hand, the statistics
$\operatorname{maj}$, $\left(\operatorname{des},\operatorname{maj}\right)$
and $\left(\operatorname{val},\operatorname{des}\right)$ (again, see
\cite{part1} for definitions) are not M-binomial.
\end{question}

\begin{verlong}
\subsection{\label{subsect.K.pfAsh}Appendix: Proof of Proposition
\ref{prop.Ast.alg} and Theorem \ref{thm.4.3}}

Let us now give proofs of Proposition \ref{prop.Ast.alg} and Theorem
\ref{thm.4.3}, which we have promised above. We will mostly rely on Lemma
\ref{lem.K.fist} and on Proposition \ref{prop.4.1.rewr}.

For the rest of Subsection \ref{subsect.K.pfAsh}, we shall make the following conventions:

\begin{convention}
Let $\operatorname{st}$ be a permutation statistic. For each permutation
$\pi$, let $\left[  \pi\right]  _{\operatorname{st}}$ denote the
$\operatorname{st}$-equivalence class of $\pi$. Let $\mathcal{A}%
_{\operatorname{st}}$ be the free $\mathbb{Q}$-vector space whose basis is
the set of all $\operatorname{st}$-equivalence classes of permutations. (This
is well-defined whether or not $\operatorname{st}$ is shuffle-compatible.)
\end{convention}

\begin{proof}
[Proof of Proposition \ref{prop.Ast.alg}.]A \textit{magmatic algebra} shall
mean a $\mathbb{Q}$-vector space equipped with a binary operation which is
written as multiplication (i.e., we write $ab$ for the image of a pair
$\left(  a,b\right)  $ under this operation), but is not required to be
associative (or have a unity). An (actual, i.e., associative unital) algebra
is thus a magmatic algebra whose multiplication is associative and has a
unity. In particular, any actual algebra is a magmatic algebra. A
\textit{magmatic algebra homomorphism} is a $\mathbb{Q}$-linear map between
two magmatic algebras that preserves the multiplication.

We make $\mathcal{A}_{\operatorname{st}}$ into a magmatic algebra by setting%
\begin{equation}
\left[  \pi\right]  _{\operatorname{st}}\left[  \sigma\right]
_{\operatorname{st}}=\sum_{\tau\in S\left(  \pi,\sigma\right)  }\left[
\tau\right]  _{\operatorname{st}} \label{pf.prop.Ast.alg.algdef}%
\end{equation}
for any two disjoint permutations $\pi$ and $\sigma$. This is well-defined,
because the right-hand side of (\ref{pf.prop.Ast.alg.algdef}) depends only on
the $\operatorname{st}$-equivalence classes $\left[  \pi\right]
_{\operatorname{st}}$ and $\left[  \sigma\right]  _{\operatorname{st}}$
rather than on the permutations $\pi$ and $\sigma$ themselves (this is because
$\operatorname{st}$ is shuffle-compatible).

Define a $\mathbb{Q}$-linear map $p:\operatorname*{QSym}\rightarrow
\mathcal{A}_{\operatorname{st}}$ by requiring that%
\begin{align*}
p\left(  F_{L}\right)   &  =\left[  \pi\right]  _{\operatorname{st}}
\ \ \ \ \ \ \ \ \ \ \text{for every composition }L\text{ and every}\\
&  \ \ \ \ \ \ \ \ \ \ \ \ \ \ \ \ \ \ \ \ \text{permutation }\pi\text{ with
}\operatorname*{Comp}\pi=L.
\end{align*}
This is well-defined, because for any given composition $L$, any two
permutations $\pi$ with $\operatorname*{Comp}\pi=L$ will have the same
$\operatorname{st}$-equivalence class $\left[  \pi\right]
_{\operatorname{st}}$ (since $\operatorname{st}$ is a descent statistic).

Thus, each permutation $\pi$ satisfies
\begin{equation}
p\left(  F_{\operatorname*{Comp}\pi}\right)  =\left[  \pi\right]
_{\operatorname{st}} \label{pf.prop.Ast.alg.pFCp}%
\end{equation}
and therefore $\left[  \pi\right]  _{\operatorname{st}}=p\left(
F_{\operatorname*{Comp}\pi}\right)  \in p\left(  \operatorname*{QSym}\right)
$. Hence, $\mathcal{A}_{\operatorname{st}}\subseteq p\left(
\operatorname*{QSym}\right)  $ (since the $\operatorname{st}$-equivalence
classes $\left[  \pi\right]  _{\operatorname{st}}$ form a basis of
$\mathcal{A}_{\operatorname{st}}$). Consequently, the map $p$ is surjective.

Moreover, we have%
\begin{equation}
p\left(  ab\right)  =p\left(  a\right)  p\left(  b\right)
\ \ \ \ \ \ \ \ \ \ \text{for all }a,b\in\operatorname*{QSym}.
\label{pf.prop.Ast.alg.pab}%
\end{equation}

[\textit{Proof of (\ref{pf.prop.Ast.alg.pab}):} Let $a,b\in
\operatorname*{QSym}$. We must prove the equality (\ref{pf.prop.Ast.alg.pab}).
Since this equality is $\mathbb{Q}$-linear in each of $a$ and $b$, we WLOG
assume that $a$ and $b$ belong to the fundamental basis of
$\operatorname*{QSym}$. That is, $a=F_{J}$ and $b=F_{K}$ for two compositions
$J$ and $K$. Consider these $J$ and $K$. Fix any two disjoint permutations
$\pi$ and $\sigma$ such that $\operatorname*{Comp}\pi=J$ and
$\operatorname*{Comp}\sigma=K$. (Such $\pi$ and $\sigma$ are easy to find.)
The definition of $p$ thus yields $p\left(  F_{J}\right)  =\left[  \pi\right]
_{\operatorname{st}}$ and $p\left(  F_{K}\right)  =\left[  \sigma\right]
_{\operatorname{st}}$. Hence,%
\begin{align}
p\left(  \underbrace{a}_{=F_{J}}\right)  p\left(  \underbrace{b}_{=F_{K}%
}\right)   &  =\underbrace{p\left(  F_{J}\right)  }_{=\left[  \pi\right]
_{\operatorname{st}}}\underbrace{p\left(  F_{K}\right)  }_{=\left[
\sigma\right]  _{\operatorname{st}}}
=\left[  \pi\right]  _{\operatorname{st}}
\left[  \sigma\right]  _{\operatorname{st}}\nonumber\\
&  =\sum_{\tau\in S\left(  \pi,\sigma\right)  }\left[  \tau\right]
_{\operatorname{st}}\ \ \ \ \ \ \ \ \ \ \left(  \text{by
(\ref{pf.prop.Ast.alg.algdef})}\right) \nonumber\\
&  =\sum_{\chi\in S\left(  \pi,\sigma\right)  }\left[  \chi\right]
_{\operatorname{st}} \label{pf.prop.Ast.alg.pab.pf.1}%
\end{align}
(here, we have renamed the summation index $\tau$ as $\chi$). On the other
hand, $a=F_{J}=F_{\operatorname*{Comp}\pi}$ (since $J=\operatorname*{Comp}\pi
$) and $b=F_{\operatorname*{Comp}\sigma}$ (similarly); multiplying these
equalities, we get%
\[
ab=F_{\operatorname*{Comp}\pi}F_{\operatorname*{Comp}\sigma}=\sum_{\chi\in
S\left(  \pi,\sigma\right)  }F_{\operatorname*{Comp}\chi}%
\]
(by Proposition \ref{prop.4.1.rewr}). Applying the map $p$ to this equality,
we find%
\begin{align*}
p\left(  ab\right)   &  =p\left(  \sum_{\chi\in S\left(  \pi,\sigma\right)
}F_{\operatorname*{Comp}\chi}\right)  =\sum_{\chi\in S\left(  \pi
,\sigma\right)  }\underbrace{p\left(  F_{\operatorname*{Comp}\chi}\right)
}_{\substack{=\left[  \chi\right]  _{\operatorname{st}}\\\text{(by
(\ref{pf.prop.Ast.alg.pFCp}))}}}=\sum_{\chi\in S\left(  \pi,\sigma\right)
}\left[  \chi\right]  _{\operatorname{st}}\\
&  =p\left(  a\right)  p\left(  b\right)  \ \ \ \ \ \ \ \ \ \ \left(  \text{by
(\ref{pf.prop.Ast.alg.pab.pf.1})}\right)  .
\end{align*}
This proves (\ref{pf.prop.Ast.alg.pab}).]

The equality (\ref{pf.prop.Ast.alg.pab}) shows that $p$ is a magmatic algebra
homomorphism (since $p$ is $\mathbb{Q}$-linear). Thus, using the surjectivity
of $p$, we can easily see that the magmatic algebra $\mathcal{A}%
_{\operatorname{st}}$ is associative\footnote{\textit{Proof.} Let
$u,v,w\in\mathcal{A}_{\operatorname{st}}$. We must show that $\left(
uv\right)  w=u\left(  vw\right)  $.
\par
There exist $a,b,c\in\operatorname*{QSym}$ such that $u=p\left(  a\right)  $,
$v=p\left(  b\right)  $ and $w=p\left(  c\right)  $ (since $p$ is surjective).
Fix such $a,b,c$. Since $\operatorname*{QSym}$ is an actual (i.e., associative
unital) algebra, we have%
\begin{align*}
p\left(  abc\right)   &  =p\left(  \left(  ab\right)  c\right)
=\underbrace{p\left(  ab\right)  }_{\substack{=p\left(  a\right)  p\left(
b\right)  \\\text{(since }p\text{ is a magmatic}\\\text{algebra homomorphism)}%
}}p\left(  c\right)  \ \ \ \ \ \ \ \ \ \ \left(  \text{since }p\text{ is a
magmatic algebra homomorphism}\right) \\
&  =\left(  \underbrace{p\left(  a\right)  }_{=u}\underbrace{p\left(
b\right)  }_{=v}\right)  \underbrace{p\left(  c\right)  }_{=w}=\left(
uv\right)  w.
\end{align*}
A similar argument shows that $p\left(  abc\right)  =u\left(  vw\right)  $.
Thus, $\left(  uv\right)  w=p\left(  abc\right)  =u\left(  vw\right)  $,
qed.}. In other words, the multiplication on $\mathcal{A}_{\operatorname{st}}$
defined in Definition \ref{def.Ast} is associative. Moreover, it is clear
that the $\operatorname{st}$-equivalence class of the $0$-permutation
$\left( \right)$
serves as a neutral element for this multiplication (because if $\varnothing$
denotes the $0$-permutation $\left( \right)$,
then $S\left(  \varnothing,\sigma\right)
=S\left(  \sigma,\varnothing\right)  =\left\{  \sigma\right\}  $ for every
permutation $\sigma$). Thus, the multiplication on $\mathcal{A}%
_{\operatorname{st}}$ defined in Definition \ref{def.Ast} is well-defined and
associative, and turns $\mathcal{A}_{\operatorname{st}}$ into a $\mathbb{Q}%
$-algebra whose unity is the $\operatorname{st}$-equivalence class of the
$0$-permutation $\left( \right)$.
This proves Proposition \ref{prop.Ast.alg} \textbf{(a)}.

\textbf{(b)} The map $p:\operatorname*{QSym}\rightarrow\mathcal{A}%
_{\operatorname{st}}$ is $\mathbb{Q}$-linear and respects multiplication (by
(\ref{pf.prop.Ast.alg.pab})). Moreover, it sends the unity of
$\operatorname*{QSym}$ to the unity of the algebra $\mathcal{A}%
_{\operatorname{st}}$\ \ \ \ \footnote{\textit{Proof.} The unity of
$\operatorname*{QSym}$ is $1=F_{\left(  {}\right)  }$, where $\left(
{}\right)  $ denotes the empty composition. Now, let $\varnothing$ denote the
$0$-permutation $\left( \right)$.
Then, the $\operatorname{st}$-equivalence class $\left[
\varnothing\right]  _{\operatorname{st}}$ is the unity of the algebra
$\mathcal{A}_{\operatorname{st}}$. But the $0$-permutation $\varnothing
= \left( \right)$
has descent composition $\operatorname*{Comp}\varnothing=\left(  {}\right)  $.
Hence, the definition of $p$ yields $p\left(  F_{\left(  {}\right)  }\right)
=\left[  \varnothing\right]  _{\operatorname{st}}$. In view of what we just
said, this equality says that $p$ sends the unity of $\operatorname*{QSym}$ to
the unity of the algebra $\mathcal{A}_{\operatorname{st}}$.}. Thus, $p$ is a
$\mathbb{Q}$-algebra homomorphism. Moreover, recall that $p$ is surjective and
satisfies $p\left(  F_{\operatorname*{Comp}\pi}\right)  =\left[  \pi\right]
_{\operatorname{st}}$ for every permutation $\pi$. Hence, there is a
surjective $\mathbb{Q}$-algebra homomorphism
$p_{\operatorname{st}}:\operatorname*{QSym}
\rightarrow\mathcal{A}_{\operatorname{st}}$ that
satisfies
\[
p_{\operatorname{st}}\left(  F_{\operatorname*{Comp}\pi}\right)  =\left[
\pi\right]  _{\operatorname{st}}\ \ \ \ \ \ \ \ \ \ \text{for every
permutation }\pi
\]
(namely, $p_{\operatorname{st}}=p$). This proves Proposition
\ref{prop.Ast.alg} \textbf{(b)}.
\end{proof}

\begin{proof}
[Proof of Theorem \ref{thm.4.3}.]\textbf{(a)} $\Longrightarrow:$ Assume that
$\operatorname{st}$ is shuffle-compatible. Proposition \ref{prop.Ast.alg}
\textbf{(b)} shows that there is a surjective $\mathbb{Q}$-algebra
homomorphism $p_{\operatorname{st}}:\operatorname*{QSym}\rightarrow
\mathcal{A}_{\operatorname{st}}$ that satisfies
\begin{equation}
p_{\operatorname{st}}\left(  F_{\operatorname*{Comp}\pi}\right)  =\left[
\pi\right]  _{\operatorname{st}}\ \ \ \ \ \ \ \ \ \ \text{for every
permutation }\pi. \label{pf.thm.4.3.a.fwd.pst}%
\end{equation}
Consider this $p_{\operatorname{st}}$.

If $\alpha$ is an $\operatorname{st}$-equivalence class of compositions, then
we let $u_{\alpha}$ denote the $\operatorname{st}$-equivalence class $\left[
\pi\right]  _{\operatorname{st}}$ of all permutations $\pi$ whose descent
composition $\operatorname*{Comp}\pi$ belongs to $\alpha$. (This is indeed a
well-defined $\operatorname{st}$-equivalence class, because
$\operatorname{st}$ is a descent statistic.) This establishes a bijection
between the $\operatorname{st}$-equivalence classes of compositions and the
$\operatorname{st}$-equivalence classes of permutations. Thus, the family
$\left(  u_{\alpha}\right)  $ (indexed by $\operatorname{st}$-equivalence
classes $\alpha$ of compositions) is just a reindexing of the basis of
$\mathcal{A}_{\operatorname{st}}$ consisting of the
$\operatorname{st}$-equivalence classes
$\left[  \pi\right]  _{\operatorname{st}}$ of
permutations. Consequently, this family is a basis of the $\mathbb{Q}$-vector
space $\mathcal{A}_{\operatorname{st}}$. Moreover, $p_{\operatorname{st}}$
is a $\mathbb{Q}$-algebra homomorphism $\operatorname*{QSym}\rightarrow
\mathcal{A}_{\operatorname{st}}$ with the property that whenever $\alpha$ is
an $\operatorname{st}$-equivalence class of compositions, we have%
\[
p_{\operatorname{st}}\left(  F_{L}\right)  =u_{\alpha}%
\ \ \ \ \ \ \ \ \ \ \text{for each }L\in\alpha.
\]
(Indeed, this follows from applying (\ref{pf.thm.4.3.a.fwd.pst}) to any
permutation $\pi$ satisfying $\operatorname*{Comp}\pi=L$.)

Thus, there exist a $\mathbb{Q}$-algebra $A$ (namely, $\mathcal{A}%
=\mathcal{A}_{\operatorname{st}}$) with basis $\left(  u_{\alpha}\right)  $
(indexed by $\operatorname{st}$-equivalence classes $\alpha$ of compositions)
and a $\mathbb{Q}$-algebra homomorphism $\phi_{\operatorname{st}}
:\operatorname*{QSym}\rightarrow A$ (namely, $\phi_{\operatorname{st}}
=p_{\operatorname{st}}$) with the property that whenever $\alpha$ is an
$\operatorname{st}$-equivalence class of compositions, we have%
\[
\phi_{\operatorname{st}}\left(  F_{L}\right)  =u_{\alpha}%
\ \ \ \ \ \ \ \ \ \ \text{for each }L\in\alpha.
\]
This proves the $\Longrightarrow$ direction of Theorem \ref{thm.4.3}
\textbf{(a)}.

$\Longleftarrow:$ Assume that there exist a $\mathbb{Q}$-algebra $A$ with
basis $\left(  u_{\alpha}\right)  $ (indexed by $\operatorname{st}%
$-equivalence classes $\alpha$ of compositions) and a $\mathbb{Q}$-algebra
homomorphism $\phi_{\operatorname{st}}:\operatorname*{QSym}\rightarrow A$
with the property that whenever $\alpha$ is an $\operatorname{st}%
$-equivalence class of compositions, we have%
\[
\phi_{\operatorname{st}}\left(  F_{L}\right)  =u_{\alpha}%
\ \ \ \ \ \ \ \ \ \ \text{for each }L\in\alpha.
\]
Consider this $A$, this $\left(  u_{\alpha}\right)  $ and this $\phi
_{\operatorname{st}}$. Lemma \ref{lem.K.fist} shows that $\operatorname*{Ker}%
\left(  \phi_{\operatorname{st}}\right)  =\mathcal{K}_{\operatorname{st}}$.
But $\operatorname*{Ker}\left(  \phi_{\operatorname{st}}\right)  $ is an
ideal of $\operatorname*{QSym}$ (since $\phi_{\operatorname{st}}$ is a
$\mathbb{Q}$-algebra homomorphism). In other words, $\mathcal{K}%
_{\operatorname{st}}$ is an ideal of $\operatorname*{QSym}$ (since
$\operatorname*{Ker}\left(  \phi_{\operatorname{st}}\right)  =\mathcal{K}%
_{\operatorname{st}}$).

Now, consider any two disjoint permutations $\pi$ and $\sigma$. Also, consider
two further disjoint permutations $\pi^{\prime}$ and $\sigma^{\prime}$
satisfying $\operatorname{st} \pi   =\operatorname{st}\left(
\pi^{\prime}\right)  $, $\operatorname{st} \sigma 
=\operatorname{st}\left(  \sigma^{\prime}\right)  $, $\left\vert
\pi\right\vert =\left\vert \pi^{\prime}\right\vert $ and $\left\vert
\sigma\right\vert =\left\vert \sigma^{\prime}\right\vert $. We shall show that
$\left\{  \operatorname{st} \tau   \ \mid\ \tau\in S\left(
\pi,\sigma\right)  \right\}  _{\operatorname*{multi}} =\left\{
\operatorname{st} \tau   \ \mid\ \tau\in S\left(  \pi^{\prime
},\sigma^{\prime}\right)  \right\}  _{\operatorname*{multi}} $ as multisets.
This will show that the multiset \newline$\left\{  \operatorname{st}
\tau \ \mid\ \tau\in S\left(  \pi,\sigma\right)  \right\}
_{\operatorname*{multi}} $ depends only on $\operatorname{st} \pi$,
$\operatorname{st} \sigma$, $\left\vert \pi\right\vert $ and
$\left\vert \sigma\right\vert $.

From $\operatorname{st} \pi   =\operatorname{st}\left(
\pi^{\prime}\right)  $ and $\left\vert \pi\right\vert =\left\vert \pi^{\prime
}\right\vert $, we conclude that $\pi$ and $\pi^{\prime}$ are
$\operatorname{st}$-equivalent. In other words, $\operatorname*{Comp}\pi$ and
$\operatorname*{Comp}\left(  \pi^{\prime}\right)  $ are $\operatorname{st}%
$-equivalent. Hence, $F_{\operatorname*{Comp}\pi}-F_{\operatorname*{Comp}%
\left(  \pi^{\prime}\right)  }\in\mathcal{K}_{\operatorname{st}}$ (by the
definition of $\mathcal{K}_{\operatorname{st}}$), so that
$F_{\operatorname*{Comp}\pi}\equiv F_{\operatorname*{Comp}\left(  \pi^{\prime
}\right)  }\operatorname{mod}\mathcal{K}_{\operatorname{st}}$. Similarly,
$F_{\operatorname*{Comp}\sigma}\equiv F_{\operatorname*{Comp}\left(
\sigma^{\prime}\right)  }\operatorname{mod}\mathcal{K}_{\operatorname{st}}$.
These two congruences, combined, yield $F_{\operatorname*{Comp}\pi
}F_{\operatorname*{Comp}\sigma}\equiv F_{\operatorname*{Comp}\left(
\pi^{\prime}\right)  }F_{\operatorname*{Comp}\left(  \sigma^{\prime}\right)
}\operatorname{mod}\mathcal{K}_{\operatorname{st}}$, because $\mathcal{K}%
_{\operatorname{st}}$ is an ideal of $\operatorname*{QSym}$.

Let $X$ be the codomain of the map $\operatorname{st}$. Let $\mathbb{Q}%
\left[  X\right]  $ be the free $\mathbb{Q}$-vector space with basis $\left(
\left[  x\right]  \right)  _{x\in X}$. Then, we can define a $\mathbb{Q}%
$-linear map $\mathbf{st}:\operatorname*{QSym}\rightarrow\mathbb{Q}\left[
X\right]  ,\ F_{J}\mapsto\left[  \operatorname{st}J\right]  $. This map
$\mathbf{st}$ sends each of the generators of $\mathcal{K}_{\operatorname{st}%
}$ to $0$ (by the definition of $\mathcal{K}_{\operatorname{st}}$), and
therefore sends the whole $\mathcal{K}_{\operatorname{st}}$ to $0$. In other
words, $\mathbf{st}\left(  \mathcal{K}_{\operatorname{st}}\right)  =0$.

We have $F_{\operatorname*{Comp}\pi}F_{\operatorname*{Comp}\sigma}\equiv
F_{\operatorname*{Comp}\left(  \pi^{\prime}\right)  }F_{\operatorname*{Comp}%
\left(  \sigma^{\prime}\right)  }\operatorname{mod}\mathcal{K}%
_{\operatorname{st}}$ and thus%
\begin{equation}
\mathbf{st}\left(  F_{\operatorname*{Comp}\pi}F_{\operatorname*{Comp}\sigma
}\right)  =\mathbf{st}\left(  F_{\operatorname*{Comp}\left(  \pi^{\prime
}\right)  }F_{\operatorname*{Comp}\left(  \sigma^{\prime}\right)  }\right)
\label{pf.thm.4.3.a.fwd.steq}%
\end{equation}
(since $\mathbf{st}\left(  \mathcal{K}_{\operatorname{st}}\right)  =0$). But
Proposition \ref{prop.4.1.rewr} yields%
\[
F_{\operatorname*{Comp}\pi}F_{\operatorname*{Comp}\sigma}=\sum_{\chi\in
S\left(  \pi,\sigma\right)  }F_{\operatorname*{Comp}\chi}.
\]
Applying the map $\mathbf{st}$ to both sides of this equality, we find%
\begin{align*}
\mathbf{st}\left(  F_{\operatorname*{Comp}\pi}F_{\operatorname*{Comp}\sigma
}\right)   &  =\mathbf{st}\left(  \sum_{\chi\in S\left(  \pi,\sigma\right)
}F_{\operatorname*{Comp}\chi}\right) \\
&  =\sum_{\chi\in S\left(  \pi,\sigma\right)  }\underbrace{\mathbf{st}\left(
F_{\operatorname*{Comp}\chi}\right)  }_{=\left[  \operatorname{st}\left(
\operatorname*{Comp}\chi\right)  \right]  =\left[  \operatorname{st}%
\chi\right]  }=\sum_{\chi\in S\left(  \pi,\sigma\right)  }\left[
\operatorname{st}\chi\right]  .
\end{align*}
Similarly,%
\[
\mathbf{st}\left(  F_{\operatorname*{Comp}\left(  \pi^{\prime}\right)
}F_{\operatorname*{Comp}\left(  \sigma^{\prime}\right)  }\right)  =\sum
_{\chi\in S\left(  \pi^{\prime},\sigma^{\prime}\right)  }\left[
\operatorname{st}\chi\right]  .
\]
But the left-hand sides of the last two equalities are equal (because of
(\ref{pf.thm.4.3.a.fwd.steq})); therefore, the right-hand sides must be equal
as well. In other words,
\[
\sum_{\chi\in S\left(  \pi,\sigma\right)  }\left[  \operatorname{st}%
\chi\right]  =\sum_{\chi\in S\left(  \pi^{\prime},\sigma^{\prime}\right)
}\left[  \operatorname{st}\chi\right]  .
\]
This shows exactly that $\left\{  \operatorname{st} \chi
\ \mid\ \chi\in S\left(  \pi,\sigma\right)  \right\}  _{\operatorname*{multi}}
=\left\{  \operatorname{st} \chi \ \mid\ \chi\in S\left(
\pi^{\prime},\sigma^{\prime}\right)  \right\}  _{\operatorname*{multi}} $. In
other words, $\left\{  \operatorname{st} \tau   \ \mid\ \tau\in
S\left(  \pi,\sigma\right)  \right\}  _{\operatorname*{multi}} =\left\{
\operatorname{st} \tau   \ \mid\ \tau\in S\left(  \pi^{\prime
},\sigma^{\prime}\right)  \right\}  _{\operatorname*{multi}} $. Thus, we have
proven that the multiset $\left\{  \operatorname{st} \tau 
\ \mid\ \tau\in S\left(  \pi,\sigma\right)  \right\}  _{\operatorname*{multi}}
$ depends only on $\operatorname{st} \pi   $,
$\operatorname{st} \sigma   $, $\left\vert \pi\right\vert $ and
$\left\vert \sigma\right\vert $. Hence, the statistic $\operatorname{st}$ is
shuffle-compatible. This proves the $\Longleftarrow$ direction of Theorem
\ref{thm.4.3} \textbf{(a)}.

\textbf{(b)} Proposition \ref{prop.Ast.alg} \textbf{(b)} shows that there is a
surjective $\mathbb{Q}$-algebra homomorphism $p_{\operatorname{st}%
}:\operatorname*{QSym}\rightarrow\mathcal{A}_{\operatorname{st}}$ that
satisfies
\begin{equation}
p_{\operatorname{st}}\left(  F_{\operatorname*{Comp}\pi}\right)  =\left[
\pi\right]  _{\operatorname{st}}\ \ \ \ \ \ \ \ \ \ \text{for every
permutation }\pi. \label{pf.thm.4.3.b.pst}%
\end{equation}
Consider this $p_{\operatorname{st}}$.

Let $\gamma$ be the $\mathbb{Q}$-linear map%
\[
\mathcal{A}_{\operatorname{st}}\rightarrow A,\ \ \ \ \ \ \ \ \ \ \left[
\pi\right]  _{\operatorname{st}}\mapsto u_{\alpha},
\]
where $\alpha$ is the $\operatorname{st}$-equivalence class of the
composition $\operatorname*{Comp}\pi$. This map $\gamma$ is clearly
well-defined (since the $\operatorname{st}$-equivalence classes $\left[
\pi\right]  _{\operatorname{st}}$ form a basis of $\mathcal{A}%
_{\operatorname{st}}$, and since the $\operatorname{st}$-equivalence class
of the composition $\operatorname*{Comp}\pi$ depends only on the
$\operatorname{st}$-equivalence class $\left[  \pi\right]
_{\operatorname{st}}$ and not on the permutation $\pi$ itself). Moreover,
$\gamma$ sends a basis of $\mathcal{A}_{\operatorname{st}}$ (the basis formed
by the $\operatorname{st}$-equivalence classes $\left[  \pi\right]
_{\operatorname{st}}$ of permutations) to a basis of $A$ (namely, to the
basis $\left(  u_{\alpha}\right)  $) bijectively; thus, $\gamma$ is an
isomorphism of $\mathbb{Q}$-vector spaces.

The diagram%
\[%
%TCIMACRO{\TeXButton{comm triangle}{\xymatrix{
%\QSym\arsurj[r]^{p_{\operatorname{st}}} \ar[dr]_{\phi_{\operatorname{st}}}
%& \mathcal{A}_{\operatorname{st}} \arsurj[d]^\gamma\\
%& A
%}}}%
%BeginExpansion
\xymatrix{
\QSym\arsurj[r]^{p_{\operatorname{st}}} \ar[dr]_{\phi_{\operatorname{st}}}
& \mathcal{A}_{\operatorname{st}} \arsurj[d]^\gamma\\
& A
}%
%EndExpansion
\]
is commutative (as one can easily check by tracing an arbitrary basis element
$F_{L}$ of $\operatorname*{QSym}$ through the diagram). Since the maps
$p_{\operatorname{st}}$ and $\phi_{\operatorname{st}}$ in this diagram are
$\mathbb{Q}$-algebra homomorphisms, and since $p_{\operatorname{st}}$ is
surjective, we thus conclude that $\gamma$ is also a $\mathbb{Q}$-algebra
homomorphism\footnote{\textit{Proof.} Let $a,b\in\mathcal{A}%
_{\operatorname{st}}$. We shall show that $\gamma\left(  ab\right)
=\gamma\left(  a\right)  \gamma\left(  b\right)  $.
\par
There exist $a^{\prime},b^{\prime}\in\operatorname*{QSym}$ such that
$a=p_{\operatorname{st}}\left(  a^{\prime}\right)  $ and
$b=p_{\operatorname{st}}\left(  b^{\prime}\right)  $ (since
$p_{\operatorname{st}}$ is surjective). Consider these $a^{\prime},b^{\prime
}$. Then, $\gamma\left(  \underbrace{a}_{=p_{\operatorname{st}}\left(
a^{\prime}\right)  }\right)  =\gamma\left(  p_{\operatorname{st}}\left(
a^{\prime}\right)  \right)  =\phi_{\operatorname{st}}\left(  a^{\prime
}\right)  $ (since the diagram is commutative) and $\gamma\left(  b\right)
=\phi_{\operatorname{st}}\left(  b^{\prime}\right)  $ (similarly). But from
$a=p_{\operatorname{st}}\left(  a^{\prime}\right)  $ and
$b=p_{\operatorname{st}}\left(  b^{\prime}\right)  $, we obtain
$ab=p_{\operatorname{st}}\left(  a^{\prime}\right)  p_{\operatorname{st}%
}\left(  b^{\prime}\right)  =p_{\operatorname{st}}\left(  a^{\prime}%
b^{\prime}\right)  $ (since $p_{\operatorname{st}}$ is a $\mathbb{Q}$-algebra
homomorphism), so that%
\begin{align*}
\gamma\left(  ab\right)   &  =\gamma\left(  p_{\operatorname{st}}\left(
a^{\prime}b^{\prime}\right)  \right)  =\phi_{\operatorname{st}}\left(
a^{\prime}b^{\prime}\right)  \ \ \ \ \ \ \ \ \ \ \left(  \text{since the
diagram is commutative}\right) \\
&  =\underbrace{\phi_{\operatorname{st}}\left(  a^{\prime}\right)  }%
_{=\gamma\left(  a\right)  }\underbrace{\phi_{\operatorname{st}}\left(
b^{\prime}\right)  }_{=\gamma\left(  b\right)  }\ \ \ \ \ \ \ \ \ \ \left(
\text{since }\phi_{\operatorname{st}}\text{ is a }\mathbb{Q}\text{-algebra
homomorphism}\right) \\
&  =\gamma\left(  a\right)  \gamma\left(  b\right)  .
\end{align*}
\par
Now, forget that we fixed $a,b$. We thus have proven that $\gamma\left(
ab\right)  =\gamma\left(  a\right)  \gamma\left(  b\right)  $ for all
$a,b\in\mathcal{A}_{\operatorname{st}}$. Similarly, $\gamma\left(  1\right)
=1$. Hence, $\gamma$ is a $\mathbb{Q}$-algebra homomorphism (since $\gamma$ is
$\mathbb{Q}$-linear).}. Since $\gamma$ is an isomorphism of $\mathbb{Q}%
$-vector spaces, we thus conclude that $\gamma$ is a $\mathbb{Q}$-algebra
isomorphism $\mathcal{A}_{\operatorname{st}}\rightarrow A$. This proves
Theorem \ref{thm.4.3} \textbf{(b)}.
\end{proof}
\end{verlong}

\section{\label{sect.dendri}Dendriform structures}

\begin{vershort}
Next, we shall recall the \textit{dendriform operations }$\left.
\prec\right.  $ and $\left.  \succeq\right.  $ on $\operatorname*{QSym}$
studied in \cite{dimcr}, and we shall connect these operations back to
LR-shuffle-compatibility. Since we consider this somewhat tangential to the
present paper, we merely summarize the main results here; more can be found in
\cite{verlong}.

\subsection{Two operations on $\operatorname*{QSym}$}

We begin with some definitions. We will use some notations from \cite{dimcr},
but we set $\mathbf{k}=\mathbb{Q}$ because we are working over the ring
$\mathbb{Q}$ in this paper. Monomials always mean formal expressions of the
form $x_{1}^{a_{1}}x_{2}^{a_{2}}x_{3}^{a_{3}}\cdots$ with $a_{1}+a_{2}%
+a_{3}+\cdots<\infty$ (see \cite[Section 2]{dimcr} for details). If
$\mathfrak{m}$ is a monomial, then $\operatorname*{Supp}\mathfrak{m}$ will
denote the finite subset
\[
\left\{  i\in\left\{  1,2,3,\ldots\right\}  \ \mid\ \text{the exponent with
which }x_{i}\text{ occurs in }\mathfrak{m}\text{ is }>0\right\}
\]
of $\left\{  1,2,3,\ldots\right\}  $. Next, we define two binary operations
\begin{align*}
&  \left.  \prec\right.  \ \left(  \text{called \textquotedblleft dendriform
less-than\textquotedblright; but it's an operation, not a relation}\right)
,\\
&  \left.  \succeq\right.  \ \left(  \text{called \textquotedblleft dendriform
greater-or-equal\textquotedblright; but it's an operation, not a
relation}\right)  ,
\end{align*}
on the ring $\mathbf{k}\left[  \left[  x_{1},x_{2},x_{3},\ldots\right]
\right]  $ of power series by first defining how they act on monomials:%
\begin{align*}
\mathfrak{m}\left.  \prec\right.  \mathfrak{n}  &  =\left\{
\begin{array}
[c]{l}%
\mathfrak{m}\cdot\mathfrak{n},\ \ \ \ \ \ \ \ \ \ \text{if }\min\left(
\operatorname*{Supp}\mathfrak{m}\right)  <\min\left(  \operatorname*{Supp}%
\mathfrak{n}\right)  ;\\
0,\ \ \ \ \ \ \ \ \ \ \text{if }\min\left(  \operatorname*{Supp}%
\mathfrak{m}\right)  \geq\min\left(  \operatorname*{Supp}\mathfrak{n}\right)
\end{array}
\right.  ;\\
\mathfrak{m}\left.  \succeq\right.  \mathfrak{n}  &  =\left\{
\begin{array}
[c]{l}%
\mathfrak{m}\cdot\mathfrak{n},\ \ \ \ \ \ \ \ \ \ \text{if }\min\left(
\operatorname*{Supp}\mathfrak{m}\right)  \geq\min\left(  \operatorname*{Supp}%
\mathfrak{n}\right)  ;\\
0,\ \ \ \ \ \ \ \ \ \ \text{if }\min\left(  \operatorname*{Supp}%
\mathfrak{m}\right)  <\min\left(  \operatorname*{Supp}\mathfrak{n}\right)
\end{array}
\right.  ;
\end{align*}
and then requiring that they both be $\mathbf{k}$-bilinear and continuous (so
their action on pairs of arbitrary power series can be computed by
\textquotedblleft opening the parentheses\textquotedblright). These operations
$\left.  \prec\right.  $ and $\left.  \succeq\right.  $ restrict to the subset
$\operatorname*{QSym}$ of $\mathbf{k}\left[  \left[  x_{1},x_{2},x_{3}%
,\ldots\right]  \right]  $ (this is proven in \cite[detailed version, Section
3]{dimcr}). They furthermore satisfy the following relations (which are easy
to verify):

\begin{itemize}
\item For all $a,b,c\in\mathbf{k}\left[  \left[  x_{1},x_{2},x_{3},\ldots\right]
\right]  $, we have
\begin{align*}
a\left.  \prec\right.  b+a\left.  \succeq\right.  b  &  =ab;\\
\left(  a\left.  \prec\right.  b\right)  \left.  \prec\right.  c  &  =a\left.
\prec\right.  \left(  bc\right)  ;\ \ \ \ \ \ \ \ \ \ \left(  a\left.
\succeq\right.  b\right)  \left.  \prec\right.  c=a\left.  \succeq\right.
\left(  b\left.  \prec\right.  c\right)  ;\\
a\left.  \succeq\right.  \left(  b\left.  \succeq\right.  c\right)   &
=\left(  ab\right)  \left.  \succeq\right.  c .
\end{align*}

\item For any $a\in\mathbf{k}\left[  \left[  x_{1},x_{2},x_{3},\ldots\right]
\right]  $, we have%
\begin{align*}
1\left.  \prec\right.  a  &  =0;\ \ \ \ \ \ \ \ \ \ a\left.  \prec\right.
1=a-\varepsilon\left(  a\right)  ; \ \ \ \ \ \ \ \ \ \ % \\
1\left.  \succeq\right.  a   =a;\ \ \ \ \ \ \ \ \ \ a\left.  \succeq\right.
1=\varepsilon\left(  a\right)  ,
\end{align*}
where $\varepsilon\left(  a\right)  $ denotes the constant term of the power
series $a$.
\end{itemize}

The operations $\left.  \prec\right.  $ and $\left.  \succeq\right.  $ are
sometimes called \textquotedblleft restricted products\textquotedblright\ due
to their similarity with the (regular) multiplication of $\operatorname*{QSym}%
$. In particular, they satisfy the following analogue of Proposition
\ref{prop.4.1.rewr}:

\begin{proposition}
Let $\pi$ and $\sigma$ be two disjoint nonempty permutations. Assume that
$\pi_{1}>\sigma_{1}$. Then,%
\begin{align*}
F_{\operatorname*{Comp}\pi}\left.  \prec\right.  F_{\operatorname*{Comp}%
\sigma}
&=\sum_{\chi\in S_{\prec}\left(  \pi,\sigma\right)  }%
F_{\operatorname*{Comp}\chi}  \qquad \text{ and } \\
F_{\operatorname*{Comp}\pi}\left.  \succeq\right.  F_{\operatorname*{Comp}%
\sigma}
&=\sum_{\chi\in S_{\succ}\left(  \pi,\sigma\right)  }%
F_{\operatorname*{Comp}\chi}.
\end{align*}

\end{proposition}

\subsection{Left- and right-shuffle-compatibility and ideals}

This proposition lets us relate the notions introduced in Definition
\ref{def.LR.left-right} to the operations $\left.  \prec\right.  $ and
$\left.  \succeq\right.  $. To state the precise connection, we need the
following notation:

\begin{definition}
Let $A$ be a $\mathbf{k}$-module equipped with some binary operation $\ast$
(written infix).

\textbf{(a)} If $B$ and $C$ are two $\mathbf{k}$-submodules of $A$, then
$B\ast C$ shall mean the $\mathbf{k}$-submodule of $A$ spanned by all elements
of the form $b\ast c$ with $b\in B$ and $c\in C$.

\textbf{(b)} A $\mathbf{k}$-submodule $M$ of $A$ is said to be a\textit{
}$\ast$\textit{-ideal} if and only if it satisfies $A\ast M\subseteq M$ and
$M\ast A\subseteq M$.
\end{definition}

Now, let us define two further variants of LR-shuffle-compatibility (to be
compared with those introduced in Definition \ref{def.LR.left-right}):

\begin{definition}
Let $\operatorname{st}$ be a permutation statistic.

\textbf{(a)} We say that $\operatorname{st}$ is \textit{weakly
left-shuffle-compatible} if for any two disjoint nonempty permutations $\pi$
and $\sigma$ having the property that%
\[
\text{each entry of }\pi\text{ is greater than each entry of }\sigma,
\]
the multiset $\left\{  \operatorname{st} \tau   \ \mid\ \tau\in
S_{\prec}\left(  \pi,\sigma\right)  \right\}  _{\operatorname*{multi}}$
depends only on $\operatorname{st} \pi$, $\operatorname{st} \sigma$,
$\left\vert \pi\right\vert $ and $\left\vert \sigma\right\vert $.

\textbf{(b)} We say that $\operatorname{st}$ is \textit{weakly
right-shuffle-compatible} if for any two disjoint nonempty permutations $\pi$
and $\sigma$ having the property that%
\[
\text{each entry of }\pi\text{ is greater than each entry of }\sigma,
\]
the multiset $\left\{  \operatorname{st} \tau   \ \mid\ \tau\in
S_{\succ}\left(  \pi,\sigma\right)  \right\}  _{\operatorname*{multi}}$
depends only on $\operatorname{st} \pi$, $\operatorname{st} \sigma$,
$\left\vert \pi\right\vert $ and $\left\vert \sigma\right\vert $.
\end{definition}

Then, the following analogues to the first part of Proposition
\ref{prop.K.ideal} hold:

\begin{theorem}
\label{thm.dendri.K.ideal}Let $\operatorname{st}$ be a descent statistic.
Then, the following three statements are equivalent:

\begin{itemize}
\item \textit{Statement A:} The statistic $\operatorname{st}$ is left-shuffle-compatible.

\item \textit{Statement B:} The statistic $\operatorname{st}$ is weakly left-shuffle-compatible.

\item \textit{Statement C:} The set $\mathcal{K}_{\operatorname{st}}$ is an
$\left.  \prec\right.  $-ideal of $\operatorname*{QSym}$.
\end{itemize}
\end{theorem}

\begin{theorem}
\label{thm.dendri.K.ideal-R}Let $\operatorname{st}$ be a descent statistic.
Then, the following three statements are equivalent:

\begin{itemize}
\item \textit{Statement A:} The statistic $\operatorname{st}$ is right-shuffle-compatible.

\item \textit{Statement B:} The statistic $\operatorname{st}$ is weakly right-shuffle-compatible.

\item \textit{Statement C:} The set $\mathcal{K}_{\operatorname{st}}$ is an
$\left.  \succeq\right.  $-ideal of $\operatorname*{QSym}$.
\end{itemize}
\end{theorem}

\begin{corollary}
\label{cor.dendri.K.ideal-LR}Let $\operatorname{st}$ be a permutation
statistic that is LR-shuffle-compatible. Then, $\operatorname{st}$ is a
shuffle-compatible descent statistic, and the set $\mathcal{K}%
_{\operatorname{st}}$ is an ideal and a $\left.  \prec\right.  $-ideal and a
$\left.  \succeq\right.  $-ideal of $\operatorname*{QSym}$.
\end{corollary}

\begin{corollary}
\label{cor.dendri.K.ideal-LRi}Let $\operatorname{st}$ be a descent statistic
such that $\mathcal{K}_{\operatorname{st}}$ is a $\left.  \prec\right.
$-ideal and a $\left.  \succeq\right.  $-ideal of $\operatorname*{QSym}$.
Then, $\operatorname{st}$ is LR-shuffle-compatible and shuffle-compatible.
\end{corollary}

Corollary \ref{cor.dendri.K.ideal-LR} can (for example) be applied to
$\operatorname{st}=\operatorname{Epk}$, which we know to be
LR-shuffle-compatible (from Theorem \ref{thm.LRcomp.Pks} \textbf{(c)}); the
result is that $\mathcal{K}_{\operatorname{Epk}}$ is an ideal and a $\left.
\prec\right.  $-ideal and a $\left.  \succeq\right.  $-ideal of
$\operatorname*{QSym}$. The same can be said about $\operatorname{Des}$ and
$\operatorname*{Lpk}$ and some other statistics.

Combining Theorem \ref{thm.dendri.K.ideal} with Theorem
\ref{thm.dendri.K.ideal-R}, we can also see that any descent statistic that is
weakly left-shuffle-compatible and weakly right-shuffle-compatible must
automatically be shuffle-compatible (because any $\left.  \prec\right.
$-ideal of $\operatorname*{QSym}$ that is also a $\left.  \succeq\right.
$-ideal of $\operatorname*{QSym}$ is an ideal of $\operatorname*{QSym}$ as
well). Note that this is only true for descent statistics! As far as arbitrary
permutation statistics are concerned, this is false; for example, the number
of inversions is weakly left-shuffle-compatible and weakly
right-shuffle-compatible but not shuffle-compatible.

Let us next define the notion of dendriform algebras:

\begin{definition}
\textbf{(a)} A \textit{dendriform algebra} over a field $\mathbf{k}$ means a
$\mathbf{k}$-algebra $A$ equipped with two further $\mathbf{k}$-bilinear
binary operations $\left.  \prec\right.  $ and $\left.  \succeq\right.  $
(these are operations, not relations, despite the symbols) from $A\times A$ to
$A$ that satisfy the four rules%
\begin{align*}
a\left.  \prec\right.  b+a\left.  \succeq\right.  b  &  =ab;\\
\left(  a\left.  \prec\right.  b\right)  \left.  \prec\right.  c  &  =a\left.
\prec\right.  \left(  bc\right)  ;\\
\left(  a\left.  \succeq\right.  b\right)  \left.  \prec\right.  c  &
=a\left.  \succeq\right.  \left(  b\left.  \prec\right.  c\right)  ;\\
a\left.  \succeq\right.  \left(  b\left.  \succeq\right.  c\right)   &
=\left(  ab\right)  \left.  \succeq\right.  c
\end{align*}
for all $a,b,c\in A$. (Depending on the situation, it is useful to also impose
a few axioms that relate the unity $1$ of the $\mathbf{k}$-algebra $A$ with
the operations $\left.  \prec\right.  $ and $\left.  \succeq\right.  $. For
example, we could require $1\left.  \prec\right.  a=0$ for each $a\in A$. For
what we are going to do, these extra axioms don't matter.)

\textbf{(b)} If $A$ and $B$ are two dendriform algebras over $\mathbf{k}$,
then a \textit{dendriform algebra homomorphism} from $A$ to $B$ means a
$\mathbf{k}$-algebra homomorphism $\phi:A\rightarrow B$ preserving the
operations $\left.  \prec\right.  $ and $\left.  \succeq\right.  $ (that is,
satisfying $\phi\left(  a\left.  \prec\right.  b\right)  =\phi\left(
a\right)  \left.  \prec\right.  \phi\left(  b\right)  $ and $\phi\left(
a\left.  \succeq\right.  b\right)  =\phi\left(  a\right)  \left.
\succeq\right.  \phi\left(  b\right)  $ for all $a,b\in A$). (Some authors
only require it to be a $\mathbf{k}$-linear map instead of being a
$\mathbf{k}$-algebra homomorphism; this boils down to the question whether
$\phi\left(  1\right)  $ must be $1$ or not. This does not make a difference
for us here.)
\end{definition}

Thus, $\operatorname*{QSym}$ (with its two operations $\left.  \prec\right.  $
and $\left.  \succeq\right.  $) becomes a dendriform algebra over $\mathbb{Q}$.

Notice that if $A$ and $B$ are two dendriform algebras over $\mathbf{k}$, then
the kernel of any dendriform algebra homomorphism $A\rightarrow B$ is an
$\left.  \prec\right.  $-ideal and a $\left.  \succeq\right.  $-ideal of $A$.
Conversely, if $A$ is a dendriform algebra over $\mathbf{k}$, and $I$ is
simultaneously a $\left.  \prec\right.  $-ideal and a $\left.  \succeq
\right.  $-ideal of $A$, then $A/I$ canonically becomes a dendriform algebra,
and the canonical projection $A\rightarrow A/I$ becomes a dendriform algebra homomorphism.

Therefore, Corollary \ref{cor.dendri.K.ideal-LR} (and the $\mathcal{A}%
_{\operatorname{st}}\cong\operatorname*{QSym}/\mathcal{K}_{\operatorname{st}%
}$ isomorphism from Proposition \ref{prop.K.ideal}) yields the following:

\begin{corollary}
If a descent statistic $\operatorname{st}$ is LR-shuffle-compatible, then its
shuffle algebra $\mathcal{A}_{\operatorname{st}}$ canonically becomes a
dendriform algebra.
\end{corollary}

We furthermore have the following analogue of Theorem~\ref{thm.4.3}, which
easily follows from Theorem \ref{thm.dendri.K.ideal} and Theorem
\ref{thm.dendri.K.ideal-R}:

\begin{theorem}
\label{thm.dendri.4.3}Let $\operatorname{st}$ be a descent statistic.

\textbf{(a)} The descent statistic $\operatorname{st}$ is
left-shuffle-compatible and right-shuffle-compatible if and only if there
exist a dendriform algebra $A$ with
basis $\left(  u_{\alpha}\right)  $ (indexed by $\operatorname{st}%
$-equivalence classes $\alpha$ of compositions)
and a dendriform algebra homomorphism
$\phi_{\operatorname{st}}:\operatorname*{QSym} \rightarrow A$
with the property that whenever $\alpha$ is an
$\operatorname{st}$-equivalence class of compositions, we have
\[
\phi_{\operatorname{st}}\left(  F_{L}\right)  =u_{\alpha}%
\ \ \ \ \ \ \ \ \ \ \text{for each }L\in\alpha.
\]

\textbf{(b)} In this case, the $\mathbb{Q}$-linear map%
\[
\mathcal{A}_{\operatorname{st}}\rightarrow A,\ \ \ \ \ \ \ \ \ \ \left[
\pi\right]  _{\operatorname{st}}\mapsto u_{\alpha},
\]
where $\alpha$ is the $\operatorname{st}$-equivalence class of the
composition $\operatorname*{Comp}\pi$, is an isomorphism of dendriform
algebras $\mathcal{A}_{\operatorname{st}}\rightarrow A$.
\end{theorem}

\begin{question}
Can the $\mathbb{Q}$-algebra $\operatorname*{Pow}\mathcal{N}$ from Definition
\ref{def.GammaZ} be endowed with two binary operations $\left.  \prec\right.
$ and $\left.  \succeq\right.  $ that make it into a dendriform algebra? Can
we then find an analogue of Proposition \ref{prop.prod1} along the following lines?

Let $\left(  P,\gamma\right)  $, $\left(  Q,\delta\right)  $ and $\left(
P\sqcup Q,\varepsilon\right)  $ be as in Proposition \ref{prop.prod1}. Assume
that each of the posets $P$ and $Q$ has a (global) minimum element; denote these
elements by $\min P$ and $\min Q$, respectively.
Let $P \left. \prec \right. Q$ be the poset obtained by adding
the relation $\min P<\min Q$ to $P\sqcup Q$.
Let $P \left. \succeq \right. Q$ be the poset obtained by adding
the relation $\min P>\min Q$ to $P\sqcup Q$.
Then, we hope to have%
\begin{align*}
\Gamma_{\mathcal{Z}}\left(  P,\gamma\right)  \left.  \prec\right.
\Gamma_{\mathcal{Z}}\left(  Q,\delta\right)   &  =\Gamma_{\mathcal{Z}}\left(
P\left.  \prec\right.  Q,\varepsilon\right)  \ \ \ \ \ \ \ \ \ \ \text{and}\\
\Gamma_{\mathcal{Z}}\left(  P,\gamma\right)  \left.  \succeq\right.
\Gamma_{\mathcal{Z}}\left(  Q,\delta\right)   &  =\Gamma_{\mathcal{Z}}\left(
P\left.  \succeq\right.  Q,\varepsilon\right)  ,
\end{align*}
assuming a simple condition on $\min P$ and $\min Q$
(say, $\gamma\left(\min P\right) <_{\ZZ} \delta\left(\min Q\right)$).
\end{question}
\end{vershort}

\begin{verlong}
Next, we shall study how the ideal $\mathcal{K}_{\operatorname{Epk}}$
interacts with some additional structure on $\operatorname*{QSym}$, viz. the
\textit{dendriform operations }$\left.  \prec\right.  $ and $\left.
\succeq\right.  $ and the \textquotedblleft runic\textquotedblright%
\ operations $\bel$ and $\tvi$. These operations were introduced in
\cite{dimcr}. Our study shall lead us back to the notions of
left-shuffle-compatibility and right-shuffle-compatibility from Section
\ref{sect.LR}. We shall reprove that $\operatorname{Epk}$ is
left-shuffle-compatible and right-shuffle-compatible; similar studies can
probably be made for other descent statistics.

\subsection{Four operations on $\operatorname*{QSym}$}

We begin with some definitions. We will use some notations from \cite{dimcr},
but we set $\mathbf{k}=\mathbb{Q}$ because we are working over the ring
$\mathbb{Q}$ in this paper. We agree that $\max\varnothing$ is understood as
$0$, while $\min\varnothing$ is understood as $\infty$ (a formal symbol that
is larger than any integer).
Monomials always mean formal expressions of the
form $x_{1}^{a_{1}}x_{2}^{a_{2}}x_{3}^{a_{3}}\cdots$ with $a_{1}+a_{2}%
+a_{3}+\cdots<\infty$ (see \cite[Section 2]{dimcr} for details). If
$\mathfrak{m}$ is a monomial, then $\operatorname*{Supp}\mathfrak{m}$ will
denote the finite subset
\[
\left\{  i\in\left\{  1,2,3,\ldots\right\}  \ \mid\ \text{the exponent with
which }x_{i}\text{ occurs in }\mathfrak{m}\text{ is }>0\right\}
\]
of $\left\{  1,2,3,\ldots\right\}  $. Next, we define four binary operations
\begin{align*}
&  \left.  \prec\right.  \ \left(  \text{called \textquotedblleft dendriform
less-than\textquotedblright; but it's an operation, not a relation}\right)
,\\
&  \left.  \succeq\right.  \ \left(  \text{called \textquotedblleft dendriform
greater-or-equal\textquotedblright; but it's an operation, not a
relation}\right)  ,\\
&  \bel\ \left(  \text{called \textquotedblleft belgthor\textquotedblright%
}\right)  ,\\
&  \tvi\ \left(  \text{called \textquotedblleft tvimadur\textquotedblright%
}\right)
\end{align*}
on the ring $\mathbf{k}\left[  \left[  x_{1},x_{2},x_{3},\ldots\right]
\right]  $ of power series by first defining how they act on monomials:%
\begin{align*}
\mathfrak{m}\left.  \prec\right.  \mathfrak{n}  &  =\left\{
\begin{array}
[c]{l}%
\mathfrak{m}\cdot\mathfrak{n},\ \ \ \ \ \ \ \ \ \ \text{if }\min\left(
\operatorname*{Supp}\mathfrak{m}\right)  <\min\left(  \operatorname*{Supp}%
\mathfrak{n}\right)  ;\\
0,\ \ \ \ \ \ \ \ \ \ \text{if }\min\left(  \operatorname*{Supp}%
\mathfrak{m}\right)  \geq\min\left(  \operatorname*{Supp}\mathfrak{n}\right)
\end{array}
\right.  ;\\
\mathfrak{m}\left.  \succeq\right.  \mathfrak{n}  &  =\left\{
\begin{array}
[c]{l}%
\mathfrak{m}\cdot\mathfrak{n},\ \ \ \ \ \ \ \ \ \ \text{if }\min\left(
\operatorname*{Supp}\mathfrak{m}\right)  \geq\min\left(  \operatorname*{Supp}%
\mathfrak{n}\right)  ;\\
0,\ \ \ \ \ \ \ \ \ \ \text{if }\min\left(  \operatorname*{Supp}%
\mathfrak{m}\right)  <\min\left(  \operatorname*{Supp}\mathfrak{n}\right)
\end{array}
\right.  ;\\
\mathfrak{m}\bel\mathfrak{n}  &  =\left\{
\begin{array}
[c]{l}%
\mathfrak{m}\cdot\mathfrak{n},\ \ \ \ \ \ \ \ \ \ \text{if }\max\left(
\operatorname*{Supp}\mathfrak{m}\right)  \leq\min\left(  \operatorname*{Supp}%
\mathfrak{n}\right)  ;\\
0,\ \ \ \ \ \ \ \ \ \ \text{if }\max\left(  \operatorname*{Supp}%
\mathfrak{m}\right)  >\min\left(  \operatorname*{Supp}\mathfrak{n}\right)
\end{array}
\right.  ;\\
\mathfrak{m}\tvi\mathfrak{n}  &  =\left\{
\begin{array}
[c]{l}%
\mathfrak{m}\cdot\mathfrak{n},\ \ \ \ \ \ \ \ \ \ \text{if }\max\left(
\operatorname*{Supp}\mathfrak{m}\right)  <\min\left(  \operatorname*{Supp}%
\mathfrak{n}\right)  ;\\
0,\ \ \ \ \ \ \ \ \ \ \text{if }\max\left(  \operatorname*{Supp}%
\mathfrak{m}\right)  \geq\min\left(  \operatorname*{Supp}\mathfrak{n}\right)
\end{array}
\right.  ;
\end{align*}
and then requiring that they all be $\mathbf{k}$-bilinear and continuous (so
their action on pairs of arbitrary power series can be computed by
\textquotedblleft opening the parentheses\textquotedblright). These operations
$\left.  \prec\right.  $, $\left.  \succeq\right.  $, $\bel$ and $\tvi$ all
restrict to the subset $\operatorname*{QSym}$ of $\mathbf{k}\left[  \left[
x_{1},x_{2},x_{3},\ldots\right]  \right]  $ (this is proven in \cite[detailed
version, Section 3]{dimcr}). They furthermore satisfy numerous
relations\footnote{These relations are all easy to prove (by linearity, it
suffices to verify them on monomials only, and this verification is
straightforward). A proof of the associativity of $\bel$ was given in
\cite[detailed version, Proposition 3.4]{dimcr}.}:

\begin{itemize}
\item The dendriform operations satisfy the four rules
\begin{align}
a\left.  \prec\right.  b+a\left.  \succeq\right.  b  &
=ab;\label{eq.dendriform.1}\\
\left(  a\left.  \prec\right.  b\right)  \left.  \prec\right.  c  &  =a\left.
\prec\right.  \left(  bc\right)  ;\nonumber\\
\left(  a\left.  \succeq\right.  b\right)  \left.  \prec\right.  c  &
=a\left.  \succeq\right.  \left(  b\left.  \prec\right.  c\right)
;\nonumber\\
a\left.  \succeq\right.  \left(  b\left.  \succeq\right.  c\right)   &
=\left(  ab\right)  \left.  \succeq\right.  c\nonumber
\end{align}
for all $a,b,c\in\mathbf{k}\left[  \left[  x_{1},x_{2},x_{3},\ldots\right]
\right]  $. (In other words, they turn $\mathbf{k}\left[  \left[  x_{1}%
,x_{2},x_{3},\ldots\right]  \right]  $ into what is called a dendriform algebra.)

\item For any $a\in\mathbf{k}\left[  \left[  x_{1},x_{2},x_{3},\ldots\right]
\right]  $, we have%
\begin{align}
1\left.  \prec\right.  a  &  =0;\label{eq.dendriform.1<a}\\
a\left.  \prec\right.  1  &  =a-\varepsilon\left(  a\right)
;\label{eq.dendriform.a<1}\\
1\left.  \succeq\right.  a  &  =a;\label{eq.dendriform.1>=a}\\
a\left.  \succeq\right.  1  &  =\varepsilon\left(  a\right)  ,
\label{eq.dendriform.a>=1}%
\end{align}
where $\varepsilon\left(  a\right)  $ denotes the constant term of the power
series $a$.

\item The binary operation $\bel$ is associative and unital (with $1$ serving
as the unity).

\item The binary operation $\tvi$ is associative and unital (with $1$ serving
as the unity).
\end{itemize}

Recall that we are using the notations $M_{\alpha}$ for the monomial
quasisymmetric functions and $F_{\alpha}$ for the fundamental quasisymmetric functions.

\begin{itemize}
\item For any two nonempty compositions $\alpha$ and $\beta$, we have
$M_{\alpha}\bel  M_{\beta}=M_{\left[  \alpha,\beta\right]  }+M_{\alpha
\odot\beta}$, where $\left[  \alpha,\beta\right]  $ and $\alpha\odot\beta$ are
two compositions defined by%
\begin{align*}
\left[  \left(  \alpha_{1},\alpha_{2},\ldots,\alpha_{\ell}\right)  ,\left(
\beta_{1},\beta_{2},\ldots,\beta_{m}\right)  \right]   &  =\left(  \alpha
_{1},\alpha_{2},\ldots,\alpha_{\ell},\beta_{1},\beta_{2},\ldots,\beta
_{m}\right)  ;\\
\left(  \alpha_{1},\alpha_{2},\ldots,\alpha_{\ell}\right)  \odot\left(
\beta_{1},\beta_{2},\ldots,\beta_{m}\right)   &  =\left(  \alpha_{1}%
,\alpha_{2},\ldots,\alpha_{\ell-1},\alpha_{\ell}+\beta_{1},\beta_{2},\beta
_{3},\ldots,\beta_{m}\right)  .
\end{align*}

\item For any two compositions $\alpha$ and $\beta$, we have $M_{\alpha
}\tvi  M_{\beta}=M_{\left[  \alpha,\beta\right]  }$.

\item For any two compositions $\alpha$ and $\beta$, we have $F_{\alpha
}\bel  F_{\beta}=F_{\alpha\odot\beta}$. (Here, $\alpha\odot\beta$ is defined
to be $\alpha$ if $\beta$ is the empty composition, and is defined to be
$\beta$ if $\alpha$ is the empty composition.)

\item For any two compositions $\alpha$ and $\beta$, we have $F_{\alpha
}\tvi  F_{\beta}=F_{\left[  \alpha,\beta\right]  }$.
\end{itemize}

Furthermore, we shall use two theorems from \cite[detailed version, Section
3]{dimcr}:

\begin{theorem}
\label{thm.beldend}Let $S$ denote the antipode of the Hopf algebra
$\operatorname*{QSym}$. Let us use Sweedler's notation $\sum_{\left(
b\right)  }b_{\left(  1\right)  }\otimes b_{\left(  2\right)  }$ for
$\Delta\left(  b\right)  $, where $b$ is any element of $\operatorname*{QSym}%
$. Then,%
\[
\sum_{\left(  b\right)  }\left(  S\left(  b_{\left(  1\right)  }\right)
\bel  a\right)  b_{\left(  2\right)  }=a\left.  \prec\right.  b
\]
for any $a\in\mathbf{k}\left[  \left[  x_{1},x_{2},x_{3},\ldots\right]
\right]  $ and $b\in\operatorname*{QSym}$ such that the constant term
of $a$ is $0$.
\end{theorem}

\begin{theorem}
\label{thm.tvidend'}Let $S$ denote the antipode of the Hopf algebra
$\operatorname*{QSym}$. Let us use Sweedler's notation $\sum_{\left(
b\right)  }b_{\left(  1\right)  }\otimes b_{\left(  2\right)  }$ for
$\Delta\left(  b\right)  $, where $b$ is any element of $\operatorname*{QSym}%
$. Then,%
\[
\sum_{\left(  b\right)  }\left(  S\left(  b_{\left(  1\right)  }\right)
\tvi  a\right)  b_{\left(  2\right)  }=b\left.  \succeq\right.  a
\]
for any $a\in\mathbf{k}\left[  \left[  x_{1},x_{2},x_{3},\ldots\right]
\right]  $ and $b\in\operatorname*{QSym}$.
\end{theorem}

(Notice that Theorem \ref{thm.tvidend'} differs from \cite[detailed version,
Theorem 3.15]{dimcr} in that we are writing $b\left.  \succeq\right.  a$
instead of $a\left.  \preceq\right.  b$. But this is the same thing, since
$a\left.  \preceq\right.  b=b\left.  \succeq\right.  a$ for all $a,b\in
\mathbf{k}\left[  \left[  x_{1},x_{2},x_{3},\ldots\right]  \right]  $.)

\begin{noncompile}
The following fact is neat, but has so far been useless:

\begin{proposition}
\label{prop.tvibel.antipode}Let $S$ denote the antipode of the Hopf algebra
$\operatorname*{QSym}$. Let $a\in\operatorname*{QSym}$ and $b\in
\operatorname*{QSym}$. Then,%
\begin{align*}
S\left(  a\tvi b\right)   &  =S\left(  b\right)  \bel S\left(  a\right)
\ \ \ \ \ \ \ \ \ \ \text{and}\\
S\left(  a\bel b\right)   &  =S\left(  b\right)  \tvi S\left(  a\right)  .
\end{align*}

\end{proposition}
\end{noncompile}

\begin{noncompile}
Let us restate Theorem \ref{thm.beldend} and Theorem \ref{thm.tvidend'} using
a slight variation on the comultiplication map $\Delta$:

\begin{theorem}
\label{thm.beldend.Del'}Let $S$ denote the antipode of the Hopf algebra
$\operatorname*{QSym}$. Let us use the Sweedler-like notation $\sum_{\left(
b\right)  }b_{\left[  1\right]  }\otimes b_{\left[  2\right]  }$ for
$\Delta\left(  b\right)  -1\otimes b$, where $b$ is any element of
$\operatorname*{QSym}$. Then,%
\[
\sum_{\left(  b\right)  }\left(  S\left(  b_{\left[  1\right]  }\right)
\bel
a\right)  b_{\left[  2\right]  }=-a\mathfrak{\left.  \succeq\right.  }b
\]
for any $a\in\mathbf{k}\left[  \left[  x_{1},x_{2},x_{3},\ldots\right]
\right]  $ and $b\in\operatorname*{QSym}$ such that the constant term
of $a$ is $0$.
\end{theorem}

\begin{theorem}
\label{thm.tvidend'.Del'}Let $S$ denote the antipode of the Hopf algebra
$\operatorname*{QSym}$. Let us use the Sweedler-like notation $\sum_{\left(
b\right)  }b_{\left[  1\right]  }\otimes b_{\left[  2\right]  }$ for
$\Delta\left(  b\right)  -1\otimes b$, where $b$ is any element of
$\operatorname*{QSym}$. Then,%
\[
\sum_{\left(  b\right)  }\left(  S\left(  b_{\left[  1\right]  }\right)
\tvi
a\right)  b_{\left[  2\right]  }=-b\left.  \prec\right.  a
\]
for any $a\in\mathbf{k}\left[  \left[  x_{1},x_{2},x_{3},\ldots\right]
\right]  $ and $b\in\operatorname*{QSym}$.
\end{theorem}

\begin{proof}
[Proof of Theorem \ref{thm.beldend.Del'}.]Let
$a\in\mathbf{k}\left[  \left[  x_{1},x_{2},x_{3},\ldots\right]
\right]  $ and $b\in\operatorname*{QSym}$ be such that the constant term
of $a$ is $0$.

Let us use Sweedler's notation, too.
From $\sum_{\left(  b\right)  }b_{\left[  1\right]  }\otimes b_{\left[
2\right]  }=\Delta\left(  b\right)  -1\otimes b=\sum_{\left(  b\right)
}b_{\left(  1\right)  }\otimes b_{\left(  2\right)  }-1\otimes b$, we have%
\begin{align*}
\sum_{\left(  b\right)  }\left(  S\left(  b_{\left[  1\right]  }\right)
\bel a\right)  b_{\left[  2\right]  }  &  =\underbrace{\sum_{\left(  b\right)
}\left(  S\left(  b_{\left(  1\right)  }\right)  \bel a\right)  b_{\left(
2\right)  }}_{\substack{=a\left.  \prec\right.  b\\\text{(by Theorem
\ref{thm.beldend})}}}-\underbrace{\left(  1\bel a\right)  }_{=a}b\\
&  =a\left.  \prec\right.  b-ab=-a\left.  \succeq\right.  b
\end{align*}
(since $a\left.  \prec\right.  b+a\left.  \succeq\right.  b=ab$). This proves
Theorem \ref{thm.beldend.Del'}.
\end{proof}

\begin{proof}
[Proof of Theorem \ref{thm.tvidend'.Del'}.]Analogous.
\end{proof}
\end{noncompile}

\subsection{The dendriform operations on the fundamental basis}

Recall Definition \ref{def.LRshuf} and Definition \ref{def.LR.pi1}. The
following theorem is analogous to Theorem~\ref{thm.4.1}:

\begin{theorem}
\label{thm.dendri.4.1}Let $\pi$ be a nonempty permutation with descent
composition $J$. Let $\sigma$ be a nonempty permutation with descent
composition $K$. Assume that the permutations $\pi$ and $\sigma$ are disjoint,
and that $\pi_{1}>\sigma_{1}$. For any composition $L$, let $c_{J,K}^{L,\prec
}$ be the number of permutations with descent composition $L$ among the left
shuffles of $\pi$ and $\sigma$, and let $c_{J,K}^{L,\succ}$ be the number of
permutations with descent composition $L$ among the right shuffles of $\pi$
and $\sigma$. Then,%
\[
F_{J}\left.  \prec\right.  F_{K}=\sum_{L}c_{J,K}^{L,\prec}F_{L}%
\]
and%
\[
F_{J}\left.  \succeq\right.  F_{K}=\sum_{L}c_{J,K}^{L,\succ}F_{L}.
\]

\end{theorem}

Note the condition $\pi_{1}>\sigma_{1}$, which is not present in
Theorem~\ref{thm.4.1}, and which makes Theorem \ref{thm.dendri.4.1} somewhat
harder to apply.

Theorem \ref{thm.dendri.4.1} can be proven similarly to
\cite[(5.2.6)]{HopfComb}, but it relies on some variants of the disjoint
union of two posets.
We shall show this proof after first establishing some auxiliary facts.

Throughout this section, the
notion of a ``labelled poset'' will be understood in the sense
of \cite[Definition 5.2.1]{HopfComb}:
Namely, a \textit{labelled poset} simply means a poset whose underlying
set is a finite subset of $\ZZ$ (but whose order is not necessarily
inherited from $\ZZ$).
Thus, a labelled poset is not equipped with a map that serves as its
labelling,
but instead its underlying set must be a finite set of integers.

Let us first recall a well-known fact about posets:

\begin{lemma}
\label{lem.poset.add-rel}
Let $R$ be a poset. Let $u$ and $v$ be two distinct elements of $R$
such that we don't have $u > v$ in $R$.
Then, there exists a unique poset $R^{\prime}$ such that
\begin{itemize}
\item we have $R^{\prime} = R$ as sets, and
\item
we have the logical equivalence
\[
\left(  x<y\text{ in }R^{\prime}\right)  \ \Longleftrightarrow\ \left(
\left(  x<y\text{ in }R\right)  \text{ or }\left(  x\leq u\text{ and }v\leq
y\text{ in }R\right)  \right)
\]
for every two elements $x$ and $y$ of $R^{\prime}$.
\end{itemize}
% Furthermore, this poset $R^{\prime}$ is the ``smallest'' extension
% of $R$ in which the inequality $u < v$ holds (i.e., it is an
% extension of $R$ in which $u < v$ holds, and every extension of
% $R$ in which $u < v$ holds must be an extension of $R^{\prime}$).
We say that this poset $R^{\prime}$ is obtained by
\textit{adding the relation $u < v$ to $R$}.
Clearly, $u < v$ holds in $R^{\prime}$.

If $R$ is a labelled poset in the sense of \cite[Section
5.2]{HopfComb}, then $R^{\prime}$ is a labelled poset in the
sense of \cite[Section 5.2]{HopfComb} as well (since $R^{\prime} = R$
as sets).
\end{lemma}

If we add a relation to a labelled poset $R$, obtaining a new labelled
poset $R^{\prime}$, then how do the $R^{\prime}$-partitions differ from
the $R$-partitions?
The following two lemmas answer this question:\footnote{We are
going to use the following notation:
If $P$ is a poset, then we let $\leq_P$, $<_P$, $\geq_P$
and $>_P$ denote the smaller-or-equal relation of $P$, the
smaller relation of $P$, the greater-or-equal relation
of $P$, and the greater relation of $P$, respectively.
Thus, in particular, $\leq_{\ZZ}$, $<_{\ZZ}$, $\geq_{\ZZ}$
and $>_{\ZZ}$ denote the usual smaller-or-equal, smaller,
greater-or-equal and greater relations of the totally
ordered set $\ZZ$ of integers.
(For example, $a \leq_{\ZZ} b$ if and only if $b - a \in \NN$.)}

\begin{lemma}
\label{lem.poset.add-rel.part1}
Let $R$ be a labelled poset. Let $u$ and $v$ be
two elements of $R$ such that $u >_{\ZZ} v$.
Assume that we don't have $u > v$ in $R$.
Let $f$ be an $R$-partition.
Let $R^{\prime}$ be the labelled poset obtained from $R$ by adding the
relation $u < v$ to $R$.
Then, $f$ is an $R^{\prime}$-partition if and only if
$f \left(u\right) < f\left(v\right)$.
\end{lemma}

\begin{proof}[Proof of Lemma~\ref{lem.poset.add-rel.part1}.]
The \textquotedblleft only
if\textquotedblright\ direction is obvious (since $u<v$ in $R^{\prime}$ and
$u>_{\ZZ}v$). It thus remains to prove the \textquotedblleft
if\textquotedblright\ direction. So let us assume that $f\left(  u\right)
<f\left(  v\right)  $. We must then show that $f$ is an $R^{\prime}$-partition.

The poset $R^{\prime}$ is obtained from $R$ by adding the relation $u<v$ to
$R$. Thus, $R^{\prime}=R$ as sets, and we have the logical equivalence%
\[
\left(  x<y\text{ in }R^{\prime}\right)  \ \Longleftrightarrow\ \left(
\left(  x<y\text{ in }R\right)  \text{ or }\left(  x\leq u\text{ and }v\leq
y\text{ in }R\right)  \right)
\]
for every two elements $x$ and $y$ of $R^{\prime}$. Thus, any pair $\left(
x,y\right)  $ of elements of $R^{\prime}$ that satisfies $x<y$ in $R^{\prime}$
must either already satisfy $x<y$ in $R$, or satisfy $x\leq u$ and $v\leq y$
in $R$.

The map $f$ is an $R$-partition, and thus is weakly increasing as a map from
$R$ to $\mathbb{Z}$.

Now, we have%
\[
\left(  \text{if }i\in R^{\prime}\text{ and }j\in R^{\prime}\text{ satisfy
}i<j\text{ in }R^{\prime}\text{ and }i<_{\ZZ}j\text{, then }f\left(
i\right)  \leq f\left(  j\right)  \right)
\]
\footnote{\textit{Proof.} Let $i\in R^{\prime}$ and $j\in R^{\prime}$ be such
that $i<j$ in $R^\prime$ and $i<_{\ZZ}j$. We must prove that $f\left(
i\right)  \leq f\left(  j\right)  $.
\par
We have $i\in R^{\prime}=R$ and $j\in R^{\prime}=R$.
\par
Recall that any pair $\left(  x,y\right)  $ of elements of $R^{\prime}$ that
satisfies $x<y$ in $R^{\prime}$ must either already satisfy $x<y$ in $R$, or
satisfy $x\leq u$ and $v\leq y$ in $R$. Applying this to $\left(  x,y\right)
=\left(  i,j\right)  $, we conclude that the pair $\left(  i,j\right)  $ must
either already satisfy $i<j$ in $R$, or satisfy $i\leq u$ and $v\leq j$ in
$R$. In the first of these two cases, we immediately obtain $f\left(
i\right)  \leq f\left(  j\right)  $, because $f$ is an $R$-partition (and
because $i<j$ in $R$ and $i<_{\ZZ}j$). Hence, we can WLOG assume that
we are in the second case. In other words, we have $i\leq u$ and $v\leq j$ in
$R$. In particular, $i\leq u$ in $R$ and therefore $f\left(  i\right)  \leq
f\left(  u\right)  $ (since the map $f$ is weakly increasing). Similarly,
$f\left(  v\right)  \leq f\left(  j\right)  $. Thus, $f\left(  i\right)  \leq
f\left(  u\right)  <f\left(  v\right)  \leq f\left(  j\right)  $, so that
$f\left(  i\right)  \leq f\left(  j\right)  $, qed.} and%
\[
\left(  \text{if }i\in R^{\prime}\text{ and }j\in R^{\prime}\text{ satisfy
}i<j\text{ in }R^{\prime}\text{ and }i>_{\ZZ}j\text{, then }f\left(
i\right)  <f\left(  j\right)  \right)
\]
\footnote{\textit{Proof.} Let $i\in R^{\prime}$ and $j\in R^{\prime}$ be such
that $i<j$ in $R^\prime$ and $i>_{\ZZ}j$. We must prove that $f\left(
i\right)  <f\left(  j\right)  $.
\par
We have $i\in R^{\prime}=R$ and $j\in R^{\prime}=R$.
\par
Recall that any pair $\left(  x,y\right)  $ of elements of $R^{\prime}$ that
satisfies $x<y$ in $R^{\prime}$ must either already satisfy $x<y$ in $R$, or
satisfy $x\leq u$ and $v\leq y$ in $R$. Applying this to $\left(  x,y\right)
=\left(  i,j\right)  $, we conclude that the pair $\left(  i,j\right)  $ must
either already satisfy $i<j$ in $R$, or satisfy $i\leq u$ and $v\leq j$ in
$R$. In the first of these two cases, we immediately obtain $f\left(
i\right)  <f\left(  j\right)  $, because $f$ is an $R$-partition (and because
$i<j$ in $R$ and $i>_{\ZZ}j$). Hence, we can WLOG assume that we are in
the second case. In other words, we have $i\leq u$ and $v\leq j$ in $R$. In
particular, $i\leq u$ in $R$ and therefore $f\left(  i\right)  \leq f\left(
u\right)  $ (since the map $f$ is weakly increasing). Similarly, $f\left(
v\right)  \leq f\left(  j\right)  $. Thus, $f\left(  i\right)  \leq f\left(
u\right)  <f\left(  v\right)  \leq f\left(  j\right)  $, so that $f\left(
i\right)  <f\left(  j\right)  $, qed.}. These two statements entail that $f$
is an $R^{\prime}$-partition (by the definition of an $R^{\prime}$-partition).
This concludes the proof of the \textquotedblleft if\textquotedblright%
\ direction of Lemma~\ref{lem.poset.add-rel.part1}.
Thus, Lemma~\ref{lem.poset.add-rel.part1} is proven.
\end{proof}

\begin{lemma}
\label{lem.poset.add-rel.part2}
Let $R$ be a labelled poset. Let $u$ and $v$ be
two elements of $R$ such that $u <_{\ZZ} v$.
Assume that we don't have $u > v$ in $R$.
Let $f$ be an $R$-partition.
Let $R^{\prime}$ be the labelled poset obtained from $R$ by adding the
relation $u < v$ to $R$.
Then, $f$ is an $R^{\prime}$-partition if and only if
$f \left(u\right) \leq f\left(v\right)$.
\end{lemma}

\begin{proof}[Proof of Lemma~\ref{lem.poset.add-rel.part2}.]
The \textquotedblleft only
if\textquotedblright\ direction is obvious (since $u<v$ in $R^{\prime}$ and
$u<_{\ZZ}v$). It thus remains to prove the \textquotedblleft
if\textquotedblright\ direction. So let us assume that $f\left(  u\right)
\leq f\left(  v\right)  $. We must then show that $f$ is an $R^{\prime}$-partition.

The poset $R^{\prime}$ is obtained from $R$ by adding the relation $u<v$ to
$R$. Thus, $R^{\prime}=R$ as sets, and we have the logical equivalence%
\[
\left(  x<y\text{ in }R^{\prime}\right)  \ \Longleftrightarrow\ \left(
\left(  x<y\text{ in }R\right)  \text{ or }\left(  x\leq u\text{ and }v\leq
y\text{ in }R\right)  \right)
\]
for every two elements $x$ and $y$ of $R^{\prime}$. Thus, any pair $\left(
x,y\right)  $ of elements of $R^{\prime}$ that satisfies $x<y$ in $R^{\prime}$
must either already satisfy $x<y$ in $R$, or satisfy $x\leq u$ and $v\leq y$
in $R$.

The map $f$ is an $R$-partition, and thus is weakly increasing as a map from
$R$ to $\mathbb{Z}$.

Now, we have%
\[
\left(  \text{if }i\in R^{\prime}\text{ and }j\in R^{\prime}\text{ satisfy
}i<j\text{ in }R^{\prime}\text{ and }i<_{\ZZ}j\text{, then }f\left(
i\right)  \leq f\left(  j\right)  \right)
\]
\footnote{\textit{Proof.} Let $i\in R^{\prime}$ and $j\in R^{\prime}$ be such
that $i<j$ in $R^\prime$ and $i<_{\ZZ}j$. We must prove that $f\left(
i\right)  \leq f\left(  j\right)  $.
\par
We have $i\in R^{\prime}=R$ and $j\in R^{\prime}=R$.
\par
Recall that any pair $\left(  x,y\right)  $ of elements of $R^{\prime}$ that
satisfies $x<y$ in $R^{\prime}$ must either already satisfy $x<y$ in $R$, or
satisfy $x\leq u$ and $v\leq y$ in $R$. Applying this to $\left(  x,y\right)
=\left(  i,j\right)  $, we conclude that the pair $\left(  i,j\right)  $ must
either already satisfy $i<j$ in $R$, or satisfy $i\leq u$ and $v\leq j$ in
$R$. In the first of these two cases, we immediately obtain $f\left(
i\right)  \leq f\left(  j\right)  $, because $f$ is an $R$-partition (and
because $i<j$ in $R$ and $i<_{\ZZ}j$). Hence, we can WLOG assume that
we are in the second case. In other words, we have $i\leq u$ and $v\leq j$ in
$R$. In particular, $i\leq u$ in $R$ and therefore $f\left(  i\right)  \leq
f\left(  u\right)  $ (since the map $f$ is weakly increasing). Similarly,
$f\left(  v\right)  \leq f\left(  j\right)  $. Thus, $f\left(  i\right)  \leq
f\left(  u\right)  \leq f\left(  v\right)  \leq f\left(  j\right)  $, so that
$f\left(  i\right)  \leq f\left(  j\right)  $, qed.} and%
\[
\left(  \text{if }i\in R^{\prime}\text{ and }j\in R^{\prime}\text{ satisfy
}i<j\text{ in }R^{\prime}\text{ and }i>_{\ZZ}j\text{, then }f\left(
i\right)  <f\left(  j\right)  \right)
\]
\footnote{\textit{Proof.} Let $i\in R^{\prime}$ and $j\in R^{\prime}$ be such
that $i<j$ in $R^\prime$ and $i>_{\ZZ}j$. We must prove that $f\left(
i\right)  <f\left(  j\right)  $.
\par
Assume the contrary. Thus, $f\left(  i\right)  \geq f\left(  j\right)  $, so
that $f\left(  j\right)  \leq f\left(  i\right)  $.
\par
We have $i\in R^{\prime}=R$ and $j\in R^{\prime}=R$.
\par
Recall that any pair $\left(  x,y\right)  $ of elements of $R^{\prime}$ that
satisfies $x<y$ in $R^{\prime}$ must either already satisfy $x<y$ in $R$, or
satisfy $x\leq u$ and $v\leq y$ in $R$. Applying this to $\left(  x,y\right)
=\left(  i,j\right)  $, we conclude that the pair $\left(  i,j\right)  $ must
either already satisfy $i<j$ in $R$, or satisfy $i\leq u$ and $v\leq j$ in
$R$. In the first of these two cases, we immediately obtain $f\left(
i\right)  <f\left(  j\right)  $, because $f$ is an $R$-partition (and because
$i<j$ in $R$ and $i>_{\ZZ}j$). Hence, we can WLOG assume that we are in
the second case. In other words, we have $i\leq u$ and $v\leq j$ in $R$. In
particular, $i\leq u$ in $R$ and therefore $f\left(  i\right)  \leq f\left(
u\right)  $ (since the map $f$ is weakly increasing). Similarly, $f\left(
v\right)  \leq f\left(  j\right)  $. Thus, $f\left(  i\right)  \leq f\left(
u\right)  \leq f\left(  v\right)  \leq f\left(  j\right)  \leq f\left(
i\right)  $. All inequality signs in this chain of inequalities must be
equalities (since its left and right hand sides are equal). In other words, we
have $f\left(  i\right)  =f\left(  u\right)  =f\left(  v\right)  =f\left(
j\right)  =f\left(  i\right)  $.
\par
If we had $i>_{\ZZ}u$, then we would have $i\neq u$ and therefore $i<u$
in $R$ (since $i\leq u$ in $R$). Therefore, if we had $i>_{\ZZ}u$, then
we would have $f\left(  i\right)  <f\left(  u\right)  $ (since $i<u$ in $R$,
but $f$ is an $R$-partition), which would contradict $f\left(  i\right)
=f\left(  u\right)  $. Hence, we cannot have $i>_{\ZZ}u$. Thus,
$i\leq_{\ZZ}u$.
\par
If we had $v>_{\ZZ}j$, then we would have $v\neq j$ and therefore $v<j$
in $R$ (since $v\leq j$ in $R$). Therefore, if we had $v>_{\ZZ}j$, then
we would have $f\left(  v\right)  <f\left(  j\right)  $ (since $v<j$ in $R$,
but $f$ is an $R$-partition), which would contradict $f\left(  v\right)
=f\left(  j\right)  $. Hence, we cannot have $v>_{\ZZ}j$. Thus,
$v\leq_{\ZZ}j$.
\par
Now, $i\leq_{\ZZ}u<_{\ZZ}v\leq_{\ZZ}j$. This contradicts
$i>_{\ZZ}j$. This contradiction shows that our assumption was wrong,
qed.}. These two statements entail that $f$ is an $R^{\prime}$-partition (by
the definition of an $R^{\prime}$-partition). This concludes the proof of the
\textquotedblleft if\textquotedblright\ direction of
Lemma~\ref{lem.poset.add-rel.part2}. Thus,
Lemma~\ref{lem.poset.add-rel.part2} is proven.
\end{proof}

For our next results, we need the following notation:

\begin{definition} \label{def.labelled-poset.min}
A \textit{minimum element} of a poset $P$ means an element
$m \in P$ such that every $p \in P$ satisfies $m \leq p$.
This is not the same as a minimal element.
A minimum element of a poset $P$ is unique if it exists;
it is denoted by $\min P$.

This notation overrides the notation $\min S$ for a nonempty
finite subset $S$ of $\ZZ$.
Thus, if $P$ is a labelled poset that has a minimum element,
then $\min P$ shall denote this minimum element,
not the smallest element of $P$ as a subset of $\ZZ$.
(For example, if $P$ is the labelled poset $\left\{6 < 2 < 3 < 4\right\}$,
then $\min P$ shall mean $6$, not $2$.)
\end{definition}

The following fact is an analogue of \cite[Lemma 5.2.17]{HopfComb}:

\begin{proposition}
\label{prop.dendri.4.1.lem}We shall use the notations of \cite[Section
5.2]{HopfComb}. Let $P$ and $Q$ be two disjoint labelled posets, each of which
has a minimum element. Assume that $\min P>_{\ZZ}\min Q$. Consider the
disjoint union $P\sqcup Q$ of $P$ and $Q$.

\textbf{(a)} Add a further relation $\min P<\min Q$ to $P\sqcup Q$; denote the
resulting labelled poset by $P\left.  \prec\right.  Q$. Then, $F_{P}\left(
\mathbf{x}\right)  \left.  \prec\right.  F_{Q}\left(  \mathbf{x}\right)
=F_{P\left.  \prec\right.  Q}\left(  \mathbf{x}\right)  $.

\textbf{(b)} Add a further relation $\min P>\min Q$ to $P\sqcup Q$; denote the
resulting labelled poset by $P\left.  \succeq\right.  Q$. Then, $F_{P}\left(
\mathbf{x}\right)  \left.  \succeq\right.  F_{Q}\left(  \mathbf{x}\right)
=F_{P\left.  \succeq\right.  Q}\left(  \mathbf{x}\right)  $.
\end{proposition}

Let us recall that (as agreed in Definition~\ref{def.labelled-poset.min})
we are using the notation $\min P$ for the minimum element
of the \textbf{poset} $P$; not the smallest element of $P$ as a subset
of $\ZZ$. The same applies to the notation $\min Q$.

\begin{proof}[Proof of Proposition~\ref{prop.dendri.4.1.lem} (sketched).]
We imitate the proof of \cite[Lemma 5.2.17]{HopfComb}:

\textbf{(a)} If $f : P \sqcup Q \to \left\{1,2,3,\ldots\right\}$ is a
$P \left. \prec \right. Q$-partition, then its restrictions $f\mid_P$
and $f\mid_Q$ are a $P$-partition and a $Q$-partition, respectively, and
have the property that
$\left(f\mid_P\right) \left( \min P \right) < \left(f\mid_Q\right) \left( \min Q \right)$
(indeed, this property must hold because $\min P < \min Q$ in
$P \left. \prec \right. Q$ but $\min P >_{\ZZ} \min Q$).
Conversely, any pair of a $P$-partition $g$ and a $Q$-partition $h$
having the property that $g \left(\min P\right) < h \left(\min Q\right)$
can be combined to form a
$P \left. \prec \right. Q$-partition\footnote{To see this, we need to
apply Lemma~\ref{lem.poset.add-rel.part1} to
$R = P \sqcup Q$, $u = \min P$, $v = \min Q$
and $R^{\prime} = P \left. \prec \right. Q$.}.

Using the notations of \cite[Definition 5.2.1]{HopfComb}, we have
\begin{equation}
\sum_{f\text{ is a }P\left.  \prec\right.  Q\text{-partition}}\mathbf{x}%
_{f}=\sum_{\substack{g\text{ is a }P\text{-partition;}\\h\text{ is a
}Q\text{-partition;}\\g\left(  \min P\right)  <h\left(  \min Q\right)
}}\mathbf{x}_{g}\mathbf{x}_{h}\label{pf.prop.dendri.4.1.lem.a.1}%
\end{equation}
(since the previous two sentences establish a bijection between the addends on
the left hand side of this equality and the addends on its right hand side).

But the definitions of $F_{P}\left(  \mathbf{x}\right)  $ and $F_{Q}\left(
\mathbf{x}\right)  $ yield
$F_{P}\left(  \mathbf{x}\right)
=  \sum_{\substack{g\text{ is a } P\text{-partition}}}\mathbf{x}_{g}$
and
$F_{Q}\left( \mathbf{x}\right)
= \sum_{\substack{h\text{ is a }Q\text{-partition}}}\mathbf{x}_{h}$.
Thus,
\begin{align*}
& F_{P}\left(  \mathbf{x}\right)  \left.  \prec\right.  F_{Q}\left(
\mathbf{x}\right)   \\
& =\left(  \sum_{\substack{g\text{ is a }%
P\text{-partition}}}\mathbf{x}_{g}\right)  \left.  \prec\right.  \left(
\sum_{\substack{h\text{ is a }Q\text{-partition}}}\mathbf{x}_{h}\right)  \\
& =\sum_{\substack{g\text{ is a }P\text{-partition;}\\h\text{ is a
}Q\text{-partition}}}\underbrace{\mathbf{x}_{g}\left.  \prec\right.
\mathbf{x}_{h}}_{\substack{=\left\{
\begin{array}
[c]{l}%
\mathbf{x}_{g}\mathbf{x}_{h},\ \ \ \ \ \ \ \ \ \ \text{if }\min\left(
\operatorname*{Supp}\left(  \mathbf{x}_{g}\right)  \right)  <\min\left(
\operatorname*{Supp}\left(  \mathbf{x}_{h}\right)  \right)  ;\\
0,\ \ \ \ \ \ \ \ \ \ \text{if }\min\left(  \operatorname*{Supp}\left(
\mathbf{x}_{g}\right)  \right)  \geq\min\left(  \operatorname*{Supp}\left(
\mathbf{x}_{h}\right)  \right)
\end{array}
\right.  \\\text{(by the definition of }\mathbf{x}_{g}\left.  \prec\right.
\mathbf{x}_{h}\text{)}}}\\
& =\sum_{\substack{g\text{ is a }P\text{-partition;}\\h\text{ is a
}Q\text{-partition}}}\underbrace{\left\{
\begin{array}
[c]{l}%
\mathbf{x}_{g}\mathbf{x}_{h},\ \ \ \ \ \ \ \ \ \ \text{if }\min\left(
\operatorname*{Supp}\left(  \mathbf{x}_{g}\right)  \right)  <\min\left(
\operatorname*{Supp}\left(  \mathbf{x}_{h}\right)  \right)  ;\\
0,\ \ \ \ \ \ \ \ \ \ \text{if }\min\left(  \operatorname*{Supp}\left(
\mathbf{x}_{g}\right)  \right)  \geq\min\left(  \operatorname*{Supp}\left(
\mathbf{x}_{h}\right)  \right)
\end{array}
\right.  }_{\substack{=\left\{
\begin{array}
[c]{l}%
\mathbf{x}_{g}\mathbf{x}_{h},\ \ \ \ \ \ \ \ \ \ \text{if }g\left(  \min
P\right)  <h\left(  \min Q\right)  ;\\
0,\ \ \ \ \ \ \ \ \ \ \text{if }g\left(  \min P\right)  \geq h\left(  \min
Q\right)
\end{array}
\right.  \\\text{(since }\min\left(  \operatorname*{Supp}\left(
\mathbf{x}_{g}\right)  \right)  =g\left(  \min P\right)  \text{ (because the
map }g\text{ is a }P\text{-partition, thus}\\\text{weakly increasing on
}P\text{, and therefore }g\left(  \min P\right)  =\min\left(  g\left(
P\right)  \right)  =\min\left(  \operatorname*{Supp}\left(  \mathbf{x}%
_{g}\right)  \right)  \text{)}\\\text{and similarly }\min\left(
\operatorname*{Supp}\left(  \mathbf{x}_{h}\right)  \right)  =h\left(  \min
Q\right)  \text{)}}}\\
& =\sum_{\substack{g\text{ is a }P\text{-partition;}\\h\text{ is a
}Q\text{-partition}}}\left\{
\begin{array}
[c]{l}%
\mathbf{x}_{g}\mathbf{x}_{h},\ \ \ \ \ \ \ \ \ \ \text{if }g\left(  \min
P\right)  <h\left(  \min Q\right)  ;\\
0,\ \ \ \ \ \ \ \ \ \ \text{if }g\left(  \min P\right)  \geq h\left(  \min
Q\right)
\end{array}
\right.  \\
& =\sum_{\substack{g\text{ is a }P\text{-partition;}\\h\text{ is a
}Q\text{-partition;}\\g\left(  \min P\right)  <h\left(  \min Q\right)
}}\mathbf{x}_{g}\mathbf{x}_{h}=\sum_{f\text{ is a }P\left.  \prec\right.
Q\text{-partition}}\mathbf{x}_{f}\ \ \ \ \ \ \ \ \ \ \left(  \text{by
(\ref{pf.prop.dendri.4.1.lem.a.1})}\right)  \\
& =F_{P\left.  \prec\right.  Q}\left(  \mathbf{x}\right)
\end{align*}
(by the definition of $F_{P\left.  \prec\right.  Q}\left(  \mathbf{x}\right)
$). This proves Proposition \ref{prop.dendri.4.1.lem} \textbf{(a)}.

\textbf{(b)} Proposition \ref{prop.dendri.4.1.lem} \textbf{(b)} is proven
similarly to Proposition \ref{prop.dendri.4.1.lem} \textbf{(a)}
(but this time we need to use Lemma~\ref{lem.poset.add-rel.part2}
instead of Lemma~\ref{lem.poset.add-rel.part1}).
\end{proof}

Also, the following simple fact is used:

\begin{lemma}
\label{lem.dendri.4.1.lem2}We shall use the notations of \cite[Section
5.2]{HopfComb}. Let $P$ and $Q$ be two disjoint posets, each of which has a
minimum element. Consider the disjoint union $P\sqcup Q$ of $P$ and $Q$ as the
set-theoretic union $P\cup Q$. Assume that $P$ and $Q$ are subsets of
$\mathbb{P}$ (the set of positive integers); thus, any linear extension of $P$
or of $Q$ or of $P\sqcup Q$ is a permutation (a word with letters in
$\mathbb{P}$).

\textbf{(a)} Add a further relation $\min P<\min Q$ to $P\sqcup Q$; denote the
resulting poset by $P\left.  \prec\right.  Q$. Then,
\[
\mathcal{L}\left(  P\left.  \prec\right.  Q\right)  =\bigsqcup_{\pi
\in\mathcal{L}\left(  P\right)  ;\ \sigma\in\mathcal{L}\left(  Q\right)
}S_{\prec}\left(  \pi,\sigma\right)  .
\]

\textbf{(b)} Add a further relation $\min P>\min Q$ to $P\sqcup Q$; denote the
resulting poset by $P\left.  \succeq\right.  Q$. Then,
\[
\mathcal{L}\left(  P\left.  \succeq\right.  Q\right)  =\bigsqcup_{\pi
\in\mathcal{L}\left(  P\right)  ;\ \sigma\in\mathcal{L}\left(  Q\right)
}S_{\succ}\left(  \pi,\sigma\right)  .
\]

\end{lemma}

\begin{proof}
[Proof of Lemma \ref{lem.dendri.4.1.lem2} (sketched).]Recall that we regard
linear extensions of a finite poset $R$ as lists of elements of $R$. Thus, if
$R$ is a finite poset, if $u$ and $v$ are two elements of $R$, and if $w$ is a
linear extension of $R$, then we have $u<v$ in $w$ if and only if $u$ appears
before $v$ in the list $w$.

We shall also use the following notation: If $w$ is a list of elements of some
set $U$, and if $V$ is a subset of $U$, then $w\mid_{V}$ means the result of
removing all entries from $w$ that don't belong to $V$. For example,
$\left. \left(
2,7,1,6,3,4\right)  \mid_{\left\{  1,2,4,5,6\right\}  } \right.
=\left( 2,1,6,4\right)  $.

Now, consider our two posets $P$ and $Q$. Recall that a linear extension of
$P\sqcup Q$ is simply a list $\left(  w_{1},w_{2},\ldots,w_{n}\right)  $ of
all elements of $P\sqcup Q$ such that no two integers $i<j$ satisfy $w_{i}\geq
w_{j}$ in $P\sqcup Q$. In other words, a linear extension of $P\sqcup Q$ is a
list $w$ of all elements of $P\sqcup Q$ such that if $x$ and $y$ are two
elements of $P\sqcup Q$ satisfying $x<y$ in $P\sqcup Q$, then $x$ appears
before\footnote{\textquotedblleft Before\textquotedblright\ doesn't imply
\textquotedblleft immediately before\textquotedblright. For example, $2$
appears before $4$ in the list $\left(  1,2,3,4\right)  $.} $y$ in the list
$w$. By the definition of $P\sqcup Q$, this rewrites as follows: A linear
extension of $P\sqcup Q$ is a list $w$ of all elements of $P\sqcup Q$ with the
following two properties:

\begin{itemize}
\item If $x$ and $y$ are two elements of $P$ satisfying $x<y$ in $P$, then $x$
appears before $y$ in the list $w$.

\item If $x$ and $y$ are two elements of $Q$ satisfying $x<y$ in $Q$, then $x$
appears before $y$ in the list $w$.
\end{itemize}

In other words, a linear extension of $P\sqcup Q$ is a list $w$ of all elements
of $P\sqcup Q$ such that $w\mid_{P}$ is a linear extension of $P$ and
$w\mid_{Q}$ is a linear extension of $Q$.

This yields the following:

\begin{itemize}
\item If $w$ is a linear extension of $P\sqcup Q$, then $w\mid_{P}$ is a
linear extension of $P$, and $w\mid_{Q}$ is a linear extension of $Q$, and we
have $w\in S\left(  w\mid_{P},w\mid_{Q}\right)  $ (since both $w\mid_{P}$ and
$w\mid_{Q}$ are subsequences of $w$, and their sizes add up to the size of
$w$). Therefore, if $w$ is a linear extension of $P\sqcup Q$, then
$w\in\bigcup_{\pi\in\mathcal{L}\left(  P\right)  ;\ \sigma\in\mathcal{L}%
\left(  Q\right)  }S\left(  \pi,\sigma\right)  $ (because $w\in S\left(
\pi,\sigma\right)  $ for
$\pi = \left. w\mid_{P}\right. \in\mathcal{L}\left(  P\right)  $ and
$\sigma = \left. w\mid_{Q} \right. \in\mathcal{L}\left(  Q\right)  $).

\item Conversely, any $w\in\bigcup_{\pi\in\mathcal{L}\left(  P\right)
;\ \sigma\in\mathcal{L}\left(  Q\right)  }S\left(  \pi,\sigma\right)  $ is a
linear extension of $P\sqcup Q$.
\end{itemize}

Combining these two facts, we conclude the following:

\begin{statement}
\textit{Observation 1:} The linear extensions of $P\sqcup Q$ are precisely the
elements of $\bigcup_{\pi\in\mathcal{L}\left(  P\right)  ;\ \sigma
\in\mathcal{L}\left(  Q\right)  }S\left(  \pi,\sigma\right)  $.
\end{statement}

Next, we notice the following:

\begin{statement}
\textit{Observation 2:} The union $\bigcup_{\pi\in\mathcal{L}\left(  P\right)
;\ \sigma\in\mathcal{L}\left(  Q\right)  }S\left(  \pi,\sigma\right)  $ is a
disjoint union (i.e., the sets $S\left(  \pi,\sigma\right)  $ for distinct
pairs $\left(  \pi,\sigma\right)  $ are disjoint).
\end{statement}

[\textit{Proof of Observation 2:} If we are given an element $w\in S\left(
\pi,\sigma\right)  $ for some $\pi\in\mathcal{L}\left(  P\right)  $ and
$\sigma\in\mathcal{L}\left(  Q\right)  $, then we can uniquely reconstruct
$\left(  \pi,\sigma\right)  $ from $w$ (namely, $\left(  \pi,\sigma\right)  $
is given by $\pi=w\mid_{P}$ and $\sigma=w\mid_{Q}$). Thus, the sets
$S\left(  \pi,\sigma\right)  $ for distinct pairs $\left(  \pi,\sigma\right)
$ are disjoint. This proves Observation 2.]

We will furthermore need a simple auxiliary claim:

\begin{statement}
\textit{Observation 3:} Let $R$ be a finite poset. Let $u$ and $v$ be two
elements of $R$. Assume that we don't have $u>v$ in $R$. Let $w$ be a linear
extension of $R$. Let $R^{\prime}$ be the poset obtained from $R$ by adding
the relation $u<v$ to $R$. Then, $w$ is a linear extension of $R^{\prime}$ if
and only if we have $u<v$ in $w$.
\end{statement}

[\textit{Proof of Observation 3:} The \textquotedblleft only
if\textquotedblright\ direction is obvious. The \textquotedblleft
if\textquotedblright\ direction is easily derived from the definition of a
linear extension, once you recall that any pair $\left(  x,y\right)  $ of
elements of $R^{\prime}$ that satisfies $x<y$ in $R^{\prime}$ must either
already satisfy $x<y$ in $R$, or satisfy $x\leq u$ and $v\leq y$ in $R$. The
details are left to the reader.]

\textbf{(a)} Observation 2 shows that the union $\bigcup_{\pi\in
\mathcal{L}\left(  P\right)  ;\ \sigma\in\mathcal{L}\left(  Q\right)
}S\left(  \pi,\sigma\right)  $ is a disjoint union. Hence, the union
$\bigcup_{\pi\in\mathcal{L}\left(  P\right)  ;\ \sigma\in\mathcal{L}\left(
Q\right)  }S_{\prec}\left(  \pi,\sigma\right)  $ is a disjoint union as well
(since $S_{\prec}\left(  \pi,\sigma\right)  $ is a subset of $S\left(
\pi,\sigma\right)  $ for all $\pi$ and $\sigma$).

We don't have $\min P>\min Q$ in $P\sqcup Q$. Thus, applying Observation 3 to
$R=P\sqcup Q$, $u=\min P$, $v=\min Q$ and $R^{\prime}=P\left.  \prec\right.
Q$, we obtain the following:

\begin{statement}
\textit{Observation 4:} Let $w$ be a linear extension of $P\sqcup Q$. Then,
$w$ is a linear extension of $P\left.  \prec\right.  Q$ if and only if we have
$\min P<\min Q$ in $w$.
\end{statement}

Next, we claim that
\begin{equation}
\mathcal{L}\left(  P\left.  \prec\right.  Q\right)  \subseteq\bigcup_{\pi
\in\mathcal{L}\left(  P\right)  ;\ \sigma\in\mathcal{L}\left(  Q\right)
}S_{\prec}\left(  \pi,\sigma\right)  .\label{pf.lem.dendri.4.1.lem2.a.5}%
\end{equation}

[\textit{Proof:} Let $w\in\mathcal{L}\left(  P\left.  \prec\right.  Q\right)
$. Thus, $w$ is a linear extension of $P\left.  \prec\right.  Q$. Hence, $w$
is a linear extension of $P\sqcup Q$, and we have $\min P<\min Q$ in $w$ (by
Observation 4). Since $w$ is a linear extension of $P\sqcup Q$, we have
$w\in\bigcup_{\pi\in\mathcal{L}\left(  P\right)  ;\ \sigma\in\mathcal{L}%
\left(  Q\right)  }S\left(  \pi,\sigma\right)  $ (by Observation 1). In other
words, $w\in S\left(  \pi,\sigma\right)  $ for some $\pi\in\mathcal{L}\left(
P\right)  $ and $\sigma\in\mathcal{L}\left(  Q\right)  $. Consider these $\pi$
and $\sigma$.

The element $\min P$ is the minimum element of the poset $P$, and thus must be
the first letter of $\pi$ (since $\pi$ is a linear extension of $P$).
Similarly, $\min Q$ must be the first letter of $\sigma$. But $w\in S\left(
\pi,\sigma\right)  $. Hence, the first letter of $w$ is either the first
letter of $\pi$ or the first letter of $\sigma$.

We have $\min P<\min Q$ in $w$. In other words, $\min P$ appears before $\min
Q$ in the list $w$. Hence, $\min Q$ cannot be the first letter of $w$. In
other words, the first letter of $w$ cannot be $\min Q$. In other words, the
first letter of $w$ cannot be the first letter of $\sigma$ (since $\min Q$ is
the first letter of $\sigma$). Thus, the first letter of $w$ is the first
letter of $\pi$ (since the first letter of $w$ is either the first letter of
$\pi$ or the first letter of $\sigma$). In other words, $w$ is a left shuffle
of $\pi$ and $\sigma$ (since $w\in S\left(  \pi,\sigma\right)  $). In other
words, $w\in S_{\prec}\left(  \pi,\sigma\right)  $ (since $S_{\prec}\left(
\pi,\sigma\right)  $ is the set of all left shuffles of $\pi$ and $\sigma$).

Now, forget that we have defined $\pi$ and $\sigma$. We thus have shown that
$w\in S_{\prec}\left(  \pi,\sigma\right)  $ for some $\pi\in\mathcal{L}\left(
P\right)  $ and $\sigma\in\mathcal{L}\left(  Q\right)  $. Thus, $w\in
\bigcup_{\pi\in\mathcal{L}\left(  P\right)  ;\ \sigma\in\mathcal{L}\left(
Q\right)  }S_{\prec}\left(  \pi,\sigma\right)  $. Since we have proven this
for any $w\in\mathcal{L}\left(  P\left.  \prec\right.  Q\right)  $, we thus
have proven (\ref{pf.lem.dendri.4.1.lem2.a.5}).]

On the other hand, we claim that
\begin{equation}
\bigcup_{\pi\in\mathcal{L}\left(  P\right)  ;\ \sigma\in\mathcal{L}\left(
Q\right)  }S_{\prec}\left(  \pi,\sigma\right)  \subseteq\mathcal{L}\left(
P\left.  \prec\right.  Q\right)  .\label{pf.lem.dendri.4.1.lem2.a.6}%
\end{equation}

[\textit{Proof:} Let $w\in\bigcup_{\pi\in\mathcal{L}\left(  P\right)
;\ \sigma\in\mathcal{L}\left(  Q\right)  }S_{\prec}\left(  \pi,\sigma\right)
$.

We have%
\[
w\in\bigcup_{\pi\in\mathcal{L}\left(  P\right)  ;\ \sigma\in\mathcal{L}\left(
Q\right)  }\underbrace{S_{\prec}\left(  \pi,\sigma\right)  }_{\subseteq
S\left(  \pi,\sigma\right)  }\subseteq\bigcup_{\pi\in\mathcal{L}\left(
P\right)  ;\ \sigma\in\mathcal{L}\left(  Q\right)  }S\left(  \pi
,\sigma\right)  .
\]
Thus, $w$ is a linear extension of $P\sqcup Q$ (by Observation 1).

Also, $w\in\bigcup_{\pi\in\mathcal{L}\left(  P\right)  ;\ \sigma
\in\mathcal{L}\left(  Q\right)  }S_{\prec}\left(  \pi,\sigma\right)  $.
Thus, $w\in S_{\prec}\left(  \pi,\sigma\right)  $ for some $\pi\in
\mathcal{L}\left(  P\right)  $ and $\sigma\in\mathcal{L}\left(  Q\right)  $.
Consider these $\pi$ and $\sigma$. From $w\in S_{\prec}\left(  \pi
,\sigma\right)  $, we conclude that $w$ is a left shuffle of $\pi$ and
$\sigma$. In other words, $w$ is a shuffle of $\pi$ and $\sigma$ such that the
first letter of $w$ is the first letter of $\pi$.

The element $\min P$ is the minimum element of the poset $P$, and thus must be
the first letter of $\pi$ (since $\pi$ is a linear extension of $P$). In other
words, $\min P$ is the first letter of $w$ (since the first letter of $w$ is
the first letter of $\pi$). Hence, $\min P$ appears before $\min Q$ in the
list $w$ (since $\min P$ and $\min Q$ are two distinct letters of $w$). In
other words, we have $\min P<\min Q$ in $w$. Hence, Observation 4 yields that
$w$ is a linear extension of $P\left.  \prec\right.  Q$ (since $w$ is a linear
extension of $P\sqcup Q$). In other words, $w\in\mathcal{L}\left(  P\left.
\prec\right.  Q\right)  $. Since we have proven this for any $w\in\bigcup
_{\pi\in\mathcal{L}\left(  P\right)  ;\ \sigma\in\mathcal{L}\left(  Q\right)
}S_{\prec}\left(  \pi,\sigma\right)  $, we thus have proven
(\ref{pf.lem.dendri.4.1.lem2.a.6}).]

Combining (\ref{pf.lem.dendri.4.1.lem2.a.5}) and
(\ref{pf.lem.dendri.4.1.lem2.a.6}), we obtain
\[
\mathcal{L}\left(  P\left.  \prec\right.  Q\right)  =\bigcup_{\pi
\in\mathcal{L}\left(  P\right)  ;\ \sigma\in\mathcal{L}\left(  Q\right)
}S_{\prec}\left(  \pi,\sigma\right)  =\bigsqcup_{\pi\in\mathcal{L}\left(
P\right)  ;\ \sigma\in\mathcal{L}\left(  Q\right)  }S_{\prec}\left(
\pi,\sigma\right)
\]
(since the union $\bigcup_{\pi\in\mathcal{L}\left(  P\right)  ;\ \sigma
\in\mathcal{L}\left(  Q\right)  }S_{\prec}\left(  \pi,\sigma\right)  $ is a
disjoint union). This proves Lemma \ref{lem.dendri.4.1.lem2} \textbf{(a)}.

\textbf{(b)} Lemma \ref{lem.dendri.4.1.lem2} \textbf{(b)} is proven similarly
to Lemma \ref{lem.dendri.4.1.lem2} \textbf{(a)}.
\end{proof}

We can now prove Theorem \ref{thm.dendri.4.1}.
First, let us rewrite Theorem \ref{thm.dendri.4.1} as follows:\footnote{Recall
that for any permutation $\varphi$, we have let $\operatorname*{Comp}\varphi$
denote the descent composition of $\varphi$.}

\begin{corollary}
\label{cor.dendri.4.1}Let $\pi$ and $\sigma$ be two disjoint nonempty
permutations. Assume that $\pi_{1}>\sigma_{1}$. Then,%
\[
F_{\operatorname*{Comp}\pi}\left.  \prec\right.  F_{\operatorname*{Comp}%
\sigma}=\sum_{\chi\in S_{\prec}\left(  \pi,\sigma\right)  }%
F_{\operatorname*{Comp}\chi}%
\]
and%
\[
F_{\operatorname*{Comp}\pi}\left.  \succeq\right.  F_{\operatorname*{Comp}%
\sigma}=\sum_{\chi\in S_{\succ}\left(  \pi,\sigma\right)  }%
F_{\operatorname*{Comp}\chi}.
\]

\end{corollary}

\begin{proof}
[Proof of Corollary \ref{cor.dendri.4.1} (sketched).] Let $n=\left\vert
\pi\right\vert $ and $m=\left\vert \sigma\right\vert $. We shall use the
notations of \cite[Section 5.2]{HopfComb}; in particular,
\textquotedblleft labelled poset\textquotedblright\ will be defined as in
\cite[Definition 5.2.1]{HopfComb}.

If $R$ is a labelled poset, and if $w$ is a linear extension of $R$,
then $w$ can be regarded as a labelled poset itself, but also as a
permutation (since $w$ is a list of distinct elements of $R$, and thus
a word over the alphabet $\mathbb{P}$ with no two equal letters).
The first interpretation (as a labelled poset)
gives rise to a quasisymmetric function
$F_w\left(\xx\right)$ (defined as in \cite[Definition 5.2.1]{HopfComb}).
The second interpretation (as a permutation) leads to a composition
$\operatorname{Comp} w$. These two objects are connected by the equality
\begin{align}
F_w \left( \xx \right) = F_{\operatorname{Comp} w} .
\label{pf.cor.dendri.4.1.Fw}
\end{align}
(This follows from \cite[Proposition 5.2.10]{HopfComb}; but keep in mind that
\cite[Proposition 5.2.10]{HopfComb} denotes $F_{\operatorname*{Comp}\pi}$
by $L_{\alpha}$ in this context.)

Let $P$ be the labelled poset whose elements are $\pi_{1},\pi_{2},\ldots
,\pi_{n}$ and whose order is the total order given by $\pi_{1}<\pi_{2}%
<\cdots<\pi_{n}$. Thus, $P$ is totally ordered, and its minimum element is
$\min P=\pi_{1}$. Also, \cite[Proposition 5.2.10]{HopfComb} yields
$F_{P}\left(  \mathbf{x}\right)  =F_{\operatorname*{Comp}\pi}$. (Keep in mind
that \cite[Proposition 5.2.10]{HopfComb} denotes $F_{\operatorname*{Comp}\pi}$
by $L_{\alpha}$ in this context.)

Let $Q$ be the labelled poset whose elements are $\sigma_{1},\sigma_{2}%
,\ldots,\sigma_{m}$ and whose order is the total order given by $\sigma
_{1}<\sigma_{2}<\cdots<\sigma_{m}$. Thus, $Q$ is totally ordered, and its
minimum element is $\min Q=\sigma_{1}$. Also, \cite[Proposition 5.2.10]%
{HopfComb} yields $F_{Q}\left(  \mathbf{x}\right)  =F_{\operatorname*{Comp}%
\sigma}$.

The posets $P$ and $Q$ are disjoint (since the permutations $\pi$ and $\sigma$
are disjoint). Define three labelled posets $P\sqcup Q$, $P\left.
\prec\right.  Q$ and $P\left.  \succeq\right.  Q$ as in Proposition
\ref{prop.dendri.4.1.lem}.

Lemma \ref{lem.dendri.4.1.lem2} \textbf{(a)} yields $\mathcal{L}\left(
P\left.  \prec\right.  Q\right)  =\bigsqcup_{\pi^{\prime}\in\mathcal{L}\left(
P\right)  ;\ \sigma^{\prime}\in\mathcal{L}\left(  Q\right)  }S_{\prec}\left(
\pi^{\prime},\sigma^{\prime}\right)  $ (where we are using the letters
$\pi^{\prime}$ and $\sigma^{\prime}$ for our subscripts, since the letters
$\pi$ and $\sigma$ are already taken). But $P$ is totally ordered; thus, there
exists only one linear extension $\pi^{\prime}\in\mathcal{L}\left(  P\right)
$, namely, $\pi^{\prime}=\pi$. In other words, $\mathcal{L}\left(  P\right)
=\left\{  \pi\right\}  $. Similarly, $\mathcal{L}\left(  Q\right)  =\left\{
\sigma\right\}  $. Hence,%
\begin{align}
\mathcal{L}\left(  P\left.  \prec\right.  Q\right)    & =\bigsqcup
_{\pi^{\prime}\in\mathcal{L}\left(  P\right)  ;\ \sigma^{\prime}\in
\mathcal{L}\left(  Q\right)  }S_{\prec}\left(  \pi^{\prime},\sigma^{\prime
}\right) \nonumber\\
&  =\bigsqcup_{\pi^{\prime}\in\left\{  \pi\right\}  ;\ \sigma^{\prime
}\in\left\{  \sigma\right\}  }S_{\prec}\left(  \pi^{\prime},\sigma^{\prime
}\right)  \ \ \ \ \ \ \ \ \ \ \left(
\begin{array}
[c]{c}%
\text{since }\mathcal{L}\left(  P\right)  =\left\{  \pi\right\}  \\
\text{and }\mathcal{L}\left(  Q\right)  =\left\{  \sigma\right\}
\end{array}
\right)  \nonumber\\
& =S_{\prec}\left(  \pi,\sigma\right)  .\label{pf.cor.dendri.4.1.L1}%
\end{align}

Recall that every labelled poset $R$ satisfies%
\begin{equation}
F_{R}\left(  \mathbf{x}\right)  =\sum_{w\in\mathcal{L}\left(  R\right)  }%
F_{w}\left(  \mathbf{x}\right)  \label{pf.cor.dendri.4.1.FR}%
\end{equation}
(by \cite[Theorem 5.2.11]{HopfComb}).

Now, $\pi_{1}>\sigma_{1}$ in $\mathbb{Z}$. In other words, $\pi_{1}%
>_{\ZZ}\sigma_{1}$. In other words, $\min P>_{\ZZ}\min Q$ (since
$\min P=\pi_{1}$ and $\min Q=\sigma_{1}$). Thus,
\begin{align*}
& \underbrace{F_{\operatorname*{Comp}\pi}}_{=F_{P}\left(  \mathbf{x}\right)
}\left.  \prec\right.  \underbrace{F_{\operatorname*{Comp}\sigma}}%
_{=F_{Q}\left(  \mathbf{x}\right)  }\\
& =F_{P}\left(  \mathbf{x}\right)  \left.  \prec\right.  F_{Q}\left(
\mathbf{x}\right)  =F_{P\left.  \prec\right.  Q}\left(  \mathbf{x}\right)
\ \ \ \ \ \ \ \ \ \ \left(  \text{by Proposition \ref{prop.dendri.4.1.lem}
\textbf{(a)}}\right)  \\
& =\sum_{w\in\mathcal{L}\left(  P\left.  \prec\right.  Q\right)  }F_{w}\left(
\mathbf{x}\right)  \ \ \ \ \ \ \ \ \ \ \left(  \text{by
(\ref{pf.cor.dendri.4.1.FR})}\right)  \\
& =\sum_{w\in S_{\prec}\left(  \pi,\sigma\right)  }\underbrace{F_{w}\left(
\mathbf{x}\right)  }_{\substack{=F_{\operatorname*{Comp}w}\\\text{(by
\eqref{pf.cor.dendri.4.1.Fw})}}}\ \ \ \ \ \ \ \ \ \ \left(  \text{by
(\ref{pf.cor.dendri.4.1.L1})}\right)  \\
& =\sum_{w\in S_{\prec}\left(  \pi,\sigma\right)  }F_{\operatorname*{Comp}%
w}=\sum_{\chi\in S_{\prec}\left(  \pi,\sigma\right)  }F_{\operatorname*{Comp}%
\chi}.
\end{align*}
A similar argument (using Proposition \ref{prop.dendri.4.1.lem} \textbf{(b)}
and Lemma \ref{lem.dendri.4.1.lem2} \textbf{(b)} instead of Proposition
\ref{prop.dendri.4.1.lem} \textbf{(a)} and Lemma \ref{lem.dendri.4.1.lem2}
\textbf{(a)}) shows that%
\[
F_{\operatorname*{Comp}\pi}\left.  \succeq\right.  F_{\operatorname*{Comp}%
\sigma}=\sum_{\chi\in S_{\succ}\left(  \pi,\sigma\right)  }%
F_{\operatorname*{Comp}\chi}.
\]
Thus, Corollary \ref{cor.dendri.4.1} is proven.
\end{proof}

Hence, Theorem \ref{thm.dendri.4.1} is proven as well (since Corollary
\ref{cor.dendri.4.1} was just a restatement of Theorem \ref{thm.dendri.4.1}).

\subsection{Ideals}

\begin{definition}
Let $A$ be a $\mathbf{k}$-module equipped with some binary operation $\ast$
(written infix).

\textbf{(a)} If $B$ and $C$ are two $\mathbf{k}$-submodules of $A$, then
$B\ast C$ shall mean the $\mathbf{k}$-submodule of $A$ spanned by all elements
of the form $b\ast c$ with $b\in B$ and $c\in C$.

\textbf{(b)} A $\mathbf{k}$-submodule $M$ of $A$ is said to be a \textit{left
}$\ast$\textit{-ideal} if and only if it satisfies $A\ast M\subseteq M$.

\textbf{(c)} A $\mathbf{k}$-submodule $M$ of $A$ is said to be a \textit{right
}$\ast$\textit{-ideal} if and only if it satisfies $M\ast A\subseteq M$.

\textbf{(d)} A $\mathbf{k}$-submodule $M$ of $A$ is said to be a\textit{
}$\ast$\textit{-ideal} if and only if it is both a left $\ast$-ideal and a
right $\ast$-ideal.
\end{definition}

\begin{theorem}
\label{thm.ideal-crit2}Let $M$ be an ideal of $\operatorname*{QSym}$. Let
$A=\operatorname*{QSym}$.
Let $A^+$ denote the $\mathbf{k}$-submodule of $A$ consisting of those
quasisymmetric functions whose constant term is $0$.

\textbf{(a)} If $A\bel  M\subseteq M$ and $M \subseteq A^+$,
then $M\left.  \prec\right.  A\subseteq M$.

\textbf{(b)} If $A\tvi  M\subseteq M$, then $A\left.  \succeq\right.
M\subseteq M$.

\textbf{(c)} If $A\tvi  M\subseteq M$ and $A\bel  M\subseteq M$
and $M \subseteq A^+$, then $M$ is a
$\left.  \prec\right.  $-ideal and a $\left.  \succeq\right.  $-ideal of
$\operatorname*{QSym}$.
\end{theorem}

\begin{proof}
[Proof of Theorem \ref{thm.ideal-crit2}.]\textbf{(a)} Assume that
$A\bel  M\subseteq M$. If $a\in M$ and $b\in A$, then the constant term
of $a$ is $0$ (since $a \in M \subseteq A^+$), and we have
\begin{align*}
a\left.  \prec\right.  b  &  =\sum_{\left(  b\right)  }\left(
\underbrace{S\left(  b_{\left(  1\right)  }\right)  }_{\in A}%
\bel \underbrace{a}_{\in M}\right)  \underbrace{b_{\left(  2\right)  }}_{\in
A}\ \ \ \ \ \ \ \ \ \ \left(  \text{by Theorem \ref{thm.beldend}}\right) \\
&  \in\underbrace{\left(  A\bel M\right)  }_{\subseteq M}A\subseteq
MA\subseteq M\ \ \ \ \ \ \ \ \ \ \left(  \text{since }M\text{ is an ideal of
}A\right)  .
\end{align*}
Thus, $M\left.  \prec\right.  A\subseteq M$. This proves Theorem
\ref{thm.ideal-crit2} \textbf{(a)}.

\textbf{(b)} Assume that $A\tvi  M\subseteq M$. If $a\in M$ and $b\in A$, then%
\begin{align*}
b\left.  \succeq\right.  a  &  =\sum_{\left(  b\right)  }\left(
\underbrace{S\left(  b_{\left(  1\right)  }\right)  }_{\in A}%
\tvi \underbrace{a}_{\in M}\right)  \underbrace{b_{\left(  2\right)  }}_{\in
A}\ \ \ \ \ \ \ \ \ \ \left(  \text{by Theorem \ref{thm.tvidend'}}\right) \\
&  \in\underbrace{\left(  A\tvi M\right)  }_{\subseteq M}A\subseteq
MA\subseteq M\ \ \ \ \ \ \ \ \ \ \left(  \text{since }M\text{ is an ideal of
}A\right)  .
\end{align*}
Thus, $A\left.  \succeq\right.  M\subseteq M$. This proves Theorem
\ref{thm.ideal-crit2} \textbf{(b)}.

\textbf{(c)} Assume that $A\tvi  M\subseteq M$ and $A\bel  M\subseteq M$.
Then, Theorem \ref{thm.ideal-crit2} \textbf{(b)} yields $A\left.
\succeq\right.  M\subseteq M$. Thus, $M$ is a left $\left.  \succeq\right.  $-ideal.

Now, any $b\in M$ and $a\in A$ satisfy%
\begin{align*}
a\left.  \prec\right.  b  &  =\underbrace{a}_{\in A}\underbrace{b}_{\in
M}-\underbrace{a}_{\in A}\left.  \succeq\right.  \underbrace{b}_{\in
M}\ \ \ \ \ \ \ \ \ \ \left(  \text{by (\ref{eq.dendriform.1})}\right) \\
&  \in\underbrace{AM}_{\substack{\subseteq M\\\text{(since }M\text{ is an
ideal of }A\text{)}}}-\underbrace{A\left.  \succeq\right.  M}_{\subseteq
M}\subseteq M-M\subseteq M.
\end{align*}
In other words, $A\left.  \prec\right.  M\subseteq M$. In other words, $M$ is
a left $\left.  \prec\right.  $-ideal.

But Theorem \ref{thm.ideal-crit2} \textbf{(a)} yields $M\left.  \prec\right.
A\subseteq M$. In other words, $M$ is a right $\left.  \prec\right.  $-ideal.

Any $a\in M$ and $b\in A$ satisfy%
\begin{align*}
a\left.  \succeq\right.  b  &  =\underbrace{a}_{\in M}\underbrace{b}_{\in
A}-\underbrace{a}_{\in M}\left.  \prec\right.  \underbrace{b}_{\in
A}\ \ \ \ \ \ \ \ \ \ \left(  \text{by (\ref{eq.dendriform.1})}\right) \\
&  \in\underbrace{MA}_{\substack{\subseteq M\\\text{(since }M\text{ is an
ideal of }A\text{)}}}-\underbrace{M\left.  \prec\right.  A}_{\subseteq
M}\subseteq M-M\subseteq M.
\end{align*}
In other words, $M\left.  \succeq\right.  A\subseteq M$. In other words, $M$
is a right $\left.  \succeq\right.  $-ideal.

Hence, $M$ is a $\left.  \prec\right.  $-ideal (since $M$ is a left $\left.
\prec\right.  $-ideal and a right $\left.  \prec\right.  $-ideal) and a
$\left.  \succeq\right.  $-ideal (since $M$ is a left $\left.  \succeq\right.
$-ideal and a right $\left.  \succeq\right.  $-ideal). This proves Theorem
\ref{thm.ideal-crit2} \textbf{(c)}.
\end{proof}

Another simple fact is the following:

\begin{proposition}
\label{prop.ideal-crit3}Let $M$ be simultaneously a $\left.  \prec\right.
$-ideal and a $\left.  \succeq\right.  $-ideal of $\operatorname*{QSym}$.
Then, $M$ is an ideal of $\operatorname*{QSym}$.
\end{proposition}

\begin{proof}
[Proof of Proposition \ref{prop.ideal-crit3}.]Any $a\in M$ and $b\in
\operatorname*{QSym}$ satisfy%
\begin{align*}
ab  &  =\underbrace{a}_{\in M}\left.  \prec\right.  \underbrace{b}%
_{\in\operatorname*{QSym}}+\underbrace{a}_{\in M}\left.  \succeq\right.
\underbrace{b}_{\in\operatorname*{QSym}}\ \ \ \ \ \ \ \ \ \ \left(  \text{by
(\ref{eq.dendriform.1})}\right) \\
&  \in\underbrace{M\left.  \prec\right.  \operatorname*{QSym}}%
_{\substack{\subseteq M\\\text{(since }M\text{ is a }\left.  \prec\right.
\text{-ideal of }\operatorname*{QSym}\text{)}}}+\underbrace{M\left.
\succeq\right.  \operatorname*{QSym}}_{\substack{\subseteq M\\\text{(since
}M\text{ is a }\left.  \succeq\right.  \text{-ideal of }\operatorname*{QSym}%
\text{)}}}\subseteq M+M\subseteq M.
\end{align*}
In other words, $M$ is an ideal of $\operatorname*{QSym}$
(since $\operatorname*{QSym}$ is commutative). This proves
Proposition \ref{prop.ideal-crit3}.
\end{proof}

\begin{question}
Proposition \ref{prop.ideal-crit3} says that if a $\mathbb{Q}$-vector subspace
$M$ of $\operatorname*{QSym}$ is simultaneously a $\left.  \prec\right.
$-ideal and a $\left.  \succeq\right.  $-ideal, then it is also an ideal.
Similarly, if $M$ is an ideal and a $\left.  \prec\right.  $-ideal, then it is
a $\left.  \succeq\right.  $-ideal. Can we state any other such criteria?
\end{question}

\subsection{Application to $\mathcal{K}_{\operatorname{Epk}}$}

We now claim the following:

\begin{theorem}
\label{thm.Epk.dend}The ideal $\mathcal{K}_{\operatorname{Epk}}$ of
$\operatorname*{QSym}$ is a $\tvi  $-ideal, a $\bel  $-ideal, a $\left.
\prec\right.  $-ideal and a $\left.  \succeq\right.  $-ideal of
$\operatorname*{QSym}$.
\end{theorem}

\begin{proof}
[Proof of Theorem \ref{thm.Epk.dend}.]Let $A=\operatorname*{QSym}$. Corollary
\ref{cor.Epk.ideal} shows that $\mathcal{K}_{\operatorname{Epk}}$ is an ideal
of $\operatorname*{QSym}$.

Let us recall the binary relation $\rightarrow$ on the set of compositions
defined in Proposition \ref{prop.K.Epk.F}.

\begin{statement}
\textit{Claim 1:} Let $J$ and $K$ be two compositions satisfying $J\rightarrow
K$. Let $G$ be a further composition. Then, $\left[  G,J\right]
\rightarrow\left[  G,K\right]  $.
\end{statement}

[\textit{Proof of Claim 1:} Write the composition $J$ in the form $J=\left(
j_{1},j_{2},\ldots,j_{m}\right)  $. Write the composition $G$ in the form
$G=\left(  g_{1},g_{2},\ldots,g_{p}\right)  $.

We have $J\rightarrow K$. In other words, there exists an $\ell\in\left\{
2,3,\ldots,m\right\}  $ such that $j_{\ell}>2$ and $K=\left(  j_{1}%
,j_{2},\ldots,j_{\ell-1},1,j_{\ell}-1,j_{\ell+1},j_{\ell+2},\ldots
,j_{m}\right)  $ (by the definition of the relation $\rightarrow$). Consider
this $\ell$. Clearly, $\ell>1$ (since $\ell\in\left\{  2,3,\ldots,m\right\}
$), so that $p+\ell>\underbrace{p}_{\geq0}+1\geq1$.

From $G=\left(  g_{1},g_{2},\ldots,g_{p}\right)  $ and $J=\left(  j_{1}%
,j_{2},\ldots,j_{m}\right)  $, we obtain%
\begin{equation}
\left[  G,J\right]  =\left(  g_{1},g_{2},\ldots,g_{p},j_{1},j_{2},\ldots
,j_{m}\right)  . \label{pf.thm.Epk.dend.c1.pf.GJ=}%
\end{equation}
From $G=\left(  g_{1},g_{2},\ldots,g_{p}\right)  $ and $K=\left(  j_{1}%
,j_{2},\ldots,j_{\ell-1},1,j_{\ell}-1,j_{\ell+1},j_{\ell+2},\ldots
,j_{m}\right)  $, we obtain%
\begin{equation}
\left[  G,K\right]  =\left(  g_{1},g_{2},\ldots,g_{p},j_{1},j_{2}%
,\ldots,j_{\ell-1},1,j_{\ell}-1,j_{\ell+1},j_{\ell+2},\ldots,j_{m}\right)  .
\label{pf.thm.Epk.dend.c1.pf.GK=}%
\end{equation}
From looking at (\ref{pf.thm.Epk.dend.c1.pf.GJ=}) and
(\ref{pf.thm.Epk.dend.c1.pf.GK=}), we conclude immediately that the
composition $\left[  G,K\right]  $ is obtained from $\left[  G,J\right]  $ by
\textquotedblleft splitting\textquotedblright\ the entry $j_{\ell}>2$ into two
consecutive entries $1$ and $j_{\ell}-1$, and that this entry $j_{\ell}$ was
not the first entry (indeed, this entry is the $\left(  p+\ell\right)  $-th
entry, but $p+\ell>1$). Hence, $\left[  G,J\right]  \rightarrow\left[
G,K\right]  $ (by the definition of the relation $\rightarrow$). This proves
Claim 1.]

\begin{statement}
\textit{Claim 2:} We have $A\tvi  \mathcal{K}_{\operatorname{Epk}}%
\subseteq\mathcal{K}_{\operatorname{Epk}}$.
\end{statement}

[\textit{Proof of Claim 2:} We must show that $a\tvi  m\in\mathcal{K}%
_{\operatorname{Epk}}$ for every $a\in A$ and $m\in\mathcal{K}%
_{\operatorname{Epk}}$. So let us fix $a\in A$ and $m\in\mathcal{K}%
_{\operatorname{Epk}}$.

Proposition \ref{prop.K.Epk.F} shows that the $\mathbb{Q}$-vector space
$\mathcal{K}_{\operatorname{Epk}}$ is spanned by all differences of the form
$F_{J}-F_{K}$, where $J$ and $K$ are two compositions satisfying $J\rightarrow
K$. Hence, we can WLOG assume that $m$ is such a difference (because the
relation $a\tvi  m\in\mathcal{K}_{\operatorname{Epk}}$, which we must prove,
is $\mathbb{Q}$-linear in $m$). Assume this. Thus, $m=F_{J}-F_{K}$ for some
two compositions $J$ and $K$ satisfying $J\rightarrow K$. Consider these $J$
and $K$.

From $J\rightarrow K$, we easily conclude that the composition $J$ is
nonempty. Thus, $\left\vert J\right\vert \neq0$. But from $J\rightarrow K$, we
also obtain $\left\vert J\right\vert =\left\vert K\right\vert $. Hence,
$\left\vert K\right\vert =\left\vert J\right\vert \neq0$. Thus, the
composition $K$ is nonempty.

Recall that the family $\left(  F_{L}\right)  _{L\text{ is a composition}}$ is
a basis of the $\mathbb{Q}$-vector space $\operatorname*{QSym}=A$. Hence, we
can WLOG assume that $a$ belongs to this family (since the relation
$a\tvi  m\in\mathcal{K}_{\operatorname{Epk}}$, which we must prove, is
$\mathbb{Q}$-linear in $a$). Assume this. Thus, $a=F_{G}$ for some composition
$G$. Consider this $G$.

If $G$ is the empty composition, then $a=F_{G}=1$, and therefore
$\underbrace{a}_{=1}\tvi  m=1\tvi  m=m\in\mathcal{K}_{\operatorname{Epk}}$
holds. Thus, for the rest of this proof, we WLOG assume that $G$ is not the
empty composition. Thus, $G$ is nonempty.

Recall that for any two compositions $\alpha$ and $\beta$, we have $F_{\alpha
}\tvi  F_{\beta}=F_{\left[  \alpha,\beta\right]  }$. Applying this to
$\alpha=G$ and $\beta=J$, we obtain $F_{G}\tvi  F_{J}=F_{\left[  G,J\right]
}$. Similarly, $F_{G}\tvi  F_{K}=F_{\left[  G,K\right]  }$.

But Claim 1 yields $\left[  G,J\right]  \rightarrow\left[  G,K\right]  $.
Hence, the difference $F_{\left[  G,J\right]  }-F_{\left[  G,K\right]  }$ is
one of the differences which span the ideal $\mathcal{K}_{\operatorname{Epk}%
}$ according to Proposition \ref{prop.K.Epk.F}. Thus, in particular, this
difference lies in $\mathcal{K}_{\operatorname{Epk}}$. In other words,
$F_{\left[  G,J\right]  }-F_{\left[  G,K\right]  }\in\mathcal{K}%
_{\operatorname{Epk}}$.

Now,%
\begin{align*}
\underbrace{a}_{=F_{G}}\tvi \underbrace{m}_{=F_{J}-F_{K}}  &  =F_{G}%
\tvi \left(  F_{J}-F_{K}\right)  =\underbrace{F_{G}\tvi F_{J}}_{=F_{\left[
G,J\right]  }}-\underbrace{F_{G}\tvi F_{K}}_{=F_{\left[  G,K\right]  }}\\
&  =F_{\left[  G,J\right]  }-F_{\left[  G,K\right]  }\in\mathcal{K}%
_{\operatorname{Epk}}.
\end{align*}
This proves Claim 2.]

\begin{statement}
\textit{Claim 3:} Let $J$ and $K$ be two compositions satisfying $J\rightarrow
K$. Let $G$ be a further composition. Then, $\left[  J,G\right]
\rightarrow\left[  K,G\right]  $.
\end{statement}

[\textit{Proof of Claim 3:} This is proven in the same way as we proved Claim
1, with the only difference that $j_{\ell}$ is now the $\ell$-th entry of
$\left[  J,G\right]  $ and not the $\left(  p+\ell\right)  $-th entry (but
this is still sufficient, since $\ell>1$).]

\begin{statement}
\textit{Claim 4:} We have $\mathcal{K}_{\operatorname{Epk}}\tvi  A\subseteq
\mathcal{K}_{\operatorname{Epk}}$.
\end{statement}

[\textit{Proof of Claim 4:} This is proven in the same way as we proved Claim
2, with the only difference that now we need to use Claim 3 instead of Claim 1.]

Combining Claim 2 and Claim 4, we conclude that $\mathcal{K}%
_{\operatorname{Epk}}$ is a $\tvi$-ideal of $A=\operatorname*{QSym}$.

\begin{statement}
\textit{Claim 5:} Let $J$ and $K$ be two nonempty compositions satisfying
$J\rightarrow K$. Let $G$ be a further nonempty composition. Then, $G\odot
J\rightarrow G\odot K$.
\end{statement}

[\textit{Proof of Claim 5:} Write the composition $J$ in the form $J=\left(
j_{1},j_{2},\ldots,j_{m}\right)  $. Write the composition $G$ in the form
$G=\left(  g_{1},g_{2},\ldots,g_{p}\right)  $. Thus, $p>0$ (since the
composition $G$ is nonempty).

We have $J\rightarrow K$. In other words, there exists an $\ell\in\left\{
2,3,\ldots,m\right\}  $ such that $j_{\ell}>2$ and $K=\left(  j_{1}%
,j_{2},\ldots,j_{\ell-1},1,j_{\ell}-1,j_{\ell+1},j_{\ell+2},\ldots
,j_{m}\right)  $ (by the definition of the relation $\rightarrow$). Consider
this $\ell$. Clearly, $\ell\geq2$ (since $\ell\in\left\{  2,3,\ldots
,m\right\}  $), so that $\underbrace{p}_{>0}+\underbrace{\ell}_{\geq
2}-1>0+2-1=1$.

From $G=\left(  g_{1},g_{2},\ldots,g_{p}\right)  $ and $J=\left(  j_{1}%
,j_{2},\ldots,j_{m}\right)  $, we obtain%
\begin{equation}
G\odot J=\left(  g_{1},g_{2},\ldots,g_{p-1},g_{p}+j_{1},j_{2},j_{3}%
,\ldots,j_{m}\right)  . \label{pf.thm.Epk.dend.c5.pf.GJ=}%
\end{equation}
From $G=\left(  g_{1},g_{2},\ldots,g_{p}\right)  $ and $K=\left(  j_{1}%
,j_{2},\ldots,j_{\ell-1},1,j_{\ell}-1,j_{\ell+1},j_{\ell+2},\ldots
,j_{m}\right)  $, we obtain%
\begin{equation}
G\odot K=\left(  g_{1},g_{2},\ldots,g_{p-1},g_{p}+j_{1},j_{2},j_{3}%
,\ldots,j_{\ell-1},1,j_{\ell}-1,j_{\ell+1},j_{\ell+2},\ldots,j_{m}\right)
\label{pf.thm.Epk.dend.c5.pf.GK=}%
\end{equation}
(notice that the $g_{p}+j_{1}$ term is \textbf{not} a $g_{p}+1$ term, because
$\ell\geq2$).

From looking at (\ref{pf.thm.Epk.dend.c5.pf.GJ=}) and
(\ref{pf.thm.Epk.dend.c5.pf.GK=}), we conclude immediately that the
composition $G\odot K$ is obtained from $G\odot J$ by \textquotedblleft
splitting\textquotedblright\ the entry $j_{\ell}>2$ into two consecutive
entries $1$ and $j_{\ell}-1$, and that this entry $j_{\ell}$ was not the first
entry (indeed, this entry is the $\left(  p+\ell-1\right)  $-th entry, but
$p+\ell-1>1$). Hence, $G\odot J\rightarrow G\odot K$ (by the definition of the
relation $\rightarrow$). This proves Claim 5.]

\begin{statement}
\textit{Claim 6:} We have $A\bel\mathcal{K}_{\operatorname{Epk}}%
\subseteq\mathcal{K}_{\operatorname{Epk}}$.
\end{statement}

[\textit{Proof of Claim 6:} This is proven in the same way as we proved Claim
2, with the only difference that now we need to use Claim 5 instead of Claim 1
and that we need to use the formula $F_{\alpha}\bel  F_{\beta}=F_{\alpha
\odot\beta}$ instead of $F_{\alpha}\tvi  F_{\beta}=F_{\left[  \alpha
,\beta\right]  }$.]

\begin{statement}
\textit{Claim 7:} Let $J$ and $K$ be two nonempty compositions satisfying
$J\rightarrow K$. Let $G$ be a further nonempty composition. Then, $J\odot
G\rightarrow K\odot G$.
\end{statement}

[\textit{Proof of Claim 7:} Write the composition $J$ in the form $J=\left(
j_{1},j_{2},\ldots,j_{m}\right)  $. Write the composition $G$ in the form
$G=\left(  g_{1},g_{2},\ldots,g_{p}\right)  $. Thus, $p>0$ (since the
composition $G$ is nonempty).

We have $J\rightarrow K$. In other words, there exists an $\ell\in\left\{
2,3,\ldots,m\right\}  $ such that $j_{\ell}>2$ and $K=\left(  j_{1}%
,j_{2},\ldots,j_{\ell-1},1,j_{\ell}-1,j_{\ell+1},j_{\ell+2},\ldots
,j_{m}\right)  $ (by the definition of the relation $\rightarrow$). Consider
this $\ell$. Clearly, $\ell\geq2$ (since $\ell\in\left\{  2,3,\ldots
,m\right\}  $), so that $\ell>1$.

From $G=\left(  g_{1},g_{2},\ldots,g_{p}\right)  $ and $J=\left(  j_{1}%
,j_{2},\ldots,j_{m}\right)  $, we obtain%
\begin{equation}
J\odot G=\left(  j_{1},j_{2},\ldots,j_{m-1},j_{m}+g_{1},g_{2},g_{3}%
,\ldots,g_{p}\right)  . \label{pf.thm.Epk.dend.c7.pf.JG=}%
\end{equation}
Now, we distinguish between the following two cases:

\textit{Case 1:} We have $\ell=m$.

\textit{Case 2:} We have $\ell\neq m$.

Let us first consider Case 1. In this case, we have $\ell=m$. Thus,
$m=\ell\geq2>1$ and $j_{m}+g_{1}=\underbrace{j_{\ell}}_{>2}+\underbrace{g_{1}%
}_{\geq0}>2$.

From $G=\left(  g_{1},g_{2},\ldots,g_{p}\right)  $ and
\begin{align*}
K  &  =\left(  j_{1},j_{2},\ldots,j_{\ell-1},1,j_{\ell}-1,j_{\ell+1}%
,j_{\ell+2},\ldots,j_{m}\right) \\
&  =\left(  j_{1},j_{2},\ldots,j_{m-1},1,j_{m}-1\right)
\ \ \ \ \ \ \ \ \ \ \left(  \text{since }\ell=m\right)  ,
\end{align*}
we obtain%
\begin{align}
K\odot G  &  =\left(  j_{1},j_{2},\ldots,j_{m-1},1,\left(  j_{m}-1\right)
+g_{1},g_{2},g_{3},\ldots,g_{p}\right) \nonumber\\
&  =\left(  j_{1},j_{2},\ldots,j_{m-1},1,j_{m}+g_{1}-1,g_{2},g_{3}%
,\ldots,g_{p}\right)  . \label{pf.thm.Epk.dend.c7.pf.KG=c1}%
\end{align}

From looking at (\ref{pf.thm.Epk.dend.c7.pf.JG=}) and
(\ref{pf.thm.Epk.dend.c7.pf.KG=c1}), we conclude immediately that the
composition $K\odot G$ is obtained from $J\odot G$ by \textquotedblleft
splitting\textquotedblright\ the entry $j_{m}+g_{1}>2$ into two consecutive
entries $1$ and $j_{m}+g_{1}-1$, and that this entry $j_{m}+g_{1}$ was not the
first entry (indeed, this entry is the $m$-th entry, but $m>1$). Hence,
$J\odot G\rightarrow K\odot G$ (by the definition of the relation
$\rightarrow$). This proves Claim 7 in Case 1.

Let us next consider Case 2. In this case, we have $\ell\neq m$. Hence,
$\ell\in\left\{  2,3,\ldots,m-1\right\}  $ (since $\ell\in\left\{
2,3,\ldots,m\right\}  $).

From $G=\left(  g_{1},g_{2},\ldots,g_{p}\right)  $ and $K=\left(  j_{1}%
,j_{2},\ldots,j_{\ell-1},1,j_{\ell}-1,j_{\ell+1},j_{\ell+2},\ldots
,j_{m}\right)  $, we obtain%
\begin{align}
&  K\odot G\nonumber\\
&  =\left(  j_{1},j_{2},\ldots,j_{\ell-1},1,j_{\ell}-1,j_{\ell+1},j_{\ell
+2},\ldots,j_{m-1},j_{m}+g_{1},g_{2},g_{3},\ldots,g_{p}\right)
\label{pf.thm.Epk.dend.c7.pf.KG=c2}%
\end{align}
(notice that the $j_{m}+g_{1}$ term is \textbf{not} a $\left(  j_{\ell
}-1\right)  +g_{1}$ term, because $\ell\neq m$).

From looking at (\ref{pf.thm.Epk.dend.c7.pf.JG=}) and
(\ref{pf.thm.Epk.dend.c7.pf.KG=c2}), we conclude immediately that the
composition $K\odot G$ is obtained from $J\odot G$ by \textquotedblleft
splitting\textquotedblright\ the entry $j_{\ell}>2$ into two consecutive
entries $1$ and $j_{\ell}-1$, and that this entry $j_{\ell}$ was not the first
entry (indeed, this entry is the $\ell$-th entry, but $\ell>1$). Hence,
$J\odot G\rightarrow K\odot G$ (by the definition of the relation
$\rightarrow$). This proves Claim 7 in Case 2.

We have now proven Claim 7 in both Cases 1 and 2. Thus, Claim 7 always holds.]

\begin{statement}
\textit{Claim 8:} We have $\mathcal{K}_{\operatorname{Epk}}\bel
A\subseteq\mathcal{K}_{\operatorname{Epk}}$.
\end{statement}

[\textit{Proof of Claim 8:} This is proven in the same way as we proved Claim
6, with the only difference that now we need to use Claim 7 instead of Claim 5.]

Combining Claim 6 and Claim 8, we conclude that $\mathcal{K}%
_{\operatorname{Epk}}$ is a $\bel$-ideal of $A=\operatorname*{QSym}$.

Let $A^+$ denote the $\mathbf{k}$-submodule of $A$ consisting of those
quasisymmetric functions whose constant term is $0$.
It is easy to see that $\mathcal{K}_{\operatorname{Epk}} \subseteq A^+$
\ \ \ \ \footnote{\textit{Proof.} It suffices to show that for any
two $\operatorname{Epk}$-equivalent compositions $J$ and $K$,
the difference $F_J - F_K$ belongs to $A^+$ (since these differences
span $\mathcal{K}_{\operatorname{Epk}}$).
But this is easy: If $J$ and $K$ are two $\operatorname{Epk}$-equivalent
compositions, then $J$ and $K$ are either both empty (in which case
$F_J - F_K = F_{\varnothing} - F_{\varnothing} = 0 \in A^+$) or both
nonempty (in which case $F_J$ and $F_K$ and therefore also $F_J - F_K$
are homogeneous of positive degree, whence $F_J - F_K \in A^+$ again).}.
Furthermore, $A \tvi \mathcal{K}_{\operatorname{Epk}} \subseteq
\mathcal{K}_{\operatorname{Epk}}$ (since $\mathcal{K}_{\operatorname{Epk}}$
is a $\tvi$-ideal of $A$) and $A \bel \mathcal{K}_{\operatorname{Epk}} \subseteq
\mathcal{K}_{\operatorname{Epk}}$ (since $\mathcal{K}_{\operatorname{Epk}}$
is a $\bel$-ideal of $A$). Thus,
Theorem \ref{thm.ideal-crit2} \textbf{(c)} (applied to $M=\mathcal{K}%
_{\operatorname{Epk}}$) shows that $\mathcal{K}_{\operatorname{Epk}}$ is a
$\left.  \prec\right.  $-ideal and a $\left.  \succeq\right.  $-ideal of
$\operatorname*{QSym}$.

Thus, altogether, we have proven that $\mathcal{K}_{\operatorname{Epk}}$ is a
$\tvi$-ideal, a $\bel$-ideal, a $\left.  \prec\right.  $-ideal and a $\left.
\succeq\right.  $-ideal of $\operatorname*{QSym}$. This proves Theorem
\ref{thm.Epk.dend}.
\end{proof}

\begin{question}
What other descent statistics $\operatorname{st}$ have the property that
$\mathcal{K}_{\operatorname{st}}$ is an $\tvi$-ideal, $\bel$-ideal, $\left.
\prec\right.  $-ideal and/or $\left.  \succeq\right.  $-ideal? We will see
some answers in Subsection \ref{subsect.dendri.other-stats}, but a more
systematic study would be interesting.
\end{question}

\subsection{Dendriform shuffle-compatibility}

We have seen (in Proposition \ref{prop.K.ideal}) that the kernel
$\mathcal{K}_{\operatorname{st}}$ of a descent statistic $\operatorname{st}$
is an ideal of $\operatorname*{QSym}$ if and only if $\operatorname{st}$ is
shuffle-compatible. It is natural to ask whether similar combinatorial
interpretations exist for when the kernel $\mathcal{K}_{\operatorname{st}}$
of a descent statistic $\operatorname{st}$ is a $\tvi$-ideal, a $\bel$-ideal,
a $\left.  \prec\right.  $-ideal or a $\left.  \succeq\right.  $-ideal. In
this section, we shall prove such interpretations.

Now, let us define two further variants of LR-shuffle-compatibility (to be
compared with those introduced in Definition \ref{def.LR.left-right}):

\begin{definition}
Let $\operatorname{st}$ be a permutation statistic.

\textbf{(a)} We say that $\operatorname{st}$ is \textit{weakly
left-shuffle-compatible} if for any two disjoint nonempty permutations $\pi$
and $\sigma$ having the property that%
\begin{equation}
\text{each entry of }\pi\text{ is greater than each entry of }\sigma,
\label{eq.def.dendri.dsc.weak-ass}%
\end{equation}
the multiset $\left\{  \operatorname{st} \tau   \ \mid\ \tau\in
S_{\prec}\left(  \pi,\sigma\right)  \right\}  _{\operatorname*{multi}}$
depends only on $\operatorname{st} \pi$, $\operatorname{st} \sigma$,
$\left\vert \pi\right\vert $ and $\left\vert \sigma\right\vert $.

\textbf{(b)} We say that $\operatorname{st}$ is \textit{weakly
right-shuffle-compatible} if for any two disjoint nonempty permutations $\pi$
and $\sigma$ having the property that%
\[
\text{each entry of }\pi\text{ is greater than each entry of }\sigma,
\]
the multiset $\left\{  \operatorname{st} \tau   \ \mid\ \tau\in
S_{\succ}\left(  \pi,\sigma\right)  \right\}  _{\operatorname*{multi}}$
depends only on $\operatorname{st} \pi$, $\operatorname{st} \sigma$,
$\left\vert \pi\right\vert $ and $\left\vert \sigma\right\vert $.
\end{definition}

Then, the following analogues to the first part of Proposition
\ref{prop.K.ideal} hold:

\begin{theorem}
\label{thm.dendri.K.ideal}Let $\operatorname{st}$ be a descent statistic.
Then, the following three statements are equivalent:

\begin{itemize}
\item \textit{Statement A:} The statistic $\operatorname{st}$ is left-shuffle-compatible.

\item \textit{Statement B:} The statistic $\operatorname{st}$ is weakly left-shuffle-compatible.

\item \textit{Statement C:} The set $\mathcal{K}_{\operatorname{st}}$ is an
$\left.  \prec\right.  $-ideal of $\operatorname*{QSym}$.
\end{itemize}
\end{theorem}

\begin{theorem}
\label{thm.dendri.K.ideal-R}Let $\operatorname{st}$ be a descent statistic.
Then, the following three statements are equivalent:

\begin{itemize}
\item \textit{Statement A:} The statistic $\operatorname{st}$ is right-shuffle-compatible.

\item \textit{Statement B:} The statistic $\operatorname{st}$ is weakly right-shuffle-compatible.

\item \textit{Statement C:} The set $\mathcal{K}_{\operatorname{st}}$ is an
$\left.  \succeq\right.  $-ideal of $\operatorname*{QSym}$.
\end{itemize}
\end{theorem}

Let us prove Theorem \ref{thm.dendri.K.ideal} directly, without using shuffle algebras:

\begin{proof}
[Proof of Theorem \ref{thm.dendri.K.ideal} (sketched).]The implication
A$\Longrightarrow$B is obvious.

\textit{Proof of the implication B}$\Longrightarrow$\textit{C:} Assume that
Statement B holds. Thus, the statistic $\operatorname{st}$ is weakly left-shuffle-compatible.

Let us show that the set $\mathcal{K}_{\operatorname{st}}$ is a $\left.
\prec\right.  $-ideal of $\operatorname*{QSym}$. Indeed, it suffices to show
that every two $\operatorname{st}$-equivalent compositions $J$ and $K$ and
every further composition $L$ satisfy%
\begin{equation}
\left(  F_{J}-F_{K}\right)  \left.  \prec\right.  F_{L}\in\mathcal{K}%
_{\operatorname{st}}\ \ \ \ \ \ \ \ \ \ \text{and}\ \ \ \ \ \ \ \ \ \ F_{L}%
\left.  \prec\right.  \left(  F_{J}-F_{K}\right)  \in\mathcal{K}%
_{\operatorname{st}} \label{pf.thm.dendri.K.ideal.two-incs}%
\end{equation}
(because of the definition of $\mathcal{K}_{\operatorname{st}}$). So let $J$
and $K$ be two $\operatorname{st}$-equivalent compositions, and let $L$ be a
further composition. If $J=K$, then (\ref{pf.thm.dendri.K.ideal.two-incs})
follows immediately from realizing that $F_{J}-F_{K}=0$; thus, we WLOG assume
that $J\neq K$. But $\left\vert J\right\vert =\left\vert K\right\vert $ (since
$J$ and $K$ are $\operatorname{st}$-equivalent). Hence, $\left\vert
J\right\vert =\left\vert K\right\vert >0$ (since otherwise, we would have
$\left\vert J\right\vert =\left\vert K\right\vert =0$, which would imply that
both $J$ and $K$ would be the empty composition, contradicting $J\neq K$).
Thus, the power series $F_{J}$ and $F_{K}$ are homogeneous of degree
$\left\vert J\right\vert =\left\vert K\right\vert >0$; consequently,
$\varepsilon\left(  F_{J}\right)  =0$ and $\varepsilon\left(  F_{K}\right)
=0$. Hence, $\varepsilon\left(  F_{J}-F_{K}\right)  =\underbrace{\varepsilon
\left(  F_{J}\right)  }_{=0}-\underbrace{\varepsilon\left(  F_{K}\right)
}_{=0}=0$.

The compositions $J$ and $K$ are nonempty (since $\left\vert J\right\vert
=\left\vert K\right\vert >0$). If $L$ is empty, then
(\ref{pf.thm.dendri.K.ideal.two-incs}) holds for easy reasons (indeed, we have
$F_{L}=1$ in this case, and therefore (\ref{eq.dendriform.a<1}) yields
\[
\left(  F_{J}-F_{K}\right)  \left.  \prec\right.  F_{L}=\left(  F_{J}%
-F_{K}\right)  -\underbrace{\varepsilon\left(  F_{J}-F_{K}\right)
}_{\substack{=0}}=F_{J}-F_{K}\in\mathcal{K}_{\operatorname{st}},
\]
and similarly (\ref{eq.dendriform.1<a}) leads to $F_{L}\left.  \prec\right.
\left(  F_{J}-F_{K}\right)  \in\mathcal{K}_{\operatorname{st}}$). Hence, we
WLOG assume that $L$ is nonempty.

Pick three disjoint permutations $\varphi$, $\psi$ and $\sigma$ having descent
compositions $J$, $K$ and $L$, respectively, and having the property that%
\[
\text{each entry of }\varphi\text{ is greater than each entry of }\sigma
\]
and%
\[
\text{each entry of }\psi\text{ is greater than each entry of }\sigma\text{.}%
\]
(Such permutations $\varphi$, $\psi$ and $\sigma$ exist, since the set
$\mathbb{P}$ is infinite.)

The permutations $\varphi$ and $\psi$ are $\operatorname{st}$-equivalent
(since their descent compositions $J$ and $K$ are $\operatorname{st}%
$-equivalent). In other words, $\left\vert \varphi\right\vert =\left\vert
\psi\right\vert $ and $\operatorname{st} \varphi
=\operatorname{st} \psi$.

The statistic $\operatorname{st}$ is weakly left-shuffle-compatible. Thus,
the multiset\newline$\left\{  \operatorname{st} \tau 
\ \mid\ \tau\in S_{\prec}\left(  \pi,\sigma\right)  \right\}
_{\operatorname*{multi}}$ (where $\pi$ is a nonempty permutation disjoint from
$\sigma$ and having the property that each entry of $\pi$ is greater than each
entry of $\sigma$) depends only on $\operatorname{st} \pi   $
and $\left\vert \pi\right\vert $ (by the definition of \textquotedblleft
weakly left-shuffle-compatible\textquotedblright)\footnote{Recall that
$\sigma$ is fixed here, which is why we don't have to say that it depends on
$\operatorname{st} \sigma   $ and $\left\vert \sigma\right\vert
$ as well.}. Therefore, the multisets $\left\{  \operatorname{st} \tau
\ \mid\ \tau\in S_{\prec}\left(  \varphi,\sigma\right)  \right\}
_{\operatorname*{multi}}$ and \newline$\left\{  \operatorname{st} \tau
\ \mid\ \tau\in S_{\prec}\left(  \psi,\sigma\right)  \right\}
_{\operatorname*{multi}}$ are equal (since $\left\vert \varphi\right\vert
=\left\vert \psi\right\vert $ and $\operatorname{st} \varphi
=\operatorname{st} \psi$). Hence, there exists a bijection
$\alpha:S_{\prec}\left(  \varphi,\sigma\right)  \rightarrow S_{\prec}\left(
\psi,\sigma\right)  $ such that each $\chi\in S_{\prec}\left(  \varphi
,\sigma\right)  $ satisfies%
\begin{equation}
\operatorname{st}\left(  \alpha\left(  \chi\right)  \right)
=\operatorname{st} \chi
. \label{pf.thm.dendri.K.ideal.alpha}
\end{equation}
Consider this $\alpha$. Clearly, each $\chi\in S_{\prec}\left(  \varphi
,\sigma\right)  $ satisfies%
\[
\left(  \chi\text{ and }\alpha\left(  \chi\right)  \text{ are }%
\operatorname{st}\text{-equivalent}\right)
\]
(because of (\ref{pf.thm.dendri.K.ideal.alpha}) and since $\left\vert
\chi\right\vert =\underbrace{\left\vert \varphi\right\vert }_{=\left\vert
\psi\right\vert }+\left\vert \sigma\right\vert =\left\vert \psi\right\vert
+\left\vert \sigma\right\vert =\left\vert \alpha\left(  \chi\right)
\right\vert $) and therefore%
\[
\left(  \operatorname*{Comp}\chi\text{ and }\operatorname*{Comp}\left(
\alpha\left(  \chi\right)  \right)  \text{ are }\operatorname{st}%
\text{-equivalent}\right)
\]
(since $\operatorname{st}$ is a descent statistic) and thus
$F_{\operatorname*{Comp}\chi}-F_{\operatorname*{Comp}\left(  \alpha\left(
\chi\right)  \right)  }\in\mathcal{K}_{\operatorname{st}}$ (by the definition
of $\mathcal{K}_{\operatorname{st}}$) and therefore%
\begin{equation}
F_{\operatorname*{Comp}\chi}\equiv F_{\operatorname*{Comp}\left(
\alpha\left(  \chi\right)  \right)  }\operatorname{mod}\mathcal{K}%
_{\operatorname{st}}. \label{pf.thm.dendri.K.ideal.alphaeq}%
\end{equation}

The first claim of Corollary \ref{cor.dendri.4.1} yields%
\begin{align*}
F_{\operatorname*{Comp}\varphi}\left.  \prec\right.  F_{\operatorname*{Comp}%
\sigma}  &  =\sum_{\chi\in S_{\prec}\left(  \varphi,\sigma\right)
}F_{\operatorname*{Comp}\chi}\ \ \ \ \ \ \ \ \ \ \text{and}\\
F_{\operatorname*{Comp}\psi}\left.  \prec\right.  F_{\operatorname*{Comp}%
\sigma}  &  =\sum_{\chi\in S_{\prec}\left(  \psi,\sigma\right)  }%
F_{\operatorname*{Comp}\chi}.
\end{align*}
Hence,%
\begin{align*}
F_{\operatorname*{Comp}\varphi}\left.  \prec\right.  F_{\operatorname*{Comp}%
\sigma}  &  =\sum_{\chi\in S_{\prec}\left(  \varphi,\sigma\right)
}\underbrace{F_{\operatorname*{Comp}\chi}}_{\substack{\equiv
F_{\operatorname*{Comp}\left(  \alpha\left(  \chi\right)  \right)
}\operatorname{mod}\mathcal{K}_{\operatorname{st}}\\\text{(by
(\ref{pf.thm.dendri.K.ideal.alphaeq}))}}}\equiv\sum_{\chi\in S_{\prec}\left(
\varphi,\sigma\right)  }F_{\operatorname*{Comp}\left(  \alpha\left(
\chi\right)  \right)  }\\
&  =\sum_{\chi\in S_{\prec}\left(  \psi,\sigma\right)  }%
F_{\operatorname*{Comp}\chi}\\
&  \ \ \ \ \ \ \ \ \ \ \left(
\begin{array}
[c]{c}%
\text{here, we have substituted }\chi\text{ for }\alpha\left(  \chi\right)
\text{ in the sum,}\\
\text{since the map }\alpha:S_{\prec}\left(  \varphi,\sigma\right)
\rightarrow S_{\prec}\left(  \psi,\sigma\right)  \text{ is a bijection}%
\end{array}
\right) \\
&  =F_{\operatorname*{Comp}\psi}\left.  \prec\right.  F_{\operatorname*{Comp}%
\sigma}\operatorname{mod}\mathcal{K}_{\operatorname{st}}.
\end{align*}
Since $\operatorname*{Comp}\varphi=J$, $\operatorname*{Comp}\psi=K$ and
$\operatorname*{Comp}\sigma=L$ (by the definition of $\varphi$, $\psi$ and
$\sigma$), this rewrites as $F_{J}\left.  \prec\right.  F_{L}\equiv
F_{K}\left.  \prec\right.  F_{L}\operatorname{mod}\mathcal{K}%
_{\operatorname{st}}$. In other words, $F_{J}\left.  \prec\right.
F_{L}-F_{K}\left.  \prec\right.  F_{L}\in\mathcal{K}_{\operatorname{st}}$. In
other words, $\left(  F_{J}-F_{K}\right)  \left.  \prec\right.  F_{L}%
\in\mathcal{K}_{\operatorname{st}}$. This proves the first claim of
(\ref{pf.thm.dendri.K.ideal.two-incs}). The second is proven similarly.
Altogether, we thus conclude that $\mathcal{K}_{\operatorname{st}}$ is a
$\left.  \prec\right.  $-ideal of $\operatorname*{QSym}$. In other words,
Statement C holds. This proves the implication B$\Longrightarrow$C.

\textit{Proof of the implication C}$\Longrightarrow$\textit{A:} Assume that
Statement C holds. Thus, the set $\mathcal{K}_{\operatorname{st}}$ is an
$\left.  \prec\right.  $-ideal of $\operatorname*{QSym}$.

Let $X$ be the codomain of the map $\operatorname{st}$. Let $\mathbb{Q}%
\left[  X\right]  $ be the free $\mathbb{Q}$-vector space with basis $\left(
\left[  x\right]  \right)  _{x\in X}$. Then, we can define a $\mathbb{Q}%
$-linear map $\mathbf{st}:\operatorname*{QSym}\rightarrow\mathbb{Q}\left[
X\right]  ,\ F_{J}\mapsto\left[  \operatorname{st}J\right]  $. This map
$\mathbf{st}$ sends each of the generators of $\mathcal{K}_{\operatorname{st}%
}$ to $0$ (by the definition of $\mathcal{K}_{\operatorname{st}}$), and
therefore sends the whole $\mathcal{K}_{\operatorname{st}}$ to $0$. In other
words, $\mathbf{st}\left(  \mathcal{K}_{\operatorname{st}}\right)  =0$.

Now, consider any two disjoint nonempty permutations $\pi$ and $\sigma$ having
the property that $\pi_{1}>\sigma_{1}$. Also, consider two further disjoint
nonempty permutations $\pi^{\prime}$ and $\sigma^{\prime}$ having the property
that $\pi_{1}^{\prime}>\sigma_{1}^{\prime}$ and satisfying $\operatorname{st}
\pi =\operatorname{st}\left(  \pi^{\prime}\right)  $,
$\operatorname{st} \sigma   =\operatorname{st}\left(
\sigma^{\prime}\right)  $, $\left\vert \pi\right\vert =\left\vert \pi^{\prime
}\right\vert $ and $\left\vert \sigma\right\vert =\left\vert \sigma^{\prime
}\right\vert $. We shall show that
\[
\left\{  \operatorname{st} \tau
\ \mid\ \tau\in S_{\prec}\left(  \pi,\sigma\right)  \right\}
_{\operatorname*{multi}}=\left\{  \operatorname{st} \tau 
\ \mid\ \tau\in S_{\prec}\left(  \pi^{\prime},\sigma^{\prime}\right)
\right\}  _{\operatorname*{multi}} .
\]
This will show that the multiset
$\left\{  \operatorname{st} \tau   \ \mid\ \tau\in S_{\prec
}\left(  \pi,\sigma\right)  \right\}  _{\operatorname*{multi}}$ depends only
on $\operatorname{st} \pi   $, $\operatorname{st} \sigma$,
$\left\vert \pi\right\vert $ and $\left\vert \sigma \right\vert $.

From $\operatorname{st} \pi   =\operatorname{st}\left(
\pi^{\prime}\right)  $ and $\left\vert \pi\right\vert =\left\vert \pi^{\prime
}\right\vert $, we conclude that $\pi$ and $\pi^{\prime}$ are
$\operatorname{st}$-equivalent. In other words, $\operatorname*{Comp}\pi$ and
$\operatorname*{Comp}\left(  \pi^{\prime}\right)  $ are $\operatorname{st}%
$-equivalent. Hence, $F_{\operatorname*{Comp}\pi}-F_{\operatorname*{Comp}%
\left(  \pi^{\prime}\right)  }\in\mathcal{K}_{\operatorname{st}}$ (by the
definition of $\mathcal{K}_{\operatorname{st}}$), so that
$F_{\operatorname*{Comp}\pi}\equiv F_{\operatorname*{Comp}\left(  \pi^{\prime
}\right)  }\operatorname{mod}\mathcal{K}_{\operatorname{st}}$. Similarly,
$F_{\operatorname*{Comp}\sigma}\equiv F_{\operatorname*{Comp}\left(
\sigma^{\prime}\right)  }\operatorname{mod}\mathcal{K}_{\operatorname{st}}$.
These two congruences, combined, yield $F_{\operatorname*{Comp}\pi}\left.
\prec\right.  F_{\operatorname*{Comp}\sigma}\equiv F_{\operatorname*{Comp}%
\left(  \pi^{\prime}\right)  }\left.  \prec\right.  F_{\operatorname*{Comp}%
\left(  \sigma^{\prime}\right)  }\operatorname{mod}\mathcal{K}%
_{\operatorname{st}}$. (Indeed, we can conclude $a\left.  \prec\right.
c\equiv b\left.  \prec\right.  d\operatorname{mod}\mathcal{K}%
_{\operatorname{st}}$ whenever we have $a\equiv b\operatorname{mod}%
\mathcal{K}_{\operatorname{st}}$ and $c\equiv d\operatorname{mod}%
\mathcal{K}_{\operatorname{st}}$; this is because we know that $\mathcal{K}%
_{\operatorname{st}}$ is a $\left.  \prec\right.  $-ideal of
$\operatorname*{QSym}$.)

From $F_{\operatorname*{Comp}\pi}\left.  \prec\right.  F_{\operatorname*{Comp}%
\sigma}\equiv F_{\operatorname*{Comp}\left(  \pi^{\prime}\right)  }\left.
\prec\right.  F_{\operatorname*{Comp}\left(  \sigma^{\prime}\right)
}\operatorname{mod}\mathcal{K}_{\operatorname{st}}$, we obtain%
\begin{equation}
\mathbf{st}\left(  F_{\operatorname*{Comp}\pi}\left.  \prec\right.
F_{\operatorname*{Comp}\sigma}\right)  =\mathbf{st}\left(
F_{\operatorname*{Comp}\left(  \pi^{\prime}\right)  }\left.  \prec\right.
F_{\operatorname*{Comp}\left(  \sigma^{\prime}\right)  }\right)
\label{pf.thm.dendri.K.ideal.steq}%
\end{equation}
(since $\mathbf{st}\left(  \mathcal{K}_{\operatorname{st}}\right)  =0$).

The first claim of Corollary \ref{cor.dendri.4.1} yields%
\[
F_{\operatorname*{Comp}\pi}\left.  \prec\right.  F_{\operatorname*{Comp}%
\sigma}=\sum_{\chi\in S_{\prec}\left(  \pi,\sigma\right)  }%
F_{\operatorname*{Comp}\chi}.
\]
Applying the map $\mathbf{st}$ to both sides of this equality, we find%
\begin{align*}
\mathbf{st}\left(  F_{\operatorname*{Comp}\pi}\left.  \prec\right.
F_{\operatorname*{Comp}\sigma}\right)   &  =\mathbf{st}\left(  \sum_{\chi\in
S_{\prec}\left(  \pi,\sigma\right)  }F_{\operatorname*{Comp}\chi}\right) \\
&  =\sum_{\chi\in S_{\prec}\left(  \pi,\sigma\right)  }\underbrace{\mathbf{st}%
\left(  F_{\operatorname*{Comp}\chi}\right)  }_{=\left[  \operatorname{st}%
\left(  \operatorname*{Comp}\chi\right)  \right]  =\left[  \operatorname{st}%
\chi\right]  }=\sum_{\chi\in S_{\prec}\left(  \pi,\sigma\right)  }\left[
\operatorname{st}\chi\right]  .
\end{align*}
Similarly,%
\[
\mathbf{st}\left(  F_{\operatorname*{Comp}\left(  \pi^{\prime}\right)
}\left.  \prec\right.  F_{\operatorname*{Comp}\left(  \sigma^{\prime}\right)
}\right)  =\sum_{\chi\in S_{\prec}\left(  \pi^{\prime},\sigma^{\prime}\right)
}\left[  \operatorname{st}\chi\right]  .
\]
But the left-hand sides of the last two equalities are equal (because of
(\ref{pf.thm.dendri.K.ideal.steq})); therefore, the right-hand sides must be
equal as well. In other words,
\[
\sum_{\chi\in S_{\prec}\left(  \pi,\sigma\right)  }\left[  \operatorname{st}%
\chi\right]  =\sum_{\chi\in S_{\prec}\left(  \pi^{\prime},\sigma^{\prime
}\right)  }\left[  \operatorname{st}\chi\right]  .
\]
This shows exactly that $\left\{  \operatorname{st} \chi
\ \mid\ \chi\in S_{\prec}\left(  \pi,\sigma\right)  \right\}
_{\operatorname*{multi}} =\left\{  \operatorname{st} \chi
\ \mid\ \chi\in S_{\prec}\left(  \pi^{\prime},\sigma^{\prime}\right)
\right\}  _{\operatorname*{multi}} $. In other words, $\left\{
\operatorname{st} \tau   \ \mid\ \tau\in S_{\prec}\left(
\pi,\sigma\right)  \right\}  _{\operatorname*{multi}} =\left\{
\operatorname{st} \tau   \ \mid\ \tau\in S_{\prec}\left(
\pi^{\prime},\sigma^{\prime}\right)  \right\}  _{\operatorname*{multi}} $.
Thus, we have proven that the multiset $\left\{  \operatorname{st} \tau
\ \mid\ \tau\in S_{\prec}\left(  \pi,\sigma\right)  \right\}
_{\operatorname*{multi}} $ depends only on $\operatorname{st} \pi$,
$\operatorname{st} \sigma$, $\left\vert \pi\right\vert $ and
$\left\vert \sigma\right\vert $. Hence, the statistic
$\operatorname{st}$ is left-shuffle-compatible. In other words, Statement A
holds. This proves the implication C$\Longrightarrow$A.

Now that we have proven all three implications A$\Longrightarrow$B,
B$\Longrightarrow$C and C$\Longrightarrow$A, the proof of Theorem
\ref{thm.dendri.K.ideal} is complete.
\end{proof}

\begin{proof}
[Proof of Theorem \ref{thm.dendri.K.ideal-R}.]The proof of Theorem
\ref{thm.dendri.K.ideal-R} is analogous to the above proof of Theorem
\ref{thm.dendri.K.ideal}.
\end{proof}

\begin{corollary}
\label{cor.dendri.K.ideal-LR}Let $\operatorname{st}$ be a permutation
statistic that is LR-shuffle-compatible. Then, $\operatorname{st}$ is a
shuffle-compatible descent statistic, and the set $\mathcal{K}%
_{\operatorname{st}}$ is an ideal and a $\left.  \prec\right.  $-ideal and a
$\left.  \succeq\right.  $-ideal of $\operatorname*{QSym}$.
\end{corollary}

\begin{proof}
[Proof of Corollary \ref{cor.dendri.K.ideal-LR} (sketched).] Proposition
\ref{prop.LRcomp.head} yields that $\operatorname{st}$ is
head-graft-compatible and shuffle-compatible. Proposition
\ref{prop.LRcomp.equivs} shows that $\operatorname{st}$ is
left-shuffle-compatible, right-shuffle-compatible and head-graft-compatible
(since $\operatorname{st}$ is LR-shuffle-compatible). Hence, Proposition
\ref{prop.des-stat.hgc} shows that $\operatorname{st}$ is a descent
statistic. Thus, Theorem \ref{thm.dendri.K.ideal} yields that $\mathcal{K}%
_{\operatorname{st}}$ is a $\left.  \prec\right.  $-ideal of
$\operatorname*{QSym}$ (since $\operatorname{st}$ is
left-shuffle-compatible). Likewise, Theorem \ref{thm.dendri.K.ideal-R} yields
that $\mathcal{K}_{\operatorname{st}}$ is a $\left.  \succeq\right.  $-ideal
of $\operatorname*{QSym}$ (since $\operatorname{st}$ is
right-shuffle-compatible). Finally, Proposition \ref{prop.K.ideal} yields that
$\mathcal{K}_{\operatorname{st}}$ is an ideal of $\operatorname*{QSym}$
(since $\operatorname{st}$ is a shuffle-compatible descent statistic). This
proves Corollary \ref{cor.dendri.K.ideal-LR}.
\end{proof}

A converse of Corollary \ref{cor.dendri.K.ideal-LR} also holds:

\begin{corollary}
\label{cor.dendri.K.ideal-LRi}Let $\operatorname{st}$ be a descent statistic
such that $\mathcal{K}_{\operatorname{st}}$ is a $\left.  \prec\right.
$-ideal and a $\left.  \succeq\right.  $-ideal of $\operatorname*{QSym}$.
Then, $\operatorname{st}$ is LR-shuffle-compatible and shuffle-compatible.
\end{corollary}

\begin{proof}
[Proof of Corollary \ref{cor.dendri.K.ideal-LRi} (sketched).] Theorem
\ref{thm.dendri.K.ideal} yields that $\operatorname{st}$ is
left-shuffle-compatible (since $\mathcal{K}_{\operatorname{st}}$ is an
$\left.  \prec\right.  $-ideal of $\operatorname*{QSym}$). Likewise, Theorem
\ref{thm.dendri.K.ideal-R} yields that $\operatorname{st}$ is
right-shuffle-compatible (since $\mathcal{K}_{\operatorname{st}}$ is an
$\left.  \succeq\right.  $-ideal of $\operatorname*{QSym}$). Hence, Corollary
\ref{cor.LRcomp.two} shows that $\operatorname{st}$ is LR-shuffle-compatible.
Thus, Proposition \ref{prop.LRcomp.head} yields that $\operatorname{st}$ is
head-graft-compatible and shuffle-compatible. This proves Corollary
\ref{cor.dendri.K.ideal-LRi}.
\end{proof}

As a consequence of Theorem \ref{thm.dendri.K.ideal} and Theorem
\ref{thm.dendri.K.ideal-R}, we can see that any descent statistic that is
weakly left-shuffle-compatible and weakly right-shuffle-compatible must
automatically be shuffle-compatible\footnote{\textit{Proof.} Let
$\operatorname{st}$ be a descent statistic that is weakly
left-shuffle-compatible and weakly right-shuffle-compatible. We must prove
that $\operatorname{st}$ is shuffle-compatible.
\par
The implication B$\Longrightarrow$C in Theorem \ref{thm.dendri.K.ideal} shows
that the set $\mathcal{K}_{\operatorname{st}}$ is a $\left.  \prec\right.
$-ideal of $\operatorname*{QSym}$. Similarly, the set $\mathcal{K}%
_{\operatorname{st}}$ is a $\left.  \succeq\right.  $-ideal of
$\operatorname*{QSym}$. Hence, Proposition \ref{prop.ideal-crit3} (applied to
$M=\mathcal{K}_{\operatorname{st}}$) yields that $\mathcal{K}%
_{\operatorname{st}}$ is an ideal of $\operatorname*{QSym}$. By Proposition
\ref{prop.K.ideal}, this shows that $\operatorname{st}$ is
shuffle-compatible.}. Note that this is only true for descent statistics! As
far as arbitrary permutation statistics are concerned, this is false; for
example, the number of inversions is weakly left-shuffle-compatible and weakly
right-shuffle-compatible but not shuffle-compatible.

Recall that every permutation statistic that is left-shuffle-compatible and
right-shuffle-compatible must automatically be LR-shuffle-compatible (by
Corollary \ref{cor.LRcomp.two}) and therefore also shuffle-compatible (by
Corollary \ref{cor.LRcomp.back}) and head-graft-compatible (again by Corollary
\ref{cor.LRcomp.back}) and therefore a descent statistic (by Proposition
\ref{prop.des-stat.hgc}).

\begin{corollary}
\label{cor.dendri.Epk}The descent statistic $\operatorname{Epk}$ is
left-shuffle-compatible and right-shuffle-compatible.
\end{corollary}

Corollary \ref{cor.dendri.Epk} follows by combining Theorem
\ref{thm.LRcomp.Pks} \textbf{(c)} with Proposition \ref{prop.LRcomp.equivs}. But
we can also give a proof using Theorem \ref{thm.dendri.K.ideal}:

\begin{proof}
[Proof of Corollary \ref{cor.dendri.Epk}.]To prove that $\operatorname{Epk}$
is left-shuffle-compatible, combine Theorem \ref{thm.dendri.K.ideal} with
Theorem \ref{thm.Epk.dend}. Similarly for right-shuffle-compatibility.
\end{proof}

Using Theorem \ref{thm.dendri.4.1}, we can state an analogue of Theorem
\ref{thm.4.3}. Let us first define the notion of dendriform algebras:

\begin{definition}
\textbf{(a)} A \textit{dendriform algebra} over a field $\mathbf{k}$ means a
$\mathbf{k}$-algebra $A$ equipped with two further $\mathbf{k}$-bilinear
binary operations $\left.  \prec\right.  $ and $\left.  \succeq\right.  $
(these are operations, not relations, despite the symbols) from $A\times A$ to
$A$ that satisfy the four rules%
\begin{align*}
a\left.  \prec\right.  b+a\left.  \succeq\right.  b  &  =ab;\\
\left(  a\left.  \prec\right.  b\right)  \left.  \prec\right.  c  &  =a\left.
\prec\right.  \left(  bc\right)  ;\\
\left(  a\left.  \succeq\right.  b\right)  \left.  \prec\right.  c  &
=a\left.  \succeq\right.  \left(  b\left.  \prec\right.  c\right)  ;\\
a\left.  \succeq\right.  \left(  b\left.  \succeq\right.  c\right)   &
=\left(  ab\right)  \left.  \succeq\right.  c
\end{align*}
for all $a,b,c\in A$. (Depending on the situation, it is useful to also impose
a few axioms that relate the unity $1$ of the $\mathbf{k}$-algebra $A$ with
the operations $\left.  \prec\right.  $ and $\left.  \succeq\right.  $. For
example, we could require $1\left.  \prec\right.  a=0$ for each $a\in A$. For
what we are going to do in the following, it does not matter whether we make
this requirement.)

\textbf{(b)} If $A$ and $B$ are two dendriform algebras over $\mathbf{k}$,
then a \textit{dendriform algebra homomorphism} from $A$ to $B$ means a
$\mathbf{k}$-algebra homomorphism $\phi:A\rightarrow B$ preserving the
operations $\left.  \prec\right.  $ and $\left.  \succeq\right.  $ (that is,
satisfying $\phi\left(  a\left.  \prec\right.  b\right)  =\phi\left(
a\right)  \left.  \prec\right.  \phi\left(  b\right)  $ and $\phi\left(
a\left.  \succeq\right.  b\right)  =\phi\left(  a\right)  \left.
\succeq\right.  \phi\left(  b\right)  $ for all $a,b\in A$). (Some authors
only require it to be a $\mathbf{k}$-linear map instead of being a
$\mathbf{k}$-algebra homomorphism; this boils down to the question whether
$\phi\left(  1\right)  $ must be $1$ or not. This does not make a difference
for us here.)
\end{definition}

Thus, $\operatorname*{QSym}$ (with its two operations $\left.  \prec\right.  $
and $\left.  \succeq\right.  $) becomes a dendriform algebra over $\mathbb{Q}$.

Notice that if $A$ and $B$ are two dendriform algebras over $\mathbf{k}$, then
the kernel of any dendriform algebra homomorphism $A\rightarrow B$ is an
$\left.  \prec\right.  $-ideal and a $\left.  \succeq\right.  $-ideal of $A$.
Conversely, if $A$ is a dendriform algebra over $\mathbf{k}$, and $I$ is
simultaneously a $\left.  \prec\right.  $-ideal and a $\left.  \succeq
\right.  $-ideal of $A$, then $A/I$ canonically becomes a dendriform algebra,
and the canonical projection $A\rightarrow A/I$ becomes a dendriform algebra homomorphism.

Therefore, Theorem \ref{thm.dendri.K.ideal} and Theorem
\ref{thm.dendri.K.ideal-R} (and the $\mathcal{A}_{\operatorname{st}}%
\cong\operatorname*{QSym}/\mathcal{K}_{\operatorname{st}}$ isomorphism from
Proposition \ref{prop.K.ideal}) yield the following:

\begin{corollary}
\label{cor.dendri.quotient-dendri}If a descent statistic $\operatorname{st}$
is left-shuffle-compatible and right-shuffle-compatible, then its shuffle
algebra $\mathcal{A}_{\operatorname{st}}$ canonically becomes a dendriform algebra.
\end{corollary}

We furthermore have the following analogue of Theorem~\ref{thm.4.3}, which
easily follows from Theorem \ref{thm.dendri.K.ideal} and Theorem
\ref{thm.dendri.K.ideal-R}:

\begin{theorem}
\label{thm.dendri.4.3}Let $\operatorname{st}$ be a descent statistic.

\textbf{(a)} The descent statistic $\operatorname{st}$ is
left-shuffle-compatible and right-shuffle-compatible if and only if there
exist a dendriform algebra $A$ with
basis $\left(  u_{\alpha}\right)  $ (indexed by $\operatorname{st}%
$-equivalence classes $\alpha$ of compositions)
and a dendriform algebra homomorphism
$\phi_{\operatorname{st}}:\operatorname*{QSym} \rightarrow A$
with the property that whenever $\alpha$ is an
$\operatorname{st}$-equivalence class of compositions, we have
\[
\phi_{\operatorname{st}}\left(  F_{L}\right)  =u_{\alpha}%
\ \ \ \ \ \ \ \ \ \ \text{for each }L\in\alpha.
\]

\textbf{(b)} In this case, the $\mathbb{Q}$-linear map%
\[
\mathcal{A}_{\operatorname{st}}\rightarrow A,\ \ \ \ \ \ \ \ \ \ \left[
\pi\right]  _{\operatorname{st}}\mapsto u_{\alpha},
\]
where $\alpha$ is the $\operatorname{st}$-equivalence class of the
composition $\operatorname*{Comp}\pi$, is an isomorphism of dendriform
algebras $\mathcal{A}_{\operatorname{st}}\rightarrow A$.
\end{theorem}

\begin{question}
Can the $\mathbb{Q}$-algebra $\operatorname*{Pow}\mathcal{N}$ from Definition
\ref{def.GammaZ} be endowed with two binary operations $\left.  \prec\right.
$ and $\left.  \succeq\right.  $ that make it into a dendriform algebra? Can
we then find an analogue of Proposition \ref{prop.prod1} along the following lines?

Let $\left(  P,\gamma\right)  $, $\left(  Q,\delta\right)  $ and $\left(
P\sqcup Q,\varepsilon\right)  $ be as in Proposition \ref{prop.prod1}. Assume
that each of the posets $P$ and $Q$ has a minimum element; denote these
elements by $\min P$ and $\min Q$, respectively. Define two posets $P\left.
\prec\right.  Q$ and $P\left.  \succeq\right.  Q$ as in Proposition
\ref{prop.dendri.4.1.lem}. Then, we hope to have%
\begin{align*}
\Gamma_{\mathcal{Z}}\left(  P,\gamma\right)  \left.  \prec\right.
\Gamma_{\mathcal{Z}}\left(  Q,\delta\right)   &  =\Gamma_{\mathcal{Z}}\left(
P\left.  \prec\right.  Q,\varepsilon\right)  \ \ \ \ \ \ \ \ \ \ \text{and}\\
\Gamma_{\mathcal{Z}}\left(  P,\gamma\right)  \left.  \succeq\right.
\Gamma_{\mathcal{Z}}\left(  Q,\delta\right)   &  =\Gamma_{\mathcal{Z}}\left(
P\left.  \succeq\right.  Q,\varepsilon\right) ,
\end{align*}
assuming a simple condition on $\min P$ and $\min Q$
(say, $\gamma\left(\min P\right) <_{\ZZ} \delta\left(\min Q\right)$).

Ideally, this would be a generalization of Proposition
\ref{prop.dendri.4.1.lem}.
\end{question}

\subsection{Criteria for $\mathcal{K}_{\operatorname{st}}$ to be a stack
ideal}

We have so far studied the combinatorial significance of when the kernel
$\mathcal{K}_{\operatorname{st}}$ of a statistic $\operatorname{st}$ is a
$\left.  \prec\right.  $-ideal or a $\left.  \succeq\right.  $-ideal of
$\operatorname*{QSym}$. What about $\bel$-ideals and $\tvi$-ideals? It turns
out that the answer to this question is given (on the level of compositions)
by the following (easily verified) proposition:

\begin{proposition}
\label{prop.bel-tvi-comp}Let $\operatorname{st}$ be a descent statistic.

\textbf{(a)} The set $\mathcal{K}_{\operatorname{st}}$ is a left $\bel$-ideal
of $\operatorname*{QSym}$ if and only if $\operatorname{st}$ has the
following property: If $J$ and $K$ are two $\operatorname{st}$-equivalent
nonempty compositions, and if $G$ is any nonempty composition, then $G\odot J$
and $G\odot K$ are $\operatorname{st}$-equivalent.

\textbf{(b)} The set $\mathcal{K}_{\operatorname{st}}$ is a right
$\bel$-ideal of $\operatorname*{QSym}$ if and only if $\operatorname{st}$ has
the following property: If $J$ and $K$ are two $\operatorname{st}$-equivalent
nonempty compositions, and if $G$ is any nonempty composition, then $J\odot G$
and $K\odot G$ are $\operatorname{st}$-equivalent.

\textbf{(c)} The set $\mathcal{K}_{\operatorname{st}}$ is a left $\tvi$-ideal
of $\operatorname*{QSym}$ if and only if $\operatorname{st}$ has the
following property: If $J$ and $K$ are two $\operatorname{st}$-equivalent
nonempty compositions, and if $G$ is any nonempty composition, then $\left[
G,J\right]  $ and $\left[  G,K\right]  $ are $\operatorname{st}$-equivalent.

\textbf{(d)} The set $\mathcal{K}_{\operatorname{st}}$ is a right
$\tvi$-ideal of $\operatorname*{QSym}$ if and only if $\operatorname{st}$ has
the following property: If $J$ and $K$ are two $\operatorname{st}$-equivalent
nonempty compositions, and if $G$ is any nonempty composition, then $\left[
J,G\right]  $ and $\left[  K,G\right]  $ are $\operatorname{st}$-equivalent.
\end{proposition}

Proposition \ref{prop.bel-tvi-comp} allows us to give a new proof of Theorem
\ref{thm.Epk.dend}, which makes no use of Proposition \ref{prop.K.Epk.F}.
Instead, it will rely on analyzing $\operatorname{Epk}\left(  \left[
A,B\right]  \right)  $ and $\operatorname{Epk}\left(  A\odot B\right)  $ when
$A$ and $B$ are two nonempty compositions.

First, we introduce a notation: If $S$ is a set of integers, and $p$ is an
integer, then $S+p$ shall denote the set $\left\{  s+p\ \mid\ s\in S\right\}
$.

We shall use the following simple lemma:

\begin{lemma}
\label{lem.Epk.AB}Let $A$ and $B$ be two nonempty compositions. Let
$n=\left\vert A\right\vert $.

\textbf{(a)} We have $\operatorname{Epk}\left(  \left[  A,B\right]  \right)
=\left(  \operatorname{Epk}A\right)  \cup\left(  \left(  \operatorname{Epk}%
B+n\right)  \setminus\left\{  n+1\right\}  \right)  $.

\textbf{(b)} We have $\operatorname{Epk}\left(  A\odot B\right)  =\left(
\left(  \operatorname{Epk}A\right)  \setminus\left\{  n\right\}  \right)
\cup\left(  \operatorname{Epk}B+n\right)  $.
\end{lemma}

\begin{proof}
[Proof of Lemma \ref{lem.Epk.AB}.]Let $m=\left\vert B\right\vert $. Consider
any $n$-permutation $\alpha=\left(  \alpha_{1},\alpha_{2},\ldots,\alpha
_{n}\right)  $ satisfying $\operatorname*{Comp}\alpha=A$. (Such $\alpha$
exists, since $n=\left\vert A\right\vert $.) Consider any $m$-permutation
$\beta=\left(  \beta_{1},\beta_{2},\ldots,\beta_{m}\right)  $ satisfying
$\operatorname*{Comp}\beta=B$. (Such $\beta$ exists, since $m=\left\vert
B\right\vert $.) From $\operatorname*{Comp}\alpha=A$, we obtain
$\operatorname{Epk}\alpha=\operatorname{Epk}A$. Similarly,
$\operatorname{Epk}\beta=\operatorname{Epk}B$.

\textbf{(a)} WLOG assume that $\alpha_{i}>\beta_{j}$ for all $i\in\left[
n\right]  $ and $j\in\left[  m\right]  $. (Indeed, we can achieve this by
choosing a positive integer $g$ that is larger than each entry of $\beta$, and
adding $g$ to each entry of $\alpha$.) Thus, in particular, the entries of
$\alpha$ are distinct from the entries of $\beta$. Also, $\alpha_{n}>\beta
_{1}$ (since $\alpha_{i}>\beta_{j}$ for all $i\in\left[  n\right]  $ and
$j\in\left[  m\right]  $).

Let $\gamma$ be the $\left(  n+m\right)  $-permutation $\left(  \alpha
_{1},\alpha_{2},\ldots,\alpha_{n},\beta_{1},\beta_{2},\ldots,\beta_{m}\right)
$. Then, the descents of $\gamma$ are obtained as follows:

\begin{itemize}
\item Each descent of $\alpha$ is a descent of $\gamma$.

\item The number $n$ is a descent of $\gamma$ (since $\alpha_{n}>\beta_{1}$).

\item Adding $n$ to each descent of $\beta$ yields a descent of $\gamma$ (that
is, if $i$ is a descent of $\beta$, then $i+n$ is a descent of $\gamma$).
\end{itemize}

These are all the descents of $\gamma$. Thus,%
\[
\operatorname{Des}\gamma=\operatorname{Des}\alpha\cup\left\{  n\right\}
\cup\left(  \operatorname{Des}\beta+n\right)  .
\]
Hence,%
\begin{align*}
\operatorname*{Comp}\left(  \operatorname{Des}\gamma\right)   &
=\operatorname*{Comp}\left(  \operatorname{Des}\alpha\cup\left\{  n\right\}
\cup\left(  \operatorname{Des}\beta+n\right)  \right) \\
&  =\left[  \operatorname*{Comp}\left(  \operatorname{Des}\alpha\right)
,\operatorname*{Comp}\left(  \operatorname{Des}\beta\right)  \right]
\end{align*}
(because of how $\operatorname*{Comp}S$ is defined for a set $S$). Since
$\operatorname*{Comp}\left(  \operatorname{Des}\pi\right)
=\operatorname*{Comp}\pi$ for any permutation $\pi$, this rewrites as%
\[
\operatorname*{Comp}\gamma=\left[  \operatorname*{Comp}\alpha
,\operatorname*{Comp}\beta\right]  .
\]
In view of $\operatorname*{Comp}\alpha=A$ and $\operatorname*{Comp}\beta=B$,
this rewrites as $\operatorname*{Comp}\gamma=\left[  A,B\right]  $. Thus,
$\operatorname{Epk}\left(  \left[  A,B\right]  \right)  =\operatorname{Epk}%
\gamma$.

On the other hand, recall again that $\gamma=\left(  \alpha_{1},\alpha
_{2},\ldots,\alpha_{n},\beta_{1},\beta_{2},\ldots,\beta_{m}\right)  $ and
$\alpha_{n}>\beta_{1}$. Thus, the exterior peaks of $\gamma$ are obtained as follows:

\begin{itemize}
\item Each exterior peak of $\alpha$ is an exterior peak of $\gamma$. (This
includes $n$, if $n$ is an exterior peak of $\alpha$, because $\alpha
_{n}>\beta_{1}$.)

\item Adding $n$ to each exterior peak of $\beta$ yields an exterior peak of
$\gamma$, except for the number $n+1$, which is \textbf{not} an exterior peak
of $\gamma$ (since $\alpha_{n}>\beta_{1}$).
\end{itemize}

These are all the exterior peaks of $\gamma$. Thus,%
\[
\operatorname{Epk}\gamma=\left(  \operatorname{Epk}\alpha\right)
\cup\left(  \left(  \operatorname{Epk}\beta+n\right)  \setminus\left\{
n+1\right\}  \right)  .
\]
In view of $\operatorname{Epk}\alpha=\operatorname{Epk}A$,
$\operatorname{Epk}\beta=\operatorname{Epk}B$ and $\operatorname{Epk}%
\gamma=\operatorname{Epk}\left(  \left[  A,B\right]  \right)  $, this
rewrites as
\[
\operatorname{Epk}\left(  \left[  A,B\right]  \right)  =\left(
\operatorname{Epk}A\right)  \cup\left(  \left(  \operatorname{Epk}%
B+n\right)  \setminus\left\{  n+1\right\}  \right)  .
\]
This proves Lemma \ref{lem.Epk.AB} \textbf{(a)}.

\textbf{(b)} WLOG assume that $\alpha_{i}<\beta_{j}$ for all $i\in\left[
n\right]  $ and $j\in\left[  m\right]  $. (Indeed, we can achieve this by
choosing a positive integer $g$ that is larger than each entry of $\alpha$,
and adding $g$ to each entry of $\beta$.) Thus, in particular, the entries of
$\alpha$ are distinct from the entries of $\beta$. Also, $\alpha_{n}<\beta
_{1}$ (since $\alpha_{i}<\beta_{j}$ for all $i\in\left[  n\right]  $ and
$j\in\left[  m\right]  $).

Let $\gamma$ be the $\left(  n+m\right)  $-permutation $\left(  \alpha
_{1},\alpha_{2},\ldots,\alpha_{n},\beta_{1},\beta_{2},\ldots,\beta_{m}\right)
$. Then, the descents of $\gamma$ are obtained as follows:

\begin{itemize}
\item Each descent of $\alpha$ is a descent of $\gamma$.

\item Adding $n$ to each descent of $\beta$ yields a descent of $\gamma$ (that
is, if $i$ is a descent of $\beta$, then $i+n$ is a descent of $\gamma$).
\end{itemize}

These are all the descents of $\gamma$ (in particular, $n$ is \textbf{not} a
descent of $\gamma$, since $\alpha_{n}<\beta_{1}$). Thus,%
\[
\operatorname{Des}\gamma=\operatorname{Des}\alpha\cup\left(
\operatorname{Des}\beta+n\right)  .
\]
Hence,%
\begin{align*}
\operatorname*{Comp}\left(  \operatorname{Des}\gamma\right)   &
=\operatorname*{Comp}\left(  \operatorname{Des}\alpha\cup\left(
\operatorname{Des}\beta+n\right)  \right) \\
&  =\operatorname*{Comp}\left(  \operatorname{Des}\alpha\right)
\odot\operatorname*{Comp}\left(  \operatorname{Des}\beta\right)
\end{align*}
(because of how $\operatorname*{Comp}S$ is defined for a set $S$). Since
$\operatorname*{Comp}\left(  \operatorname{Des}\pi\right)
=\operatorname*{Comp}\pi$ for any permutation $\pi$, this rewrites as%
\[
\operatorname*{Comp}\gamma=\operatorname*{Comp}\alpha\odot\operatorname*{Comp}%
\beta.
\]
In view of $\operatorname*{Comp}\alpha=A$ and $\operatorname*{Comp}\beta=B$,
this rewrites as $\operatorname*{Comp}\gamma=A\odot B$. Thus,
$\operatorname{Epk}\left(  A\odot B\right)  =\operatorname{Epk}\gamma$.

On the other hand, recall again that $\gamma=\left(  \alpha_{1},\alpha
_{2},\ldots,\alpha_{n},\beta_{1},\beta_{2},\ldots,\beta_{m}\right)  $ and
$\alpha_{n}<\beta_{1}$. Thus, the exterior peaks of $\gamma$ are obtained as follows:

\begin{itemize}
\item Each exterior peak of $\alpha$ is an exterior peak of $\gamma$, except
for the number $n$, which is \textbf{not} an exterior peak of $\gamma$ (since
$\alpha_{n}<\beta_{1}$).

\item Adding $n$ to each exterior peak of $\beta$ yields an exterior peak of
$\gamma$. (This includes $n+1$, if $1$ is an exterior peak of $\beta$, because
$\alpha_{n}<\beta_{1}$.)
\end{itemize}

These are all the exterior peaks of $\gamma$. Thus,%
\[
\operatorname{Epk}\gamma=\left(  \left(  \operatorname{Epk}\alpha\right)
\setminus\left\{  n\right\}  \right)  \cup\left(  \operatorname{Epk}%
\beta+n\right)  .
\]
In view of $\operatorname{Epk}\alpha=\operatorname{Epk}A$,
$\operatorname{Epk}\beta=\operatorname{Epk}B$ and $\operatorname{Epk}%
\gamma=\operatorname{Epk}\left(  A\odot B\right)  $, this rewrites as
\[
\operatorname{Epk}\left(  A\odot B\right)  =\left(  \left(
\operatorname{Epk}A\right)  \setminus\left\{  n\right\}  \right)  \cup\left(
\operatorname{Epk}B+n\right)  .
\]
This proves Lemma \ref{lem.Epk.AB} \textbf{(b)}.
\end{proof}

We can now easily prove Theorem \ref{thm.Epk.dend} again:

\begin{proof}
[Second proof of Theorem \ref{thm.Epk.dend} (sketched).]Let
$A=\operatorname*{QSym}$. Corollary \ref{cor.Epk.ideal} shows that
$\mathcal{K}_{\operatorname{Epk}}$ is an ideal of $\operatorname*{QSym}$.

We now argue the following claim:\footnote{These claims are numbered Claim 2,
Claim 4, Claim 6 and Claim 8, in order to match the numbering of the
corresponding claims in the first proof of Theorem \ref{thm.Epk.dend} above.}

\begin{statement}
\textit{Claim 2:} We have $A\tvi\mathcal{K}_{\operatorname{Epk}}%
\subseteq\mathcal{K}_{\operatorname{Epk}}$.
\end{statement}

[\textit{Proof of Claim 2:} We must show that $\mathcal{K}%
_{\operatorname{Epk}}$ is a left $\tvi$-ideal of $\operatorname*{QSym}$.
According to Proposition \ref{prop.bel-tvi-comp} \textbf{(c)}, this boils down
to proving that if $J$ and $K$ are two $\operatorname{Epk}$-equivalent
nonempty compositions, and if $G$ is any nonempty composition, then $\left[
G,J\right]  $ and $\left[  G,K\right]  $ are $\operatorname{Epk}$-equivalent.

So let $J$ and $K$ be two $\operatorname{Epk}$-equivalent nonempty
compositions. Thus, $\left\vert J\right\vert =\left\vert K\right\vert >0$ and
$\operatorname{Epk}J=\operatorname{Epk}K$.
Let $G$ be any nonempty composition.

Define a positive integer $n$ by $n=\left\vert G\right\vert $. Lemma
\ref{lem.Epk.AB} \textbf{(a)} (applied to $A=G$ and $B=J$) yields
\begin{equation}
\operatorname{Epk}\left(  \left[  G,J\right]  \right)  =\left(
\operatorname{Epk}G\right)  \cup\left(  \left(  \operatorname{Epk}%
J+n\right)  \setminus\left\{  n+1\right\}  \right)  .
\label{pf.thm.Epk.dend.2nd.c2.pf.4}%
\end{equation}
Similarly,%
\begin{equation}
\operatorname{Epk}\left(  \left[  G,K\right]  \right)  =\left(
\operatorname{Epk}G\right)  \cup\left(  \left(  \operatorname{Epk}%
K+n\right)  \setminus\left\{  n+1\right\}  \right)  .
\label{pf.thm.Epk.dend.2nd.c2.pf.5}%
\end{equation}
The right hand sides of (\ref{pf.thm.Epk.dend.2nd.c2.pf.4}) and
(\ref{pf.thm.Epk.dend.2nd.c2.pf.5}) are equal (since $\operatorname{Epk}%
J=\operatorname{Epk}K$). Hence, the left hand sides are equal as well. In
other words, $\operatorname{Epk}\left(  \left[  G,J\right]  \right)
=\operatorname{Epk}\left(  \left[  G,K\right]  \right)  $. Combining this
with%
\[
\left\vert \left[  G,J\right]  \right\vert =\left\vert G\right\vert
+\underbrace{\left\vert J\right\vert }_{=\left\vert K\right\vert }=\left\vert
G\right\vert +\left\vert K\right\vert =\left\vert \left[  G,K\right]
\right\vert ,
\]
we conclude that $\left[  G,J\right]  $ and $\left[  G,K\right]  $ are
$\operatorname{Epk}$-equivalent. As we have said, this concludes the proof of
Claim 2.]

Similarly to Claim 2, we can show the following three claims:

\begin{statement}
\textit{Claim 4:} We have $\mathcal{K}_{\operatorname{Epk}}\tvi
A\subseteq\mathcal{K}_{\operatorname{Epk}}$.
\end{statement}

\begin{statement}
\textit{Claim 6:} We have $A\bel\mathcal{K}_{\operatorname{Epk}}%
\subseteq\mathcal{K}_{\operatorname{Epk}}$.
\end{statement}

\begin{statement}
\textit{Claim 8:} We have $\mathcal{K}_{\operatorname{Epk}}\bel
A\subseteq\mathcal{K}_{\operatorname{Epk}}$.
\end{statement}

(Of course, in proving Claims 4, 6 and 8, we need to use the other three parts
of Proposition \ref{prop.bel-tvi-comp} instead of Proposition
\ref{prop.bel-tvi-comp} \textbf{(c)}, and we occasionally need to use Lemma
\ref{lem.Epk.AB} \textbf{(b)} instead of Lemma \ref{lem.Epk.AB} \textbf{(a)}.)

Combining Claim 2 and Claim 4, we conclude that $\mathcal{K}%
_{\operatorname{Epk}}$ is a $\tvi$-ideal of $A$.

Combining Claim 6 and Claim 8, we conclude that $\mathcal{K}%
_{\operatorname{Epk}}$ is a $\bel$-ideal of $A$.

As in the first proof above, we can show that
$\mathcal{K}_{\operatorname{Epk}} \subseteq A^+$.

Thus, Theorem \ref{thm.ideal-crit2} \textbf{(c)} (applied to $M=\mathcal{K}%
_{\operatorname{Epk}}$) shows that $\mathcal{K}_{\operatorname{Epk}}$ is a
$\left.  \prec\right.  $-ideal and a $\left.  \succeq\right.  $-ideal of
$\operatorname*{QSym}$. This proves Theorem \ref{thm.Epk.dend} again.
\end{proof}

\subsection{\label{subsect.dendri.other-stats}Left/right-shuffle-compatibility
of other statistics}

Let us now briefly analyze the kernels $\mathcal{K}_{\operatorname{st}}$ of
some other descent statistics, following the same approach that we took in our
above second proof of Theorem \ref{thm.Epk.dend} again. Much of what follows
will merely reproduce results from Section \ref{sect.LR}.

\subsubsection{The descent set $\operatorname{Des}$}

First of all, the following is obvious:

\begin{proposition}
\label{prop.Des-set.dend}The ideal $\mathcal{K}_{\operatorname{Des}}$ of
$\operatorname*{QSym}$ is the trivial ideal $0$, and is a $\tvi$-ideal, a
$\bel$-ideal, a $\left.  \prec\right.  $-ideal and a $\left.  \succeq\right.
$-ideal of $\operatorname*{QSym}$.
\end{proposition}

\begin{corollary}
\label{cor.dendri.Des-set}The descent statistic $\operatorname{Des}$ is
left-shuffle-compatible and right-shuffle-compatible.
\end{corollary}

\begin{proof}
[Proof of Corollary \ref{cor.dendri.Des-set}.]Corollary
\ref{cor.dendri.Des-set} can be derived from Proposition
\ref{prop.Des-set.dend} in the same way as Corollary \ref{cor.dendri.Epk} was
derived from Theorem \ref{thm.Epk.dend}.
\end{proof}

\subsubsection{The descent number $\operatorname{des}$}

The permutation statistic $\operatorname{des}$ (called the \textit{descent
number}) is defined as follows: For each permutation $\pi$, we set
$\operatorname{des}\pi=\left\vert \operatorname{Des}\pi\right\vert $ (that
is, $\operatorname{des}\pi$ is the number of all descents of $\pi$). It was
proven in \cite[Theorem 4.6 \textbf{(a)}]{part1} that this statistic
$\operatorname{des}$ is shuffle-compatible. Furthermore, $\operatorname{des}%
$ is clearly a descent statistic. Hence, Proposition \ref{prop.K.ideal}
(applied to $\operatorname{st}=\operatorname{des}$) shows that
$\mathcal{K}_{\operatorname{des}}$ is an ideal of $\operatorname*{QSym}$. We
now claim the following:

\begin{proposition}
\label{prop.des.dend}The ideal $\mathcal{K}_{\operatorname{des}}$ of
$\operatorname*{QSym}$ is a $\tvi$-ideal, a $\bel$-ideal, a $\left.
\prec\right.  $-ideal and a $\left.  \succeq\right.  $-ideal of
$\operatorname*{QSym}$.
\end{proposition}

\begin{corollary}
\label{cor.dendri.des}The descent statistic $\operatorname{des}$ is
left-shuffle-compatible and right-shuffle-compatible.
\end{corollary}

The proofs rely on the following fact (similar to Lemma \ref{lem.Epk.AB}):

\begin{lemma}
\label{lem.des.AB}Let $A$ and $B$ be two nonempty compositions. Let
$n=\left\vert A\right\vert $.

\textbf{(a)} We have $\operatorname{des}\left(  \left[  A,B\right]  \right)
=\operatorname{des}A+\operatorname{des}B+1$.

\textbf{(b)} We have $\operatorname{des}\left(  A\odot B\right)
=\operatorname{des}A+\operatorname{des}B$.
\end{lemma}

\begin{proof}
[Proof of Lemma \ref{lem.des.AB}.]If $I$ is a nonempty composition, then
$\operatorname{des}I$ equals the length of $I$ minus $1$. Lemma
\ref{lem.des.AB} follows easily from this.
\end{proof}

\begin{proof}
[Proof of Proposition \ref{prop.des.dend}.]Analogous to the above second proof
of Theorem \ref{thm.Epk.dend}, but using Lemma \ref{lem.des.AB} instead of
Lemma \ref{lem.Epk.AB}.
\end{proof}

\begin{proof}
[Proof of Corollary \ref{cor.dendri.des}.]Corollary \ref{cor.dendri.des} can
be derived from Proposition \ref{prop.des.dend} in the same way as Corollary
\ref{cor.dendri.Epk} was derived from Theorem \ref{thm.Epk.dend}.
\end{proof}

\subsubsection{The major index $\operatorname*{maj}$}

The permutation statistic $\operatorname*{maj}$ (called the \textit{major
index}) is defined as follows: For each permutation $\pi$, we set
$\operatorname*{maj}\pi=\sum_{i\in\operatorname{Des}\pi}i$ (that is,
$\operatorname*{maj}\pi$ is the sum of all descents of $\pi$). It was proven
in \cite[Theorem 3.1 \textbf{(a)}]{part1} that this statistic
$\operatorname*{maj}$ is shuffle-compatible. Furthermore, $\operatorname*{maj}%
$ is clearly a descent statistic. Hence, Proposition \ref{prop.K.ideal}
(applied to $\operatorname{st}=\operatorname*{maj}$) shows that
$\mathcal{K}_{\operatorname*{maj}}$ is an ideal of $\operatorname*{QSym}$. We
now claim the following:

\begin{proposition}
\label{prop.maj.dend}The ideal $\mathcal{K}_{\operatorname*{maj}}$ of
$\operatorname*{QSym}$ is a right $\tvi$-ideal and a right $\bel$-ideal, but
neither a $\left.  \prec\right.  $-ideal nor a $\left.  \succeq\right.
$-ideal of $\operatorname*{QSym}$.
\end{proposition}

\begin{corollary}
\label{cor.dendri.maj}The descent statistic $\operatorname*{maj}$ is neither
left-shuffle-compatible nor right-shuffle-compatible.
\end{corollary}

The proofs rely on the following fact (similar to Lemma \ref{lem.Epk.AB}):

\begin{lemma}
\label{lem.maj.AB}Let $A$ and $B$ be two nonempty compositions. Let
$n=\left\vert A\right\vert $.

\textbf{(a)} We have $\operatorname*{maj}\left(  \left[  A,B\right]  \right)
=\operatorname*{maj}A+\operatorname*{maj}B+n\cdot\left(  \operatorname{des}%
B+1\right)  $.

\textbf{(b)} We have $\operatorname*{maj}\left(  A\odot B\right)
=\operatorname*{maj}A+\operatorname*{maj}B+n\cdot\operatorname{des}B$.
\end{lemma}

\begin{proof}
[Proof of Lemma \ref{lem.maj.AB}.]If $I=\left(  i_{1},i_{2},\ldots
,i_{k}\right)  $ is a nonempty composition, then
\begin{align*}
\operatorname*{maj}I  &  =i_{1}+\left(  i_{1}+i_{2}\right)  +\left(
i_{1}+i_{2}+i_{3}\right)  +\cdots+\left(  i_{1}+i_{2}+\cdots+i_{k-1}\right) \\
&  =\left(  k-1\right)  i_{1}+\left(  k-2\right)  i_{2}+\cdots+\left(
k-k\right)  i_{k}.
\end{align*}
Lemma \ref{lem.maj.AB} follows easily from this.
\end{proof}

\begin{proof}
[Proof of Proposition \ref{prop.maj.dend}.]To prove that $\mathcal{K}%
_{\operatorname*{maj}}$ is a right $\tvi$-ideal of $\operatorname*{QSym}$, we
proceed as in the proof of Claim 2 in the second proof of Theorem
\ref{thm.Epk.dend}, but using Lemma \ref{lem.maj.AB} instead of Lemma
\ref{lem.Epk.AB}. Similarly, we can show that $\mathcal{K}%
_{\operatorname*{maj}}$ is a right $\bel$-ideal of $\operatorname*{QSym}$.

To prove that $\mathcal{K}_{\operatorname*{maj}}$ is not a $\left.
\prec\right.  $-ideal of $\operatorname*{QSym}$ (and not even a left $\left.
\prec\right.  $-ideal of $\operatorname*{QSym}$), it suffices to find some
$m\in\mathcal{K}_{\operatorname*{maj}}$ and some $a\in\operatorname*{QSym}$
such that $a\left.  \prec\right.  m\notin\mathcal{K}_{\operatorname*{maj}}$.
For example, we can take $m=F_{\left(  1,1,2\right)  }-F_{\left(  3,1\right)
}$ and $a=F_{\left(  1\right)  }$; then, $a\left.  \prec\right.  m=F_{\left(
1,1,1,2\right)  }-F_{\left(  1,3,1\right)  }\notin\mathcal{K}%
_{\operatorname*{maj}}$. The same values of $m$ and $a$ also satisfy $a\left.
\succeq\right.  m\notin\mathcal{K}_{\operatorname*{maj}}$, $m\left.
\prec\right.  a\notin\mathcal{K}_{\operatorname*{maj}}$ and $m\left.
\succeq\right.  a\notin\mathcal{K}_{\operatorname*{maj}}$; thus,
$\mathcal{K}_{\operatorname*{maj}}$ is not a $\left.  \succeq\right.  $-ideal
of $\operatorname*{QSym}$ either. Proposition \ref{prop.maj.dend} is now proven.
\end{proof}

\begin{proof}
[Proof of Corollary \ref{cor.dendri.maj}.]Again, this follows from Proposition
\ref{prop.maj.dend}.
\end{proof}

\subsubsection{The joint statistic $\left(  \operatorname{des}%
,\operatorname*{maj}\right)  $}

The next permutation statistic we shall study is the so-called joint statistic
$\left(  \operatorname{des},\operatorname*{maj}\right)  $. This statistic is
defined as the permutation statistic that sends each permutation $\pi$ to the
ordered pair $\left(  \operatorname{des}\pi,\operatorname*{maj}\pi\right)  $.
(Calling it $\left(  \operatorname{des},\operatorname*{maj}\right)  $ is thus
a slight abuse of notation.) It was proven in \cite[Theorem 4.5 \textbf{(a)}%
]{part1} that this statistic $\left(  \operatorname{des},\operatorname*{maj}%
\right)  $ is shuffle-compatible. Furthermore, $\left(  \operatorname{des}%
,\operatorname*{maj}\right)  $ is clearly a descent statistic. Hence,
Proposition \ref{prop.K.ideal} (applied to $\operatorname{st}=\left(
\operatorname{des},\operatorname*{maj}\right)  $) shows that $\mathcal{K}%
_{\left(  \operatorname{des},\operatorname*{maj}\right)  }$ is an ideal of
$\operatorname*{QSym}$. We now claim the following:

\begin{proposition}
\label{prop.(des,maj).dend}The ideal $\mathcal{K}_{\left(  \operatorname{des}%
,\operatorname*{maj}\right)  }$ of $\operatorname*{QSym}$ is a $\tvi$-ideal, a
$\bel$-ideal, a $\left.  \prec\right.  $-ideal and a $\left.  \succeq\right.
$-ideal of $\operatorname*{QSym}$.
\end{proposition}

\begin{corollary}
\label{cor.dendri.(des,maj)}The descent statistic $\left(  \operatorname{des}%
,\operatorname*{maj}\right)  $ is left-shuffle-compatible and right-shuffle-compatible.
\end{corollary}

\begin{proof}
[Proof of Proposition \ref{prop.(des,maj).dend}.]Analogous to the above second
proof of Theorem \ref{thm.Epk.dend}, but using Lemma \ref{lem.des.AB} together
with Lemma \ref{lem.maj.AB} instead of Lemma \ref{lem.Epk.AB}.
\end{proof}

\begin{proof}
[Proof of Corollary \ref{cor.dendri.(des,maj)}.]Corollary
\ref{cor.dendri.(des,maj)} can be derived from Proposition
\ref{prop.(des,maj).dend} in the same way as Corollary \ref{cor.dendri.Epk}
was derived from Theorem \ref{thm.Epk.dend}.
\end{proof}

\subsubsection{The left peak set $\operatorname*{Lpk}$}

Recall the permutation statistic $\operatorname*{Lpk}$ (the left peak set)
defined in Definition \ref{def.Des-et-al}. It was proven in \cite[Theorem 4.9
\textbf{(a)}]{part1} that this statistic $\operatorname*{Lpk}$ is
shuffle-compatible. Furthermore, $\operatorname*{Lpk}$ is clearly a descent
statistic. Hence, Proposition \ref{prop.K.ideal} (applied to
$\operatorname{st}=\operatorname*{Lpk}$) shows that $\mathcal{K}%
_{\operatorname*{Lpk}}$ is an ideal of $\operatorname*{QSym}$. We now claim
the following:

\begin{proposition}
\label{prop.Lpk.dend}The ideal $\mathcal{K}_{\operatorname*{Lpk}}$ of
$\operatorname*{QSym}$ is a left $\tvi$-ideal, a $\bel$-ideal, a $\left.
\prec\right.  $-ideal and a $\left.  \succeq\right.  $-ideal of
$\operatorname*{QSym}$.
\end{proposition}

\begin{corollary}
\label{cor.dendri.Lpk}The descent statistic $\operatorname*{Lpk}$ is
left-shuffle-compatible and right-shuffle-compatible.
\end{corollary}

The proofs rely on the following fact (similar to Lemma \ref{lem.Epk.AB}):

\begin{lemma}
\label{lem.Lpk.AB}Let $A$ and $B$ be two nonempty compositions. Let
$n=\left\vert A\right\vert $.

\textbf{(a)} We have $\operatorname*{Lpk}\left(  \left[  A,B\right]  \right)
=\left(  \operatorname*{Lpk}A\right)  \cup\left(  \left(  \operatorname*{Lpk}%
B+n\right)  \setminus\left\{  n+1\right\}  \right)  \cup\left\{
n\ \mid\ n-1\notin\operatorname{Des}A\right\}  $.

\textbf{(b)} We have $\operatorname*{Lpk}\left(  A\odot B\right)  =\left(
\operatorname*{Lpk}A\right)  \cup\left(  \operatorname*{Lpk}B+n\right)  $.
\end{lemma}

\begin{proof}
[Proof of Lemma \ref{lem.Lpk.AB}.]Not unlike the proof of Lemma
\ref{lem.Epk.AB} (but left to the reader).
\end{proof}

\begin{proof}
[Proof of Proposition \ref{prop.Lpk.dend}.]Analogous to the above second proof
of Theorem \ref{thm.Epk.dend}, but using Lemma \ref{lem.Lpk.AB} instead of
Lemma \ref{lem.Epk.AB}. This time, however, the analogue of Claim 4 will be
false (i.e., we don't have $\mathcal{K}_{\operatorname*{Lpk}}\tvi
A\subseteq\mathcal{K}_{\operatorname*{Lpk}}$), because the formula for
$\operatorname*{Lpk}\left(  \left[  A,B\right]  \right)  $ in Lemma
\ref{lem.Lpk.AB} \textbf{(a)} depends on $\operatorname{Des}A$. Thus,
$\mathcal{K}_{\operatorname*{Lpk}}$ is merely a left $\tvi$-ideal, not a
$\tvi$-ideal. (But this does not prevent us from applying Theorem
\ref{thm.ideal-crit2} \textbf{(c)}, because that theorem does not require
$M\tvi A\subseteq M$.)
\end{proof}

\begin{proof}
[Proof of Corollary \ref{cor.dendri.Lpk}.]Corollary \ref{cor.dendri.Lpk} can
be derived from Proposition \ref{prop.Lpk.dend} in the same way as Corollary
\ref{cor.dendri.Epk} was derived from Theorem \ref{thm.Epk.dend}.
\end{proof}

\subsubsection{The right peak set $\operatorname*{Rpk}$}

Recall the permutation statistic $\operatorname*{Rpk}$ (the right peak set)
defined in Definition \ref{def.Des-et-al}. It follows from \cite[Theorem 4.9
\textbf{(a)} and Theorem 3.5]{part1} that this statistic $\operatorname*{Rpk}$
is shuffle-compatible (since $\operatorname{Lpk}$ and $\operatorname{Rpk}$ are
$r$-equivalent, using the terminology of \cite{part1}). Furthermore,
$\operatorname*{Rpk}$ is clearly a descent statistic. Hence, Proposition
\ref{prop.K.ideal} (applied to $\operatorname{st}=\operatorname*{Rpk}$) shows
that $\mathcal{K}_{\operatorname*{Rpk}}$ is an ideal of $\operatorname*{QSym}%
$. We now claim the following:

\begin{proposition}
\label{prop.Rpk.dend}The ideal $\mathcal{K}_{\operatorname*{Rpk}}$ of
$\operatorname*{QSym}$ is a $\tvi$-ideal, a right $\bel$-ideal, a left
$\left.  \prec\right.  $-ideal and a left $\left.  \succeq\right.  $-ideal,
but neither a $\left.  \prec\right.  $-ideal nor a $\left.  \succeq\right.
$-ideal of $\operatorname*{QSym}$.
\end{proposition}

\begin{corollary}
\label{cor.dendri.Rpk}The descent statistic $\operatorname*{Rpk}$ is neither
left-shuffle-compatible nor right-shuffle-compatible.
\end{corollary}

The proofs rely on the following fact (similar to Lemma \ref{lem.Epk.AB}):

\begin{lemma}
\label{lem.Rpk.AB}Let $A$ and $B$ be two nonempty compositions. Let
$n=\left\vert A\right\vert $ and $m=\left\vert B\right\vert $.

\textbf{(a)} We have $\operatorname*{Rpk}\left(  \left[  A,B\right]  \right)
=\left(  \operatorname*{Rpk}A\right)  \cup\left(  \operatorname*{Rpk}%
B+n\right)  $.

\textbf{(b)} We have $\operatorname*{Rpk}\left(  A\odot B\right)  =\left(
\left(  \operatorname*{Rpk}A\right)  \setminus\left\{  n\right\}  \right)
\cup\left(  \operatorname*{Rpk}B+n\right)  \cup\left\{  n+1\ \mid
\ 1\in\operatorname{Des}B\text{ or }m=1\right\}  $.
\end{lemma}

\begin{proof}
[Proof of Lemma \ref{lem.Rpk.AB}.]Not unlike the proof of Lemma
\ref{lem.Epk.AB} (but left to the reader).
\end{proof}

\begin{proof}
[Proof of Proposition \ref{prop.Rpk.dend}.]To prove that $\mathcal{K}%
_{\operatorname*{Rpk}}$ is a $\tvi$-ideal and a right $\bel$-ideal, we proceed
as in the above second proof of Theorem \ref{thm.Epk.dend}, but using Lemma
\ref{lem.Rpk.AB} instead of Lemma \ref{lem.Epk.AB}. This time, however, the
analogue of Claim 6 will be false (i.e., we don't have $A\bel\mathcal{K}%
_{\operatorname*{Rpk}}\subseteq\mathcal{K}_{\operatorname*{Rpk}}$), because
the formula for $\operatorname*{Rpk}\left(  A\odot B\right)  $ in Lemma
\ref{lem.Rpk.AB} \textbf{(b)} depends on $\operatorname{Des}B$. Thus,
$\mathcal{K}_{\operatorname*{Rpk}}$ is merely a right $\bel$-ideal, not a
$\bel$-ideal. This prevents us from applying Theorem \ref{thm.ideal-crit2}
\textbf{(c)}. However, we can apply Theorem \ref{thm.ideal-crit2} \textbf{(b)}
instead, and obtain $\operatorname*{QSym}\left.  \succeq\right.
\mathcal{K}_{\operatorname*{Rpk}}\subseteq\mathcal{K}_{\operatorname*{Rpk}}$.
In other words, $\mathcal{K}_{\operatorname*{Rpk}}$ is a left $\left.
\succeq\right.  $-ideal of $\operatorname*{QSym}$. Using
(\ref{eq.dendriform.1}), we thus easily see that $\mathcal{K}%
_{\operatorname*{Rpk}}$ is a left $\left.  \prec\right.  $-ideal of
$\operatorname*{QSym}$ as well.

To prove that $\mathcal{K}_{\operatorname*{Rpk}}$ is not a $\left.
\prec\right.  $-ideal of $\operatorname*{QSym}$ (and not even a right $\left.
\prec\right.  $-ideal of $\operatorname*{QSym}$), it suffices to find some
$m\in\mathcal{K}_{\operatorname*{Rpk}}$ and some $a\in\operatorname*{QSym}$
such that $m\left.  \prec\right.  a\notin\mathcal{K}_{\operatorname*{Rpk}}$.
For example, we can take $m=F_{\left(  1,2\right)  }-F_{\left(  3\right)  }$
and $a=F_{\left(  2\right)  }$; then,
\[
m\left.  \prec\right.  a=-F_{\left(  3,2\right)  }-F_{\left(  2,3\right)
}-F_{\left(  2,2,1\right)  }+F_{\left(  1,2,2\right)  }+F_{\left(
1,1,3\right)  }+F_{\left(  1,1,2,1\right)  }\notin\mathcal{K}%
_{\operatorname*{Rpk}}.
\]
The same values of $m$ and $a$ also satisfy $m\left.  \succeq\right.
a\notin\mathcal{K}_{\operatorname*{Rpk}}$; thus, $\mathcal{K}%
_{\operatorname*{Rpk}}$ is not a $\left.  \succeq\right.  $-ideal of
$\operatorname*{QSym}$ either. Proposition \ref{prop.Rpk.dend} is now proven.
\end{proof}

\begin{proof}
[Proof of Corollary \ref{cor.dendri.Rpk}.]Follows from Proposition
\ref{prop.Rpk.dend}.
\end{proof}

\subsubsection{The peak set $\operatorname*{Pk}$}

Recall the permutation statistic $\operatorname*{Pk}$ (the peak set) defined
in Definition \ref{def.Des-et-al}. It was proven in \cite[Theorem 4.7
\textbf{(a)}]{part1} that this statistic $\operatorname*{Pk}$ is
shuffle-compatible. Furthermore, $\operatorname*{Pk}$ is clearly a descent
statistic. Hence, Proposition \ref{prop.K.ideal} (applied to
$\operatorname{st}=\operatorname*{Pk}$) shows that $\mathcal{K}%
_{\operatorname*{Pk}}$ is an ideal of $\operatorname*{QSym}$. We now claim the following:

\begin{proposition}
\label{prop.Pk.dend}The ideal $\mathcal{K}_{\operatorname*{Pk}}$ of
$\operatorname*{QSym}$ is a left $\tvi$-ideal, a right $\bel$-ideal, a left
$\left.  \prec\right.  $-ideal and a left $\left.  \succeq\right.  $-ideal,
but neither a $\left.  \prec\right.  $-ideal nor a $\left.  \succeq\right.
$-ideal of $\operatorname*{QSym}$.
\end{proposition}

\begin{corollary}
\label{cor.dendri.Pk}The descent statistic $\operatorname*{Pk}$ is neither
left-shuffle-compatible nor right-shuffle-compatible.
\end{corollary}

The proofs rely on the following fact (similar to Lemma \ref{lem.Epk.AB}):

\begin{lemma}
\label{lem.Pk.AB}Let $A$ and $B$ be two nonempty compositions. Let
$n=\left\vert A\right\vert $ and $m=\left\vert B\right\vert $.

\textbf{(a)} We have $\operatorname*{Pk}\left(  \left[  A,B\right]  \right)
=\left(  \operatorname*{Pk}A\right)  \cup\left(  \operatorname*{Pk}B+n\right)
\cup\left\{  n\ \mid\ n-1\notin\operatorname{Des}A\text{ and }n>1\right\}  $.

\textbf{(b)} We have $\operatorname*{Pk}\left(  A\odot B\right)  =\left(
\operatorname*{Pk}A\right)  \cup\left(  \operatorname*{Pk}B+n\right)
\cup\left\{  n+1\ \mid\ 1\in\operatorname{Des}B\right\}  $.
\end{lemma}

\begin{proof}
[Proof of Lemma \ref{lem.Pk.AB}.]Not unlike the proof of Lemma
\ref{lem.Epk.AB} (but left to the reader).
\end{proof}

\begin{proof}
[Proof of Proposition \ref{prop.Pk.dend}.]To prove that $\mathcal{K}%
_{\operatorname*{Pk}}$ is a left $\tvi$-ideal and a right $\bel$-ideal, we
proceed as in the above second proof of Theorem \ref{thm.Epk.dend}, but using
Lemma \ref{lem.Pk.AB} instead of Lemma \ref{lem.Epk.AB}. This time, however,
the analogues of Claim 4 and Claim 6 will be false (i.e., neither
$\mathcal{K}_{\operatorname*{Pk}}\tvi A\subseteq\mathcal{K}%
_{\operatorname*{Pk}}$ nor $A\bel\mathcal{K}_{\operatorname*{Pk}}%
\subseteq\mathcal{K}_{\operatorname*{Pk}}$ will hold), because the formula for
$\operatorname*{Pk}\left(  \left[  A,B\right]  \right)  $ in Lemma
\ref{lem.Pk.AB} \textbf{(a)} depends on $\operatorname{Des}A$ whereas the
formula for $\operatorname*{Pk}\left(  A\odot B\right)  $ in Lemma
\ref{lem.Pk.AB} \textbf{(b)} depends on $\operatorname{Des}B$. Again, this
prevents us from applying Theorem \ref{thm.ideal-crit2} \textbf{(c)}. However,
we can apply Theorem \ref{thm.ideal-crit2} \textbf{(b)} instead, and obtain
$\operatorname*{QSym}\left.  \succeq\right.  \mathcal{K}_{\operatorname*{Pk}%
}\subseteq\mathcal{K}_{\operatorname*{Pk}}$. In other words, $\mathcal{K}%
_{\operatorname*{Pk}}$ is a left $\left.  \succeq\right.  $-ideal of
$\operatorname*{QSym}$. Using (\ref{eq.dendriform.1}), we thus easily see that
$\mathcal{K}_{\operatorname*{Pk}}$ is a left $\left.  \prec\right.  $-ideal of
$\operatorname*{QSym}$ as well.

To prove that $\mathcal{K}_{\operatorname*{Pk}}$ is not a $\left.
\prec\right.  $-ideal of $\operatorname*{QSym}$ (and not even a right $\left.
\prec\right.  $-ideal of $\operatorname*{QSym}$), it suffices to find some
$m\in\mathcal{K}_{\operatorname*{Pk}}$ and some $a\in\operatorname*{QSym}$
such that $m\left.  \prec\right.  a\notin\mathcal{K}_{\operatorname*{Pk}}$.
For example, we can take $m=F_{\left(  1,2\right)  }-F_{\left(  3\right)  }$
and $a=F_{\left(  2\right)  }$; then,
\[
m\left.  \prec\right.  a=-F_{\left(  3,2\right)  }-F_{\left(  2,3\right)
}-F_{\left(  2,2,1\right)  }+F_{\left(  1,2,2\right)  }+F_{\left(
1,1,3\right)  }+F_{\left(  1,1,2,1\right)  }\notin\mathcal{K}%
_{\operatorname*{Pk}}.
\]
The same values of $m$ and $a$ also satisfy $m\left.  \succeq\right.
a\notin\mathcal{K}_{\operatorname*{Pk}}$; thus, $\mathcal{K}%
_{\operatorname*{Pk}}$ is not a $\left.  \succeq\right.  $-ideal of
$\operatorname*{QSym}$ either. Proposition \ref{prop.Pk.dend} is now proven.
\end{proof}

\begin{proof}
[Proof of Corollary \ref{cor.dendri.Pk}.]Follows from Proposition
\ref{prop.Pk.dend}.
\end{proof}
\end{verlong}

\end{document}